%% file: SKM_arxiv.tex
\documentclass[envcountsame, envcountsect,openany]{svmonoalt}

\usepackage{type1cm}         

\usepackage{makeidx}         
\usepackage{graphicx}        
\usepackage{multicol}        
\usepackage[bottom]{footmisc}

\usepackage{newtxtext}       %
\usepackage{newtxmath}       

\makeindex             
                       
\usepackage{mathtools}
\usepackage{amscd}
\usepackage{pdflscape}
\usepackage{graphicx}
\usepackage{bm}
\usepackage{multirow}
\usepackage[normalem]{ulem}
\usepackage{enumitem}
\setlist[enumerate,1]{label={\rm (\arabic*)}} 
\setlist[enumerate,2]{label={\rm (\alph*)}} 
\usepackage{chngcntr}
\counterwithout{figure}{chapter}
\usepackage{booktabs}
\usepackage{dcolumn}
\newcolumntype{d}[1]{D{.}{.}{#1}}
\newcommand\mc[1]{\multicolumn{1}{c}{#1}}
\usepackage{tabu}
\usepackage{hhline}
\usepackage{color}
\usepackage{colortbl}
\usepackage{tikz}
\usetikzlibrary{shapes,arrows,positioning, decorations.pathreplacing, calc, patterns, matrix, fit,cd}
\newcommand{\vcenteredbox}[1]{\begingroup
\setbox0=\hbox{#1}\parbox{\wd0}{\box0}\endgroup}

\newcommand{\bracetwo}{%
  \vcenteredbox{\tikz{%
      \draw[decorate,decoration={brace,amplitude=4pt},xshift=0pt, yshift=0pt] (0,0) -- (0,-2em);
  }}
}
\newcommand{\bracethree}{%
  \vcenteredbox{\tikz{%
      \draw[decorate,decoration={brace,amplitude=4pt},xshift=0pt, yshift=0pt] (0,0) -- (0,-4.5em);
  }}
}

\usepackage[pdfpagelabels, hyperfootnotes=false, colorlinks, linkcolor={blue}, citecolor={magenta}, filecolor={blue}, urlcolor={blue}]{hyperref}

\newcommand{\A}{{\mathbb A}}
\newcommand{\Q}{{\mathbb Q}}
\newcommand{\Z}{{\mathbb Z}}
\newcommand{\R}{{\mathbb R}}
\newcommand{\C}{{\mathbb C}}
\newcommand{\N}{{\mathbb N}}

\newcommand{\idm}{{1}}
\newcommand{\p}{\mathfrak p}
\newcommand{\OF}{{\mathfrak o}}
\newcommand{\GL}{{\rm GL}}

\newcommand{\SL}{{\rm SL}}

\newcommand{\GSp}{{\rm GSp}}

\newcommand{\St}{{\rm St}}
\newcommand{\triv}{1}

\newcommand{\K}[1]{{\rm K}(\p^{#1})}
\newcommand{\Ks}[1]{{\mathrm{K}_s}(\p^{#1})}
\newcommand{\Kl}[1]{{\rm Kl}(\p^{#1})}

\newcommand{\Normalizer}{{\rm N}}
\newcommand{\Mat}{{\rm M}}

\newcommand{\Ind}{{\rm Ind}}

\newcommand{\cInd}{\text{c-Ind}}

\newcommand{\Real}{{\rm Real}}
\newcommand{\Hom}{{\rm Hom}}
\newcommand{\SSp}{{\rm Sp}}

\newcommand{\vl}{{\rm vol}}

\newcommand{\Left}{{\rm Left}}
\newcommand{\Right}{{\rm Right}}
\newcommand{\setchar}{{\rm char}}
\newcommand{\Sym}{{\rm Sym}}
\newcommand{\catone}{{category~1\ }}
\newcommand{\cattwo}{{category~2\ }}

\newcommand*{\transpose}[2][-3mu]{\ensuremath{\mskip1mu\prescript{\smash{\mathrm t\mkern#1}}{}{\mathstrut#2}}}%
\newcommand{\Mod}[1]{\ (\textnormal{mod}\ #1)}
\newcommand{\catonebox}{\fbox{\oldstylenums{1}}}
\newcommand{\cattwobox}{\fbox{\oldstylenums{2}}}

\spnewtheorem*{introtheorem}{Theorem}{\bf}{\it}
\spnewtheorem*{introproposition}{Proposition}{\bf}{\it}
\spnewtheorem*{introcorollary}{Corollary}{\bf}{\it}
\newcommand{\intronumber}[1]{{\bf\hspace{-1mm}{#1}\hspace{1mm}}} 

\allowdisplaybreaks
\begin{document}

\hypersetup{pageanchor=false} 
\author{Jennifer Johnson-Leung\\ Brooks Roberts\\ Ralf Schmidt}
\title{Stable Klingen Vectors and Paramodular Newforms}
\maketitle
\hypersetup{pageanchor=true} 

\frontmatter

\include{SKMS_acknowledgment}

\tableofcontents

\mainmatter

\include{SKMS_chapter1} 
\include{SKMS_part1} 
\include{SKMS_chapter2} 
\include{SKMS_chapter3} 
\include{SKMS_chapter4} 
\include{SKMS_chapter5} 
\include{SKMS_chapter6} 
\include{SKMS_chapter7} 
\include{SKMS_chapter8} 
\include{SKMS_chapter9} 
\include{SKMS_part2} 
\include{SKMS_chapter10} 
\include{SKMS_chapter11} 
\include{SKMS_chapter12} 
\include{SKMS_appendix} 
\include{SKMS_symbols} 

\backmatter

\bibliographystyle{spmpsci}
\bibliography{SKMS_references}
\printindex




\end{document}

%% file: SKMS_acknowledgment.tex
%
%

\extrachap{Acknowledgments}

This work was supported by a Collaboration Grant for Mathematicians from the Simons Foundation (PI Ralf Schmidt).

%% file: SKMS_chapter1.tex
\chapter{Introduction}
\label{introchap}

Siegel modular forms of degree two defined with respect
to the paramodular groups 
$$
\mathrm{K}(N)
=
\SSp(4,\Q) \cap
\begin{bsmallmatrix}
\Z & N\Z & \Z & \Z \vphantom{N^{-1}}\\
\Z & \Z & \Z & N^{-1} \Z \\
\Z & N \Z & \Z & \Z \vphantom{N^{-1}}\\
N\Z & N\Z & \N\Z & \Z \vphantom{N^{-1}}
\end{bsmallmatrix}, \qquad N \in \Z_{>0},
$$
are of significant current interest, due in part to 
the conjecture of Brumer and Kramer \cite{BK1,BK2}. This
conjecture asserts that for every  abelian surface $A$ defined
over $\Q$ with trivial endomorphism ring and conductor $N$
there exists a suitable Siegel modular form $F$ of degree and weight two defined with respect to $\mathrm{K}(N)$ such that $L(s,A)=L(s,F)$,
where $L(s,F)$ is the spin $L$-function associated to $F$. 
This is an attractive conjecture because it is specific 
enough to be checked for examples (see \cite{PY1}, \cite{BDPSM}, \cite{C}, \cite{DK}, \cite{BPY}, \cite{PSY3}, \cite{BK3}, \cite{BPPTGVY},  
\cite{PSY1},  \cite{BKK}, \cite{GPY}, \cite{CCG}).
On the modular forms side, such examples $F$ are often investigated via their Fourier coefficients, and the local factors $L_p(s,F)$ are calculated
by determining the Hecke eigenvalues of $F$ at the prime $p$.

Representation theory can be used to study Siegel modular forms of degree two defined with respect to the paramodular groups.
The work \cite{NF} presented a theory of local new- and oldforms
for irreducible, admissible representations $(\pi,V)$ of $\GSp(4,F)$ with trivial
central character, where $F$ is a 
nonarchimedean local field  of characteristic zero and  we define 
$\GSp(4,F)$ with respect to 
\begin{equation}
\label{tradJeq}
\begin{bsmallmatrix}
&&1&\\
&&&1\\
-1&&&\\
&-1&&
\end{bsmallmatrix}.
\end{equation}
In this theory one considers the subspace of vectors $V(n)$ in $V$
fixed by the local paramodular groups $\K{n}$ for non-negative integers $n$;
here, $\p$ is the prime ideal of the ring of integers $\OF$ of $F$,
$$
\K{n} = \{g \in \GSp(4,F): \lambda(g) \in \OF^\times \}
\cap 
\begin{bsmallmatrix}
\OF & \p^n & \OF & \OF \\
\OF & \OF & \OF & \p^{-n} \\
\OF & \p^n & \OF & \OF\\
\p^n & \p^n  & \p^n & \OF
\end{bsmallmatrix},
$$
and $\lambda(g)$ is the similitude factor of $g$. 
A representation $\pi$ is said to be paramodular if $V(n)$ is non-zero for 
some $n$. If $\pi$ is paramodular and $N_\pi$ is the smallest $n$ such that $V(n) \neq 0$,
then it is proven in \cite{NF} that $V(N_\pi)$ is one-dimensional; the integer
$N_\pi$ is called the paramodular level of $\pi$. 
Moreover, the vectors in $V(n)$ for $n \geq N_\pi$ are obtained from 
$V(N_\pi)$ via three level raising operators.
The work \cite{NF} also introduced two Hecke operators $T_{0,1}$
and $T_{1,0}$ that act on the spaces $V(n)$. If $\lambda_\pi$ and $\mu_\pi$ are the 
eigenvalues of $T_{0,1}$ and $T_{1,0}$, respectively, on $V(N_\pi)$, then 
the spin $L$-factor $L(s,\pi)$ of $\pi$ may be expressed in terms of $\lambda_\pi$ and $\mu_\pi$. 
These results are  ingredients
for the theory of new- and oldforms for Siegel modular forms of degree two
that are defined with respect the paramodular groups $\mathrm{K}(N)$ (see \cite{RSA}, \cite{S2}, \cite{S3}). 

While a satisfactory theory of new- and oldforms exists, open problems remain
concerning basic theory. One such problem concerns Fourier coefficients and the calculation
of Hecke eigenvalues.
The local Hecke operators $T_{0,1}$ and $T_{1,0}$ induce, for each prime $p$ of $\Z$, 
 Hecke operators $T(1,1,p,p)$ and $T(1,p,p,p^2)$ that act
on the vector space $S_k(\mathrm{K}(N))_{\mathrm{new}}$ of cuspidal
newforms of degree two and weight $k$  defined with respect to $\mathrm{K}(N)$.
If $F \in S_k(\mathrm{K}(N))_{\mathrm{new}}$ is an eigenform
for $T(1,1,p,p)$ and $T(1,p,p,p^2)$ for all $p$, then the local factor of the 
spin $L$-function $L_p(s,F)$ at $p$ is expressible in terms
of the eigenvalues $p^{k-3}\lambda_p$ and $p^{2(k-3)}\mu_p$ of these Hecke operators. If $p^2 \nmid N$,
then $T(1,1,p,p)$ and $T(1,p,p,p^2)$ admit upper block expressions, i.e., are given
as sums of the form $\sum_{i=1}^m F|_k g_i$ where the $g_i$ lie in the Siegel
parabolic subgroup of $\GSp(4,\Q)$ (if $p\nmid N$, this is well known; if 
$p\mid\mid N$, see Sect.~2.3 of \cite{S2}). Consequently, when $p^2 \nmid N$, 
there 
are formulas relating the  $\lambda_p$ and $\mu_p$ to 
the Fourier coefficients of $F$ (for an example of this when $N=1$ 
see \cite{MR3541708}).
However, if $p^2 \mid N$, then no such upper block expressions can exist for 
$T(1,1,p,p)$
and $T(1,p,p,p^2)$. Thus, when $p^2 \mid N$, 
it is an open problem to find formulas relating the $\lambda_p$ and $\mu_p$ to 
the Fourier coefficients of $F$; an application of such formulas would provide a way to calculate the local factors of the spin $L$-function from Fourier coefficients.

In this work we present new results about paramodular representations of
$\GSp(4,F)$, and we solve the just mentioned problem using  this local theory. Briefly summarized, our main idea
is to consider certain new local spaces $V_s(n)$ of invariant vectors and the 
corresponding spaces $S_k(\mathrm{K}_s(N))$ of Siegel modular forms of 
degree two; the spaces $V_s(n)$ and $S_k(\mathrm{K}_s(N))$ contain 
$V(n)$ and $S_k(\mathrm{K}(N))$, respectively. We completely determine the 
structure of the local spaces $V_s(n)$, and their relation to the $V(n)$,
for all irreducible, admissible representations of $\GSp(4,F)$ with trivial 
central character. While the spaces $V_s(n)$ are closely related to the $V(n)$,
we find that some problems that are difficult in the setting of the $V(n)$
become easier in the context of the $V_s(n)$. Our main application uses
this observation to solve the above mentioned problem.
As part of this new theory, we introduce novel upper block operators
that effectively serve as a substitute for the paramodular Hecke operators
mentioned above.
From these, we obtain formulas relating the eigenvalues
of a paramodular Hecke eigenform $F$ in $S_k(\mathrm{K}(N))_{\mathrm{new}}$ to the 
Fourier coefficients of $F$ when $p^2 \mid N$. Besides answering a natural question,
our formulas can also be used to efficiently compute Hecke eigenvalues at $p$ from Fourier coefficients when $p^2\mid N$.
Using our formulas, we were able to quickly determine $\lambda_2$ and 
$\mu_2$ for the examples in $S_k(\mathrm{K}(16))_{\mathrm{new}}$ for 
$k = 6,\dots, 14$ presented in \cite{PSY} and  \cite{PSYW}, and the values
of $\lambda_2$ and $\mu_2$ we obtained agree with those from 
\cite{PSY} and  \cite{PSYW}, which were computed using a different method. 
In addition, for a paramodular Hecke eigenform $F$ in $S_k(\mathrm{K}(N))_{\mathrm{new}}$ and $p^2 \mid N$,
we  apply our formulas to prove that the Fourier coefficients
$a(p^t S)$ of $F$ for $t \geq 0$ satisfy a recurrence relation determined
by $L_p(s,F)$ and $a(S)$, $a(pS)$, and $a(p^2 S)$. 
Our results are described in more detail in the remainder of this introduction.
This work is divided into two parts. The first part is concerned with 
local theory, and second part covers applications to Siegel modular forms. 

Finally, to provide a bit more orientation, we mention some other active areas
of investigation concerning Siegel modular forms defined with respect to paramodular groups.
These include  
Eichler-Jacquet-Langlands type
correspondences (\cite{I1}, \cite{I2}, \cite{I3}, \cite{I4}, \cite{IK}, \cite{LW}, \cite{RT}, \cite{vH}, \cite{DPRT}), 
new- and oldforms (\cite{RSA}, \cite{NF}, \cite{S2}, \cite{S3}),
Borcherds products  and lifting
(\cite{G2}, \cite{G3}, \cite{G1}, \cite{D1}, \cite{S4}, \cite{Kr}, \cite{GPY2}, \cite{PSY5}, 
\cite{GV}, \cite{Wi}, \cite{PSY4}, \cite{MiT}, \cite{GPY}, \cite{He}), 
twisting (\cite{JLR1}, \cite{JLR2}, \cite{JLR}), 
congruences (\cite{F1}, \cite{BL}, \cite{DPRT}),
theta and Eisenstein series (\cite{Tak}, \cite{SS}, \cite{Shu}, \cite{BSP1}),
the B\"ocherer conjecture (\cite{MR2825979} and \cite{MR3434888}),
Fourier coefficients and Bessel models (\cite{marzec2016bessel}, \cite{MR3703944}, \cite{MR4244304}), 
and relations to physics (\cite{Na}, \cite{BCGK}, \cite{BCGK2}, \cite{BH}).

\section{Local results}

In this section we summarize our local results. We begin with some definitions.  As above, let 
$F$ be a nonarchimedean local field of characteristic zero with ring
of integers $\OF$. Let $\p$ be the prime ideal of $\OF$, let $\varpi$ be a 
generator of $\p$, and let $q$
be the order of $\OF/\p$. For this introduction we continue to define $\GSp(4,F)$
with respect to \eqref{tradJeq} (though we will use a more convenient 
form in Part~1 of this work). 
Let $n$ be a non-negative integer. 
We define the \emph{Klingen congruence subgroup of level $\p^n$} to be
$$
\mathrm{Kl}(\p^n) = \GSp(4,\OF) \cap 
\begin{bsmallmatrix}
\OF & \p^n & \OF & \OF \\
\OF & \OF & \OF & \OF \\
\OF & \p^n & \OF & \OF\\
\p^n & \p^n & \p^n & \OF
\end{bsmallmatrix},
$$
and we define the \emph{stable Klingen congruence subgroup of level $\p^n$}
to be the subgroup $\mathrm{K}_s(\p^n)$ of $\GSp(4,F)$ generated by $\mathrm{Kl}(\p^n)$ and 
the subgroup 
$$
\begin{bsmallmatrix}
1&&&\\
&1&&\p^{-n+1}\\
&&1&\\
&&&1
\end{bsmallmatrix}.
$$
We have $\mathrm{K}_s(\p^0)=\GSp(4,\OF)$,  $\mathrm{K}_s(\p^1)=\mathrm{Kl}(\p)$, and if $n$ is positive,
then 
$$
\mathrm{K}_s(\p^n) = \{g \in \GSp(4,F)\mid \lambda(g) \in \OF^\times \}
\cap 
\begin{bsmallmatrix}
\OF & \p^n & \OF & \OF \\
\OF & \OF & \OF & \p^{-n+1} \\
\OF & \p^n & \OF & \OF\\
\p^n & \p^n  & \p^n & \OF
\end{bsmallmatrix}.
$$
Evidently, the stable Klingen subgroup $\mathrm{K}_s(\p^n)$ is contained in the 
paramodular subgroup $\mathrm{K}(\p^n)$. Let $(\pi,V)$ be a smooth representation
of $\GSp(4,F)$ for which the center acts trivially. We define 
$$
 V_s(n)=\{v\in V \mid \pi(g)v=v\text{ for all }g\in\Ks{n}\},
$$
and refer to the elements of $V_s(n)$ as \emph{stable Klingen vectors of level $\p^n$}.
If $V_s(n)$ is non-zero for some $n \geq 0$, then we let $N_{\pi,s}$ be the smallest
such $n$ and say that $\pi$ admits non-zero stable Klingen vectors and 
refer to $N_{\pi,s}$ as the \emph{stable Klingen level} of $\pi$. 
As mentioned above, the main idea of this work is to investigate stable
Klingen vectors with a view toward applications to the paramodular theory. 
We were motivated to consider  stable Klingen vectors because of the shadow of a newform. 
The shadow of a newform is contained in $V_s(N_\pi-1)$ when $\pi$ is paramodular and
was introduced, but not deeply studied, in \cite{NF}. The shadow of a newform, 
 defined below, plays a significant role in this work. 

\subsection*{A partition} Our first major local result establishes a fundamental connection between the spaces
$V_s(n)$ and $V(n)$ and a consequent partition 
of the set of paramodular representations.  Since  $V(n-1)+V(n)$ is contained $ V_s(n)$ for all $n$, we see that if $\pi$
is paramodular, then $\pi$ admits non-zero stable Klingen vectors. The 
following theorem implies that the converse also holds when $\pi$ is irreducible. 

\begin{introtheorem}\intronumber{\ref{upperboundtheorem}}
Let $n$ be an integer such that $n \geq 0$. 
Let $(\pi,V)$ be an irreducible, admissible representation of $\GSp(4,F)$ with trivial central character. 
Then
 \begin{equation*}
  \dim V_s(n)\leq\dim V(n)+\dim V(n+1).
 \end{equation*}
\end{introtheorem}

The proof of this theorem uses the $P_3$-quotient of $\pi$ (see Sect.~\ref{repsec}) along with detailed
results about stable Klingen vectors in representations that are induced from the Siegel
parabolic subgroup of $\GSp(4,F)$ (see Chap.~\ref{inducedrepchap}). Theorem \ref{upperboundtheorem}
has the following corollary. 

\theoremstyle{plain}
\begin{introcorollary}\intronumber{\ref{upperboundpropcor}}
Let $(\pi,V)$ be an irreducible, admissible representation of $\GSp(4,F)$ with trivial central character. 
Then $\pi$ admits non-zero paramodular vectors if and only if $\pi$ admits non-zero stable Klingen vectors. If $\pi$ is paramodular, then 
 $N_{\pi,s}=N_\pi-1$ or $N_{\pi,s}=N_\pi$.
\end{introcorollary}

This corollary
divides the set of paramodular representations into two classes. Assume
that $\pi$ is a paramodular, irreducible, admissible representation of $\GSp(4,F)$
with trivial central character. We will say that $\pi$ is a \emph{\catone paramodular representation}
if $N_{\pi,s} = N_\pi -1$; if $N_{\pi,s}=N_\pi$, then we will say that $\pi$
is a \emph{\cattwo paramodular representation}. Interestingly, one can 
prove that this dichotomy admits a characterization independent of the concept
of stable Klingen vectors. We prove the following result.

\begin{introcorollary}\intronumber{\ref{Lparametercharaterizationcor}}
Let $(\pi,V)$ be an irreducible, admissible representation of $\GSp(4,F)$ with trivial central character. 
Assume that $\pi$ is paramodular. 
Then the following are equivalent.
 \begin{enumerate}
  \item  $\pi$ is a \catone paramodular representation, i.e., $N_{\pi,s}=N_\pi-1$.
  \item  The decomposition of the $L$-parameter of $\pi$ into indecomposable representations contains no unramified one-dimensional factors.
 \end{enumerate}
\end{introcorollary}

The partition of the set of paramodular representations
into \catone and \cattwo representations plays an important
role in the results described below.

\subsection*{Structure of the $V_s(n)$} Our second major local result is the determination of the structure of the spaces $V_s(n)$. A cornerstone
for this is the calculation of  the dimensions of the spaces $V_s(n)$
and  certain related vector spaces for all irreducible, admissible representations of $\GSp(4,F)$
with trivial central character and all non-negative integers $n$. Again let
$(\pi,V)$ be a smooth representation of $\GSp(4,F)$ for which the center acts trivially. 
For $n$ a non-negative integer the vector space $V_s(n)$ contains $V(n-1)+V(n)$,
and we define the quotient of $V_s(n)$ by this subspace to be
$$
\bar V_s(n) = V_s(n)/(V(n-1)+V(n)).
$$
The quotient $\bar V_s(n)$ measures the difference between $V_s(n)$
and the subspace of paramodular vectors $V(n-1)+V(n)$; when $V$ is 
infinite-dimensional this subspace
 is actually the direct sum of $V(n-1)$ and $V(n)$, and has structure as 
determined in \cite{NF}. 
If $\bar V_s(n)$ is non-zero for some $n$, then we let $\bar N_{\pi,s}$ be the smallest such $n$
and refer to $\bar N_{\pi,s}$ as the \emph{quotient stable Klingen level} of $\pi$. 
We prove the following theorem about the dimensions of the $V_s(n)$ and $\bar V_s(n)$. 

\begin{introtheorem}\intronumber{\ref{dimensionstheorem}}
For every irreducible, admissible representations $(\pi,V)$ of the group $\GSp(4,F)$ with trivial central character, 
the stable Klingen level $N_{\pi,s}$, 
the dimensions of the spaces $V_s(n)$,
the quotient stable Klingen level $\bar N_{\pi,s}$, and
the dimensions of the spaces $\bar V_s(n)$
are given in Table \ref{dimensionstable}. 
\end{introtheorem}

The proof of this theorem uses the inequality of Theorem \ref{upperboundtheorem}, the already
mentioned results about stable Klingen vectors in representations
induced from the Siegel parabolic subgroup from Chap.~\ref{inducedrepchap},
zeta integrals, and results from \cite{NF}. 

Theorem \ref{dimensionstheorem}
completely determines the partition of paramodular representations into
\catone and \cattwo representations; additionally, 
 Theorem \ref{dimensionstheorem} reveals dimensional growth patterns 
 that are explained by the following result. To state the 
theorem we need to introduce two level raising operators. Let $(\pi,V)$
be a smooth representation for which the center acts trivially, and let 
$n$ be a non-negative integer. We define
$$
\tau, \theta: V_s(n) \longrightarrow V_s(n+1)
$$
by 
$$
\tau v = q^{-1} \sum_{z \in \OF/\p} \pi (
\begin{bsmallmatrix}
1&&&\\
&1&&z\varpi^{-n}\\
&&1&\\
&&&1
\end{bsmallmatrix})v
$$
and 
$$
\theta v = 
\pi(
\begin{bsmallmatrix}
1&&&\\
&1&&\\
&&\varpi&\\
&&&\varpi
\end{bsmallmatrix})v
+
\sum_{x \in \OF /\p}
\pi(
\begin{bsmallmatrix}
1&&x&\\
&1&&\\
&&1&\\
&&&1
\end{bsmallmatrix}
\begin{bsmallmatrix}
\varpi&&&\\
&1&&\\
&&1&\\
&&&\varpi
\end{bsmallmatrix})
$$
for $v \in V_s(n)$. Then $\tau$ and $\theta$ are well-defined, commute,
and induce maps from $\bar V_s(n)$ to $\bar V_s(n+1)$ (see Sect.~\ref{opersksec} and 
Sect.~\ref{levelraisingsec}). The following theorem determines
the structure of the spaces $V_s(n)$ for $n \geq N_{\pi,s}$. 

\begin{introtheorem}\intronumber{\ref{gendecompcor}}
Let $(\pi,V)$ be an infinite-dimensional, irreducible, admissible representation
of $\GSp(4,F)$ with central character. Assume that $\pi$ is paramodular.  Define 
$v_0 \in V$ as follows. If $\bar N_{\pi,s}$ is not defined, set $v_0=0$; if $\bar N_{\pi,s}$
is defined, let $v_0$ be an element of $V_s(\bar N_{\pi,s})$ that is not contained
in $V(\bar N_{\pi,s}-1) + V(\bar N_{\pi,s})$. For integers $n$ such that $n \geq N_{\pi,s}$
let $E(n)$ be the subspace of $V_s(n)$ spanned by the vectors
$$
\tau^i \theta^j v_0, \qquad i,j \geq 0, \quad i+j = n-\bar N_{\pi,s}.
$$
Then 
$$
V_s(n) = V(n-1) \oplus V(n) \oplus E(n)
$$
for integers $n$ such that $n \geq N_{\pi,s}$. 
\end{introtheorem}

The main ingredients for the proof of this theorem are zeta integrals and 
Theorem~\ref{dimensionstheorem}. If $\pi$ is generic, then the integer $\bar N_{\pi,s}$
is defined and $V(n-1)+V(n)$ is a proper subspace of $V_s(n)$ for $n \geq \max(N_\pi,1)$; in
this case there is a natural choice for the vector $v_0$  and the $\tau^i \theta^j v_0$ 
form a basis for $E(n)$ (see Sect.~\ref{structureSKsec} and Fig.~\ref{genericfig} on p.~\pageref{genericfig}
for a visualization). 
If $\pi$ is non-generic, and does not belong to subgroup  IVb or IVc, then 
$V_s(n)=V(n-1)+V(n)$ for $n \geq 0$. Thus, for almost all non-generic $\pi$,
including all Saito-Kurokawa representations, we have $V_s(n) =V(n-1)+V(n)$
for $n \geq 0$ so that $\bar N_{\pi,s}$
is not defined.

As an interesting corollary of Theorem \ref{gendecompcor} we prove that the quotients
$\bar V_s(n)$ exhibit a regular structure reminiscent of $\Gamma_0(\p^n)$-invariant
vectors in irreducible, admissible representations of $\GL(2,F)$ with trivial
central character.

\begin{introcorollary}\intronumber{\ref{genquocor}}
Let $(\pi,V)$ be an irreducible, admissible representation
of $\GSp(4,F)$ with trivial central character. 
Assume that $\pi$ is paramodular.
Assume that $\bar N_{\pi,s}$ is defined. Then 
\begin{equation*}
\bar N_{\pi,s} = 
\begin{cases}
1&\text{if $N_\pi=0$ or $N_\pi=1$,}\\
N_\pi-1&\text{if $N_\pi \geq 2$ and $\pi$ is a \catone representation,}\\
N_\pi&\text{if $N_\pi \geq 2$ and $\pi$ is a \cattwo representation}
\end{cases}
\end{equation*}
and
\begin{equation*}
\dim \bar V_s(\bar N_{\pi,s}) =1.
\end{equation*}
Let $v_{s,\mathrm{new}}$ be a non-zero element of the one-dimensional
vector space $\bar V_s(\bar N_{\pi,s})$. Then the vectors 
\begin{equation*}
\tau^i \theta^j v_{s,\mathrm{new}}, \qquad i,j \geq 0, \quad i+j=n-\bar N_{\pi,s}
\end{equation*}
span $\bar V_s(n)$ for $n \geq \bar N_{\pi,s}$.
\end{introcorollary}

\subsection*{Paramodular Hecke eigenvalues}
Our third major local result relates stable Klingen vectors to paramodular Hecke eigenvalues.
To explain the connection we need some definitions. Let $(\pi,V)$
be a smooth representation of $\GSp(4,F)$ for which the center acts trivially,
and let $n$ be a non-negative integer. In analogy to the paramodular Hecke
operators $T_{0,1}$ and $T_{1,0}$, let
\begin{align*}
\mathrm{K}_s(\p^n) 
\begin{bsmallmatrix}
\varpi&&&\\
&\varpi\vphantom{\varpi^2}&&\\
&&1&\\
&&&1
\end{bsmallmatrix}
\mathrm{K}_s(\p^n) 
&= 
\sqcup_i g_i \mathrm{K}_s(\p^n),\\
\mathrm{K}_s(\p^n) 
\begin{bsmallmatrix}
\varpi&&&\\
&\varpi^2&&\\
&&\varpi&\\
&&&1
\end{bsmallmatrix}
\mathrm{K}_s(\p^n) 
&= 
\sqcup_i h_j \mathrm{K}_s(\p^n)
\end{align*}
be disjoint decompositions, and define 
$$
T_{0,1}^s, T_{1,0}^s: V_s(n) \longrightarrow V_s(n)
$$
by 
$$
T_{0,1}^s v = \sum_i \pi(g_i)v, \qquad T_{1,0}^s v = \sum_j \pi(h_j) v
$$
for $v \in V_s(n)$. 
We refer to $T_{0,1}^s$ and $T_{1,0}^s$ as \emph{stable Klingen Hecke operators}.\index{stable Hecke operators}
\index{Hecke operator!stable}
In contrast to the paramodular case, if $n \geq 1$,
then $T_{0,1}^s$ and $T_{1,0}^s$ admit simple upper block formulas. We prove
that if $n \geq 1$, then 
\begin{align*}
 T^s_{0,1}(v)&=\sum_{y,z\in\OF/\p}\pi(\begin{bsmallmatrix}1&&&\\y&1&&z\varpi^{-n+1}\\&&1&-y\\&&&1\end{bsmallmatrix}
 \begin{bsmallmatrix}
 1&&&\vphantom{\varpi^{-n+1}}\\
 &\varpi&&\vphantom{\varpi^{-n+1}}\\
 &&\varpi&\vphantom{\varpi^{-n+1}}\\
 &&&1\vphantom{\varpi^{-n+1}}
 \end{bsmallmatrix})v\\
 &\quad+\sum_{c,y,z\in\OF/\p}\pi(\begin{bsmallmatrix}1&&c&y\\&1&y&z\varpi^{-n+1}\\&&1\\&&&1\end{bsmallmatrix}
 \begin{bsmallmatrix}
 \varpi&&&\vphantom{\varpi^{-n+1}}\\
 &\varpi&&\vphantom{\varpi^{-n+1}}\\
 &&1&\vphantom{\varpi^{-n+1}}\\
 &&&1\vphantom{\varpi^{-n+1}}\end{bsmallmatrix})v,\\
  T^s_{1,0}(v)&=\sum_{\substack{x,y\in\OF/\p\\z\in\OF/\p^2}} \pi(\begin{bsmallmatrix}1&&&y\\x&1&y&z\varpi^{-n+1}\\&&1&-x\\&&&1\end{bsmallmatrix}
 \begin{bsmallmatrix}
 \varpi&&&\vphantom{\varpi^{-n+1}}\\
 &\varpi^2&&\vphantom{\varpi^{-n+1}}\\
 &&\varpi&\vphantom{\varpi^{-n+1}}\\
 &&&1\vphantom{\varpi^{-n+1}}
 \end{bsmallmatrix})v
\end{align*}
for $v \in V_s(n)$ (see Sect.~\ref{stableheckesec}). If $n \geq 2$, then we also introduce a level lowering operator
$$
\sigma_n: V_s(n+1)  \longrightarrow V_s(n)
$$
that is given by the upper block formula
$$
\sigma_n v = q^{-3} 
\sum_{x,y,z \in \OF/\p}
\pi(
 \begin{bsmallmatrix}1&&&y\\x&1&y&z\varpi^{-n+1}\\&&1&-x\\&&&1\end{bsmallmatrix}
\begin{bsmallmatrix}
1&&&\vphantom{\varpi^{-n+1}}\\
&\varpi&&\vphantom{\varpi^{-n+1}}\\
&&1&\vphantom{\varpi^{-n+1}}\\
&&&\varpi^{-1}\end{bsmallmatrix})v
$$
for $v \in V_s(n+1)$ (see Sect.~\ref{sigmasec}). Finally, assume that $\pi$ is irreducible, admissible, and paramodular.
Let $v_{\mathrm{new}} \in V(N_\pi)$ be a newform, i.e., a non-zero element
of $V(N_\pi)$. When $N_\pi\geq2$, we define the \emph{shadow} of $v_{\mathrm{new}}$ to be the vector 
$$
v_{s}=  \sum_{x,y,z \in \OF/\p}
\pi(
\begin{bsmallmatrix}
1\vphantom{a_b^c}&x\varpi^{N_\pi-1}&&\\
&1&&\\
&y\varpi^{N_\pi-1}&1\\
y\varpi^{N_\pi-1}&z\varpi^{N_\pi-1}&-x\varpi^{N_\pi-1}&1
\end{bsmallmatrix})v_{\rm new}.
$$
The vector $v_s$ is contained in $V_s(N_\pi-1)$. The 
next proposition proves that $\pi$ being a \catone paramodular representation
is equivalent to a number of conditions involving the just introduced concepts.  
These equivalences are used in subsequent parts of this work.
For \ref{shadowsigmavnewpropitem3}
and \ref{shadowsigmavnewpropitem4} note that $v_{\mathrm{new}}$ is contained
in both $V_s(N_\pi)$ and $V_s(N_\pi+1)$. 

\begin{introproposition}\intronumber{\ref{shadowsigmavnewprop}}
Let $(\pi,V)$ be an irreducible, admissible representation of $\GSp(4,F)$ with trivial central character. 
Assume that $\pi$ is paramodular, and that $N_\pi\geq2$. 
Let $v_{\mathrm{new}} \in V(N_\pi)$ be a newform, i.e., a non-zero element of $V(N_\pi)$.
Let $v_s \in V(N_\pi-1)$ be the shadow of $v_{\mathrm{new}}$ as defined in Lemma \ref{Wsgenericlemma1}.
Then the following are equivalent:
\begin{enumerate}
\item  $\pi$ is a \catone paramodular representation, i.e., $N_{\pi,s}=N_\pi-1$.
\item   The $T_{1,0}$-eigenvalue $\mu_\pi$ on $v_{\rm new}$ is non-zero.
\item  $\sigma_{N_\pi-1}v_{\rm new}\neq0$.
\item  $\sigma_{N_\pi}v_{\rm new}\neq0$.
\item  $T^s_{1,0} v_{\rm new}\neq0$.
\item  $v_s\neq0$.
\end{enumerate}
\end{introproposition}

The next theorem ties together stable Klingen operators  and paramodular Hecke eigenvalues 
for irreducible $\pi$. As a  consequence of this theorem the paramodular
Hecke eigenvalues $\lambda_\pi$ and $\mu_\pi$ may be calculated using only upper block
operators. We apply this theorem to Siegel modular forms in the second part of this work.

\begin{introtheorem}\intronumber{\ref{TSKleveltheorem}}
Let $(\pi,V)$ be an irreducible, admissible representation of the group $\GSp(4,F)$ with trivial central character. Assume that $\pi$ is paramodular and that $N_\pi\geq2$. Let $v_{\mathrm{new}} \in V(N_\pi)$ be a newform, and let $v_s$
be the \index{shadow of a newform} shadow of $v_{\mathrm{new}}$. 
\begin{enumerate}
\item  Assume that $\pi$ is a \catone representation, so that $N_{\pi,s}=N_\pi-1$ and $\mu_\pi \neq 0$.  Then $V_s(N_{\pi,s})=V_s(N_\pi-1)$ is one-dimensional and
\begin{equation*}
T^s_{0,1}v_s=\lambda_\pi v_s\quad\text{and}\quad T^s_{1,0}v_s=(\mu_\pi+q^2)v_s.
\end{equation*}
\item  Assume that $\pi$ is a \cattwo representation, so that $N_{\pi,s}=N_\pi$ and $\mu_\pi=0$. Let $v_{\mathrm{new}} \in V(N_\pi)$ be a newform. The vector space $V_s(N_{\pi,s})$ is spanned by the vectors $v_{\mathrm{new}}$ and $T_{0,1}^s v_{\mathrm{new}}$. 
If $v_{\mathrm{new}}$ is not an eigenvector for $T_{0,1}^s$, so that $V_s(N_{\pi,s})$ is two-dimensional, then $\pi$
is generic, and the matrix of $T_{0,1}^s$ in the ordered basis $v_{\mathrm{new}}$, $T_{0,1}^s v_{\mathrm{new}}$ is
\begin{equation*}
\begin{bsmallmatrix}
0&-q^3\\
1&\lambda_\pi
\end{bsmallmatrix},
\end{equation*}
so that 
\begin{equation*}
q^3 v_{\mathrm{new}} +(T_{0,1}^s)^2 v_{\mathrm{new}} = \lambda_\pi T_{0,1}^s v_{\mathrm{new}}.
\end{equation*}
If $v_{\mathrm{new}}$ is an eigenvector for $T_{0,1}^s$, so that $V_s(N_{\pi,s})$ is one-dimensional, then $\pi$ is non-generic, and 
\begin{equation*}
T_{0,1}^s v_{\mathrm{new}} = (1+q^{-1})^{-1} \lambda_{\pi} v_{\mathrm{new}}.
\end{equation*}
\end{enumerate}
\end{introtheorem}

\subsection*{Further local results}
In this work we also prove a number of further local results about
stable Klingen vectors. 
The first collection of results is about the case
when $\pi$ is generic. Let $(\pi,V)$ be a generic, irreducible, admissible
representation of $\GSp(4,F)$ with trivial central character. For 
the first result, assume that $V=\mathcal{W}(\pi,\psi_{c_1,c_2})$
is a Whittaker model of $\pi$ with $c_1, c_2 \in \OF^\times$ (see
Sect.~\ref{repsec}). 
We prove that if $n \geq 0$ and $W \in V_s (n)$, 
then $W\neq 0$ if and only if $W$ does not vanish on the diagonal
subgroup of $\GSp(4,F)$. This generalizes a foundational result from
the paramodular theory. 
Our second result about generic $\pi$ is about the kernel of $\sigma_{n-1}:V_s(n) \to V_s(n-1)$. 
We prove that 
$$
V_s(n) = V(n-1) \oplus V(n) + \ker(\sigma_{n-1})
$$
for integers $n \geq \max(N_{\pi,s},2)$. It follows that $\sigma_{n-1}$ induces
an isomorphism
$$
(V(n-1) \oplus V(n))/\mathcal{K}_{n} \stackrel{\sim}{\longrightarrow} V_s(n-1)
$$
for integers $n \geq \max(N_{\pi,s},2)$, where $\mathcal{K}_n$ is the intersection
of $V(n-1) \oplus V(n)$ with $\ker(\sigma_{n-1})$. 
We determine $\mathcal{K}_n$ for $n \geq \max(N_{\pi,s},2)$ and show that 
this subspace is at most two-dimensional. To prove these statements we
use the  result about the non-vanishing of stable Klingen vectors
on the diagonal of $\GSp(4,F)$. 
Third, assuming that $L(s,\pi)=1$ (this includes all supercuspidal $\pi$), we define a useful graphical model
for $V_s(n)$ for integers $n \geq N_{\pi,s}$. The existence of this 
model also uses  the just mentioned nonvanishing result. In this model, our level changing operators have simple
and visual interpretations. Lastly, still under the hypothesis that $L(s,\pi)=1$,
we prove that $N_\pi \geq 4$. 
For these results see Chap.~\ref{genericchap}.

Finally, for all irreducible, admissible representations $(\pi,V)$ of $\GSp(4,F)$
with trivial central character for which $V_s(1)$ is non-zero, we determine the 
characteristic polynomials of $T_{0,1}^s$ and $T_{1,0}^s$ on $V_s(1)$. If $\pi$
is such a representation, then $\pi$ is Iwahori spherical, and thus may be realized
as a subquotient of a representation induced from an unramified character of the 
Borel subgroup of $\GSp(4,F)$; moreover, the action of the Iwahori-Hecke algebra
of $\GSp(4,F)$ on such induced representations may be explicitly realized. 
We use this observation, via the expression of
$T_{0,1}^s$ and $T_{1,0}^s$ in the Iwahori-Hecke algebra of $\GSp(4,F)$, 
to calculate the desired characteristic polynomials. These results appear
in Chap.~\ref{iwahorichap}. 

\section{Results about Siegel modular forms}

Let $F \in S_k(\mathrm{K}(N))_{\mathrm{new}}$ be an eigenvector of $T(1,1,q,q)$ and $T(1,q,q,q^2)$ for all but finitely many primes $q$ of $\Z$
with $q \nmid N$; by \cite{S3} $F$ is an eigenvector of $T(1,1,q,q)$ and $T(1,q,q,q^2)$ for all primes 
$q$ of $\Z$. Let $T(1,1,q,q)F=q^{k-3}\lambda_q F$ and $T(1,q,q,q^2)F=q^{2(k-3)}\mu_q F$ for each prime $q$ of $\Z$. 
In the second part of this work we apply some of the above local results
to solve the problem of finding formulas relating
$\lambda_p$ and $\mu_p$ to the Fourier coefficients of $F$ for the primes $p$ of $\Z$ such that $p^2 \mid N$.

To make the application, we consider Siegel modular forms defined with respect
to  the global analogues of the local stable Klingen congruence subgroups
from the previous section. We define
$$
\mathrm{K}_s(N)
=
\SSp(4,\Q) \cap
\begin{bsmallmatrix}
\Z & N\Z & \Z & \Z \vphantom{N_s^{-1}}\\
\Z & \Z & \Z & N_s^{-1} \Z \\
\Z & N \Z & \Z & \Z \vphantom{N^{-1}}\\
N\Z & N\Z & \N\Z & \Z \vphantom{N^{-1}}
\end{bsmallmatrix}, \qquad N \in \Z_{>0},
$$
where $N_s = N \prod_{p \mid N} p^{-1}$ and $p$ runs over the primes dividing $N$. 
We refer to $\mathrm{K}_s(N)$ as the \emph{stable Klingen congruence subgroup}\index{stable Klingen congruence subgroup} of level $N$. Evidently, $\mathrm{K}_s(N)$ is contained in $\mathrm{K}(N)$. 
As usual,
we let $M_k(\mathrm{K}_s(N))$  denote the vector space
of Siegel modular forms of degree two and weight $k$ defined with respect to $\mathrm{K}_s(N)$;
the subspace of cusp forms is denoted by $S_k(\mathrm{K}_s(N))$. If $F \in M_k(\mathrm{K}_s(N))$,
then $F$ has a Fourier expansion of the form
$$
F(Z) = \sum_{S \in B(N)} a(S) \E^{2 \pi \I \mathrm{Tr}(SZ)}.
$$
Here, $B(N)$ consists of the half-integral positive semi-definite matrices $S =
\begin{bsmallmatrix} \alpha & \beta \\ \beta & \gamma \end{bsmallmatrix}$ such that $N_s \mid \gamma$;
if $F \in S_k(\mathrm{K}_s(N))$, then the sum runs over the subset $B(N)^+$ of positive definite $S$. 
We note that the space $M_k(\mathrm{K}(N))$ is contained in $M_k(\mathrm{K}_s(N))$. 
If $p$ is a prime dividing $N$, then the local stable Klingen Hecke operators from the previous section
induce operators
$$
T_{0,1}^s(p), T_{1,0}^s (p): M_k(\mathrm{K}_s(N)) \longrightarrow M_k(\mathrm{K}_s(N)),
$$
and if $p^2$ divides $N$, then the local map $\sigma$ induces a level lowering operator
$$
\sigma_p:M_k(\mathrm{K}_s(N)) \longrightarrow M_k(\mathrm{K}_s(Np^{-1})).
$$
These maps take cusp forms to cusp forms. 
We prove that these operators are given by slash formulas involving only upper block matrices.
As a consequence, we show that
if $F$ has the above Fourier expansion, then the 
Fourier expansions of $T_{0,1}^s(p)F$, $T_{1,0}^s (p)F$, and $\sigma_pF$ are given by
\begin{align*}
(T_{0,1}^{s}(p)F)(Z)
&=\sum_{S=\begin{bsmallmatrix}\alpha&\beta\\ \beta&\gamma \end{bsmallmatrix}\in B(N)} 
\Big( p^{3-k}  a(pS)\\
&\quad+  \sum_{\substack{y \in \Z/p\Z \\p\mid(\alpha+2\beta y+\gamma y^2)}} 
p  a(p^{-1}S[\begin{bsmallmatrix}1&\\y&p\end{bsmallmatrix}]) \Big) \E^{2\pi \I \mathrm{Tr}(SZ)},\\
(T_{1,0}^{s}(p)F)(Z)&= \sum_{\substack{S\in B(N)}} \sum_{\substack{a\in\Z/p\Z}} p^{3-k}
 a(S[\begin{bsmallmatrix}1&\\a&p\end{bsmallmatrix}]) \E^{2\pi \I \mathrm{Tr} (SZ)},\\
(\sigma_{p}F)(Z)&= \sum_{\substack{S\in B(Np^{-1})}}\sum_{\substack{a\in\Z/p\Z}} p^{-k-1}
 a(S[\begin{bsmallmatrix}1&\\a&p\end{bsmallmatrix}]) \E^{2\pi \I \mathrm{Tr}{(SZ)}}.
\end{align*}
In these formulas if $S$ and $A$ are $2 \times 2$ matrices, then $S[A]=\transpose{A}SA$. 
We also calculate the Fourier-Jacobi expansions of $T_{0,1}^s(p)F$, $T_{1,0}^s (p)F$, and $\sigma_pF$.
These formulas involve  new operators on Jacobi forms that may be of independent interest.
Similarly, we derive formulas for the operators induced by the other local maps discussed in the previous section; see Chap.~\ref{operchap}. 

\subsection*{The main theorem} Let $F \in S_k(\mathrm{K}(N))_{\mathrm{new}}$ be an eigenvector of $T(1,1,q,q)$ and $T(1,q,q,q^2)$ for all but finitely primes $q$ of $\Z$
with $q \nmid N$; as mentioned above, by \cite{S3} $F$ is an eigenvector of $T(1,1,q,q)$ and $T(1,q,q,q^2)$ for all primes 
$q$ of $\Z$. Let $T(1,1,q,q)F=q^{k-3}\lambda_q F$ and $T(1,q,q,q^2)F=q^{2(k-3)}\mu_q F$ for each prime $q$ of $\Z$. 
Using Theorem \ref{TSKleveltheorem} and other local results we obtain the following result, which relates 
$\lambda_p$ and $\mu_p$  to the action of the upper block operators $T^s_{0,1}(p)$, $T^s_{1,0}(p)$, 
and $\sigma_p$ on $F$ for primes $p$ such that $p^2 \mid N$. This theorem mentions the representation $\pi$ 
generated by the adelization of $F$; see Theorem \ref{Ralftheorem} (based on \cite{S2} and \cite{S3}) for a summary of the properties of $\pi$.

\begin{introtheorem}\intronumber{\ref{globalalgthm}}
Let $N$ and $k$ be integers such that $N>0$ and $k>0$, and let $F \in S_k(\mathrm{K}(N))$. 
Assume that $F$ is a newform and is an eigenvector for the Hecke operators $T(1,1,q,q)$ and $T(1,q,q,q^2)$
for all but finitely many of the primes $q$ of $\Z$ such that $q \nmid N$; then by Theorem
\ref{Ralftheorem}, $F$ is an eigenvector for $T(1,1,q,q)$ and $T(1,q,q,q^2)$ for all primes $q$ of $\Z$.
Let $\otimes_{v \leq \infty} \pi_v$ be as in Theorem \ref{Ralftheorem}. 
For every prime $q$ of $\Z$ let $\lambda_q, \mu_q \in \C$ be such that 
\begin{align*}
\tag{\ref{globallambdaeq}} T(1,1,q,q)F&=q^{k-3}\lambda_qF,\\
\tag{\ref{globalmueq}} T(1,q,q,q^2)F&=q^{2(k-3)}\mu_qF.
\end{align*}
Let $p$ be a prime of $\Z$ with $v_p(N)\geq2$. Then 
\begin{equation*}
\tag{\ref{globalalgthmeq01}}
\mu_p = 0 \Longleftrightarrow \sigma_p F = 0 \Longleftrightarrow T_{1,0}^s(p) F = 0.
\end{equation*}
Moreover:
\begin{enumerate}
\item  If $v_p(N) \geq 3$, then $\sigma_p^2 F =0$.
\item  We have
\begin{equation*}
\tag{\ref{globalnewformreleq}}
\mu_p F = p^4 \tau_p^2 \sigma_p F - p^2 \eta_p \sigma_p F.
\end{equation*}
\item  Assume that $\mu_p\neq 0$. Then  $\sigma_p F \neq 0$, and 
\begin{align*}
T^{s}_{0,1}(p)(\sigma_p F)&=\lambda_p(\sigma_pF),\tag{\ref{stablelambdacase1eq}}\\
T^{s}_{1,0}(p)(\sigma_p F)&=(\mu_p+p^2)(\sigma_pF)\tag{\ref{stablemucase1eq}},
\end{align*}
and the representation $\pi_p$ is generic.
\item  Assume that $\mu_p=0$. Then
\begin{equation*}
\tag{\ref{sigmazeromuzeroeq}}
T_{1,0}^s(p) F =0.
\end{equation*}
If $F$ is not an eigenvector for $T^{s}_{0,1}(p)$, then $\pi_p$ is generic and 
\begin{align*}
T^{s}_{0,1}(p)^2F=-p^3F+\lambda_p T^{s}_{0,1}(p)F. \tag{\ref{stablelambdacase3eq}}
\end{align*}
The newform $F$ is an eigenvector for $T^{s}_{0,1}(p)$ if and only if 
\begin{align*}
T^{s}_{0,1}(p)F=(1+p^{-1})^{-1}\lambda_pF;\tag{\ref{stablelambdacase2eq}}
\end{align*}
in this case $\pi_p$ is non-generic. 
\end{enumerate}
\end{introtheorem}

\subsection*{Fourier coefficients} The next corollary translates the assertions of Theorem~\ref{globalalgthm} into identities
involving the Fourier coefficients of $F$. This is a natural application of this theorem since the operators
$T_{0,1}^s(p)$, $T_{1,0}^s(p)$, and $\sigma_p$ are upper block operators. These relations in turn allow
for efficient numerical calculations, at least in the examples that we discuss subsequent to the statement of
the corollary. 

\begin{introcorollary}\intronumber{\ref{fourierheckethm}}
Let $N$ and $k$ be integers such that $N>0$ and $k>0$, and let $F \in S_k(\mathrm{K}(N))$. 
Assume that $F$ is a newform and is an eigenvector for the Hecke operators $T(1,1,q,q)$ and $T(1,q,q,q^2)$
for all but finitely many of the primes $q$ of $\Z$ such that $q \nmid N$;  by Theorem
\ref{Ralftheorem}, $F$ is an eigenvector for $T(1,1,q,q)$, $T(1,q,q,q^2)$ and $U(q)$ for all primes $q$ of $\Z$.
Let $\lambda_q, \mu_q \in \C$ and $\varepsilon_q \in \{\pm 1\}$ be such that 
\begin{align*}
T(1,1,q,q)F&=q^{k-3}\lambda_qF,\\
T(1,q,q,q^2)F&=q^{2(k-3)}\mu_qF,\\
U(q)F&=\varepsilon_q F
\end{align*}
for all primes $q$ of $\Z$. Regard $F$ as an element of $S_k(\mathrm{K}_s(N))$, and let 
\begin{align*}
F(Z)=\sum_{S\in B(N)^+}a(S) \E^{2\pi \I \mathrm{Tr}(SZ)}
\end{align*}
be the Fourier expansion of $F$. Let $\pi\cong \otimes_{v \leq \infty} \pi_v$ be as in Theorem \ref{globalalgthm}. 
Let $p$ be a prime of $\Z$ with $v_p(N)\geq2$. Then
\begin{gather*}
\sum_{a \in \Z/p\Z} a(S[ \begin{bsmallmatrix} 1& \\ a&p \end{bsmallmatrix} ] ) =0 \qquad
\text{for $S \in B(Np^{-1})^+$}  \\
\Updownarrow\tag{\ref{fourierheckethmeq1000}}\\
\mu_p = 0 \\
\Updownarrow\tag{\ref{fourierheckethmeq1001}}\\
\sum_{a \in \Z/p\Z} a(S[ \begin{bsmallmatrix} 1& \\ a&p \end{bsmallmatrix} ] ) =0 \qquad 
\text{for $S \in B(N)^+$}. 
\end{gather*}
Moreover:
\begin{enumerate}
\item 
If $v_p(N) \geq 3$ and $S\in B(Np^{-2})^+$, then 
\begin{equation*}
\tag{\ref{fourierheckethmsigmassquareitemeq}}
\sum_{z \in \Z/ p^2 \Z}  a(S[\begin{bsmallmatrix} 1& \\ z&p^2 \end{bsmallmatrix}]) = 0.
\end{equation*}
\item 
If $S= \begin{bsmallmatrix} \alpha&\beta\\ \beta&\gamma \end{bsmallmatrix} \in B(Np)^+$, then 
\begin{equation*}
\tag{\ref{mufouriereq}}
\mu_p  a(S)=
\begin{cases}
\sum\limits_{x\in\Z/p\Z} p^{3-k}  a(S[\begin{bsmallmatrix}1&\\x&p\end{bsmallmatrix}])&\text{if $p\nmid 2\beta$,}\\
\sum\limits_{x\in\Z/p\Z} p^{3-k}  a(S[\begin{bsmallmatrix}1&\\x&p\end{bsmallmatrix}])
-\sum\limits_{x\in\Z/p\Z} p  a(S[\begin{bsmallmatrix}1&\\xp^{-1}&1\end{bsmallmatrix}])
&\text{if $p\mid 2\beta$}.
\end{cases}
\end{equation*}
\item 
Assume that $\mu_p \neq 0$. 
If $S= \begin{bsmallmatrix} \alpha&\beta\\ \beta&\gamma \end{bsmallmatrix}\in B(Np^{-1})^+$, then 
\begin{align*}
\lambda_p \sum_{x\in\Z/p\Z} a(S[\begin{bsmallmatrix}1&\\x&p\end{bsmallmatrix}])
&=\sum_{x\in\Z/p\Z}p^{3-k}   a(pS[\begin{bsmallmatrix}1&\\x&p\end{bsmallmatrix}])\nonumber \\
&\quad+\sum_{\substack{z\in\Z/p^2\Z\\p\mid(\alpha+2\beta z+\gamma z^2)}}p  a(p^{-1}S[\begin{bsmallmatrix}1&\\z&p^2\end{bsmallmatrix}])
\tag{\ref{munozerolambdaeq}}
\end{align*}
and 
\begin{equation*}
\sum_{y\in\Z/p^2\Z}  a(S[\begin{bsmallmatrix}1&\\ y&p^2\end{bsmallmatrix}])
=\begin{cases}\varepsilon_p \sum\limits_{x\in\Z/p\Z}p^{k-2}  a(S[\begin{bsmallmatrix}1&\\x&p\end{bsmallmatrix}])&\text{if } v_p(N)=2,\\
0&\text{if } v_p(N)>2.
\end{cases}\tag{\ref{epsilonfourierformulaeq}}
\end{equation*}
The representation $\pi_p$ is generic.
\item Assume that $\mu_p=0$. 
Then 
\begin{equation*}
\sum_{a \in \Z/p\Z} a(S[ \begin{bsmallmatrix} 1& \\ a&p \end{bsmallmatrix} ] ) =0\qquad \text{for $S \in B(Np^{-1})^+$.}
\tag{\ref{fourierheckethmeq1002}}
\end{equation*}
We have $T_{0,1}^s(p)F = \sum_{S \in B(N)^+} c(S) \E^{2\pi \I \mathrm{Tr}(SZ)}$ where
\begin{equation*}
\tag{\ref{fourierheckethmeq10021}}
c(S)= p^{3-k}  a(pS)+ \sum_{\substack{x\in\Z/p\Z\\ p\mid (\alpha+2\beta x)}}
p  a(p^{-1}S[\begin{bsmallmatrix}1&\\ x&p\end{bsmallmatrix}])
\end{equation*}
for $S=\begin{bsmallmatrix} \alpha&\beta\\ \beta&\gamma \end{bsmallmatrix} \in B(N)^+$. 
If $F$ is not an eigenvector for  $T^{s}_{0,1}(p)$, then $\pi_p$ is generic, 
\begin{align*}
\lambda_p  c(S)
&=p^3 a(S)+p^{6-2k}  a(p^2S)+\sum_{\substack{y\in\Z/p\Z\\ p\mid(\alpha+2\beta y)}}
p^{4-k}  a(S[\begin{bsmallmatrix}1&\\y&p\end{bsmallmatrix}]) \\
&\quad+\sum_{\substack{z\in\Z/p^2\Z\\ p^2 \mid (\alpha+ 2\beta z+\gamma z^2)}}
p^{2}  a(p^{-2}S[\begin{bsmallmatrix}1&\\z&p^2\end{bsmallmatrix}])\tag{\ref{muzerosecondcaseeq}}
\end{align*}
for $S= \begin{bsmallmatrix} \alpha&\beta\\ \beta&\gamma \end{bsmallmatrix} \in B(N)^+$,
and $c(S) \neq 0$ for some $S \in B(N)^+$. 
The newform $F$ is an eigenvector for $T^{s}_{0,1}(p)$ if and only if
\begin{equation*} 
\tag{\ref{muzerofirstcaseeq}}
\lambda_p  a(S) = (1+p^{-1}) c(S) 
=(1+p) p^{2-k}  a(pS)+\sum_{\substack{x\in\Z/p\Z\\ p \mid (\alpha +  2\beta x)}}
(1+p)  a(p^{-1}S[\begin{bsmallmatrix}1&\\x&p\end{bsmallmatrix}])
\end{equation*}
for $S= \begin{bsmallmatrix} \alpha&\beta\\ \beta&\gamma \end{bsmallmatrix} \in B(N)^+$; in this case $\pi_p$ is non-generic.
\end{enumerate}
\end{introcorollary}

The formulas of Corollary  \ref{fourierheckethm} provide a solution to the problem
of finding formulas that relate $\lambda_p$ and $\mu_p$ to the Fourier coefficients of $F$
when $F$ is as in Corollary  \ref{fourierheckethm} and $p^2 \mid N$. 
As mentioned earlier, when $p^2 \mid N$, the Hecke operators $T(1,1,p,p)$ and $T(1,p,p,p^2)$
do not admit formulas involving only upper block matrices; consequently, a direct calculation is not 
possible. Instead, our solution uses the new theory of stable Klingen vectors to relate
these Hecke eigenvalues to the new upper block operators $T_{0,1}^s(p), T_{1,0}^s(p)$, and $\sigma_p$. 

The reader may recognize the sums occurring in this corollary as involving $\GL(2)$ Hecke operators. 
In Sect.~\ref{altformsec} we use this observation 
to reformulate the results of Corollary \ref{fourierheckethm} into  briefer conceptual statements (see 
Theorem \ref{abstractformulathm}).

\subsection*{Calculations} The formulas of Corollary  \ref{fourierheckethm} seem to be of practical value. Using a computer algebra program, we verified that the identities
of Corollary~\ref{fourierheckethm} do indeed hold for the case 
$N=16$, $p=2$,  $k = 6,\dots, 14$ for the examples presented in \cite{PSY} and  \cite{PSYW} (see Sect.~\ref{examplessec} for 
a discussion of this verification). Besides providing a welcome check, our calculations also showed that the formulas of Corollary  \ref{fourierheckethm}
can be used to quickly determine $\lambda_2$ and $\mu_2$ from the Fourier coefficients of the examples
of \cite{PSY} and  \cite{PSYW}. See Sect.~\ref{eigenvaluescompsec} for a description of how to use Corollary  \ref{fourierheckethm}
to compute Hecke eigenvalues, along with some illustrative examples.

\subsection*{A recurrence relation}
Finally, as an application of Corollary \ref{fourierheckethm}
we prove that the Fourier coefficients of $F$ as in Corollary \ref{fourierheckethm}
satisfy a recurrence relation. The following theorem shows that the Fourier
coefficients $a(p^t S)$ of $F$ are determined by $a(S)$, $a(pS)$, $a(p^2S)$
and the Hecke eigenvalues $\lambda_p$ and $\mu_p$ for $p$ such that $v_p(N) \geq 2$. 

\begin{introtheorem}\intronumber{\ref{recurrencetheorem}}
 Let $N$ and $k$ be integers such that $N>0$ and $k>0$, and let $F \in S_k(\mathrm{K}(N))$. 
Assume that $F$ is a newform and is an eigenvector for the Hecke operators $T(1,1,q,q)$ and $T(1,q,q,q^2)$
for all but finitely many of the primes $q$ of $\Z$ such that $q \nmid N$;  by Theorem
\ref{Ralftheorem}, $F$ is an eigenvector for $T(1,1,q,q)$, $T(1,q,q,q^2)$ and $U(q)$ for all primes $q$ of $\Z$.
Let $\lambda_q, \mu_q \in \C$  be such that 
\begin{align*}
T(1,1,q,q)F&=q^{k-3}\lambda_qF,\\
T(1,q,q,q^2)F&=q^{2(k-3)}\mu_qF
\end{align*}
for all primes $q$ of $\Z$. Regard $F$ as an element of $S_k(\mathrm{K}_s(N))$, and let 
\begin{align*}
F(Z)=\sum_{S\in B(N)^+}a(S) \E^{2\pi \I \mathrm{Tr}(SZ)}
\end{align*}
be the Fourier expansion of $F$. 
Let $p$ be a prime of $\Z$ with $v_p(N)\geq2$. If $S \in B(N)^+$, then there is formal identity
of power series in $p^{-s}$
\begin{equation*}
\tag{\ref{recurrencetheoremeq1}}
\sum\limits_{t=0}^\infty \dfrac{a(p^t S)}{p^{ts}} = N(p^{-s},S) L_p(s,F)
\end{equation*}
where 
\begin{align*}
N(p^{-s},S) &= a(S) + \big( a(pS) - p^{k-3} \lambda_p a(S) \big) p^{-s}\nonumber \\
\tag{\ref{recurrencetheoremeq2}} &\quad + \big( a(p^2 S) - p^{k-3} \lambda_p a(pS) + p^{2k-5} (\mu_p +p^2) a(S) \big) p^{-2s}
\end{align*}
and 
\begin{equation*}
\tag{\ref{recurrencetheoremeq3}}
L_p(s,F) = \dfrac{1}{1 - p^{k-3} \lambda_p p^{-s} +p^{2k-5} (\mu_p + p^2 ) p^{-2s} }
\end{equation*}
is the spin $L$-factor of $F$ at $p$ (e.g., see \cite{MR2887605}, p.~547).
\end{introtheorem}

A similar statement when $\mathrm{K}(N)$ is replaced by $\Gamma_0(N)$ and $p \nmid N$
can be found in Sect.~4.3.2 of \cite{MR884891}. We also note that the examples of \cite{PSY} and  \cite{PSYW}
can be used to exhibit $F$ and $S$ such that the coefficients of $p^{-ts}$ for $t=0,1$ and $2$ in $N(p^{-s},S)$ are non-zero.

%% file: SKMS_part1.tex
\begin{partbacktext}
\part{Local Theory}
\end{partbacktext}

%% file: SKMS_chapter2.tex
\chapter{Background}
\label{backgroundchap}

In this chapter we recall some
essential definitions and results concerning
the group $\GSp(4)$ over a nonarchimedean
local field of characteristic zero
and its representation theory. We also
review the theory of paramodular vectors.
\section{Some definitions}
\label{somedefsec}
The following objects and 
definitions will be fixed
for Part~1 of this work.
\subsection*{The base field and matrices}
\label{basefieldsubsecstar}
We let $F$ be a nonarchimedean local field of
characteristic zero. Let 
$\OF$ denote the ring of integers of $F$, and let $\p$
be the maximal proper ideal of $\OF$. Once and
for all, we fix a generator $\varpi$ of $\p$,
so that $\p = \OF \varpi$. 
If $x \in F^\times$, then there exists a unique
integer $v(x)$ and $u \in \OF^\times$ such
that $x = u \varpi^{v(x)}$. The order of $\OF/\p$
is denoted by $q$, and we denote by $|\cdot|$
the absolute value on $F$ such that $|\varpi|=q^{-1}$.
We will use the Haar measure $dx$ on $F$ that gives $\OF$
volume $1$. Occasionally, we will use the Haar measure
on $F^\times$ that is given by $d^\times x = dx/|\cdot|$. 
We let $\psi: F \to \C^\times$ denote a fixed continuous
homomorphism such that $\psi(\OF) =1$ but $\psi(\p^{-1}) \neq 1$. 
If $g$ is a matrix,
then we denote the \emph{transpose} of $g$ by $\transpose{g}$.  \label{transposeg}
\index{transpose of a matrix}
In this work, if a matrix is specified, then a blank entry
means a zero entry. 
If $A=\begin{bsmallmatrix}a_1&a_2\\a_3&a_4\end{bsmallmatrix}$ 
is in $\GL(2,F)$, then we define
\begin{equation}
\label{Aprimedefeq}
A'=
\begin{bsmallmatrix}
&1\\
1&
\end{bsmallmatrix}
\transpose{A}^{-1}
\begin{bsmallmatrix}
&1\\
1&
\end{bsmallmatrix} 
=
\det(A)^{-1}
\begin{bsmallmatrix}
a_1&-a_2\\
-a_3&a_4
\end{bsmallmatrix}.
\end{equation}
It is useful to use the asterisk symbol $*$ as
a wild card to define a subset of a given set of matrices. For example,
when working with $\GL(2,F)$, the phrase ``let $B=\begin{bsmallmatrix} *&* \\ & * \end{bsmallmatrix}$'' 
would define  the subgroup of all upper triangular elements of $\GL(2,F)$.
It is useful to note that if $x \in F^\times$, then 
\begin{equation}
\label{usefuleq}
\begin{bsmallmatrix}
\vphantom{-}1^{\vphantom{-1}}&\\
x&\vphantom{-}1 
\end{bsmallmatrix}
=
\begin{bsmallmatrix}
\vphantom{-}1^{\vphantom{-1}}&x^{-1}\\
&\vphantom{-}1
\end{bsmallmatrix}
\begin{bsmallmatrix}
-x^{-1}&\\
&-x
\end{bsmallmatrix}
\begin{bsmallmatrix}
&\vphantom{-}1^{\vphantom{-1}}\\
-1&
\end{bsmallmatrix}
\begin{bsmallmatrix}
\vphantom{-}1^{\vphantom{-1}}&x^{-1}\\
&\vphantom{-}1
\end{bsmallmatrix}.
\end{equation}

\subsection*{The symplectic similitude group}
\index{GSp@$\GSp(4,F)$}
\label{GSp4Fdef}
In this first part we define $\GSp(4,F)$ to be the set of 
$g$ in $\GL(4,F)$ such that $\transpose{g} J g = \lambda J$ for some 
$\lambda \in F^\times$, where
\begin{equation}
\label{Jdefeq}
J=
\begin{bsmallmatrix}
&&&1\\
&&1&\\
&-1&&\\
-1&&&
\end{bsmallmatrix}.
\end{equation}
In the second part of this work, which concerns Siegel modular forms,
we will use the traditional form of $J$ to define $\GSp(4)$.
If $g \in \GSp(4,F)$, then the unit $\lambda$
such that $\transpose{g} J g = \lambda J$ is unique, and will be denoted by $\lambda(g)$. \label{simdef}
It is easy to see that $\GSp(4,F)$ is a subgroup of $\GL(4,F)$. 
If  $g=\begin{bsmallmatrix} A&B \\ C&D \end{bsmallmatrix}
\in \GSp(4,F)$ with $A=\begin{bsmallmatrix}a_1&a_2\\a_3&a_4\end{bsmallmatrix},
B=\begin{bsmallmatrix}b_1&b_2\\b_3&b_4\end{bsmallmatrix},
C=\begin{bsmallmatrix}c_1&c_2\\c_3&c_4\end{bsmallmatrix}$, and 
$D=\begin{bsmallmatrix}d_1&d_2\\d_3&d_4\end{bsmallmatrix}$ in $\Mat(2,F)$, then
\begin{equation}
\label{ginveq}
g^{-1} 
=\lambda(g)^{-1}
\begin{bsmallmatrix}
\hfill d_4&\hfill d_2&\hfill -b_4&\hfill -b_2\\
\hfill d_3&\hfill d_1&\hfill  -b_3&\hfill -b_1\\
\hfill -c_4&\hfill -c_2&\hfill  a_4&\hfill a_2\\
\hfill -c_3&\hfill -c_1&\hfill a_3&\hfill a_1
\end{bsmallmatrix}.
\end{equation}
We define $\SSp(4,F)$ to be the subgroup of $g$ in $\GSp(4,F)$ such that \label{Sp4Fdef}
$\lambda(g)=1$, i.e., $\transpose{g} J g = J$. 
The center $Z$ of the group $\GSp(4,F)$ consists of all the matrices of \label{Zcenterdef}
the form
$$
\begin{bsmallmatrix}
z&&&\\
&z&&\\
&&z&\\
&&&z
\end{bsmallmatrix}, \qquad z \in F^\times.
$$
We give $\GSp(4,F)\subset \GL(4,F)$ the relative topology.
The groups $\GSp(4,F)$ and $\SSp(4,F)$ are unimodular. 
We let 
\begin{equation}\label{s1s2defeq}
 s_1=\begin{bsmallmatrix}&1\\1\\&&&1\\&&1\end{bsmallmatrix}\quad\text{and}\quad
 s_2=\begin{bsmallmatrix}1\\&&1\\&-1\\&&&1\end{bsmallmatrix},
\end{equation}
and 
\begin{equation}
\label{deltaijdefeq}
\Delta_{i,j}
=
\begin{bsmallmatrix}
\varpi^{2i+j}&&&\\
&\varpi^{i+j}&&\\
&&\varpi^i&\\
&&&1
\end{bsmallmatrix}
\end{equation}
for integers $i$ and $j$.

\subsection*{Characters and representations} Let $G$ be a group
of $td$-type, as defined in \cite{Car},
and assume that the topology of $G$ has a countable basis. 

\runinhead{Characters} A \emph{character}\index{character} 
of $G$ is a continuous homomorphism $G \to \C^\times$. The \emph{trivial character} of $G$ is denoted
by $1_G$. If $\chi:F^\times \to \C^\times$ is a character,
then we let $a(\chi)$ be the smallest non-negative integer $n$
such that $a(1+\p^{n})=1$; here, we use the convention that
$1+\p^0=\OF^\times$. We refer to $a(\chi)$ \label{conductorpageref} as the 
\emph{conductor} of $\chi$. 
\index{conductor} 
\index{Dirichlet character!conductor}
We will use the continuous character $\nu: F^\times \to \C^\times$
defined by $\nu(x) = |x|$
for $x \in F^\times$; we
have $a(\nu)=0$.  \label{absolutevalchardef}
\runinhead{Schwartz functions} If $X$ is a closed subset of $G$, then the complex vector
space of locally constant and compactly supported functions
on $X$ will be denoted by $\mathcal{S}(X)$. \index{Schwartz functions}
\label{schwartzfndef}

\runinhead{Representations} A \emph{representation}\index{representation} of $G$
is a pair $(\pi,V)$, where $V$ is a complex vector space,
and $\pi:G \to \GL_{\C}(V)$ is a homomorphism. 
Let $(\pi,V)$ be a representation of $G$. We say that $\pi$
is \emph{smooth}\index{smooth representation}\index{representation!smooth}
if for every $v \in V$, the stabilizer of $v$ in $G$ is open. 
Evidently, characters of $G$ can be identified with one-dimensional
smooth representations of $G$.
We
say that $\pi$ is \emph{admissible}\index{admissible representation}\index{representation!admissible}
if $\pi$ is smooth and, for every compact open subgroup $K$ of $G$, the subspace of vectors
in $V$ fixed by every element of $K$ is finite-dimensional. Let $(\pi,V)$
be a smooth representation of $G$. We say that $\pi$ is 
\emph{irreducible}\index{irreducible representation}\index{representation!irreducible}
if the only $G$ subspaces of $V$ are $0$ and $V$. An 
\emph{irreducible constituent}\index{irreducible constituent}\index{representation!irreducible constituent},
or \emph{irreducible subquotient},\index{irreducible subquotient}\index{representation!irreducible subquotient}
of $\pi$ is an irreducible representation of $G$ that is isomorphic to $W/W'$, where
$W$ and $W'$ are $G$ subspaces of $V$ with $W'\subset W$. 
  The group $G$ acts on the complex
vector space of  linear functionals $\lambda: V \to \C$ via the formula $(g \cdot \lambda)(v)
=\lambda (\pi(g^{-1})v)$ for $v\in V$ and $g\in G$; we say that $\lambda$ is  \emph{smooth}
if there exists a compact, open subgroup $K$ of $G$ such that $k\cdot \lambda = \lambda$
for $k\in K$. The complex vector space of all smooth linear functionals on $V$ will be
denoted by $V^\vee$
and is a smooth representation $\pi^\vee$ of $G$ called the 
\emph{contragredient}\index{contragredient representation}\index{representation!contragredient}
of~$\pi$.  \label{contradef}
  If  $\pi$ admits a \emph{central character}\index{central character}\index{representation!central character},
then we denote it by $\omega_\pi$. 
\label{centralchardef}

\runinhead{Induction} Let $H$ be a closed subgroup of $G$, and let $(\pi,V)$ be a smooth representation
of $H$. Let $\cInd_H^G \pi$ be the complex vector space of all functions $f:G \to V$
such that $f(hg)=\pi(h)f(g)$ for $h\in H$ and $g\in G$, $f(gk)=f(g)$
for all $k$ in a compact, open subgroup of $G$ and $g\in G$, and there exists
a compact subset $X$ of $G$ such that $f(g)=0$ for $g$ not in $HX$. Then with 
the right translation action
defined by $(g'\cdot f)(g) = f(gg')$ for $g,g' \in G$ and $f$ in 
$\cInd_H^G\pi$, the space $\cInd_H^G \pi$ is a smooth representation of $G$. 
Next, assume that
$G$ is unimodular, and let $M$ and $U$ be closed subgroups of $G$ such that
$M\cap U=1$, $U$ is unimodular, $M$ normalizes $U$, the group $P=MU$ is closed in $G$, 
and $P \backslash G$ is compact. We define the 
\emph{modular character}\index{modular character}
\label{modchardef}
$\delta_P: P \to \C^\times$ of $P$ in the following way. Fix a Haar
measure on $U$. For $p$ in $P$ we let $\delta_P(p)$ be the unique
positive number such that 
\begin{equation}
\label{modchareq}
\int\limits_{U}f(p^{-1}up)\, du = \delta_P(p) \int\limits_U f(u)\, du
\end{equation}
for $f$ in $\mathcal{S}(U)$. Further assume that $(\pi,V)$ is a smooth representation
of $M$. The complex vector space of all functions $f:G \to V$ such that
$f(mug)=\delta_P(m)^{1/2}\pi(m)f(g)$ for $m\in M$, $u\in U$ and $g\in G$, and $f(gk)=f(g)$
for all $k$ in some compact, open subgroup of $G$ and $g\in G$, will be denoted
by $\Ind_P^G \pi$. With the right translation action, $\Ind_P^G \pi$ is a smooth
representation of $G$ called the representation obtained from $\pi$ by
\emph{normalized induction}\index{normalized induction}\index{representation!normalized induction}.

\runinhead{Jacquet modules} Let $N$ be a closed subgroup of $G$,  let $\chi$ be a character of $N$,
and let $(\pi,V)$ be a smooth representation of $G$. We define $V(N,\chi)$
to be the $\C$ subspace of $V$ spanned by all the vectors of the form $\chi (n)v - \pi (n) v$ for 
$n \in N$ and $v \in V$, 
and we set $V_{N,\chi} = V/V(N,\chi)$. The subspace $V(N,\chi)$ is a $N$ subspace of
$V$, and $N$ acts on the quotient $V_{N,\chi}$ by the character $\chi$. If $H$ is a closed
subgroup of $G$ that normalizes $N$ and is such that $\chi(hnh^{-1}) = \chi(n)$ for $n \in N$
and $h \in H$, then $V(N,\chi)$ is also an $H$ subspace, and $H$ acts on $V_{N,\chi}$.
We refer to $V_{N,\chi}$ as a \emph{Jacquet module}\index{Jacquet module}. \label{Jacquetmoddef}

\section{Representations}
\label{repsec}
In this section we review some needed definitions and results about 
representations of $\GSp(4,F)$. 

\subsection*{Parabolic induction}
In this work we will often consider representations of $\GSp(4,F)$
that are constructed via parabolic induction from the three standard
proper parabolic subgroups. 

The first, called the \emph{Borel subgroup}\index{Borel subgroup}\index{GSp@$\GSp(4,F)$!Borel subgroup}
of $\GSp(4,F)$, is the subgroup of
all upper triangular elements, i.e.
$$
B=\begin{bsmallmatrix}
*&*&*&*\\
&*&*&*\\
&&*&*\\
&&&*
\end{bsmallmatrix}. \label{BorelsubgroupGSp4}
$$
We also define the following subgroups of $B$:
\begin{equation}
\label{Udefeq}
T = \begin{bsmallmatrix}
*&&&\\
&*&&\\
&&*&\\
&&&*
\end{bsmallmatrix}, \qquad
U = 
\begin{bsmallmatrix}
1&*&*&*\\
&1&*&*\\
&&1&*\\
&&&1
\end{bsmallmatrix}.
\end{equation}
Then $T$ normalizes $U$, and $B=TU=UT$. 
Every element $h$ of $B$ can be written in the 
form 
$$
h=tu, \qquad \text{where} \quad
t=\begin{bsmallmatrix}
a&&&\\
&b&&\\
&&c b^{-1}&\\
&&&c a^{-1}
\end{bsmallmatrix}
\quad\text{and}\quad
u=
\begin{bsmallmatrix}
1&\vphantom{x}&\vphantom{y}&\vphantom{z}\\
&1&w&\vphantom{y}\\
&&1&\vphantom{-x}\\
&&&1
\end{bsmallmatrix}
\begin{bsmallmatrix}
1&x&y&z\\
&1&&y&\\
&&1&-x\\
&&&1
\end{bsmallmatrix}, 
$$
with $a,b,c \in F^\times$ and $w,x,y,z \in F$. 
With this notation,
$\lambda(h)=c$. Let $\chi_1$, $\chi_2$ and $\sigma$
be characters of $F^\times$. The representation of
$\GSp(4,F)$ obtained from the character of $T$ defined
by $t \mapsto \chi_1(a) \chi_2(b) \sigma (c)$ for 
$t$ as above via
normalized  induction as in Sect.~\ref{somedefsec} is denoted by 
$\chi_1 \times \chi_2 \rtimes \sigma$. The standard
model for $\chi_1 \times \chi_2 \rtimes \sigma$ is the complex vector space all the 
locally constant functions $f: \GSp(4,F) \to \C$
that satisfy
$$
f(hg) = |a^2b||c|^{-\frac{3}{2}} \chi_1(a) \chi_2(b) \sigma (c) f(g)
$$
for $g$ in $\GSp(4,F)$ and $h$ in $B$ as above;
the action of $\GSp(4,F)$ is by right 
translation. We remark
that the modular character (see \eqref{somedefsec}) of $B$ is given by
$\delta_B(h) = |a|^4 |b|^2 |c|^{-3}$ for $h$ as above, 
and that $\chi_1 \times \chi_2 \rtimes \sigma$
admits a central character, which is $\chi_1 \chi_2 \sigma^2$.

The 
\emph{Siegel parabolic subgroup}\index{Siegel parabolic subgroup}\index{GSp@$\GSp(4,F)$!Siegel parabolic subgroup}
of $\GSp(4,F)$ is defined to be 
$$
P=
\begin{bsmallmatrix}
*&*&*&*\\
*&*&*&*\\
&&*&*\\
&&*&*\\
\end{bsmallmatrix}. \label{Siegelparadef}
$$
We also define the following subgroups of $P$:
\begin{equation}
\label{NPdefeq}
M_P = 
\begin{bsmallmatrix}
*&*&&\\
*&*&&\\
&&*&*\\
&&*&*
\end{bsmallmatrix}, \qquad
N_P = 
\begin{bsmallmatrix}
1&&*&*\\
&1&*&*\\
&&1&\\
&&&1
\end{bsmallmatrix}.
\end{equation}
Then $M_P$ normalizes $N_P$, and we have $P = M_P N_P=N_P M_P$. 
An element $p$ of $P$ has a decomposition of the form
$$
p = mn, \qquad\text{where}\quad
m=\begin{bsmallmatrix}
a_1&a_2&&\\
a_3&a_4&&\\
&&c d_1 & c d_2\\
&&c d_3 & c d_4
\end{bsmallmatrix}\quad\text{and}\quad
n=\begin{bsmallmatrix}
1&&y&z\\
&1&x&y\\
&&1&\\
&&&1
\end{bsmallmatrix}
$$
with $A=\begin{bsmallmatrix} a_1 & a_2 \\ a_3 & a_4 \end{bsmallmatrix} \in \GL(2,F)$,
$\begin{bsmallmatrix} d_1 & d_2 \\ d_3 & d_4 \end{bsmallmatrix}
= A'$ (see 
\eqref{Aprimedefeq}), $c \in F^\times$, and $x,y,z \in F$. 
We note that $\lambda (p) = c$.
Let $(\pi,V)$ be an admissible representation of $\GL(2,F)$,
and let $\sigma$ be a character of $F^\times$. 
The representation of $\GSp(4,F)$ obtained from the 
representation of $M_P$ on $V$ defined by $m \mapsto \sigma (c) \pi(A)$
for $m$ as above 
via normalized  induction as in Sect.~\ref{somedefsec} is denoted
by $\pi \rtimes \sigma$. The modular character of $P$ (see \eqref{somedefsec})
is given by $\delta_P( p ) = |\det(A)|^3 |c|^{-3}$
for $p$ as above. Thus,
the standard model of $\pi \rtimes \sigma$ is comprised
of all the locally constant functions $f: \GSp(4,F) \to V$
such that
$$
f(pg) =  |c^{-1}\det(A) |^{\frac{3}{2}} \sigma (c) \pi (A) f(g)
$$
for $p$ as above and $g \in \GSp(4,F)$. If the representation
$\pi$ admits a central character $\omega_\pi$, then $\pi \rtimes \sigma$
admits $\omega_\pi \sigma^2$ as a central character. 

Finally, we
define the 
\emph{Klingen parabolic subgroup}\index{Klingen parabolic subgroup}\index{GSp@$\GSp(4,F)$!Klingen parabolic subgroup}
of $\GSp(4,F)$ to be
$$ \label{KlingenQdef}
Q = 
\begin{bsmallmatrix}
*&*&*&*\\
&*&*&*\\
&*&*&*\\
&&&*
\end{bsmallmatrix}.
$$
We also define the following subgroups of $Q$:
\begin{equation}
\label{NQdefeq}
M_Q=
\begin{bsmallmatrix}
*&&&\\
&*&*&\\
&*&*&\\
&&&*
\end{bsmallmatrix}, \qquad
N_Q = 
\begin{bsmallmatrix}
1&*&*&*\\
&1&&*\\
&&1&*\\
&&&1
\end{bsmallmatrix}.
\end{equation}
If $q$ is in $Q$, then $q$ has a decomposition
$$
q = mn, \qquad\text{where}\quad 
m=\begin{bsmallmatrix}
t&&&\\
&a&b&\\
&c&d&\\
&&&\lambda  t^{-1}
\end{bsmallmatrix} \quad\text{and}\quad
\begin{bsmallmatrix}
1&x&y&z\\
&1&&y\\
&&1&-x\\
&&&1
\end{bsmallmatrix},
$$
with $\begin{bsmallmatrix} a&b \\ c&d \end{bsmallmatrix} \in \GL(2,F)$, $t \in F^\times$,  $x,y,z \in F$,
and  $\lambda = \det(\begin{bsmallmatrix} a&b \\ c&d \end{bsmallmatrix})$. We note
that $\lambda(q)=\lambda$. 
Let $(\pi,V)$ be an admissible
representation of $\GSp(2,F)=\GL(2,F)$, and let $\chi$ be a character of $F^\times$.
The representation of $\GSp(4,F)$ obtained
by normalized induction as in Sect.~\ref{somedefsec} from the representation of $M_Q$ on $V$ defined
by $m \mapsto \chi(t) \pi (\begin{bsmallmatrix} a & b \\ c & d \end{bsmallmatrix})$
for $m$ as above is denoted by $\chi \rtimes \pi$. 
The modular character of $Q$ is given by $\delta_Q(q) =
|t|^4 |ad-bc|^{-2}$ for $q$ as above. The standard model of $\chi \rtimes \pi$
consists of all the locally constant functions $f: \GSp(4,F) \to V$ such that
$$
f(qg) = |t^2(ad-bc)^{-1}| \chi(t) \pi ( 
\begin{bsmallmatrix}
a&b \\ c&d
\end{bsmallmatrix}
)f(g)
$$
for $q$ as above and $g \in \GSp(4,F)$. 

\subsection*{Generic representations}
Let $U$ be the subgroup of $B$ defined in \eqref{Udefeq}.
Let $c_1,c_2 \in F$, and let $\psi:F\to \C^\times$
be the character fixed in Sect.~\ref{somedefsec}. We define
a character $\psi_{c_1,c_2}:U \to \C^\times$ by 
\begin{equation}
\label{psic1c2defeq}
\psi_{c_1,c_2} (
\begin{bsmallmatrix}
1&x&*&*\\
&1&y&*\\
&&1&-x\\
&&&1
\end{bsmallmatrix} ) = \psi(c_1 x + c_2 y), \qquad x,y \in F.
\end{equation}
Let $(\pi,V)$ be an irreducible, admissible representation
of $\GSp(4,F)$. We say that $\pi$ is 
\emph{generic}\index{generic representation}\index{representation!generic}
if $\Hom_{U}(V,\psi_{c_1,c_2})$ is non-zero for some $c_1,c_2 \in F^\times$; we note that this
does not depend on the choice of $c_1$ and $c_2$ in $F^\times$. Assume
that $\pi$ is generic. Then $\pi$ has a 
Whittaker model.
This means that $\pi$ is isomorphic to a subspace 
of the complex vector space of all functions $W: \GSp(4,F) \to \C$
such that
$$
W ( \begin{bsmallmatrix}
1&x&*&*\\
&1&y&*\\
&&1&-x\\
&&&1
\end{bsmallmatrix} g) = \psi (c_1 x + c_2 y) W(g), \qquad x,y \in F, \quad g \in \GSp(4,F).
$$
Here, $\GSp(4,F)$ acts on the space of all such functions by right translation. 
The subspace isomorphic to $\pi$ is unique, and will be denoted by $\mathcal{W}(\pi,\psi_{c_1,c_2})$.
We call $\mathcal{W}(\pi,\psi_{c_1,c_2})$ the \emph{Whittaker model} of $\pi$
with respect to $\psi_{c_1,c_2}$. 
\index{Whittaker model}
In the remainder of this work, if we consider a generic representation $\pi$ and its
Whittaker model $\mathcal{W}(\pi,\psi_{c_1,c_2})$, then we always assume that $c_1,c_2 \in \OF^\times$. \label{whitmodeldef}

\subsection*{The list of non-supercuspidal representations} Thanks to \cite{ST}, it is possible
to provide a convenient list of all the non-supercuspidal, irreducible, admissible representations
of $\GSp(4,F)$. This list is described in \cite{NF}; here, we provide a brief summary. The list
also appears in Table \ref{maintable}. The list consists of eleven groups of representations, denoted
by capital roman numerals. These groups
are further divided into (sometimes only one, and at most four) subgroups, denoted by small roman letters (e.g., IIa and IIb);
in the case that a group admits only one subgroup, the letter is omitted.
Every non-supercuspidal, irreducible, admissible representation
of $\GSp(4,F)$ is isomorphic to a representation that is a member of a group and a subgroup. 
Representations from distinct subgroups are not isomorphic, with the exception that subgroup Vb equals subgroup Vc if the parameters are allowed to run through all possibilities.
The notation is such that 
representations of an ``a'' subgroup (e.g., IIa) are generic (if a group consists only
of one subgroup, then all the representations in that group are generic). 
\runinhead{Group I} The members of this group are the irreducible representations
of the form $\chi_1 \times \chi_2 \rtimes \sigma$, where $\chi_1,\chi_2$ and $\sigma$ are
characters of $F^\times$. A representation of this form is irreducible if and only if
$\chi_1 \neq \nu^{\pm 1}$, $\chi_2 \neq \nu^{\pm 1}$, and $\chi_2 \neq \chi_1^{\pm 1} \nu^{\pm 1}$. 
This group of representations  has only one subgroup. 
\runinhead{Group II} Let $\chi$ be a character of $F^\times$ such that $\chi \neq \nu^{\pm \frac{3}{2}}$
and $\chi^2 \neq \nu^{\pm 1}$,
and let $\sigma$ be a character of $F^\times$.
The members of Group II are the irreducible constituents of the representations of the form
$\nu^{\frac{1}{2}} \chi \times \nu^{-\frac{1}{2}} \chi \rtimes \sigma$. 
The representation $\nu^{\frac{1}{2}} \chi \times \nu^{-\frac{1}{2}} \chi \rtimes \sigma$
has two irreducible constituents. 
One of these irreducible constituents, the subrepresentation, is $\chi \St_{\GL(2)} \rtimes \sigma$;
here $\St_{\GL(2)}$ is the Steinberg representation of $\GL(2,F)$. The representations
of the form $\chi \St_{\GL(2)} \rtimes \sigma$ comprise subgroup IIa. The other irreducible constituent, the quotient,
of $\nu^{\frac{1}{2}} \chi \times \nu^{-\frac{1}{2}} \chi \rtimes \sigma$
is $\chi 1_{\GL(2)}  \rtimes \sigma$. The representations of the 
form $\chi 1_{\GL(2)} \rtimes \sigma$
make up subgroup  IIb.
\runinhead{Group III} Let $\chi$ be a character of $F^\times$ such that
$\chi \neq 1_{F^\times}$ and $\chi \neq \nu^{\pm 2}$, and let $\sigma$
be a character of $F^\times$. The elements of Group III
are the irreducible constituents of the representations of the form $\chi \times \nu \rtimes \nu^{-\frac{1}{2}} \sigma$. 
The representation $\chi \times \nu \rtimes \nu^{-\frac{1}{2}} \sigma$ has two irreducible
constituents. The irreducible representation $\chi \rtimes \sigma \St_{\GSp(2)}$ is
a subrepresentation of $\chi \times \nu \rtimes \nu^{-\frac{1}{2}} \sigma$, and representations
of the form $\chi \rtimes \sigma \St_{\GSp(2)}$ make up subgroup IIIa. The irreducible 
representation $\chi \rtimes \sigma 1_{\GSp(2)}$ is a quotient of $\chi \times \nu \rtimes \nu^{-\frac{1}{2}}\sigma$,
and representations of the form $\chi \rtimes \sigma 1_{\GSp(2)}$ comprise subgroup IIIb. 
\runinhead{Group IV} Let $\sigma$ be a character of $F^\times$. The irreducible
constituents of the representations of the form $\nu^2 \times \nu \rtimes \nu^{-\frac{3}{2}} \sigma$ make
up Group IV. The representation $\nu^2 \times \nu \rtimes \nu^{-\frac{3}{2}} \sigma$ has four
irreducible constituents, which can be described as follows. There are two
exact sequences
\begin{equation}
\label{groupIVseq1eq}
0 \to \nu^{\frac{3}{2}} \St_{\GL(2)} \rtimes \nu^{-\frac{3}{2}}\sigma 
\to \nu^2 \times \nu \rtimes \nu^{-\frac{3}{2}} \sigma
\to \nu^{\frac{3}{2}} 1_{\GL(2)} \rtimes \nu^{-\frac{3}{2}} \sigma \to 0
\end{equation}
and
\begin{equation}
\label{groupIVseq2eq}
0 \to \nu^2 \times \nu^{-1}\sigma \St_{\GSp(2)}
\to  \nu^2 \times \nu \rtimes \nu^{-\frac{3}{2}} \sigma
\to \nu^2 \rtimes \nu^{-1} \sigma 1_{\GSp(2)} \to 0.
\end{equation}
In these exact sequences, the second and fourth representations 
are reducible, and each has two irreducible constituents as in the following
table. The subrepresentations are on the top and left, and the quotients are on the bottom and right.
\begin{equation}
\label{groupIVtableeq}
\begin{array}{ccc}
\toprule
& \nu^{\frac{3}{2}}\St_{\GL(2)} \rtimes \nu^{-\frac{3}{2}}\sigma & \nu^{\frac{3}{2}} 1_{\GL(2)} \rtimes \nu^{-\frac{3}{2}}\sigma\\
\cmidrule{2-3}
\nu^2 \rtimes \nu^{-1}\sigma \St_{\GSp(2)} & \sigma \St_{\GSp(4)} & L(\nu^2,\nu^{-1}\sigma\St_{\GSp(2)})\\
\nu^2 \rtimes \nu^{-1} \sigma 1_{\GSp(2)} & L(\nu^{\frac{3}{2}}\St_{\GL(2)},\nu^{-\frac{3}{2}}\sigma) & \sigma 1_{\GSp(4)}\\
\bottomrule
\end{array}
\end{equation}
Thus,  $\sigma \St_{\GSp(4)}$, $L(\nu^2,\nu^{-1}\sigma\St_{\GSp(2)})$, $L(\nu^{\frac{3}{2}}\St_{\GL(2)},\nu^{-\frac{3}{2}}\sigma)$,
and $\sigma 1_{\GSp(4)}$ are the irreducible constituents of $\nu^2 \times \nu \rtimes \nu^{-\frac{3}{2}} \sigma$,
and representations of this form constitute, respectively, subgroups IVa, IVb, IVc and IVd. 
\runinhead{Group V} Let $\xi$ be a non-trivial quadratic character of $F^\times$, and
let $\sigma$ be a character of $F^\times$. Then the irreducible
constituents of representations of the form $\nu\xi \rtimes \xi \rtimes \nu^{-\frac{1}{2}}\sigma$ 
make up Group V. The representation $\nu\xi \rtimes \xi \rtimes \nu^{-\frac{1}{2}}\sigma$  has four
irreducible constituents which are delineated as follows. There are two exact sequences
\begin{equation}
\label{groupVseq1eq}
0\to
\nu^{\frac{1}{2}}\xi\St_{\GL(2)} \rtimes \nu^{-\frac{1}{2}}\sigma
\to
\nu\xi \times \xi \rtimes \nu^{-\frac{1}{2}}\sigma
\to
\nu^{\frac{1}{2}} \xi 1_{\GL(2)} \rtimes \nu^{-\frac{1}{2}} \sigma
\to
0
\end{equation}
and
\begin{equation}
\label{groupVseq2eq}
0
\to
\nu^{\frac{1}{2}}\xi \St_{\GL(2)} \rtimes \xi \nu^{-\frac{1}{2}} \sigma
\to
\nu\xi \times \xi \rtimes \nu^{-\frac{1}{2}}\sigma
\to
\nu^{\frac{1}{2}}\xi 1_{\GL(2)} \rtimes \xi \nu^{-\frac{1}{2}} \sigma
\to
0.
\end{equation}
In these exact sequences, the second and fourth representations 
are reducible, and each has two irreducible constituents as in the following
table. The subrepresentations are on the top and left, and the quotients are on the bottom and right.
\begin{equation}
\label{groupVtableeq}
\begin{array}{ccc}
\toprule
&\nu^{\frac{1}{2}}\xi\St_{\GL(2)}\rtimes \nu^{-\frac{1}{2}}\xi\sigma& \nu^{\frac{1}{2}}\xi 1_{\GL(2)} \rtimes \xi\nu^{-\frac{1}{2}} \sigma\\
\cmidrule{2-3}
\nu^{\frac{1}{2}}\xi\St_{\GL(2)}\rtimes \nu^{-\frac{1}{2}}\sigma & \delta([\xi,\nu\xi],\nu^{-\frac{1}{2}}\sigma) & L(\nu^{\frac{1}{2}}\xi\St_{\GL(2)},\nu^{-\frac{1}{2}}\sigma)\\
\nu^{\frac{1}{2}}\xi1_{\GL(2)}\rtimes \nu^{-\frac{1}{2}}\sigma & L(\nu^{\frac{1}{2}}\xi \St_{\GL(2)}, \nu^{-\frac{1}{2}}\xi \sigma)
&L(\nu\xi,\xi \rtimes \nu^{-\frac{1}{2}}\sigma)\\
\bottomrule
\end{array}
\end{equation}
The representation $\nu\xi \rtimes \xi \rtimes \nu^{-\frac{1}{2}}\sigma$ has as irreducible constituents
 $\delta([\xi,\nu\xi],\nu^{-\frac{1}{2}}\sigma)$, $L(\nu^{\frac{1}{2}}\xi\St_{\GL(2)},\nu^{-\frac{1}{2}}\sigma)$,
$L(\nu^{\frac{1}{2}}\xi \St_{\GL(2)}, \nu^{-\frac{1}{2}}\xi \sigma)$, and $L(\nu\xi,\xi \rtimes \nu^{-\frac{1}{2}}\sigma)$, and  
representations of this form
make up subgroups Va, Vb, Vc and Vd, respectively. We note that the Va representation $\delta([\xi,\nu\xi],\nu^{-\frac{1}{2}}\sigma)$
is essentially square integrable, and while the Vb representation $L(\nu^{\frac{1}{2}}\xi\St_{\GL(2)},\nu^{-\frac{1}{2}}\sigma)$
and  the Vc representation $L(\nu^{\frac{1}{2}}\xi \St_{\GL(2)}, \nu^{-\frac{1}{2}}\xi \sigma)$ are twists of each other, these representations
are not isomorphic.
\runinhead{Group VI} Let $\sigma$ be a character of $F^\times$. 
The irreducible constituents of representations of the form 
$\nu \times 1_{F^\times} \rtimes \nu^{-\frac{1}{2}} \sigma$ make up Group VI.
The representation $\nu \times 1_{F^\times} \rtimes \nu^{-\frac{1}{2}} \sigma$
has four irreducible constituents which can
be described as follows. There are two exact sequences
\begin{equation}
\label{groupVIseq1eq}
0
\to
\nu^{\frac{1}{2}} \St_{\GL(2)} \rtimes\nu^{-\frac{1}{2}} \sigma
\to
\nu \times 1_{F^\times} \rtimes \nu^{-\frac{1}{2}} \sigma 
\to
\nu^{\frac{1}{2}} 1_{\GL(2)} \rtimes \nu^{-\frac{1}{2}} \sigma
\to
0
\end{equation}
and
\begin{equation}
\label{groupVIseq2eq}
0
\to
1_{F^\times} \rtimes \sigma \St_{\GSp(2)}
\to 
\nu \times 1_{F^\times} \rtimes \nu^{-\frac{1}{2}} \sigma
\to
1_{F^\times} \rtimes \sigma 1_{\GSp(2)}
\to
0.
\end{equation}
In these exact sequences, the second and fourth representations 
are reducible, and each has two irreducible constituents as in the following
table. The subrepresentations are on the top and left, and the quotients are on the bottom and right.
\begin{equation}
\label{groupVItableeq}
\begin{array}{ccc}
\toprule
&\nu^{\frac{1}{2}}\St_{\GL(2)} \rtimes \nu^{-\frac{1}{2}}\sigma & \nu^{\frac{1}{2}}1_{\GL(2)} \rtimes \nu^{-\frac{1}{2}}\sigma\\
\cmidrule{2-3}
1_{F^\times} \rtimes \sigma \St_{\GSp(2)} & \tau (S,\nu^{-\frac{1}{2}}\sigma) & \tau (T,\nu^{-\frac{1}{2}} \sigma)\\
1_{F^\times} \rtimes \sigma 1_{\GSp(2)} & L(\nu^{\frac{1}{2}} \St_{\GL(2)}, \nu^{-\frac{1}{2}} \sigma) & L(\nu, 1_{F^\times} \rtimes 
\nu^{-\frac{1}{2}} \sigma)\\
\bottomrule
\end{array}
\end{equation}
The representation $\nu \times 1_{F^\times} \rtimes \nu^{-\frac{1}{2}} \sigma$ has the representations
$\tau (S,\nu^{-\frac{1}{2}}\sigma)$, 
$\tau (T,\nu^{-\frac{1}{2}} \sigma)$, 
$L(\nu^{\frac{1}{2}} \St_{\GL(2)}, \nu^{-\frac{1}{2}} \sigma)$,
and $L(\nu, 1_{F^\times} \rtimes \nu^{-\frac{1}{2}} \sigma)$ as irreducible constituents,
and representations of this form constitute subgroups VIa, VIb, VIc, and VId, respectively.
We note that the representations $\tau (S,\nu^{-\frac{1}{2}}\sigma)$ and $\tau (T,\nu^{-\frac{1}{2}} \sigma)$ are essentially
tempered but not essentially square integrable. 
\runinhead{Group VII} The members of this group are the irreducible representations of the form $\chi \rtimes \pi$,
where $\chi$ is a character of $F^\times$, and $\pi$ is a supercuspidal, irreducible, admissible representation
of $\GSp(2,F) = \GL(2,F)$. A representation of this form is irreducible if and only if $\chi \neq 1_{F^\times}$
and $\chi \neq \xi \nu^{\pm 1}$ for every character $\xi$ of $F^\times$ of order two such that $\xi \pi \cong \pi$. 
This group of representations has  only one subgroup. 
\runinhead{Group VIII} Let $\pi$ be a supercuspidal, irreducible, admissible representation
of $\GSp(2,F) = \GL(2,F)$. The members
of Group VIII are the irreducible constituents 
of representations of the form  $1 \rtimes \pi$.
The representation $1 \rtimes \pi$ is the direct sum of two 
essentially tempered, irreducible, admissible representations $\tau(S,\pi)$ and $\tau(T,\pi)$.
Representations of the form $\tau(S,\pi)$ and $\tau(T,\pi)$ make up subgroups VIIIa and VIIIb,
respectively. 
\runinhead{Group IX} Let $\xi$ be a non-trivial quadratic character of $F^\times$, and let 
$\pi$ be a supercuspidal, irreducible, admissible representation of $\GSp(2,F) = \GL(2,F)$
such that $\xi \pi \cong \pi$. The members of group IX are the irreducible constituents
of the representations of the form $\nu \xi \rtimes \nu^{-\frac{1}{2}} \pi$. The representation
$\nu \xi \rtimes \nu^{-\frac{1}{2}} \pi$ has two irreducible constituents. One of these irreducible
constituents, the subrepresentation $\delta (\nu \xi, \nu^{-\frac{1}{2}} \pi)$, is essentially
square-integrable; the quotient $L(\nu \xi, \nu^{-\frac{1}{2}} \pi)$ is non-tempered. Representations
of the form $\delta (\nu \xi, \nu^{-\frac{1}{2}} \pi)$ make up subgroup IXa, and
representations of the form $L(\nu \xi, \nu^{-\frac{1}{2}} \pi)$ make up subgroup IXb.
\runinhead{Group X} The members of this group are the irreducible representations of
the form $\pi \rtimes \sigma$ where $\pi$ is a supercuspidal, irreducible, admissible
representation of $\GL(2,F)$, and $\sigma$ is a character of $F^\times$. A representation
of this form is irreducible if and only if $\pi$ is not of the form $\nu^{\pm \frac{1}{2}} \rho$
where $\rho$ is a supercuspidal, irreducible, admissible
representation of $\GL(2,F)$ with trivial central character. This group has only one
subgroup.
\runinhead{Group XI} Let $\pi$ be a supercuspidal, irreducible, admissible
representation of $\GL(2,F)$ with trivial central character, and let $\sigma$ be 
a character of $F^\times$. The members of group XI are the irreducible constituents
of representations of the form $\nu^{\frac{1}{2}}\pi \rtimes \nu^{-\frac{1}{2}} \sigma$. 
The representation $\nu^{\frac{1}{2}}\pi \rtimes \nu^{-\frac{1}{2}} \sigma$ has two irreducible
constituents. One of these constituents, the subrepresentation $\delta(\nu^{\frac{1}{2}}\pi,\nu^{-\frac{1}{2}}\sigma)$,
is essentially square integrable; the quotient $L(\nu^{\frac{1}{2}}\pi,\nu^{-\frac{1}{2}}\sigma)$ is
non-tempered. Representations of the form $\delta(\nu^{\frac{1}{2}}\pi,\nu^{-\frac{1}{2}}\sigma)$
make up subgroup XIa, and representations of the form $L(\nu^{\frac{1}{2}}\pi,\nu^{-\frac{1}{2}}\sigma)$
make up subgroup XIb. 

\subsection*{Saito-Kurokawa representations}
\label{SKsubsec}
Let $\pi$ be an infinite-dimensional, irreducible, admissible
representation of $\GL(2,F)$ with trivial central character,
and let $\sigma$ be a character of $F^\times$. Assume
that $\pi \ncong \nu^{\frac{3}{2}} \times \nu^{-\frac{3}{2}}$. 
Then the representation $\nu^{\frac{1}{2}}\pi \rtimes \nu^{-\frac{1}{2}}\sigma$
of $\GSp(4,F)$ has a unique irreducible quotient $Q(\nu^{\frac{1}{2}}\pi,\nu^{-\frac{1}{2}}\sigma)$;
we say that $Q(\nu^{\frac{1}{2}}\pi,\nu^{-\frac{1}{2}}\sigma)$
is a
\emph{Saito-Kurokawa} representation.
\index{Saito-Kurokawa representation}
\index{representation!Saito-Kurokawa}
The central character of $Q(\nu^{\frac{1}{2}}\pi,\nu^{-\frac{1}{2}}\sigma)$ is $\sigma^2$. 
The Saito-Kurokawa representations are all the representations
in groups IIb, Vb, Vc, VIc, and XIb, and form an important
family of non-generic representations of $\GSp(4,F)$. See Sect.~5.5 of \cite{NF}
for more information. 

\subsection*{The $P_3$-quotient}
\label{jacobisubsec}
As developed in \cite{NF} and \cite{RS2}, there is a useful map of representations of $\GSp(4,F)$ 
 to representations of the group $P_3$. Here, $P_3$ \index{Paa@$P_3$}
is the subgroup of $\GL(3,F)$ defined as
\begin{equation}
\label{P3defeq}
P_3 = 
\begin{bsmallmatrix}
*&*&*\\
*&*&*\\
&&1
\end{bsmallmatrix}.
\end{equation}
To explain the map we require some more definitions. 
The \emph{Jacobi subgroup}
of $\GSp(4,F)$ is defined to be
$$
G^J = 
\begin{bsmallmatrix}
1&*&*&*\\
&*&*&*\\
&*&*&*\\
&&&1
\end{bsmallmatrix}. \label{Jacobigroupdef}
$$
The Jacobi subgroup $G^J$\index{Jacobi subgroup}\index{GSp@$\GSp(4,F)$!Jacobi subgroup}
is a normal subgroup of the Klingen parabolic subgroup $Q$.
The center of the Jacobi group is
\begin{equation}
\label{ZJdefeq}
Z^J = 
\begin{bsmallmatrix}
1&&&* \\
&1&&\\
&&1&\\
&&&1
\end{bsmallmatrix}.
\end{equation}
The group $Z^J$\index{center of the Jacobi subgroup}
\index{GSp@$\GSp(4,F)$!center of the Jacobi subgroup} is also a normal subgroup of $Q$. Next, let $\bar Q$ be the following subgroup of $Q$:
$$ \label{Qbardef}
\bar Q = 
\begin{bsmallmatrix}
*&*&*&*\\
&*&*&*\\
&*&*&*\\
&&&1
\end{bsmallmatrix}.
$$
Then $Z^J$ is a normal subgroup of $\bar Q$. Every element of $\bar Q$
can be uniquely written in the form 
$$
q = 
\begin{bsmallmatrix}
ad-bc &\vphantom{-y}&\vphantom{x}&\vphantom{z}\\
&a&b&\vphantom{x} \\
&c&d&\vphantom{y}\\
&&&1
\end{bsmallmatrix}
\begin{bsmallmatrix}
1&-y&x&z\\
&1&&x\\
&&1&y\\
&&&1
\end{bsmallmatrix}
$$
with $\begin{bsmallmatrix}a&b \\ c&d \end{bsmallmatrix}\in \GL(2,F)$ and $x,y,z \in F$. Define
$i: \bar Q \to P_3$ by 
$$
i(q) = 
\begin{bsmallmatrix} 
a&b&\\
c&d&\\
&&1
\end{bsmallmatrix}
\begin{bsmallmatrix}
1&&x\\
&1&y\\
&&1
\end{bsmallmatrix}
$$
for $q$ as above. Then $i$ is a surjective homomorphism with kernel
$Z^J$, and thus induces an isomorphism $\bar Q/ Z^J \cong P_3$. Now suppose
that $(\pi,V)$ is a smooth representation of $\GSp(4,F)$. Let $V(Z^J)$
be the $\C$ subspace of $V$ spanned by the vectors $v - \pi (z)v$
for $v \in V$ and $z \in Z^J$. Since $\bar Q$ normalizes $Z^J$, the subspace 
$V(Z^J)$ is a $\bar Q$ subspace of $V$. 
It follows that the quotient $V_{Z^J} = V/V(Z^J)$
is also a representation
of $\bar Q$. Since $Z^J$ acts trivially on $V_{Z^J}$, the group $\bar Q/Z^J \cong P_3$
acts on $V_{Z^J}$, so that $V_{Z^J}$ may be regarded as a representation
of $P_3$. 
We refer to $V_{Z^J}$ as the \emph{$P_3$-quotient}\index{Paaa@$P_3$-quotient} of \label{VZJdef}
$V$. We will write $p: V \to V_{Z^J}$ for the projection from $V$ to $V_{Z^J}$.
We have $p(\pi(q)v)=i(q)v$ for $v \in V$ and $q \in \bar Q$.

To state the main result about the $P_3$-quotient of an irreducible,
admissible representation of $\GSp(4,F)$ we need to first recall some facts
about representations of $P_3$, as explained in \cite{BZ}.  Define a subgroup
$U_3$ 
of $P_3$ by \label{U3def}
$$
U_3 = 
\begin{bsmallmatrix}
1&*&*\\
&1&*\\
&&1
\end{bsmallmatrix}.
$$
Define characters $\Theta$ and $\Theta'$ of $U_3$ by
$$
\Theta(
\begin{bsmallmatrix}
1&u_{12}&*\\
&1&u_{23}\\
&&1
\end{bsmallmatrix}) = \psi(u_{12}+u_{23}), \qquad
\Theta'(
\begin{bsmallmatrix}
1&u_{12}&*\\
&1&u_{23}\\
&&1
\end{bsmallmatrix}) = \psi(u_{23})
$$
where $u_{12}, u_{23} \in F$. There are three families of
representations of $P_3$, obtained from representations
of the groups $\GL(0,F)=1$,
$\GL(1,F)=F^\times$ and $\GL(2,F)$
by an induction. First of all, we define
$$
\tau_{\GL(0)}^{P_3} (1)= \cInd_{U_3}^{P_3} (\Theta). 
$$
The representation $\tau_{\GL(0)}^{P_3} (1)$ is smooth
and irreducible. 
Second, let $(\chi,V)$ be smooth representation of $F^\times$. 
Define a smooth representation $ \Theta' \otimes \chi$ of the
subgroup $\begin{bsmallmatrix} *&*&*\\&1&*\\&&1 \end{bsmallmatrix}$
of $P_3$ on the space $V$ of $\chi$ by 
$$
(\Theta'\otimes \chi) (
\begin{bsmallmatrix}
a&*&*\\
&1&y\\
&&1
\end{bsmallmatrix})v = \Theta'( \begin{bsmallmatrix}
a&*&*\\
&1&y\\
&&1
\end{bsmallmatrix}) \chi (a) v = \psi(y) \chi (a) v
$$
for $a \in F^\times$, $y \in F$, and $v \in V$. Then we define
$$
\tau_{\GL(1)}^{P_3} (\chi) = \cInd_{\begin{bsmallmatrix} *&*&*\\&1&*\\&&1 \end{bsmallmatrix}}^{P_3} (\Theta' \otimes \chi)
$$
The representation $\tau_{\GL(1)}^{P_3} (\chi)$ is smooth, and if $\chi$ is irreducible (i.e., a character of $F^\times$), 
then $\tau_{\GL(1)}^{P_3} (\chi)$ is irreducible. 
Third, let $\rho$ be a smooth representation of $\GL(2,F)$. We define
$$
\tau_{\GL(2)}^{P_3} (\rho)
$$
to have the same space as $\rho$, with action given by 
$$
\tau_{\GL(2)}^{P_3} (\rho) ( \begin{bsmallmatrix} a &b & * \\ c&d&* \\ &&1 \end{bsmallmatrix}) = 
\rho(\begin{bsmallmatrix} a&b \\ c&d \end{bsmallmatrix}),
$$
where $\begin{bsmallmatrix} a&b \\ c&d \end{bsmallmatrix} \in \GL(2,F)$. Clearly, $\tau_{\GL(2)}^{P_3} (\rho)$
is smooth, and if $\rho$ is irreducible, then $\tau_{\GL(2)}^{P_3} (\rho)$ is irreducible. Let $\eta$
be an irreducible,
smooth representation  of $P_3$. Then $\eta$ is isomorphic to a representation of the form $\tau_{\GL(0)}^{P_3} (1)$,
$\tau_{\GL(1)}^{P_3} (\chi)$, where $\chi$ is a character of $F^\times$, or $\tau_{\GL(2)}^{P_3} (\rho)$
where $\rho$ is an irreducible, admissible representation of $P_3$; moreover, if $\eta \cong \tau_{\GL(k_1)}^{P_3}(\pi_1)$
and $\eta \cong \tau_{\GL(k_2)}^{P_3}(\pi_2)$ where $k_1,k_2 \in \{0,1,2\}$ and $\pi_1$ and $\pi_2$ are of the types
just described, then $k_1=k_2$ and $\pi_1 \cong \pi_2$. 

\begin{theorem}
\label{mainP3theorem}
Let $(\pi,V)$ be an irreducible, admissible representation of the group $\GSp(4,F)$. 
The quotient $V_{Z^J} = V/V(Z^J)$ is a smooth representation of $\bar Q/Z^J$,
and hence defines a smooth representation of $P_3 \cong \bar Q/Z^J$. As
a representation of $P_3$, $V_{Z^J}$ has a finite filtration by $P_3$
subspaces such that the successive quotients are irreducible and of the form
$\tau_{\GL(0)}^{P_3}(1)$, $\tau_{\GL(1)}^{P_3}(\chi)$, or $\tau_{\GL(2)}^{P_3}(\rho)$,
where $\chi$ is a character of $F^\times$, and $\rho$ is an irreducible,
admissible representation of $\GL(2,F)$. Moreover:
\begin{enumerate}
\item \label{mainP3theoremitem1} There exists a chain of $P_3$ subspaces
$$
0 \subset V_2 \subset V_1 \subset V_0=V_{Z^J}
$$
such that
\begin{align*}
V_2 &\cong \dim \Hom_{U}(V, \psi_{-1,1}) \cdot \tau_{\GL(0)}^{P_3}(1),\\
V_1/V_2 & \cong \tau_{\GL(1)}^{P_3} (V_{U,\psi_{-1,0}}),\\
V_0/V_1 & \cong \tau_{\GL(2)}^{P_3} (V_{N_Q}).
\end{align*}
Here, the vector space $V_{U,\psi_{-1,0}}$ admits a smooth action of $\GL(1,F) \cong F^\times$
induced by the operators
$$
\pi(\begin{bsmallmatrix} a&&& \\ &a&& \\ &&1& \\ &&&1 \end{bsmallmatrix}), \qquad a \in F^\times,
$$
and $V_{N_Q}$ admits a smooth action of $\GL(2,F)$ induced by the operators
$$
\pi(\begin{bsmallmatrix} \det(g) && \\ & g & \\ &&1 \end{bsmallmatrix}, \qquad g \in \GL(2,F).
$$
\item \label{mainP3theoremitem2} The representation $\pi$ is generic if and only if $V_2 \neq 0$, and if $\pi$ is generic,
then $V_2 \cong \tau_{\GL(0)}^{P_3} (1)$. 
\item \label{mainP3theoremitem3} We have $V_2 \cong V_{Z^J}$ if and only if $\pi$ is supercuspidal. If $\pi$
is supercuspidal and generic, then $V_{Z^J} = V_2 \cong \tau_{\GL(0)}^{P_3} (1)$ is non-zero
and irreducible. If $\pi$ is supercuspidal and non-generic, then $V_{Z^J} = V_2 =0$. 
\end{enumerate}
\end{theorem}
\begin{proof}
See Sect.~2.5 of \cite{NF} and Sect.~3 of \cite{RS2}.
\end{proof}

The semisimplifications of the quotients $V_0/V_1$ and $V_1/V_2$ from Theorem \ref{mainP3theorem}
for all irreducible, admissible representations $(\pi,V)$ of $\GSp(4,F)$ are listed in Appendix A.4
of \cite{NF}. Note that this reference assumed that $\pi$ has trivial central character; however,
as pointed out in \cite{RS2}, the results are exactly the same without this assumption. Note also that there is typo in Table
A.5 of \cite{NF}: The entry for Vd in the ``s.s.($V_0/V_1$)'' column should be 
$\tau_{\GL(2)}^{P_3}(\nu (\nu^{-\frac{1}{2}} \sigma \times \nu^{-\frac{1}{2}}\xi \sigma))$. 

\subsection*{Zeta integrals}
\label{zetaintegralssubsecstar}
Let $\pi$ be a generic, irreducible, admissible representation of 
$\GSp(4,F)$ with trivial central character. Let $c_1, c_2 \in \OF^\times$, and let $\mathcal{W}(\pi,\psi_{c_1,c_2})$
be the Whittaker model of $\pi$ defined with respect to $\psi_{c_1,c_2}$, as defined above.
For $W \in \mathcal{W}(\pi,\psi_{c_1,c_2})$ and $s \in \C$, we define the \emph{zeta integral}
\index{zeta integral} $Z(s,W)$ by 
$$
Z(s,W) = \int\limits_{F^\times} \int\limits_{F} 
W(
\begin{bsmallmatrix}
a&&&\\
&a&&\\
&x&1&\\
&&&1
\end{bsmallmatrix}
)|a|^{s-\frac{3}{2}}\, dx\,d^\times a. \label{zetaidef}
$$
The following results  were proven in \cite{NF}. There exists a real number $s_0$
such that, for all $W \in \mathcal{W}(\pi,\psi_{c_1,c_2})$,  $Z(s,W)$ converges for  $s \in \C$ with $\Real(s)>s_0$  
to an element of $\C(q^{-s})$. Thus,  $Z(s,W)$ for $W \in \mathcal{W}(\pi,\psi_{c_1,c_2})$ has a 
meromorphic continuation to $\C$. Let $I(\pi)$ be the $\C$ vector subspace of $\C(q^{-s})$
spanned by the $Z(s,W)$ for $W \in \mathcal{W}(\pi,\psi_{c_1,c_2})$. Then $I(\pi)$ is independent
of the choice of $c_1,c_2 \in \OF^\times$. Moreover, $I(\pi)$ is a non-zero $\C[q^{-s},q^s]$ module
containing $\C$, and there exists $R(X) \in \C[X]$ such that $R(q^{-s}) I(\pi) \subset \C[q^{-s},q^s]$,
so that $I(\pi)$ is a fractional ideal of the principal ideal domain $\C[q^{-s},q^s]$ whose
quotient field is $\C(q^{-s})$. The fractional ideal $I(\pi)$ admits a generator of the form
$1/Q(q^{-s})$ with $Q(0)=1$, where $Q(X) \in \C[X]$. We define
\begin{equation}
\label{Lfunctiondefeq}
L(s,\pi) = \dfrac{1}{Q(q^{-s})},
\end{equation}
and call $L(s,\pi)$ the \emph{$L$-function}\index{Lfunction@$L$-function of a generic representation}
\index{generic representation!$L$-function of} of $\pi$. Define
\begin{equation}
\label{wdefeq}
w=\begin{bsmallmatrix}
&&1&\\
&&&-1\\
1&&&\\
&-1&&
\end{bsmallmatrix}.
\end{equation}
Then there exists $\gamma(s, \pi, \psi_{c_1,c_2}) \in \C(q^{-s})$ such that 
$$
Z(1-s,\pi(w) W) = \gamma(s, \pi, \psi_{c_1,c_2}) Z(s,W)
$$
for $W \in \mathcal{W}(\pi,\psi_{c_1,c_2})$. The factor $\gamma(s, \pi, \psi_{c_1,c_2})$
does not depend on the choice of $c_1,c_2\in \OF^\times$. Further define
$$
\varepsilon(s, \pi, \psi_{c_1,c_2})  =\gamma(s, \pi, \psi_{c_1,c_2}) \dfrac{L(s,\pi)}{L(1-s,\pi)}.
$$
Then there exists $\varepsilon \in \{\pm 1\}$ and 
an integer $N$ such that $\varepsilon(s, \pi, \psi_{c_1,c_2}) = \varepsilon q^{-N(s-\frac{1}{2})}$. 
We will write $\gamma(s,\pi)$ 
for $\gamma(s, \pi, \psi_{c_1,c_2})$ and $\varepsilon(s,\pi)$ 
for $\varepsilon(s, \pi, \psi_{c_1,c_2})$, it being understood
that we have fixed $\psi$ as on p.~\pageref{basefieldsubsecstar} and that 
$c_1, c_2 \in \OF^\times$. \label{gammafactordef} \label{epsilonfactordef}

\section{The paramodular theory}
\label{parasec}

In this section we recall some facts about paramodular vectors from \cite{NF}.
Let $n$ be a non-negative integer. We define 
\begin{equation}
\label{paradefeq}
\K{n} = \{ g \in \GSp(4,F)\mid\lambda(g) \in \OF^\times \} \cap
\begin{bsmallmatrix}
\OF&\OF&\OF&\p^{-n}\\
\p^n&\OF&\OF&\OF\\
\p^n&\OF&\OF&\OF\\
\p^n&\p^n&\p^n&\OF
\end{bsmallmatrix}.
\end{equation}
Using \eqref{ginveq}, it is easy to verify that $\K{n}$ is a compact subgroup of $\GSp(4,F)$. 
We refer to $\K{n}$ as the \emph{paramodular subgroup of level $\p^n$}.
\index{paramodular subgroup of level $\p^n$}
\index{GSp@$\GSp(4,F)$!paramodular group}
We have $\K{0}=\GSp(4,\OF)$. Let $(\pi,V)$ be a smooth representation of the group $\GSp(4,F)$ for 
which the center of $\GSp(4,F)$ acts trivially. We define 
$$
V(n)=\{ v \in V\mid \text{$\pi(k)v=v$ for all $k \in \K{n}$}\}. \label{nparadef}
$$
We refer to the non-zero elements of $V(n)$ as \emph{paramodular vectors}.\index{paramodular vectors}
Paramodular vectors have the following important property.

\begin{theorem}
\label{linindparatheorem}
Let $(\pi,V)$ be a smooth representation of $\GSp(4,F)$
for which the center acts trivially. Assume that the 
subspace of vectors fixed by $\SSp(4,F)$ is trivial. 
Then paramodular vectors from different levels are linearly
independent, i.e., 
$v_1\in V(n_1),\dots, v_t \in V(n_t)$
with $0 \leq n_1 < \cdots < n_t$ and $v_1 + \cdots + v_t=0$
implies that $v_1=\cdots =v_t =0$
\end{theorem}
\begin{proof}
This is Theorem 3.1.3 of \cite{NF}.
\end{proof}

We note that the assumption in Theorem \ref{linindparatheorem}
is satisfied if $\pi$ is an infinite-dimensional, irreducible,
admissible representation of $\GSp(4,F)$ with trivial central
character.

Let $(\pi,V)$ be a smooth representation of $\GSp(4,F)$
for which the center acts trivially. 
We define three operators 
$$
\eta,\theta, \theta': V \longrightarrow V \label{levelraisdef}
$$
by
\begin{align}
\eta v&= 
\pi(\begin{bsmallmatrix}
\varpi^{-1}&&&\\
&1&&\\
&&1&\\
&&&\varpi
\end{bsmallmatrix})v, \label{eta1tdefeq}\\
\theta v
& = \pi(
\begin{bsmallmatrix}
1&&&\\
&1&&\\
&&\varpi&\\
&&&\varpi
\end{bsmallmatrix} )v
+
q \int\limits_{\OF}
\pi(
\begin{bsmallmatrix}
1&&&\\
&1&c&\\
&&1&\\
&&&1
\end{bsmallmatrix}
\begin{bsmallmatrix}
1&&&\\
&\varpi&&\\
&&1&\\
&&&\varpi
\end{bsmallmatrix}
)v\, dc, \label{thetadefeq}\\
\theta'v
&=\eta v + q \int\limits_{\OF} \pi ( 
\begin{bsmallmatrix}
1&&& c \varpi^{-(n+1)}\\
&1&&\\
&&1&\\
&&&1
\end{bsmallmatrix})v\,dc \label{thetap1defeq}
\end{align}
for $v \in V$. If $v \in V(n)$, then $\eta v \in V(n+2)$ and $\theta v, \theta' v \in V(n+1)$. 
Further assume that $\pi$ is admissible and irreducible. If $V(n) \neq 0$ for some non-negative integer
$n$, then we say that $\pi$ is \emph{paramodular}, \index{paramodular representation}\index{representation!paramodular}
and we define $N_\pi$ 
\label{paraleveldef} to be the smallest non-negative integer $n$ such that $V(n) \neq 0$; 
we call $N_\pi$ the \emph{paramodular level}\index{paramodular level} of $\pi$. The 
following theorem is the main structural result about paramodular vectors.
See Table \ref{noklingentable} for the list of irreducible, admissible representations
of $\GSp(4,F)$ with trivial central character that are \emph{not} paramodular.
\begin{theorem}
\label{basicparatheorem}
Let $(\pi,V)$ be an irreducible, admissible representation of the group $\GSp(4,F)$
with trivial central character. Assume that $\pi$ is paramodular. Then $\dim V(N_\pi)=1$.
Let $v$ be a non-zero element of the one-dimensional space $V(N_\pi)$.  If $n \geq N_\pi$, then the space
$V(n)$ is spanned by the vectors
$$
\theta'^{i} \theta^j \eta^k v, \qquad i,j,k \geq 0, \quad i+j+2k=n-N_\pi.
$$
\end{theorem}
\begin{proof}
This is Theorem 7.5.1 and Theorem 7.5.7 of \cite{NF}.
\end{proof}

Any non-zero element of the one-dimensional
vector space $V(N_\pi)$ is called a \emph{paramodular newform}\index{paramodular newform}. 

\begin{figure}
\caption{Paramodular representations}
\label{pararepsfig}
\bigskip
\begin{tikzpicture}[scale=0.7]
\draw (0,0) rectangle (4,4);
\draw (0,3) -- (4,3);
\draw (0,2) -- (4,2);
\draw (2,0) -- (2,4);
\draw [decorate,decoration={brace,amplitude=5pt,raise=4pt},yshift=0pt]
(0,2) -- (0,4) node [black,midway,xshift=-1.2cm] {\scriptsize tempered};
\draw [decorate,decoration={brace,amplitude=5pt,raise=4pt},yshift=0pt]
(0,0) -- (0,2) node [black,midway,xshift=-1.4cm] {\scriptsize non-tempered};
\draw [decorate,decoration={brace,amplitude=5pt,raise=4pt},yshift=0pt]
(4,4) -- (4,3) node [black,midway,xshift=1.4cm] {\scriptsize supercuspidal};
\draw [decorate,decoration={brace,amplitude=5pt,raise=4pt},yshift=0pt]
(4,3) -- (4,0) node [black,midway,xshift=1.7cm] {\scriptsize non-supercuspidal};
\draw [decorate,decoration={brace,amplitude=5pt,raise=4pt},yshift=0pt]
(2,0) -- (0,0) node [black,midway,yshift=-0.5cm] {\scriptsize generic};
\draw [decorate,decoration={brace,amplitude=5pt,raise=4pt},yshift=0pt]
(4,0) -- (2,0) node [black,midway,yshift=-0.5cm] {\scriptsize non-generic};
\draw[fill=gray, fill opacity=0.3] (0,2) rectangle (2,4);
\draw[fill=gray, fill opacity=0.3] (2,0) rectangle (3,2);
\draw[fill=gray, fill opacity=0.3] (0,0) rectangle (2,2);
\draw[pattern=crosshatch dots] (2.4,0) rectangle (3,2);
\end{tikzpicture}\\
A schematic diagram of the paramodular representations among
 all irreducible, admissible representations of $\GSp(4,F)$ with trivial central
character. 
Paramodular representations are in gray, and  paramodular Saito-Kurokawa
representations are  dotted.
\end{figure}
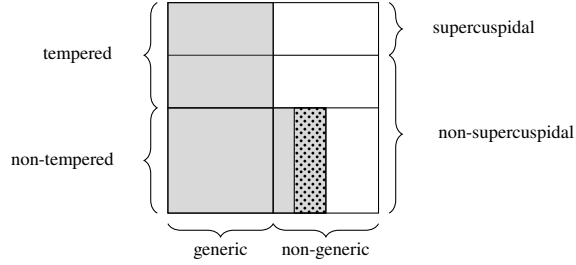

The next theorem shows that all generic representations are paramodular, and gives some sense
of the place of paramodular representations among all representations. For a schematic summary of
this result see Fig.~\ref{pararepsfig}. More generally, the work 
\cite{NF} determines exactly which representations $\pi$ of $\GSp(4,F)$ with trivial 
central character are paramodular and calculates $N_\pi$
for all paramodular representations.  

\begin{theorem}
\label{gentemptheorem}
Let $(\pi,V)$ be an irreducible, admissible representation of the group $\GSp(4,F)$ with trivial
central character. If $\pi$ is generic, then $\pi$ is paramodular. If $\pi$ is tempered, then 
 $\pi$ is paramodular if and only if $\pi$ is generic. 
\end{theorem}
\begin{proof}
This follows from Theorem 7.5.4 and Theorem 7.5.8 of \cite{NF}.
\end{proof}

As regards non-generic representations, ``most'' of the 
non-generic, irreducible, admissible representations of the group $\GSp(4,F)$ with trivial
central character that are paramodular are Saito-Kurokawa representations (see Sect.~\ref{repsec}).
If $Q(\nu^{\frac{1}{2}} \pi, \nu^{-\frac{1}{2}}\sigma)$ 
is a  Saito-Kurokawa representation with $\sigma^2=1$,
then $Q(\nu^{\frac{1}{2}} \pi, \nu^{-\frac{1}{2}}\sigma)$ is paramodular
if and only if $\sigma$ is unramified.

Again let $(\pi,V)$ be a smooth representation of $\GSp(4,F)$ for which
the center of $\GSp(4,F)$ acts trivially. Let $n$ be a non-negative integer. 
Let 
\begin{align*}
\K{n}
\begin{bsmallmatrix}
\varpi&&&\\
&\varpi&&\\
&&1&\\
&&&1
\end{bsmallmatrix}
\K{n}
&=\bigsqcup_{i \in I} g_i \K{n}, \\
\K{n}
\begin{bsmallmatrix}
\varpi^2&&&\\
&\varpi&&\\
&&\varpi&\\
&&&1
\end{bsmallmatrix}
\K{n}
&=\bigsqcup_{j \in J} h_j \K{n}
\end{align*}
be disjoint decompositions. We define endomorphisms
\begin{equation}
\label{Heckeopdefeq}
T_{0,1}: V(n) \longrightarrow V(n) \quad \text{and} \quad T_{1,0}: V(n) \longrightarrow V(n)
\end{equation}
by 
$$
T_{0,1} v = \sum_{i \in I} \pi (g_i) v \quad \text{and}\quad  T_{1,0} v = \sum_{j \in J} \pi (h_j) v
$$
for $v \in V(n)$. We refer to $T_{0,1}$ and $T_{1,0}$ as the \emph{paramodular Hecke operators of level $n$}.\index{paramodular Hecke operators!local}
\index{Hecke operator!paramodular}
Also, let 
\begin{equation}
\label{ALeq}
u_n =
\begin{bsmallmatrix}
&&1&\\
&&&-1\\
\varpi^n&&&\\
&-\varpi^n&&&
\end{bsmallmatrix}.
\end{equation}
The element $u_n$ normalizes $\K{n}$, and 
the operator $\pi(u_n)$ maps the space $V(n)$ into itself; we call $\pi(u_n)$ the \emph{paramodular Atkin-Lehner operator of level $n$}.
\index{paramodular Atkin-Lehner operator!local}\index{Atkin-Lehner operator, paramodular!local}
If $\pi$ is a paramodular, irreducible, admissible representation of $\GSp(4,F)$
with trivial central character, then any non-zero element $v_\mathrm{new}$ of the one-dimensional
vector space $V(N_\pi)$ is an eigenvector for $\pi(u_{N_\pi})$, $T_{0,1}$, and $T_{1,0}$;
these eigenvalues will be denoted by $\varepsilon_\pi$, $\lambda_\pi$, and $\mu_\pi$,
respectively, so that
\begin{equation}
\label{epsilonlambdamudefeq}
\pi(u_{N_\pi})v_{\mathrm{new}} = \varepsilon_\pi v_{\mathrm{new}},\quad
T_{0,1} v_{\mathrm{new}} = \lambda_\pi v_{\mathrm{new}},\quad
T_{1,0} v_{\mathrm{new}} = \mu_\pi v_{\mathrm{new}}.
\end{equation}
The eigenvalues $\varepsilon_\pi$, $\lambda_\pi$, and $\mu_\pi$ were determined 
for all paramodular, irreducible, admissible representations $\pi$ of $\GSp(4,F)$
with trivial central character in \cite{NF}; for the convenience of the reader, this data
is reproduced in Table \ref{levelsandeigenvaluestable}.
The meaning of these eigenvalues for generic representations is given by the following theorem.
\begin{theorem}
\label{basicgenparatheorem}
Let $(\pi,V)$ be a generic, irreducible, admissible
representation of $\GSp(4,F)$ with trivial central character; 
$\pi$ is paramodular by Theorem \ref{gentemptheorem}.
 Let $W \in V(N_\pi)$ be a newform. Let $\varepsilon_\pi$
be the Atkin-Lehner eigenvalue of $W$, i.e., $\pi(u_n) W = \varepsilon_\pi W$, and
let $\lambda_\pi$ and $\mu_\pi$ be the Hecke eigenvalues of $W$,
defined by $T_{0,1} W = \lambda_\pi W$ and $T_{1,0} W = \mu_\pi W$. We may choose
$W$ such that
$$
Z(s,W) = L(s,\pi),
$$
and we have
$$
\varepsilon(s,\pi) =\varepsilon_\pi q^{-N_\pi (s-1/2)}.
$$

Here, the zeta integral $Z(s,W)$, $L(s,\pi)$ and $\varepsilon(s,\pi)$  are defined as on p.~\pageref{zetaintegralssubsecstar}.
Moreover:
\begin{enumerate}
\item \label{basicgenparatheoremitem1} If $N_\pi=0$, so that $\pi$ is unramified, then
$$
L(s,\pi)
=
\dfrac{1}{1-q^{-\frac{3}{2}}\lambda_\pi q^{-s} 
+(q^{-2} \mu_\pi +1+q^{-2}) q^{-2s} -q^{-\frac{3}{2}}\lambda_\pi q^{-3s} +q^{-4s}}.
$$
\item \label{basicgenparatheoremitem2} If $N_\pi=1$, then,
$$
L(s,\pi)
=
\dfrac{1}{ 1-q^{-\frac{3}{2}}(\lambda_\pi+\varepsilon_\pi)q^{-s}+(q^{-2}\mu_\pi+1)q^{-2s} +\varepsilon_\pi q^{-\frac{1}{2}}q^{-3s}}.
$$
\item \label{basicgenparatheoremitem3} If $N_\pi \geq 1$, then 
$$
L(s,\pi)
=
\dfrac{1}{1-q^{-\frac{3}{2}}\lambda_\pi q^{-s} + (q^{-2} \mu_\pi +1) q^{-2s}}.
$$
\end{enumerate}
\end{theorem}
\begin{proof}
See Theorem 7.5.4, Corollary 7.5.5, and Theorem 7.5.3 of \cite{NF}. 
\end{proof}
A similar result holds for all paramodular representations.
\begin{theorem}
\label{basictwoparatheorem}
Let $(\pi,V)$ be a paramodular, irreducible, admissible
representation of $\GSp(4,F)$ with trivial central character. 
Let $\varphi_\pi: W_F' \to \GSp(4,\C)$ be the $L$-parameter
of $\pi$. Let $v \in V(N_\pi)$ be a non-zero vector. Let $\varepsilon_\pi$
be the Atkin-Lehner eigenvalue of $v$, i.e., $\pi(u_n) v = \varepsilon_\pi v$, and
let $\lambda_\pi$ and $\mu_\pi$ be the Hecke eigenvalues of~$v$,
defined by $T_{0,1} v = \lambda_\pi v$ and $T_{1,0} v = \mu_\pi v$. Then  
$$
\varepsilon(s,\varphi_\pi) = \varepsilon_\pi q^{-N_\pi (s-1/2)}.
$$
Moreover, $L(s,\varphi_\pi)$ is given by exactly the same formulas
as in (1), (2), and (3) of Theorem \ref{basicgenparatheorem}.
\end{theorem}
\begin{proof}
The work \cite{GT} assigns an $L$-parameter $\varphi_\pi$ to  $\pi$. As remarked in \cite{GT},
the $L$-parameter assigned to $\pi$ by \cite{GT} coincides with
the $L$-parameter assigned to $\pi$ in \cite{NF} when $\pi$
is non-supercuspidal; in this case, the $L$-parameter of $\pi$ is 
determined by the desiderata of the local Langlands conjecture. The 
statement of the theorem when $\pi$ is non-supercuspidal is Theorem 7.5.9
of \cite{NF}. Assume that $\pi$ is supercuspidal. By Theorem \ref{gentemptheorem},
$\pi$ is generic. By \cite{Soudry1993}, \cite{JiangSoudry2004} and \cite{GT} one has $L(s,\varphi_\pi)=L(s,\pi)$
and $\varepsilon(s,\varphi_\pi)=\varepsilon(s,\pi)$, where $L(s,\pi)$
and $\varepsilon(s,\pi)$ are as on p.~\pageref{zetaintegralssubsecstar}. 
The theorem follows now from Theorem \ref{basicgenparatheorem}.
\end{proof} 

In Chap.~\ref{backSMFchap} we will use the following result about oldforms 
in smooth representations. Let 
$(\pi,V)$ be a smooth representation of $\GSp(4,F)$ for which the center acts trivially, 
and let $n$ be a non-negative integer. We let $V(n)_{\mathrm{old}}$ be the subspace
of $V(n)$ spanned by the vectors of the form $\theta v$ and $\theta' v$ for $v \in V(n-1)$ 
 and the vectors of the form $\eta v$ for $v \in V(n-2)$ (for this definition we take $V(k)=0$
 if $k<0$). We note that the next lema follows trivially from Theorem~\ref{basicparatheorem} in case $\pi$ is an irreducible, admissible representation.

\begin{lemma}\label{localoldpreservinglemma}
 Let $(\pi,V)$ be a smooth representation of $\GSp(4,F)$ for which the center acts trivially, 
and let $n$ be a non-negative integer. The endomorphisms $T_{0,1},T_{1,0},$ and $ \pi(u_n)$ of $V(n)$
map $V(n)_{\mathrm{old}}$ into $V(n)_{\mathrm{old}}$. 
\end{lemma}
\begin{proof}
The statement is trivial if $n=0$ because $V(0)_{\mathrm{old}}=0$. Assume that $n \geq 1$. 
Let $w \in V(n)_{\mathrm{old}}$. We need to prove that $T_{0,1}w, T_{1,0}w, \pi(u_n)w \in V(n)_{\mathrm{old}}$. 
One can prove that
$$
T_{0,1} \eta v = q \theta' \theta v,\qquad
T_{1,0} \eta v = q^2 \theta'^2 v -q^2 (q+1)\eta v
$$
for $v \in V(n-1)$, and
$$
\pi(u_n) \eta v = \eta \pi (u_{n-2}) v
$$
for $v \in V(n-2)$ (again, if $n-2<0$ we take $V(n-2)=0$). It follows that we may assume that $w$ is of the form $\theta v$ or $\theta' v$
for some $v \in V(n-1)$. 
If $n \geq 2$, then $T_{0,1}w, T_{1,0}w, \pi(u_n)w \in V(n)_{\mathrm{old}}$ by Proposition 6.3 of \cite{NF}. 
If $n=1$, then $n-1=0$. One may prove that if $v \in V(0)$ then
\begin{align*}
 T_{0,1} \theta  v&= \theta  T_{0,1} v+(q^2-1)\theta' v ,\\
T_{0,1}\theta' v &= \theta' T_{0,1} v+(q^2-1)\theta v ,\\
T_{1,0}\theta v&=q \theta' T_{0,1} v -q(q+1) \theta v,\\
T_{1,0} \theta' v &=   \theta' T_{1,0} v+(q^3+1)\theta' v -\theta T_{0,1} v.
\end{align*}
This proves the lemma.
\end{proof}

%% file: SKMS_chapter3.tex
\chapter{Stable Klingen Vectors}
\label{basicfactschapter}

Let $(\pi,V)$ be a smooth representation of $\GSp(4,F)$ with trivial central character.
In this chapter, for non-negative integers $n$, we introduce the stable Klingen
congruence subgroups $\Ks{n}$ of $\GSp(4,F)$ and begin the consideration
of the subspaces $V_s(n)$ of vectors in $V$ fixed by $\Ks{n}$. We refer
to such vectors as stable Klingen vectors. Since $\Ks{n}$ is contained
in the paramodular group $\K{n}$, every paramodular vector is a stable Klingen vector,
and we will use this connection throughout this work. In this chapter
we will also introduce and study various operators between the spaces $V_s(n)$ and $V_s(m)$
for non-negative integers $n$ and $m$. Similar level changing operators are important in the 
paramodular theory. However, in contrast to the paramodular theory, we will find
that level changing operators in the stable Klingen theory often admit upper block
formulas. As we will see in the second part of this work, such formulas are useful for effective calculations and applications.
The results of this chapter are algebraic in the sense that we only assume
that $\pi$ is smooth; subsequent chapters will assume that $\pi$ is admissible and
irreducible.

\section{The stable Klingen subgroup}
In this section we define the stable Klingen subgroups of $\GSp(4,F)$ and describe some of their
basic properties. We will  make reference to the paramodular theory; see Sect.~\ref{parasec} for a summary. 
Let $n$ be an integer such that $n\geq 0$. We define  
\begin{equation}\label{Klndefeq}
 \Kl{n}=\GSp(4,\OF)\cap\begin{bsmallmatrix}\OF&\OF&\OF&\OF\\\p^n&\OF&\OF&\OF\\\p^n&\OF&\OF&\OF\\\p^n&\p^n&\p^n&\OF\end{bsmallmatrix}.
\end{equation}
Using \eqref{ginveq}, it is easy to see that $\Kl{n}$ is a subgroup of $\GSp(4,\OF)$. We
refer to $\Kl{n}$ as 
the \emph{Klingen congruence subgroup
of level $\p^n$}\index{Klingen congruence subgroup}.
Evidently, $\Kl{n}$ is contained in the paramodular subgroup $\K{n}$. Also, $\Kl{0}=\GSp(4,\OF)$. Next,
we let $\Ks{n}$
to be the subgroup of $\GSp(4,F)$ generated by $\Kl{n}$ and
the subgroup 
$$
\begin{bsmallmatrix} 1&&&\p^{-n+1}\\&1&&\\&&1&\\&&&1 \end{bsmallmatrix}.
$$
We refer to $\Ks{n}$ as the \emph{stable Klingen subgroup of level $\p^n$}\index{stable Klingen subgroup}. Thus,
\begin{equation}\label{Ksndefeq}
 \Ks{n}=\langle\Kl{n},\begin{bsmallmatrix}1&&&\p^{-n+1}\\&1\\&&1\\&&&1\end{bsmallmatrix}\rangle.
\end{equation}
We note that $\Ks{0}=\GSp(4,\OF)$ and $\Ks{1}=\Kl{}$.  It is useful to note that the element
$s_2$ from \eqref{s1s2defeq} is contained in $\Kl{n}$, and hence in $\Ks{n}$ and $\K{n}$, for all
integers $n$ such that $n \geq 0$. 
\begin{lemma}
\label{Ksshapelemma}
Let $n$ be an integer such that $n \geq 1$. Then
\begin{equation}\label{Ksshapeeq}
\Ks{n}=\{g\in\GSp(4,F)\mid\lambda(g)\in\OF^\times\}\cap
\begin{bsmallmatrix}
\OF&\OF&\OF&\p^{-n+1}\\
\p^n&\OF&\OF&\OF\\
\p^n&\OF&\OF&\OF\\
\p^n&\p^n&\p^n&\OF
\end{bsmallmatrix}.
\end{equation}
Moreover, if $g=(g_{ij}) \in \Ks{n}$, then 
\begin{equation}\label{1144unitseq}
g_{11},g_{44} \in \OF^\times.
\end{equation}
\end{lemma}
\begin{proof}
It is obvious that the left-hand side is contained in the right-hand side. Conversely, let $g=(g_{ij})$ be an element of the right-hand side. Then $\det(g)$ is in $\OF^\times$. This implies that $g_{11}$ and $g_{44}$ are also in $\OF^\times$. Consequently,
\begin{equation}
\label{upgeq}
 \begin{bsmallmatrix}1&&&-g_{14}g_{44}^{-1}\\&1\\&&1\\&&&1\end{bsmallmatrix}g\in\Kl{n}.
\end{equation}
This implies that $g\in\Ks{n}$.
\end{proof}

Let $n$ be an integer such that $n\geq 1$. Then it is possible to conjugate $\Ks{n}$ into $\GSp(4,\OF)$. If $n$ is even, then
\begin{equation}\label{Ksconjeq1}
\begin{bsmallmatrix}\varpi^{\frac{n}{2}}\\&1\\&&1\\&&&\varpi^{-\frac{n}{2}}\end{bsmallmatrix}\Ks{n}
\begin{bsmallmatrix}\varpi^{-\frac{n}{2}}\\&1\\&&1\\&&&\varpi^{\frac{n}{2}}\end{bsmallmatrix}\subset
\begin{bsmallmatrix}\OF&\p^{\frac{n}{2}}&\p^{\frac{n}{2}}&\p\\\p^{\frac{n}{2}}&\OF&\OF&\p^{\frac{n}{2}}\\\p^{\frac{n}{2}}&\OF&\OF&\p^{\frac{n}{2}}\\\OF&\p^{\frac{n}{2}}&\p^{\frac{n}{2}}&\OF\end{bsmallmatrix},
\end{equation}
and if $n$ is odd, then
\begin{equation}\label{Ksconjeq2}
 \begin{bsmallmatrix}\varpi^{\frac{n-1}2}\\&1\\&&1\\&&&\varpi^{\frac{1-n}2}\end{bsmallmatrix}\Ks{n}\begin{bsmallmatrix}\varpi^{\frac{1-n}2}\\&1\\&&1\\&&&\varpi^{\frac{n-1}2}\end{bsmallmatrix}\subset\begin{bsmallmatrix}\OF&\p^{\frac{n-1}2}&\p^{\frac{n-1}2}&\OF\\\p^{\frac{n+1}2}&\OF&\OF&\p^{\frac{n-1}2}\\\p^{\frac{n+1}2}&\OF&\OF&\p^{\frac{n-1}2}\\\p&\p^{\frac{n+1}2}&\p^{\frac{n+1}2}&\OF\end{bsmallmatrix}.
\end{equation}
For the paramodular group, such a conjugation is only possible if $n$ is even. 

In the next lemma we describe the coset decompositions of $\Kl{n}$ in $\Ks{n}$ and of $\Ks{n}$ in $\K{n}$. 
Let $n$ be an integer such that $n\geq 0$. We define 
\begin{equation}
\label{tndefeq}
t_n=
\begin{bsmallmatrix}
&&&-\varpi^{-n}\\
&1&&\\
&&1&\\
\varpi^n
\end{bsmallmatrix}.
\end{equation}
Clearly, $t_n$ is contained in the paramodular group $\K{n}$ of level $\p^n$. 

\begin{lemma}
\label{Ksdecomplemma}
Let $n$ be  an integer such that $n\geq 1$. Then
\begin{equation}\label{KsKlcosetseq}
 \Ks{n}=\bigsqcup\limits_{x\in\OF/\p^{n-1}}\begin{bsmallmatrix}1&&&x\varpi^{-n+1}\\&1\\&&1\\&&&1\end{bsmallmatrix}\Kl{n}
\end{equation}
and
\begin{align}\label{paramodularKsdecompositioneq2}
\K{n}=\bigsqcup\limits_{x\in\OF/\p}
\begin{bsmallmatrix}1&&&x\varpi^{-n}\\&1\\&&1\\&&&1\end{bsmallmatrix}\Ks{n}\sqcup t_n\Ks{n}.
\end{align}
Moreover,
\begin{equation}\label{KsnKn1Kneq}
 \Ks{n}=\K{n-1}\cap\K{n}.
\end{equation}
\end{lemma}
\begin{proof}
It is clear that the right-hand side of \eqref{KsKlcosetseq} is contained in the left-hand side. The opposite
inclusion follows from the proof of Lemma \ref{Ksshapelemma}. It is straightforward to verify
that the decomposition is disjoint. 
Next, by Lemma 3.3.1 of \cite{NF},
\begin{align}\label{paramodularKlingendecompositioneq2}
  \K{n}=\bigsqcup\limits_{x\in\OF/\p^n}
  \begin{bsmallmatrix}1&&&x\varpi^{-n}\\&1\\&&1\\&&&1\end{bsmallmatrix}\Kl{n}\sqcup\bigsqcup\limits_{x\in\OF/\p^{n-1}}
  t_n\begin{bsmallmatrix}1&&&x\varpi^{-n+1}\\&1\\&&1\\&&&1\end{bsmallmatrix}\Kl{n}.
\end{align}
Multiplying \eqref{paramodularKlingendecompositioneq2} from the right by $\Ks{n}$, we get
\begin{align}\label{paramodularKsdecompositioneq1}
  \K{n}=\bigcup\limits_{x\in\OF/\p}
  \begin{bsmallmatrix}1&&&x\varpi^{-n}\\&1\\&&1\\&&&1\end{bsmallmatrix}\Ks{n}\cup
  t_n\Ks{n}.
\end{align}
It is easy to see that this decomposition is disjoint. To prove \eqref{KsnKn1Kneq}, we
note first that $\Ks{n}$ is contained in $\K{n-1} \cap \K{n}$. Conversely, let $g \in \K{n-1}\cap \K{n}$.
Then $g$ is in the right-hand side of \eqref{Ksshapeeq}. By \eqref{Ksshapeeq}, $g \in \Ks{n}$.
\end{proof}

Since the Klingen subgroup has an Iwahori decomposition, it follows from \eqref{KsKlcosetseq} that the stable Klingen subgroup has as well. Let $n$ be an integer such that $n\geq 1$. Then $\Ks{n}$ is equal to 
\begin{equation}\label{Ksniwahorieq1}
(\begin{bsmallmatrix}
 1&\OF&\OF&\p^{-n+1}\vphantom{\varpi^{-n+1}}\\
 &1&&\OF\vphantom{\varpi^{-n+1}}\\
 &&1&\OF\vphantom{\varpi^{-n+1}}\\
 &&&1\vphantom{\varpi^{-n+1}}
 \end{bsmallmatrix} \cap \Ks{n})
 (\begin{bsmallmatrix}
 \OF^\times&&&\vphantom{\varpi^{-n+1}}\\
 &\OF&\OF\\&\OF&\OF\vphantom{\varpi^{-n+1}}\\
 &&&\OF^\times\vphantom{\varpi^{-n+1}}
 \end{bsmallmatrix}\cap \Ks{n})
(\begin{bsmallmatrix}1\\
\p^n&1\\\p^n&&1\\
\p^n&\p^n&\p^n&1
\end{bsmallmatrix}\cap \Ks{n})
\end{equation}
and is also equal to
\begin{equation}\label{Ksniwahorieq2}
(\begin{bsmallmatrix}1\\
\p^n&1\\\p^n&&1\\
\p^n&\p^n&\p^n&1
\end{bsmallmatrix} \cap \Ks{n})
(\begin{bsmallmatrix}
\OF^\times&&&\vphantom{\varpi^{-n+1}}\\
&\OF&\OF\vphantom{\varpi^{-n+1}}\\
&\OF&\OF\vphantom{\varpi^{-n+1}}\\
&&&\OF^\times\vphantom{\varpi^{-n+1}}
\end{bsmallmatrix} \cap \Ks{n})
(\begin{bsmallmatrix}
1&\OF&\OF&\p^{-n+1}\vphantom{\varpi^{-n+1}}\\
&1&&\OF\vphantom{\varpi^{-n+1}}\\
&&1&\OF\vphantom{\varpi^{-n+1}}\\
&&&1\end{bsmallmatrix} \cap \Ks{n}).
\end{equation}
\section{Stable Klingen vectors}
\label{stableklingenvecsec}
We now introduce the spaces of vectors which play a central role in this work. 
Let $(\pi,V)$ be a smooth representation of $\GSp(4,F)$ for which the center acts trivially,
and let $n$ be an integer such that $n \geq 0$. 
Let $\Ks{n}$ be the stable Klingen subgroup defined in \eqref{Ksndefeq}, and let
\begin{equation}\label{Vsndefeq}
 V_s(n)=\{v\in V\:|\:\pi(g)v=v\text{ for all }g\in\Ks{n}\}
\end{equation}
be the corresponding space of vectors in $V$ fixed by $\Ks{n}$. We refer to the elements of $V_s(n)$
as \emph{stable Klingen vectors of level $\p^n$}.\index{stable Klingen vector}
If  $n\geq 1$, then since $\Ks{n} = \K{n-1} \cap \K{n}$ by Lemma \ref{Ksdecomplemma},
we have $V(n-1)+V(n) \subset V_s(n)$. In this work we define $V(-1)=0$; with this definition, 
we also have $V(n-1)+V(n) \subset V_s(n)$ if $n=0$. We refer to $V(n-1)+V(n)$ as 
the \emph{paramodular subspace}
\index{paramodular subspace}
of $V_s(n)$. We define
\begin{equation}\label{Vsbardefeq}
 \bar V_s(n)=V_s(n)/(V(n-1)+ V(n)).
\end{equation}
Since $V_s(0)=V(0)$, we have $\bar V_s(0)=0$. We refer
to the elements of $\bar V_s(n)$ as \emph{quotient stable Klingen vectors}.
\index{quotient stable Klingen vectors}
If $V_s(n) \neq 0$ for some integer such that $n\geq 0$, then we define $N_{\pi,s}$ 
to be the smallest such integer, and we call $N_{\pi,s}$ the 
\emph{stable Klingen level} \label{stableKlingenleveldef} of $\pi$.\index{stable Klingen level}
If there exists no integer $n\geq 0$ such that $V_s(n) \neq 0$, then we will 
say that $N_{\pi,s}$ is not defined. If $\bar V_s(n) \neq 0$ for some integer such that $n \geq 0$, then
we define $\bar N_{\pi,s}$ to be the smallest such integer, and we call $\bar N_{\pi,s}$
the \emph{quotient stable Klingen level} of $\pi$;
\index{quotient stable Klingen level}
if there exists no integer $n\geq 0$ such that $\bar V_s(n) \neq 0$, then we
will say that $\bar N_{\pi,s}$ is not defined.
As a point of orientation, we note that it is trivial that if $\pi$ admits non-zero
paramodular vectors, then $\pi$ admits non-zero stable Klingen vectors; in Corollary~\ref{upperboundpropcor}
we will see that the converse also holds when $\pi$ is irreducible.
\section{Paramodularization}
\label{paramodularizationsec}
  
Let $(\pi,V)$ be a smooth representation of $\GSp(4,F)$ for which the center acts trivially.
Let $n$ be an integer such that $n\geq 0$, and let $dk$ be a Haar measure on $\GSp(4,F)$.
We define the \emph{paramodularization map}\index{paramodularization map}
$p_n: V_s(n) \to V(n)$ by 
\begin{equation}\label{pndefeq}
 p_nv=\frac1{{\rm vol}(\K{n})}\int\limits_{\K{n}}\pi(k)v\,dk,\qquad \text{for $v\in V_s(n)$}.
\end{equation}
\begin{lemma}
\label{paraformlemma}
Let $(\pi,V)$ be a smooth representation of $\GSp(4,F)$ for which the center acts trivially.
Let $n$ be an integer such that $n\geq 0$. Then 
\begin{align}\label{pnVsneq}
 p_nv&=\frac1{q+1}\bigg(\pi(t_n)v
 +\sum_{x\in\OF/\p}\pi(\begin{bsmallmatrix}1&&&x\varpi^{-n}\\&1\\&&1\\&&&1\end{bsmallmatrix})v\bigg)\\
 &=\frac1{q+1}\bigg(\pi(t_n)v
 +q\int\limits_{\OF} \pi(\begin{bsmallmatrix}1&&&x\varpi^{-n}\\&1\\&&1\\&&&1\end{bsmallmatrix})v\, dx\bigg)
\label{pnVsneq3}
\end{align}
for $v \in V_s(n)$. The map $p_n:V_s(n) \to V(n)$ is a linear projection, so that $p_n^2 =p_n$ and
\begin{equation}
\label{pnVsneq4}
V_s(n) = V(n) \oplus \ker(p_n).
\end{equation}

\end{lemma}
\begin{proof}
The assertion \eqref{pnVsneq} is clear if $n=0$. If $n$ is an integer such that $n\geq 1$, then \eqref{pnVsneq}
follows from \eqref{paramodularKsdecompositioneq2}. It is clear that $p_n$ is a linear
projection. 
\end{proof}

\section{Operators on stable Klingen vectors}
\label{opersksec}

In this work we will use a variety of linear operators to investigate
stable Klingen vectors. These linear maps are instances of a general definition.
To describe this concept, let 
$(\pi,V)$ be a smooth representation of $\GSp(4,F)$ 
for which the center acts trivially, let
$n$ and $m$ be integers such that $n,m \geq 0$,
and let $g \in \GSp(4,F)$. We define a linear map 
$$
T_g: V_s(n) \longrightarrow V_s(m)
$$
by
\begin{equation}
\label{Tgformeq}
T_g v = \frac{1}{\vl(\Ks{m} \cap g \Ks{n} g^{-1})} \int\limits_{K_s(\p^m)} \pi(k) \pi(g) v\, dk
\end{equation}
for $v \in V_s(n)$. Here, $dk$ is a fixed Haar measure on $\GSp(4,F)$; we note that the 
definition of $T_g$ does not depend on the choice of Haar measure. The operator $T_g$
has a simpler expression. Let $J=\Ks{m} \cap g\Ks{n}g^{-1}$, and let
\begin{equation}
\label{Tgdecompeq}
\Ks{m} = \bigsqcup_{i \in I} g_i J
\end{equation}
be a disjoint decomposition, where $g_i \in \Ks{m}$ for $i \in I$. The set $I$ is finite since $J$
is open and $\Ks{m}$ is compact. A calculation shows that 
\begin{equation}\label{Tgvsumeq}
T_g v = \sum_{i \in I} \pi(g_i) \pi(g) v
\end{equation}
for $v \in V_s(n)$. We will say that $T_g:V_s(n) \to V_s(m)$ is a 
\emph{level raising operator}\index{level raising operator}
if $n<m$, that $T_g$ is a 
\emph{level lowering operator}\index{level lowering operator}
if $m<n$, and that $T_g$ is a 
\emph{Hecke operator}\index{Hecke operator}
if $n=m$. 
We will also say that $T_g$ is an 
\emph{upper block operator}\index{upper block operator}
if there exists a disjoint decomposition \eqref{Tgdecompeq} such 
that $g_ig$ lie in the Siegel parabolic subgroup $P$ for all $i \in I$. 
In the following sections we will define some important
examples of level raising, level lowering, and Hecke operators and investigate
the basic properties of these linear maps. 

We close this section with a lemma showing that the above
definition of a Hecke operator agrees with the traditional definition.

\begin{lemma}
\label{heckedefeqlemma}
Let $(\pi,V)$ be a smooth representation
of $\GSp(4,F)$ for which the center acts trivially,
and let $n$ be an integer such that $n \geq 0$. Let
$g \in \GSp(4,F)$, and let
$$
\Ks{n} g \Ks{n} = \bigsqcup_{i \in I} g_i \Ks{n}
$$
be a disjoint decomposition. Then 
$$
T_g v = \sum_{i \in I} \pi(g_i)v
$$
for $v \in V_s(n)$.
\end{lemma}
\begin{proof}
For each $i \in I$, let $k_i,k_i' \in \Ks{n}$ be such that $g_i = k_i g k_i'$.
Define $g_i' = g_ik_i'^{-1}$ for $i \in I$. Then $g_i'=k_ig$ for $i \in I$, 
$\Ks{n} g\Ks{n} = \sqcup_{i\in I}g_i' \Ks{n}$, and $\sum_{i \in I} \pi(g_i')v =
\sum_{i \in I} \pi(g_i)v$ for $v \in V_s(n)$. An argument shows that
$\Ks{n} = \sqcup_{i \in I} g_i'g^{-1} J$, where $J = \Ks{n} \cap g\Ks{n}g^{-1}$.
Hence, if $v \in V_s(n)$, then $T_g v = \sum_{i\in I} \pi(g_i'g^{-1}) \pi(g) v
=\sum_{i\in I} \pi(g_i') v = \sum_{i \in I} \pi(g_i) v$. 
\end{proof}

\section{Level raising operators}
\label{levelraisingsec}
Let $(\pi,V)$ be a smooth representation of $\GSp(4,F)$ for which the center acts trivially, and let $n$ be an integer such that $n\geq 0$. 
In this section we will consider the three level raising operators
$T_g:V_s(n) \to V_s(n+1)$ determined by $g=t_n$, $g=1$, and 
\begin{equation}
\label{01elementeq}
g = 
\begin{bsmallmatrix}
1&&&\\
&1&&\\
&&\varpi&\\
&&&\varpi
\end{bsmallmatrix}.
\end{equation}
The level raising operator determined by $t_n$ is particularly simple.

\begin{lemma}
\label{tncosetlemma}
Let $(\pi,V)$ be a smooth representation of $\GSp(4,F)$ for which the center acts trivially, and let $n$ be an integer such that $n\geq 0$. We have 
\begin{equation}
\label{tnKsneq}
\Ks{n+1} \cap t_n \Ks{n} t_n^{-1} = \Ks{n+1}.
\end{equation}
Thus, the level raising operator $T_g:V_s(n) \to V_s(n+1)$ determined by $g = t_n$ is given by 
$T_g(v) = \pi(t_n) v$ for $v \in V_s(n)$.
\end{lemma}
\begin{proof}
The equation \eqref{tnKsneq} is equivalent to 
$\Ks{n+1} \subset t_n \Ks{n} t_n^{-1}$, which is equivalent to
$t_n^{-1} \Ks{n+1} t_n \subset \Ks{n}$. 
This is clear if $n=0$. Assume that $n>0$. Then $t_n^{-1} \Ks{n+1} t_n \subset \Ks{n}$
 follows by a calculation using
the criterion \eqref{Ksshapeeq}.
\end{proof}

With the notation as in Lemma \ref{tncosetlemma}, this lemma asserts that the level raising operator $T_{t_n}:V_s(n) \to V_s(n+1)$
is given by applying $\pi(t_n)$; consequently, we will write 
   \begin{equation}\label{tnlevelraisinglemmaeq1}
    t_n:\:V_s(n)\longrightarrow V_s(n+1)
   \end{equation}
for this level raising operator. Using, for example, \eqref{1144unitseq}, it can be
verified that $t_n:V_s(n) \to V_s(n+1)$ is not an upper block operator. Further
properties of the operator $t_n$ are given in the following lemma.

\begin{lemma}\label{tnlevelraisinglemma}
Let $(\pi,V)$ be a smooth representation of $\GSp(4,F)$ for which the center acts trivially, and let $n$ be an integer such that $n\geq 0$.
 \begin{enumerate}
  \item\label{tnlevelraisinglemmaitem1}  The level raising operator $t_n: V_s(n) \to V_s(n+1)$ is injective.
  \item\label{tnlevelraisinglemmaitem2}  The level raising operator $t_n: V_s(n) \to V_s(n+1)$ induces a linear map
   \begin{equation}\label{tnlevelraisinglemmaeq2}
    t_n:\:\bar V_s(n)\longrightarrow\bar V_s(n+1)
   \end{equation}
   which is also injective.
  \item\label{tnlevelraisinglemmaitem3}  The map $t_{n+1}\circ t_n:V_s(n)\to V_s(n+2)$ is given by applying 
\begin{equation}
\label{etadefeq}
\eta=\pi(\begin{bsmallmatrix} \varpi^{-1} &&& \\ &1&& \\ &&1& \\ &&&\varpi \end{bsmallmatrix}).
\end{equation}
 \end{enumerate}
\end{lemma}
\begin{proof}
\ref{tnlevelraisinglemmaitem1}. This is obvious. 

\ref{tnlevelraisinglemmaitem2}. Because the element $t_n$ is in $\K{n}$, we have $t_nV(n) =V(n)$. This proves \ref{tnlevelraisinglemmaitem2} in 
the case $n=0$; assume that $n$ is an integer such that $n\geq 1$. We have $t_n V(n-1) = \eta t_{n-1} V(n-1) = \eta V(n-1)
\subset V(n+1)$ (see Sect.~\ref{parasec}). 
It follows that $t_n$ maps $V(n-1)+V(n)$ into $V(n)+V(n+1)$, and 
therefore induces a map $\bar V_s(n)\longrightarrow\bar V_s(n+1)$. To prove injectivity, 
assume that $v\in V_s(n)$, and that $t_nv$ lies in the subspace $V(n)+V(n+1)$ of $V_s(n+1)$. We will show 
that $v$ lies in the subspace $V(n-1)+V(n)$ of $V_s(n)$, thus proving the injectivity of \eqref{tnlevelraisinglemmaeq2}.
Write $t_nv=v_1+v_2$ with $v_1\in V(n)$ and $v_2\in V(n+1)$. We have
\begin{equation}\label{tnlevelraisinglemmaeq3}
 v_2\text{ is invariant under }\begin{bsmallmatrix}1\\\p^n&1\\\p^n&&1\\\p^{n+1}&\p^n&\p^n&1\end{bsmallmatrix},
\end{equation}
because both $t_nv$ and $v_1$ have this property. Then also
\begin{equation}\label{tnlevelraisinglemmaeq4}
 v_2=t_{n+1}v_2\text{ is invariant under }\begin{bsmallmatrix}1&\p^{-1}&\p^{-1}&\p^{-n-1}\\&1&&\p^{-1}\\&&1&\p^{-1}\\&&&1\end{bsmallmatrix}.
\end{equation}
We see from \eqref{tnlevelraisinglemmaeq3} and \eqref{tnlevelraisinglemmaeq4} that $\eta^{-1}v_2\in V(n-1)$. Applying $t_n$ to $t_nv=v_1+v_2$, we get
$$
 v=t_nv_1+t_nv_2=v_1+\eta^{-1} t_{n+1} v_2=v_1+\eta^{-1}v_2\in V(n)+V(n-1),
$$
as claimed.

\ref{tnlevelraisinglemmaitem3}. This is obvious.
\end{proof}

Let $(\pi,V)$ be a smooth representation of $\GSp(4,F)$ for which the center acts trivially, and let $n$ be an integer such that $n\geq 0$. We describe the level raising operator $T_1:V_s(n) \to V_s(n+1)$ induced by $g=1$. 
If $n=0$, then the map $T_1:V_s(0) \to V_s(1)$ is just inclusion, since $\Ks{0}=\GSp(4,\OF)$ and $\Ks{}=\Kl{}$.
We need the following lemma for the case $n \geq 1$.

\begin{lemma}
\label{g1decomplemma}
Let $n$ be an integer such that $n \geq 1$. There is a disjoint decomposition
\begin{equation}
\label{nnplusonedecompeq}
\Ks{n+1}=
\bigsqcup_{x \in \OF/\p} 
\begin{bsmallmatrix}
1&&&x\varpi^{-n}\\
&1&&\\
&&1&\\
&&&1
\end{bsmallmatrix}
\left( \Ks{n+1} \cap \Ks{n} \right).
\end{equation}
\end{lemma}
\begin{proof}
It is clear that the union of the cosets on the right-hand side of \eqref{nnplusonedecompeq}
is contained in $\Ks{n+1}$. 
Using \eqref{Ksniwahorieq1} (with $n$ replaced by $n+1$), an 
argument shows that $\Ks{n+1}$ is contained in this union. 
Finally, it is easy to see that the cosets on the right-hand side of \eqref{nnplusonedecompeq} are  disjoint.
\end{proof}

Again let the notation be as in the paragraph preceding Lemma \ref{g1decomplemma},
and assume that $n  \geq 1$. 
Then by Lemma \ref{g1decomplemma} the level raising operator 
$T_1: V_s(n) \to V_s(n+1)$ is given by
$$
T_1(v) = 
\sum_{x \in \OF/\p}
\pi(
\begin{bsmallmatrix}
1&&&x\varpi^{-n}\\
&1&&\\
&&1&\\
&&&1
\end{bsmallmatrix})v
$$
for $v \in V_s(n)$. In fact, we will work with 
a slight modification of $T_1$ when $n \geq 1$. For 
any integer $n$ such that $n \geq 0$, define 
\begin{equation}\label{taunlevelraisinglemmaeq1}
\tau_n:\:V_s(n)\longrightarrow V_s(n+1)
\end{equation}
by
\begin{equation}\label{taundefeq}
 \tau_n v=\int\limits_\OF\pi(\begin{bsmallmatrix}1&&&z\varpi^{-n}\\&1\\&&1\\&&&1\end{bsmallmatrix})v\,dz
\end{equation}
for $v \in V_s(n)$. We see that if $n=0$, then $\tau_n=T_1$;
however, if $n \geq 1$, then $\tau_n = q^{-1} T_1$. The operator $\tau_n:V_s(n) \to V_s(n+1)$
is clearly an upper block operator. 
We also note that 
using $\tau_n$, equation \eqref{pnVsneq3} for the paramodularization operator $p_n$ can  be written as
\begin{equation}\label{pnVsneq2}
 (q+1)p_nv=q\tau_nv+t_nv\qquad\text{for }v\in V_s(n).
\end{equation}

\begin{lemma}
\label{taunlevelraisinglemma}
 Let $(\pi,V)$ be a smooth representation of $\GSp(4,F)$ for which the center acts trivially, and let $n$ be an integer such that $n\geq 0$.
 \begin{enumerate}
  \item \label{taunlevelraisinglemmaitem1} The level raising operator $\tau_n:V_s(n) \to V_s(n+1)$ is injective.
  \item\label{taunlevelraisinglemmaitem2} The level raising operator $\tau_n:V_s(n) \to V_s(n+1)$ induces a linear map
   \begin{equation}\label{taunlevelraisinglemmaeq2}
    \tau_n:\:\bar V_s(n)\longrightarrow\bar V_s(n+1),
    \end{equation}
    which is also injective.
  \item\label{taunlevelraisinglemmaitem3} We have $\tau_n=-q^{-1}t_n$ as operators from $\bar V_s(n)$ to $\bar V_s(n+1)$.
 \end{enumerate}
\end{lemma}
\begin{proof}
\ref{taunlevelraisinglemmaitem1}.  To prove injectivity, let $v\in V_s(n)$ and assume that $\tau_nv=0$. Then, from \eqref{pnVsneq2} we get
\begin{equation}\label{taunlevelraisinglemmaeq3}
 p_nv=\frac1{q+1}t_nv\qquad\text{for }v\in V_s(n).
\end{equation}
Applying $t_n$ to both sides of this equation, we see that $v\in V(n)$. It now follows that $v=\tau_nv=0$.

\ref{taunlevelraisinglemmaitem2} and \ref{taunlevelraisinglemmaitem3}. It follows from \eqref{pnVsneq2} that
\begin{equation}\label{tauntneq}
 q\,\tau_nv+V(n)=-t_nv+V(n)\qquad\text{for }v\in V_s(n).
\end{equation}
Since $t_n$ maps $V(n-1)+V(n)$ to $V(n)+V(n+1)$, so does $\tau_n$. Thus $\tau_n$ induces a map $\bar V_s(n)\longrightarrow\bar V_s(n+1)$. It follows from \eqref{tauntneq} that this map coincides with $-q^{-1}t_n$. Since \eqref{tnlevelraisinglemmaeq2} is injective, so is \eqref{taunlevelraisinglemmaeq2}.
\end{proof}

Let the assumptions and notation be as in Lemma \ref{taunlevelraisinglemma}. 
Two other expressions for $\tau_n$ will be useful. 
First, if $v \in V_s(n)$, then $\tau_n$ also can  be written as a sum 
$$
\tau_n v =q^{-1} \sum_{z \in \OF/\p} \pi (
\begin{bsmallmatrix}
1&&&z\varpi^{-n}\\
&1&&\\
&&1&\\
&&&1
\end{bsmallmatrix})v.
$$
Second, 
let 
$dk$ be the Haar measure on $\Ks{n+1}$ that assigns $\Ks{n+1}$ volume $1$.
If $v \in V_s(n)$, then $\tau_n v$ can be simply written as
\begin{equation}
\label{tauinteq}
\tau_n v = \int\limits_{\Ks{n+1}} \pi(k)v\, dk.
\end{equation}
We  extend $\tau_n: V_s(n) \to V_s(n+1)$ to all of $V$ by 
defining $\tau_n v$ for $v \in V$ to be as in \eqref{taundefeq}.
This extensions simplifies formulas involving iterations. That is,
let $m$ be an integer such that $m \geq n$. Then we have
$$
\tau_m \tau_{m-1} \cdots \tau_n v = \tau_m v
$$
for $v \in V_s(n)$. 
Finally, if $n$ is an integer such that $n<0$, then we will also define $\tau_n:V\to V$
by \eqref{taundefeq}.

Again let
$(\pi,V)$ be a smooth representation of $\GSp(4,F)$ for which
the center acts trivially, and let $n$ be an integer such that $n \geq 0$.
Our third level raising operator 
is  $T_g:V_s(n) \to V_s(n+1)$ with $g$ as in \eqref{01elementeq}.
The following lemma provides the required coset decomposition.

\begin{lemma}
\label{thetagdecomplemma}
Let $n$ be an integer such that $n \geq 0$, and let
$$
J = \Ks{n+1} \cap
\begin{bsmallmatrix}
1&&&\\
&1&&\\
&&\varpi&\\
&&&\varpi
\end{bsmallmatrix}
\Ks{n}
\begin{bsmallmatrix}
1&&&\\
&1&&\\
&&\varpi&\\
&&&\varpi
\end{bsmallmatrix}^{-1}.
$$
There is a disjoint decomposition
\begin{equation}
\label{thetagdecomplemmaeq1}
\Ks{n+1} = s_2 J \sqcup \bigsqcup_{x \in \OF/\p}
\begin{bsmallmatrix}
1&&&\\
&1&&\\
&x&1&\\
&&&1
\end{bsmallmatrix} J.
\end{equation}
\end{lemma}
\begin{proof}
Let $k \in \GSp(4,F)$. 
Using \eqref{Ksshapeeq}, a calculation shows that
$k \in J$ if and only if $\lambda(k) \in \OF^\times$
and 
$$
k \in 
\begin{bsmallmatrix}
\OF&\OF&\OF&\p^{-n}\\
\p^{n+1}&\OF&\OF&\OF\\
\p^{n+1}&\p&\OF&\OF\\
\p^{n+1}&\p^{n+1}&\p^{n+1}&\OF
\end{bsmallmatrix}.
$$
Again using \eqref{Ksshapeeq}, it follows that the union of the cosets on
the right-hand side of \eqref{thetagdecomplemmaeq1} is contained 
in $\Ks{n+1}$. Conversely, let $k \in \Ks{n+1}$, and write
$k=(k_{ij})_{1 \leq i,j \leq 4}$. Assume that $k_{22} \in \OF^\times$. 
Then 
$$
\begin{bsmallmatrix}
1&&&\\
&1&&\\
&-k_{32}k_{22}^{-1} &1& \\
&&&1
\end{bsmallmatrix}k 
\in 
\begin{bsmallmatrix}
\OF&\OF&\OF&\p^{-n}\\
\p^{n+1}&\OF&\OF&\OF\\
\p^{n+1}&\p&\OF&\OF\\
\p^{n+1}&\p^{n+1}&\p^{n+1}&\OF
\end{bsmallmatrix}.
$$
This implies that $k$ is in the union. Assume that $k_{22} \in \p$. Then
$$
s_2 k
\in 
\begin{bsmallmatrix}
\OF&\OF&\OF&\p^{-n}\\
\p^{n+1}&\OF&\OF&\OF\\
\p^{n+1}&\p&\OF&\OF\\
\p^{n+1}&\p^{n+1}&\p^{n+1}&\OF
\end{bsmallmatrix}.
$$
This again implies that $k$ is in the union. It is straightforward
to verify that the cosets on the right-hand side of \eqref{thetagdecomplemmaeq1}
are disjoint. 
\end{proof}

By Lemma \ref{thetagdecomplemma}, with the notation preceding this lemma,
we see that $T_g$ sends $v \in V_s(n)$ to the following element of $V_s(n+1)$:
$$
\pi(s_2)\pi(\begin{bsmallmatrix}
1&&&\\
&1&&\\
&&\varpi&\\
&&&\varpi
\end{bsmallmatrix}) v+ \sum_{x \in \OF/\p} \pi(
\begin{bsmallmatrix}
1&&&\\
&1&&\\
&x&1&\\
&&&1
\end{bsmallmatrix}) \pi(\begin{bsmallmatrix}
1&&&\\
&1&&\\
&&\varpi&\\
&&&\varpi
\end{bsmallmatrix}) v.
$$
Since $s_2 \in \Ks{m}$ for all $m \geq 0$, we see that
 this vector is equal to 
$$
\pi(\begin{bsmallmatrix}
1&&&\\
&1&&\\
&&\varpi&\\
&&&\varpi
\end{bsmallmatrix}) v+ \sum_{x \in \OF/\p} \pi(
\begin{bsmallmatrix}
1&&&\\
&1&x&\\
&&1&\\
&&&1
\end{bsmallmatrix}) \pi(\begin{bsmallmatrix}
1&&&\\
&\varpi&&\\
&&1&\\
&&&\varpi
\end{bsmallmatrix}) v.
$$
Thus, if $v \in V(n)$, then $T_g(v) =\theta(v)$,
where $\theta:V(n) \to V(n+1)$ is the paramodular level raising operator from Sect.~\ref{parasec}.
Consequently, we also will write $\theta$ for $T_g$, so that 
\begin{equation}\label{thetadefeq1}
 \theta:\:V_s(n)\longrightarrow V_s(n+1)
\end{equation}
is given by 
$$
\theta v= \pi(s_2)\pi(\begin{bsmallmatrix}
1&&&\\
&1&&\\
&&\varpi&\\
&&&\varpi
\end{bsmallmatrix})v + \sum_{x \in \OF/\p} \pi(
\begin{bsmallmatrix}
1&&&\\
&1&&\\
&x&1&\\
&&&1
\end{bsmallmatrix}) \pi(\begin{bsmallmatrix}
1&&&\\
&1&&\\
&&\varpi&\\
&&&\varpi
\end{bsmallmatrix}) v
$$
for $v \in V_s(n)$. 
Some alternative expressions for $\theta$ are
\begin{align*}
\theta v&= \pi(
\begin{bsmallmatrix}
1&&&\\
&1&&\\
&&\varpi&\\
&&&\varpi
\end{bsmallmatrix})v
+
\sum_{x \in \OF/\p}
\pi(
\begin{bsmallmatrix}
1&&&\\
&1&x&\\
&&1&\\
&&&1
\end{bsmallmatrix}
\begin{bsmallmatrix}
1&&&\\
&\varpi&&\\
&&1&\\
&&&\varpi
\end{bsmallmatrix})v,\\
\theta v&= \pi(
\begin{bsmallmatrix}
1&&&\\
&\varpi&&\\
&&1&\\
&&&\varpi
\end{bsmallmatrix})v
+
\sum_{x \in \OF/\p}
\pi(
\begin{bsmallmatrix}
1&&&\\
&1&&\\
&x&1&\\
&&&1
\end{bsmallmatrix}
\begin{bsmallmatrix}
1&&&\\
&1&&\\
&&\varpi&\\
&&&\varpi
\end{bsmallmatrix})v
\end{align*}
for $v \in V_s(n)$. 
Just as for $\tau_n$, we may express $\theta$ using an integral over $\OF$. Thus,
\begin{align}
\label{thetadefeq3}
\theta v
& = \pi(
\begin{bsmallmatrix}
1&&&\\
&1&&\\
&&\varpi&\\
&&&\varpi
\end{bsmallmatrix} )v
+
q \int\limits_{\OF}
\pi(
\begin{bsmallmatrix}
1&&&\\
&1&c&\\
&&1&\\
&&&1
\end{bsmallmatrix}
\begin{bsmallmatrix}
1&&&\\
&\varpi&&\\
&&1&\\
&&&\varpi
\end{bsmallmatrix}
)v\, dc,\\
\theta v & = \pi(
\begin{bsmallmatrix}
1&&&\\
&\varpi&&\\
&&1&\\
&&&\varpi
\end{bsmallmatrix} )v
+
q \int\limits_{\OF}
\pi(
\begin{bsmallmatrix}
1&&&\\
&1&&\\
&c&1&\\
&&&1
\end{bsmallmatrix}
\begin{bsmallmatrix}
1&&&\\
&1&&\\
&&\varpi&\\
&&&\varpi
\end{bsmallmatrix}
)v\, dc
\end{align}
for $v \in V_s(n)$. The operator $\theta:V_s(n) \to V_s(n+1)$ is clearly
an upper block operator. 
\begin{lemma}
\label{thetastableklingenlemma}
 Let $(\pi,V)$ be a smooth representation of $\GSp(4,F)$ for which the center acts trivially, and let $n$ be an integer such that $n\geq 0$.
The level raising operator $\theta:V_s(n) \to V_s(n+1)$ induces a linear map
\begin{equation}
\label{thetadefeq23}
 \theta:\:\bar V_s(n)\longrightarrow \bar V_s(n+1).
\end{equation}
\end{lemma}
\begin{proof}
 This follows from $\theta V(m) \subset V(m+1)$ for integers $m$ such that $m \geq 0$ (see Sect.~\ref{parasec}).
\end{proof}

Finally, we point out yet another useful expression for $\theta$.
Let the notation be as in Lemma \ref{thetastableklingenlemma}. 
Let $dk$ be the Haar measure on $\GSp(4,F)$ that assigns $\Ks{n+1}$ volume $1$. 
Then 
\begin{equation}
\label{thetaineq}
\theta v
= (q+1)\int\limits_{\Ks{n+1}} \pi (k 
\begin{bsmallmatrix}
1&&&\\
&1&&\\
&&\varpi&\\
&&&\varpi
\end{bsmallmatrix})v\, dk
=
(q+1)\int\limits_{\Ks{n+1}} \pi (k
\begin{bsmallmatrix}
1&&&\\
&\varpi&&\\
&&1&\\
&&&\varpi
\end{bsmallmatrix})v\, dk
\end{equation}
for $v \in V_s(n)$.

\section{Four conditions}\label{impcondsec}
In this section we prove the equivalence
of four conditions on a stable Klingen vector.
In Chap.~\ref{dimchap}, these conditions will play an important role
in proving that an irreducible, admissible representation
$\pi$ of $\GSp(4,F)$ with trivial central character
admits non-zero stable Klingen vectors if and only
if $\pi$ is paramodular, and in establishing the partition 
of paramodular representations into two classes. We begin
by introducing a variant of paramodularization.

\begin{lemma}
\label{rhopoplemma}
Let $(\pi,V)$ be a smooth representation of $\GSp(4,F)$ for which the center acts trivially, and let $n$ be an integer such that $n\geq 0$. Define 
\begin{equation}
\label{rhopoplemmaeq01}
\rho_n' : V_s(n) \longrightarrow V(n+1)
\end{equation}
by 
\begin{equation}
\label{rhoprimedefeq}
\rho_n' v = (1+q^{-1}) p_{n+1} t_n v
\end{equation}
for $v \in V_s(n)$.
\begin{enumerate}
\item \label{rhopoplemmaitem1} The linear map $\rho_n'$ has the alternative form
\begin{equation}
\label{rhopoplemmaeq2}
\rho_n' v = q^{-1} \eta v + \tau_{n+1} t_n v 
\end{equation}
for $v \in V_s(n)$.
\item \label{rhopoplemmaitem2} Assume that $v \in V(n)$. Then 
\begin{equation}
\label{rhopoplemmaeq4}
\rho_n' v = q^{-1} \eta v + \tau_{n+1} v = q^{-1} \theta ' v.
\end{equation}
\item \label{rhopoplemmaitem3} Assume that $n\geq 1$ and $v \in V(n-1)$. Then
\begin{equation}
\label{rhopoplemmaeq5}
\rho_n' v = (1+q^{-1}) \eta v.
\end{equation}
\item \label{rhopoplemmaitem4} Assume that $v \in V_s(n)$ and that the paramodularization $p_nv$ is zero. Then
\begin{equation}
\label{rhopoplemmaeq3}
\rho_n' v = q^{-1} \eta v - q \tau_{n+1} v. 
\end{equation}
\end{enumerate} 
\end{lemma}
\begin{proof}
\ref{rhopoplemmaitem1}. Let $v \in V_s(n)$. We have, using \eqref{pnVsneq2} and \ref{tnlevelraisinglemmaitem3} of Lemma \ref{tnlevelraisinglemma},
\begin{align*}
\rho_n' v &= q^{-1}(q+1) p_{n+1}t_n v \\
&= q^{-1} (t_{n+1} t_n v + q \tau_{n+1} t_nv) \\
&= q^{-1}\eta v +  \tau_{n+1} t_n v.
\end{align*}

\ref{rhopoplemmaitem2}. Let $v \in V(n)$. Then $t_n v =v$, and hence by \ref{rhopoplemmaitem1}, $\rho_n' v = q^{-1}\eta v +  \tau_{n+1}  v$. This
is $q^{-1} \theta ' v$ by \eqref{thetap1defeq}.

\ref{rhopoplemmaitem3}. Let $v \in V(n-1)$. Then by \ref{rhopoplemmaitem1},
\begin{align*}
\rho_n'v &=q^{-1} \eta v + \tau_{n+1} t_n v\\
&= q^{-1} \eta v +\eta \tau_{n-1} t_{n-1} v \\
&= q^{-1} \eta v +\eta  v \\
& = (1+q^{-1}) \eta v.
\end{align*}

\ref{rhopoplemmaitem4}. Let $v \in V_s(n)$ and assume that $p_nv =0$. By \eqref{pnVsneq2} we have $t_n v = -q \tau_n v$. Hence,
by \ref{rhopoplemmaitem1},
\begin{align*}
\rho_n' v &= q^{-1} \eta v +\tau_{n+1}(-q \tau_n v) \\
&= q^{-1} \eta v -q\tau_{n+1} \tau_n v \\
&= q^{-1} \eta v -q\tau_{n+1} v.
\end{align*}
This completes the proof. 
\end{proof}

We will now prove the main result of this section using $\rho_n'$. 

\begin{lemma}\label{funnyidentitylemma5}
 Let $(\pi,V)$ be a smooth representation of $\GSp(4,F)$ for which the center acts trivially. Assume that the subspace of vectors of $V$ fixed by ${\rm Sp}(4,F)$ is trivial. Let $n$ be an integer such that $n \geq 0$, and let $v \in V_s(n)$. The following are equivalent:
\begin{enumerate}
\item \label{funnyidentitylemma5item1} The vector $v$ satisfies
 \begin{equation}\label{funnyidentitylemma5eq1}
  \tau_{n+1}v=q^{-2}\eta v.
 \end{equation}
\item \label{funnyidentitylemma5item2} For all integers $m \geq n$, the paramodular projection $p_m(\tau_{m-1}v)$ from \eqref{pnVsneq} is zero, or equivalently,
 \begin{equation}\label{funnyidentitylemma5eq2}
-q\tau_m v =t_m \tau_{m-1} v.
 \end{equation} 
\item\label{funnyidentitylemma5item3}The paramodular projections $p_n(v)$ and $p_{n+1}(\tau_{n}v )$ from \eqref{pnVsneq} are zero, or equivalently, 
 \begin{equation}\label{funnyidentitylemma5eq21}
  -q\tau_n v=t_nv \quad \text{and} \quad -q\tau_{n+1} v=t_{n+1} \tau_n v.
 \end{equation}
\item \label{funnyidentitylemma5item4} The paramodular projection  $p_n(v)$ from \eqref{pnVsneq} is zero and $\rho_n'(v)=0$. 
\end{enumerate}
\end{lemma}
\begin{proof}
\ref{funnyidentitylemma5item1} $\Rightarrow$ \ref{funnyidentitylemma5item2}. Write $v=v_1+v_2$ with $v_1\in V(n)$ and $v_2\in V_s(n)$ such that $p_n(v_2)=0$ (see \eqref{pnVsneq4}). Then \begin{align*}
 \rho'_nv&=\rho'_nv_1+\rho'_nv_2\\
 &=q^{-1}\eta v_1+\tau_{n+1}v_1+q^{-1}\eta v_2-q\tau_{n+1}v_2
 \quad\text{(see \eqref{rhopoplemmaeq4} and \eqref{rhopoplemmaeq3})}\\
 &=q^{-1}\eta(v_1+v_2)+(q+1)\tau_{n+1}v_1-q\tau_{n+1}(v_1+v_2)\\
 &=q^{-1}\eta v-q\tau_{n+1}v+(q+1)\tau_{n+1}v_1\\
 &=(q+1)\tau_{n+1}v_1
\quad\text{(see \eqref{funnyidentitylemma5eq1})}\\
 &=(q+1)q^{-1}\sum_{z\in\OF/\p}\pi(\begin{bsmallmatrix}1&&&z\varpi^{-n-1}\\&1\\&&1\\&&&1\end{bsmallmatrix})v_1,
\end{align*}
so that, by \eqref{thetap1defeq},
\begin{equation*}
 \rho'_nv=(q+1)q^{-1}\theta'v_1-(q+1)q^{-1}\eta v_1.
\end{equation*}
Observe that $\rho'_nv\in V(n+1)$ by \eqref{rhopoplemmaeq01}. By the linear independence of paramodular vectors at different levels, Theorem \ref{linindparatheorem} (this theorem uses the hypothesis that the subspace of vectors fixed by $\SSp(4,F)$ is trivial), we get $\eta v_1=0$, and hence $v_1=0$. Thus $v=v_2$, so that $p_nv=0$. Using \eqref{pnVsneq2}, we obtain \eqref{funnyidentitylemma5eq2} for $m=n$. 

Now let $m$ be an integer such that $m \geq n+1$ and set $v' =\tau_{m-1} v$. Then $v' \in V_s(m)$. Moreover, 
\begin{align*}
\tau_{n+1}v&=q^{-2}\eta v\\
\tau_{m+1} \tau_{n+1}v & = q^{-2} \tau_{m+1} \eta v\\
\tau_{m+1} v & = q^{-2}  \eta \tau_{m-1} v\\
\tau_{m+1} \tau_{m-1}v & = q^{-2}  \eta \tau_{m-1} v\\
\tau_{m+1}v' & = q^{-2}  \eta v'.
\end{align*}
Arguing as in the last paragraph with $v'$ replacing $v$ shows that $-q\tau_{m} v'=t_m v'$, or equivalently, $-q\tau_{m} v=t_m  \tau_{m-1}v$; this is  \eqref{funnyidentitylemma5eq2}.

\ref{funnyidentitylemma5item2} $\Rightarrow$ \ref{funnyidentitylemma5item3}. This is trivial.

\ref{funnyidentitylemma5item3} $\Rightarrow$ \ref{funnyidentitylemma5item1}.  We have $\tau_{n+1}v = -q^{-1} t_{n+1} \tau_n v = q^{-2} t_{n+1} t_n v= q^{-2} \eta v$.

\ref{funnyidentitylemma5item1} $\Rightarrow$ \ref{funnyidentitylemma5item4}. See the proof of \ref{funnyidentitylemma5item1} $\Rightarrow$ \ref{funnyidentitylemma5item2}.

\ref{funnyidentitylemma5item4} $\Rightarrow$ \ref{funnyidentitylemma5item1}. This follows from \ref{rhopoplemmaitem4} of Lemma \ref{rhopoplemma}.
\end{proof}

\section{Level lowering operators}\label{sigmasec}

Let $(\pi,V)$ be a smooth representation of $\GSp(4,F)$ for which the center
acts trivially, and let $n$ be an integer such that $n \geq 1$. In this work
we will use the two level lowering operators $T_g:V_s(n+1) \to V_s(n)$ determined
by $g=1$ and $g= \eta^{-1}$. 

\begin{lemma}
\label{tracedecomplemma}
Let $n$ be an integer such that $n \geq 1$. There is a disjoint decomposition
\begin{equation}
\label{tracedecompeq}
\Ks{n}
=
\bigsqcup_{x,y,z \in \OF/\p}
\begin{bsmallmatrix}
1&&&\\
x\varpi^{n}&1&&\\
y\varpi^{n}&&1&\\
z\varpi^{n}&y\varpi^{n}&-x\varpi^{n}&1
\end{bsmallmatrix}
\left(
\Ks{n} \cap \Ks{n+1} \right).
\end{equation}
\end{lemma}
\begin{proof}
It is clear that the union of the cosets on 
the right-hand side of \eqref{tracedecompeq}
is contained in $\Ks{n}$. Conversely, let
$k \in \Ks{n}$. By \eqref{Ksniwahorieq2}, we may write $k=k_1k_2k_3$
where $k_1,k_2,k_3 \in \Ks{n}$ and 
$$
k_1 \in 
\begin{bsmallmatrix}
1&&&\\
\p^{n}&1&&\\
\p^{n}&&1&\\
\p^{n}&\p^{n}&\p^{n}&1
\end{bsmallmatrix},\quad 
k_2 \in
\begin{bsmallmatrix}
\vphantom{\p^{n}}\OF^\times&&&\\
&\OF\vphantom{\p^{n}}&\OF&\\
&\OF\vphantom{\p^{n}}&\OF&\\
&&&\vphantom{\p^{n}}\OF^\times
\end{bsmallmatrix},\quad
\text{and}\quad
k_3 \in
\begin{bsmallmatrix}
1\vphantom{\p^{n}}&\OF&\OF&\p^{-n+1}\\
&1\vphantom{\p^{n}}&&\OF\\
&&1\vphantom{\p^{n}}&\OF\\
&&&1\vphantom{\p^{n}}
\end{bsmallmatrix}.
$$
Since $k_2k_3 \in \Ks{n} \cap \Ks{n+1}$, we see that $k$ is in the union.
It is straightforward to see that the cosets in \eqref{tracedecompeq}
are disjoint.
\end{proof}

Let the notation be as in the paragraph preceding Lemma \ref{tracedecomplemma}. By \eqref{Tgvsumeq} and Lemma \ref{tracedecomplemma}, the level lowering operator $T_1:V_s(n+1) \to V_s(n)$
sends $v \in V_s(n+1)$ to 
$$
\sum_{x,y,z \in \OF/\p} 
\pi(
\begin{bsmallmatrix}
1&&&\\
x\varpi^n&1&&\\
y\varpi^n&&1&\\
z\varpi^n&y\varpi^n&-x\varpi^n&1
\end{bsmallmatrix})v.
$$
In fact, we will work with a multiple of $T_1$. Define
\begin{equation}
\label{snsetupdefeq}
s_n: V_s(n+1) \longrightarrow V_s(n)
\end{equation}
by
\begin{equation}
\label{sndefeq}
s_n v = \int\limits_{\OF}\int\limits_{\OF}\int\limits_{\OF}
\pi( 
\begin{bsmallmatrix}
1&&&\\
x \varpi^{n}&1&&\\
y\varpi^{n}&&1&\\
z\varpi^{n}&y\varpi^{n}&-x\varpi^{n}&1
\end{bsmallmatrix})v\, dx\, dy\, dz
\end{equation}
for $v \in V_s(n+1)$. Evidently, $s_n v=q^{-3} T_1 v$ for $v \in V_s(n+1)$. 
It is easy to see that $s_n:V_s(n+1) \to V_s(n)$ is not an upper block operator.

Next, we consider the level lowering operator $T_{\eta^{-1}}:V_s(n+1) \to V_s(n)$.
As usual, we need a coset decomposition.

\begin{lemma}
\label{sigmadecomplemma}
Let $n$ be an integer such that $n \geq 1$, and let
$$
J=\Ks{n} \cap \eta^{-1} \Ks{n+1} \eta.
$$
There is a disjoint decomposition
\begin{equation}
\label{sigmadecompeq}
\Ks{n}
=
\bigsqcup_{x,y,z \in \OF/\p}
\begin{bsmallmatrix}
1&x&y&z\varpi^{-n+1}\\
&1&&y\\
&&1&-x\\
&&&1
\end{bsmallmatrix} J.
\end{equation}
\end{lemma}
\begin{proof}
Using \eqref{Ksshapeeq}, we see that
if $k \in \GSp(4,F)$, then $k \in J$
if and only if $\lambda(k) \in \OF^\times$
and 
$$
k \in
\begin{bsmallmatrix}
\OF&\p&\p&\p^{-n+2}\\
\p^n&\OF&\OF&\p\\
\p^n&\OF&\OF&\p\\
\p^n&\p^n&\p^n&\OF
\end{bsmallmatrix}.
$$
It follows that the union of the cosets
on the right-hand side of \eqref{sigmadecompeq}
is contained in $\Ks{n}$. 
Conversely, let
$k \in \Ks{n}$. By \eqref{Ksniwahorieq1}, we may write $k=k_1k_2k_3$
where $k_1,k_2,k_3 \in \Ks{n}$ and 
$$
k_1 \in 
\begin{bsmallmatrix}
1\vphantom{\p^{n}}&\OF&\OF&\p^{-n+1}\\
&1\vphantom{\p^{n}}&&\OF\\
&&1\vphantom{\p^{n}}&\OF\\
&&&1\vphantom{\p^{n}}
\end{bsmallmatrix},\quad 
k_2 \in
\begin{bsmallmatrix}
\vphantom{\p^{n}}\OF^\times&&&\\
&\OF\vphantom{\p^{n}}&\OF&\\
&\OF\vphantom{\p^{n}}&\OF&\\
&&&\vphantom{\p^{n}}\OF^\times
\end{bsmallmatrix},\quad
\text{and}\quad
k_3 \in
\begin{bsmallmatrix}
1&&&\\
\p^{n}&1&&\\
\p^{n}&&1&\\
\p^{n}&\p^{n}&\p^{n}&1
\end{bsmallmatrix}.
$$
Since $k_2k_3 \in J$, we see that $k$ is in the union.
It is straightforward to see that the cosets in \eqref{sigmadecompeq}
are disjoint.
\end{proof}

Let the notation be as at the beginning of this section. By \eqref{Tgvsumeq} and Lemma~\ref{sigmadecomplemma},
the level lowering operator $T_{\eta^{-1}}: V_s(n+1) \to V_s(n)$ sends $v \in V_s(n+1)$
to 
$$
\sum_{x,y,z \in \OF/\p}
\pi(
 \begin{bsmallmatrix}1&x&y&z\varpi^{-n+1}\\&1&&y\\&&1&-x\\&&&1\end{bsmallmatrix}
\begin{bsmallmatrix}
\varpi&&&\vphantom{\varpi^{-n+1}}\\
&1&&\vphantom{\varpi^{-n+1}}\\
&&1&\vphantom{\varpi^{-n+1}}\\
&&&\varpi^{-1}\end{bsmallmatrix})v.
$$
Again, we prefer to work with a multiple of $T_{\eta^{-1}}$. 
 Define 
\begin{equation}
\label{sigmandefformeq}
\sigma_n : V_s(n+1) \longrightarrow V_s(n)
\end{equation}
by 
\begin{equation}\label{sigmaopseq}
 \sigma_n v =\int\limits_\OF\int\limits_\OF\int\limits_{\OF}\pi(
 \begin{bsmallmatrix}1&x&y&z\varpi^{-n+1}\\&1&&y\\&&1&-x\\&&&1\end{bsmallmatrix}
\begin{bsmallmatrix}
\varpi&&&\vphantom{\varpi^{-n+1}}\\
&1&&\vphantom{\varpi^{-n+1}}\\
&&1&\vphantom{\varpi^{-n+1}}\\
&&&\varpi^{-1}\end{bsmallmatrix})v\,dx\,dy\,dz\end{equation}
for $v \in V_s(n)$. Evidently, $\sigma_n v = q^{-3}T_{\eta^{-1}} v$ for $v \in V_s(n+1)$. 
Clearly, the map $\sigma_n:V_s(n+1) \to V_s(n)$ is an upper block operator.

\begin{lemma}\label{sigmaopslemma}
Let $(\pi,V)$ be a smooth representation of $\GSp(4,F)$ for which the center acts trivially, and 
let $n$ be an integer such that $n \geq 1$. 
 \begin{enumerate}
  \item \label{sigmaopslemmaitem1} If $w\in V(n-1)$, then $\sigma_n\eta w=w$.
  \item \label{sigmaopslemmaitem2} If $n \geq 2$ and $w\in V(n-2)$, then $\sigma_n\eta w=\tau_{n-1}w$.
  \item \label{sigmaopslemmaitem3} If $w\in V(n)$, then $\sigma_n\theta'w=w+q\sigma_nw$, so that $w=\sigma_n(\theta'w-qw)$.
  \item \label{sigmaopslemmaitem4} If $w\in V_s(n)$ with $p_n(w)=0$, then $w=q\sigma_n(\rho_n'w+q\tau_nw)$.
  \item \label{sigmaopslemmaitem5} If $n\geq 2$ and $w \in V_s(n)$, then $t_{n-1}s_{n-1} w = \sigma_n t_n w$. 
  \item \label{sigmaopslemmaitem6} If $n \geq 2$ and $w \in V_s(n)$, then $\sigma_n \tau_n w = \tau_{n-1} \sigma_{n-1} w$.
 \end{enumerate}
\end{lemma}
\begin{proof}

\ref{sigmaopslemmaitem1} and \ref{sigmaopslemmaitem2}. These statements are straightforward.

\ref{sigmaopslemmaitem3}. This follows by a calculation using the formula \eqref{thetap1defeq} for $\theta'$.

\ref{sigmaopslemmaitem4}. Let $w\in V_s(n)$ with $p_n(w)=0$. By \ref{rhopoplemmaitem4} of Lemma \ref{rhopoplemma} we have 
$\rho'_nw=q^{-1}\eta w-q\tau_{n+1}w$.
We calculate
\begin{align*}
\sigma_n\rho'_nw&=\int\limits_\OF\int\limits_\OF\int\limits_{\OF}\pi(
\begin{bsmallmatrix}1&x&y&z\varpi^{-n+1}\\&1&&y\\&&1&-x\\&&&1\end{bsmallmatrix}
\begin{bsmallmatrix}
\varpi&&&\vphantom{\varpi^{-n+1}}\\
&1&&\vphantom{\varpi^{-n+1}}\\
&&1&\vphantom{\varpi^{-n+1}}\\
&&&\varpi^{-1}
\end{bsmallmatrix}
)(\rho'_nw)\,dz\,dy\,dx\\
&=\int\limits_\OF\int\limits_\OF\int\limits_{\OF}\pi(
\begin{bsmallmatrix}1&x&y&z\varpi^{-n+1}\\&1&&y\\&&1&-x\\&&&1\end{bsmallmatrix}
\begin{bsmallmatrix}
\varpi&&&\vphantom{\varpi^{-n+1}}\\
&1&&\vphantom{\varpi^{-n+1}}\\
&&1&\vphantom{\varpi^{-n+1}}\\
&&&\varpi^{-1}
\end{bsmallmatrix}
)(q^{-1}\eta w)\,dz\,dy\,dx\\
 &\quad-\int\limits_\OF\int\limits_\OF\int\limits_{\OF}\pi(
\begin{bsmallmatrix}1&x&y&z\varpi^{-n+1}\\&1&&y\\&&1&-x\\&&&1\end{bsmallmatrix}
\begin{bsmallmatrix}
\varpi&&&\vphantom{\varpi^{-n+1}}\\
&1&&\vphantom{\varpi^{-n+1}}\\
&&1&\vphantom{\varpi^{-n+1}}\\
&&&\varpi^{-1}
\end{bsmallmatrix}
)(q\tau_{n+1}w)\,dz\,dy\,dx\\
&=q^{-1}\int\limits_\OF\int\limits_\OF\int\limits_{\OF}\pi(
\begin{bsmallmatrix}1&x&y&z\varpi^{-n+1}\\&1&&y\\&&1&-x\\&&&1\end{bsmallmatrix})w\,dz\,dy\,dx\\
&\quad-q\int\limits_\OF\int\limits_\OF\int\limits_{\OF}\pi(
\begin{bsmallmatrix}1&x&y&z\varpi^{-n+1}\\&1&&y\\&&1&-x\\&&&1\end{bsmallmatrix}
\begin{bsmallmatrix}
\varpi&&&\vphantom{\varpi^{-n+1}}\\
&1&&\vphantom{\varpi^{-n+1}}\\
&&1&\vphantom{\varpi^{-n+1}}\\
&&&\varpi^{-1}
\end{bsmallmatrix}
)w\,dz\,dy\,dx\\
&=q^{-1}w-q\int\limits_\OF\int\limits_\OF\int\limits_{\OF}\pi(
\begin{bsmallmatrix}1&x&y&z\varpi^{-n+1}\\&1&&y\\&&1&-x\\&&&1\end{bsmallmatrix}
\begin{bsmallmatrix}
\varpi&&&\vphantom{\varpi^{-n+1}}\\
&1&&\vphantom{\varpi^{-n+1}}\\
&&1&\vphantom{\varpi^{-n+1}}\\
&&&\varpi^{-1}
\end{bsmallmatrix}
)(\tau_nw)\,dz\,dy\,dx\\
  &=q^{-1}w-q\sigma_n\tau_nw.
\end{align*}
Hence $w=q\sigma_n(\rho'_nw+q\tau_nw)$, as claimed.

\ref{sigmaopslemmaitem5} and \ref{sigmaopslemmaitem6}. These statements follow by direct calculations.
\end{proof}

Again, other expressions for $\sigma_n$ will be useful. Let the notation be
as in Lemma~\ref{sigmaopslemma}. If $v \in V_s(n+1)$, then $\sigma_n v$
can be written as a sum:
\begin{equation}
\label{sigmaalteq}
\sigma_n v = q^{-3} 
\sum_{x,y,z \in \OF/\p}
\pi(
 \begin{bsmallmatrix}1&x&y&z\varpi^{-n+1}\\&1&&y\\&&1&-x\\&&&1\end{bsmallmatrix}
\begin{bsmallmatrix}
\varpi&&&\vphantom{\varpi^{-n+1}}\\
&1&&\vphantom{\varpi^{-n+1}}\\
&&1&\vphantom{\varpi^{-n+1}}\\
&&&\varpi^{-1}\end{bsmallmatrix})v.
\end{equation}
Also, let $dk$ be the Haar measure on $\GSp(4,F)$ that assigns $\Ks{n}$
volume 1. If $v \in V_s(n+1)$, then 
\begin{equation}
\label{sigmaKsneq}
\sigma_n v 
=
\int\limits_{\Ks{n}} \pi(k 
\begin{bsmallmatrix}
\varpi&&&\\
&1&&\\
&&1&\\
&&&\varpi^{-1}
\end{bsmallmatrix})v\, dk.
\end{equation}

\begin{proposition}\label{sigmasurjectiveprop}
 Let $(\pi,V)$ be a smooth representation of $\GSp(4,F)$ for which the center acts trivially. Let $n$ be an integer such that $n\geq 1$. Then the level lowering operator
 \begin{equation}\label{sigmasurjectivepropeq1}
  \sigma_n:V_s(n+1)\longrightarrow V_s(n)
 \end{equation}
 is surjective. Moreover,
 \begin{equation}\label{sigmasurjectivepropeq2}
  \sigma_n(V(n)+ V(n+1))\supset V(n-1)+ V(n).
 \end{equation}
\end{proposition}
\begin{proof}
By Lemma \ref{paraformlemma} we have $V_s(n)=V(n)\oplus{\rm ker}(p_n)$. If $v\in V(n)$, then $v$ lies in 
the image of $\sigma_n$ by \ref{sigmaopslemmaitem3} of Lemma \ref{sigmaopslemma}. If $v\in{\rm ker}(p_n)$, 
then $v$ lies in the image of $\sigma_n$ by \ref{sigmaopslemmaitem4} of Lemma \ref{sigmaopslemma}. 
Hence the map \eqref{sigmasurjectivepropeq1} is surjective.
Next, \ref{sigmaopslemmaitem1} of Lemma~\ref{sigmaopslemma}  shows that $V(n-1)\subset\sigma_n(V(n+1))$, 
and \ref{sigmaopslemmaitem3} of Lemma \ref{sigmaopslemma} shows that $V(n)\subset\sigma_n(V(n)+ V(n+1))$. 
Hence we get \eqref{sigmasurjectivepropeq2}.
\end{proof}
Let $(\pi,V)$ be a smooth representation of $\GSp(4,F)$ for which the center acts trivially. Let $n$ be an integer such that $n\geq 1$. 
We may consider the composition of $\tau_n:V_s(n) \to V_s(n+1)$ and $\sigma_n: V_s(n+1) \to V_s(n)$, which 
 defines an endomorphism
\begin{equation}\label{sigmaendoeq}
  \sigma_n\circ \tau_n:V_s(n)\longrightarrow V_s(n).
\end{equation}
\begin{lemma}\label{sigmaendolemma}
 Let $(\pi,V)$ be a smooth representation of $\GSp(4,F)$ for which the center acts trivially. For $n\geq2$, the kernel of the endomorphism $\sigma_n \circ \tau_n$ of $V_s(n)$ coincides with the kernel of the map $\sigma_{n-1}:V_s(n)\to V_s(n-1)$.
\end{lemma}
\begin{proof}
By \ref{sigmaopslemmaitem6} of Lemma \ref{sigmaopslemma} we have $(\sigma_n\circ \tau_n) v=\tau_{n-1}(\sigma_{n-1}v)$ for $v \in V_s(n)$. Since the map $\tau_{n-1}:V_s(n-1)\to V_s(n)$ is injective by \ref{taunlevelraisinglemmaitem1} of Lemma \ref{taunlevelraisinglemma}, the assertion follows.
\end{proof}

\section{Stable Hecke operators}
\label{stableheckesec}

Let $(\pi,V)$ be a smooth representation for which the center acts trivially,
and let $n$ be an integer such that $n \geq 1$. In this section we will consider
the Hecke operators $T_g:V_s(n) \to V_s(n)$ determined by
$$
g=\Delta_{0,1}=
\begin{bsmallmatrix}
\varpi&&&\\
&\varpi&&\\
&&1&\\
&&&1
\end{bsmallmatrix}
\quad 
\text{and}
\quad
g=\Delta_{1,0}=
\begin{bsmallmatrix}
\varpi^2&&&\\
&\varpi&&\\
&&\varpi&\\
&&&1
\end{bsmallmatrix}.
$$
The following lemma provides the required coset decompositions; see also Sect.~\ref{opersksec}
and Lemma \ref{heckedefeqlemma}.

\begin{lemma}
Let $n$ be an integer such that $n \geq 1$. 
\label{stable01lemma}
 \begin{enumerate}
  \item \label{stable01lemmaitem1} There is a disjoint decomposition
   \begin{align}\label{stable01lemmaeq0}
     \Ks{n}\begin{bsmallmatrix}
     \varpi&&&\vphantom{\varpi^{-n+1}}\\
     &\varpi&&\vphantom{\varpi^{-n+1}}\\
     &&1&\vphantom{\varpi^{-n+1}}\\
     &&&1\vphantom{\varpi^{-n+1}}
     \end{bsmallmatrix}\Ks{n}&=\bigsqcup_{y,z\in\OF/\p} \begin{bsmallmatrix}1&y&&z\varpi^{-n+1}\\&1&&\\&&1&-y\\&&&1\end{bsmallmatrix}
     \begin{bsmallmatrix}
     \varpi&&&\vphantom{\varpi^{-n+1}}\\
     &1&&\vphantom{\varpi^{-n+1}}\\
     &&\varpi&\vphantom{\varpi^{-n+1}}\\
     &&&1\vphantom{\varpi^{-n+1}}
     \end{bsmallmatrix}\Ks{n}\nonumber\\
     &\quad\sqcup\bigsqcup_{c,y,z\in\OF/\p}\begin{bsmallmatrix}1&&y&z\varpi^{-n+1}\\&1&c&y\\&&1\\&&&1\end{bsmallmatrix}
     \begin{bsmallmatrix}
     \varpi&&&\vphantom{\varpi^{-n+1}}\\
     &\varpi&&\vphantom{\varpi^{-n+1}}\\
     &&1&\vphantom{\varpi^{-n+1}}\\
     &&&1\vphantom{\varpi^{-n+1}}
     \end{bsmallmatrix}\Ks{n}.
   \end{align}
  \item \label{stable01lemmaitem2} There is a disjoint decomposition
   \begin{align}\label{stable01lemmaeq10}
     \Ks{n}
     \begin{bsmallmatrix}
     \varpi^2&&&\vphantom{\varpi^{-n+1}}\\
     &\varpi&&\vphantom{\varpi^{-n+1}}\\
     &&\varpi&\vphantom{\varpi^{-n+1}}\\
     &&&1\vphantom{\varpi^{-n+1}}
     \end{bsmallmatrix}\Ks{n}&=\bigsqcup_{\substack{x,y\in\OF/\p\\z\in\OF/\p^2}} 
     \begin{bsmallmatrix}1&x&y&z\varpi^{-n+1}\\&1&&y\\&&1&-x\\&&&1\end{bsmallmatrix}
     \begin{bsmallmatrix}
     \varpi^2&&&\vphantom{\varpi^{-n+1}}\\
     &\varpi&&\vphantom{\varpi^{-n+1}}\\
     &&\varpi&\vphantom{\varpi^{-n+1}}\\
     &&&1\vphantom{\varpi^{-n+1}}
     \end{bsmallmatrix}\Ks{n}.
   \end{align}
 \end{enumerate}
\end{lemma}
\begin{proof}
\ref{stable01lemmaitem1}. The asserted identity is equivalent to
\begin{align}\label{stable01lemmaeq1}
 \Ks{n}\begin{bsmallmatrix}
 \varpi&&&\vphantom{\varpi^{-n+1}}\\
 &\varpi&&\vphantom{\varpi^{-n+1}}\\
 &&1&\vphantom{\varpi^{-n+1}}\\
 &&&1\vphantom{\varpi^{-n+1}}
 \end{bsmallmatrix}\Ks{n}
 &=\bigsqcup_{y,z\in\OF/\p} s_2
 \begin{bsmallmatrix}1&&y&z\varpi^{-n+1}\\&1&&y\\&&1\\&&&1\end{bsmallmatrix}
 \begin{bsmallmatrix}
 \varpi&&&\vphantom{\varpi^{-n+1}}\\
 &\varpi&&\vphantom{\varpi^{-n+1}}\\
 &&1&\vphantom{\varpi^{-n+1}}\\
 &&&1\vphantom{\varpi^{-n+1}}
 \end{bsmallmatrix}\Ks{n}\nonumber\\
 &\quad\sqcup\bigsqcup_{c,y,z\in\OF/\p}\begin{bsmallmatrix}1&&y&z\varpi^{-n+1}\\&1&c&y\\&&1\\&&&1\end{bsmallmatrix}
 \begin{bsmallmatrix}
 \varpi&&&\vphantom{\varpi^{-n+1}}\\
 &\varpi&&\vphantom{\varpi^{-n+1}}\\
 &&1&\vphantom{\varpi^{-n+1}}\\
 &&&1\vphantom{\varpi^{-n+1}}
 \end{bsmallmatrix}\Ks{n},
\end{align}
where $s_2$ is defined in \eqref{s1s2defeq}. 
Clearly, the right hand side is contained in the left hand side. By the Iwahori factorization \eqref{Ksniwahorieq1}, the left hand side is contained in
$$
 \begin{bsmallmatrix}1\\&\OF&\OF\\&\OF&\OF\\&&&1\end{bsmallmatrix}\begin{bsmallmatrix}1&&\OF&\p^{-n+1}\\&1&&\OF\\&&1&\\&&&1\end{bsmallmatrix}\begin{bsmallmatrix}\varpi\\&\varpi\\&&1\\&&&1\end{bsmallmatrix}\Ks{n}.
$$
Here, and below, when we write, for example, $ \begin{bsmallmatrix}1\\&\OF&\OF\\&\OF&\OF\\&&&1\end{bsmallmatrix}$,
we actually mean the intersection of this set with $\GSp(4,F)$. 
We have
$$
 \SL(2,\OF)=\bigsqcup_{c\in\OF/\p}\begin{bsmallmatrix}1&c\\&1\end{bsmallmatrix}\transpose{\,\Gamma_0(\p)}\sqcup
 \begin{bsmallmatrix}&1\\-1&\end{bsmallmatrix}\transpose{\,\Gamma_0(\p)},
$$
where $\Gamma_0(\p)=\{\begin{bsmallmatrix}a&b\\c&d\end{bsmallmatrix} \in \SL(2,\OF)\mid c \equiv 0\Mod{\p}\}$.
Hence the left hand side of \eqref{stable01lemmaeq1} is contained in
\begin{align*}
&\bigcup_{c\in\OF/\p}
\begin{bsmallmatrix}1\\&1&c\\&&1\\&&&1\end{bsmallmatrix}
\begin{bsmallmatrix}1&&\OF&\p^{-n+1}\\&1&&\OF\\&&1&\\&&&1\end{bsmallmatrix}
\begin{bsmallmatrix}\varpi\\&\varpi\\&&1\\&&&1\end{bsmallmatrix}\Ks{n}\\
&\qquad\quad\cup s_2
\begin{bsmallmatrix}1&&\OF&\p^{-n+1}\\&1&&\OF\\&&1&\\&&&1\end{bsmallmatrix}
\begin{bsmallmatrix}\varpi\\&\varpi\\&&1\\&&&1\end{bsmallmatrix}\Ks{n}.
\end{align*}
It follows that the left hand side of \eqref{stable01lemmaeq1} is contained in the right hand side of \eqref{stable01lemmaeq1}. The disjointness of the right hand side is easy to see.

\ref{stable01lemmaitem2}. This follows by a similar argument.
\end{proof}

Let $(\pi,V)$ be a smooth representation of $\GSp(4,F)$ for which the center acts trivially, and
let $n$ be an integer such that $n \geq 1$. By Lemma \ref{stable01lemma}, there are endomorphisms $T^s_{0,1}$ and $T^s_{1,0}$ of $V_s(n)$ given by
\begin{align}\label{Ts01eq}
 T^s_{0,1}v&=\sum_{y,z\in\OF/\p}\pi(\begin{bsmallmatrix}1&y&&z\varpi^{-n+1}\\&1&&\\&&1&-y\\&&&1\end{bsmallmatrix}
 \begin{bsmallmatrix}
 \varpi&&&\vphantom{\varpi^{-n+1}}\\
 &1&&\vphantom{\varpi^{-n+1}}\\
 &&\varpi&\vphantom{\varpi^{-n+1}}\\
 &&&1\vphantom{\varpi^{-n+1}}
 \end{bsmallmatrix})v\nonumber\\
 &\quad+\sum_{c,y,z\in\OF/\p}\pi(\begin{bsmallmatrix}1&&y&z\varpi^{-n+1}\\&1&c&y\\&&1\\&&&1\end{bsmallmatrix}
 \begin{bsmallmatrix}
 \varpi&&&\vphantom{\varpi^{-n+1}}\\
 &\varpi&&\vphantom{\varpi^{-n+1}}\\
 &&1&\vphantom{\varpi^{-n+1}}\\
 &&&1\vphantom{\varpi^{-n+1}}\end{bsmallmatrix})v
\end{align}
and
\begin{align}\label{Ts10eq}
 T^s_{1,0}v&=\sum_{\substack{x,y\in\OF/\p\\z\in\OF/\p^2}} \pi(\begin{bsmallmatrix}1&x&y&z\varpi^{-n+1}\\&1&&y\\&&1&-x\\&&&1\end{bsmallmatrix}
 \begin{bsmallmatrix}
 \varpi^2&&&\vphantom{\varpi^{-n+1}}\\
 &\varpi&&\vphantom{\varpi^{-n+1}}\\
 &&\varpi&\vphantom{\varpi^{-n+1}}\\
 &&&1\vphantom{\varpi^{-n+1}}
 \end{bsmallmatrix})v
\end{align}
for $v \in V_s(n)$. We refer to the endomorphisms $T_{0,1}^s$ and $T_{1,0}^s$ of 
$V_s(n)$ as the \emph{stable Klingen Hecke operators}, or simply \emph{stable Hecke operators}, of level $n$. \index{stable Hecke operators}
\index{Hecke operator!stable}We note that  $T_{0,1}^s$ and $T_{1,0}^s$ are upper block operators.

\begin{lemma}\label{Heckecommutelemma}
Let $(\pi,V)$ be a smooth representation of $\GSp(4,F)$ for which the center acts trivially, and let $n$ be
an integer such that $n \geq 1$. The endomorphisms $T^s_{0,1}$ and $T^s_{1,0}$ of $V_s(n)$ commute.
\end{lemma}
\begin{proof}
Let $K=\Ks{n}$, and $dk$ be the Haar measure on $\GSp(4,F)$ that
assigns $K$ volume $1$.  Let $v \in V_s(n)$. We have
\begin{align*}
T_{0,1}^s v& = 
 \int\limits_{K} \pi (k) T_{0,1}^s v \, dk\\
&=\sum_{y,z\in\OF/\p}\int\limits_{K}
\pi(k\begin{bsmallmatrix}1&y&&z\varpi^{-n+1}\\&1&&\\&&1&-y\\&&&1\end{bsmallmatrix}
\begin{bsmallmatrix}
\varpi&&&\vphantom{\varpi^{-n+1}}\\
&1&&\vphantom{\varpi^{-n+1}}\\
&&\varpi&\vphantom{\varpi^{-n+1}}\\
&&&1\vphantom{\varpi^{-n+1}}
\end{bsmallmatrix})v\, dk\\
&\quad +\sum_{c,y,z\in\OF/\p} \int\limits_{K}
\pi(k
\begin{bsmallmatrix}1&&y&z\varpi^{-n+1}\\&1&c&y\\&&1\\&&&1\end{bsmallmatrix}
\begin{bsmallmatrix}
\varpi&&&\vphantom{\varpi^{-n+1}}\\
&\varpi&&\vphantom{\varpi^{-n+1}}\\
&&1&\vphantom{\varpi^{-n+1}}\\
&&&1\vphantom{\varpi^{-n+1}}
\end{bsmallmatrix})v\,dk\\
&=q^2 \int\limits_{K}
\pi(k
\begin{bsmallmatrix}\varpi\\&1\\&&\varpi\\&&&1\end{bsmallmatrix})v\, dk
+q^3 \int\limits_{K}
\pi(k
\begin{bsmallmatrix}\varpi\\&\varpi\\&&1\\&&&1\end{bsmallmatrix})v\,dk,\\
\end{align*}
so that
\begin{equation}\label{Heckecommutelemmaeq20}
 T_{0,1}^s v=(q^2+q^3) \int\limits_{K}
\pi(k
\begin{bsmallmatrix}\varpi\\&\varpi\\&&1\\&&&1\end{bsmallmatrix})v\,dk.
\end{equation}
For the last step we recall that $s_2 \in K=\Ks{n}$. Similarly,
\begin{equation}
\label{Heckecommutelemmaeq21}
T_{1,0}^s v
=
q^4 \int\limits_{K}
\pi(k
\begin{bsmallmatrix}\varpi^2\\&\varpi\\&&\varpi\\&&&1\end{bsmallmatrix})v\,dk.
\end{equation}
By \eqref{Heckecommutelemmaeq21} and \eqref{Ts01eq} we now have
\begin{align*}
&T_{1,0}^s(T_{0,1}^s v)
 = q^4\int\limits_{K}
\pi(k
\begin{bsmallmatrix}\varpi^2\\&\varpi\\&&\varpi\\&&&1\end{bsmallmatrix})(T_{0,1}^s v)\,dk\\
&\qquad = \sum_{y,z\in\OF/\p}q^4\int\limits_{K}
\pi(k
\begin{bsmallmatrix}
\varpi^2&&&\vphantom{\varpi^{-n+1}}\\
&\varpi&&\vphantom{\varpi^{-n+1}}\\
&&\varpi&\vphantom{\varpi^{-n+1}}\\
&&&1\vphantom{\varpi^{-n+1}}
\end{bsmallmatrix}
\begin{bsmallmatrix}1&y&&z\varpi^{-n+1}\\&1&&\\&&1&-y\\&&&1\end{bsmallmatrix}
\begin{bsmallmatrix}
\varpi&&&\vphantom{\varpi^{-n+1}}\\
&1&&\vphantom{\varpi^{-n+1}}\\
&&\varpi&\vphantom{\varpi^{-n+1}}\\
&&&1\vphantom{\varpi^{-n+1}}
\end{bsmallmatrix})v\,dk\\
&\qquad\quad+\sum_{c,y,z\in\OF/\p} q^4 \int\limits_{K}
\pi(k
\begin{bsmallmatrix}
\varpi^2&&&\vphantom{\varpi^{-n+1}}\\
&\varpi&&\vphantom{\varpi^{-n+1}}\\
&&\varpi&\vphantom{\varpi^{-n+1}}\\
&&&1\vphantom{\varpi^{-n+1}}
\end{bsmallmatrix}
\begin{bsmallmatrix}1&&y&z\varpi^{-n+1}\\&1&c&y\\&&1\\&&&1\end{bsmallmatrix}
\begin{bsmallmatrix}
\varpi&&\vphantom{\varpi^{-n+1}}\\
&\varpi&&\vphantom{\varpi^{-n+1}}\\
&&1&\vphantom{\varpi^{-n+1}}\\
&&&1\vphantom{\varpi^{-n+1}}
\end{bsmallmatrix})v\,dk\\
&\qquad = \sum_{y,z\in\OF/\p} q^4 \int\limits_{K}
\pi(k
\begin{bsmallmatrix}1&y\varpi&&z\varpi^{-n+3}\\&1&&\\&&1&-y\varpi\\&&&1\end{bsmallmatrix}
\begin{bsmallmatrix}
\varpi^3&&&\vphantom{\varpi^{-n+1}}\\
&\varpi&&\vphantom{\varpi^{-n+1}}\\
&&\varpi^2&\vphantom{\varpi^{-n+1}}\\
&&&1\vphantom{\varpi^{-n+1}}
\end{bsmallmatrix})v\,dk\\
&\qquad\quad+\sum_{c,y,z\in\OF/\p} q^4 \int\limits_{K}
\pi(k
\begin{bsmallmatrix}1&&y\varpi&z\varpi^{-n+3}\\&1&c&y\varpi\\&&1\\&&&1\end{bsmallmatrix}
\begin{bsmallmatrix}
\varpi^3&&&\vphantom{\varpi^{-n+1}}\\
&\varpi^2&&\vphantom{\varpi^{-n+1}}\\
&&\varpi&\vphantom{\varpi^{-n+1}}\\
&&&1\vphantom{\varpi^{-n+1}}
\end{bsmallmatrix})v\,dk\\
&\qquad = q^6\int\limits_{K}
\pi(k
\begin{bsmallmatrix}\varpi^3\\&\varpi\\&&\varpi^2\\&&&1\end{bsmallmatrix})v\,dk
+q^7 \int\limits_{K}
\pi(k
\begin{bsmallmatrix}\varpi^3\\&\varpi^2\\&&\varpi\\&&&1\end{bsmallmatrix})v\,dk\\
&\qquad =(q^6+q^7) \int\limits_{K}
\pi(k
\begin{bsmallmatrix}\varpi^3\\&\varpi^2\\&&\varpi\\&&&1\end{bsmallmatrix})v\,dk.
\end{align*}
Similarly, by \eqref{Heckecommutelemmaeq20} and \eqref{Ts10eq},
\begin{align*}
&T_{0,1}^s(T_{1,0}^s v)= (q^2+q^3)\int\limits_{K}
\pi(k
\begin{bsmallmatrix}\varpi\\&\varpi\\&&1\\&&&1\end{bsmallmatrix})T_{1,0}^s v\,dk\\
&\qquad =\sum_{\substack{x,y\in\OF/\p\\z\in\OF/\p^2}} (q^2+q^3) \int\limits_{K}
\pi(k
\begin{bsmallmatrix}
\varpi&&&\vphantom{\varpi^{-n+1}}\\
&\varpi&&\vphantom{\varpi^{-n+1}}\\
&&1&\vphantom{\varpi^{-n+1}}\\
&&&1\vphantom{\varpi^{-n+1}}
\end{bsmallmatrix}
\begin{bsmallmatrix}1&x&y&z\varpi^{-n+1}\\&1&&y\\&&1&-x\\&&&1\end{bsmallmatrix}
\begin{bsmallmatrix}
\varpi^2&&&\vphantom{\varpi^{-n+1}}\\
&\varpi&&\vphantom{\varpi^{-n+1}}\\
&&\varpi&\vphantom{\varpi^{-n+1}}\\
&&&1\vphantom{\varpi^{-n+1}}
\end{bsmallmatrix})v\,dk\\
&\qquad =\sum_{\substack{x,y\in\OF/\p\\z\in\OF/\p^2}}(q^2+q^3)\int\limits_{K}
\pi(k
\begin{bsmallmatrix}1&x&y\varpi&z\varpi^{-n+2}\\&1&&y\varpi\\&&1&-x\\&&&1\end{bsmallmatrix}
\begin{bsmallmatrix}
\varpi^3&&&\vphantom{\varpi^{-n+1}}\\
&\varpi^2&&\vphantom{\varpi^{-n+1}}\\
&&\varpi&\vphantom{\varpi^{-n+1}}\\
&&&1\vphantom{\varpi^{-n+1}}
\end{bsmallmatrix})v\,dk\\
&\qquad =(q^6+q^7) \int\limits_{K}
\pi(k
\begin{bsmallmatrix}\varpi^3\\&\varpi^2\\&&\varpi\\&&&1\end{bsmallmatrix})v\,dk.
\end{align*}
It now follows that $T_{1,0}^s(T_{0,1}^s v)=T_{0,1}^s(T_{1,0}^s v)$. 
\end{proof}

\begin{lemma}\label{T01Ts01lemma}
Let $(\pi,V)$ be a smooth representation of $\GSp(4,F)$ for which the center acts trivially,
and let $n$ be an integer such that $n \geq 1$. Let $v \in V(n)$. Then
\begin{equation}\label{T01Ts01lemmaeq1}
  (1+q^{-1})p_nT^s_{0,1}v=T_{0,1}v
\end{equation}
and
\begin{equation}\label{T01Ts01lemmaeq2}
  (1+q^{-1})p_nT^s_{1,0}v=T_{1,0}v.
 \end{equation}
\end{lemma}
\begin{proof}
By \eqref{pnVsneq2}, \eqref{Ts01eq}, and the assumption that $v \in V(n)$,
\begin{align*}
&(q+1)p_nT^s_{0,1}v=q\tau_nT^s_{0,1}v+t_nT^s_{0,1}v\\
 &\qquad=q\tau_n\sum_{y,z\in\OF/\p}\pi(\begin{bsmallmatrix}1&y&&z\varpi^{-n+1}\\&1&&\\&&1&-y\\&&&1\end{bsmallmatrix}
 \begin{bsmallmatrix}
 \varpi&&&\vphantom{\varpi^{-n+1}}\\
 &1&&\vphantom{\varpi^{-n+1}}\\
 &&\varpi&\vphantom{\varpi^{-n+1}}\\
 &&&1\vphantom{\varpi^{-n+1}}
 \end{bsmallmatrix})v\\
 &\qquad\quad+q\tau_n\sum_{c,y,z\in\OF/\p}\pi(\begin{bsmallmatrix}1&&y&z\varpi^{-n+1}\\&1&c&y\\&&1\\&&&1\end{bsmallmatrix}
 \begin{bsmallmatrix}
 \varpi&&&\vphantom{\varpi^{-n+1}}\\
 &\varpi&&\vphantom{\varpi^{-n+1}}\\
 &&1&\vphantom{\varpi^{-n+1}}\\
 &&&1\vphantom{\varpi^{-n+1}}
 \end{bsmallmatrix})v\\
 &\qquad\quad+t_n\sum_{y,z\in\OF/\p}\pi(\begin{bsmallmatrix}1&y&&z\varpi^{-n+1}\\&1&&\\&&1&-y\\&&&1\end{bsmallmatrix}
 \begin{bsmallmatrix}
 \varpi&&&\vphantom{\varpi^{-n+1}}\\
 &1&&\vphantom{\varpi^{-n+1}}\\
 &&\varpi&\vphantom{\varpi^{-n+1}}\\
 &&&1\vphantom{\varpi^{-n+1}}
 \end{bsmallmatrix})v\\
 &\qquad\quad+t_n\sum_{c,y,z\in\OF/\p}\pi(\begin{bsmallmatrix}1&&y&z\varpi^{-n+1}\\&1&c&y\\&&1\\&&&1\end{bsmallmatrix}
 \begin{bsmallmatrix}
 \varpi&&&\vphantom{\varpi^{-n+1}}\\
 &\varpi&&\vphantom{\varpi^{-n+1}}\\
 &&1&\vphantom{\varpi^{-n+1}}\\
 &&&1\vphantom{\varpi^{-n+1}}
 \end{bsmallmatrix})v\\
 &\qquad=q\sum_{y,z\in\OF/\p}\pi(\begin{bsmallmatrix}1&y&&z\varpi^{-n}\\&1&&\\&&1&-y\\&&&1\end{bsmallmatrix}
 \begin{bsmallmatrix}
 \varpi&&&\vphantom{\varpi^{-n+1}}\\
 &1&&\vphantom{\varpi^{-n+1}}\\
 &&\varpi&\vphantom{\varpi^{-n+1}}\\
 &&&1\vphantom{\varpi^{-n+1}}
 \end{bsmallmatrix})v\\
 &\qquad\quad+q\sum_{c,y,z\in\OF/\p}\pi(\begin{bsmallmatrix}1&&y&z\varpi^{-n}\\&1&c&y\\&&1\\&&&1\end{bsmallmatrix}
 \begin{bsmallmatrix}
 \varpi&&&\vphantom{\varpi^{-n+1}}\\
 &\varpi&&\vphantom{\varpi^{-n+1}}\\
 &&1&\vphantom{\varpi^{-n+1}}\\
 &&&1\vphantom{\varpi^{-n+1}}
 \end{bsmallmatrix})v\\
 &\qquad\quad+q\sum_{y\in\OF/\p}\pi(t_n\begin{bsmallmatrix}1&y&&\\&1&&\\&&1&-y\\&&&1\end{bsmallmatrix}
 \begin{bsmallmatrix}
 \varpi\vphantom{1}&&\vphantom{y}\\
 &1&\vphantom{c}&\vphantom{y}\\
 &&\varpi\vphantom{1}\\
 &&&1
 \end{bsmallmatrix})v\\
 &\qquad\quad+q\sum_{c,y\in\OF/\p}\pi(t_n\begin{bsmallmatrix}1&&y&\\&1&c&y\\&&1\\&&&1\end{bsmallmatrix}
 \begin{bsmallmatrix}\varpi\vphantom{1}&&\vphantom{y}\\&\varpi\vphantom{1}&\vphantom{c}&\vphantom{y}\\&&1\\&&&1\end{bsmallmatrix})v\\
 &\qquad=qT_{0,1}v.
\end{align*}
The last equality follows from Lemma 6.1.2 of \cite{NF}. This proves \eqref{T01Ts01lemmaeq1}.
Next, by \eqref{pnVsneq2}, \eqref{Ts10eq}, and the assumption that $v \in V(n)$,
\begin{align*}
&(q+1)p_nT^s_{1,0}v=q\tau_nT^s_{1,0}v+t_nT^s_{1,0}v\\
&\qquad=q\sum_{\substack{x,y\in\OF/\p\\z\in\OF/\p^2}} \pi(
\begin{bsmallmatrix}1&x&y&z\varpi^{-n}\\&1&&y\\&&1&-x\\&&&1\end{bsmallmatrix}
\begin{bsmallmatrix}
\varpi^2&&&\vphantom{\varpi^{-n+1}}\\
&\varpi&&\vphantom{\varpi^{-n+1}}\\
&&\varpi&\vphantom{\varpi^{-n+1}}\\
&&&1\vphantom{\varpi^{-n+1}}
\end{bsmallmatrix})v\\
&\qquad\quad+q\sum_{x,y,z\in\OF/\p} \pi(t_n
\begin{bsmallmatrix}1&x&y&z\varpi^{-n+1}\\&1&&y\\&&1&-x\\&&&1\end{bsmallmatrix}
\begin{bsmallmatrix}
\varpi^2&&&\vphantom{\varpi^{-n+1}}\\
&\varpi&&\vphantom{\varpi^{-n+1}}\\
&&\varpi&\vphantom{\varpi^{-n+1}}\\
&&&1\vphantom{\varpi^{-n+1}}\end{bsmallmatrix})v\\
&\qquad=qT_{1,0}v.
\end{align*}
The last equality again follows from Lemma 6.1.2 of \cite{NF}.
\end{proof}

\begin{lemma}
\label{heckeupdownlemma}
Let $(\pi,V)$ be a smooth representation of $\GSp(4,F)$ for 
which the center acts trivially. Let $n$ be a integer such that $n \geq 1$.
Let $v \in V_s(n)$. Then 
\begin{equation}
\label{sigmaT10eq}
T_{0,1}^s v = q^2 \sigma_n \theta v\quad \text{and} \quad 
T_{1,0}^s v = q^4 \sigma_n \tau_n v.
\end{equation}
Assume that $n \geq 2$. Then we also have 
\begin{equation}
\label{heckeupdownlemmaeq1}
T_{1,0}^s v= q^4 \tau_{n-1} \sigma_{n-1} v, 
\end{equation}
and $T_{1,0}^s v =0$ if and only if $\sigma_{n-1} v =0$. 
\end{lemma}
\begin{proof}
Let $K =\Ks{n}$. 
Let $dk$ be the Haar measure on $K$ that assigns $K$ volume 1. 
Let $v \in V_s(n)$.
Then
\begin{align*}
&q^2 \sigma_n \theta v\\
&\qquad = q^2 \int\limits_K \pi(k) \pi(
\begin{bsmallmatrix}
\varpi&&&\\
&1&&\\
&&1&\\
&&&\varpi^{-1} 
\end{bsmallmatrix})
\theta v \, dk\qquad \text{(by \eqref{sigmaKsneq})}\\
&\qquad = q^2 \int\limits_K \pi(k) \pi(
\begin{bsmallmatrix}
\varpi&&&\\
&1&&\\
&&1&\\
&&&\varpi^{-1} 
\end{bsmallmatrix}
\begin{bsmallmatrix}
1&&&\\
&1&&\\
&&\varpi&\\
&&&\varpi
\end{bsmallmatrix} 
) v \, dk\\
&\qquad\quad + q^3 \int\limits_K \int\limits_{\OF} \pi(k) \pi(
\begin{bsmallmatrix}
\varpi&&&\\
&1&&\\
&&1&\\
&&&\varpi^{-1} 
\end{bsmallmatrix}
\begin{bsmallmatrix}
1&&&\\
&1&c&\\
&&1&\\
&&&1
\end{bsmallmatrix}
\begin{bsmallmatrix}
1&&&\\
&\varpi&&\\
&&1&\\
&&&\varpi
\end{bsmallmatrix}
)v\, dc\, dk\qquad\!\text{(by \eqref{thetadefeq3})}\\
&\qquad = q^2 \int\limits_K \pi(k) \pi(
\begin{bsmallmatrix}
\varpi&&&\\
&1&&\\
&&\varpi&\\
&&&1
\end{bsmallmatrix} 
) v \, dk
 + q^3 \int\limits_K  \pi(k) \pi(
\begin{bsmallmatrix}
\varpi&&&\\
&\varpi&&\\
&&1&\\
&&&1
\end{bsmallmatrix}
)v\,  dk\\
&\qquad = ( q^2 +q^3)\int\limits_K  \pi(k) \pi(
\begin{bsmallmatrix}
\varpi&&&\\
&\varpi&&\\
&&1&\\
&&&1
\end{bsmallmatrix}
)v\,  dk\\
&\qquad = T_{0,1}^s v \qquad \text{(by \eqref{Heckecommutelemmaeq20})}.
\end{align*}
A similar proof shows that $T_{1,0}^s v = q^4 \sigma_n \tau_n v$. 
Assume that $n \geq 2$. 
Then
$\sigma_n \tau_n v=\tau_{n-1} \sigma_{n-1} v$ by \ref{sigmaopslemmaitem6} of Lemma \ref{sigmaopslemma}.
Finally, we have $T_{1,0}^s v =0$ if and only if ${\sigma_{n-1} v=0}$ by $T_{1,0}^s v = q^4 \sigma_n \tau_n v$ and Lemma \ref{sigmaendolemma}.
\end{proof}

\section{Commutation relations}

In this section we derive commutation relations for the operators we defined in the previous sections. These
relations will play an important role in some of the subsequent chapters. 

\begin{lemma}\label{commlemma}
 Let $(\pi,V)$ be a smooth representation of $\GSp(4,F)$ for which the center acts trivially. Let $n$ be an integer such that $n\geq 0$, and $v\in V_s(n)$. Then
 \begin{align}
  \label{commreleq1}\theta\tau_nv&=\tau_{n+1}\theta v,\\
  \label{commreleq2}\theta t_nv&=t_{n+1}\theta v,\\
  \label{commreleq4}\theta\rho'_nv&=\rho'_{n+1}\theta v,\\
  \label{commreleq5}\theta p_nv&=p_{n+1}\theta v.
 \end{align}
\end{lemma}
\begin{proof}
It is easy to verify these statements using the defining formulas.
\end{proof}

Let $(\pi,V)$ be a smooth representation of $\GSp(4,F)$ for which the center acts trivially. For $v \in V_s(1)$ define
\begin{equation}\label{eoperatoreq}
 e(v)=\int\limits_\OF\int\limits_\OF\int\limits_\OF\pi(\begin{bsmallmatrix}1&x\varpi^{-1}&y\varpi^{-1}&z\varpi^{-1}\\&1&&y\varpi^{-1}\\&&1&-x\varpi^{-1}\\&&&1\end{bsmallmatrix})v\,dx\,dy\,dz.
\end{equation}
It is easy to verify that $e(v) \in V_s(2)$  for $v \in V_s(1)$, so that this defines a linear map $e:V_s(1)\to V_s(2)$.

\begin{lemma}\label{eTs10identitylemma}
 Let $(\pi,V)$ be a smooth representation for which the center acts trivially.
 Assume that $v\in V_s(1)$ is such that the paramodularization $p_1(T^s_{1,0}v)$ is zero. Then
 \begin{equation}\label{eTs10identityeq}
  e(v)=q^{-4}\eta T^s_{1,0}v-q^{-2}\tau_2T^s_{1,0}v+q^{-2}\tau_1T^s_{1,0}v.
 \end{equation}
\end{lemma}
\begin{proof}
By \eqref{rhopoplemmaeq3},
$$
 \rho'_1T^s_{1,0}v=q^{-1}\eta T^s_{1,0}v-q\tau_2T^s_{1,0}v
$$
is in $V(2)$, and in particular in $V_s(2)$. Hence both sides of \eqref{eTs10identityeq} are in $V_s(2)$. By Lemma \ref{taunlevelraisinglemma}, to prove \eqref{eTs10identityeq}, it is enough to prove that
\begin{equation}\label{eTs10identityeq2}
 \tau_2e(v)=\tau_2\big(q^{-4}\eta T^s_{1,0}v-q^{-2}\tau_2T^s_{1,0}v+q^{-2}\tau_1T^s_{1,0}v\big).
\end{equation}
This is equivalent to
\begin{equation}\label{eTs10identityeq3}
 \tau_2e(v)=q^{-4}\tau_2\eta T^s_{1,0}v.
\end{equation}
It follows from the definitions \eqref{eoperatoreq} and \eqref{Ts10eq} that \eqref{eTs10identityeq3} holds.
\end{proof}

\begin{lemma}\label{Ts01thetalemma}
Let $(\pi,V)$ be a smooth representation of $\GSp(4,F)$ for which the center acts trivially,
and let $n$ be an integer. Let $v \in V_s(n)$. 
Then
 \begin{align}
  \label{T01taueq}T^s_{0,1}\tau_nv&=\tau_nT^s_{0,1}v\qquad\text{if }n\geq1,\\
  \label{T10taueq}T^s_{1,0}\tau_nv&=\tau_nT^s_{1,0}v\qquad\text{if }n\geq1,\\
  \label{Ts10thetaeq}T^s_{1,0}\theta v&=q^2T^s_{0,1}\tau_nv\qquad\text{if }n\geq0,\\
  \label{Ts01thetalemmaeq1}T^s_{0,1}\theta v&=\theta T^s_{0,1}v+q^3\tau_nv-q^3\eta\sigma_{n-1}v\qquad\text{if }n\geq2,\\
  \label{Ts01thetalemmaeq1b}T^s_{0,1}\theta v&=\theta T^s_{0,1}v+q^3\tau_1 v-q^3e(v)\qquad\text{if }n=1,\\
  \label{Ts10thetapeq}T^s_{1,0}\theta'v&=q^4v+qT^s_{1,0}\tau_nv,\qquad\text{if }n\geq1\text{ and }v\in V(n),\\
  \label{T01rhoeq}T^s_{0,1}\rho_n'v&=q\theta v-q\tau_n T^s_{0,1}v\qquad\text{if }n\geq1\text{ and }p_n(v)=0,\\
  \label{T01rhoeqb}T^s_{0,1}\rho_n'v&=q\theta v+\tau_n T^s_{0,1}v\qquad\text{if }n\geq1\text{ and }v\in V(n),\\
  \label{T10rhoeq}T^s_{1,0}\rho_n'v&=q^3\tau_nv-q\tau_n T^s_{1,0}v\qquad\text{if }n\geq1\text{ and }p_n(v)=0,\\
  \label{T10rhoeqb}T^s_{1,0}\rho_n'v&=q^3\tau_nv+\tau_n T^s_{1,0}v\qquad\text{if }n\geq1\text{ and }v\in V(n).
 \end{align}
In these equations, if $T_{0,1}^s$ or $T_{1,0}^s$ appears on the left-hand side, then
this operator is an endomorphism of $V_s(n+1)$, and if $T_{0,1}^s$ or $T_{1,0}^s$ appears on the right-hand side,
then this operator is an endomorphism of $V_s(n)$.
\end{lemma}
\begin{proof}
Equations \eqref{T01taueq}, \eqref{T10taueq}, and \eqref{Ts10thetaeq}  follow from straightforward calculations. 
To prove \eqref{Ts01thetalemmaeq1} and \eqref{Ts01thetalemmaeq1b}, we begin by proving a preliminary identity.
Let $K=\Ks{n+1}$, and let $dk$ be the Haar measure on $K$ that assigns $K$ volume $1$. Let $v \in V_s(n)$; as before,
$n$ is an integer such that $n \geq 1$. 
Then:
\begin{align*}
&
\int\limits_{K}\int\limits_{\OF}
\pi(k
\begin{bsmallmatrix}
1&y\varpi^{-1}&&\\
&1&&\\
&&1&-y\varpi^{-1}\\
&&&1
\end{bsmallmatrix})v\, dy\, dk\\
&\qquad=
\int\limits_{K}\int\limits_{\OF}\int\limits_{\OF}
\pi(k
\begin{bsmallmatrix}
1&&&\vphantom{\varpi^{-n+1}}\\
&1&b&\vphantom{\varpi^{-n+1}}\\
&&1&\vphantom{\varpi^{-n+1}}\\
&&&1\vphantom{\varpi^{-n+1}}
\end{bsmallmatrix}
\begin{bsmallmatrix}
1&y\varpi^{-1}&&\\
&1&&\\
&&1&-y\varpi^{-1}\\
&&&1
\end{bsmallmatrix})v\, dy\, db\, dk\\
&\qquad=
\int\limits_{K}\int\limits_{\OF}\int\limits_{\OF}
\pi(k
\begin{bsmallmatrix}
1&y\varpi^{-1}&&\\
&1&&\\
&&1&-y\varpi^{-1}\\
&&&1
\end{bsmallmatrix}
\begin{bsmallmatrix}
1&&-yb\varpi^{-1}&y^2b\varpi^{-2}\\
&1&b&-yb\varpi^{-1}\\
&&1&\\
&&&1
\end{bsmallmatrix}
)v\, dy\, db\, dk\\
&\qquad=
\int\limits_{K}\int\limits_{\OF}\int\limits_{\OF^\times}
\pi(k
\begin{bsmallmatrix}
1&y\varpi^{-1}&&\\
&1&&\\
&&1&-y\varpi^{-1}\\
&&&1
\end{bsmallmatrix}
\begin{bsmallmatrix}
1&&-yb\varpi^{-1}&y^2b\varpi^{-2}\\
&1&&-yb\varpi^{-1}\\
&&1&\\
&&&1
\end{bsmallmatrix}
)v\, dy\, db\, dk\\
&\qquad\quad+
\int\limits_{K}\int\limits_{\OF}\int\limits_{\p}
\pi(k
\begin{bsmallmatrix}
1&y\varpi^{-1}&&\\
&1&&\\
&&1&-y\varpi^{-1}\\
&&&1
\end{bsmallmatrix}
\begin{bsmallmatrix}
1&&-yb\varpi^{-1}&y^2b\varpi^{-2}\\
&1&&-yb\varpi^{-1}\\
&&1&\\
&&&1
\end{bsmallmatrix}
)v\,  dy\, db\, dk\\
&\qquad=
\int\limits_{K}\int\limits_{\OF}\int\limits_{\OF^\times}
\pi(k
\begin{bsmallmatrix}
1&y\varpi^{-1}&&\\
&1&&\\
&&1&-y\varpi^{-1}\\
&&&1
\end{bsmallmatrix}
\begin{bsmallmatrix}
1&&-b\varpi^{-1}&yb\varpi^{-2}\\
&1&&-b\varpi^{-1}\\
&&1&\\
&&&1
\end{bsmallmatrix}
)v\, dy\, db\, dk\\
&\qquad\quad+q^{-1}
\int\limits_{K}
\pi(k
)v\, dk\\
&\qquad=
\int\limits_{K}\int\limits_{\OF}\int\limits_{\OF}
\pi(k
\begin{bsmallmatrix}
1&y\varpi^{-1}&&\\
&1&&\\
&&1&-y\varpi^{-1}\\
&&&1
\end{bsmallmatrix}
\begin{bsmallmatrix}
1&&-b\varpi^{-1}&yb\varpi^{-2}\\
&1&&-b\varpi^{-1}\\
&&1&\\
&&&1
\end{bsmallmatrix}
)v\, dy\, db\, dk\\
&\qquad\quad-
\int\limits_{K}\int\limits_{\OF}\int\limits_{\p}
\pi(k
\begin{bsmallmatrix}
1&y\varpi^{-1}&&\\
&1&&\\
&&1&-y\varpi^{-1}\\
&&&1
\end{bsmallmatrix}
\begin{bsmallmatrix}
1&&-b\varpi^{-1}&yb\varpi^{-2}\\
&1&&-b\varpi^{-1}\\
&&1&\\
&&&1
\end{bsmallmatrix}
)v\, dy\, db\, dk\\
&\qquad\quad+q^{-1}
\int\limits_{K}
\pi(k
)v\,dk\\
&\qquad=
\int\limits_{K}\int\limits_{\OF}\int\limits_{\OF}
\pi(k
\begin{bsmallmatrix}
1&y\varpi^{-1}&-b\varpi^{-1}&\\
&1&&-b\varpi^{-1}\\
&&1&-y\varpi^{-1}\\
&&&1
\end{bsmallmatrix}
)v\, dy\, db\, dk\\
&\qquad\quad-q^{-1}
\int\limits_{K}\int\limits_{\OF}
\pi(k
\begin{bsmallmatrix}
1&&-b\varpi^{-1}&\\
&1&&-b\varpi^{-1}\\
&&1&\\
&&&1
\end{bsmallmatrix}
)v\,  db\, dk
+q^{-1}
\int\limits_{K}
\pi(k
)v\,  dk\\
&\qquad=
\int\limits_{K}\int\limits_{\OF}\int\limits_{\OF}\int\limits_{\OF}
\pi(k
\begin{bsmallmatrix}
1&y\varpi^{-1}&b\varpi^{-1}&z\varpi^{-n}\\
&1&&b\varpi^{-1}\\
&&1&-y\varpi^{-1}\\
&&&1
\end{bsmallmatrix}
)v\, dy\, db\, dk\\
&\qquad\quad-q^{-1}
\int\limits_{K}\int\limits_{\OF}
\pi(k
\begin{bsmallmatrix}
1&b\varpi^{-1}&&\\
&1&&\\
&&1&-b\varpi^{-1}\\
&&&1
\end{bsmallmatrix}
)v\,  db\, dk
+q^{-1}
\int\limits_{K}
\pi(k
)v\,  dk.
\end{align*}
We have proven that
\begin{multline}
\label{intideq}
(q+1)
\int\limits_{K}\int\limits_{\OF}
\pi(k
\begin{bsmallmatrix}
1&y\varpi^{-1}&&\\
&1&&\\
&&1&-y\varpi^{-1}\\
&&&1
\end{bsmallmatrix}
)v\,  dy\, dk\\
=\int\limits_{K}
\pi(k
)v\,  dk
+q\int\limits_{K}\int\limits_{\OF}\int\limits_{\OF}\int\limits_{\OF}
\pi(k
\begin{bsmallmatrix}
1&y\varpi^{-1}&b\varpi^{-1}&z\varpi^{-n}\\
&1&&b\varpi^{-1}\\
&&1&-y\varpi^{-1}\\
&&&1
\end{bsmallmatrix}
)v\, dy\, db\, dk.
\end{multline}
Next, as in the proof of Lemma \ref{Heckecommutelemma}, we have
$$
T^s_{0,1}\theta v = (q^2+q^3)\int\limits_K\pi(k
\begin{bsmallmatrix}\varpi\\&\varpi\\&&1\\&&&1\end{bsmallmatrix})\theta v\, dk.
$$
Hence, by \eqref{thetadefeq3},
\begin{align*}
T^s_{0,1}\theta v
&=(q^2+q^3)\int\limits_K\pi(k
\begin{bsmallmatrix}\varpi\\&\varpi\\&&1\\&&&1\end{bsmallmatrix}
\begin{bsmallmatrix}1&&&\\&1&&\\&&\varpi&\\&&&\varpi\end{bsmallmatrix})v\, dk\\
&\quad+q (q^2+q^3)\int\limits_K\int\limits_{\OF}\pi(k 
\begin{bsmallmatrix}\varpi\\&\varpi\\&&1\\&&&1\end{bsmallmatrix}
\begin{bsmallmatrix}1&&&\\&1&c&\\&&1&\\&&&1\end{bsmallmatrix}
\begin{bsmallmatrix}1&&&\\&\varpi&&\\&&1&\\&&&\varpi\end{bsmallmatrix})v\,dc\,dk\\
&=(q^2+q^3)\int\limits_K\pi(k)v\, dk
+(q^3+q^4)\int\limits_K\pi(k 
\begin{bsmallmatrix}\varpi&&&\\&\varpi^2&&\\&&1&\\&&&\varpi\end{bsmallmatrix})v\,dk.
\end{align*}
We also have, by \eqref{thetaineq}, and then \eqref{Ts01eq}, 
\begin{align*}
&\theta T_{0,1}^s v 
= (q+1)\int\limits_K 
\pi(k
\begin{bsmallmatrix}
1&&&\\
&\varpi&&\\
&&1&\\
&&&\varpi
\end{bsmallmatrix})T_{0,1}^s v \, dk\\
&\qquad=\sum_{y,z\in\OF/\p}
(q+1)\int\limits_K 
\pi(k
\begin{bsmallmatrix}
1&&&\vphantom{\varpi^{-n+1}}\\
&\varpi&&\vphantom{\varpi^{-n+1}}\\
&&1&\vphantom{\varpi^{-n+1}}\\
&&&\varpi\vphantom{\varpi^{-n+1}}
\end{bsmallmatrix}
\begin{bsmallmatrix}1&y&&z\varpi^{-n+1}\\&1&&\\&&1&-y\\&&&1\end{bsmallmatrix}
\begin{bsmallmatrix}
\varpi&&&\vphantom{\varpi^{-n+1}}\\
&1&&\vphantom{\varpi^{-n+1}}\\
&&\varpi&\vphantom{\varpi^{-n+1}}\\
&&&1\vphantom{\varpi^{-n+1}}
\end{bsmallmatrix})v\,dk\\
&\qquad\quad+\sum_{c,y,z\in\OF/\p}(q+1)\int\limits_K 
\pi(k
\begin{bsmallmatrix}
1&&&\vphantom{\varpi^{-n+1}}\\
&\varpi&&\vphantom{\varpi^{-n+1}}\\
&&1&\vphantom{\varpi^{-n+1}}\\
&&&\varpi\vphantom{\varpi^{-n+1}}
\end{bsmallmatrix}
\begin{bsmallmatrix}1&&y&z\varpi^{-n+1}\\&1&c&y\\&&1\\&&&1\end{bsmallmatrix}
\begin{bsmallmatrix}
\varpi&&&\vphantom{\varpi^{-n+1}}\\
&\varpi&&\vphantom{\varpi^{-n+1}}\\
&&1&\vphantom{\varpi^{-n+1}}\\
&&&1\vphantom{\varpi^{-n+1}}
\end{bsmallmatrix})v\, dk\\
&\qquad=
(q^3+q^2)\int\limits_K \int\limits_\OF
\pi(k
\begin{bsmallmatrix}1&y\varpi^{-1}&&\\&1&&\\&&1&-y\varpi^{-1}\\&&&1\end{bsmallmatrix}
)v\,dy\,dk\\
&\qquad\quad+(q^4+q^3)\int\limits_K 
\pi(k
\begin{bsmallmatrix}\varpi\\&\varpi^2\\&&1\\&&&\varpi\end{bsmallmatrix})v\, dk.
\end{align*}
Hence, by these formulas for $T^s_{0,1}\theta v$, $\theta T^s_{0,1}v$, and \eqref{intideq},
\begin{align*}
&T^s_{0,1}\theta v -\theta T^s_{0,1} v\\
&\qquad=(q^2+q^3)\int\limits_K\pi(k)v\, dk
-(q^3+q^2)\int\limits_K \int\limits_\OF
\pi(k
\begin{bsmallmatrix}1&y\varpi^{-1}&&\\&1&&\\&&1&-y\varpi^{-1}\\&&&1\end{bsmallmatrix}
)v\,dy\,dk\\
&\qquad=q^3\int\limits_K\pi(k)v\, dk
-q^3
\int\limits_{K}\int\limits_{\OF}\int\limits_{\OF}\int\limits_{\OF}
\pi(k
\begin{bsmallmatrix}
1&y\varpi^{-1}&b\varpi^{-1}&z\varpi^{-n}\\
&1&&b\varpi^{-1}\\
&&1&-y\varpi^{-1}\\
&&&1
\end{bsmallmatrix}
)v\, dy\, db\, dk.
\end{align*}
The first term on the right hand side equals $q^3\tau_nv$ by \eqref{tauinteq}.
This proves \eqref{Ts01thetalemmaeq1} and \eqref{Ts01thetalemmaeq1b} since
$$
\int\limits_{K}\int\limits_{\OF}\int\limits_{\OF}\int\limits_{\OF}
\pi(k
\begin{bsmallmatrix}
1&y\varpi^{-1}&b\varpi^{-1}&z\varpi^{-n}\\
&1&&b\varpi^{-1}\\
&&1&-y\varpi^{-1}\\
&&&1
\end{bsmallmatrix}
)v\, dy\, db\, dk
=
\begin{cases}
\eta\sigma_{n-1} v &\text{if $n\geq 2$},\\
e(v)&\text{if $n=1$},
\end{cases}
$$
by \eqref{sigmaopseq} and \eqref{eoperatoreq}.

Finally, to complete the proof, we note that:  \eqref{Ts10thetapeq} follows  from \eqref{Ts10eq} and \eqref{thetap1defeq};
 \eqref{T01rhoeq} follows  from \eqref{rhopoplemmaeq3};
 \eqref{T01rhoeqb} follows  from \eqref{rhopoplemmaeq4};
 \eqref{T10rhoeq} follows  from \eqref{rhopoplemmaeq3};
and  \eqref{T10rhoeqb} follows  from \eqref{rhopoplemmaeq4}.
\end{proof}

\section{A result about eigenvalues}
Let $(\pi,V)$ be a smooth representation of $\GSp(4,F)$ for which the center acts trivially.
In the following lemma we prove results that show that the spectrum of the endomorphism $T_{1,0}^s$ of $V_s(n)$ for $n\geq1$ is contained in the union of $\{0\}$ and the spectrum of the endomorphism $T_{1,0}^s$ of $V_s(1)$.
\begin{lemma}\label{Ts10eigenvalueslemma}
 Let $(\pi,V)$ be a smooth representation of $\GSp(4,F)$ for which the center acts trivially.
Let $n$ be an integer such that $n \geq 1$.
 \begin{enumerate}
  \item \label{Ts10eigenvalueslemmaitem1} Assume that $V_s(1)=0$.  
  Then the operator $(T^s_{1,0})^{n-1}$ is zero on $V_s(n)$.
  \item \label{Ts10eigenvalueslemmaitem2} Assume that $T^s_{1,0}$ acts on $V_s(1)$ by multiplication by $c\in\C$. Then the operator $(T^s_{1,0}-c\cdot{\rm id})(T^s_{1,0})^{n-1}$ is zero on $V_s(n)$. Moreover, if $n \geq 2$, then 
   \begin{equation}\label{Ts10eigenvalueslemmaeq5}
    {\rm im}((T^s_{1,0})^k)=\tau_{n-1}V_s(n-k)
   \end{equation}
   and $\dim({\rm im}((T^s_{1,0})^k))=\dim(V_s(n-k))$ for $1 \leq k \leq n-1$. 
  \item \label{Ts10eigenvalueslemmaitem3} Let $\{\mu_1,\ldots,\mu_r\}$ be the set of non-zero eigenvalues of $T^s_{1,0}$ on $V_s(1)$, and let $\{\mu_1',\ldots,\mu'_{r'}\}$ be the set of non-zero eigenvalues of $T^s_{1,0}$ on $V_s(n)$. Then
   \begin{equation}\label{Ts10eigenvalueslemmaeq6}
    \{\mu_1,\ldots,\mu_r\}=\{\mu_1',\ldots,\mu'_{r'}\}.
   \end{equation}
 \end{enumerate}
\end{lemma}
\begin{proof}
\ref{Ts10eigenvalueslemmaitem1}. We may assume that $n\geq2$. Let $v \in V_s(n)$. It follows from \eqref{sigmaopseq} that
\begin{equation}
 \sigma_{n-k}\sigma_{n-k+1}\ldots\sigma_{n-1}v=\int\limits_\OF\int\limits_\OF\int\limits_\OF\pi(\begin{bsmallmatrix}1&x&y&z\varpi^{-n+k+1}\\&1&&y\\&&1&-x\\&&&1\end{bsmallmatrix}\eta^{-k})v\,dx\,dy\,dz
\end{equation}
for $1\leq k\leq n-1$. Also, by \eqref{Ts10eq},
\begin{equation}
 (T^s_{1,0})^kv=q^{4k}\int\limits_\OF\int\limits_\OF\int\limits_\OF\pi(\begin{bsmallmatrix}1&x&y&z\varpi^{-n+1}\\&1\\&&1\\&&&1\end{bsmallmatrix}\eta^{-k})v\,dx\,dy\,dz.
\end{equation}
It follows that
\begin{equation}\label{Ts10eigenvalueslemmaeq1}
 q^{4k}\tau_{n-1}\sigma_{n-k}\sigma_{n-k+1}\ldots\sigma_{n-1}v=(T^s_{1,0})^kv
\end{equation}
for $1\leq k \leq n-1$. 
For $k=n-1$,
\begin{equation}\label{Ts10eigenvalueslemmaeq2}
 q^{4n-4}\tau_{n-1}\sigma_1\sigma_2\ldots\sigma_{n-1}v=(T^s_{1,0})^{n-1}v.
\end{equation}
Since $\sigma_1\sigma_2\ldots\sigma_{n-1}v\in V_s(1)=0$, we conclude that $(T^s_{1,0})^{n-1}v=0$. 

\ref{Ts10eigenvalueslemmaitem2}. We may assume that $n\geq 2$. Let $v \in V_s(n)$. By applying $T^s_{1,0}-c\cdot{\rm id}$ to \eqref{Ts10eigenvalueslemmaeq2} we obtain
\begin{equation}\label{Ts10eigenvalueslemmaeq3}
  q^{4n-4}(T^s_{1,0}-c\cdot{\rm id})\tau_{n-1}\sigma_1\sigma_2\ldots\sigma_{n-1}v=(T^s_{1,0}-c\cdot{\rm id})(T^s_{1,0})^{n-1}v.
\end{equation}
By \eqref{T10taueq},
\begin{equation}\label{Ts10eigenvalueslemmaeq4}
 q^{4n-4} \tau_{n-1}(T^s_{1,0}-c\cdot{\rm id})\sigma_1\sigma_2\ldots\sigma_{n-1}v=(T^s_{1,0}-c\cdot{\rm id})(T^s_{1,0})^{n-1}v.
\end{equation}
By hypothesis $T^s_{1,0}-c\cdot{\rm id}$ acts by zero on $V_s(1)$. 
This proves that $(T^s_{1,0}-c\cdot{\rm id})(T^s_{1,0})^{n-1}$ is zero on $V_s(n)$. 
Equation \eqref{Ts10eigenvalueslemmaeq5} follows from \eqref{Ts10eigenvalueslemmaeq1} and 
Proposition~\ref{sigmasurjectiveprop}. Lemma~\ref{taunlevelraisinglemma} \ref{taunlevelraisinglemmaitem1} implies the last assertion.

\ref{Ts10eigenvalueslemmaitem3}. We may assume that $n\geq 2$. Let $\mu \in \{\mu'_1,\dots,\mu'_{r'}\}$. Let $v \in V_s(n)$ be a non-zero vector
such that  $T^s_{1,0}v=\mu v$. By \eqref{Ts10eigenvalueslemmaeq2}, the vector
$$
 v_1:=\sigma_1\sigma_2\ldots\sigma_{n-1}v\in V_s(1)
$$
is non-zero. We have
\begin{align*}
 \tau_{n-1}(T^s_{1,0}-\mu\cdot{\rm id}) v_1&=(T^s_{1,0}-\mu\cdot{\rm id})\tau_{n-1}\tau_{n-2}\cdots \tau_1
 v_1\qquad \text{(by \eqref{T10taueq})}\\
 &=(T^s_{1,0}-\mu\cdot{\rm id})\tau_{n-1}\tau_{n-2}\cdots \tau_1\sigma_1\sigma_2\ldots\sigma_{n-1}v\\
 &=q^{4-4n}(T^s_{1,0}-\mu\cdot{\rm id})(T^s_{1,0})^{n-1}v\qquad \text{(by \eqref{Ts10eigenvalueslemmaeq2})}\\
 &=q^{4-4n}(T^s_{1,0})^{n-1}(T^s_{1,0}-\mu\cdot{\rm id})v\\
 &=0.
\end{align*}
Lemma \ref{taunlevelraisinglemma} \ref{Ts10eigenvalueslemmaitem1} now implies that $(T^s_{1,0}-\mu\cdot{\rm id}) v_1=0$. Hence $\mu$ is a non-zero eigenvalue for $T^s_{1,0}$ on $V_s(1)$, so that $\mu \in \{\mu_1,\dots,\mu_r\}$. 

Assume now that $\mu \in \{\mu_1,\dots,\mu_r\}$. Let $v \in V_s(1)$ be a non-zero vector
such that  $T^s_{0,1}v=\mu v$. Then by \eqref{T10taueq},
\begin{align*}
 T^s_{1,0}\tau_{n-1}v=\tau_{n-1}T^s_{1,0}v=\mu\tau_{n-1}v.
\end{align*}
The vector $\tau_{n-1}v\in V_s(n)$ is non-zero by \ref{taunlevelraisinglemmaitem1} of Lemma \ref{taunlevelraisinglemma}. Hence $\mu$ is a non-zero eigenvalue for $T^s_{1,0}$ on $V_s(n)$, so that $\mu \in \{\mu'_1,\dots,\mu'_{r'}\}$. 
\end{proof}

%% file: SKMS_chapter4.tex
\chapter{Some Induced Representations}
\label{inducedrepchap}
In this chapter we prove a number of important results concerning stable
Klingen vectors in representations of $\GSp(4,F)$ induced from the Siegel parabolic
subgroup $P$ of $\GSp(4,F)$. These results are key ingredients for the dimension formulas proved in Chap.~\ref{dimchap}.

\section{Double coset representatives}
We begin by determining a system of representatives for $P \backslash \GSp(4,F)/\Ks{n}$
for an integer $n$ such that $n \geq 1$.  
For an integer $i$, let
\begin{equation}\label{Midefeq}
 M_i=\begin{bsmallmatrix}1\\&1\\\varpi^i&&1\\&\varpi^i&&1\end{bsmallmatrix}.
\end{equation}
Let $t_n$ be the element of $\GSp(4,F)$ defined in \eqref{tndefeq}. The 
following lemma is the main result of this section. 
\begin{lemma}\label{PKsdisjointlemma1}
Let $n$ be an integer such that $n\geq 1$. Then
 \begin{align}\label{PKsdisjointeq6}
  \GSp(4,F)&=P\Ks{n}\;\sqcup\bigsqcup_{0<i\leq n/2}P M_i\Ks{n}\nonumber\\
  &\quad\sqcup\;Pt_n\Ks{n}\;\sqcup\bigsqcup_{0<i<n/2}P M_it_n\Ks{n}.
 \end{align}
\end{lemma}
\begin{proof}
By Proposition 5.1.2 of \cite{NF} we have
\begin{equation}\label{NF51eq1}
 \GSp(4,F)=P\K{n}\sqcup\bigsqcup_{0<i\leq n/2}P M_i\K{n}.
\end{equation}
Using \eqref{paramodularKsdecompositioneq2}, it follows that
\begin{align*}
 \GSp(4,F)&=\bigcup\limits_{v\in\OF/\p}P\begin{bsmallmatrix}1&&&v\varpi^{-n}\\&1\\&&1\\&&&1\end{bsmallmatrix}\Ks{n}\;\cup\;P t_n\Ks{n}\\
&\quad\cup\bigcup_{0<i\leq n/2}\,\bigcup\limits_{v\in\OF/\p}P M_i\begin{bsmallmatrix}1&&&v\varpi^{-n}\\&1\\&&1\\&&&1\end{bsmallmatrix}\Ks{n}\\
&\quad\cup\bigcup_{0<i\leq n/2}P  M_it_n\Ks{n}\\
&=P \Ks{n}\;\cup\;P t_n\Ks{n}\\
&\quad\cup\bigcup_{0<i\leq n/2}\,\bigcup\limits_{v\in\OF/\p}P M_i\begin{bsmallmatrix}1&&&v\varpi^{-n}\\&1\\&&1\\&&&1\end{bsmallmatrix}\Ks{n}\\
&\quad\cup\bigcup_{0<i\leq n/2}P M_it_n\Ks{n}.
\end{align*}
Since 
$
\begin{bsmallmatrix}
v&&&\\
&1&&\\
&&v&\\
&&&1
\end{bsmallmatrix}
$
for $v \in \OF^\times$ commutes with $M_i$, we have
\begin{align*}
 \GSp(4,F)&=\bigcup_{0<i\leq n/2}P M_i\begin{bsmallmatrix}1&&&\varpi^{-n}\\&1\\&&1\\&&&1\end{bsmallmatrix}\Ks{n}\\
  &\quad\cup\bigg(P\Ks{n}\;\cup\bigcup_{0<i\leq n/2}P M_i\Ks{n}\bigg)\\
  &\quad\cup\bigg(Pt_n\Ks{n}\;\cup\bigcup_{0<i\leq n/2}P M_it_n\Ks{n}\bigg).
\end{align*}
From the disjointness of \eqref{NF51eq1} it follows that
\begin{align}\label{PKsdisjointeq3}
 \GSp(4,F)&=\bigsqcup_{0<i\leq n/2}P M_i\begin{bsmallmatrix}1&&&\varpi^{-n}\\&1\\&&1\\&&&1\end{bsmallmatrix}\Ks{n}\nonumber\\
  &\quad\cup\bigg(P\Ks{n}\;\sqcup\bigsqcup_{0<i\leq n/2}P M_i\Ks{n}\bigg)\nonumber\\
  &\quad\cup\bigg(Pt_n\Ks{n}\;\sqcup\bigsqcup_{0<i\leq n/2}P M_it_n\Ks{n}\bigg).
\end{align}
If $n$ is even, then there are additional identities
\begin{align}\label{PKsdisjointeq11}
&M_{n/2}
\begin{bsmallmatrix}
\vphantom{-\varpi^n}1&&&\varpi^{-n}\\
&\vphantom{-\varpi^n}1\\
&&\vphantom{-\varpi^n}1\\
&&&\vphantom{-\varpi^n}1
\end{bsmallmatrix}
\begin{bsmallmatrix}
\vphantom{\varpi^n}-1&&&\\&\vphantom{\varpi^n}1&&\\&&\vphantom{\varpi^n}-1&\\&&&\vphantom{\varpi^n}1
\end{bsmallmatrix}
\begin{bsmallmatrix}1\\\varpi^n&1&1+\varpi^{n/2}&\\-\varpi^n&-1&-\varpi^{n/2}\\\varpi^n(1-\varpi^{n/2})&&-\varpi^n&1
\end{bsmallmatrix}\nonumber\\
 &\qquad=
\begin{bsmallmatrix}
\vphantom{\varpi^n}-1&&&\\&\vphantom{\varpi^n}1&&\\&&\vphantom{\varpi^n}-1&\\&&&\vphantom{\varpi^n}1
\end{bsmallmatrix}
 \begin{bsmallmatrix}&\varpi^{-n/2}\\-\varpi^{n/2}&1\\&&&-\varpi^{-n/2}\\&&\varpi^{n/2}&1\end{bsmallmatrix}
 \begin{bsmallmatrix}1&&-\varpi^{-n/2}&-\varpi^{-n}\\&1&\varpi^{n/2}&-\varpi^{-n/2}\\&&1\\&&&1\end{bsmallmatrix}M_{n/2}
\end{align}
and
\begin{align}\label{PKsdisjointeq12}
&\begin{bsmallmatrix}&-\varpi^{-n/2}\\\varpi^{n/2}\\&&&\varpi^{-n/2}\\&&-\varpi^{n/2}\end{bsmallmatrix}
 \begin{bsmallmatrix}1&&-\varpi^{-n/2}\\&1&&-\varpi^{-n/2}\\&&1\\&&&1\end{bsmallmatrix}M_{n/2}t_n\nonumber\\
 &\qquad=M_{n/2}\begin{bsmallmatrix}1\\&&-1\\&1\\&&&1\end{bsmallmatrix}.
\end{align}
They show that if $n$ is even, then
\begin{equation}\label{PKsdisjointeq8b}
 P M_{n/2}\begin{bsmallmatrix}1&&&\varpi^{-n}\\&1\\&&1\\&&&1\end{bsmallmatrix}\Kl{n}=PM_{n/2}\Kl{n}=PM_{n/2}t_n\Kl{n}.
\end{equation}
Hence, if $n$ is even, then we get from \eqref{PKsdisjointeq3} that
\begin{align}\label{PKsdisjointeq4}
 \GSp(4,F)&=\bigsqcup_{0<i<n/2}P M_i\begin{bsmallmatrix}1&&&\varpi^{-n}\\&1\\&&1\\&&&1\end{bsmallmatrix}\Ks{n}\nonumber\\
  &\quad\cup\bigg(P\Ks{n}\;\sqcup\bigsqcup_{0<i\leq n/2}P M_i\Ks{n}\bigg)\nonumber\\
  &\quad\cup\bigg(Pt_n\Ks{n}\;\sqcup\bigsqcup_{0<i<n/2}P M_it_n\Ks{n}\bigg).
\end{align}
If $n$ is odd, then we can never have $i=n/2$, so that \eqref{PKsdisjointeq4} holds for both $n$ even and $n$ odd. The identity
\begin{align}\label{PKsdisjointeq13}
&M_i\begin{bsmallmatrix}
\vphantom{-\varpi^{-i}}1&&&\varpi^{-n}\\
&\vphantom{-\varpi^{-i}}1\\
&&\vphantom{-\varpi^{-i}}1\\
&&&\vphantom{-\varpi^{-i}}1
\end{bsmallmatrix}
\begin{bsmallmatrix}
\vphantom{\varpi^{-i}}-1&&&\\&\vphantom{\varpi^{-i}}1&&\\&&\vphantom{\varpi^{-i}}-1&\\&&&\vphantom{\varpi^{-i}}1
\end{bsmallmatrix}
\begin{bsmallmatrix}
1\\
&1+\varpi^{n-2i}&-\varpi^{n-2i}\\
&-1&1\\
\varpi^n+\varpi^{2n-2i}&&&1
\end{bsmallmatrix} \nonumber\\
&=
\begin{bsmallmatrix}
\vphantom{\varpi^{-i}}-1&&&\\&\vphantom{\varpi^{-i}}1&&\\&&\vphantom{\varpi^{-i}}-1&\\&&&\vphantom{\varpi^{-i}}1
\end{bsmallmatrix}
\begin{bsmallmatrix}
1&\varpi^{-i}\\
\varpi^{n-i}&1+\varpi^{n-2i}\\
&&1&-\varpi^{-i}\\
&&-\varpi^{n-i}&1+\varpi^{n-2i}
\end{bsmallmatrix}
\begin{bsmallmatrix}
1&&\varpi^{n-3i}&-\varpi^{n-4i}-\varpi^{-2i}\\&1&-\varpi^{n-2i}&\varpi^{n-3i}\\&&1\\&&&1
\end{bsmallmatrix}M_it_n
\end{align}
shows that
\begin{equation}\label{PKsdisjointeq7}
 PM_i\begin{bsmallmatrix}1&&&\varpi^{-n}\\&1\\&&1\\&&&1\end{bsmallmatrix}\Kl{n}=PM_it_n\Kl{n}
\end{equation}
for $0<i\leq n/2$. Hence we can omit the first line on the right hand side of \eqref{PKsdisjointeq4}, and get
\begin{align}\label{PKsdisjointeq5}
 \GSp(4,F)&=\bigg(P\Ks{n}\;\sqcup\bigsqcup_{0<i\leq n/2}P M_i\Ks{n}\bigg)\nonumber\\
  &\quad\cup\bigg(Pt_n\Ks{n}\;\sqcup\bigsqcup_{0<i<n/2}P M_it_n\Ks{n}\bigg).
\end{align}
It follows from \eqref{NF51eq1} and \eqref{PKsdisjointeq5} that the only possible equalities between the involved double cosets are
\begin{equation}\label{PKsdisjointlemma1eq0}
 PM_i\Ks{n}=P M_it_n\Ks{n}
\end{equation}
for some $i$ with $0<i<n/2$, or
\begin{equation}\label{PKsdisjointlemma1eq0b}
 P\Ks{n}=Pt_n\Ks{n}.
\end{equation}
Assume that $0<i<n/2$, and that
\begin{equation}\label{PKsdisjointlemma1eq1}
 M_ik=pM_it_n.
\end{equation}
for some $k\in\Ks{n}$ and $p\in P$.
Then
\begin{equation}\label{PKsdisjointlemma1eq2}
 M_i^{-1}pM_it_n=k\in\Ks{n}.
\end{equation}
Write
$$
 p=\begin{bsmallmatrix}A& \\ & uA'\end{bsmallmatrix} \begin{bsmallmatrix} 1& B\\ & \vphantom{u'}1\end{bsmallmatrix}
\quad
\mathrm{with}
\quad
A=\begin{bsmallmatrix} \vphantom{b_1}a_1& a_2\\ \vphantom{b_3}a_3 & a_4\end{bsmallmatrix} 
\quad \text{and}\quad
B=\begin{bsmallmatrix} b_1& b_2 \\ b_3 & b_1 \end{bsmallmatrix},
$$
where $a_1,a_2,a_3,a_4,b_1,b_2,b_3 \in F$ and $u \in F^\times$. 
Let $k_{ij}$, $1 \leq i,j \leq 4$, be the entries of the matrix $k$ in  \eqref{PKsdisjointlemma1eq2}. Calculations show that
\begin{equation}\label{PKsdisjointlemma1eq3}
 k_{14}+\varpi^{i-n}k_{13}=-a_1\varpi^{-n}.
\end{equation}
Since $k$ is in $\Ks{n}$, we obtain
\begin{equation}\label{PKsdisjointlemma1eq4}
 -a_1\varpi^{-n}\in\p^{-n+1}+\p^{i-n}=\p^{-n+1},
\end{equation}
which implies that
\begin{equation}\label{PKsdisjointlemma1eq4b}
 a_1\in\p.
\end{equation}
It follows from \eqref{PKsdisjointlemma1eq2} that
\begin{equation}\label{PKsdisjointlemma1eq5}
 p=\begin{bsmallmatrix} A& \\ & uA'\end{bsmallmatrix} \begin{bsmallmatrix} 1 & B\\ & \vphantom{uA'}1\end{bsmallmatrix} \in P\cap M_i\K{n}M_i^{-1}.
\end{equation}
By Lemma 5.2.1 of \cite{NF}, we have
\begin{equation}\label{PKsdisjointlemma1eq6}
 A=\begin{bsmallmatrix} a_1 & a_2 \\ a_3 & a_4\end{bsmallmatrix} \in\begin{bsmallmatrix} \OF& \p^{-i}\\ 
\p^{n-i}&\OF\end{bsmallmatrix} \quad\text{and}\quad\det(A)\in\OF^\times.
\end{equation}
Since $\p^{n-2i}\subset\p$ by assumption, it follows that $a_1,a_4\in\OF^\times$. Now we have a contradiction with \eqref{PKsdisjointlemma1eq4b}, proving that \eqref{PKsdisjointlemma1eq0} is impossible.

Now assume that \eqref{PKsdisjointlemma1eq0b} holds. Then $pt_n\in\Ks{n}$ for some $p\in P$. Hence also $pt_ns_2\in\Ks{n}$. Calculations show that $pt_ns_2$ is of the form
$
 \begin{bsmallmatrix} A& B\\ C& 0\end{bsmallmatrix}
$
for some $2\times 2$ matrices $A$, $B$, and $C$ with entries from $F$. 
We claim that $\Ks{n}$ does not contain any matrices of this kind. If $\begin{bsmallmatrix} A&B\\ C & 0\end{bsmallmatrix}\in\Ks{n}$, then $\det(B)\det(C)\in\OF^\times$. On the other hand, $\det(B)\in\p^{-n+1}$ and $\det(C)\in\p^n$ by \eqref{Ksshapeeq}. Hence $\det(B)\det(C)\in\p$, a contradiction. This proves our claim, and concludes the proof.
\end{proof}

\begin{remark}\label{KnKsnremark}
{\rm
It follows from \eqref{paramodularKsdecompositioneq1} that
\begin{equation}\label{PKsdisjointeq8}
 P\K{n}=P\Ks{n}\sqcup Pt_n\Ks{n}
\end{equation}
(the disjointness follows from Lemma \ref{PKsdisjointlemma1}). It follows from \eqref{paramodularKsdecompositioneq1} and \eqref{PKsdisjointeq7} that
\begin{equation}\label{PKsdisjointeq9}
 PM_i\K{n}=PM_i\Ks{n}\sqcup PM_it_n\Ks{n}\quad\text{for }0<i<n/2.
\end{equation}
It follows from \eqref{paramodularKsdecompositioneq1} and \eqref{PKsdisjointeq8b} that if $n$ is even, then 
\begin{equation}\label{PKsdisjointeq9b}
 PM_{n/2}\K{n}=PM_{n/2}\Ks{n}.
\end{equation}
}
\end{remark}
\section{Stable Klingen vectors in Siegel induced representations}
Let $(\pi,W)$ be an admissible representation of $\GL(2,F)$
admitting a central character~$\omega_\pi$, and let $\sigma$
be a character of $F^\times$. Assume that $\omega_\pi \sigma^2 =1$,
so that the induced representation $\pi \rtimes \sigma$ has 
trivial central character. Let $V$ be the standard model of $\pi \rtimes \sigma$.
In this section we will determine the stable Klingen level
of $\pi \rtimes \sigma$, and, in the case $\pi$ is irreducible,
the dimensions of the vector spaces $V_s(n)$ of stable Klingen
vectors for all integers $n \geq 0$. Explicitly, 
$V$ consists of the locally constant functions $f:\GSp(4,F)\rightarrow W$ such that
\begin{equation}\label{indreptrafoeq}
 f(\begin{bsmallmatrix} A&*\\ & uA'\end{bsmallmatrix}g)=|u^{-1}\det(A)|^{3/2}
 \sigma(u)\pi(A)f(g)
\end{equation}
for $\begin{bsmallmatrix} A&*\\ & uA'\end{bsmallmatrix} \in P$ and $g \in \GSp(4,F)$ (see Sect.~\ref{repsec}).
Consequently, if $n$ is an integer such that $n \geq 1$, and $f \in V_s(n)$, 
then by Lemma \ref{PKsdisjointlemma1} and 
\eqref{indreptrafoeq}, $f$ is determined by its values on the set
\begin{equation}
\label{indreptrafoeq021}
X=\{1,t_n\} \cup \{M_i \mid 0 <i\leq n/2\} \cup \{M_it_n \mid 0< i<n/2\}.
\end{equation}
The main result of this section is Theorem \ref{Siegelindleveltheorem}. 
We begin with some lemmas that will be useful in determining the elements
of $V_s(n)$ supported on a particular element of $X$. 

\begin{lemma}
\label{siegelprelemma}
Let $n$ be an integer such that $n \geq 0$. 
Let $g \in \GSp(4,F)$, $u \in F^\times$, and $A \in \GL(2,F)$. If
there exists $X \in \Mat(2,F)$ of the form
$X=\begin{bsmallmatrix}x&y\\z&x\end{bsmallmatrix}$
such that 
$\begin{bsmallmatrix} A& \\ & uA'\end{bsmallmatrix}
\begin{bsmallmatrix} \idm & X\\ & \vphantom{uA'}\idm\end{bsmallmatrix}
\in P \cap g \Ks{n} g^{-1}$, then $u \in \OF^\times$ and $\det(A) \in \OF^\times$.
\end{lemma}
\begin{proof}
Consider the homomorphism
$$
 P\cap g\Ks{n} g^{-1}\longrightarrow F^\times\times F^\times 
\quad
\text{defined by}
\quad
 \begin{bsmallmatrix} B& *\\ & vB'\end{bsmallmatrix} \longmapsto (v,\det(B)).
$$
Its image is a compact subgroup of $F^\times\times F^\times$, and is therefore
contained in $\OF^\times\times\OF^\times$. This implies that $u\in \OF^\times$ and
$\det(A)\in\OF^\times$.
\end{proof}

\begin{lemma}
\label{siegelpre2lemma}
Let $n$ be an integer such that $n\geq 0$. Let $u \in F^\times$
and $A \in \GL(2,F)$. The following statements are equivalent
\begin{enumerate}
\item \label{siegelpre2lemmaitem1} There exists an $X\in \Mat(2,F)$ of the form $X=\begin{bsmallmatrix} x& y \\ z & x\end{bsmallmatrix}$
such that
\begin{equation}
\begin{bsmallmatrix} A& \\ & uA'\end{bsmallmatrix} 
\begin{bsmallmatrix} \idm & X\\ & \vphantom{uA'}\idm\end{bsmallmatrix}
\in P\cap  \Ks{n}.
\end{equation}
\item \label{siegelpre2lemmaitem2} The following conditions are satisfied:
\subitem $\bullet$ $u\in\OF^\times$,
\subitem $\bullet$ $\det(A)\in\OF^\times$,
\subitem $\bullet$ $A\in\begin{bsmallmatrix} \OF& \OF\\ \p^{n} & \OF\end{bsmallmatrix}$.
\end{enumerate}
\end{lemma}
\begin{proof}
Assume that \ref{siegelpre2lemmaitem1} holds. Then $u \in \OF^\times$
and $\det(A) \in \OF^\times$ by Lemma \ref{siegelprelemma}. 
The assumption  $\begin{bsmallmatrix} A& \\ & uA'\end{bsmallmatrix} 
\begin{bsmallmatrix} \idm & X\\ & \vphantom{uA'}\idm\end{bsmallmatrix}
\in P\cap  \Ks{n}$, along with \eqref{Ksshapeeq}, implies that
$A\in\begin{bsmallmatrix} \OF& \OF\\ \p^{n} & \OF\end{bsmallmatrix}$.
Thus, \ref{siegelpre2lemmaitem2} holds. If \ref{siegelpre2lemmaitem2} holds, then \ref{siegelpre2lemmaitem1} holds with $X=0$.
\end{proof}

\begin{lemma}
\label{siegelpre3lemma}
Let $n$ be an integer such that $n\geq 0$. Let $u \in F^\times$
and $A \in \GL(2,F)$. The following statements are equivalent
\begin{enumerate}
\item \label{siegelpre3lemmaitem1} There exists an $X\in \Mat(2,F)$ of the form $X=\begin{bsmallmatrix} x& y \\ z & x\end{bsmallmatrix}$
such that
\begin{equation}
\begin{bsmallmatrix} A& \\ & uA'\end{bsmallmatrix} 
\begin{bsmallmatrix} \idm & X\\ & \vphantom{uA'}\idm\end{bsmallmatrix}
\in P\cap  t_n\Ks{n}t_n^{-1}.
\end{equation}
\item \label{siegelpre3lemmaitem2} The following conditions are satisfied:
\subitem $\bullet$ $u\in\OF^\times$,
\subitem $\bullet$ $\det(A)\in\OF^\times$,
\subitem $\bullet$ $A\in\begin{bsmallmatrix} \OF& \OF\\ \p^{n} & \OF\end{bsmallmatrix}$.
\end{enumerate}
\end{lemma}
\begin{proof}
The proof is similar to the proof of Lemma \ref{siegelpre2lemma}.
\end{proof}

\begin{lemma}\label{Siegelindlevellemma1}
Let $n$ be an integer such that $n\geq 1$. Let $u\in F^\times$ and $A=\begin{bsmallmatrix} a& b\\ c & d \end{bsmallmatrix}\in\GL(2,F)$ be given.
 For $1\leq i<n/2$ the following statements are equivalent:
 \begin{enumerate}
  \item \label{Siegelindlevellemma1item1} There exists an $X\in \Mat(2,F)$ of the form $X=\begin{bsmallmatrix} x& y \\ z & x\end{bsmallmatrix}$
   such that
   \begin{equation}\label{Siegelindlevellemmaeq0}
    \begin{bsmallmatrix} A& \\ & uA'\end{bsmallmatrix} \begin{bsmallmatrix} \idm & X\\ & \vphantom{uA'}\idm\end{bsmallmatrix}\in P\cap M_i \Ks{n} M_i^{-1}.
   \end{equation}
   Here $M_i$ is as in \eqref{Midefeq}.
  \item \label{Siegelindlevellemma1item2} The following conditions are satisfied:
   \subitem $\bullet$ $u\in\OF^\times$,
   \subitem $\bullet$ $\det(A)\in\OF^\times$,
   \subitem $\bullet$ $u\det(A)^{-1}\in 1+\p^i$,
   \subitem $\bullet$ $A\in\begin{bsmallmatrix} \OF& \p^{-i}\\ \p^{n-i} & \OF\end{bsmallmatrix}$.
 \end{enumerate}
\end{lemma}
\begin{proof} 
We begin with a preliminary observation. Let $X =\begin{bsmallmatrix} x&y\\z&x \end{bsmallmatrix}
\in \Mat(2,F)$. Then using \eqref{Ksshapeeq}, a calculation shows that 
$\begin{bsmallmatrix} A& \\ & uA'\end{bsmallmatrix} \begin{bsmallmatrix} \idm & X\\ & \vphantom{uA'}\idm\end{bsmallmatrix}\in P\cap M_i \Ks{n} M_i^{-1}$, or equivalently, 
$M_i^{-1}\begin{bsmallmatrix} A& \\ & uA'\end{bsmallmatrix}
\begin{bsmallmatrix} \idm& X\\ & \vphantom{uA'}\idm\end{bsmallmatrix} M_i\in \Ks{n}$, if and only if $u \in \OF^\times $ and the following
(\ref{Siegelindlevellemmaeq1}) to (\ref{Siegelindlevellemmaeq4})
are fulfilled:
\begin{align}\label{Siegelindlevellemmaeq1}
 &A(\idm+\varpi^iX)\in\begin{bsmallmatrix} \OF& \OF\\ \p^n & \OF\end{bsmallmatrix},\\
\label{Siegelindlevellemmaeq2}
 &AX\in\begin{bsmallmatrix} \OF & \p^{-n+1}\\ \OF & \OF \end{bsmallmatrix},\\
\label{Siegelindlevellemmaeq3}
 &uA'-A-\varpi^iAX\in\begin{bsmallmatrix} \p^{n-i}&\p^{-i}\\ \p^{n-i} & \p^{n-i}\end{bsmallmatrix},\\
\label{Siegelindlevellemmaeq4}
 &uA'-\varpi^iAX\in\begin{bsmallmatrix} \OF& \OF\\ \p^n & \OF\end{bsmallmatrix}.
\end{align}

Now assume that \ref{Siegelindlevellemma1item1} holds; we will prove \ref{Siegelindlevellemma1item2}. 
By the previous paragraph,  (\ref{Siegelindlevellemmaeq1}) to (\ref{Siegelindlevellemmaeq4}) hold.
By Lemma \ref{siegelprelemma} we have $u\in \OF^\times$ and
$\det(A)\in\OF^\times$. By (\ref{Siegelindlevellemmaeq3}) and
(\ref{Siegelindlevellemmaeq4}) we get
$$
A=\begin{bsmallmatrix} a & b \\ c & d\end{bsmallmatrix} \in\begin{bsmallmatrix} \p^{n-i} & \p^{-i}\\ \p^{n-i} & \p^{n-i}\end{bsmallmatrix}
+\begin{bsmallmatrix} \OF\vphantom{\p^{n-i}} & \vphantom{\p^{n-i}}\OF\vphantom{\p^{n-i}} \\ \p^n & \OF\vphantom{\p^{n-i}}\end{bsmallmatrix} 
=\begin{bsmallmatrix} \OF & \p^{-i} \\ \p^{n-i} & \OF\end{bsmallmatrix},
$$
hence $a,d\in\OF$, $b\in\p^{-i}$ and $c\in\p^{n-i}$.
Since $i<n/2$, we have $bc\in\p^{n-2i}\subset\p$, and it follows from $ad-bc\in\OF^\times$
that $a,d\in\OF^\times$. By (\ref{Siegelindlevellemmaeq2}) and
(\ref{Siegelindlevellemmaeq3}) we get
$$
 uA'-A\in\begin{bsmallmatrix} \p^{n-i}& \p^{-i} \\ \p^{n-i} & \p^{n-i}\end{bsmallmatrix}+
 \begin{bsmallmatrix} \p^i & \p^{i-n+1} \\ \p^i & \p^i\end{bsmallmatrix}=\begin{bsmallmatrix} \p^i & \p^{i-n+1} \\ \p^i & \p^i\end{bsmallmatrix}.
$$
Since $A'=\frac1{ad-bc}\begin{bsmallmatrix} a & -b \\ -c & d\end{bsmallmatrix}$, it follows that
$(\frac u{ad-bc}-1)a\in\p^i$, hence $\frac u{ad-bc}\in1+\p^i$.
Therefore all the conditions in \ref{Siegelindlevellemma1item2} are fulfilled.

Assume that \ref{Siegelindlevellemma1item2} holds.
We define
$$
 X:=\varpi^{-i}\Big(\frac u{ad-bc}A^{-1}\begin{bsmallmatrix} a & \varpi^ib \\
 \varpi^ic & d\end{bsmallmatrix}-\idm\Big).
$$
It is then easy to verify that $X$ has the required form and
that (\ref{Siegelindlevellemmaeq1}) to (\ref{Siegelindlevellemmaeq4})
are fulfilled. By the first paragraph, \ref{Siegelindlevellemma1item1} holds.
\end{proof}

\begin{lemma}\label{Siegelindlevellemma3}
Let $n$ be an integer such that $n\geq 1$.  
Let $u\in F^\times$ and $A=\begin{bsmallmatrix}a&b \\ c&d \end{bsmallmatrix} \in\GL(2,F)$ be given.
 For $1\leq i<n/2$ the following statements are equivalent:
 \begin{enumerate}
  \item \label{Siegelindlevellemma3item1} There exists an $X\in \Mat(2,F)$ of the form $X=\begin{bsmallmatrix} x& y \\ z & x\end{bsmallmatrix}$
   such that
   \begin{equation}\label{Siegelindlevellemma3eq0}
    \begin{bsmallmatrix} A& \\ & uA'\end{bsmallmatrix} \begin{bsmallmatrix} \idm & \vphantom{A'}X\\ & \idm\end{bsmallmatrix}
    \in P\cap M_it_n \Ks{n}t_n^{-1} M_i^{-1}.
   \end{equation}
   Here $M_i$ is as in \eqref{Midefeq}.
  \item \label{Siegelindlevellemma3item2} The following conditions are satisfied:
   \subitem $\bullet$ $u\in\OF^\times$,
   \subitem $\bullet$ $\det(A)\in\OF^\times$,
   \subitem $\bullet$ $u\det(A)^{-1}\in 1+\p^i$,
   \subitem $\bullet$ $A\in\begin{bsmallmatrix} \OF & \p^{-i}\\ \p^{n-i+1} & \OF\end{bsmallmatrix}$.
 \end{enumerate}
\end{lemma}
\begin{proof} 
We begin with a preliminary observation.
Let $X =\begin{bsmallmatrix} x&y\\z&x \end{bsmallmatrix}
\in \Mat(2,F)$.
Then using \eqref{Ksshapeeq}, a calculation shows that 
$\begin{bsmallmatrix} A& \\ & uA'\end{bsmallmatrix} \begin{bsmallmatrix} \idm & \vphantom{A'}X\\ & \idm\end{bsmallmatrix}
\in P\cap M_it_n \Ks{n}t_n^{-1} M_i^{-1}$, or equivalently,
$M_i^{-1}\begin{bsmallmatrix} A& \\ & uA'\end{bsmallmatrix} 
\begin{bsmallmatrix} \idm & \vphantom{A'}X\\ & \idm\end{bsmallmatrix} M_i\in t_n\Ks{n}t_n^{-1}$,
if and only if $u \in \OF^\times $ and the following
(\ref{Siegelindlevellemma3eq1}) to (\ref{Siegelindlevellemma3eq4})
are fulfilled:
\begin{align}\label{Siegelindlevellemma3eq1}
 &A(\idm+\varpi^iX)\in\begin{bsmallmatrix}\OF& \OF\\ \p^n & \OF\end{bsmallmatrix},\\
\label{Siegelindlevellemma3eq2}
 &AX\in\begin{bsmallmatrix} \OF & \p^{-n} \\ \OF & \OF\end{bsmallmatrix},\\
\label{Siegelindlevellemma3eq3}
 &uA'-A-\varpi^iAX\in\begin{bsmallmatrix} \p^{n-i} & \p^{-i} \\ \p^{n-i+1} & \p^{n-i}\end{bsmallmatrix},\\
\label{Siegelindlevellemma3eq4}
 &uA'-\varpi^iAX\in\begin{bsmallmatrix} \OF & \OF \\ \p^n & \OF\end{bsmallmatrix}.
\end{align}

Now assume that \ref{Siegelindlevellemma3item1} holds; we will prove \ref{Siegelindlevellemma3item2}. 
By the previous paragraph, (\ref{Siegelindlevellemma3eq1}) to (\ref{Siegelindlevellemma3eq4})
hold. 
By Lemma \ref{siegelprelemma} we have $u\in \OF^\times$ and
$\det(A)\in\OF^\times$. By (\ref{Siegelindlevellemma3eq3}) and
(\ref{Siegelindlevellemma3eq4}) we get
$$
A=\begin{bsmallmatrix} a & b \\ c & d\end{bsmallmatrix} 
\in\begin{bsmallmatrix} \p^{n-i} & \p^{-i}\\ \p^{n-i+1} & \p^{n-i} \end{bsmallmatrix}
+\begin{bsmallmatrix} \vphantom{\p^n} \OF &\vphantom{\p^{n-i}} \OF\\ \vphantom{\p^{n-i}} \p^n & \vphantom{\p^{n-i}} \OF\end{bsmallmatrix} 
=\begin{bsmallmatrix}\OF & \p^{-i} \\ \p^{n-i+1} & \OF\end{bsmallmatrix},
$$
hence $a,d\in\OF$, $b\in\p^{-i}$, and $c\in\p^{n-i+1}$.
Since $i<n/2$, we have $bc\in\p^{n-2i+1}\subset\p$, and it follows from $ad-bc\in\OF^\times$
that $a,d\in\OF^\times$. By (\ref{Siegelindlevellemma3eq2}) and
(\ref{Siegelindlevellemma3eq3}) we get
$$
 uA'-A\in\begin{bsmallmatrix} \p^{n-i} & \p^{-i} \\ \p^{n-i+1} & \p^{n-i}\end{bsmallmatrix} +
 \begin{bsmallmatrix} \p^i & \p^{i-n}\\ \p^i & \p^i \end{bsmallmatrix}=\begin{bsmallmatrix} \p^i & \p^{i-n} \\ \p^i & \p^i\end{bsmallmatrix}.
$$
Since $A'=\frac1{ad-bc}\begin{bsmallmatrix} a & -b\\ -c & d\end{bsmallmatrix}$, it follows that
$(\frac u{ad-bc}-1)a\in\p^i$, hence $\frac u{ad-bc}\in1+\p^i$.
Therefore all the conditions in \ref{Siegelindlevellemma3item2} are fulfilled.

Assume that  \ref{Siegelindlevellemma3item2} holds.
We define
$$
 X:=\varpi^{-i}\Big(\frac u{ad-bc}A^{-1}\begin{bsmallmatrix} a & \varpi^ib\\ 
 \varpi^ic & d\end{bsmallmatrix} -\idm\Big).
$$
It is then easy to verify that $X$ has the required form and
that (\ref{Siegelindlevellemma3eq1}) to (\ref{Siegelindlevellemma3eq4})
are fulfilled. By the first paragraph, \ref{Siegelindlevellemma3item1} holds.
\end{proof}

\begin{lemma}\label{Siegelindlevellemma2}
Let $n$ be an even integer such that $n\geq 2$. Let $u\in F^\times$ and $A=\begin{bsmallmatrix} a & b \\ c & d \end{bsmallmatrix} \in\GL(2,F)$ be given.
 The following statements are equivalent:
 \begin{enumerate}
  \item \label{Siegelindlevellemma2item1} There exists an $X\in \Mat(2,F)$ of the form $X=\begin{bsmallmatrix} x& y\\ z & x \end{bsmallmatrix}$
   such that
   \begin{equation}\label{Siegelindlevellemma2eq0}
    \begin{bsmallmatrix} A& \\ & uA'\end{bsmallmatrix} \begin{bsmallmatrix} \idm& X\\ &\vphantom{A'} \idm\end{bsmallmatrix} \in P\cap M_{n/2} \Ks{n} M_{n/2}^{-1}.
   \end{equation}
  \item \label{Siegelindlevellemma2item2} The following conditions are satisfied:
   \subitem $\bullet$ $u\in\OF^\times$,
   \subitem $\bullet$ $\det(A)\in\OF^\times$,
   \subitem $\bullet$ $u\det(A)^{-1}\in 1+\p^{n/2}$,
   \subitem $\bullet$ $A\in\begin{bsmallmatrix} \OF& \p^{-n/2+1}\\ \p^{n/2} & \OF\end{bsmallmatrix}$.
 \end{enumerate}
\end{lemma}
\begin{proof} We begin with a preliminary observation.
Let $X =\begin{bsmallmatrix} x&y\\z&x \end{bsmallmatrix}
\in \Mat(2,F)$.
Then using \eqref{Ksshapeeq},
a calculation shows that 
$$
\begin{bsmallmatrix} A& \\ & uA'\end{bsmallmatrix} 
\begin{bsmallmatrix} \idm& X\\ &\vphantom{A'} \idm\end{bsmallmatrix} 
\in P\cap M_{n/2} \Ks{n} M_{n/2}^{-1},
$$
or equivalently, 
$M_{n/2}^{-1}
\begin{bsmallmatrix} A & \\ & uA'\end{bsmallmatrix}
\begin{bsmallmatrix} \idm & X\\ & \vphantom{A'}\idm\end{bsmallmatrix} M_{n/2}\in \Ks{n}$,
if and only if $u \in \OF^\times $ and the following
(\ref{Siegelindlevellemma2eq1}) to (\ref{Siegelindlevellemma2eq4})
are fulfilled:
\begin{align}\label{Siegelindlevellemma2eq1}
 &A(\idm+\varpi^{n/2}X)\in\begin{bsmallmatrix} \OF & \OF \\ \p^n & \OF\end{bsmallmatrix},\\
\label{Siegelindlevellemma2eq2}
 &AX\in\begin{bsmallmatrix} \OF & \p^{-n+1} \\ \OF & \OF\end{bsmallmatrix},\\
\label{Siegelindlevellemma2eq3}
 &uA'-A-\varpi^{n/2}AX\in\begin{bsmallmatrix} \p^{n/2} & \p^{-n/2} \\ \p^{n/2} & \p^{n/2} \end{bsmallmatrix},\\
\label{Siegelindlevellemma2eq4}
 &uA'-\varpi^{n/2}AX\in\begin{bsmallmatrix} \OF & \OF \\ \p^n & \OF\end{bsmallmatrix}.
\end{align}

Assume that \ref{Siegelindlevellemma2item1} holds; we will prove \ref{Siegelindlevellemma2item2}.
By the previous paragraph, (\ref{Siegelindlevellemma2eq1}) to (\ref{Siegelindlevellemma2eq4})
hold.
By Lemma \ref{siegelprelemma} we have $u\in \OF^\times$ and
$\det(A)\in\OF^\times$. From (\ref{Siegelindlevellemma2eq2}) and
(\ref{Siegelindlevellemma2eq4}) we get
$$
 uA'=\frac u{ad-bc}\begin{bsmallmatrix}a&-b\\-c&d\end{bsmallmatrix}\in\begin{bsmallmatrix}\p^{n/2}& \p^{-n/2+1} \\ 
\p^{n/2}& \p^{n/2} \end{bsmallmatrix} 
 +\begin{bsmallmatrix} \vphantom{\p^{n-i}}\OF & \vphantom{\p^{n-i}}\OF \\ \vphantom{\p^{n-i}}\p^n & \vphantom{\p^{n-i}}\OF \end{bsmallmatrix} 
 =\begin{bsmallmatrix} \OF & \p^{-n/2+1} \\ \p^{n/2} & \OF\end{bsmallmatrix},
$$
hence $a,d\in\OF$, $b\in\p^{-n/2+1}$, and $c\in\p^{n/2}$.
Since $bc\in\p$ and $ad-bc\in\OF^\times$, we have $a,d\in\OF^\times$. By (\ref{Siegelindlevellemma2eq2}) and
(\ref{Siegelindlevellemma2eq3}) we get
$$
 uA'-A\in\begin{bsmallmatrix} \p^{n/2} & \p^{-n/2} \\ \p^{n/2} & \p^{n/2} \end{bsmallmatrix} +
 \begin{bsmallmatrix} \p^{n/2} & \p^{-n/2+1} \\ \p^{n/2} & \p^{n/2} \end{bsmallmatrix}
 =\begin{bsmallmatrix} \p^{n/2} & \p^{-n/2} \\ \p^{n/2} & \p^{n/2}\end{bsmallmatrix}.
$$
Since $A'=\frac1{ad-bc}\begin{bsmallmatrix} a & -b \\ -c & d \end{bsmallmatrix}$, it follows that
$(\frac u{ad-bc}-1)a\in\p^{n/2}$, and hence $\frac u{ad-bc}\in1+\p^{n/2}$. Therefore, all the conditions in \ref{Siegelindlevellemma2item2} are satisfied.

Assume that \ref{Siegelindlevellemma2item2} holds. We define
$$
 X:=\varpi^{-n/2}\Big(\frac u{ad-bc}A^{-1}\begin{bsmallmatrix} a & \varpi^{n/2}b \\ 
 \varpi^{n/2}c & d \end{bsmallmatrix}-\idm\Big).
$$
It is then easy to verify that $X$ has the required form and
that (\ref{Siegelindlevellemma2eq1}) to (\ref{Siegelindlevellemma2eq4})
are fulfilled. By the first paragraph, \ref{Siegelindlevellemma2item1} holds.
\end{proof}

Let $n$ be an integer such that $n\geq 0$.
If $n=0$, then we define $\Gamma_0(\p^n)=\GL(2,\OF)$.
If $n>0$, then we define 
\begin{equation}\label{K1defeq}
 \Gamma_0(\p^n):=\left\{\begin{bsmallmatrix} a & b \\ c & d \end{bsmallmatrix} \in\GL(2,\OF)\mid c\in\p^n
 \right\}.
\end{equation}
Evidently, $\Gamma_0(\p^n)$ is a subgroup of $\GL(2,\OF)$.
Let $(\tau,W)$ be an admissible representation of $\GL(2,F)$. Assume that $\tau$ admits
a central character, and that this character is trivial. We define
$\tau^{\Gamma_0(\p^n)} = \{w \in W\mid\tau(k)w=w\ \text{for $k \in \Gamma_0(\p^n)$}\}$.
If $\tau^{\Gamma_0(\p^m)} \neq 0$ for some $m$, then we let $N_{\tau}$
be the smallest such $m$; for convenience, if $\tau^{\Gamma_0(\p^m)} =0$
for all $m \geq 0$, then we also define $N_\tau=\infty$. We call $N_\tau$
the \emph{level} of $\tau$.
\label{gl2leveldef}

\begin{theorem}\label{Siegelindleveltheorem}
 Let $(\pi,W)$ be an admissible representation of $\GL(2,F)$ admitting a central character $\omega_\pi$, and let $\sigma$ be a character of $F^\times$. We assume that $\omega_\pi\sigma^2=1$,
 so that the induced representation $\pi\rtimes\sigma$ has trivial
 central character. 
Let $V$ be the space of $\pi \rtimes \sigma$. 
 \begin{enumerate}
  \item \label{Siegelindleveltheoremitem1} The minimal stable Klingen level of $\pi\rtimes\sigma$ is
   \begin{equation}\label{Siegelindleveltheoremeq0}
    N_{\pi\rtimes\sigma,s}=\begin{cases}
     N_{\sigma\pi}&\text{if }a(\sigma)=0,\\
     2a(\sigma)&\text{if }a(\sigma)>0\text{ and }N_{\sigma\pi}=0,\\
     N_{\sigma\pi}+2a(\sigma)-1&\text{if }a(\sigma)>0\text{ and }N_{\sigma\pi}>0.
    \end{cases}
   \end{equation}
   If this number is finite, then
   $$
    \dim V_s(N_{\pi\rtimes\sigma,s})=
    \begin{cases}
     \dim\big(V_{\sigma\pi}^{\GL(2,\OF)}\big)&\text{if $a(\sigma)=0$ and $N_{\sigma\pi}=0$},\\[2ex]
     2\dim\big(V_{\sigma\pi}^{\Gamma_0(\p^{N_{\sigma\pi}})}\big)&\text{if $a(\sigma)=0$ and $N_{\sigma\pi}>0$},\\[2ex]
     \dim\big(V_{\sigma\pi}^{\Gamma_0(\p)}\big)&\text{if $a(\sigma)>0$ and  $N_{\sigma\pi}=0$},\\[2ex]
     \dim\big(V_{\sigma\pi}^{\Gamma_0(\p^{N_{\sigma\pi}})}\big)&\text{if $a(\sigma)>0$ and  $N_{\sigma\pi}>0$}.
    \end{cases}
   $$
  \item \label{Siegelindleveltheoremitem2} Assume that $\pi$ is irreducible and infinite-dimensional.
   If $a(\sigma)=0$ and $N_{\sigma\pi}=0$, then
   $$\renewcommand{\arraystretch}{1.3}
    \dim V_s(n)=\frac{n^2+5n+2}2\qquad\text{for all integers }n\geq0.
   $$  
   If $a(\sigma)=0$ and $N_{\sigma\pi}>0$, then for all integers $n \geq 0$
   $$\renewcommand{\arraystretch}{1.3}
    \dim V_s(n)=\left\{\begin{array}{l@{\qquad\mbox{if }}l}
    \displaystyle\frac{(n-N_{\pi\rtimes\sigma,s}+1)(n-N_{\pi\rtimes\sigma,s}+4)}2&n\geq N_{\pi\rtimes\sigma,s},\\
    \;0&n<N_{\pi\rtimes\sigma,s}.\end{array}\right.
   $$
   If $a(\sigma)>0$ and $N_{\sigma\pi}=0$, then for all integers $n \geq 0$
   $$\renewcommand{\arraystretch}{1.3}
    \dim V_s(n)=\left\{\begin{array}{l@{\qquad\mbox{if }}l}
    \displaystyle\frac{(n-N_{\pi\rtimes\sigma,s}+1)(n-N_{\pi\rtimes\sigma,s}+4)}2&n\geq N_{\pi\rtimes\sigma,s},\\[2ex]
    \;0&n<N_{\pi\rtimes\sigma,s}.\end{array}\right.
   $$
   If $a(\sigma)>0$ and $N_{\sigma\pi}>0$, then for all integers $n \geq 0$
   $$\renewcommand{\arraystretch}{1.3}
    \dim V_s(n)=\left\{\begin{array}{l@{\qquad\mbox{if }}l}
    \displaystyle\frac{(n-N_{\pi\rtimes\sigma,s}+1)(n-N_{\pi\rtimes\sigma,s}+2)}2&n\geq N_{\pi\rtimes\sigma,s},\\
    \;0&n<N_{\pi\rtimes\sigma,s}.\end{array}\right.
   $$
  \item \label{Siegelindleveltheoremitem3} Assume that $\pi=\chi\triv_{\GL(2)}$ for a character $\chi$ of $F^\times$ and that $\sigma\chi$ is unramified. Then
   $$\renewcommand{\arraystretch}{1.3}
    \dim V_s(n)=n-N_{\pi\rtimes\sigma,s}+1\qquad\text{for all }n\geq N_{\pi\rtimes\sigma,s}.
   $$
  \item \label{Siegelindleveltheoremitem4} If $\pi=\chi\triv_{\GL(2)}$ for a character $\chi$ of $F^\times$ and $\sigma\chi$ is ramified, then
$V_s(n)=0$ for all integers $n \geq 0$. 
 \end{enumerate}
\end{theorem}
\begin{proof}
Let $W$ be the space of $\pi$. We
note that since $\omega_\pi \sigma^2=1$, the representation
$\sigma \pi$, which also has space $W$, has trivial central character.
Thus, the remarks preceding the statement of the theorem apply to $\sigma \pi$.  
Let $n$ be an integer such that $n\geq1$. 
Our first task is to find a formula for $\dim V_s(n)$ in terms of $\sigma$
and $\pi$.
By Lemma \ref{PKsdisjointlemma1} and 
\eqref{indreptrafoeq}, the elements of $V_s(n)$ are determined by their values on the set
$X$ from \eqref{indreptrafoeq021}. 
It follows that
$$
\dim V_s(n) = \sum_{g \in X} \dim V_s(n)_g,
$$
where $V_s(n)_g = \{f \in V_s(n)\mid\text{$f$ is supported on $Pg\Ks{n}$}\}$ for $g \in X$. Let $g \in X$. 
Define $C_g: V_s(n)_g \to W$ by $C_g(f)=f(g)$ for $f \in V_s(n)_g$. 
This linear map is injective; therefore, $\dim V_s(n)_g=\dim \mathrm{im}(C_g)$. It is straightforward
to verify that $\mathrm{im}(C_g) = W_g$, where
\begin{equation}
\label{indreptrafoeq21}
W_g =
\{w \in W\mid \sigma (u) \pi(A) w = w\ 
\text{for $\begin{bsmallmatrix} A&*\\ & uA'\end{bsmallmatrix} \in P \cap g\Ks{n}g^{-1}$}\}.
\end{equation}
For this, note that if $\begin{bsmallmatrix} A&*\\ & uA'\end{bsmallmatrix} \in P \cap g\Ks{n}g^{-1}$,
then $u,\det(A) \in \OF^\times$ by Lemma~\ref{siegelprelemma}. Hence,
\begin{equation}
\label{indreptrafoeq22}
\dim V_s(n) = \sum_{g \in X} \dim W_g.
\end{equation}
Next, we calculate the dimensions of the spaces $W_g$ for $g \in X$. 

Assume that $g =1$. By \eqref{indreptrafoeq21} and Lemma \ref{siegelpre2lemma}, we have
$$
W_1=\{w \in W\mid \sigma (u) \pi(A) w = w\ \mathrm{for}\ u \in \OF^\times, A \in \Gamma_0(\p^n)\}. 
$$
It follows that if $\sigma$ is ramified, then $W_1=0$, and if $\sigma$ is unramified,
then 
\begin{equation}
\label{indreptrafoeqf1}
W_1 = (\sigma \pi)^{\Gamma_0(\p^n)}.
\end{equation}
Thus,
\begin{equation}
\label{indreptrafoeq23}
\dim W_1 = \begin{cases}
0&\text{ if $\sigma$ is ramified},\\
\dim (\sigma \pi)^{\Gamma_0(\p^n)}&\text{ if $\sigma$ is unramified}.
\end{cases}
\end{equation}
A similar argument using \eqref{indreptrafoeq21} and Lemma \ref{siegelpre3lemma} proves that
\begin{equation}
\label{indreptrafoeq24}
\dim W_{t_n} = \begin{cases}
0&\text{ if $\sigma$ is ramified},\\
\dim (\sigma \pi)^{\Gamma_0(\p^n)}&\text{ if $\sigma$ is unramified}.
\end{cases}
\end{equation}
Assume that $g=M_i$ for some $0<i<n/2$. 
By \eqref{indreptrafoeq21} and Lemma \ref{Siegelindlevellemma1} we have
$$
W_{M_i}
=
\left\{
w \in W\mid
\sigma(u)\pi(A)w = w\quad\text{for}\quad
\begin{array}{l}
u \in \OF^\times,\\
A \in \begin{bsmallmatrix} \OF & \p^{-i} \\ \p^{n-i} & \OF \end{bsmallmatrix},\\
\det(A) \in \OF^\times,\\
u\det(A)^{-1} \in 1+\p^i
\end{array}
\right\}
$$
We see that if $W_{M_i}\neq 0$, then $\sigma(1+\p^i)=1$, i.e., $i \geq a(\sigma)$. 
Assume that $i \geq a(\sigma)$. Then
\begin{align}
W_{M_i}&=\{w \in W\mid (\sigma \pi) (A)w =w\ \mathrm{for}\ A \in \begin{bsmallmatrix} \OF & \p^{-i} \\ \p^{n-i} & \OF \end{bsmallmatrix}\ \mathrm{with}\ \det(A) \in \OF^\times \}\nonumber \\
&=(\sigma \pi)(\begin{bsmallmatrix}1&\\& \varpi^i\end{bsmallmatrix})\left((\sigma\pi)^{\Gamma_0(\p^{n-2i})}\right).
\label{indreptrafoeqf2}
\end{align}
Hence, for $0<i<n/2$, 
\begin{equation}
\label{indreptrafoeq25}
\dim W_{M_i}
=
\begin{cases}
0&\text{ if $i<a(\sigma)$},\\
\dim (\sigma\pi)^{\Gamma_0(\p^{n-2i})}&\text{ if $i \geq a(\sigma)$}.
\end{cases}
\end{equation}
A similar argument using  Lemma \ref{Siegelindlevellemma3}
proves that for $0<i<n/2$, 
\begin{equation}
\label{indreptrafoeq26}
\dim W_{M_it_n}
=
\begin{cases}
0&\text{ if $i<a(\sigma)$},\\
\dim (\sigma\pi)^{\Gamma_0(\p^{n-2i+1})}&\text{ if $i \geq a(\sigma)$}.
\end{cases}
\end{equation}
Finally, assume that $n$ is even. By 
\eqref{indreptrafoeq21} and Lemma \ref{Siegelindlevellemma2}
we have
$$
W_{M_{n/2}}
=
\left\{
w \in W\mid
\sigma(u)\pi(A)w = w\quad\text{for}\quad
\begin{array}{l}
u \in \OF^\times,\\
A \in \begin{bsmallmatrix} \OF & \p^{-n/2+1} \\ \p^{n/2} & \OF \end{bsmallmatrix},\\
\det(A) \in \OF^\times,\\
u\det(A)^{-1} \in 1+\p^{n/2}
\end{array}
\right\}.
$$
Evidently, if $W_{M_{n/2}}\neq 0$, then $\sigma(1+\p^{n/2})=1$,
i.e., $n/2\geq a(\sigma)$. Assume that $n/2\geq a(\sigma)$.
Then
\begin{align}
W_{M_{n/2}}&=\{w \in W\mid (\sigma \pi) (A)w =w\ \mathrm{for}\ A \in \begin{bsmallmatrix} \OF & \p^{-n/2+1} \\ \p^{n/2} & \OF \end{bsmallmatrix}\ \mathrm{with}\ \det(A) \in \OF^\times \}\nonumber\\
&=(\sigma \pi)(\begin{bsmallmatrix}1&\\& \varpi^{n/2-1}\end{bsmallmatrix})\left((\sigma\pi)^{\Gamma_0(\p)}\right).
\label{indreptrafoeqf3}
\end{align}
It follows that if $n$ is even, then 
\begin{equation}
\label{indreptrafoeq27}
\dim W_{M_{n/2}}
=
\begin{cases}
0&\text{ if $n/2<a(\sigma)$},\\
\dim (\sigma\pi)^{\Gamma_0(\p)}&\text{ if $n/2 \geq a(\sigma)$}.
\end{cases}
\end{equation}

By \eqref{indreptrafoeq22}, \eqref{indreptrafoeq23}, \eqref{indreptrafoeq24}, \eqref{indreptrafoeq25}, \eqref{indreptrafoeq26}, and \eqref{indreptrafoeq27} we now see that
 if $a(\sigma)=0$, then
\begin{align}
\label{Siegelindleveltheoremeq10}
\dim V_s(n)&=2\dim (\sigma\pi)^{\Gamma_0(\p^n)}+\sum_{0<i<\frac n2}\Big((\sigma\pi)^{\Gamma_0(\p^{n-2i})}+\dim (\sigma\pi)^{\Gamma_0(\p^{n-2i+1})}\Big)\nonumber\\
 &\quad+\left\{
  \begin{array}{cl}
   \!\dim (\sigma\pi)^{\Gamma_0(\p)}&\!\text{if $n$ even}\\[1ex]0&\!\text{if $n$ odd}
  \end{array}
 \right\},
\end{align}
and if $a(\sigma)>0$, then
\begin{align}\label{Siegelindleveltheoremeq11}
 \dim V_s(n)&=\sum_{a(\sigma)\leq i<\frac n2}\Big(\dim (\sigma\pi)^{\Gamma_0(\p^{n-2i})}+\dim (\sigma\pi)^{\Gamma_0(\p^{n-2i+1})}\Big)\nonumber\\
 &\quad+\left\{
  \begin{array}{cl}
   \!\dim (\sigma\pi)^{\Gamma_0(\p)}&\!\text{if $n$ even and }a(\sigma)\leq\frac n2\!\\[1ex]0&\!\text{if $n$ odd or }a(\sigma)>\frac n2\!
  \end{array}
 \right\}.
\end{align}
The formulas \eqref{Siegelindleveltheoremeq10} and \eqref{Siegelindleveltheoremeq11} hold for all $n\geq1$.

\ref{Siegelindleveltheoremitem1}. The assertion is easy to prove if $a(\sigma)=0$ and $N_{\sigma\pi}=0$. Assume therefore that either $a(\sigma)>0$ or $N_{\sigma\pi}>0$.
Then $V_s(0)=0$. It follows from \eqref{Siegelindleveltheoremeq10} and \eqref{Siegelindleveltheoremeq11} that the minimal stable Klingen level is
$$
N_{\pi\rtimes\sigma,s}=\begin{cases}
N_{\sigma\pi}&\text{if }a(\sigma)=0,\\
N_{\sigma\pi}+2a(\sigma)-1&\text{if }a(\sigma)>0\text{ and }N_{\sigma\pi}\geq 1,\\
2a(\sigma)&\text{if }a(\sigma)>0\text{ and }N_{\sigma\pi}=0.
\end{cases}
$$
The dimension of $V_s(N_{\pi\rtimes\sigma,s})$ follows by substituting $N_{\pi\rtimes\sigma,s}$ into \eqref{Siegelindleveltheoremeq10} resp.\ \eqref{Siegelindleveltheoremeq11}.

\ref{Siegelindleveltheoremitem2} follows from evaluating the formulas \eqref{Siegelindleveltheoremeq10} and \eqref{Siegelindleveltheoremeq11}, keeping in mind that
\begin{equation}\label{GL2dimformulaeq}
 \dim V_{\sigma\pi}^{\Gamma_0(\p^m)}=\begin{cases}
                              m-N_{\sigma\pi}+1&\text{if }m\geq N_{\sigma\pi},\\
                              0&\text{if }m<N_{\sigma\pi}.
                             \end{cases}
\end{equation}

\ref{Siegelindleveltheoremitem3} and \ref{Siegelindleveltheoremitem4}  follow similarly from \eqref{Siegelindleveltheoremeq10} and \eqref{Siegelindleveltheoremeq11}.
\end{proof}

For future use, we point out some facts that are 
evident from the proof of Theorem~\ref{Siegelindleveltheorem}.
If the notation is as in Theorem~\ref{Siegelindleveltheorem},  $n$ is an
integer such that $n \geq 1$, and $f \in V_s(n)$, then
\begin{align}
f(1), f(t_n) & \in 
\begin{cases}
\{0\}&\text{if $\sigma$ is ramified},\\
(\sigma \pi)^{\Gamma_0(\p^n)}&\text{if $\sigma$ is unramified},
\end{cases}
\label{Xextraeq1} \\
(\sigma \pi)(\begin{bsmallmatrix}1&\\&\varpi^{-i} \end{bsmallmatrix})f(M_i)
& \in 
\begin{cases}
\{0\}&\text{if $0<i<a(\sigma)$},\\
(\sigma \pi)^{\Gamma_0(\p^{n-2i})}&\text{if $a(\sigma)\leq  i < n/2$} 
\end{cases}
\label{Xextraeq2} \\
(\sigma \pi)(\begin{bsmallmatrix}1&\\&\varpi^{-i} \end{bsmallmatrix})f(M_it_n)
& \in 
\begin{cases}
\{0\}&\text{if $0<i<a(\sigma)$},\\
(\sigma \pi)^{\Gamma_0(\p^{n-2i+1})}&\text{if $a(\sigma)\leq  i < n/2$},
\end{cases}
\label{Xextraeq3} \\
(\sigma \pi)(\begin{bsmallmatrix}1&\\&\varpi^{-(n/2-1)} \end{bsmallmatrix})f(M_{n/2})
& \in (\sigma \pi)^{\Gamma_0(\p)}\quad \text{if $n$ is even}. \label{Xextraeq4}
\end{align}
\section{Non-existence of certain vectors}\label{nonexistencesec}
The following result will be used in the proof of Proposition \ref{funnyidentityprop}.
\begin{lemma}\label{rhopantilemma}
 Let $\pi$ be an admissible representation of $\GL(2,F)$ admitting a central character $\omega_\pi$. Assume that either
 \begin{equation}\label{rhopantilemmaeq10}
  \text{$\pi$ does not contain any non-zero vectors fixed by $\SL(2,F)$},
 \end{equation}
 or
 \begin{equation}\label{rhopantilemmaeq11}
  \text{$\pi=\chi1_{\GL(2)}$ for $\chi$ a character of $F^\times$ with $\chi\neq\nu^{1/2}$}.
 \end{equation}
 Let $\sigma$ be a character of $F^\times$. We assume that $\omega_\pi\sigma^2=1$, so that the induced representation $\pi\rtimes\sigma$ of $\GSp(4,F)$ has trivial
 central character. 
Let $V$ be the space of $\pi \rtimes \sigma$, and
let $n$ be an integer such that $n\geq 0$.  Suppose that $f\in V_s(n)$ satisfies
 \begin{equation}\label{rhopantilemmaeq2}
  \tau_{n+1}f=q^{-2}\eta f.
 \end{equation}
 Then $f=0$.
\end{lemma}
\begin{proof}
We assume that $V$ is the standard model of $\pi\rtimes \sigma$ (see Sect.~\ref{repsec}).
We note that Lemma \ref{funnyidentitylemma5} implies that
 \begin{equation}\label{rhopantilemmaeq1}
  -q\tau_nf=t_nf
 \end{equation}
since \eqref{rhopantilemmaeq2} holds by assumption.
If $n=0$, then \eqref{rhopantilemmaeq1} is $-f=f$, so that $f=0$.
Assume that $n\geq1$.  By Lemma \ref{PKsdisjointlemma1},
\begin{align}\label{PKsdisjointeq30}
\GSp(4,F)&=P\Ks{n}\;\sqcup\bigsqcup_{0<i\leq n/2}P M_i\Ks{n}\nonumber\\
&\quad\sqcup\;Pt_n\Ks{n}\;\sqcup\bigsqcup_{0<i<n/2}P M_it_n\Ks{n}.
\end{align}
To prove that $f=0$ it will suffice to 
prove that $f$ vanishes on $1$, $t_n$, $M_i$ and $M_it_n$ for $0<i\leq n/2$. 
Evaluating \eqref{rhopantilemmaeq2} at $1$, and using \eqref{indreptrafoeq}, gives
\begin{equation}\label{rhopantilemmaeq3}
 f(1)=q^{-1/2}\pi(\begin{bsmallmatrix} \varpi^{-1}& \\ & 1 \end{bsmallmatrix} )f(1).
\end{equation}
Since $f \in V_s(n)$,  we also deduce from \eqref{indreptrafoeq} that 
\begin{align}
\pi(\begin{bsmallmatrix}1& x \\ &1 \end{bsmallmatrix})f(1)&=f(1)\quad\text{for $x\in \OF$},\label{rhopantilemmaeq31}\\
\pi(\begin{bsmallmatrix}1&\\ y&1 \end{bsmallmatrix})f(1)&=f(1)\quad\text{for $y \in \p^n$},\label{rhopantilemmaeq32}\\
\pi(\begin{bsmallmatrix}v&\\ &v^{-1} \end{bsmallmatrix})f(1)&=f(1)\quad\text{for $v \in \OF^\times$}.\label{rhopantilemmaeq34}
\end{align}
By \eqref{rhopantilemmaeq3} and \eqref{rhopantilemmaeq32} we have 
\begin{equation}
\label{rhopantilemmaeq33}
\pi(\begin{bsmallmatrix}1&\\ y&1 \end{bsmallmatrix})f(1)=f(1)\quad\text{for $y \in F$}.
\end{equation}
The identity
$$
\begin{bsmallmatrix}
&1\\
-1&
\end{bsmallmatrix}
=
\begin{bsmallmatrix}
-1&\\
&-1
\end{bsmallmatrix}
\begin{bsmallmatrix}
1&-1\\
&1
\end{bsmallmatrix}
\begin{bsmallmatrix}
1&\vphantom{-1}\\1&1
\end{bsmallmatrix}
\begin{bsmallmatrix}
1&-1\\
&1
\end{bsmallmatrix}
$$
now implies that 
$\pi(\begin{bsmallmatrix}
&1\\
-1&
\end{bsmallmatrix})f(1)=f(1)$. This, along with \eqref{rhopantilemmaeq33}, proves that 
$f(1)$ is fixed by $\SL(2,F)$ (see Satz A 5.4 of \cite{Fr}). If \eqref{rhopantilemmaeq10} is satisfied, then it follows that $f(1)=0$. Assume that \eqref{rhopantilemmaeq11} is satisfied. Then \eqref{rhopantilemmaeq3} becomes
\begin{equation}\label{rhopantilemmaeq3b}
 f(1)=q^{-1/2}\chi(\varpi)^{-1}f(1).
\end{equation}
By \eqref{indreptrafoeq} and the assumption $f \in V_s(n)$, 
\begin{equation}\label{rhopantilemmaeq3c}
 f(1)=f(\begin{bsmallmatrix}u\\&1\\&&1\\&&&u^{-1}\end{bsmallmatrix})=\chi(u)f(1)\quad\text{for $u\in\OF^\times$}.
\end{equation}
If $f(1)\neq0$, then \eqref{rhopantilemmaeq3b} and \eqref{rhopantilemmaeq3c} imply that $\chi=\nu^{1/2}$, contradiction. Hence $f(1)=0$. Evaluating \eqref{rhopantilemmaeq1} at $1$ yields $f(t_n)=-q f(1)=0$.

To complete the proof we will prove that 
$f(M_i)=f(M_it_n)=0$ for $0\leq i\leq n/2$ by induction on $i$. There are matrix identities
\begin{align}
M_0
&=
\begin{bsmallmatrix}
&1&&\varpi^{-n}\\
\varpi^n&&1&\\
&&&\varpi^{-n}\\
&&1&
\end{bsmallmatrix}
t_n 
\begin{bsmallmatrix}
1&&1&\vphantom{\varpi^{-n}}\\
&&-1&\vphantom{\varpi^{-n}}\\
&1&&1\vphantom{\varpi^{-n}}\\
&&&1\vphantom{\varpi^{-n}}
\end{bsmallmatrix},\label{rhopantilemmaeq01} \\
M_0t_n
&=
\begin{bsmallmatrix}
&1&&\varpi^{-n}\\
\varpi^n&&1&\\
&&&\varpi^{-n}\\
&&1&
\end{bsmallmatrix}
\begin{bsmallmatrix}
-1&&&\vphantom{\varpi^{-n}}\\
&&-1&\vphantom{\varpi^{-n}}\\
\varpi^n&1&&\vphantom{\varpi^{-n}}\\
&&\varpi^n&-1\vphantom{\varpi^{-n}}
\end{bsmallmatrix}.\label{rhopantilemmaeq02}
\end{align}
These identities imply that 
$$
PM_0\Ks{n}=Pt_n\Ks{n}\quad\text{and}\quad PM_0t_n\Ks{n}=P\Ks{n}.
$$
Since $f(t_n)=0$ and $f(1)=0$, we conclude that $f(M_0)=0$
and $f(M_0t_n)=0$. This proves our assertion for $i=0$. Assume 
that $i>0$, and that the assertion has been proven for $i-1$. 
We have the following identity for $z \in \OF^\times$:
\begin{align}\label{rhopantilemmaeq4}
 &\begin{bsmallmatrix}1&z\varpi^{-i}\\z^{-1}\varpi^{n-i}&1+\varpi^{n-2i}\\&&1&-z\varpi^{-i}\\&&-z^{-1}\varpi^{n-i}&1+\varpi^{n-2i}\end{bsmallmatrix}\begin{bsmallmatrix}1&&\varpi^{n-3i}&-z\varpi^{n-4i}-z\varpi^{-2i}\\&1&-z^{-1}\varpi^{n-2i}&\varpi^{n-3i}\\&&1\\&&&1\end{bsmallmatrix}M_it_n\nonumber\\
 &\qquad=M_i\begin{bsmallmatrix}1&&&-z\varpi^{-n}\\\vphantom{\varpi^n}&1\\&&1\\&\vphantom{z^{-1}}&&1\end{bsmallmatrix}\begin{bsmallmatrix}z\\&1+\varpi^{n-2i}&-z^{-1}\varpi^{n-2i}\\&-z&1\\\varpi^n+\varpi^{2n-2i}&&&z^{-1}\end{bsmallmatrix}.
\end{align}
Now
\begin{align}\label{rhopantilemmaeq5}
(\tau_nf)(M_i)&=\int\limits_\OF f(M_i\begin{bsmallmatrix}1&&&z\varpi^{-n}\\&1\\&&1\\&&&1\end{bsmallmatrix})\,dz\nonumber\\
&=\int\limits_\p f(M_i\begin{bsmallmatrix}1&&&z\varpi^{-n}\\&1\\&&1\\&&&1\end{bsmallmatrix})\,dz+\int\limits_{\OF^\times} f(M_i\begin{bsmallmatrix}1&&&z\varpi^{-n}\\&1\\&&1\\&&&1\end{bsmallmatrix})\,dz\nonumber\\
&=q^{-1}f(M_i)+\int\limits_{\OF^\times} f(M_i\begin{bsmallmatrix}1&&&-z\varpi^{-n}\\&1\\&&1\\&&&1\end{bsmallmatrix})\,dz\nonumber\\
&=q^{-1}f(M_i)+\int\limits_{\OF^\times}\pi(
 \begin{bsmallmatrix} 1& z\varpi^{-i}\\ z^{-1}\varpi^{n-i} & 1+\varpi^{n-2i}\end{bsmallmatrix}) f(M_it_n)\,dz.
\end{align}
For the last equality we used \eqref{rhopantilemmaeq4}.
Also,
\begin{align}
(\tau_{n+1}f)(M_i)&=\int\limits_\OF f(M_i\begin{bsmallmatrix}1&&&z\varpi^{-n-1}\\&1\\&&1\\&&&1\end{bsmallmatrix})\,dz\nonumber\\
&=\int\limits_\p f(M_i
\begin{bsmallmatrix}1&&&z\varpi^{-n-1}\\&1\\&&1\\&&&1\end{bsmallmatrix})\,dz+\int\limits_{\OF^\times} f(M_i
\begin{bsmallmatrix}1&&&z\varpi^{-n-1}\\&1\\&&1\\&&&1\end{bsmallmatrix})\,dz\nonumber\\
&=q^{-1}(\tau_nf)(M_i)+\int\limits_{\OF^\times} f(M_i
\begin{bsmallmatrix}1&&&-z\varpi^{-n-1}\\&1\\&&1\\&&&1\end{bsmallmatrix})\,dz\nonumber\\
&=q^{-1}(\tau_nf)(M_i)+\int\limits_{\OF^\times}\pi(
 \begin{bsmallmatrix} 1 & z\varpi^{-i} \\ z^{-1}\varpi^{n+1-i} & 1+\varpi^{n+1-2i} \end{bsmallmatrix}) f(M_it_{n+1})\,dz\nonumber\\
&=q^{-1}(\tau_nf)(M_i)+\int\limits_{\OF^\times}\pi(
 \begin{bsmallmatrix} 1 & z\varpi^{-i} \\ z^{-1}\varpi^{n+1-i} & 1+\varpi^{n+1-2i} \end{bsmallmatrix}) f(M_i\eta t_n)\,dz\nonumber\\
&=q^{-1}(\tau_nf)(M_i)+\int\limits_{\OF^\times}\pi(
 \begin{bsmallmatrix} 1 & z\varpi^{-i} \\ z^{-1}\varpi^{n+1-i} & 1+\varpi^{n+1-2i} \end{bsmallmatrix}) f(\eta M_{i-1} t_n)\,dz\nonumber\\
\label{rhopantilemmaeq6}&=q^{-1}(\tau_nf)(M_i)\\
&\quad+q^{3/2}\int\limits_{\OF^\times}\pi(
\begin{bsmallmatrix} 1 & z\varpi^{-i} \\ z^{-1}\varpi^{n+1-i} & 1+\varpi^{n+1-2i}\end{bsmallmatrix}
\begin{bsmallmatrix} \varpi^{-1}& \\ \vphantom{z^{-1}}& 1 \end{bsmallmatrix}) f(M_{i-1} t_n)\,dz.\nonumber
\end{align}
Here, we again used \eqref{rhopantilemmaeq4} in the fourth equality.
Evaluating \eqref{rhopantilemmaeq1} and \eqref{rhopantilemmaeq2} at $M_i$ gives
 \begin{equation}\label{rhopantilemmaeq7}
  -q(\tau_nf)(M_i)=f(M_it_n)
 \end{equation}
 and
 \begin{equation}\label{rhopantilemmaeq8}
  (\tau_{n+1}f)(M_i)=q^{-2}f(M_i\eta)=q^{-1/2}\pi(
  \begin{bsmallmatrix} \varpi^{-1}& \\ & 1\end{bsmallmatrix})f(M_{i-1})=0;
 \end{equation}
note that $f(M_{i-1})=0$ by the induction hypothesis.
It follows from \eqref{rhopantilemmaeq6}, \eqref{rhopantilemmaeq8}, and the induction hypothesis that $(\tau_nf)(M_i)=0$. From \eqref{rhopantilemmaeq7} we have $f(M_it_n)=0$. 
Since $f(M_it_n)=0$ and $(\tau_nf)(M_i)=0$, 
we also deduce from \eqref{rhopantilemmaeq5} that $f(M_i)=0$. This concludes the proof.
\end{proof}
\section{Characterization of paramodular vectors}\label{paracharsec}
To end this chapter, we prove a result which characterizes the paramodular vectors among the stable Klingen vectors in a Siegel induced representation.

\begin{lemma}\label{paracharlemma}
Let $(\pi,W)$ be an infinite-dimensional, irreducible, admissible representation of $\GL(2,F)$ with central character $\omega_\pi$, and let $\sigma$ be a character of $F^\times$. We assume that $\omega_\pi\sigma^2=1$,
so that the induced representation $\pi\rtimes\sigma$ has trivial
central character. 
Let $V$ be the standard model of $\pi \rtimes \sigma$ (see Sect.~\ref{repsec}).  
 \begin{enumerate}
\item \label{paracharlemmaitem1} Let $n$ be an integer such that $n\geq N_{\pi\rtimes \sigma,s}$ and let $f \in V_s(n)$.
Then $f$ lies in $V(n)$ if and only if the following conditions are satisfied.
   \begin{align}
    \label{paracharlemma1}f(t_n)&=f(1),\\
    \label{paracharlemma2}f(M_it_n)&=f(M_i)\quad\text{for }0<i<n/2,\\
    \label{paracharlemma3}f(M_{n/2})&=\pi(
    \begin{bsmallmatrix} &-\varpi^{-n/2}\\ \varpi^{n/2}& \end{bsmallmatrix})f(M_{n/2})\quad\text{if $n$ is even}.
   \end{align}
  \item \label{paracharlemmaitem2} Let $n$ be an integer such that $n\geq N_{\pi\rtimes \sigma,s}$ and $n \geq 1$,  and let $f \in V_s(n)$.
Then $f$ lies in $V(n-1)$ if and only if the following conditions are satisfied.
   \begin{align}
    \label{paracharlemma4}f(t_n)&=q^{3/2}\pi(
    \begin{bsmallmatrix} \varpi^{-1}& \\ & 1 \end{bsmallmatrix})f(1),\\
    \label{paracharlemma5}f(M_it_n)&=q^{3/2}\pi(
    \begin{bsmallmatrix} \varpi^{-1}& \\ & 1\end{bsmallmatrix})f(M_{i-1})\quad\text{for }1<i<n/2,\\
    \label{paracharlemma6}f(M_1t_n)&=q^{3-(3/2)n}\pi(
    \begin{bsmallmatrix} & \varpi^{-1}\\ \varpi^{n-1} & \end{bsmallmatrix})f(1),\\     \label{paracharlemma7}f(M_{n/2})&=q^{3/2}\pi(
    \begin{bsmallmatrix} & -\varpi^{-n/2} \\ \varpi^{n/2-1} & \end{bsmallmatrix} )f(M_{n/2-1})\quad\text{if $n$ is even},\\
    \label{paracharlemma8}f(M_{(n-1)/2})&=\pi(
    \begin{bsmallmatrix} & -\varpi^{-(n-1)/2} \\ \varpi^{(n-1)/2} & \end{bsmallmatrix} )f(M_{(n-1)/2})\quad\text{if $n$ is odd}.
   \end{align}
 \end{enumerate}
\end{lemma}
\begin{proof}
We begin by noting  that by \ref{Siegelindleveltheoremitem1} of Theorem \ref{Siegelindleveltheorem},
\begin{equation}
\label{parachareq0010}
N_{\pi \rtimes \sigma,s}
=
\begin{cases}
N_{\sigma\pi}+2a(\sigma)&\text{if $a(\sigma)=0$ or $N_{\sigma\pi}=0$},\\
N_{\sigma\pi}+2a(\sigma)-1&\text{if $a(\sigma)>0$ and $N_{\sigma\pi}>0$}.
\end{cases}
\end{equation}
We see that if $n \geq N_{\pi\rtimes \sigma,s}$, then $a(\sigma) \leq n/2$. We also note that the elements of $V$ satisfy \eqref{indreptrafoeq};
this will be used repeatedly in the following proof.

\ref{paracharlemmaitem1}. Let $n$ be an integer such that $n\geq N_{\pi\rtimes \sigma,s}$.
Let $A(n)$ be the subspace of all $f\in V_s(n)$ satisfying \eqref{paracharlemma1}--\eqref{paracharlemma3}. 
We need to prove that $V(n)=A(n)$. 
Let $f\in V(n)$. Then $f(gt_n)=f(g)$ for all $g\in\GSp(4,F)$. Substituting $g=1$ and $g=M_i$, we obtain \eqref{paracharlemma1}--\eqref{paracharlemma3}; for \eqref{paracharlemma3}, observe the matrix identity \eqref{PKsdisjointeq12}. Thus, $V(n) \subset A(n)$. 
To prove $V(n)=A(n)$, we will show that $\dim V(n)=\dim A(n)$. 
The dimension of $V(n)$ is known: by Theorem 5.2.2 of \cite{NF}, the minimal paramodular level of $\pi\rtimes\sigma$ is $N_{\pi\rtimes\sigma}=N_{\sigma\pi}+2a(\sigma)$, and
\begin{equation}
\label{parachareq002}
\renewcommand{\arraystretch}{1.3}
\dim V(n)
=\left\{\begin{array}{l@{\qquad\mbox{if }}l}
\displaystyle\Big\lfloor\frac{(n-N_{\pi\rtimes\sigma}+2)^2}4\Big\rfloor&n\geq N_{\pi\rtimes\sigma},\\
\;0&n<N_{\pi\rtimes\sigma}.\end{array}\right.
\end{equation}
To estimate the dimension of $A(n)$ we need a preliminary result. 
For each integer $i$ such that $a(\sigma)\leq i \leq n/2$ define
 $g_i \in \GSp(4,F)$ by 
$$
g_i = 
\begin{cases}
1&\text{if $i=0$},\\
M_i&\text{if $1<i\leq n/2$}.
\end{cases}
$$
Let $f \in A(n)$. We claim that 
\begin{equation}
\label{parachareq001}
(\sigma\pi)(
\begin{bsmallmatrix} 
1& \\ &\varpi^{-i}
\end{bsmallmatrix})
f(g_i)
\in (\sigma\pi)^{\Gamma_0(\p^{n-2i})}
\quad 
\text{for $a(\sigma)\leq i \leq n/2$}.
\end{equation}
Let $i$ be an integer such that $a(\sigma) \leq i \leq n/2$. If $i=0$, so that
necessarily $a(\sigma)=0$, then \eqref{parachareq001} follows from
\eqref{Xextraeq1}. 
If  $0<i<n/2$,
then \eqref{parachareq001} follows from \eqref{Xextraeq2}.
Assume that $i=n/2$, so that $n$ is necessarily even.
By 
\eqref{Xextraeq4} we have 
\begin{align}
(\sigma \pi)(\begin{bsmallmatrix} 1&\\ &\varpi^{-(n/2-1)}\end{bsmallmatrix})f(M_{n/2})
&\in (\sigma \pi)^{\Gamma_0(\p)}.\label{parachareq03}
\end{align}
Define $w=(\sigma \pi)(\begin{bsmallmatrix} 1&\\ &\varpi^{-n/2}\end{bsmallmatrix})f(M_{n/2})$.
By \eqref{parachareq03}, we have $(\sigma \pi) (\begin{bsmallmatrix}1&\\&\varpi\end{bsmallmatrix})w
\in (\sigma \pi)^{\Gamma_0(\p)}$. This implies that
\begin{align}
\label{parachareq04}
(\sigma\pi)(\begin{bsmallmatrix}1&\\y&1\end{bsmallmatrix})w &=w\quad\text{for}\quad y \in \OF,\\
\label{parachareq06}
(\sigma\pi)(\begin{bsmallmatrix}1&\\&u\end{bsmallmatrix})w &=w\quad\text{for}\quad u \in \OF^\times.
\end{align}
Since $f \in A(n)$, \eqref{paracharlemma3} holds; therefore,
\begin{equation}
\label{parachareq05}
(\sigma\pi)(\begin{bsmallmatrix}&-1\\1&\end{bsmallmatrix})w =w.
\end{equation}
By \eqref{parachareq04} and \eqref{parachareq05} we have $(\sigma\pi)(k)w=w$ for $k \in \SL(2,\OF)$
(see Satz A 5.4 of \cite{Fr}); by \eqref{parachareq06}, recalling the definition of $w$, we conclude that 
\begin{equation}
\label{parachareq07}
(\sigma \pi)(\begin{bsmallmatrix} 1&\\ &\varpi^{-n/2}\end{bsmallmatrix})f(M_{n/2})
\in (\sigma\pi)^{\GL(2,\OF)}.
\end{equation}
This completes the argument for \eqref{parachareq001}.
Next, let 
$$
B(n)=\bigoplus_{a(\sigma)\leq i \leq n/2} (\sigma \pi)^{\Gamma_0(\p^{n-2i})}.
$$
Define $\varphi: A(n)\to B(n)$
 by 
\begin{equation*}
\varphi(f)= \bigoplus_{a(\sigma)\leq i \leq n/2}
(\sigma\pi)(
\begin{bsmallmatrix} 
1& \\ &\varpi^{-i}
\end{bsmallmatrix})
f(g_i).
\end{equation*}
Then $\varphi$ is a well-defined linear map by \eqref{parachareq001}.
We claim that $\varphi$ is injective. Suppose that $f \in A(n)$
is such that $\varphi(f)=0$. By the definition of $\varphi$,
we have $f(g_i)=0$ for  all integers $i$
such that $a(\sigma) \leq i \leq n/2$. Since 
\eqref{Xextraeq1} and \eqref{Xextraeq2}
hold, and since
$f$ satisfies
\eqref{paracharlemma1} and \eqref{paracharlemma2}, it follows
$f$ vanishes on the set $X$ from \eqref{indreptrafoeq021}, and
hence on $\GSp(4,F)$ by Lemma \ref{PKsdisjointlemma1} and 
\eqref{indreptrafoeq}. 
Thus, $\varphi$ is injective.
It follows
that
$
\dim A(n) \leq \dim B(n).
$
We will now prove that  $\dim B(n)=\dim V(n)$. Evidently,
$$
\dim B(n) = 
\sum_{\substack{a(\sigma)\leq i \leq n/2}} \dim (\sigma \pi)^{\Gamma_0(\p^{n-2i})}.
$$
Hence, using \eqref{GL2dimformulaeq},
\begin{align*}
\dim B(n) 
&= 
\sum_{\substack{a(\sigma)\leq i \leq  (n-N_{\sigma \pi})/2}} n-N_{\sigma \pi}-2i+1.
\end{align*}
Assume first that this sum is empty. Then  $\dim B(n)=0$ and $n<N_{\sigma \pi}+2a(\sigma)$.
Since 
 $\dim V(n)=0$ by \eqref{parachareq002}, it follows that $\dim B(n)=\dim V(n)=0$ in this
case. Now assume that the sum is not empty. 
A calculation shows that the sum is 
$$
\displaystyle\Big\lfloor\frac{(n-N_{\sigma\pi}-2a(\sigma)+2)^2}4\Big\rfloor.
$$
By \eqref{parachareq002} this is $\dim V(n)$, so that 
$\dim B(n)=\dim V(n)$, in all cases. Since
$$
 \dim V(n)\leq \dim A(n)\leq\dim B(n)=\dim V(n),
$$
we get $\dim V(n)=\dim A(n)$; hence, $V(n)=A(n)$, completing the proof of \ref{paracharlemmaitem1}.

\ref{paracharlemmaitem2}. Let $n$ be an integer such that $n\geq N_{\pi\rtimes \sigma,s}$ and $n \geq 1$. 
Let $C(n)$ be the subspace of all $f \in V_s(n)$ satisfying \eqref{paracharlemma4}--\eqref{paracharlemma8}.
We need to prove that $V(n-1)=C(n)$. Let $f \in V(n-1)$; we will prove
that $f$ satisfies \eqref{paracharlemma4}--\eqref{paracharlemma8}.
Since $f \in V(n-1)$ we have that $f(gt_{n-1})=f(g)$ for all $g\in\GSp(4,F)$. 
Equation \eqref{paracharlemma4} follows by substituting $g=1$ and
using 
\begin{equation}
\label{parachareq009}
t_{n-1}=
\begin{bsmallmatrix}
\varpi&&&\\
&1&&\\
&&1&\\
&&&\varpi^{-1}
\end{bsmallmatrix}
t_n.
\end{equation}
Equation \eqref{paracharlemma5} follows by substituting $g=M_{i-1}$
and using \eqref{parachareq009} again. 
Equation \eqref{paracharlemma6} follows by substituting $g=M_0$
and using \eqref{parachareq009} and \eqref{rhopantilemmaeq01}.
Equation \eqref{paracharlemma7} follows by evaluating $f$ at both sides of the identity \eqref{PKsdisjointeq12}.
Equation \eqref{paracharlemma8} follows by evaluating $f$ at both sides of the identity
$$
 M_{(n-1)/2}=\begin{bsmallmatrix}&-\varpi^{(1-n)/2}&&\varpi^{1-n}\\\varpi^{(n-1)/2}&&-1\\&&&\varpi^{-(n-1)/2}\\&&-\varpi^{(n-1)/2}\end{bsmallmatrix}M_{(n-1)/2}t_{n-1}s_2.
$$
Since $f$ satisfies \eqref{paracharlemma4}--\eqref{paracharlemma8} it follows that $f \in C(n)$. Thus,
$V(n-1) \subset C(n)$. To prove that $V(n-1)=C(n)$ we will prove that $\dim V(n-1)=\dim C(n)$. The dimension
of $V(n-1)$ is known: by Theorem 5.2.2 of \cite{NF},  the minimal paramodular level of $\pi\rtimes\sigma$ is $N_{\pi\rtimes\sigma}=N_{\sigma\pi}+2a(\sigma)$, and
\begin{equation}
\label{paracharlemma4010}
\renewcommand{\arraystretch}{1.3}
    \dim V(n-1)
    =\left\{\begin{array}{l@{\qquad\mbox{if }}l}
    \displaystyle\Big\lfloor\frac{(n-N_{\pi\rtimes\sigma}+1)^2}4\Big\rfloor&n\geq N_{\pi\rtimes\sigma}+1,\\
    \;0&n<N_{\pi\rtimes\sigma}+1.\end{array}\right.
\end{equation}
We estimate the dimension of $C(n)$ as follows. For each integer $i$ such that
$a(\sigma) \leq i <n/2$ define $h_i \in \GSp(4,F)$ by
$$
h_i=
\begin{cases}
1&\text{if $i=0$},\\
M_i &\text{if $1<i<n/2$}.
\end{cases}
$$
Let $f \in C(n)$. An argument as in the proof of \ref{paracharlemmaitem1} shows that
\begin{equation}
\label{paracharlemma401}
(\sigma\pi)(
\begin{bsmallmatrix}
1&\\
&\varpi^{-i}
\end{bsmallmatrix})
f(h_i) \in
(\sigma \pi)^{\Gamma_0(\p^{n-1-2i})}
\quad\text{for $a(\sigma) \leq i < n/2$}.
\end{equation}
For this, note since $f \in V(n-1)$ we
have $f \in V_s(n-1)$.  Now let 
$$
D(n)=\bigoplus_{a(\sigma)\leq i < n/2} (\sigma \pi)^{\Gamma_0(\p^{n-1-2i})}.
$$
Define $\rho:C(n) \to D(n)$ by 
\begin{equation*}
\rho(f)= \bigoplus_{a(\sigma)\leq i < n/2}
(\sigma\pi)(
\begin{bsmallmatrix} 
1& \\ &\varpi^{-i}
\end{bsmallmatrix})
f(h_i).
\end{equation*}
Then $\rho$ is a well-defined linear map by \eqref{paracharlemma401}.
As in the proof of \ref{paracharlemmaitem1}, $\rho$ is injective. It follows that $\dim C(n) \leq \dim D(n)$. 
We will prove that $\dim D(n) = \dim V(n-1)$. 
As in the proof of \ref{paracharlemmaitem1}, we find that
$$
\dim D(n)
=
\sum_{a(\sigma) \leq i \leq (n-N_{\sigma \pi})/2 -1/2} n-N_{\sigma \pi} -2i.
$$
Assume first that this sum is empty. Then $\dim D(n)=0$, and $n-1<N_{\sigma\pi}+2a(\sigma)$. 
By \eqref{paracharlemma4010} we have $\dim V(n-1)=0$, so that $\dim D(n)=\dim V(n-1)=0$
in this case. Now assume that the sum is not empty.
A calculation shows that the sum is 
$$
\displaystyle\Big\lfloor\frac{(n-N_{\sigma\pi}-2a(\sigma)+1)^2}4\Big\rfloor.
$$
By \eqref{paracharlemma4010} this is $\dim V(n-1)$, so that 
$\dim D(n)=\dim V(n-1)$, in all cases. Since
$$
 \dim V(n-1)\leq \dim C(n)\leq\dim D(n)=\dim V(n-1),
$$
we get $\dim V(n-1)=\dim C(n)$; hence, $V(n-1)=C(n)$, completing the proof of~\ref{paracharlemmaitem2}.
\end{proof}

%% file: SKMS_chapter5.tex
\chapter{Dimensions}
\label{dimchap}

Let $(\pi,V)$ be an irreducible, admissible representation of $\GSp(4,F)$
with trivial central character. In this chapter we will compute the dimensions
of $V_s(n)$ and $\bar V_s(n)$ for all non-negative integers $n$;
the result appears in Theorem~\ref{dimensionstheorem}. 
To prove Theorem~\ref{dimensionstheorem} we will use three tools. The first
tool is the theory 
about induced representations from Chap.~\ref{inducedrepchap}. The second 
tool is an important inequality. For 
non-negative integers $n$ 
we will prove that
\begin{equation}\label{upperbound}
\dim V_s(n) \leq \dim V(n) + \dim V(n+1).
\end{equation}
As a corollary of this upper bound on $\dim V_s(n)$, we will show  that
$\pi$ admits non-zero stable Klingen vectors
if and only if $\pi$ is paramodular, and if
$\pi$ is paramodular, then 
$N_{\pi,s}$ is either $N_\pi-1$ or $N_\pi$.
Thus, paramodular representations fall into one of two classes. 
If $N_{\pi,s} = N_\pi-1$, then 
we will say that $\pi$ is a \index{category 1 representation}\catone representation; if $N_{\pi,s}=N_\pi$, then we will
say that $\pi$ is \index{category 2 representation}\cattwo representation. Our third and final 
tool for proving Theorem \ref{dimensionstheorem} is a certain stable Klingen
vector which we refer to as the shadow of a newform.

Theorem \ref{dimensionstheorem} has some significant immediate consequences. 
For example, in Corollary \ref{Lparametercharaterizationcor}
we will prove that $\pi$ is a \catone paramodular representation
if and only if the decomposition of the $L$-parameter of $\pi$
into indecomposable representations contains no unramified one-dimensional
factors.
In this chapter we will also deduce results indicating that if $\pi$ is a paramodular
representation and $N_\pi \geq 2$, then the stable Klingen vectors\index{stable Klingen vector}
in $\pi$ exhibit a simple structure.

\section{The upper bound}
In this section we prove the upper bound \eqref{upperbound}. 
Our approach uses the $P_3$-quotient recalled in Sect.~\ref{repsec} as well as the non-existence result from Sect.~\ref{nonexistencesec}.

\begin{lemma}\label{vsZJinjlemma}
 Let $(\pi,V)$ be an admissible representation of $\GSp(4,F)$ for which the center  acts trivially. Let $n$ be an integer such that $n \geq 0$. Let $p:V \to V_{Z^J}$ be the projection map. If $v \in V_s(n)$ is such that $p(v)=0$, then $v=0$.
\end{lemma}
\begin{proof}
Assume that $v \in V_s(n)$ and $p(v)=0$. Then there exists a positive integer $M \geq n$ such that 
$$
\int\limits_{\p^{-M}} \pi(\begin{bsmallmatrix} 1&&&x \\ &1&& \\ &&1& \\ &&&1 \end{bsmallmatrix}) v\, dx =0.
$$
This implies that $(\tau_{M} \tau_{M-1} \dots \tau_n)(v) =0$. Applying (1) of Lemma \ref{taunlevelraisinglemma}  repeatedly proves that $v=0$. 
\end{proof}

\begin{lemma}\label{pnzerolemmaeq1lemma}
 Let $(\pi,V)$ be an irreducible, admissible representation of the group $\GSp(4,F)$ with trivial central character. Assume that $\pi$ is not of one of the following types:
 \begin{itemize}
  \item\label{pnzerolemmaeq1lemmaitem1} $\chi_1\times\chi_2\rtimes\sigma$ of type I with $\chi_1,\chi_2,\sigma$ unramified and ($\chi_1=1$ or $\chi_2=1$);
  \item \label{pnzerolemmaeq1lemmaitem2} Type VIc with unramified $\sigma$;
  \item \label{pnzerolemmaeq1lemmaitem3} Type VId with unramified $\sigma$.
 \end{itemize}
 Let $n$ be an integer such that $n \geq 0$. Let $v \in V_s(n)$. If $\tau_{n+1} v = q^{-2} \eta v$, then $v=0$.
\end{lemma}
\begin{proof}
Assume that  $\tau_{n+1}v= q^{-2} \eta v$. Assume that $v\neq0$; we will obtain a contradiction. Applying the projection $p:V\to V_{Z^J}$, we get
\begin{equation}\label{pnzerolemmaeq1lemmaeq1}
 p(v)=q^{-2}\begin{bsmallmatrix}\varpi^{-1}\\&\varpi^{-1}\\&&1\end{bsmallmatrix}p(v).
\end{equation}
The vector $p(v)$ is also non-zero by Lemma \ref{vsZJinjlemma}. 
Hence $p(v)$ defines a non-zero vector $w$ in an irreducible subquotient $\tau$ of
the $P_3$-filtration of $V_{Z^J}$ (see Sect.~\ref{repsec}). This vector has the property that
\begin{equation}\label{pnzerolemmaeq1lemmaeq3}
 \tau(P_3(\OF))w=w
\end{equation}
and
\begin{equation}\label{pnzerolemmaeq1lemmaeq2}
 w=q^{-2}\tau(\begin{bsmallmatrix}\varpi^{-1}\\&\varpi^{-1}\\&&1\end{bsmallmatrix})w.
\end{equation}
Assume that $\tau = \tau_{\GL(0)}^{P_3}(1)$ or $\tau=\tau_{\GL(1)}^{P_3}(\chi)$
for some character $\chi$ of $F^\times$; we will obtain a contradiction. 
By the definition of $\tau_{\GL(0)}^{P_3}(1)$ or $\tau_{\GL(1)}^{P_3}(\chi)$
there exists a compact set $X \subset P_3$ such that if $y \in P_3$ is such 
that $w(y) \neq 0$, then $y \in 
\begin{bsmallmatrix}
*&*&*\\
&1&*\\
&&1
\end{bsmallmatrix}
X$. Fix $y \in P_3$  such that $w(y) \neq 0$. By \eqref{pnzerolemmaeq1lemmaeq2},
we also have $w(y\begin{bsmallmatrix} \varpi^{-k}&&\\&\varpi^{-k}&\\&&1 \end{bsmallmatrix})
\neq 0$ for all $k \in \Z$. It follows that 
$y\begin{bsmallmatrix} \varpi^{-k}&&\\&\varpi^{-k}&\\&&1 \end{bsmallmatrix}
\in 
\begin{bsmallmatrix}
*&*&*\\
&1&*\\
&&1
\end{bsmallmatrix} X$
for all $k \in \Z$. In particular, we see that there exists a compact subset
$K$ of $F$ such that the $(2,1)$ and $(2,2)$ entries of 
$y\begin{bsmallmatrix} \varpi^{-k}&&\\&\varpi^{-k}&\\&&1 \end{bsmallmatrix}$
lie in $K$ for all $k \in \Z$. This implies that $y_{21}=y_{22}=0$, a contradiction. 
It follows that
$\tau=\tau_{\GL(2)}^{P_3}(\rho)$ for some irreducible, admissible representation $\rho$ of $\GL(2,F)$. 
By \eqref{pnzerolemmaeq1lemmaeq3} the 
representation $\rho$ is unramified. Condition \eqref{pnzerolemmaeq1lemmaeq2} translates to
$\omega_\rho=\nu^2$.
Using Table A.5 of \cite{NF}, we see that 
the only representations that admit $\tau_{\GL(2)}^{P_3}(\rho)$ with such $\rho$ in their $P_3$-filtration are $\chi_1\times\chi_2\rtimes\sigma$ of type I with $\chi_1,\chi_2,\sigma$ unramified and ($\chi_1=1$ or $\chi_2=1$), or VIc or VId with unramified $\sigma$ (note  that there is typo in Table
A.5 of \cite{NF}: The entry for Vd in the ``s.s.($V_0/V_1$)'' column should be 
$\tau_{\GL(2)}^{P_3}(\nu (\nu^{-\frac{1}{2}} \sigma \times \nu^{-\frac{1}{2}}\xi \sigma))$) 
; by assumption, $\pi$ is not one of these representations. This is a contradiction.
\end{proof}

\begin{proposition}\label{funnyidentityprop}
Let $(\pi,V)$ be an irreducible, admissible representation of the group 
$\GSp(4,F)$ with trivial central character. Let $n$ be an integer such that $n\geq 0$. Assume that $v\in V_s(n)$ satisfies
 \begin{equation}\label{funnyidentitypropeq1}
  \tau_{n+1}v=q^{-2}\eta v.
 \end{equation}
 Then $v=0$.
\end{proposition}
\begin{proof}
By Lemma \ref{pnzerolemmaeq1lemma}, we may assume that $\pi$ is one of the following types:
\begin{itemize}
  \item $\chi_1\times\chi_2\rtimes\sigma$ of type I with $\chi_1,\chi_2,\sigma$ unramified and ($\chi_1=1$ or $\chi_2=1$);
  \item Type VIc with unramified $\sigma$;
  \item Type VId with unramified $\sigma$.
\end{itemize}
We claim that $\pi$ is a subrepresentation of  $\pi'\rtimes\sigma'$ for some admissible representation $\pi'$ of $\GL(2,F)$ satisfying the hypotheses of Lemma \ref{rhopantilemma} and character $\sigma'$ of $F^\times$ such that $\omega_{\pi'}\sigma'^2=1$.
This is clear in the first case because $\chi_1\times\chi_2$ is irreducible. Assume that $\pi$ is as in the second case, so that
$\pi = L(\nu^{\frac{1}{2}}\St_{\GL(2)},\nu^{-\frac{1}{2}}\sigma)$ for some
unramified character $\sigma$ of $F^\times$ with $\sigma^2=1$.
By \eqref{groupVItableeq}, $\pi$ is a quotient of $\nu^{\frac{1}{2}}\St_{\GL(2)}\rtimes\nu^{-\frac{1}{2}}\sigma$.
Taking contragredients and using $\pi^\vee \cong \pi$, we see that 
$\pi$ is a subrepresentation of $\nu^{-\frac{1}{2}}\St_{\GL(2)}\rtimes\nu^{\frac{1}{2}}\sigma^{-1}$.
Assume that $\pi$ is as in the third case, so 
that ${\pi \cong L(\nu,1_{F^\times} \rtimes \nu^{-\frac{1}{2}}\sigma)}$
for some unramified character of $F^\times$ such that $\sigma^2=1$. 
Then $\pi$ is a quotient of $\nu^{\frac{1}{2}}\triv_{\GL(2)}\rtimes\nu^{-\frac{1}{2}}\sigma$ by \eqref{groupVItableeq}. Taking contragredients, we see that $\pi$ is a subrepresentation of $\nu^{-\frac{1}{2}}\triv_{\GL(2)}\rtimes\nu^{\frac{1}{2}}\sigma^{-1}$. This proves our claim in all cases.

We now apply Lemma \ref{rhopantilemma} to conclude that $v=0$. 
\end{proof}

For the following result recall the paramodularization operator $p_n:V_s(n)\to V(n)$ defined in~\eqref{pndefeq}.
\begin{lemma}\label{pnzerolemma}
Let $(\pi,V)$ be a smooth representation of $\GSp(4,F)$ for which the center of $\GSp(4,F)$ acts trivially, and 
assume that the subspace of vectors of $V$ fixed by $\SSp(4,F)$ is trivial.
Let $n$ be an integer such that $n \geq0$. 
Define 
 $$
  c_n: V_s(n) \longrightarrow V(n) \oplus V(n+1)
 $$
 by
 $$
  c_n(v) = (p_n v , p_{n+1} \tau_nv)
 $$
for $v \in V_s(n)$. 
 Let $v \in V_s(n)$. If $c_n(v)=0$, then
 \begin{equation}\label{pnzerolemmaeq1}
  \tau_{n+1} v = q^{-2} \eta v.
 \end{equation}
\end{lemma}
\begin{proof}
Assume that $c_n(v)=0$. Then $p_n v =0$ and $p_{n+1}\tau_nv=0$. The result follows now from Lemma \ref{funnyidentitylemma5}.
\end{proof}

\begin{theorem}\label{upperboundtheorem}
Let $n$ be an integer such that $n \geq 0$. 
Let $(\pi,V)$ be an irreducible, admissible representation of $\GSp(4,F)$ with trivial central character. 
Then
 \begin{equation}\label{upperboundpropeq1}
  \dim V_s(n)\leq\dim V(n)+\dim V(n+1).
 \end{equation}
\end{theorem}
\begin{proof}
It suffices to show that the map $c_n$ defined in Lemma \ref{pnzerolemma} is injective. 
Assume first that $\pi$ is one-dimensional, so that $\pi=\chi 1_{\GSp(4,F)}$ for some
character $\chi$ of $F^\times$. If $\chi$ is ramified, then $V_s(n)=V(n)=V(n+1)=0$ so that
$c_n$ is injective. If $\chi$ is unramified, then $c_n(v)=(v,v)$ for $v \in V_s(n)$, so that
$c_n$ is injective. Assume that $\pi$ is not one-dimensional, so that $\pi$ is infinite-dimensional.
Then the subspace of vectors in $V$ fixed by $\SSp(4,F)$ is trivial.
Let $v\in V_s(n)$ be such that  $c_n(v)=0$. By   Lemma~\ref{pnzerolemma}, $\tau_{n+1}v=q^{-2}\eta v$. Hence $v=0$ by 
Proposition~\ref{funnyidentityprop}.
\end{proof}

\begin{corollary}\label{upperboundpropcor}
Let $(\pi,V)$ be an irreducible, admissible representation of $\GSp(4,F)$ with trivial central character. 
Then $\pi$ admits non-zero paramodular vectors if and only if $\pi$ admits non-zero stable Klingen vectors. If $\pi$ is paramodular, then 
 $N_{\pi,s}=N_\pi-1$ or $N_{\pi,s}=N_\pi$.
\end{corollary}
\begin{proof}
If $\pi$ is paramodular, then $\pi$ admits non-zero stable Klingen vectors because $V(n) \subset V_s(n)$
for all integers $n \geq 0$. Assume that $\pi$ admits non-zero stable Klingen vectors. Then there exists
an integer $n$ such that $n \geq 0$ and $V_s(n) \neq 0$. By \eqref{upperboundpropeq1}, $0< \dim V_s(n) \leq 
V(n)+\dim V(n+1)$, so that $\pi$ is paramodular. For the final assertion, assume that $\pi$ is paramodular.
Then
$$
0 < \dim V_s(N_{\pi,s}) \leq \dim V(N_{\pi,s}) + \dim V(N_{\pi,s}+1)
$$
by \eqref{upperboundpropeq1}. It follows that $N_\pi \leq N_{\pi,s}$ or $N_\pi \leq N_{\pi,s}+1$. Hence,
$N_\pi-1 \leq N_{\pi,s}$. We also have 
$$
0< \dim V(N_\pi) \leq \dim V_s(N_\pi)
$$
because $V(N_\pi) \subset V_s(N_\pi)$. This implies that $N_{\pi,s} \leq N_{\pi}$. We now have 
$N_\pi-1 \leq N_{\pi,s} \leq N_\pi$ so that $N_{\pi,s}=N_\pi-1$ or $N_{\pi,s}=N_\pi$. 
\end{proof}

\label{catdefpageanchor}
Let $(\pi,V)$ be an irreducible, admissible representation of $\GSp(4,F)$ with trivial
central character. Assume that $\pi$ is paramodular. By Corollary \ref{upperboundpropcor}
we have $N_{\pi,s} = N_\pi-1$ or $N_{\pi,s}=N_\pi$. If $N_{\pi,s} = N_\pi-1$, then 
we will say that $\pi$ is a \emph{\catone representation}; if $N_{\pi,s}=N_\pi$, then we will
say that $\pi$ is \emph{\cattwo representation}. \index{category 1 representation} \index{category 2 representation}
We note that if $\pi$ is a \catone representation, then $\bar N_{\pi,s}$ is also equal to $N_\pi-1$. 
Also, it is evident that if $\pi$ is a \catone representation, then $N_\pi  \geq 2$. Thus, we have
the following diagram:
$$
\begin{array}{lll}
 \!N_\pi\\
 0&\multirow{2}{*}{\bracetwo}&\multirow{2}{*}{\text{only \cattwo representations}}\\
 1\\
 2&\multirow{4}{*}{\bracethree}&\multirow{4}{*}{\text{\catone or \cattwo representations}}\\
 3\\
 \vdots\\
\end{array}
$$

\vspace{2ex}
\section{The shadow of a newform}

The following lemma defines the shadow of a newform and proves
some basic properties. As we shall see, the shadow of 
a newform is a key example of a stable Klingen vector. This vector
was originally introduced in Sect.~7.4 of \cite{NF}. In the following lemma, 
in \eqref{shadowgeneraldefeq3} and \eqref{shadowgeneraldefeq4}, the vector $v_{\mathrm{new}}$
is regarded as an element of $V_s(N_\pi+1)$, while in \eqref{Wsgenericlemma1eq1},
$v_{\mathrm{new}}$ is regarded as an element of $V_s(N_\pi)$.
We recall the level lowering operators $s_n,\sigma_n:V_s(n+1)\to V_s(n)$ defined in \eqref{sndefeq} and~\eqref{sigmaopseq}.

\begin{lemma}
\label{Wsgenericlemma1}
Let $(\pi,V)$ be an irreducible, admissible representation of the group $\GSp(4,F)$
with trivial central character. Assume that $\pi$ is paramodular and 
that $N_\pi \geq 2$. Let $v_{\mathrm{new}}$ be a newform for $\pi$, i.e.,
a non-zero element of the one-dimensional space $V(N_\pi)$. Define
\begin{equation}\label{shadowdefeq2}
v_{s}=q^3 s_{N_\pi-1} v_{\mathrm{new}},
\end{equation}
so that
\begin{equation}\label{shadowdefeq200}
v_{s}= q^3 \int\limits_{\OF}\int\limits_{\OF}\int\limits_{\OF}
\pi(
\begin{bsmallmatrix}
1\vphantom{a_b^c}&&&\\
x\varpi^{N_\pi-1}&1\\
y\varpi^{N_\pi-1}&&1\\
z\varpi^{N_\pi-1}&y\varpi^{N_\pi-1}&-x\varpi^{N_\pi-1}&1
\end{bsmallmatrix})v_{\rm new}\, dx\,dy\,dz.
\end{equation}
We refer to $v_s$ as the \emph{shadow}\index{shadow of a newform}
of the newform $v_{\mathrm{new}}$. 
The vector $v_{s}$ is contained in the space $V_s(N_\pi-1)$. Moreover, 
\begin{align}
t_{N_\pi-1}v_{s}&=q^3\sigma_{N_\pi}v_{\rm new}\label{shadowgeneraldefeq3},\\
\tau_{N_\pi-1}v_{s}&=-q^2\sigma_{N_\pi}v_{\rm new}\label{shadowgeneraldefeq4},\\
v_s&=-q^2\sigma_{N_\pi-1}v_{\rm new}, \label{Wsgenericlemma1eq1}\\
p_{N_\pi} \tau_{N_\pi-1} v_s &= -(q+q^2)^{-1} \mu_\pi v_{\mathrm{new}},\label{vsmulemmaeq1}\\
\rho_{N_\pi-1}'v_s&=q^{-1}\mu_\pi v_{\mathrm{new}},\label{vsmulemmaeq2}\\
-q^2\tau_{N_\pi-1}v_s & = T_{1,0}^s v_{\mathrm{new}}, \label{Wsgenericlemma100}\\
q^{-1} \mu_\pi v_{\mathrm{new}} &= q^3 \tau_{N_\pi} \tau_{N_\pi-1}\sigma_{N_\pi-1} v_{\mathrm{new}} - q \eta \sigma_{N_\pi-1} v_{\mathrm{new}}.
\label{vnewalluppereq}
\end{align}
\end{lemma}
\begin{proof}
The vector $v_s$ is contained in $V_s(N_\pi-1)$ by the definition of $s_{N_\pi-1}$. 
To prove \eqref{shadowgeneraldefeq3} we calculate:
\begin{align*}
t_{N_\pi-1} v_s & = q^3t_{N_\pi-1} s_{N_\pi-1} v_{\mathrm{new}}\\
&= q^3 \sigma_{N_\pi} t_{N_\pi} v_{\mathrm{new}}\qquad \text{(by \ref{sigmaopslemmaitem5} of Lemma \ref{sigmaopslemma})}\\
& = q^3 \sigma_{N_\pi}  v_{\mathrm{new}}\qquad \text{(since $v_{\mathrm{new}} \in V(N_\pi)$)}.
\end{align*}
Equation \eqref{shadowgeneraldefeq4} follows from \eqref{shadowgeneraldefeq3}
because $-q \tau_{N_\pi-1} v_s=t_{N_\pi-1} v_s$ (by \eqref{pnVsneq2}, since $p_{N_\pi-1}v_s=0$).
For \eqref{Wsgenericlemma1eq1}, we have:
\begin{align*}
 \tau_{N_\pi-1} v_s
& = -q^2 \sigma_{N_\pi}  v_{\mathrm{new}}\qquad \text{(by \eqref{shadowgeneraldefeq4})}\\
& =- q^2 \sigma_{N_\pi} \tau_{N_\pi} v_{\mathrm{new}}\qquad \text{(since $v_{\mathrm{new}} \in V(N_\pi)$)}\\
& = -q^2 \tau_{N_\pi-1} \sigma_{N_\pi-1} v_{\mathrm{new}}\qquad \text{(by \ref{sigmaopslemmaitem6} of Lemma \ref{sigmaopslemma})}.
\end{align*}
Hence, $\tau_{N_\pi-1}( q^{2}\sigma_{N_\pi-1} v_{\mathrm{new}} + v_s)=0$. Since $\tau_{N_\pi-1}: V_s(N_\pi-1)
\to V(N_\pi)$ is injective by Lemma \ref{taunlevelraisinglemma}, we get 
$q^2\sigma_{N_\pi-1} v_{\mathrm{new}} +  v_s=0$; this is \eqref{Wsgenericlemma1eq1}. 
To prove \eqref{vsmulemmaeq1},
let $dk$ be the Haar measure on $\GSp(4,F)$ that assigns $\K{n}$ volume $1$.
We have
\begin{align*}
&p_{N_\pi} \tau_{N_\pi-1} v_s
=-q^{-1} p_{N_\pi} t_{N_\pi-1} v_s \qquad \text{(by \eqref{pnVsneq2})}\\
&\qquad=-q^{2} p_{N_\pi} \sigma_{N_\pi} v_{\mathrm{new}}\qquad \text{(by \eqref{shadowgeneraldefeq3})}\\
&\qquad=-q^2 \int\limits_{\K{n}} \int\limits_{\OF}\int\limits_{\OF}\int\limits_{\OF}
\pi(k
\begin{bsmallmatrix}
1&x&y&z\varpi^{-N_\pi+1}\\
&1&&y\\
&&1&-x\\
&&&1
\end{bsmallmatrix}
\begin{bsmallmatrix}
\varpi\vphantom{\varpi^{-N_\pi+1}}&&&\\
&1\vphantom{\varpi^{-N_\pi+1}}&&\\
&&1\vphantom{\varpi^{-N_\pi+1}}&\\
&&&\varpi^{-1}\vphantom{\varpi^{-N_\pi+1}}
\end{bsmallmatrix})v_{\mathrm{new}}\, dx\,dy\,dz\,dk\\
&\qquad=-q^2 \int\limits_{\K{n}} 
\pi(k
\begin{bsmallmatrix}
\varpi^2&&&\\
&\varpi&&\\
&&\varpi&\\
&&&1
\end{bsmallmatrix})v_{\mathrm{new}}\,dk\\
&\qquad=-q^2(q^3+q^4)^{-1} T_{1,0}v_{\mathrm{new}}\qquad \text{(by Lemma 6.1.2 of \cite{NF})}\\
&\qquad=-(q+q^2)^{-1} \mu_\pi v_{\mathrm{new}}.
\end{align*}
This is \eqref{vsmulemmaeq1}. And
\begin{align*}
\rho_{N_\pi-1}'v_s&= (1+q^{-1})p_{N_\pi}t_{N_\pi-1}v_s \qquad \text{(by \eqref{rhoprimedefeq})}\\
&= -q(1+q^{-1})p_{N_\pi} \tau_{N_\pi-1} v_s \qquad \text{(by \eqref{pnVsneq2}, since $p_{N_\pi-1}v_s=0$)}\\
&= q^{-1} \mu_\pi v_{\mathrm{new}} \qquad \text{(by \eqref{vsmulemmaeq1})}.
\end{align*}
This proves \eqref{vsmulemmaeq2}. Next,
\begin{align*}
-q^2\tau_{N_\pi-1}v_s 
&=-q^2\tau_{N_\pi-1}(-q^2\sigma_{N_\pi-1}v_{\mathrm{new}})\qquad\text{(by \eqref{Wsgenericlemma1eq1})}\\
&=q^4 \tau_{N_\pi-1}\sigma_{N_\pi-1} v_{\mathrm{new}}\\
&=T_{1,0}^s v_{\mathrm{new}} \qquad \text{(by \eqref{heckeupdownlemmaeq1})}.
\end{align*}
This is \eqref{Wsgenericlemma100}. Finally,
\begin{align*}
q^{-1} \mu_\pi v_{\mathrm{new}}
&=\rho_{N_\pi-1}' v_s \qquad \text{(by \eqref{vsmulemmaeq2})} \\
&= q^{-1} \eta v_s + \tau_{N_\pi} t_{N_\pi-1} v_s \qquad \text{(by \eqref{rhopoplemmaeq2})}\\
&=-q\eta \sigma_{N_\pi-1} v_{\mathrm{new}} + q^3 \tau_{N_\pi} \sigma_{N_\pi} v_{\mathrm{new}} \qquad 
\text{(by \eqref{Wsgenericlemma1eq1} and \eqref{shadowgeneraldefeq3})}\\
&=-q\eta \sigma_{N_\pi-1} v_{\mathrm{new}} + q^3 \tau_{N_\pi} \sigma_{N_\pi}\tau_{N_\pi} v_{\mathrm{new}} \qquad 
\text{($v_{\mathrm{new}} \in V(N_\pi)$)}\\
&=-q\eta \sigma_{N_\pi-1} v_{\mathrm{new}} + q^3 \tau_{N_\pi} \tau_{N_\pi-1} \sigma_{N_\pi-1}v_{\mathrm{new}} \qquad 
\text{(\ref{sigmaopslemmaitem6} of Lemma \ref{sigmaopslemma})}.
\end{align*}
This is \eqref{vnewalluppereq}.
\end{proof}

With the notation and assumptions as in Lemma \ref{Wsgenericlemma1}, a basic observation
is that the shadow of a newform provides an example of a vector in $V_s(N_\pi-1)$. By 
Corollary \ref{upperboundpropcor} the stable Klingen level $N_{\pi,s}$ of $\pi$ is either $N_\pi-1$
or $N_\pi$. Thus, if the shadow of a newform is non-zero, then $N_{\pi,s} = N_\pi-1$.

\section{Zeta integrals and diagonal evaluation}
\label{zetaintsec}
We will sometimes use zeta integrals to investigate stable Klingen vectors in generic representations.
Let $(\pi,V)$ be a generic, irreducible, admissible representation of $\GSp(4,F)$ with trivial central character. 
Let $c_1,c_2 \in \OF^\times$, and let $V=\mathcal{W}(\pi,\psi_{c_1,c_2})$ be the Whittaker model of $\pi$
with respect to $\psi_{c_1,c_2}$ (see Sect.~\ref{repsec}). 
Let $n$ be an integer such that $n \geq 0$, and let $W \in V_s(n)$. Since $W$ is invariant under the elements
$$
\begin{bsmallmatrix}
1&x&&\\
&1&&\\
&&1&-x\\
&&&1
\end{bsmallmatrix}\quad\text{and}\quad
\begin{bsmallmatrix}
1&&&\\
&1&y&\\
&&1&\\
&&&1
\end{bsmallmatrix}
$$
for $x,y \in \OF$, it follows that for $a,b,c \in F^\times$, 
\begin{equation}
\label{Wabcvaneq}
W(
\begin{bsmallmatrix}
a&&&\\
&b&&\\
&&cb^{-1}&\\
&&&ca^{-1}
\end{bsmallmatrix})=0 \qquad \text{if $v(a)<v(b)$ or $2v(b)<v(c)$}.
\end{equation}
Let $i$ and $j$ be integers, and recall the element $\Delta_{i,j}$
defined in \eqref{deltaijdefeq}. By \eqref{Wabcvaneq} we see that
\begin{equation}
\label{delltavanisheq}
W(\Delta_{i,j})=0\qquad\text{if $i<0$ or $j<0$.}
\end{equation}
We recall from Sect.~\ref{repsec} that the zeta integral of 
$W$ is
\begin{equation}\label{zetaintdef}
Z(s,W)=\int\limits_{F^\times}\int\limits_FW
(\begin{bsmallmatrix}a\\&a\\&x&1\\&&&1\end{bsmallmatrix})|a|^{s-\frac{3}{2}}\,dx\,d^\times a.
\end{equation}
By Lemma 4.1.1 of \cite{NF}, the zeta integral of $W$
is given by the simplified formula
\begin{equation}\label{Vszetaintdef}
Z(s,W)=\int\limits_{F^\times}W(\begin{bsmallmatrix}a\\&a\\&&1\\&&&1\end{bsmallmatrix})|a|^{s-\frac{3}{2}}\,
d^\times a.
\end{equation}
Moreover, by \eqref{Wabcvaneq}, we even have that
\begin{equation}\label{Vszetaintposdef}
Z(s,W)=\int\limits_{\substack{F^\times\\v(a)\geq 0}}W(\begin{bsmallmatrix}a\\&a\\&&1\\&&&1\end{bsmallmatrix})|a|^{s-\frac{3}{2}}\,
d^\times a.
\end{equation}
The following lemma describes the zeta integrals of stable Klingen vectors obtained
by applying the operators of Sect.~\ref{levelraisingsec}.

\begin{lemma}
\label{zetarelationslemma}
Let $(\pi,V)$ be a generic, irreducible, admissible representation of $\GSp(4,F)$ with trivial central character. 
Let $c_1,c_2 \in \OF^\times$, and let $V=\mathcal{W}(\pi,\psi_{c_1,c_2})$ be the Whittaker model of $\pi$
with respect to $\psi_{c_1,c_2}$. Let $n$ be an integer such that $n \geq 0$ and let $W \in V_s(n)$. Then
\begin{align}
 \label{zetarelationseq1}Z(s,\tau_n W)&=Z(s,W),\\
 \label{zetarelationseq2}Z(s,\theta W)&=q^{-s+3/2}Z(s,W),\\
 \label{zetarelationseq3}Z(s,\eta W)&=0,\\
 \label{zetarelationseq5}Z(s,\rho'_n W)&=\begin{cases}
                                       Z(s,W)&\text{if }W\in V(n),\\
                                       -qZ(s,W)&\text{if }p_n(W)=0.
                                      \end{cases}
\end{align}
\end{lemma}
\begin{proof}
It is straightforward to verify \eqref{zetarelationseq1}.  For \eqref{zetarelationseq2},
we have by \eqref{thetadefeq3},
\begin{align*}
Z(s,\theta W)
&= \int\limits_{F^\times}W(\begin{bsmallmatrix}a\\&a\\&&1\\&&&1\end{bsmallmatrix}
\begin{bsmallmatrix}
1&&&\\
&1&&\\
&&\varpi&\\
&&&\varpi
\end{bsmallmatrix}
)|a|^{s-\frac{3}{2}}\,
d^\times a\\
&\quad+
q\int\limits_{\OF}\int\limits_{F^\times}W(\begin{bsmallmatrix}a\\&a\\&&1\\&&&1\end{bsmallmatrix}
\begin{bsmallmatrix}
1&&&\\
&1&c&\\
&&1&\\
&&&1
\end{bsmallmatrix}
\begin{bsmallmatrix}
1&&&\\
&\varpi&&\\
&&1&\\
&&&\varpi
\end{bsmallmatrix}
)|a|^{s-\frac{3}{2}}\,
d^\times a\\
&= q^{-s+\frac{3}{2}} Z(s,W)
+
q\int\limits_{\OF}\int\limits_{F^\times}\psi(c_2 a c) W(\begin{bsmallmatrix}a\\&a\varpi\\&&1\\&&&\varpi\end{bsmallmatrix}
)|a|^{s-\frac{3}{2}}\,
d^\times a.
\end{align*}
It follows from \eqref{Wabcvaneq} that
$$
W(\begin{bsmallmatrix}a\\&a\varpi\\&&1\\&&&\varpi\end{bsmallmatrix}
) =0
$$
for all $a \in F^\times$. 
This proves \eqref{zetarelationseq2}. The assertion \eqref{zetarelationseq3}
has a similar proof. Finally, 
for \eqref{zetarelationseq5}, observe \eqref{rhopoplemmaeq4} and \eqref{rhopoplemmaeq3}.
\end{proof}

To close this section we make some basic observations about stable Hecke operators
for generic representations.  Let $(\pi,V)$ be a generic, irreducible, admissible representation of the group $\GSp(4,F)$ with trivial central character. Let $V=\mathcal{W}(\pi,\psi_{c_1,c_2})$ be the Whittaker model of $\pi$; as usual, $c_1,c_2 \in \OF^\times$. Let
$i,j \in \Z$.
At various points in this work we will evaluate elements $W$ of $V_s(n)$ at $\Delta_{i,j}$ when
$n$ is an integer such that $n \geq 0$. As an initial observation we note that, by \eqref{delltavanisheq},
if $i<0$ or $j<0$, then $W(\Delta_{i,j})=0$. 
As concerns the stable Hecke operators we have the following lemma.

\begin{lemma}
\label{deltaWlemma}
Let $(\pi,V)$ be a generic, irreducible, admissible representation of the group $\GSp(4,F)$ with trivial central character. Let $V=\mathcal{W}(\pi,\psi_{c_1,c_2})$ be the Whittaker model of $\pi$; as usual, $c_1,c_2 \in \OF^\times$. Let $n$ be an integer such that $n \geq 0$, and let $W \in V_s(n)$. Let $i,j \in \Z$.  If $i,j \geq 0$, then
\begin{align}
(T_{0,1}^s W)(\Delta_{i,j}) &= q^2 W(\Delta_{i+1,j-1}) + q^3 W(\Delta_{i,j+1}),\label{deltaWlemmaeq1}\\
(T_{1,0}^s W)(\Delta_{i,j}) & = q^4 W(\Delta_{i+j,j}).\label{deltaWlemmaeq2}
\end{align}
\end{lemma}
\begin{proof}
These formulas follow from \eqref{Ts01eq} and \eqref{Ts10eq}.
\end{proof}

\section{Dimensions for some generic representations}
Let $(\pi,V)$ be a generic, irreducible, admissible representation of $\GSp(4,F)$
with trivial central character, and assume that $N_\pi \geq 2$. Under the 
assumption that $\mu_\pi\neq 0$, we can prove the following result about the structure
of $V_s(n)$ when $n \geq N_{\pi,s}=N_{\pi}-1$.  This theorem will be used in the proof
of the main result of this chapter. 

To prove this theorem we will need a result from \cite{NF} about the shadow of a newform in $\pi$ which we now recall. 
Let $V$ be the Whittaker model $\mathcal{W}(\pi,\psi_{c_1,c_2})$ of $\pi$,
let $W_{\mathrm{new}}$ be a newform for $V$, and 
let $W_s$ be the shadow of $W_{\mathrm{new}}$; we have $T_{1,0} W_{\mathrm{new}} = \mu_\pi W_{\mathrm{new}}$.
In Proposition 7.4.8 of \cite{NF} it is proven that 
\begin{equation}
\label{WsWnewzetaeq}
Z(s,W_s) = -q^{-2} \mu_\pi Z(s, W_{\mathrm{new}}). 
\end{equation}
Since $Z(s,W_{\mathrm{new}})$ is a non-zero multiple of $L(s,\pi)$ by Theorem 7.5.4 of \cite{NF},
it follows that if $\mu_\pi \neq 0$, then $W_s \neq 0$, so that $N_{\pi,s} = N_\pi -1$ by 
Corollary \ref{upperboundpropcor}, and $\pi$ is a \catone representation. 
In \eqref{genericdimensionstheoremeq1} of the following theorem, to simplify notation, we write $\tau$ for the level
raising operator $\tau_k: V_s(k)\to V_s(k+1)$ for each integer $k \geq 0$; we also note that by \eqref{commreleq1}
the operators $\tau$ and $\theta$ commute.

\begin{theorem}\label{genericdimensionstheorem}
Let $\pi$ be a generic, irreducible, admissible representation of the group $\GSp(4,F)$ with trivial central character. 
Assume that $N_\pi\geq2$. 
Let $W_{\mathrm{new}}$ be a newform for $\pi$, and let $W_s$ be the shadow of $W_{\mathrm{new}}$.
Assume that $\mu_\pi \neq 0$ so that $W_s \neq 0$, and $\pi$ is a \catone \index{category 1 representation} representation. 
Then
\begin{equation}\label{genericdimensionstheoremeq1}
V_s(n)=V(n-1)\oplus V(n)\oplus\bigoplus_{\substack{i,j\geq0\\i+j=n-N_\pi+1}}\C\tau^i\theta^{j}W_s
\end{equation}
 and
 \begin{equation}\label{genericdimensionstheoremeq5}
  \dim V_s(n)=\frac{(n-N_\pi+2)(n-N_\pi+3)}2
 \end{equation}
 and
 \begin{equation}\label{genericdimensionstheoremeq6}
  \dim \bar V_s(n)=n-N_\pi+2
 \end{equation}
 for all $n\geq N_{\pi,s}=\bar N_{\pi,s}=N_\pi-1$.
\end{theorem}
\begin{proof}
First we will prove that the sum on the right hand side of \eqref{genericdimensionstheoremeq1} is direct by induction on $n$. The case $n=N_{\pi,s}$ is clear.
Assume that $n\geq N_{\pi,s}+1$, and that the sum on the right-hand side of \eqref{genericdimensionstheoremeq1} is direct for $n-1$. Suppose that
\begin{equation}\label{genericdimensionstheoremeq2}
 W_1+W_2+\sum_{\substack{i,j\geq0\\i+j=n-N_\pi+1}}c_i\tau^i\theta^{j}W_s=0
\end{equation}
for some $W_1\in V(n-1)$, $W_2\in V(n)$ and complex numbers $c_0,\ldots, c_{n-N_\pi+1}$. Taking zeta integrals,
and using Lemma \ref{zetarelationslemma} and  \eqref{WsWnewzetaeq}, we get
\begin{equation}\label{genericdimensionstheoremeq3}
 Z(s,W_1)+Z(s,W_2)=q^{-2}\mu_\pi\sum_{\substack{i,j\geq0\\i+j=n-N_\pi+1}}(q^{-s+3/2})^jc_i Z(s,W_{\rm new}).
\end{equation}
By Theorem 7.5.6 of \cite{NF} we have
$$
 Z(s,W_1)\in(\C+\C q^{-s}+\cdots+\C(q^{-s})^{n-N_\pi-1})Z(s,W_{\rm new})
$$
and
$$
 Z(s,W_2)\in(\C+\C q^{-s}+\cdots+\C(q^{-s})^{n-N_\pi})Z(s,W_{\rm new})
$$
It therefore follows from 
\eqref{genericdimensionstheoremeq3} that $c_0=0$. 
From \eqref{genericdimensionstheoremeq2} we thus obtain
\begin{equation}\label{genericdimensionstheoremeq4}
 W_1+W_2+\tau_{n-1}\sum_{\substack{i,j\geq0\\i+j=n-N_\pi}}c_{i+1}\tau^i\theta^{j}W_s=0.
\end{equation}
The map $\tau_{n-1}:\bar V_s(n-1)\to\bar V_s(n)$ is injective
by Lemma \ref{taunlevelraisinglemma}. It follows that
$$
 \sum_{\substack{i,j\geq0\\i+j=n-N_\pi}}c_{i+1}\tau^i\theta^{j}W_s\in V(n-2)\oplus V(n-1).
$$
By the induction hypothesis, $c_1=\ldots=c_{n-N_\pi+1}=0$. Hence $W_1+W_2=0$, which implies $W_1=W_2=0$
by Theorem \ref{linindparatheorem}. 
We have proven that the sum on the right-hand side of \eqref{genericdimensionstheoremeq1} is direct.

Next, since the sum is direct, since $\dim V(m) = \lfloor (m-N_\pi+2)^2/4 \rfloor$ for integers
$m$ such that $m \geq N_\pi$ by Theorem 7.5.6 of \cite{NF}, and since $\lfloor k^2/4 \rfloor +
\lfloor (k+1)^2/4 \rfloor = k(k+1)/2$ for any integer $k$, 
\begin{align*}
 \dim V_s(n)&\geq \dim V(n-1)+\dim V(n)+n-N_\pi+2\\
 &=\Big\lfloor\frac{(n-N_\pi+1)^2}4\Big\rfloor+\Big\lfloor\frac{(n-N_\pi+2)^2}4\Big\rfloor+n-N_\pi+2\\
 &=\frac{(n-N_\pi+2)(n-N_\pi+3)}2.
\end{align*}
On the other hand, by Theorem \ref{upperboundtheorem}, 
\begin{align*}
 \dim V_s(n)&\leq \dim V(n)+\dim V(n+1)\\
 &=\Big\lfloor\frac{(n-N_\pi+2)^2}4\Big\rfloor+\Big\lfloor\frac{(n-N_\pi+3)^2}4\Big\rfloor\\
 &=\frac{(n-N_\pi+2)(n-N_\pi+3)}2.
\end{align*}
The equality \eqref{genericdimensionstheoremeq5} follows. It is also now clear
that we have equality in \eqref{genericdimensionstheoremeq1} and that
\eqref{genericdimensionstheoremeq6} holds. 
\end{proof}

A similar result for arbitrary generic representations will be proven in Theorem~\ref{W0theorem}.
\section{Dimensions for some non-generic representations}

In this section we will prove a number of lemmas about the 
dimensions of the spaces of stable Klingen vectors in various
non-generic representations. These results will be used in the 
proof of Theorem \ref{dimensionstheorem}, the main result of this
chapter.

\begin{lemma}\label{IIIbdimlemma}
Let $(\pi,V)$ be a representation of type IIIb, so that there exist
characters $\chi$ and $\sigma$ of $F^\times$ with
$\chi\notin\{1,\nu^{\pm2}\}$ such that $\pi=\chi\rtimes\sigma\triv_{\GSp(2)}$. Assume that $\pi$ has trivial central
character, i.e., $\chi\sigma^2=1$. 
Then, for all integers $n$ such that $n\geq0$,
$$
\dim V_s(n)=
\begin{cases}
2n+1&\text{if $\sigma$ is unramified},\\
0&\text{if $\sigma$ is ramified}.
\end{cases}
$$
\end{lemma}
\begin{proof}
If $\sigma$ is ramified, then $\pi$ is not paramodular by Theorem 3.4.3 of \cite{NF};
hence $V_s(n)=0$ for all integers  $n \geq 0$ by Corollary \ref{upperboundpropcor}.
Assume that $\sigma$ is unramified. 
Then $\chi$ is also unramified. 
By Sect.~\ref{repsec} there is  an exact sequence
$$
 0\longrightarrow \chi \rtimes \sigma \St_{\GSp(2)} \longrightarrow\chi\times\nu\rtimes\nu^{-\frac{1}{2}}\sigma\longrightarrow
\pi  \longrightarrow0.
$$
By Theorem \ref{Siegelindleveltheorem}, $\chi\times\nu\rtimes\nu^{-\frac{1}{2}}\sigma$ has  stable Klingen level $0$, and 
the dimension of the space of stable Klingen vectors of level $\p^n$ in this representation is
$(n^2+5n+2)/2$ for integers $n \geq 0$. 
By Theorem \ref{genericdimensionstheorem} and Table A.14 of \cite{NF}, the IIIa representation $\chi \rtimes \sigma \St_{\GSp(2)}$ has stable Klingen level $1$, and
the dimension of the space of stable Klingen vectors of level $\p^n$ in this representation is
$n(n+1)/2$ for integers $n \geq 1$. 
The statement of the lemma now follows from the exact sequence. 
\end{proof}

\begin{lemma}\label{IVbdimlemma}
Let $(\pi,V)$ be a representation of type IVb, so that
there exists a character $\sigma$ of $F^\times$ such that
$\pi=L(\nu^2,\nu^{-1}\sigma \St_{\GSp(2)})$.
Assume that $\pi$ has trivial central character, i.e., $\sigma^2=1$.
Then, for all integers $n \geq 0$,
$$
\dim V_s(n)=
\begin{cases}
n&\text{if $\sigma$ is unramified},\\
0&\text{if $\sigma$ is ramified}.
\end{cases}
 $$
\end{lemma}
\begin{proof}
If $\sigma$ is ramified, then $\pi$ is not paramodular by Theorem 3.4.3 of \cite{NF};
hence $V_s(n)=0$ for all integers  $n \geq 0$ by Corollary \ref{upperboundpropcor}.
Assume that $\sigma$ is unramified.
By \eqref{groupIVtableeq} there is an exact sequence
$$
0
\longrightarrow
\pi
\longrightarrow
\nu^{\frac{3}{2}}1_{\GL(2)}\rtimes\nu^{-\frac{3}{2}}\sigma
\longrightarrow 
\sigma 1_{\GSp(4)}
\longrightarrow0.
$$
By Theorem \ref{Siegelindleveltheorem}, $\nu^{\frac{3}{2}}1_{\GL(2)}\rtimes \nu^{-\frac{3}{2}}\sigma$ has  stable Klingen level $0$, and 
the dimension of the space of stable Klingen vectors of level $\p^n$ in this representation is
$n+1$ for integers $n \geq 0$. It is clear that $\sigma 1_{\GSp(4)}$
has stable Klingen level $0$, and the dimension of the space of stable Klingen vectors of level $\p^n$ in this representation is $1$ for integers $n \geq 0$. The statement of the lemma now follows from the exact sequence.
\end{proof}

\begin{lemma}\label{IVcdimlemma}
Let $(\pi,V)$ be a representation of type IVc, so that 
there exists a character $\sigma$ of $F^\times$ such that
$\pi=L(\nu^{\frac{3}{2}}\St_{\GL(2)},\nu^{-\frac{3}{2}}\sigma)$. 
Assume that $\pi$ has trivial central character, i.e., $\sigma^2=1$. 
Then, for all integers $n$ such that $n\geq0$,
$$
\dim V_s(n)=
\begin{cases}
2n&\text{if $\sigma$ is unramified},\\
0&\text{if $\sigma$ is ramified}.
\end{cases}
 $$
\end{lemma}
\begin{proof} 
If $\sigma$ is ramified, then $\pi$ is not paramodular by Theorem 3.4.3 of \cite{NF};
hence $V_s(n)=0$ for all integers  $n \geq 0$ by Corollary \ref{upperboundpropcor}.
Assume that $\sigma$ is unramified.
By \eqref{groupIVtableeq} there is an exact sequence
$$
0
\longrightarrow
\sigma\St_{\GSp(4)} 
\longrightarrow
\nu^{\frac{3}{2}}\St_{\GL(2)}\rtimes\nu^{-\frac{3}{2}}\sigma
\longrightarrow 
\pi 
\longrightarrow0.
$$
By Theorem \ref{Siegelindleveltheorem}, $\nu^{\frac{3}{2}}\St_{\GL(2)}\rtimes \nu^{-\frac{3}{2}}\sigma$ has  stable Klingen level $1$, and 
the dimension of the space of stable Klingen vectors of level $\p^n$ in this representation is
$n(n+3)/2$ for integers $n \geq 1$.
By Theorem \ref{genericdimensionstheorem} and Table A.14 of \cite{NF}, the IVa representation 
$\sigma\St_{\GSp(4)}$ has  stable Klingen level $2$, and 
the dimension of the space of stable Klingen vectors of level $\p^n$ in this representation is
$(n-1)n/2$ for integers $n \geq 2$.
The statement of the lemma now follows from the exact sequence. 
\end{proof}

\begin{lemma}\label{Vddimlemma}
Let $(\pi,V)$ be a representation of type Vd,
so that there exist a non-trivial quadratic character $\xi$ of $F^\times$ and  a
character $\sigma$ of $F^\times$ such that 
$\pi=L(\nu\xi,\xi\rtimes\nu^{-\frac{1}{2}}\sigma)$. 
Assume that $\pi$ as trivial central character, i.e., $\sigma^2=1$. 
Then, for all integers $n$ such that $n\geq0$,
 $$
  \dim V_s(n)=\begin{cases}
              1&\text{ if $\xi$ and $\sigma$ are unramified},\\
              0&\text{ if $\xi$ or $\sigma$ is ramified}.
             \end{cases}
 $$
\end{lemma}
\begin{proof}
If $\xi$ or $\sigma$ is ramified, then $\pi$ is not paramodular by Theorem 3.4.3 of \cite{NF};
hence $V_s(n)=0$ for all integers  $n \geq 0$ by Corollary \ref{upperboundpropcor}.
Assume that $\xi$ and $\sigma$ are unramified. By Table A.12 of \cite{NF}
$$
  \dim V(n)=\begin{cases}
              1&\text{ if $n$ is even},\\
              0&\text{ if $n$ is odd},
             \end{cases}
$$
for all integers $n\geq0$. 
Since $V_s(0)=V(0)$ and $V_s(n)$ contains $V(n-1)\oplus V(n)$ for integers $n \geq 1$, 
it follows that $\dim V_s(n)\geq1$ for all integers $n\geq0$. 
By Theorem~\ref{upperboundtheorem}, we also have $\dim V_s(n)\leq1$ for all integers $n\geq0$. 
This completes the proof.
\end{proof}
\begin{lemma}\label{Vbdimlemma}
Let $(\pi,V)$ be a representation of type Vb, 
so that there exist a non-trivial quadratic character $\xi$ of $F^\times$ and  a
character $\sigma$ of $F^\times$ such that
$\pi=L(\nu^{\frac{1}{2}}\xi\St_{\GL(2)},\nu^{-\frac{1}{2}}\sigma)$. 
Assume that $\pi$ has trivial central character, i.e., $\sigma^2=1$. 
Then, for all integers $n$ such that $n\geq0$,
 $$
  \dim V_s(n)=\begin{cases}
              n&\text{ if $\sigma$ and $\xi$ are unramified},\\
              n-2a(\xi)+1&\text{ if $n\geq2a(\xi)$, $\sigma$ is unramified, and $\xi$ is ramified},\\
              0&\text{ if $n<2a(\xi)$, $\sigma$ is unramified, and $\xi$ is ramified},\\
              0&\text{ if $\sigma$ is ramified}.
             \end{cases}
 $$
\end{lemma}
\begin{proof}
By \eqref{groupVtableeq} there is an exact sequence
$$
0
\longrightarrow 
\pi 
\longrightarrow
\nu^{\frac{1}{2}}\xi1_{\GL(2)}\rtimes\xi\nu^{-\frac{1}{2}}\sigma
\longrightarrow 
L(\nu \xi, \xi \rtimes \nu^{-\frac{1}{2}}\sigma)
\longrightarrow
0.
$$
The dimensions of the spaces of stable Klingen vectors of level $\p^n$ for integers $n \geq 0$ in 
 $L(\nu \xi, \xi \rtimes \nu^{-\frac{1}{2}}\sigma)$  and in $\nu^{\frac{1}{2}}\xi1_{\GL(2)}\rtimes\nu^{-\frac{1}{2}}\xi\sigma$ are
 given by Lemma \ref{Vddimlemma} and  Theorem \ref{Siegelindleveltheorem}, respectively. The statement of the lemma
now follows  from the exact sequence.
\end{proof}

\begin{lemma}\label{Vcdimlemma}
Let $(\pi,V)$ be a representation of type Vc, 
so that there exist a non-trivial quadratic character $\xi$ of $F^\times$ and  a
character $\sigma$ of $F^\times$ such that
$\pi=L(\nu^{\frac{1}{2}}\xi\St_{\GL(2)},\xi\nu^{-\frac{1}{2}}\sigma)$. 
Assume that $\pi$ has trivial central character, i.e., $\sigma^2=1$. 
Then, for all integers $n$ such that $n\geq0$,
 $$
  \dim V_s(n)=\begin{cases}
              n&\text{ if $\sigma$ and $\xi$ are unramified},\\
              n-2a(\xi)+1&\text{ if $n\geq2a(\xi)$, $\xi\sigma$ is unramified, and $\xi$ is ramified},\\
              0&\text{ if $n<2a(\xi)$, $\xi\sigma$ is unramified, and $\xi$ is ramified},\\
              0&\text{ if $\xi\sigma$ is ramified}.
             \end{cases}
 $$
\end{lemma}
\begin{proof}
This follows from Lemma \ref{Vbdimlemma}.
\end{proof}

\begin{lemma}\label{VIcdimlemma}
Let $(\pi,V)$ be a representation of type VIc, so that 
there exists a character $\sigma$ of $F^\times$ such that 
$\pi=L(\nu^{\frac{1}{2}}\St_{\GL(2)},\nu^{-\frac{1}{2}}\sigma)$. 
Assume that $\pi$ has trivial central character, i.e., $\sigma^2=1$.
Then, for all integers $n$ such that $n\geq0$,
$$
\dim V_s(n)=
\begin{cases}
n&\text{if $\sigma$ is unramified},\\
0&\text{if $\sigma$ is ramified}.
\end{cases}
$$
\end{lemma}
\begin{proof}
If $\sigma$ is ramified, then $\pi$ is not paramodular by Theorem 3.4.3 of \cite{NF};
hence $V_s(n)=0$ for all integers  $n \geq 0$ by Corollary \ref{upperboundpropcor}.
Assume that $\sigma$ is unramified.
By \eqref{groupVItableeq} there is an exact sequence
$$
0\longrightarrow
\tau(S,\nu^{-\frac{1}{2}}\sigma) 
\longrightarrow
\nu^{\frac{1}{2}}\St_{\GL(2)}\rtimes\nu^{-\frac{1}{2}}\sigma
\longrightarrow
\pi
\longrightarrow
0.
$$
By Theorem \ref{Siegelindleveltheorem}, 
$\nu^{\frac{1}{2}}\St_{\GL(2)}\rtimes\nu^{-\frac{1}{2}}\sigma$ 
has stable Klingen level $1$,
and the dimension of the space of stable Klingen vectors of level $\p^n$ in this representation is
$n(n+3)/2$ for integers $n \geq 1$.
By Theorem \ref{genericdimensionstheorem} and Table A.14 of \cite{NF}, the VIa representation 
$\tau(S,\nu^{-\frac{1}{2}}\sigma)$ has  stable Klingen level $1$, and 
the dimension of the space of stable Klingen vectors of level $\p^n$ in this representation is
$n(n+1)/2$ for integers $n \geq 1$.
The statement of the lemma follows now from the exact sequence.
\end{proof}

\begin{lemma}\label{VIddimlemma}
Let $(\pi,V)$ be a representation of type VId, so that 
there exists a character $\sigma$ of $F^\times$ such that 
$\pi=L(\nu,1_{F^\times}\rtimes\nu^{-\frac{1}{2}}\sigma)$. 
Assume that $\pi$ has trivial central character, i.e., $\sigma^2=1$.
Then, for all integers $n$ such that $n\geq0$,
$$
\dim V_s(n)=\begin{cases}
              n+1&\text{ if $\sigma$ is unramified},\\
              0&\text{ if $\sigma$ is ramified}.
             \end{cases}
 $$
\end{lemma}
\begin{proof}
If $\sigma$ is ramified, then $\pi$ is not paramodular by Theorem 3.4.3 of \cite{NF};
hence $V_s(n)=0$ for all integers  $n \geq 0$ by Corollary \ref{upperboundpropcor}.
Assume that $\sigma$ is unramified.
By \eqref{groupVItableeq} there is an exact sequence
$$
0
\longrightarrow
\tau(T,\nu^{-\frac{1}{2}}\sigma)
\longrightarrow
\nu^{\frac{1}{2}}1_{\GL(2)}\rtimes\nu^{-\frac{1}{2}}\sigma
\longrightarrow
\pi
\longrightarrow
0.
$$
By Theorem 3.4.3 of \cite{NF} the IVb representation is not paramodular,
and hence has no non-zero stable Klingen vectors of any level. 
It follows from the exact sequence that, for integers $n \geq 0$, the dimension of $V_s(n)$ 
is the same as the dimension
of the space of stable Klingen vectors of level $\p^n$ in 
$\nu^{\frac{1}{2}}1_{\GL(2)}\rtimes\nu^{-\frac{1}{2}}\sigma $;
by Theorem \ref{Siegelindleveltheorem} this dimension
is $n+1$. 
\end{proof}

\begin{lemma}\label{XIbdimlemma}
Let $(\pi,V)$ be a representation of type XIb, so that 
there exists an irreducible, admissible, supercuspidal representation $\tau$ of $\GL(2,F)$ 
with trivial central character and a character $\sigma$ of $F^\times$
such that $\pi=L(\nu^{\frac{1}{2}}\tau,\nu^{-\frac{1}{2}}\sigma)$.
Assume that $\pi$ has trivial central character, i.e., $\sigma^2=1$.
Then, for all integers $n$ such that $n\geq0$,
$$
  \dim V_s(n)=\begin{cases}
               n-N_\tau+1&\text{if $n\geq N_\tau$ and $\sigma$ is unramified},\\
               0&\text{if $n<N_\tau$ or $\sigma$ is ramified}.
              \end{cases}
$$
\end{lemma}
\begin{proof}
If $\sigma$ is ramified, then $\pi$ is not paramodular by Theorem 3.4.3 of \cite{NF};
hence $V_s(n)=0$ for all integers  $n \geq 0$ by Corollary \ref{upperboundpropcor}.
Assume that $\sigma$ is unramified.
By Sect.~\ref{repsec} there is an exact sequence
$$
0
\longrightarrow
\delta(\nu^{\frac{1}{2}}\tau,\nu^{-\frac{1}{2}}\sigma)
\longrightarrow
\nu^{\frac{1}{2}}\tau\rtimes\nu^{-\frac{1}{2}}\sigma
\longrightarrow
\pi
\longrightarrow
0.
$$
By Theorem \ref{Siegelindleveltheorem}, 
$\nu^{\frac{1}{2}}\tau\rtimes\nu^{-\frac{1}{2}}\sigma$ has stable Klingen level $N_\tau$,
and the dimension of the space of stable Klingen vectors of level $\p^n$ in this representation is
$ (n-N_\tau+1)(n-N_\tau+4)/2$ for integers $n \geq N_\tau$.
By Theorem \ref{genericdimensionstheorem} and Table A.14 of \cite{NF}, the XIa representation 
$\delta(\nu^{\frac{1}{2}}\tau,\nu^{-\frac{1}{2}}\sigma)$ has stable Klingen level $N_\tau$, and 
the dimension of the space of stable Klingen vectors of level $\p^n$ in this representation is
$(n-N_\tau+1)(n-N_\tau+2)/2$ for integers $n \geq N_\tau$.
The statement of the lemma follows now from the exact sequence.
\end{proof}
\section{The table of dimensions}

We will now prove the main result of this chapter. 

\begin{theorem}\label{dimensionstheorem}
For every irreducible, admissible representations $(\pi,V)$ of the group $\GSp(4,F)$ with trivial central character, 
the stable Klingen level $N_{\pi,s}$, 
the dimensions of the spaces $V_s(n)$,
the quotient stable Klingen level $\bar N_{\pi,s}$, and
the dimensions of the spaces $\bar V_s(n)$
are given in Table \ref{dimensionstable}. 
\end{theorem}
\begin{proof}
We begin by noting that Table A.12 of \cite{NF}
indicates, for every $\pi$, whether or not $\pi$ is paramodular,
and lists $N_\pi$ in the case that $\pi$ is paramodular. If $\pi$
is not paramodular, then this is indicated as such in Table \ref{dimensionstable}; if $\pi$ is paramodular and $\bar N_{\pi,s}$ is not defined, then the entry for $\bar N_{\pi,s}$ is $-$.
Also, in this proof we will sometimes use Theorem \ref{genericdimensionstheorem};
the hypothesis $\mu_\pi \neq 0$ 
used in this theorem can be verified using Table A.14 of \cite{NF} 
(this information also appears in Table \ref{levelsandeigenvaluestable} of this work).

The entries in the table are now verified as follows.
For Groups I and II, the entries follow from Theorem~\ref{Siegelindleveltheorem}.
For Group IIIa, the entries  follow from Theorem~\ref{genericdimensionstheorem}.
For Group IIIb, the entries follow from Lemma~\ref{IIIbdimlemma}.
For Group IVa, the entries  follow from Theorem~\ref{genericdimensionstheorem}.
For Group IVb, the entries follow Lemma~\ref{IVbdimlemma}.
For Group IVc, the entries follow Lemma~\ref{IVcdimlemma}.
For Group IVd, the entries are easily verified as $\pi$ is a twist of the trivial representation.
For Group Va, the entries  follow from Theorem~\ref{genericdimensionstheorem}.
For Group Vb, the entries follow Lemma~\ref{Vbdimlemma}.
For Group Vc, the entries follow Lemma~\ref{Vcdimlemma}.
For Group Vd, the entries follow Lemma~\ref{Vddimlemma}.
For Group VIa, the entries  follow from Theorem~\ref{genericdimensionstheorem}.
For Group VIb, we note that these representations are never paramodular.
For Group VIc, the entries follow Lemma~\ref{VIcdimlemma}.
For Group VId, the entries follow Lemma~\ref{VIddimlemma}.
For Groups VII and VIIIa, the entries  follow from Theorem~\ref{genericdimensionstheorem}.
For Group VIIIb, we note that these representations are never paramodular.
For Group IXa, the entries  follow from Theorem~\ref{genericdimensionstheorem}.
For Group IXb, we note that these representations are never paramodular.
For Group X, the entries follow from Theorem~\ref{Siegelindleveltheorem}.
For Group XIa, the entries follow from Theorem~\ref{genericdimensionstheorem}.
For Group XIb, the entries follow Lemma~\ref{XIbdimlemma}. 
If $\pi$ is generic and supercuspidal, the entries follow from Theorem~\ref{genericdimensionstheorem};
the assertion that $a=N_\pi \geq 4$ will be proven in Theorem~\ref{Npi4theorem}.
Finally, if $\pi$ is non-generic and supercuspidal, then $\pi$ is never paramodular.
\end{proof}

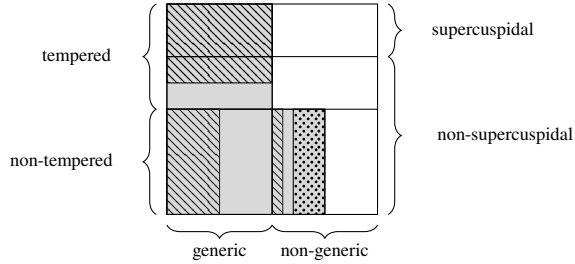
\begin{figure}
\caption{Paramodular representations with categories}
\label{pararepcatsfig}
\bigskip
\begin{tikzpicture}[scale=0.7]
\draw (0,0) rectangle (4,4);
\draw (0,3) -- (4,3);
\draw (0,2) -- (4,2);
\draw (2,0) -- (2,4);
\draw [decorate,decoration={brace,amplitude=5pt,raise=4pt},yshift=0pt]
(0,2) -- (0,4) node [black,midway,xshift=-1.2cm] {\scriptsize tempered};
\draw [decorate,decoration={brace,amplitude=5pt,raise=4pt},yshift=0pt]
(0,0) -- (0,2) node [black,midway,xshift=-1.4cm] {\scriptsize non-tempered};
\draw [decorate,decoration={brace,amplitude=5pt,raise=4pt},yshift=0pt]
(4,4) -- (4,3) node [black,midway,xshift=1.4cm] {\scriptsize supercuspidal};
\draw [decorate,decoration={brace,amplitude=5pt,raise=4pt},yshift=0pt]
(4,3) -- (4,0) node [black,midway,xshift=1.7cm] {\scriptsize non-supercuspidal};
\draw [decorate,decoration={brace,amplitude=5pt,raise=4pt},yshift=0pt]
(2,0) -- (0,0) node [black,midway,yshift=-0.5cm] {\scriptsize generic};
\draw [decorate,decoration={brace,amplitude=5pt,raise=4pt},yshift=0pt]
(4,0) -- (2,0) node [black,midway,yshift=-0.5cm] {\scriptsize non-generic};
\draw[fill=gray, fill opacity=0.3] (0,2) rectangle (2,4);
\draw[pattern=north west lines, pattern] (0,3) rectangle (2,4);
\draw[pattern=north west lines, pattern] (0,2.5) rectangle (2,3);
\draw[pattern=north west lines, pattern] (0,0) rectangle (1,2);
\draw[pattern=north west lines, pattern] (2,0) rectangle (2.2,2);
\draw[fill=gray, fill opacity=0.3] (2,0) rectangle (3,2);
\draw[fill=gray, fill opacity=0.3] (0,0) rectangle (2,2);
\draw[pattern=crosshatch dots] (2.4,0) rectangle (3,2);
\end{tikzpicture}\\
A schematic diagram of the paramodular representations among
all irreducible, admissible representations of $\GSp(4,F)$ with trivial central
character. Paramodular representations are in gray and \catone paramodular representations \index{category 1 representation}are hatched.
The paramodular Saito-Kurokawa
representations are  dotted; all of these are \index{category 2 representation} \cattwo paramodular representations.
\end{figure}

Using Theorem \ref{dimensionstheorem} we can update the schematic diagram of 
paramodular representations from Fig.~\ref{pararepsfig} on p.~\pageref{pararepsfig}. The 
result is shown in Fig.~\ref{pararepcatsfig} on p.~\pageref{pararepcatsfig}.
We note that every paramodular Saito-Kurokawa representation is a \cattwo paramodular
representation. Also, the only non-generic \catone representations are the IVb 
representations with $\sigma$ unramified. 
For the convenience of the reader we also include Table \ref{levelsandeigenvaluestable} 
on p.~\pageref{levelsandeigenvaluestable}. This table lists the irreducible, admissible
representations $\pi$ of $\GSp(4,F)$ with trivial central character; if $\pi$ is paramodular,
then the table lists the paramodular level $N_\pi$, the Atkin-Lehner eigenvalue
$\varepsilon_\pi$, the paramodular Hecke eigenvalues $\lambda_\pi$ and $\mu_\pi$, whether 
$\pi$ is cagegory 1 or category 2, and any comments about $\pi$. 

\section{Some consequences}

To conclude this chapter we will derive some consequences 
from Theorem \ref{dimensionstheorem}. 
Our first result provides a remarkable characterization of \catone paramodular representations.

\begin{corollary}\label{Lparametercharaterizationcor}
Let $(\pi,V)$ be an irreducible, admissible representation of $\GSp(4,F)$ with trivial central character. 
Assume that $\pi$ is paramodular. 
Then the following are equivalent.
 \begin{enumerate}
  \item \label{Lparametercharaterizationcoritem1} $\pi$ is a \catone paramodular representation, i.e., $N_{\pi,s}=N_\pi-1$.
  \item \label{Lparametercharaterizationcoritem2} The decomposition of the $L$-parameter of $\pi$ into indecomposable representations contains no unramified one-dimensional factors.
 \end{enumerate}
\end{corollary}
\begin{proof}
For non-supercuspidal representations, the equivalence of \ref{Lparametercharaterizationcoritem1} and \ref{Lparametercharaterizationcoritem2} follows by inspecting Table \ref{dimensionstable} from Theorem \ref{dimensionstheorem} and Table A.7 of \cite{NF}. Non-generic supercuspidals do not admit paramodular vectors by Theorem 3.4.3 of \cite{NF}. Assume that $\pi$ is a generic and supercuspidal. Then $N_{\pi,s}=N_\pi-1$ by Table \ref{dimensionstable}. The $L$-parameter of $\pi$ is a discrete series parameter, and thus does not factor through the Levi of any proper parabolic of $\SSp(4,\C)$. Suppose that the decomposition of the $L$-parameter $(\varphi,W)$ contains a one-dimensional factor, so that $W=V_1\oplus V_2$ for some $W_F'$-subspaces $V_1,V_2$ of $W$ with $\dim(V_1)=1$; we will obtain a contradiction. Let $\langle\cdot,\cdot\rangle$ be the non-degenerate symplectic form on $W$ preserved by $\varphi$. Let $x\in V_1$ be non-zero, so that $V_1=\C x$. Choose a non-zero $y\in V_2$ such that $\langle y,v\rangle=0$ for all $v\in V_2$; such a $y$ exists because $V_2$ is odd-dimensional. Since $W$ is non-degenerate, we may assume that $\langle x,y\rangle=1$. Let $U=(\C x)^\perp\cap V_2$. We have
\begin{align*}
 \dim(U)&=\dim((\C x)^\perp)+\dim(V_2)-\dim\big((\C x)^\perp+V_2\big)\\
 &=3+3-4=2.
\end{align*}
Clearly, $y\notin U$, so that $W=\C x\oplus U\oplus\C y$ with $(\C x+\C y)\perp U$. Let $x',y'$ be a basis of $U$ such $\langle x',y'\rangle=1$. Then $x,x',y',y$ is an ordererd basis of $W$, with respect to which the symplectic form is $J$, the matrix defined in \eqref{Jdefeq}. We now have that
$$
 \varphi(W_F')\subset\begin{bsmallmatrix} *\\&*&*&*\\&*&*&*\\&&&*\end{bsmallmatrix}\cap\SSp(4,\C).
$$
This implies that $\varphi(W_F')$ is contained in the Levi component of the Klingen parabolic, a contradiction.
\end{proof}

An examination of Table \ref{dimensionstable} also establishes the following corollary. 

\begin{corollary}
\label{atNpiscor}
Let $(\pi,V)$ be an irreducible, admissible
representation of $\GSp(4,F)$ with trivial
central character. 
Assume that $\pi$ is paramodular.
Then $\dim V_s(N_{\pi,s})$ is $1$ or $2$.
Further assume that $N_\pi \geq 2$. 
If $\pi$ is generic,
then
\begin{equation}
\label{atNpiscoreq1}
\dim V_s(N_{\pi,s})
=
\begin{cases}
1&\text{if $N_{\pi,s}=N_\pi-1$},\\
2&\text{if $N_{\pi,s}=N_\pi$},
\end{cases}
\end{equation}
and if $\pi$ is non-generic, then
\begin{equation}
\label{atNpiscoreq2}
\dim V_s(N_{\pi,s})=1.
\end{equation}
\end{corollary}

The previous corollary suggests that the theory of stable Klingen vectors is  regular
if the minimal paramodular level is least two, and this observation is reinforced by the next
proposition. The following result provides practical criteria for identifying \catone 
paramodular representations under the assumption that the minimal paramodular level is at least two. 
We will use this proposition in the next chapter.

\begin{proposition}\label{shadowsigmavnewprop}
Let $(\pi,V)$ be an irreducible, admissible representation of the group $\GSp(4,F)$ with trivial central character. 
Assume that $\pi$ is paramodular, and that $N_\pi\geq2$. 
Let $v_{\mathrm{new}} \in V(N_\pi)$ be a newform, i.e., a non-zero element of $V(N_\pi)$.
Let $v_s \in V(N_\pi-1)$ be the shadow of $v_{\mathrm{new}}$ as defined in Lemma \ref{Wsgenericlemma1}.
Then the following are equivalent:
 \begin{enumerate}
  \item \label{shadowsigmavnewpropitem1} $\pi$ is a \catone paramodular representation, i.e., $N_{\pi,s}=N_\pi-1$.
  \item  \label{shadowsigmavnewpropitem2} The $T_{1,0}$-eigenvalue $\mu_\pi$ on $v_{\rm new}$ is non-zero.
  \item \label{shadowsigmavnewpropitem3} $\sigma_{N_\pi-1}v_{\rm new}\neq0$.
    \item \label{shadowsigmavnewpropitem4} $\sigma_{N_\pi}v_{\rm new}\neq0$.
  \item \label{shadowsigmavnewpropitem5} $T^s_{1,0} v_{\rm new}\neq0$.
  \item \label{shadowsigmavnewpropitem6} $v_s\neq0$.
 \end{enumerate}
\end{proposition}
\begin{proof}
\ref{shadowsigmavnewpropitem1} $\Leftrightarrow$ \ref{shadowsigmavnewpropitem2} follows by inspecting Table \ref{dimensionstable} from Theorem \ref{dimensionstheorem} and Table A.14 of \cite{NF}.

\ref{shadowsigmavnewpropitem3} $\Leftrightarrow$ \ref{shadowsigmavnewpropitem4}. By \ref{sigmaopslemmaitem6} of Lemma \ref{sigmaopslemma} we have
$\sigma_{N_\pi} \tau_{N_\pi} v_{\mathrm{new}} = \tau_{N_\pi-1} \sigma_{N_\pi-1} v_{\mathrm{new}}$. 
Since $v_{\mathrm{new}} \in V(N_\pi)$ we have $ \tau_{N_\pi} v_{\mathrm{new}} = v_{\mathrm{new}}$,
so that $\sigma_{N_\pi}  v_{\mathrm{new}} = \tau_{N_\pi-1} \sigma_{N_\pi-1} v_{\mathrm{new}}$. The 
injectivity of $\tau_{N_\pi-1}: V_s(N_\pi-1) \to V_s(N_\pi)$ from \ref{taunlevelraisinglemmaitem1} of Lemma \ref{taunlevelraisinglemma}
now implies that \ref{shadowsigmavnewpropitem3} and \ref{shadowsigmavnewpropitem4} are equivalent.

\ref{shadowsigmavnewpropitem3} $\Leftrightarrow$ \ref{shadowsigmavnewpropitem5} follows from Lemma \ref{heckeupdownlemma}.

\ref{shadowsigmavnewpropitem4} $\Leftrightarrow$ \ref{shadowsigmavnewpropitem6} follows from \eqref{shadowgeneraldefeq3}.

\ref{shadowsigmavnewpropitem6} $\Rightarrow$ \ref{shadowsigmavnewpropitem1} follows from Corollary \ref{upperboundpropcor}.

\ref{shadowsigmavnewpropitem2} $\Rightarrow$ \ref{shadowsigmavnewpropitem5}. By \eqref{T01Ts01lemmaeq2} we have $(1+q^{-1})p_{N_\pi}T^s_{1,0}v_{\mathrm{new}}=T_{1,0}v_{\mathrm{new}}$.
Hence, if $T_{1,0}v_{\mathrm{new}} \neq 0$, then $T^s_{1,0}v_{\mathrm{new}} \neq 0$. 
\end{proof}

%% file: SKMS_chapter6.tex
\chapter{Hecke Eigenvalues and Minimal Levels}
Let $\pi$ be a paramodular, irreducible, admissible representation
of $\GSp(4,F)$ with trivial central character such that $N_\pi\geq 2$.
In this chapter we will consider the action of the stable Klingen
operators $T_{0,1}^s$ and $T_{1,0}^s$\index{stable Hecke operators} on the vector spaces $V_s(N_{\pi,s})$
and $V_s(N_\pi)$. We will show that the actions of $T_{0,1}^s$ and
$T_{1,0}^s$ can be described in terms of the paramodular Hecke eigenvalues
$\lambda_\pi$ and $\mu_\pi$, and
in the final section we will explain how these results can be used 
to compute the paramodular Hecke eigenvalues $\lambda_\pi$
and $\mu_\pi$, and determine whether or not $\pi$ is non-generic,  using explicit upper block operators 
(see Sect.~\ref{opersksec} for the definition). As outlined in 
the introduction of this work, this has useful applications to 
Siegel paramodular newforms. Before beginning, we mention that 
throughout this chapter we will often use the equivalences of 
Proposition \ref{shadowsigmavnewprop}.

\section{At the minimal stable Klingen level}

Let $\pi$ be a paramodular, irreducible, admissible representation
of $\GSp(4,F)$ with trivial central character such that $N_\pi\geq 2$.
In the following theorem we describe the action of the stable Klingen
operators $T_{0,1}^s$ and $T_{1,0}^s$ on $V_s(N_{\pi,s})$ via the 
paramodular Hecke eigenvalues $\lambda_\pi$ and $\mu_\pi$. The proof
of this theorem uses the theory developed in the preceding chapters
as well as some results from \cite{NF}.

\begin{theorem}
\label{TSKleveltheorem}
Let $(\pi,V)$ be an irreducible, admissible representation of the group $\GSp(4,F)$ with trivial central character. Assume that $\pi$ is paramodular and that $N_\pi\geq2$. Let $v_{\mathrm{new}} \in V(N_\pi)$ be a newform, and let $v_s$
be the \index{shadow of a newform} shadow of $v_{\mathrm{new}}$. 
\begin{enumerate}
\item \label{TSKleveltheoremitem1} Assume that $\pi$ is a \catone representation, so that $N_{\pi,s}=N_\pi-1$ and $\mu_\pi \neq 0$.  Then $V_s(N_{\pi,s})=V_s(N_\pi-1)$ is one-dimensional and
\begin{equation}\label{T01generictheoremeq1}
T^s_{0,1}v_s=\lambda_\pi v_s\quad\text{and}\quad T^s_{1,0}v_s=(\mu_\pi+q^2)v_s.
\end{equation}
\item \label{TSKleveltheoremitem2} Assume that $\pi$ is a \cattwo representation, so that $N_{\pi,s}=N_\pi$ and $\mu_\pi=0$. Let $v_{\mathrm{new}} \in V(N_\pi)$ be a newform. The vector space $V_s(N_{\pi,s})$ is spanned by the vectors $v_{\mathrm{new}}$ and $T_{0,1}^s v_{\mathrm{new}}$. 
If $v_{\mathrm{new}}$ is not an eigenvector for $T_{0,1}^s$, so that $V_s(N_{\pi,s})$ is two-dimensional, then $\pi$
is generic, and the matrix of $T_{0,1}^s$ in the ordered basis $v_{\mathrm{new}}$, $T_{0,1}^s v_{\mathrm{new}}$ is
\begin{equation}
\label{TSKleveltheoremeq1}
\begin{bsmallmatrix}
0&-q^3\\
1&\lambda_\pi
\end{bsmallmatrix},
\end{equation}
so that 
\begin{equation}
\label{TSKleveltheoremeq11}
q^3 v_{\mathrm{new}} +(T_{0,1}^s)^2 v_{\mathrm{new}} = \lambda_\pi T_{0,1}^s v_{\mathrm{new}}.
\end{equation}
If $v_{\mathrm{new}}$ is an eigenvector for $T_{0,1}^s$, so that $V_s(N_{\pi,s})$ is one-dimensional, then $\pi$ is non-generic, and 
\begin{equation}
\label{TSKleveltheoremeq2}
T_{0,1}^s v_{\mathrm{new}} = (1+q^{-1})^{-1} \lambda_{\pi} v_{\mathrm{new}}.
\end{equation}
\end{enumerate}
\end{theorem}
\begin{proof}
\ref{TSKleveltheoremitem1}.  By Proposition \ref{shadowsigmavnewprop}
$v_s$ is non-zero, and by Corollary \ref{atNpiscor} $V(N_{\pi,s})$ is one-dimensional. 
Let $\mu_{\pi,s}, \lambda_{\pi,s} \in \C$ be such that $T^s_{0,1}v_s = \mu_{\pi,s} v_s$
and $T^s_{1,0}v_s = \lambda_{\pi,s} v_s$. We need to prove that $\lambda_{\pi,s}=\lambda_\pi$
and $\mu_{\pi,s}=\mu_\pi+q^2$. 
We have
\begin{align*}
T^s_{0,1}\rho_{N_\pi-1}'v_s
&=q\theta v_s-q\tau_{N_\pi-1} T^s_{0,1}v_s\qquad\text{(by \eqref{T01rhoeq})}\\
q^{-1}\mu_\pi T^s_{0,1}v_{\mathrm{new}}&=q\theta v_s-q\lambda_{\pi,s}\tau_{N_\pi-1} v_s
\qquad\text{(by \eqref{vsmulemmaeq2})}.
\end{align*}
Applying $p_{N_\pi}$ to this equation, we obtain:
\begin{align*}
q^{-1} \mu_\pi p_{N_\pi} T^s_{0,1}  v_{\mathrm{new}}&=qp_{N_\pi}\theta v_s-q\lambda_{\pi,s}p_{N_\pi}\tau_{N_\pi-1} v_s\\
(1+q)^{-1} \mu_\pi  T_{0,1}  v_{\mathrm{new}}&=qp_{N_\pi}\theta v_s-q\lambda_{\pi,s}p_{N_\pi}\tau_{N_\pi-1} v_s\qquad\text{(by \eqref{T01Ts01lemmaeq1})}\\
(1+q)^{-1} \mu_\pi \lambda_\pi v_{\mathrm{new}}&=qp_{N_\pi}\theta v_s+(1+q)^{-1}\lambda_{\pi,s}  \mu_\pi v_{\mathrm{new}}\qquad\text{(by \eqref{vsmulemmaeq1})}\\
(1+q)^{-1}\mu_\pi \lambda_\pi v_{\mathrm{new}}&=q\theta p_{N_\pi-1} v_s+(1+q)^{-1}\lambda_{\pi,s}  \mu_\pi v_{\mathrm{new}}\qquad\text{(by \eqref{commreleq5})}\\
(1+q)^{-1} \mu_\pi \lambda_\pi v_{\mathrm{new}} &=(1+q)^{-1}\lambda_{\pi,s}  \mu_\pi v_{\mathrm{new}}\qquad\text{(since $p_{N_\pi-1} v_s=0$)}.
\end{align*}
Since $\mu_\pi \neq 0$  we obtain $\lambda_{\pi,s}=\lambda$, as desired.

Next, 
\begin{align*}
T_{1,0}^s \rho'_{N_\pi-1}v_s
&=q^3\tau_{N_\pi-1}v_s -q\tau_{N_\pi-1}T_{1,0}^sv_s\qquad \text{(by \eqref{T10rhoeq})}\\
q^{-1}\mu_\pi T_{1,0}^s v_{\mathrm{new}}
&=q^3\tau_{N_\pi-1}v_s -q\mu_{\pi,s}\tau_{N_\pi-1} v_s\qquad\text{(by \eqref{vsmulemmaeq2})}.
\end{align*}
Applying $p_{N_\pi}$ to this equation, we have:
\begin{align*}
q^{-1}\mu_\pi p_{N_\pi} T_{1,0}^s v_{\mathrm{new}}
&=(q^3-q\mu_{\pi,s})p_{N_\pi}\tau_{N_\pi-1} v_s\\
(1+q)^{-1}\mu_\pi T_{1,0} v_{\mathrm{new}}
&=(q^3-q\mu_{\pi,s})p_{N_\pi}\tau_{N_\pi-1}v_s\qquad\text{(by \eqref{T01Ts01lemmaeq2})}\\
(1+q)^{-1}\mu_\pi^2 v_{\mathrm{new}}
&=-(q+q^2)^{-1}(q^3-q\mu_{\pi,s})\mu_\pi v_{\mathrm{new}}\qquad\text{(by \eqref{vsmulemmaeq1})}.
\end{align*}
Since $\mu_\pi \neq 0$ by Proposition \ref{shadowsigmavnewprop} we obtain $\mu_{\pi,s}=\mu_\pi+q^2$. This
completes the proof of \ref{TSKleveltheoremitem1}. 

\ref{TSKleveltheoremitem2}. By Corollary \ref{atNpiscor}, since $N_{\pi,s}=N_\pi$,  the vector space $V_s(N_{\pi,s})$
is two-dimensional if $\pi$ is generic, and $V_s(N_{\pi,s})$ is one-dimensional if $\pi$ is
non-generic. Thus, to prove that $v_{\mathrm{new}}$ and $T_{0,1}^s v_{\mathrm{new}}$
span $V_s(N_{\pi,s})$ it will suffice to prove that $v_{\mathrm{new}}$ and $T_{0,1}^s v_{\mathrm{new}}$
are linearly independent if $\pi$ is generic. Assume that $\pi$ is generic, and that $v_{\mathrm{new}}$
and $T_{0,1}^sv_{\mathrm{new}}$ are linearly dependent; we will obtain a contradiction. Let 
$c \in \C$ be such that $T_{0,1}^s v_{\mathrm{new}}= c v_{\mathrm{new}}$. We have:
\begin{align*}
T^s_{0,1}\rho_{N_\pi}'v_{\mathrm{new}}&=q\theta v_{\mathrm{new}}+\tau_{N_\pi} T^s_{0,1}v_{\mathrm{new}} \qquad \text{(by \eqref{T01rhoeqb})}\\
q^{-1} T^s_{0,1}\theta' v_{\mathrm{new}}&=q\theta v_{\mathrm{new}}+c  v_{\mathrm{new}} \qquad \text{(by \eqref{rhopoplemmaeq4} and $v_{\mathrm{new}} \in V(N_\pi)$)}.
\end{align*}
Applying $p_{N_\pi+1}$ to this equation, we obtain:
\begin{align*}
q^{-1}p_{N_\pi+1} T^s_{0,1}\theta' v_{\mathrm{new}}&=qp_{N_\pi+1}\theta v_{\mathrm{new}}+c p_{N_\pi+1} v_{\mathrm{new}} \\
  T_{0,1}\theta' v_{\mathrm{new}}&=(q+q^2)\theta v_{\mathrm{new}}+c\theta' v_{\mathrm{new}}
\qquad\text{(by \eqref{T01Ts01lemmaeq1}, \eqref{pnVsneq2}, and \eqref{thetap1defeq})}\\
 \theta' T_{0,1} v_{\mathrm{new}}+ q^2\theta v_{\mathrm{new}}&=(q+q^2)\theta v_{\mathrm{new}}+c\theta' v_{\mathrm{new}}
\qquad\text{(by (6.17) of \cite{NF})}\\
 \lambda_\pi\theta'  v_{\mathrm{new}}+ q^2\theta v_{\mathrm{new}}&=(q+q^2)\theta v_{\mathrm{new}}+c\theta' v_{\mathrm{new}}\\
(\lambda_\pi-c)\theta'  v_{\mathrm{new}}&=q\theta v_{\mathrm{new}}.
\end{align*}
Since $\pi$ is generic, the vectors $\theta'v_{\mathrm{new}}$ and $\theta v_{\mathrm{new}}$ are linearly independent
by Theorem~7.5.6 of \cite{NF}. It follows that $q=0$, a contradiction. This completes the proof that 
$v_{\mathrm{new}}$ and $T_{0,1}^s v_{\mathrm{new}}$ span $V_s(N_{\pi,s})$. 

Next, assume that $v_{\mathrm{new}}$ is not an eigenvector for $T_{0,1}^s$, so that $v_{\mathrm{new}}$ and $T_{0,1}^s v_{\mathrm{new}}$ are linearly independent. Then by Corollary \ref{atNpiscor} the vector space $V_s(N_{\pi,s})$ is two-dimensional and $\pi$ is generic. 
We need to prove that the matrix of $T_{0,1}^s$ in the ordered basis $v_{\mathrm{new}}$, $T_{0,1}^s v_{\mathrm{new}}$
is as in \eqref{TSKleveltheoremeq1}. It is obvious that the first column of the matrix of $T_{0,1}^s$ is as in \eqref{TSKleveltheoremeq1}. To determine the second column of $T_{0,1}^s$, let $a,b \in \C$ be such
that 
\begin{equation}
\label{TSKleveltheoremeq3}
(T_{0,1}^s)^2 v_{\mathrm{new}} = a v_{\mathrm{new}} + b T_{0,1}^s v_{\mathrm{new}}.
\end{equation}
We need to prove that $a=-q^3$ and $b=\lambda_\pi$.
To do this, we will assume that $V=\mathcal{W}(\pi,\psi_{c_1,c_2})$, 
the \index{Whittaker model} Whittaker model of $\pi$ with respect to $\psi_{c_1,c_2}$; as usual, we assume that $c_1,c_2 \in \OF^\times$. 
As in (7.15) of \cite{NF}, we define $c_{i,j} = v_{\mathrm{new}}(\Delta_{i,j})$ for $i,j \in \Z$; here
$\Delta_{i,j}$ is defined as in \eqref{deltaijdefeq}.  
Let $j \in \Z$ be such that $j \geq 0$. 
Evaluating the left-hand side of \eqref{TSKleveltheoremeq3} at $\Delta_{0,j}$, we have by \eqref{deltaWlemmaeq1},
\begin{align}
((T_{0,1}^s)^2 v_{\mathrm{new}})(\Delta_{0,j})&=q^4 v_{\mathrm{new}}(\Delta_{2,j-2}) +2 q^5  v_{\mathrm{new}} (\Delta_{1,j})
+q^6 W(\Delta_{0,j+2}) \nonumber \\
&=q^4c_{2,j-2}+2q^5 c_{1,j}+q^6 c_{0,j+2}.\label{TSKleveltheoremeq5}
\end{align}
By Lemma 7.4.4 of \cite{NF}, 
$$
c_{2,j-2}=q^{-4}(\mu_\pi+q^2) c_{1,j-2}=q^{-4}(\mu_\pi+q^2)q^{-4}\mu_\pi c_{0,j-2}=0
$$ 
 and also $c_{1,j}=q^{-4}\mu_\pi c_{0,j}=0$ because $\mu_\pi=0$. Hence
 \begin{equation}\label{noshadowT01matrixlemmaeq4}
 ((T^s_{0,1})^2 v_{\rm new}))(\Delta_{0,j})=q^6c_{0,j+2}.
\end{equation}
Evaluating the right-hand side of \eqref{TSKleveltheoremeq3} at $\Delta_{0,j}$, we have by \eqref{deltaWlemmaeq1},
\begin{align}
(a v_{\mathrm{new}} + b T_{0,1}^s v_{\mathrm{new}})(\Delta_{0,j})
&=a v_{\mathrm{new}} (\Delta_{0,j}) +b q^2 v_{\mathrm{new}}(\Delta_{1,j-1}) 
+ bq^3 v_{\mathrm{new}}(\Delta_{0,j+1}) \nonumber \\
&= ac_{0,j}+b q^2 c_{1,j-1} + bq^3 c_{0,j+1}\nonumber\\
&= ac_{0,j}+ bq^3 c_{0,j+1}, \label{TSKleveltheoremeq6}
\end{align}
because again by Lemma 7.4.4 of \cite{NF}, $c_{1,j-1}=q^{-4}\mu_\pi c_{0,j-1}=0$ since $\mu_\pi=0$.
Since \eqref{noshadowT01matrixlemmaeq4} and \eqref{TSKleveltheoremeq6} are equal, 
$$
 ac_{0,j}+bq^3c_{0,j+1}=q^6c_{0,j+2}\qquad\text{for all }j\geq0.
$$
By Lemma 7.4.4 of \cite{NF}, $c_{0,j+2}=q^{-3}\lambda_\pi c_{0,j+1}-q^{-3}c_{0,j}$. Substituting, we find
$$
 ac_{0,j}+bq^3c_{0,j+1}=q^3\lambda_\pi c_{0,j+1}-q^3c_{0,j}\qquad\text{for all }j\geq0,
$$
or
$$
 (a+q^3)c_{0,j}=q^3(\lambda_\pi-b)c_{0,j+1}\qquad\text{for all }j\geq0.
$$
Now assume that $b \neq \lambda_\pi$. We then have 
$$
c_{0,j} = M^j c_{0,0}\qquad\text{for all $j \geq 0$},
$$
where 
$$
M=\dfrac{a+q^3}{q^3(\lambda_\pi-b)}.
$$
A standard calculation using Lemma 4.1.1 and Lemma 4.1.2 of \cite{NF}
shows that
\begin{align}
Z(s,v_\mathrm{new}) &= (1-q^{-1}) \sum_{j=0}^\infty c_{0,j} ( q^{\frac{3}{2}}q^{-s})^j. \label{TSKleveltheoremeq4}
\end{align}
By \eqref{TSKleveltheoremeq4} we thus obtain
\begin{equation}
\label{TSKleveltheoremeq100}
Z(s,v_{\mathrm{new}}) = \dfrac{(1-q^{-1})c_{0,0}}{1-M q^{\frac{3}{2}}q^{-s}}.
\end{equation}
On the other hand, by Proposition 7.4.5 of \cite{NF}, since $\mu_\pi=0$, we have
\begin{align}
Z(s,v_\mathrm{new})&=\dfrac{(1-q^{-1})c_{0,0}}{1-q^{-\frac{3}{2}}\lambda_\pi q^{-s} +q^{-2s}}. \label{TSKleveltheoremeq41}
\end{align}
This contradicts \eqref{TSKleveltheoremeq100}.
Hence $b=\lambda_\pi$, which then implies $a=-q^3$, as desired. 

Finally, assume that $v_{\mathrm{new}}$ is an eigenvector for $T_{0,1}^s$, so that
$V_s(N_{\pi,s})$ is one-dimensional and spanned by $v_{\mathrm{new}}$. Then, as we have
already indicated above, $\pi$ is non-generic by Corollary \ref{atNpiscor}. Let $\lambda_{\pi,s} \in \C$
be such that $T_{0,1}^s v_{\mathrm{new}} = \lambda_{\pi,s} v_{\mathrm{new}}$. 
By~\eqref{T01Ts01lemmaeq1},
\begin{align*}
(1+q^{-1})\lambda_{\pi,s}v_{\rm new}
&=(1+q^{-1})p_{N_\pi}(\lambda_{\pi,s}v_{\rm new})\\
&=(1+q^{-1})p_{N_\pi}(T^s_{0,1}v_{\rm new})\\
&=T_{0,1}v_{\rm new}\\
&=\lambda_\pi v_{\rm new}.
\end{align*}
Hence, $\lambda_{\pi,s}=(1+q^{-1})^{-1}\lambda_\pi$.
This completes the proof. 
\end{proof}

As the next lemma shows, in case \ref{TSKleveltheoremitem1} of Theorem \ref{TSKleveltheorem} the eigenvalue $\mu_\pi$
has a simple form, so that the action of $T_{1,0}^s$ on the shadow vector is also explicit.

\begin{corollary}
\label{muvalueslemma}
Let $\pi$ be 
an  irreducible, admissible representation of the group $\GSp(4,F)$ with trivial central character.
Assume that $\pi$ is  paramodular, $N_\pi \geq 2$, and $\pi$ is a \index{category 1 representation}
\catone representation, so that $N_{\pi,s}=N_\pi-1$ and $\mu_\pi \neq 0$. Let $v_{\mathrm{new}} \in V(N_\pi)$ be a newform, 
and let $v_s \in V_s(N_\pi-1)$
be the \index{shadow of a newform} shadow of $v_{\mathrm{new}}$.
Then
\begin{equation}
\label{muvalueseq}
\mu_\pi =
\begin{cases}
-q^2 +\varepsilon_\pi q & \text{if $N_\pi=2$,}\\
-q^2 &\text{if $N_\pi>2$},
\end{cases}
\end{equation}
so that
\begin{equation}
\label{muvalueseq2}
T_{1,0}^s v_s = 
\begin{cases}
\varepsilon_\pi q v_s &\text{if $N_\pi = 2$,}\\
0 &\text{if $N_\pi >2$}.
\end{cases}
\end{equation}
\end{corollary}
\begin{proof}
The assertion \eqref{muvalueseq} follows from an 
inspection of Table \ref{levelsandeigenvaluestable},
and \eqref{muvalueseq2} now follows from \ref{TSKleveltheoremitem1} of Theorem \ref{TSKleveltheorem}. 
\end{proof}

We mention that Corollary \ref{muvalueslemma} is  
consistent with \ref{Ts10eigenvalueslemmaitem3} of Lemma \ref{Ts10eigenvalueslemma}.
To explain this, let the notation be as in Corollary \ref{muvalueslemma}, and assume that $N_\pi > 2$. Since $N_{\pi,s}=N_\pi-1>1$, $V_s(1)=0$, and so by 
\ref{Ts10eigenvalueslemmaitem3} of Lemma \ref{Ts10eigenvalueslemma} the only eigenvalue of $T_{1,0}^s$
on $V_s(N_\pi-1)$ is $0$. This is indeed the case by \eqref{muvalueseq2}.

\section{Non-generic paramodular representations}
\label{nongenparasec}
Let $(\pi,V)$ be a paramodular, irreducible, admissible representation
of the group $\GSp(4,F)$ with trivial central character such that $N_\pi \geq 2$.
In this section we determine under what conditions $\pi$ is non-generic. This information
will be used in the next two sections. 

\begin{lemma}
\label{IVbcharlemma}
Let $(\pi,V)$ be an irreducible, admissible representation
of the group $\GSp(4,F)$ with trivial central character.
Assume that $\pi$ is paramodular and $N_\pi \geq 2$. Then
$\pi$ is a IVb representation if and only if
$N_\pi=2$, $\mu_\pi = -q^2+q$, and $\lambda_\pi=\pm(1+q^2)$.
\end{lemma}
\begin{proof}
If $\pi$ is a IVb representation, then $N_\pi=2$, $\mu_\pi=-q^2+q$,
and $\lambda_\pi=\pm (1+q^2)$ by 
Table \ref{levelsandeigenvaluestable} on  p.~\pageref{levelsandeigenvaluestable}.
Conversely, assume that $N_\pi=2$, $\mu_\pi=-q^2+q$,
and $\lambda_\pi=\pm (1+q^2)$.
Table \ref{levelsandeigenvaluestable}
along with $N_\pi=2$ and $\mu_\pi=-q^2+q$ then imply that $\pi$ is either a IIIa representation of the
form $\chi \rtimes \sigma \St_{\GSp(2)}$ for unramified characters $\chi$ and $\sigma$ of $F^\times$
such that $\chi \sigma^2 =1$ and $\chi \notin \{1,\nu^{\pm 2}\}$, or $\pi$ is VIa representation of the form $\tau(S,\nu^{-1/2}\sigma)$ for
an unramified character $\sigma$ of $F^\times$ such that $\sigma^2=1$, or $\pi$ is a IVb representation of the form $\pi=L(\nu^2,\nu^{-1}\sigma \St_{\GSp(2)})$ for an unramified character $\sigma$ of $F^\times$ such that $\sigma^2=1$.
Assume that the first possibility holds. Then (recalling the assumption $\lambda_\pi=\pm(1+q^2)$) we have $q(\sigma(\varpi)+\sigma(\varpi)^{-1})=\lambda_\pi=\pm(1+q^2)$, so that $\sigma(\varpi)=\pm q^{\pm 1}$;
this implies that $\chi = \nu^{\pm 2}$, a contradiction. The second possibility is similarly seen to lead
to a contradiction. It follows that $\pi$ is a IVb representation.
\end{proof}

\begin{lemma}
\label{nongenericfindlemma}
Let $(\pi,V)$ be an irreducible, admissible representation of
the group $\GSp(4,F)$ with trivial central character.
Assume that $\pi$ is paramodular and $N_\pi \geq 2$.
Let $v_{\mathrm{new}} \in V(N_\pi)$ be a newform.
Then $\pi$ is non-generic if and only if 
\begin{equation}
\label{nongenericfindcoreq1}
\text{$\mu_\pi =0$ and $v_{\mathrm{new}} \in V(N_\pi)=V(N_{\pi,s})$ is an eigenvector for $T_{0,1}^s$}
\end{equation}
or
\begin{equation}
\label{nongenericfindcoreq2}
\text{$N_\pi=2$, $\mu_\pi=-q^2+q$, and $\lambda_\pi=\pm(1+q^2)$.}
\end{equation}
\end{lemma}
\begin{proof}
Assume that $\pi$ is non-generic. Then by inspection of Table \ref{levelsandeigenvaluestable} on 
p.~\pageref{levelsandeigenvaluestable}, either $\mu_\pi=0$, so that \eqref{nongenericfindcoreq1} holds
by \ref{TSKleveltheoremitem2} of Theorem \ref{TSKleveltheorem}, or $\pi$ is a IVb representation of the 
form $\pi=L(\nu^2,\nu^{-1}\sigma \St_{\GSp(2)})$ for an unramified character $\sigma$ of $F^\times$ such that $\sigma^2=1$, $N_\pi=2$,
$\lambda_\pi = \sigma(\varpi)(1+q^2)$, and $\mu_\pi=-q^2+q$, so that \eqref{nongenericfindcoreq2} holds.
In the converse direction, if \eqref{nongenericfindcoreq1} holds, then $\pi$ is non-generic by \ref{TSKleveltheoremitem2} of
Theorem~\ref{TSKleveltheorem}, and if \eqref{nongenericfindcoreq2} holds, then $\pi$ is a IVb representation by 
Lemma \ref{IVbcharlemma}, and is hence non-generic. 
\end{proof}

\begin{lemma}
\label{nogenorIVblemma}
Let $(\pi,V)$ be an irreducible, admissible representation of the group $\GSp(4,F)$
with trivial central character. Assume that $\pi$ is paramodular, $N_\pi \geq 2$,
and $\pi$ is a \catone representation. Then the following are equivalent:
\begin{enumerate}
\item \label{nogenorIVblemmaitem1} $\pi$ is non-generic;
\item \label{nogenorIVblemmaitem2} $\pi$ is a IVb representation;
\item \label{nogenorIVblemmaitem3} $N_\pi=2$, $\mu_\pi=-q^2+q$, and $\lambda_\pi=\pm(1+q^2)$.
\end{enumerate}
\end{lemma}
\begin{proof}
\ref{nogenorIVblemmaitem1} $\Leftrightarrow$ \ref{nogenorIVblemmaitem3}. Assume \ref{nogenorIVblemmaitem1}. By Lemma \ref{nongenericfindlemma} either 
\eqref{nongenericfindcoreq1} or \eqref{nongenericfindcoreq2} holds. Since $\mu_\pi \neq 0$
because $\pi$ is a \catone representation (see Proposition \ref{shadowsigmavnewprop}),
\ref{nogenorIVblemmaitem3} holds. Assume \ref{nogenorIVblemmaitem3}. Then $\pi$ is non-generic by Lemma \ref{nongenericfindlemma}.

\ref{nogenorIVblemmaitem2} $\Leftrightarrow$ \ref{nogenorIVblemmaitem3}. This follows from Lemma \ref{IVbcharlemma}.
\end{proof}

\begin{lemma}
\label{cat2nongencharlemma}
Let $(\pi,V)$ be an irreducible, admissible representation of the group $\GSp(4,F)$
with trivial central character. Assume that $\pi$ is paramodular, $N_\pi \geq 2$,
and $\pi$ is a \cattwo representation. Then $\pi$ is non-generic if and only if
$\pi$ is a Saito-Kurokawa representation.
\end{lemma}
\begin{proof}
This follows from an inspection of Table \ref{levelsandeigenvaluestable}.
\end{proof}

\section{At the minimal paramodular level} 

Let $\pi$ be a paramodular, irreducible, admissible
representation of $\GSp(4,F)$ with trivial central 
character such that $N_\pi \geq 2$. In this section
we describe the action of the \index{stable Hecke operators} stable Klingen Hecke operators
$T_{0,1}^s$ and $T_{1,0}^s$ on $V_s(N_\pi)$ in terms
of the paramodular Hecke eigenvalues $\lambda_\pi$ and $\mu_\pi$.
Since the case when
$\pi$ is a \cattwo representation, so that $N_\pi=N_{\pi,s}$, was dealt with in \ref{TSKleveltheoremitem2} of 
Theorem \ref{TSKleveltheorem}, we need only consider \index{category 1 representation} \catone representations, i.e.,
$\pi$ such that $N_\pi=N_{\pi,s}+1$.  By Lemma~\ref{nogenorIVblemma}, there are two cases
to consider: $\pi$ is generic, or $\pi$ is a IVb representation. We first consider
the case when $\pi$ is generic.

\begin{theorem}\label{Wsgenericytheorem2}
Let $(\pi,V)$ be a generic, irreducible, admissible representation of $\GSp(4,F)$ with trivial central character. Assume that $N_\pi\geq2$ and that $\pi$ is a \catone representation, so that $N_{\pi,s}=N_\pi-1$. 
Let $v_{\rm new} \in V(N_\pi)$ be a newform, and let $v_s$ be the shadow of $v_{\mathrm{new}}$. The vector
space $V_s(N_\pi)$ is three-dimensional and has ordered basis
\begin{equation}\label{Wsgenericlemma2eq1}
v_{\rm new},\quad \tau_{N_\pi-1}v_s,\quad \theta v_s.
\end{equation}
The matrix of the endomorphism $T_{0,1}^s$ of $V_s(N_\pi)$ with respect to \eqref{Wsgenericlemma2eq1}
is 
\begin{equation}\label{Wsgenericlemma2eq6}
\begin{bsmallmatrix}0&0&-\mu_\pi q^{-1}(\mu_\pi+q^2)\\-q^2\mu_\pi^{-1}\lambda_\pi&\lambda_\pi&-q\mu_\pi\\q^2\mu_\pi^{-1}&0&\lambda_\pi\end{bsmallmatrix},
\end{equation}
and the matrix of the endomorphism $T_{1,0}^s$ of $V_s(N_\pi)$ with respect to \eqref{Wsgenericlemma2eq1} is
\begin{equation}\label{Wsgenericlemma2eq2}
\begin{bsmallmatrix}0&0&0&\\-q^2&\mu_\pi+q^2&q^2\lambda_\pi\\0&0&0\end{bsmallmatrix}.
\end{equation}
The characteristic polynomial of $T_{0,1}^s$ on $V_s(N_\pi)$ is
\begin{equation}\label{Wsgenericlemma2eq11}
p(T_{0,1}^s, V_s(N_\pi),X)=(X-\lambda_\pi)(X^2-\lambda_\pi X+q(\mu_\pi+q^2)), 
\end{equation}
and the  characteristic polynomial of $T_{1,0}^s$ on $V_s(N_\pi)$ is 
\begin{equation}
p(T_{1,0}^s, V_s(N_\pi),X)=(X -(\mu_\pi +q^2))X^2.
\end{equation}
\end{theorem}
\begin{proof}
It follows from Theorem \ref{genericdimensionstheorem} that $\dim V_s(N_\pi)=3$, and that the vectors in \eqref{Wsgenericlemma2eq1} are linearly independent.

To prove that the matrix of $T_{0,1}^s$ of $V_s(N_\pi)$ with respect to \eqref{Wsgenericlemma2eq1}
is as in \eqref{Wsgenericlemma2eq6} we proceed as follows. 
We have
\begin{align}
 T^s_{0,1}v_{\rm new}&=q\mu_\pi^{-1}T^s_{0,1}\rho'_{N_\pi-1}v_s\qquad\text{(by \eqref{vsmulemmaeq2})}\nonumber\\
 &=q\mu_\pi^{-1}(q\theta v_s-q\tau_{N_\pi-1}T^s_{0,1}v_s)\qquad \text{(by \eqref{T01rhoeq})}\nonumber\\
 &=q^2\mu_\pi^{-1}\theta v_s-q^2\mu_\pi^{-1}\tau_{N_\pi-1}\lambda_\pi v_s \qquad\text{(by \eqref{T01generictheoremeq1})}.\label{Wsgenericlemma2eq3}
\end{align}
This verifies the first column of \eqref{Wsgenericlemma2eq6}.
For the second column we have:
\begin{align*}
T_{0,1}^s \tau_{N_\pi-1}v_s
&=\tau_{N_\pi-1}T_{0,1}^s v_s\qquad \text{(by \eqref{T01taueq})}\\
&=\lambda_\pi \tau_{N_\pi-1} v_s\qquad \text{(by \eqref{T01generictheoremeq1})}.
\end{align*}
This verifies the second column of \eqref{Wsgenericlemma2eq6}.
To obtain the third column we will first prove that
\begin{equation}\label{Wsgenericlemma2eq4}
(T^s_{0,1})^2 v_{\rm new}
=-q(\mu_\pi+q^2)v_{\rm new}-(q^2\lambda_\pi^2\mu_\pi^{-1}+q^3)\tau_{N_\pi-1}v_s+q^2\lambda_\pi\mu_\pi^{-1}\theta v_s.
\end{equation}
To prove \eqref{Wsgenericlemma2eq4}, we begin with \eqref{Wsgenericlemma2eq3} and find that:
\begin{align}
(T^s_{0,1})^2 v_{\rm new} &=  T^s_{0,1}(T^s_{0,1}v_{\rm new})\nonumber \\
&=T^s_{0,1}\big(-q^2\mu_\pi^{-1}\lambda_\pi\tau_{N_\pi-1}v_s+q^2\mu_\pi^{-1}\theta v_s\big)\nonumber \\
 &=-q^2\mu_\pi^{-1}\lambda_\pi T^s_{0,1}\tau_{N_\pi-1}v_s+q^2\mu_\pi^{-1}T^s_{0,1}\theta v_s\nonumber \\
 &=-q^2\mu_\pi^{-1}\lambda_\pi \tau_{N_\pi-1}T^s_{0,1}W_s+q^2\mu_\pi^{-1}T^s_{0,1}\theta v_s
\qquad \text{(by \eqref{T01taueq})}\nonumber \\
 &=-q^2\mu_\pi^{-1}\lambda_\pi^2 \tau_{N_\pi-1}v_s+q^2\mu_\pi^{-1}T^s_{0,1}\theta v_s
\qquad \text{(by \eqref{T01generictheoremeq1})}. \label{Ts01squareintereq1}
\end{align}
To proceed further we will need to look at cases.
Assume first that $N_\pi\geq3$. Then  \eqref{Ts01squareintereq1} becomes, by \eqref{Ts01thetalemmaeq1} with $n=N_\pi-1$, and since $\sigma_{N_{\pi-2}}v_s=0$ as $V_s(N_{\pi}-2)=0$,
\begin{align*}
&(T^s_{0,1})^2 v_{\rm new}
 =-q^2\mu_\pi^{-1}\lambda_\pi^2 \tau_{N_\pi-1}v_s+q^2\mu_\pi^{-1}\big(\theta T^s_{0,1}v_s+q^3\tau_{N_\pi-1}v_s-q^3\eta\sigma_{N_\pi-2}v_s\big)\\
 &\qquad=-q^2\mu_\pi^{-1}\lambda_\pi^2 \tau_{N_\pi-1}v_s+q^2\mu_\pi^{-1}\big(\theta T^s_{0,1}v_s+q^3\tau_{N_\pi-1}v_s\big)\\
 &\qquad=-q^2\mu_\pi^{-1}\lambda_\pi^2 \tau_{N_\pi-1}v_s+q^2\mu_\pi^{-1}\big(\lambda_\pi\theta v_s+q^3\tau_{N_\pi-1}v_s\big)
\qquad \text{(by \eqref{T01generictheoremeq1})}\\
 &\qquad=(-q^2\mu_\pi^{-1}\lambda_\pi^2+q^5\mu_\pi^{-1}) \tau_{N_\pi-1}v_s+q^2\mu_\pi^{-1}\lambda_\pi\theta v_s.
\end{align*}
Since $N_\pi \geq 3$ and $\pi$ is generic, we have  $\mu_\pi=-q^2$ by Corollary~\ref{muvalueslemma}. Hence
$$
(T^s_{0,1})^2 v_{\rm new}=(\lambda_\pi^2-q^3) \tau_{N_\pi-1}v_s-\lambda_\pi\theta v_s.
$$
This is \eqref{Wsgenericlemma2eq4}.
Now assume that $N_\pi=2$. Then continuing from \eqref{Ts01squareintereq1}, 
\begin{align*}
&(T^s_{0,1})^2 v_{\rm new}=-q^2\mu_\pi^{-1}\lambda_\pi^2 \tau_1v_s+q^2\mu_\pi^{-1}\big(\theta T^s_{0,1}v_s+q^3\tau_1 v_s-q^3e(v_s)\big)\qquad\text{(by \eqref{Ts01thetalemmaeq1b})}\\
&\qquad=-q^2\mu_\pi^{-1}\lambda_\pi^2 \tau_1 v_s+q^2\mu_\pi^{-1}\big(\theta T^s_{0,1}v_s+q^3\tau_1v_s\\
&\quad\qquad-q^3\big(q^{-4}\eta T^s_{1,0}v_s-q^{-2}\tau_2T^s_{1,0}v_s+q^{-2}\tau_1T^s_{1,0}v_s\big)\big)
\qquad\text{(by \eqref{eTs10identityeq})}\\
&\qquad=-q^2\mu_\pi^{-1}\lambda_\pi^2 \tau_1v_s+q^2\mu_\pi^{-1}\big(\lambda_\pi\theta v_s+q^3\tau_1v_s\\
&\quad\qquad-q^3\big(q^{-4}( \mu_\pi+q^2)\eta v_s-q^{-2}(\mu_\pi+q^2)\tau_2v_s+q^{-2}(\mu_\pi+q^2)\tau_1v_s\big)\big)\\
&\quad\qquad \text{(by \eqref{T01generictheoremeq1})}\\
&\qquad=-q^2\mu_\pi^{-1}\lambda_\pi^2 \tau_1v_s+q^2\mu_\pi^{-1}\big(\lambda_\pi\theta v_s+q^3\tau_1v_s\\
&\quad\qquad-q^3(\mu_\pi+q^2)\big(q^{-3}(q^{-1}\eta v_s-q\tau_2v_s)+q^{-2}\tau_1v_s\big)\big)\\
&\qquad=-q^2\mu_\pi^{-1}\lambda_\pi^2 \tau_1v_s+q^2\mu_\pi^{-1}\big(\lambda_\pi\theta v_s+q^3\tau_1v_s\\
&\quad\qquad-q^3(\mu_\pi+q^2)\big(q^{-3}\rho_1'v_s+q^{-2}\tau_1v_s\big)\big)\qquad \text{(by \eqref{rhopoplemmaeq3})}\\
&\qquad=-q^2\mu_\pi^{-1}\lambda_\pi^2 \tau_1v_s+q^2\mu_\pi^{-1}\big(\lambda_\pi\theta v_s+q^3\tau_1v_s\\
&\quad\qquad-q^3(\mu_\pi+q^2)\big(q^{-4}\mu_\pi v_{\rm new}+q^{-2}\tau_1v_s\big)\big)\qquad\text{(by \eqref{vsmulemmaeq2})}\\
&\qquad=-q(\mu_\pi+q^2)v_{\rm new}-q^2(\mu_\pi^{-1}\lambda_\pi^2+q)\tau_1v_s+q^2\mu_\pi^{-1}\lambda_\pi\theta v_s.
\end{align*}
This is \eqref{Wsgenericlemma2eq4}, and concludes the proof of \eqref{Wsgenericlemma2eq4} in all cases. 
We now return to the verification of the third column of \eqref{Wsgenericlemma2eq6}. 
Equating the right-hand sides of \eqref{Ts01squareintereq1} and \eqref{Wsgenericlemma2eq4} gives
\begin{align*}
 &-q^2\mu_\pi^{-1}\lambda_\pi^2 \tau_{N_\pi-1}v_s+q^2\mu_\pi^{-1}T^s_{0,1}(\theta v_s)\\
 &\qquad=-(q\mu_\pi+q^3)v_{\rm new}-(q^2\lambda_\pi^2\mu_\pi^{-1}+q^3)\tau_{N_\pi-1}v_s+q^2\lambda_\pi\mu_\pi^{-1}\theta v_s.
\end{align*}
Hence
\begin{align*}
 q^2\mu_\pi^{-1}T^s_{0,1}(\theta v_s)&=q^2\mu_\pi^{-1}\lambda_\pi^2 \tau_{N_\pi-1}v_s\\
   &\qquad-(q\mu_\pi+q^3)v_{\rm new}-(q^2\lambda_\pi^2\mu_\pi^{-1}+q^3)\tau_{N_\pi-1}v_s+q^2\lambda_\pi\mu_\pi^{-1}\theta v_s\\
 &=-(q\mu_\pi+q^3)v_{\rm new}-q^3\tau_{N_\pi-1}v_s+q^2\lambda_\pi\mu_\pi^{-1}\theta v_s.
\end{align*}
It follows that
$$
 T^s_{0,1}(\theta v_s)=-\mu_\pi(q^{-1}\mu_\pi+q)v_{\rm new}-q\mu_\pi\tau_{N_\pi-1}v_s+\lambda_\pi\theta v_s.
$$
This verifies the third column of \eqref{Wsgenericlemma2eq6}. 

To prove that the matrix of  $T_{1,0}^s$ of $V_s(N_\pi)$ with respect to \eqref{Wsgenericlemma2eq1} is
as in \eqref{Wsgenericlemma2eq2}, we note that
by \eqref{Wsgenericlemma100}
$$
 T^s_{1,0}v_{\rm new}=-q^{2}\tau_{N_\pi-1}v_s.
$$
This verifies the first column of the matrix in \eqref{Wsgenericlemma2eq2}. The second column in \eqref{Wsgenericlemma2eq2} follows from \eqref{T10taueq} and \eqref{T01generictheoremeq1}. The third column in \eqref{Wsgenericlemma2eq2} follows from:
\begin{align*}
 T^s_{1,0}\theta v_s&=q^2T^s_{0,1}\tau_{N_\pi-1}v_s\qquad\text{(by \eqref{Ts10thetaeq})}\\
 &=q^2\tau_{N_\pi-1}T^s_{0,1}v_s\qquad \text{(by \eqref{T01taueq})}\\
 &\stackrel{}{=}q^2\tau_{N_\pi-1}\lambda_\pi v_s \qquad\text{(by \eqref{T01generictheoremeq1})}.
\end{align*}
The final assertions about characteristic polynomials follow by standard calculations.
\end{proof}
Let $\pi$ be a paramodular, irreducible, admissible representation of $\GSp(4,F)$
with trivial central character such that $N_\pi \geq 2$. 
The next proposition considers 
the actions
of $T_{0,1}^s$ and $T_{1,0}^s$ on $V(N_\pi)$ when $\pi$ is a IVb representation. 
This proposition completes our analysis of the actions of $T_{0,1}^s$ and
$T_{1,0}^s$ on $V(N_\pi)$ when $\pi$ is a \catone representation.

\begin{proposition}
\label{IVbmatprop}
Let $\sigma$ be an unramified character of $F^\times$
such that $\sigma^2=1$, and let $\pi$ be the IVb representation
$L(\nu^2,\nu^{-1}\sigma \St_{\GSp(2)})$. Then $\pi$ has trivial
central character, $\pi$ is paramodular, $N_\pi=2$, $\mu_\pi=-q^2+q$, and $\lambda_\pi=\sigma(\varpi)(1+q^2)$,
so that $\pi$ is a \catone representation. Let $v_{\mathrm{new}} \in V(2)$
be a newform, and let $v_s$ be the shadow of $v_{\mathrm{new}}$. The 
vector space $V_s(2)$ is two-dimensional and has ordered basis
\begin{equation}
 \label{IVbmatpropeq1}
v_{\mathrm{new}}, \qquad \tau_1 v_s.
\end{equation}
The matrix of the endomorphism $T_{0,1}^s$ on $V_s(2)$ with respect to \eqref{IVbmatpropeq1}
is 
\begin{equation}
\label{IVbmatpropeq2}
\begin{bsmallmatrix}
q^2(1+q^2)^{-1}\lambda_\pi &0\\
-q(1+q^2)^{-1}\lambda_\pi & \lambda_\pi
\end{bsmallmatrix}
=
\begin{bsmallmatrix}
\sigma(\varpi)q^2&0\\
-\sigma(\varpi)q&\sigma(\varpi)(1+q^2)
\end{bsmallmatrix},
\end{equation}
and the matrix of the endomorphism $T_{1,0}^s$ on $V_s(2)$ with respect to \eqref{IVbmatpropeq1}
is
\begin{equation}
\label{IVbmatpropeq3}
\begin{bsmallmatrix}
0&0\\
-q^2&\mu_\pi +q^2
\end{bsmallmatrix}
=
\begin{bsmallmatrix}
0&0\\
-q^2&q
\end{bsmallmatrix}.
\end{equation}
The characteristic polynomial of $T_{0,1}^s$ on $V_s(2)$ is
\begin{equation}
 \label{IVbmatpropeq31}
p(T_{0,1}^s, V_s(2),X)= (X-\lambda_\pi)(X-q^2(1+q^2)^{-1} \lambda_\pi),
\end{equation}
and the characteristic polynomial of $T_{1,0}^s$ on $V_s(2)$ is
\begin{equation}
p(T_{1,0}^s, V_s(2),X)= (X-(\mu_\pi+q^2))X.
\end{equation}
\end{proposition}
\begin{proof}
That $\pi$ is paramodular, $N_\pi=2$, $\mu_\pi=-q^2+q$, $\lambda_\pi=\sigma(\varpi)(1+q^2)$, and
$V_s(2)$ is two-dimensional follow from Table \ref{dimensionstable} on p.~\pageref{dimensionstable} and 
Table \ref{levelsandeigenvaluestable} on p.~\pageref{levelsandeigenvaluestable}.
To prove that \eqref{IVbmatpropeq1} is a basis for $V_s(2)$ it suffices 
to prove that  $v_{\mathrm{new}}$ and $\tau_1 v_s$ are linearly independent; suppose otherwise.
Then $v_s$ is in the kernel of the map $\tau_1:\bar V_s(1) \to \bar V_s(2)$; this contradicts
the injectivity of this map from \ref{tnlevelraisinglemmaitem2} of Lemma~\ref{taunlevelraisinglemma}.
Next, we prove that the matrix of $T_{1,0}^s$ on $V_s(2)$ with respect to \eqref{IVbmatpropeq1}
is as in \eqref{IVbmatpropeq3}. By \eqref{Wsgenericlemma100} we have $T_{1,0}^s v_{\mathrm{new}}=-q^2 \tau_{1}v_s$.
This verifies the first column of \eqref{IVbmatpropeq3}. 
For the second column, we have
\begin{align}
T_{1,0}^s \tau_1 v_s
&= \tau_1 T_{1,0}^s v_s\qquad \text{(by  \eqref{T10taueq})}\nonumber\\
&= (\mu_\pi+q^2)\tau_1 v_s\qquad \text{(by \eqref{T01generictheoremeq1})}.
\end{align}
Next, we prove that the matrix of $T_{0,1}^s$ on $V_s(2)$ with respect to \eqref{IVbmatpropeq1} is as in 
\eqref{IVbmatpropeq2}. We have
\begin{align}
 T^s_{0,1}v_{\rm new}&=q\mu_\pi^{-1}T^s_{0,1}\rho'_{1}v_s\qquad\text{(by \eqref{vsmulemmaeq2})}\nonumber\\
 &=q\mu_\pi^{-1}(q\theta v_s-q\tau_{N_\pi-1}T^s_{0,1}v_s)\qquad \text{(by \eqref{T01rhoeq})}\nonumber\\
 &=q^2\mu_\pi^{-1}\theta v_s-q^2\lambda_\pi\mu_\pi^{-1}\tau_{1} v_s \qquad\text{(by \eqref{T01generictheoremeq1})}.\label{IVbmatpropeq7}
\end{align}
To proceed further we need to express $\theta v_s$ in the basis \eqref{IVbmatpropeq1}.
Let  $a,b \in \C$ be such that 
\begin{equation}
\label{IVbmatpropeq6}
\theta v_s = a v_{\mathrm{new}} + b \tau_1 v_s.
\end{equation}
Applying $T_{1,0}^s$ to this equation, we obtain
\begin{align*}
T_{1,0}^s \theta v_s &= a T_{1,0}^s v_{\mathrm{new}} + b T_{1,0}^s \tau_1 v_s\\
&=-q^2 a \tau_1v_s + (\mu_\pi+q^2)b  \tau_1 v_s \qquad \text{(by \eqref{IVbmatpropeq3})}\\
q^2T_{0,1}^s \tau_1 v_s&=-q^2 a \tau_1v_s + (\mu_\pi+q^2)b  \tau_1 v_s \qquad \text{(by \eqref{Ts10thetaeq})}\\
q^2\tau_1 T_{0,1}^s  v_s&=-q^2 a \tau_1v_s + (\mu_\pi+q^2)b  \tau_1 v_s \qquad \text{(by \eqref{T01taueq})}\\
q^2\lambda_\pi\tau_1  v_s&=(-q^2 a  + (\mu_\pi+q^2)b  )\tau_1 v_s \qquad \text{(by \eqref{T01generictheoremeq1})}.
\end{align*}
It follows that 
\begin{equation}
\label{IVbmatpropeq21}
q^2\lambda_\pi=-q^2 a  + (\mu_\pi+q^2)b.
\end{equation}
Applying $p_2$ to \eqref{IVbmatpropeq6}, we have:
\begin{align*}
p_2 \theta v_s 
&= a p_2 v_{\mathrm{new}} + b p_2 \tau_1 v_s\\
\theta p_1 v_s 
&= a  v_{\mathrm{new}}  -(q+q^2)^{-1} \mu_\pi b v_{\mathrm{new}}\qquad \text{(by \eqref{commreleq5} and \eqref{vsmulemmaeq1})}\\
0
&= (a    -(q+q^2)^{-1} \mu_\pi b )v_{\mathrm{new}}\qquad \text{(since $V_s(1)=0$)}.
\end{align*}
Hence, 
\begin{equation}
\label{IVbmatpropeq22}
0=a    -(q+q^2)^{-1} \mu_\pi b.
\end{equation}
Solving \eqref{IVbmatpropeq21} and \eqref{IVbmatpropeq22} for $a$ and $b$, we obtain
$$
a =\sigma(\varpi)q(1-q), \qquad b=\sigma(\varpi) q(1+q),
$$
so that
\begin{equation}
\label{IVbmatpropeq23}
\theta v_s = \sigma(\varpi)q(1-q) v_{\mathrm{new}} + \sigma(\varpi)q(1+q) \tau_1 v_s.
\end{equation}
Returning to \eqref{IVbmatpropeq7}, and substituting for $\theta v_s$, we now have:
\begin{align*}
T^s_{0,1}v_{\rm new}
&=q^2\mu_\pi^{-1}\sigma(\varpi)q(1-q) v_{\mathrm{new}} + (q^2\mu_\pi^{-1}\sigma(\varpi)q(1+q)
-q^2\lambda_\pi\mu_\pi^{-1})\tau_{1} v_s\\
&=q^2 (1+q^2)^{-1}\lambda_\pi v_{\mathrm{new}}  - q(1+q^2)^{-1}\lambda_\pi \tau_1 v_s.
\end{align*}
This verifies the first column of \eqref{IVbmatpropeq2}.
The second column of \eqref{IVbmatpropeq2} follows from \eqref{T01taueq} and \eqref{T01generictheoremeq1}.
The final assertions about characteristic polynomials follow by standard computations.
\end{proof}

\section{An upper block algorithm}
\label{algorithmsec}

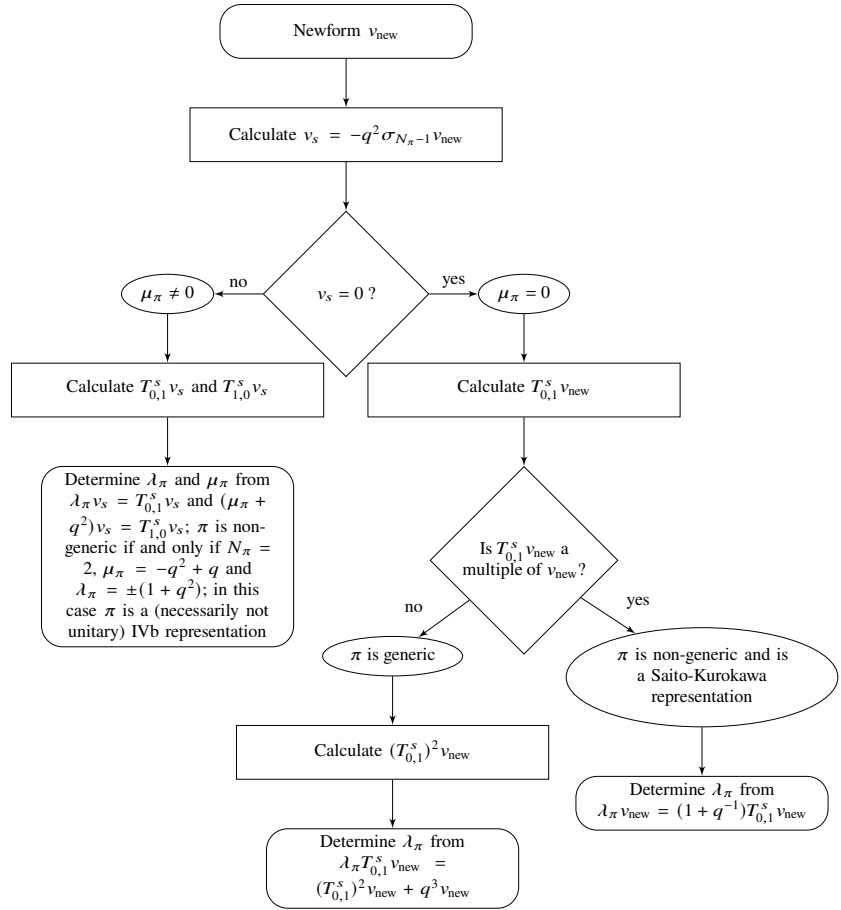
\begin{figure}
\caption{Algorithm for computing paramodular Hecke eigenvalues}
\label{algorithmfig}
\bigskip
\noindent\resizebox{\textwidth}{!}{
\tikzstyle{terminator} = [rectangle, draw,  
    text width=4cm, text centered, rounded corners=1em, minimum height=3em]
\tikzstyle{process} = [rectangle, draw,  
    text width=5cm, text centered, minimum height=3em]
\tikzstyle{decision} = [diamond, draw,  
    text width=8.4em, text badly centered, inner sep=0pt]
\tikzstyle{block} = [rectangle, draw,  
    text width=4.5cm, text centered, rounded corners]
\tikzstyle{action} = [rectangle, draw,  
    text width=4cm, text centered, rounded corners, minimum height=3em]
\tikzstyle{line} = [draw, -latex']
\tikzstyle{longfact} = [draw, ellipse, minimum height=2em, text centered, text width=3cm]  
\tikzstyle{fact} = [draw, ellipse, minimum height=2em] 
\begin{tikzpicture}[->,scale=0.85, every node/.style={scale=0.85},node distance =0.69cm, auto]
\node [terminator] (Wnew) {Newform $v_{\rm new}$};
\node [process, below =of Wnew] (sigma) {Calculate $v_s=-q^2\sigma_{N_\pi-1}v_{\rm new}$};
\node [decision, below =of sigma] (sigmazero) {$v_s=0$ ?};
\node [fact, right =of sigmazero] (muzero) {$\mu_\pi=0$};
\node [fact, left =of sigmazero] (munonzero) {$\mu_\pi\neq0$};
\node [process, below =of munonzero)] (munonzerocalc) {Calculate $T_{0,1}^s v_s$ and $T_{1,0}^s v_s$};
\node [terminator, below =of munonzerocalc)] (getlambdamunonzero) {Determine $\lambda_\pi$ and $\mu_\pi$ from
    $\lambda_\pi v_s=T_{0,1}^s v_s$ and $(\mu_\pi+q^2)v_s=T_{1,0}^s v_s$; $\pi$ is non-generic
    if and only if $N_\pi=2$, $\mu_\pi=-q^2+q$ and $\lambda_\pi=\pm(1+q^2)$; in this case
    $\pi$ is a (necessarily not unitary) IVb representation};
\node [process, below =of muzero)] (muzerocalc) {Calculate $T_{0,1}^s v_{\mathrm{new}}$};
\node [decision, below =of muzerocalc] (iseigen) {Is $T_{0,1}^s v_{\mathrm{new}}$ a multiple of $v_{\mathrm{new}}$?};
\node [longfact,below right =of iseigen] (nongeneric) {$\pi$ is non-generic and is a Saito-Kurokawa representation};
\node [terminator, below  =of nongeneric] (iseigenyes) {Determine $\lambda_\pi$ from $\lambda_\pi v_{\mathrm{new}}
=(1+q^{-1})T_{0,1}^sv_{\mathrm{new}}$};
\node [fact,below left =of iseigen] (generic) {$\pi$ is generic};
\node [process, below  =of generic] (iseigenno) {Calculate $(T_{0,1}^s)^2v_{\mathrm{new}}$};
\node [terminator, below  =of iseigenno] (iseigennogetlambda) {Determine $\lambda_\pi$ from
$\lambda_\pi T_{0,1}^s v_{\mathrm{new}} = (T_{0,1}^s)^2v_{\mathrm{new}}+q^3v_{\mathrm{new}}$};
\path [line] (Wnew) -- (sigma);
\path [line] (sigma) -- (sigmazero);
\path [line] (sigmazero)  --node[anchor=south, above] {yes} (muzero);
\path [line] (sigmazero) --node[anchor=south] {no} (munonzero);
\path [line] (munonzero) -- (munonzerocalc);
\path [line] (muzero) -- (muzerocalc);
\path [line] (munonzerocalc) -- (getlambdamunonzero);
\path [line] (muzerocalc) -- (iseigen);
\path [line] (iseigen) --node[anchor=south,xshift=5mm] {yes} (nongeneric);
\path [line] (iseigen) --node[anchor=south,xshift=-5mm] {no} (generic);
\path [line] (generic) -- (iseigenno);
\path [line] (nongeneric) -- (iseigenyes);
\path [line] (iseigenno) -- (iseigennogetlambda);
\end{tikzpicture}
}
\end{figure}

Let $\pi$ be a paramodular, irreducible, admissible representation of $\GSp(4,F)$
with trivial central character such that $N_\pi \geq 2$, and let $v_{\mathrm{new}} 
\in V(N_\pi)$ be a newform.
In this final section we describe how $\lambda_\pi$ and $\mu_\pi$ can be computed
from $v_{\mathrm{new}}$
using the \index{upper block operator} upper block operators $\sigma_{N_\pi-1}$, $T_{0,1}^s$, and $T_{1,0}^s$; at the same
time, this method also determines whether or not $\pi$ is non-generic.
Our algorithm is based on 
Proposition~\ref{shadowsigmavnewprop} and
Theorem \ref{TSKleveltheorem}, and proceeds as follows. A flow chart version of this algorithm appears in Fig.~\ref{algorithmfig}.
\begin{enumerate}
\item\label{line1} Calculate $v_s=-q^2 \sigma_{N_\pi-1} v_{\mathrm{new}}$ (see \eqref{Wsgenericlemma1eq1}). If $v_s \neq 0$, then
proceed to \ref{line2}. If $v_s =0$, then proceed to \ref{line3}.
\item\label{line2} Since $v_s \neq 0$, $\pi$ is a \catone representation so that
$N_{\pi,s}=N_\pi-1$ and $\mu_\pi \neq 0$ by Proposition \ref{shadowsigmavnewprop}. Calculate
$T_{0,1}^s v_s$ and $T_{1,0}^s v_s$. By Theorem~\ref{TSKleveltheorem} we have that
$$
\lambda_\pi v_s=T_{0,1}^s v_s \quad\text{and}\quad(\mu_\pi+q^2)v_s=T_{1,0}^s v_s.
$$
Solve for $\lambda_\pi$ and $\mu_\pi$. By Lemma \ref{nogenorIVblemma} $\pi$ is non-generic
if and only if $N_\pi=2$, $\mu_\pi=-q^2+q$, and $\lambda_\pi=\pm(1+q^2)$; in this case $\pi$
is a IVb representation.
\item\label{line3} Since $v_s = 0$, $\pi$ is a \cattwo representation so that
$N_{\pi,s}=N_\pi$ and $\mu_\pi = 0$ by Proposition \ref{shadowsigmavnewprop}. Calculate
$T_{0,1}^s v_{\mathrm{new}}$. If $T_{0,1}^s v_{\mathrm{new}}$ is a multiple of $v_{\mathrm{new}}$,
then $\pi$ is non-generic by Theorem \ref{TSKleveltheorem} and is
a Saito-Kurokawa representation by Lemma \ref{cat2nongencharlemma}; proceed to \ref{line4}. If $T_{0,1}^s v_{\mathrm{new}}$ 
is not a multiple of $v_{\mathrm{new}}$, then $\pi$ is generic by Theorem \ref{TSKleveltheorem};
proceed to \ref{line5}. 
\item\label{line4}  By Theorem \ref{TSKleveltheorem}, we have that
$$ 
\lambda_\pi v_{\mathrm{new}}=(1+q^{-1}) T_{0,1}^s v_{\mathrm{new}}.
$$
Solve for $\lambda_\pi$. 
\item\label{line5} Calculate $(T_{0,1})^2 v_{\mathrm{new}}$. By Theorem \ref{TSKleveltheorem}
we have that 
$$
\lambda_\pi T_{0,1}^s v_{\mathrm{new}} = (T_{0,1})^2 v_{\mathrm{new}} + q^3 v_{\mathrm{new}}.
$$
Solve for  $\lambda_\pi$. 
\end{enumerate}

%% file: SKMS_chapter7.tex
\chapter{The Paramodular Subspace}
\label{structuralchap}

Let $(\pi,V)$ be an irreducible, admissible representation of $\GSp(4,F)$
with trivial central character, and let $n$ be an integer such that $n \geq 0$. 
Assume that $\pi$ is paramodular. In this chapter we investigate 
the relationship between $V_s(n)$ and its subspace $V(n-1)+V(n)$ of paramodular vectors. 
From Table \ref{dimensionstable}, it is already evident that for almost
all non-generic\index{non-generic representation}\index{representation!non-generic} $\pi$, 
including all Saito-Kurokawa representations, we have
$V(n-1)+V(n)=V_s(n)$ for all $n \geq 0$; also, from Table \ref{dimensionstable} we
see that if $\pi$ is generic, then $V(n-1)+V(n)$ is a proper subspace of 
$V_s(n)$ if $n \geq \max(N_\pi,1)$. Thus, except for a few non-generic $\pi$,
the investigation of the relation between $V_s(n)$ and its subspace $V(n-1)+V(n)$
is reduced to considering $\pi$ that are generic. In Theorem \ref{genericdimensionstheorem} we 
proved that if $\pi$ is generic and a
\catone representation,\index{category 1 representation}\index{representation!category 1} then 
$$
V_s(n)=V(n-1) \oplus V(n) \oplus \bigoplus_{\substack{i,j \geq 0\\ i+j=n-N_\pi+1}} \C \tau^i \theta^j W_s
$$
where $W_s$ is the shadow of the newform\index{shadow of a newform} $W_{\mathrm{new}}$ in $V(N_\pi)$. 
Most of this chapter is devoted to generalizing this result to all generic $\pi$,
though we also give a complete account for all non-generic $\pi$.

\section{Calculation of certain zeta integrals}
\label{certainzetasec}
Let $(\pi,V)$ be a generic, irreducible, admissible
representation of $\GSp(4,F)$ with trivial central
character, and let $V=\mathcal{W}(\pi,\psi_{c_1,c_2})$,
the Whittaker model of $\pi$ with respect to $\psi_{c_1,c_2}$;
as usual, we assume that $c_1,c_2 \in \OF^\times$ (see Sect.~\ref{repsec}).
In  this section we calculate certain zeta integrals\index{zeta integral} for later use.
As in earlier chapters, we let $\lambda_\pi$, $\mu_\pi$, 
and $\varepsilon_\pi$ be the eigenvalues of $T_{0,1}$, $T_{1,0}$,
and $u_{N_\pi}$, respectively, on the one-dimensional
space $V(N_\pi)$. Since $\pi$ is infinite-dimensional,  paramodular
vectors\index{paramodular vectors} with distinct levels are linearly independent (see Theorem 
\ref{linindparatheorem}); thus, the subspace $V(n-1) + V(n)$
of $V_s(n)$ is a direct sum $V(n-1) \oplus V(n)$ for integers $n \geq 1$.
Also, as in previous chapters, we will use the equivalences from
Proposition \ref{shadowsigmavnewprop}.

We begin with some notation and a summary of some results from
\cite{NF}. 
Let $n$ be an integer such that $n \geq 0$, and let $W \in V(n)$. By Proposition 4.1.3
of \cite{NF} we have
\begin{equation}
\label{zetanewspanlemmaeq31}
Z(s,\theta' W)=q Z(s,W), \qquad Z(s,\theta W)= q^{-s+\frac{3}{2}} Z(s,W)
\end{equation}
and 
\begin{equation}
\label{zetanewspanlemmaeq32}
Z(s,\eta W)=0.
\end{equation}
For all integers $i$ and $j$ we recall that
\begin{equation}\label{Deltaijdefeq}
\Delta_{i,j}= 
\begin{bsmallmatrix}
\varpi^{2i+j}&&&\\ 
& \varpi^{i+j} && \\ 
&& \varpi^i & \\ 
&&& 1
\end{bsmallmatrix}.
\end{equation}
We have $\lambda(\Delta_{i,j})=\varpi^{2i+j}$ for integers $i$ and $j$.
Let $W_{\mathrm{new}} \in V(N_\pi)$ be a newform, i.e., a non-zero
element of the one-dimensional space $V(N_\pi)$. 
For all integers $i$ and $j$ we set
\begin{equation}\label{cijdefeq}
 c_{i,j}=W_{\rm new}(\Delta_{i,j}).
\end{equation}
By Corollary 4.3.8 of \cite{NF} the vector $W_{\rm new}$ is determined by the numbers $c_{i,j}$, and by \eqref{delltavanisheq} we have $c_{i,j}=0$ if $i<0$ or $j<0$. As indicated in  Theorem \ref{basicgenparatheorem}, $Z(s,W_{\mathrm{new}})$ is a
non-zero multiple of $L(s,\pi)$.
\begin{lemma}
\label{cijlemma}
Let $(\pi,V)$ be a generic, irreducible, admissible
representation of $\GSp(4,F)$ with trivial central
character, let $V=\mathcal{W}(\pi,\psi_{c_1,c_2})$,
and let $W_{\mathrm{new}} \in V(N_\pi)$ be a newform. Then
\begin{equation}\label{c00ZsWneweq}
  c_{0,0}=\frac{Z(s,W_{\rm new})}{1-q^{-1}}
  \begin{cases}
   (1-q^{-\frac{3}{2}}\lambda_\pi q^{-s}+(\mu q^{-2}+1+q^{-2})q^{-2s}\\
   \qquad-q^{-\frac{3}{2}}\lambda_\pi q^{-3s}+q^{-4s})&\text{if }N_\pi=0,\\
   (1-q^{-\frac{3}{2}}(\lambda_\pi+\varepsilon_\pi)q^{-s}+(\mu_\pi q^{-2}+1)q^{-2s}\\
   \qquad+\varepsilon_\pi q^{-1/2}q^{-3s})&\text{if }N_\pi=1,\\
   (1-q^{-\frac{3}{2}}\lambda_\pi q^{-s}
  +(\mu_\pi q^{-2}+1)q^{-2s})&\text{if }N_\pi\geq2,
  \end{cases}
\end{equation}
and
\begin{equation}\label{c1jZsWneweq}
 \sum_{j=0}^\infty c_{1,j}(q^{-s+\frac{3}{2}})^j=\frac{Z(s,W_{\rm new})}{1-q^{-1}}
                 \begin{cases}
                  q^{-4}(\mu_\pi+1&\\
                  \qquad-\lambda_\pi q^{\frac{1}{2}}q^{-s}+q^2q^{-2s})&\text{if }N_\pi=0,\\
                  q^{-4}(\mu_\pi+\varepsilon_\pi q^{\frac{3}{2}}q^{-s})&\text{if }N_\pi=1,\\
                  q^{-4}\mu_\pi&\text{if }N_\pi\geq2.
                 \end{cases}
\end{equation}
\end{lemma}
\begin{proof}
Equation \eqref{c00ZsWneweq} follows from  
Proposition~7.1.4, Proposition~7.2.5 and Proposition~7.4.5 of \cite{NF}.
Equation \eqref{c1jZsWneweq} follows from 
Lemma~7.1.3 of \cite{NF} (combined with \eqref{c00ZsWneweq}), 
the first display in the proof of Proposition~7.2.5 of \cite{NF}, and Lemma~7.4.4 of~\cite{NF}.
\end{proof}

\begin{lemma}\label{Ts01Wnewlemma}
 Let $(\pi,V)$ be a generic, irreducible, admissible representation of $\GSp(4,F)$ with trivial central character. Let $W_{\rm new}\in V(N_\pi)$ be a newform. We consider $W_{\rm new}$ as an element of $V_s(n)$, where $n=1$ if $N_\pi=0$ and $n=N_\pi$ if $N_\pi\geq1$. Then
 \begin{align}\label{Ts01Wnewlemmaeq1}
  &Z(s,T^s_{0,1}W_{\rm new})=\begin{cases}
                  (\lambda_\pi-q^{-s+\frac{3}{2}})Z(s,W_{\rm new})&\text{if }N_\pi=0,\\
                  (\lambda_\pi+\varepsilon_\pi-q^{-s+\frac{3}{2}})Z(s,W_{\rm new})&\text{if }N_\pi=1,\\
                  (\lambda_\pi-q^{-s+\frac{3}{2}})Z(s,W_{\rm new})&\text{if }N_\pi\geq2,
                 \end{cases}
 \end{align}
 and
 \begin{align}\label{Ts01Wnewlemmaeq10}
  &Z(s,T^s_{1,0}W_{\rm new})=
                 \begin{cases}
                  (\mu_\pi+1-\lambda_\pi q^{\frac{1}{2}}q^{-s}+q^2q^{-2s})Z(s,W_{\rm new})&\text{if }N_\pi=0,\\
                  (\mu_\pi+\varepsilon_\pi q^{\frac{3}{2}}q^{-s})Z(s,W_{\rm new})&\text{if }N_\pi=1,\\
                  \mu_\pi Z(s,W_{\rm new})&\text{if }N_\pi\geq2.
                 \end{cases}
 \end{align}
\end{lemma}
\begin{proof}
By \eqref{Vszetaintposdef}, 
$$
Z(s,T^s_{0,1}W_{\rm new})=\int\limits_{\substack{F^\times\\v(a)\geq0}}(T^s_{0,1}W_{\rm new})(\begin{bsmallmatrix}a\\&a\\&&1\\&&&1\end{bsmallmatrix})|a|^{s-\frac{3}{2}}\,d^\times a.
$$
Using \eqref{Ts01eq}, we have:
\begin{align*}
&Z(s,T^s_{0,1}W_{\rm new})\\
&\qquad=\sum_{y,z\in\OF/\p}\int\limits_{\substack{F^\times\\v(a)\geq0}}W_{\rm new}(\begin{bsmallmatrix}
a\vphantom{\varpi^{-n+1}}\\
&a\vphantom{\varpi^{-n+1}}\\
&&1\vphantom{\varpi^{-n+1}}\\
&&&1\vphantom{\varpi^{-n+1}}
\end{bsmallmatrix}
\begin{bsmallmatrix}
1&y&&z\varpi^{-n+1}\\
&1&&\\
&&1&-y\\
&&&1
\end{bsmallmatrix}
\begin{bsmallmatrix}
\varpi\vphantom{\varpi^{-n+1}}\\
&1\vphantom{\varpi^{-n+1}}\\
&&\varpi\vphantom{\varpi^{-n+1}}\\
&&&1\vphantom{\varpi^{-n+1}}
\end{bsmallmatrix})|a|^{s-\frac{3}{2}}\,d^\times a\\
&\qquad\quad+\sum_{c,y,z\in\OF/\p}\int\limits_{\substack{F^\times\\v(a)\geq0}}W_{\rm new}(
\begin{bsmallmatrix}
a\vphantom{\varpi^{-n+1}}\\
&a\vphantom{\varpi^{-n+1}}\\
&&1\vphantom{\varpi^{-n+1}}\\
&&&1\vphantom{\varpi^{-n+1}}
\end{bsmallmatrix}
\begin{bsmallmatrix}
1&&y&z\varpi^{-n+1}\vphantom{\varpi^{-n+1}}\\
&1&c&y\vphantom{\varpi^{-n+1}}\\
&&1\vphantom{\varpi^{-n+1}}\\
&&&1\vphantom{\varpi^{-n+1}}
\end{bsmallmatrix}
\begin{bsmallmatrix}
\varpi\vphantom{\varpi^{-n+1}}\\
&\varpi\vphantom{\varpi^{-n+1}}\\
&&1\vphantom{\varpi^{-n+1}}\\
&&&1\vphantom{\varpi^{-n+1}}
\end{bsmallmatrix})|a|^{s-\frac{3}{2}}\,d^\times a\\
&\qquad=q^2\int\limits_{\substack{F^\times\\v(a)\geq0}}W_{\rm new}(
\begin{bsmallmatrix}a\\&a\\&&1\\&&&1\end{bsmallmatrix}
\begin{bsmallmatrix}\varpi\\&1\\&&\varpi\\&&&1\end{bsmallmatrix})
|a|^{s-\frac{3}{2}}\,d^\times a\\
&\qquad\quad+q^3\int\limits_{\substack{F^\times\\v(a)\geq0}}W_{\rm new}(
\begin{bsmallmatrix}a\\&a\\&&1\\&&&1\end{bsmallmatrix}
\begin{bsmallmatrix}\varpi\\&\varpi\\&&1\\&&&1\end{bsmallmatrix})
|a|^{s-\frac{3}{2}}\,d^\times a\\
&\qquad=q^2(1-q^{-1})\sum_{j=0}^\infty W_{\rm new}(\begin{bsmallmatrix}\varpi^{j+1}\\&\varpi^j\\&&\varpi\\&&&1\end{bsmallmatrix}
)q^{-j(s-\frac{3}{2})}\\
&\qquad\quad+q^3(1-q^{-1})\sum_{j=0}^\infty W_{\rm new}(\begin{bsmallmatrix}\varpi^{j+1}\\&\varpi^{j+1}\\&&1\\&&&1\end{bsmallmatrix})q^{-j(s-\frac{3}{2})}\\
&\qquad=q^2(1-q^{-1})\sum_{j=0}^\infty c_{1,j-1}q^{-j(s-\frac{3}{2})}+q^3(1-q^{-1})\sum_{j=0}^\infty c_{0,j+1}q^{-j(s-\frac{3}{2})}\\
&\qquad=q^2(1-q^{-1})\sum_{j=-1}^\infty c_{1,j}q^{-(j+1)(s-\frac{3}{2})}+q^3(1-q^{-1})\sum_{j=1}^\infty c_{0,j}q^{-(j-1)(s-\frac{3}{2})}\\
&\qquad=q^{-s+\frac{7}{2}}(1-q^{-1})\sum_{j=-1}^\infty c_{1,j}q^{-j(s-\frac{3}{2})}+q^{s+\frac{3}{2}}(1-q^{-1})\sum_{j=1}^\infty c_{0,j}q^{-j(s-\frac{3}{2})}\\
&\qquad=q^{-s+\frac{7}{2}}(1-q^{-1})\sum_{j=0}^\infty c_{1,j}q^{-j(s-\frac{3}{2})}\\
&\qquad\quad+q^{s+\frac{3}{2}}(1-q^{-1})\!\Big(\sum_{j=0}^\infty c_{0,j}q^{-j(s-\frac{3}{2})}-c_{0,0}\Big)\\
&\qquad=q^{-s+\frac{7}{2}}(1-q^{-1})\sum_{j=0}^\infty c_{1,j}q^{-j(s-\frac{3}{2})}\\
&\qquad\quad+q^{s+\frac{3}{2}}Z(s,W_{\rm new})-q^{s+\frac{3}{2}}(1-q^{-1})c_{0,0}.
\end{align*}
Substituting from \eqref{c00ZsWneweq} and \eqref{c1jZsWneweq} proves \eqref{Ts01Wnewlemmaeq1}.
Similarly, by \eqref{Ts10eq}, 
\begin{align*}
&Z(s,T^s_{1,0}W_{\rm new})=\int\limits_{\substack{F^\times\\v(a)\geq0}}(T^s_{1,0}W_{\rm new})(\begin{bsmallmatrix}a\\&a\\&&1\\&&&1\end{bsmallmatrix})|a|^{s-\frac{3}{2}}\,d^\times a\\
&\qquad=\int\limits_{\substack{F^\times\\v(a)\geq0}}\sum_{\substack{x,y\in\OF/\p\\z\in\OF/\p^2}}W_{\rm new}(\begin{bsmallmatrix}
a\vphantom{\varpi^{-n+1}}\\
&a\vphantom{\varpi^{-n}}\\
&&1\vphantom{\varpi^{-n}}\\
&&&1\vphantom{\varpi^{-n}}
\end{bsmallmatrix}
\begin{bsmallmatrix}
1&x&y&z\varpi^{-n+1}\\
&1&&y\\
&&1&-x\\
&&&1
\end{bsmallmatrix}
\begin{bsmallmatrix}
\varpi^2\vphantom{\varpi^{-n}}\\
&\varpi\vphantom{\varpi^{-n}}\\
&&\varpi\vphantom{\varpi^{-n}}\\
&&&1\vphantom{\varpi^{-n}}
\end{bsmallmatrix})|a|^{s-\frac{3}{2}}\,d^\times a\\
&\qquad=q^4\int\limits_{\substack{F^\times\\v(a)\geq0}}W_{\rm new}(
\begin{bsmallmatrix}a\\&a\\&&1\\&&&1\end{bsmallmatrix}
\begin{bsmallmatrix}\varpi^2\\&\varpi\\&&\varpi\\&&&1\end{bsmallmatrix}
)|a|^{s-\frac{3}{2}}\,d^\times a\\
&\qquad=q^4(1-q^{-1})\sum_{j=0}^\infty W_{\rm new}(\begin{bsmallmatrix}\varpi^{2+j}\\&\varpi^{1+j}\\&&\varpi\\&&&1\end{bsmallmatrix})q^{-j(s-\frac{3}{2})}\\
&\qquad=q^4(1-q^{-1})\sum_{j=0}^\infty c_{1,j}q^{-j(s-\frac{3}{2})}.
\end{align*}
Hence \eqref{Ts01Wnewlemmaeq10} follows from \eqref{c1jZsWneweq}.
\end{proof}

\section{Generic representations}
\label{structureSKsec}
The main result of this section is Theorem \ref{W0theorem}, giving the structure of the spaces $V_s(n)$ for all \index{generic representation}\index{representation!generic}generic, irreducible, admissible representations of $\GSp(4,F)$ with trivial central character. As a consequence, we will show that the linear dimension growth for the spaces $\bar V_s(n)$ proven in \eqref{genericdimensionstheoremeq6} for $N_\pi\geq2$ and $\mu_\pi\neq0$ holds in all cases. We also show that the maps $\theta:V_s(n)\to V_s(n+1)$ and $\theta:\bar V_s(n)\to\bar V_s(n+1)$ are injective.

\begin{lemma}\label{zetanewspanlemma}
 Let $\pi$ be a generic, irreducible, admissible generic representation of $\GSp(4,F)$ with trivial central character. Let $W_{\rm new}\in V(N_\pi)$ be a newform.
 \begin{enumerate}
  \item \label{zetanewspanlemmaitem1} Assume that $N_\pi=0$. Then
   \begin{align}
    T^s_{0,1}W_{\rm new}&\in\langle W_{\rm new},\theta W_{\rm new},\theta'W_{\rm new}\rangle,\label{zetanewspanlemmaeq1}\\
    T^s_{1,0}W_{\rm new}&\notin\langle W_{\rm new},\theta W_{\rm new},\theta'W_{\rm new}\rangle.\label{zetanewspanlemmaeq11}
   \end{align}
   Here, we view $W_{\rm new}$ as an element of $V_s(1)$ and apply the Hecke operators $T^s_{0,1},T^s_{1,0}$ at this level.
  \item \label{zetanewspanlemmaitem2} Assume that $N_\pi=1$. Then
   \begin{equation}\label{zetanewspanlemmaeq2}
    T^s_{0,1}W_{\rm new}\notin\C W_{\rm new}\quad\text{and}\quad T^s_{1,0}W_{\rm new}\notin\C W_{\rm new}.
   \end{equation}
  \item \label{zetanewspanlemmaitem3} Assume that $N_\pi\geq2$. Then
   \begin{equation}\label{zetanewspanlemmaeq3}
    T^s_{0,1}W_{\rm new}\notin\C W_{\rm new}.
   \end{equation}
 \end{enumerate}
\end{lemma}
\begin{proof}
\ref{zetanewspanlemmaitem1}.  By  \eqref{Ts01Wnewlemmaeq10} we see that \eqref{zetanewspanlemmaeq11} holds. 
Next assume that $T^s_{0,1}W_{\rm new}\notin\langle W_{\rm new},\theta W_{\rm new},\theta'W_{\rm new}\rangle$.
By \eqref{zetanewspanlemmaeq31} the vectors
$\theta W_{\rm new}$ and $\theta'W_{\rm new}$ are linearly independent elements of $V(N_\pi+1)$; since 
$W_{\rm new} \in V(N_\pi)$, the  vectors $W_{\rm new},\theta W_{\rm new}$, and $\theta'W_{\rm new}$ are linearly independent by Theorem \ref{linindparatheorem}. By our assumption, we now have that 
$W_{\rm new},\theta W_{\rm new},\theta'W_{\rm new}$ and $T^s_{0,1}W_{\rm new}$ span the four-dimensional space $V_s(1)$ (see Table \ref{dimensionstable}). By \eqref{Ts01Wnewlemmaeq1} we then have $Z(s,V_s(1))=(\C+\C q^{-s})Z(s,W_{\rm new})$, contradicting \eqref{Ts01Wnewlemmaeq10}. This verifies \eqref{zetanewspanlemmaeq1}.

\ref{zetanewspanlemmaitem2} is immediate from \eqref{Ts01Wnewlemmaeq1} and \eqref{Ts01Wnewlemmaeq10}.

\ref{zetanewspanlemmaitem3} is immediate from \eqref{Ts01Wnewlemmaeq1}.
\end{proof}

Let $(\pi,V)$ be an irreducible, admissible representation of $\GSp(4,F)$ with trivial central character, and let 
$n$ be an integer such that $n \geq 0$. 
We recall some definitions from Sect.~\ref{stableklingenvecsec}.
By definition, $\bar V_s(n)=V_s(n)/(V(n-1)\oplus V(n))$ if $n \geq 1$, and by definition $V_s(0)=0$.  
If $\bar V_s(n)\neq0$ for some $n \geq 0$, then we let $\bar N_{\pi,s}$ be the smallest such integer $n$;
if $\bar V_s(n)=0$ for all $n \geq 0$, then we say that $\bar N_{\pi,s}$ is not defined.
Clearly, $\bar N_{\pi,s} \geq 1$.
The numbers $\bar N_{\pi,s}$  are listed in Table \ref{dimensionstable} for all $\pi$. We note that if $\pi$ 
is generic, then $\bar N_{\pi,s}$ is defined, and $\bar N_{\pi,s}=N_{\pi,s}$  except if $\pi$ is an unramified type I representation; in the latter case, $N_{\pi,s}=0$ and $\bar N_{\pi,s}=1$. We can now prove the following generalization of Theorem \ref{genericdimensionstheorem}.
In \eqref{W0theoremeq10} of the following theorem, to simplify notation, we write $\tau$
for the level
raising operator $\tau_k: V_s(k)\to V_s(k+1)$ for each integer $k \geq 0$; we also note that by \eqref{commreleq1}
the operators $\tau$ and $\theta$ commute. \label{taurefdef}

\begin{theorem}\label{W0theorem}
 Let $(\pi,V)$ be a generic, irreducible, admissible representation of $\GSp(4,F)$ with trivial central character. Then $\bar V_s(n)\neq0$ for some integer $n\geq 0$. Let $W_0$ be an element of $V_s(\bar N_{\pi,s})$ that is not in $V(\bar N_{\pi,s}-1)+V(\bar N_{\pi,s})$, so that $W_0$ represents a non-zero element of $\bar V_s(\bar N_{\pi,s})$. Then
 \begin{equation}\label{W0theoremeq1}
  Z(s,W_0)=P_0(q^{-s})Z(s,W_{\rm new}),
 \end{equation}
 where $P_0(X)\in\C[X]$ is a non-zero polynomial of degree \label{P0def}
 \begin{equation}\label{W0theoremeq2}
  \deg P_0=\begin{cases}
            2&\text{if }N_\pi=0,\\
            1&\text{if }N_\pi=1,\\
            1&\text{if $N_\pi\geq2$ and $\mu_\pi=0$},\\
            0&\text{if $N_\pi\geq2$ and $\mu_\pi \neq0$}.
           \end{cases}
 \end{equation}
 We have
 \begin{equation}\label{W0theoremeq3}
  \bar N_{\pi,s}=N_\pi-1+\deg P_0.
 \end{equation}
 Furthermore,
 \begin{equation}\label{W0theoremeq10}
  V_s(n)=V(n-1)\oplus V(n)\oplus\bigoplus_{\substack{i,j\geq0\\i+j=n-\bar N_{\pi,s}}}\C\tau^i\theta^{j}W_0
 \end{equation}
 and
 \begin{equation}\label{W0theoremeq15}
  \dim V_s(n)=\frac{(n-N_\pi+2)(n-N_\pi+3)}2-\deg(P_0)
 \end{equation}
 for integers $n\geq N_{\pi,s}$. 
In  \eqref{W0theoremeq10}, if $N_{\pi,s}=0$ and $n=0$, then we take $V(n-1)=V(-1)$ to be the zero subspace,
and the direct sum over $i,j \geq 0$, $i+j=n-\bar N_{\pi,s}=-1$ is empty.
 In addition,
 \begin{equation}\label{W0theoremeq20}
  Z(s,V_s(n))=(\C+\C q^{-s}+\ldots+\C(q^{-s})^{n-N_\pi+1})Z(s,W_{\rm new})
 \end{equation}
 for integers $n\geq\bar N_{\pi,s}$.
\end{theorem}
\begin{proof}
It follows from Table \ref{dimensionstable} that $\bar V_s(n)\neq0$ for some $n$.

Assume that $N_\pi=0$. Then $\bar N_{\pi,s}=1$ and $\dim V_s(1)=4$ by Table \ref{dimensionstable}. Lemma~\ref{zetanewspanlemma} implies that $V_s(1)=\langle W_{\rm new},\theta W_{\rm new},\theta'W_{\rm new},T^s_{1,0}W_{\rm new}\rangle$. The vector $W_0$, when expressed as a linear combination of $W_{\rm new},\theta W_{\rm new},\theta'W_{\rm new},T^s_{1,0}W_{\rm new}$, must have a non-zero $T^s_{1,0}W_{\rm new}$-component. Hence \eqref{W0theoremeq1} and \eqref{W0theoremeq2} follow from \eqref{zetanewspanlemmaeq31} and~\eqref{Ts01Wnewlemmaeq10}.

Assume that $N_\pi=1$, or that $N_\pi\geq2$ and $\mu_\pi=0$. Then $\bar N_{\pi,s}=N_\pi$ and $\dim V_s(N_\pi)=2$ by Table \ref{dimensionstable}. Lemma \ref{zetanewspanlemma} implies  $V_s(N_\pi)=\langle W_{\rm new},T^s_{0,1}W_{\rm new}\rangle$. The vector $W_0$, when expressed as a linear combination of $W_{\rm new}$ and $T^s_{0,1}W_{\rm new}$, must have a non-zero $T^s_{0,1}W_{\rm new}$-component. Hence \eqref{W0theoremeq1}  and \eqref{W0theoremeq2} follow from \eqref{zetanewspanlemmaeq31} and~\eqref{Ts01Wnewlemmaeq1}.

Assume that $N_\pi\geq2$ and $\mu_\pi\neq0$. Then $\bar N_{\pi,s}=N_\pi-1$ by Table \ref{dimensionstable}. Let $W_s$ be the shadow of $W_{\mathrm{new}}$\index{shadow of a newform} as defined in \eqref{shadowdefeq2}. By Theorem \ref{genericdimensionstheorem}, $W_s$ spans the one-dimensional space $V_s(N_\pi-1)$. Hence $W_0$ is a multiple of $W_s$. The assertions \eqref{W0theoremeq1} and \eqref{W0theoremeq2} follow from \eqref{WsWnewzetaeq}.

It is now easily verified that \eqref{W0theoremeq3} holds in all cases.

Next we prove \eqref{W0theoremeq10}. Evidently, \eqref{W0theoremeq10} holds for $N_\pi=0$ and $n=0$. We may therefore assume that $n>0$. A calculation using Table \ref{dimensionstable} shows that
$$
\dim V_s(n) = \dim V(n-1) + \dim V(n) + (n-\bar N_{\pi,s}+1).
$$
For this, it is useful to note that $\dim V(m)=\lfloor (m-N_\pi+2)^2/4\rfloor$ for integers
$m$ such that $m \geq N_\pi$ by Theorem 7.5.6 of \cite{NF} and that $\lfloor k^2/4 \rfloor + \lfloor (k+1)^2/4\rfloor
=k(k+1)/2$ for any integer $k$.
It follows that we only need to prove that the sum on the right hand side of \eqref{W0theoremeq10} is direct. For this we use induction on $n$. By the choice of $W_0$, the sum is direct for $n=\max(\bar N_{\pi,s},1)=\max(N_{\pi,s},1)$. (The only case where $\bar N_{\pi,s}$ differs from $N_{\pi,s}$ is for $N_\pi=0$, in which case $\bar N_{\pi,s}=1$.) Assume that $n>\max(N_{\pi,s},1)$, and that the statement is true for $n-1$. Suppose that
\begin{equation}\label{W0theoremeq11}
 W_1+W_2+\sum_{\substack{i,j\geq0\\i+j=n-\bar N_{\pi,s}}}c_i\tau^i\theta^{j}W_0=0
\end{equation}
for some $W_1\in V(n-1)$, $W_2\in V(n)$ and complex numbers $c_0,\ldots, c_{n-\bar N_{\pi,s}}$. Taking zeta integrals and observing \eqref{W0theoremeq1} and Lemma \ref{zetarelationslemma}, we get
\begin{equation}\label{W0theoremeq12}
 Z(s,W_1)+Z(s,W_2)=-P_0(q^{-s})\sum_{\substack{i,j\geq0\\i+j=n-\bar N_{\pi,s}}}(q^{-s+\frac{3}{2}})^jc_i Z(s,W_{\rm new}).
\end{equation}
By \eqref{zetanewspanlemmaeq31} and Theorem 7.5.6 of \cite{NF} we have
\begin{equation}\label{W0theoremeq14}
 Z(s,V(n-1)+V(n))= (\C+\C q^{-s}+\ldots+\C(q^{-s})^{n-N_\pi})Z(s,W_{\rm new}).
\end{equation}
It therefore follows from \eqref{W0theoremeq12} and \eqref{W0theoremeq3} that $c_0=0$. From \eqref{W0theoremeq11} we hence obtain
\begin{equation}\label{W0theoremeq13}
 W_1+W_2+\tau_{n-1}\sum_{\substack{i,j\geq0\\i+j=n-\bar N_{\pi,s}-1}}c_{i+1}\tau^i\theta^{j}W_s=0.
\end{equation}
The map $\tau_{n-1}:\bar V_s(n-1)\to\bar V_s(n)$ is injective by Lemma \ref{taunlevelraisinglemma}. It follows that
$$
 \sum_{\substack{i,j\geq0\\i+j=n-\bar N_{\pi,s}-1}}c_{i+1}\tau^i\theta^{j}W_s\in V(n-2)\oplus V(n-1).
$$
By the induction hypothesis, $c_1=\ldots=c_{n-\bar N_{\pi,s}}=0$. Hence $W_1+W_2=0$, which implies $W_1=W_2=0$
by Theorem \ref{linindparatheorem}. This completes the proof of \eqref{W0theoremeq10}. 

Equation \eqref{W0theoremeq15} follows by taking dimensions on both sides of \eqref{W0theoremeq10}, observing \eqref{W0theoremeq3} and using
 $\dim V(m)=\lfloor (m-N_\pi+2)^2/4\rfloor$ for integers
$m$ such that $m \geq N_\pi$ by Theorem 7.5.6 of \cite{NF}, and  $\lfloor k^2/4 \rfloor + \lfloor (k+1)^2/4\rfloor
=k(k+1)/2$ for any integer~$k$.
(Alternatively, \eqref{W0theoremeq15} can be verified from Table \ref{dimensionstable}.) 

Finally,  \eqref{W0theoremeq20} follows by taking zeta integrals of both sides of \eqref{W0theoremeq10}, observing Lemma \ref{zetarelationslemma}, \eqref{W0theoremeq3} and \eqref{W0theoremeq14}.
\end{proof}

We note that, by the proof of Theorem \ref{W0theorem}, possible choices for the vector $W_0$ are as follows:
\begin{equation}\label{W0choiceeq}
 W_0=\begin{cases}
      T^s_{1,0}W_{\rm new}&\text{if $N_\pi=0$},\\
      T^s_{0,1}W_{\rm new}\text{ or }T^s_{1,0}W_{\rm new}&\text{if $N_\pi=1$},\\
      T^s_{0,1}W_{\rm new}&\text{if $N_\pi\geq2$ and $\mu_\pi=0$},\\
      W_s&\text{if $N_\pi\geq2$ and $\mu_\pi\neq0$}.
     \end{cases}
\end{equation}
Here, $W_s$ is the shadow of $W_{\mathrm{new}}$ as defined in \eqref{shadowdefeq2}.

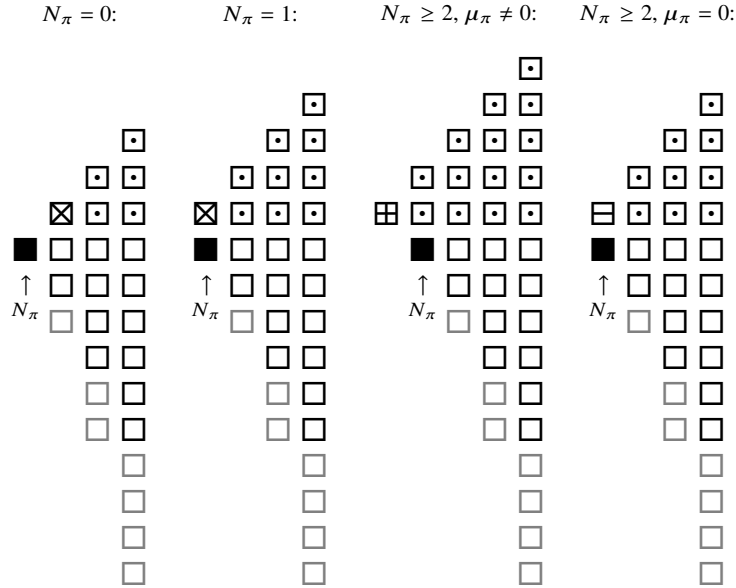
\begin{figure}
\caption{Visualization of the generic cases}\label{genericfig}
\bigskip
\newcommand{\sqr}{\scalebox{2}{$\square$}}
\newcommand{\nsqr}{\scalebox{2}{$\blacksquare$}}
\newcommand{\esqr}{\textcolor{gray}{\scalebox{2}{$\square$}}}
\newcommand{\psqr}{\textcolor{white}{\scalebox{2}{$\square$}}}
\newcommand{\rsqr}{\scalebox{2}{$\boxdot$}}
\newcommand{\nrsqr}{\scalebox{2}{$\boxtimes$}}
\newcommand{\gsqr}{\scalebox{2}{$\boxdot$}}
\newcommand{\ngsqr}{\scalebox{2}{$\boxminus$}}
\newcommand{\bsqr}{\scalebox{2}{$\boxdot$}}
\newcommand{\nbsqr}{\scalebox{2}{$\boxplus$}}
\begin{tikzpicture}[scale=0.48]
\node at (1,10) {\nsqr};
\node at (2,11) {\nrsqr};
\node at (2,10) {\sqr};
\node at (2,9) {\sqr};
\node at (2,8) {\esqr};
\node at (3,12) {\rsqr};
\node at (3,11) {\rsqr};
\node at (3,10) {\sqr};
\node at (3,9) {\sqr};
\node at (3,8) {\sqr};
\node at (3,7) {\sqr};
\node at (3,6) {\esqr};
\node at (3,5) {\esqr};
\node at (4,13) {\rsqr};
\node at (4,12) {\rsqr};
\node at (4,11) {\rsqr};
\node at (4,10) {\sqr};
\node at (4,9) {\sqr};
\node at (4,8) {\sqr};
\node at (4,7) {\sqr};
\node at (4,6) {\sqr};
\node at (4,5) {\sqr};
\node at (4,4) {\esqr};
\node at (4,3) {\esqr};
\node at (4,2) {\esqr};
\node at (4,1) {\esqr};
\node at (6,11) {\nrsqr};
\node at (6,10) {\nsqr};
\node at (7,12) {\rsqr};
\node at (7,11) {\rsqr};
\node at (7,10) {\sqr};
\node at (7,9) {\sqr};
\node at (7,8) {\esqr};
\node at (8,13) {\rsqr};
\node at (8,12) {\rsqr};
\node at (8,11) {\rsqr};
\node at (8,10) {\sqr};
\node at (8,9) {\sqr};
\node at (8,8) {\sqr};
\node at (8,7) {\sqr};
\node at (8,6) {\esqr};
\node at (8,5) {\esqr};
\node at (9,14) {\rsqr};
\node at (9,13) {\rsqr};
\node at (9,12) {\rsqr};
\node at (9,11) {\rsqr};
\node at (9,10) {\sqr};
\node at (9,9) {\sqr};
\node at (9,8) {\sqr};
\node at (9,7) {\sqr};
\node at (9,6) {\sqr};
\node at (9,5) {\sqr};
\node at (9,4) {\esqr};
\node at (9,3) {\esqr};
\node at (9,2) {\esqr};
\node at (9,1) {\esqr};
\node at (11,11) {\nbsqr};
\node at (12,12) {\rsqr};
\node at (12,11) {\rsqr};
\node at (12,10) {\nsqr};
\node at (13,13) {\rsqr};
\node at (13,12) {\rsqr};
\node at (13,11) {\rsqr};
\node at (13,10) {\sqr};
\node at (13,9) {\sqr};
\node at (13,8) {\esqr};
\node at (14,14) {\rsqr};
\node at (14,13) {\rsqr};
\node at (14,12) {\rsqr};
\node at (14,11) {\rsqr};
\node at (14,10) {\sqr};
\node at (14,9) {\sqr};
\node at (14,8) {\sqr};
\node at (14,7) {\sqr};
\node at (14,6) {\esqr};
\node at (14,5) {\esqr};
\node at (15,15) {\rsqr};
\node at (15,14) {\rsqr};
\node at (15,13) {\rsqr};
\node at (15,12) {\rsqr};
\node at (15,11) {\rsqr};
\node at (15,10) {\sqr};
\node at (15,9) {\sqr};
\node at (15,8) {\sqr};
\node at (15,7) {\sqr};
\node at (15,6) {\sqr};
\node at (15,5) {\sqr};
\node at (15,4) {\esqr};
\node at (15,3) {\esqr};
\node at (15,2) {\esqr};
\node at (15,1) {\esqr};
\node at (17,11) {\ngsqr};
\node at (17,10) {\nsqr};
\node at (18,12) {\rsqr};
\node at (18,11) {\rsqr};
\node at (18,10) {\sqr};
\node at (18,9) {\sqr};
\node at (18,8) {\esqr};
\node at (19,13) {\rsqr};
\node at (19,12) {\rsqr};
\node at (19,11) {\rsqr};
\node at (19,10) {\sqr};
\node at (19,9) {\sqr};
\node at (19,8) {\sqr};
\node at (19,7) {\sqr};
\node at (19,6) {\esqr};
\node at (19,5) {\esqr};
\node at (20,14) {\rsqr};
\node at (20,13) {\rsqr};
\node at (20,12) {\rsqr};
\node at (20,11) {\rsqr};
\node at (20,10) {\sqr};
\node at (20,9) {\sqr};
\node at (20,8) {\sqr};
\node at (20,7) {\sqr};
\node at (20,6) {\sqr};
\node at (20,5) {\sqr};
\node at (20,4) {\esqr};
\node at (20,3) {\esqr};
\node at (20,2) {\esqr};
\node at (20,1) {\esqr};
\node at (1,9) {$\uparrow$};
\node at (6,9) {$\uparrow$};
\node at (12,9) {$\uparrow$};
\node at (17,9) {$\uparrow$};
\node[scale=0.9] at (1,8.3) {$N_\pi$};
\node[scale=0.9] at (6,8.3) {$N_\pi$};
\node[scale=0.9] at (12,8.3) {$N_\pi$};
\node[scale=0.9] at (17,8.3) {$N_\pi$};
\node[scale=1] at (2.5,16.5) {$N_\pi=0$:};
\node[scale=1] at (7.5,16.5) {$N_\pi=1$:};
\node[scale=1] at (13,16.5) {$N_\pi\geq 2$, $\mu_\pi \neq 0$:};
\node[scale=1] at (18.5,16.5) {$N_\pi\geq 2$, $\mu_\pi = 0$:};
\end{tikzpicture}\\
In each of the four diagrams, a vertical column of squares
represents a basis for $V_s(n)$, beginning with $n=N_{\pi,s}$.
The squares are as follows:\\
\parbox{0cm}{\begin{tabbing}
 \raisebox{-0.2ex}[0.2ex]{\scalebox{0.6}{\nsqr}}\quad\=The paramodular newform.\\
 \raisebox{-0.2ex}[0.2ex]{\scalebox{0.6}{\sqr}}\quad\>Descendants of \raisebox{-0.2ex}[0.2ex]{\scalebox{0.6}{\nsqr}} via $\theta,\theta',\eta$, making up the spaces $V(n)$ inside $V_s(n)$.\\
 \raisebox{-0.2ex}[0.2ex]{\scalebox{0.6}{\esqr}}\quad\>The spaces $V(n-1)$ inside $V_s(n)$.\\
 \raisebox{-0.2ex}[0.2ex]{\scalebox{0.6}{\nrsqr}}\quad\>The vector $W_0=T^s_{1,0}(\text{\raisebox{-0.2ex}[0.2ex]{\scalebox{0.6}{\nsqr}}})$ in the cases $N_\pi=0$ and $N_\pi=1$.\\
 \raisebox{-0.2ex}[0.2ex]{\scalebox{0.6}{\nbsqr}}\quad\>The shadow vector.\\
 \raisebox{-0.2ex}[0.2ex]{\scalebox{0.6}{\ngsqr}}\quad\>The vector $W_0=T^s_{0,1}(\text{\raisebox{-0.2ex}[0.2ex]{\scalebox{0.6}{\nsqr}}})$.\\
 \raisebox{-0.2ex}[0.2ex]{\scalebox{0.6}{\rsqr}}\quad\>Descendants of \raisebox{-0.2ex}[0.2ex]{\scalebox{0.6}{\nrsqr}},
 \raisebox{-0.2ex}[0.2ex]{\scalebox{0.6}{\nbsqr}}, or \raisebox{-0.2ex}[0.2ex]{\scalebox{0.6}{\ngsqr}} via $\theta$ and $\tau$.
\end{tabbing}}
\end{figure}

The various cases of Theorem~\ref{W0theorem} are illustrated in Fig.~\ref{genericfig}.

\begin{corollary}\label{quotientstructurecor}
Let $(\pi,V)$ be a generic, irreducible, admissible representation of $\GSp(4,F)$ with trivial central character. 
Then
  \begin{equation}\label{quotientstructurecoreq1}
   \bar N_{\pi,s}=\begin{cases}
                   1&\text{if }N_\pi=0\text{ or }N_\pi=1,\\
                   N_\pi-1&\text{if $N_\pi\geq2$ and $\mu_\pi\neq0$},\\
                   N_\pi&\text{if $N_\pi\geq2$ and $\mu_\pi=0$},
                  \end{cases}
  \end{equation}
  and
  \begin{equation}\label{quotientstructurecoreq1b}
   \dim\bar V_s(n)=
   \begin{cases}
    n-\bar N_{\pi,s}+1&\text{for }n\geq \bar N_{\pi,s},\\
    0&\text{for }n<\bar N_{\pi,s}.
   \end{cases}
  \end{equation}
Let $n$ be an integer such that $n \geq 0$. 
  The operators $\tau,\theta:V_s(n)\to V_s(n+1)$, and also $\tau,\theta:\bar V_s(n)\to\bar V_s(n+1)$, are injective and satisfy $\theta\tau=\tau\theta$. Let $W_{s,{\rm new}}$ be a non-zero element of the one-dimensional space $\bar V_s(\bar N_{\pi,s})$. Then the vectors
  \begin{equation}\label{quotientstructurecoreq1c}
   \tau^i\theta^jW_{s,{\rm new}},\qquad i,j\geq0,\quad i+j=n-\bar N_{\pi,s},
  \end{equation}
  are a basis of $\bar V_s(n)$ for  $n\geq \bar N_{\pi,s}$.
\end{corollary}
\begin{proof}
Statements \eqref{quotientstructurecoreq1} and \eqref{quotientstructurecoreq1b} follow from Table \ref{dimensionstable}. We already noticed in \eqref{commreleq1} that $\theta\tau=\tau\theta$. The injectivity of $\tau:V_s(n)\to V_s(n+1)$ and $\tau:\bar V_s(n)\to\bar V_s(n+1)$ was proven in Lemma \ref{taunlevelraisinglemma}. The injectivity of $\theta:V_s(n)\to V_s(n+1)$ and of $\theta:\bar V_s(n)\to\bar V_s(n+1)$ follows easily from \eqref{W0theoremeq10}. The last statement also follows from \eqref{W0theoremeq10}.
\end{proof}

\section{Non-generic representations}
\label{nongenrepstructsec}
In this section we prove the analogue of Theorem \ref{W0theorem}
for non-generic representations\index{non-generic representation}\index{representation!non-generic}. 

\begin{theorem}
\label{nongenstructtheorem}
Let $(\pi,V)$ be an infinite-dimensional, non-generic, 
irreducible, admissible representation of $\GSp(4,F)$
with trivial central character. Assume that $\pi$ is
paramodular.
\begin{enumerate}
\item \label{nongenstructtheoremitem1}
The number $\bar N_{\pi,s}$ is defined if and only if 
$\pi$ belongs to subgroup IVb or IVc. Assume that
$\pi$ is such a representation. Then $N_{\pi,s}=\bar N_{\pi,s}=1$.
Let $v_0$ be an element of the one-dimensional space
$V_s(\bar N_{\pi,s})$ that is not in $V(\bar N_{\pi,s}-1) +V(\bar N_{\pi,s})$,
so that $v_0$ represents a non-zero element of $\bar V_s(\bar N_{\pi,s})$.
Then 
\begin{equation}
\label{nongenstructtheoremeq1}
V_s(n) = V(n-1) \oplus V(n) \oplus \C \tau^{n-\bar N_{\pi,s}}v_0
\end{equation}
for $n \geq N_{\pi,s}=\bar N_{\pi,s} =1$.
\item \label{nongenstructtheoremitem2}
The number $\bar N_{\pi,s}$ is not defined if and only if 
$\pi$ is a Saito-Kurokawa representation\index{Saito-Kurokawa representation}\index{representation!Saito-Kurokawa} or $\pi$ belongs to 
subgroup IIIb, Vd, or VId. Assume that $\pi$
is such a representation. 
Then 
\begin{equation}
\label{nongenstructtheoremeq2}
V_s(n) = V(n-1) \oplus V(n)
\end{equation}
for $n \geq N_{\pi,s}$
In \eqref{nongenstructtheoremeq2}, if $N_{\pi,s}=0$ and $n=0$,
then we take $V(n-1)=V(-1)$ to be the zero subspace.
\end{enumerate}
\end{theorem}
\begin{proof}
\ref{nongenstructtheoremitem1}. All of the assertions in \ref{nongenstructtheoremitem1}
follow from an inspection of Table \ref{dimensionstable}, except \eqref{nongenstructtheoremeq1}.
By Table \ref{dimensionstable} $\dim V_s(n) - \dim V(n-1) - \dim V(n)=1$ for integers $n$
such that $n \geq N_{\pi,s}$. Therefore, to verify \eqref{nongenstructtheoremeq1}, it suffices
to prove that $\tau^{n-\bar N_{\pi,s}} v_0$ is not contained in $V(n-1) \oplus V(n)$. This 
follows from \ref{taunlevelraisinglemmaitem2} of Lemma \ref{taunlevelraisinglemma}.

\ref{nongenstructtheoremitem2}. The assertions in \ref{nongenstructtheoremitem2} follow
from an inspection of Table \ref{dimensionstable}.
\end{proof}

\section{Summary statements}
\label{sumstatesec}
In this final section we provide statements that apply to all irreducible, admissible
representations of $\GSp(4,F)$ with trivial central character. These results are
corollaries of results from the previous sections. 

\begin{corollary}
\label{barNnondefcor}
Let $(\pi,V)$ be an  irreducible, admissible representation
of $\GSp(4,F)$ with central character. Assume that $\pi$ is paramodular. The number
$\bar N_{\pi,s}$ is not defined, so that $V(n-1)+V(n) = V_s(n)$ for 
all integers $n \geq N_{\pi,s}$, if and only if $\pi$ is  
a Saito-Kurokawa representation or $\pi$ belongs to subgroup IIIb, IVd, Vd, or VId.
\end{corollary}
\begin{proof}
Assume that $\bar N_{\pi,s}$ is not defined. 
By Sect.~\ref{structureSKsec} we see that $\pi$ is non-generic. Theorem \ref{nongenstructtheorem}
now implies that $\pi$ is  
a Saito-Kurokawa representation or $\pi$ belongs to subgroup IIIb, IVd, Vd, or VId (note that 
representations belonging to subgroup IVd are one-dimensional). The converse also follows from
Sect.~\ref{structureSKsec} and Theorem~\ref{nongenstructtheorem}.
\end{proof}

\begin{corollary}
\label{gendecompcor}
Let $(\pi,V)$ be an infinite-dimensional, irreducible, admissible representation
of $\GSp(4,F)$ with central character. Assume that $\pi$ is paramodular.  Define 
$v_0 \in V$ as follows. If $\bar N_{\pi,s}$ is not defined, set $v_0=0$; if $\bar N_{\pi,s}$
is defined, let $v_0$ be an element of $V_s(\bar N_{\pi,s})$ that is not contained
in $V(\bar N_{\pi,s}-1) + V(\bar N_{\pi,s})$. For integers $n$ such that $n \geq N_{\pi,s}$
let $E(n)$ be the subspace of $V_s(n)$ spanned by the vectors
$$
\tau^i \theta^j v_0, \qquad i,j \geq 0, \quad i+j = n-\bar N_{\pi,s}.
$$
Then 
$$
V_s(n) = V(n-1) \oplus V(n) \oplus E(n)
$$
for integers $n$ such that $n \geq N_{\pi,s}$. 
\end{corollary}
\begin{proof}
This follows from Theorem \ref{W0theorem} and Theorem \ref{nongenstructtheorem}.
\end{proof}

\begin{corollary}
\label{genquocor}
Let $(\pi,V)$ be an irreducible, admissible representation
of $\GSp(4,F)$ with trivial central character. 
Assume that $\pi$ is paramodular.
Assume that $\bar N_{\pi,s}$ is defined. Then 
\begin{equation}
\label{genquocoreq2}
\bar N_{\pi,s} = 
\begin{cases}
1&\text{if $N_\pi=0$ or $N_\pi=1$,}\\
N_\pi-1&\text{if $N_\pi \geq 2$ and $\pi$ is a \catone representation,}\\
N_\pi&\text{if $N_\pi \geq 2$ and $\pi$ is a \cattwo representation}
\end{cases}
\end{equation}
and
\begin{equation}
\label{genquocoreq3}
\dim \bar V_s(\bar N_{\pi,s}) =1.
\end{equation}
Let $v_{s,\mathrm{new}}$ be a non-zero element of the one-dimensional
vector space $\bar V_s(\bar N_{\pi,s})$. Then the vectors 
\begin{equation}
\label{genquocoreq1}
\tau^i \theta^j v_{s,\mathrm{new}}, \qquad i,j \geq 0, \quad i+j=n-\bar N_{\pi,s}
\end{equation}
span $\bar V_s(n)$ for $n \geq \bar N_{\pi,s}$.
\end{corollary}
\begin{proof}
By Corollary \ref{barNnondefcor}, the representation $\pi$ is not
a Saito-Kurokawa representation and  
$\pi$ does not belong  to subgroup IIIb, IVd, Vd, or VId; in particular, since $\pi$ does
not belong to subgroup IVd, $\pi$ is infinite-dimensional. 
If $\pi$ is generic, then 
\eqref{genquocoreq2}  follows from \eqref{quotientstructurecoreq1} and
Proposition \ref{shadowsigmavnewprop}. 
Assume that $\pi$ is non-generic. By \ref{nongenstructtheoremitem1}
of Theorem \ref{nongenstructtheorem}, $\pi$ belongs to subgroup IVb or IVc. If
$\pi$ belongs to subgroup IVb, then $N_\pi=2$, $N_{\pi,s}=1$, and $\bar N_{\pi,s}=1$ by
Table \ref{dimensionstable}; this verifies \eqref{genquocoreq2}.
If $\pi$ belongs to subgroup IVc, then $N_\pi=1$, $N_{\pi,s}=1$, and $\bar N_{\pi,s}=1$ by
Table \ref{dimensionstable}; this again verifies \eqref{genquocoreq2}.
Next, the assertion \eqref{genquocoreq3} follows from \eqref{quotientstructurecoreq1}
if $\pi$ is generic and Table \ref{dimensionstable} if $\pi$ is non-generic.
Finally, the vectors in \eqref{genquocoreq1} span $\bar V_s(n)$ for $n \geq \bar N_{\pi,s}$
by Corollary \ref{quotientstructurecor} if $\pi$ is generic, and by 
\ref{nongenstructtheoremitem1} of Theorem \ref{nongenstructtheorem} if $\pi$ is non-generic.
\end{proof}

%% file: SKMS_chapter8.tex
\chapter{Further Results about Generic Representations\index{generic representation}\index{representation!generic}}
\label{genericchap}
Let $(\pi,V)$ be a generic, irreducible, admissible
representation of $\GSp(4,F)$ with trivial central character.
In this chapter we prove several additional results about 
stable Klingen vectors \index{stable Klingen vector} in $\pi$. Our first result generalizes
a fundamental theorem from the paramodular theory. Assume
that $V$ is the Whittaker model of $\mathcal{W}(\pi,\psi_{c_1,c_2})$
of~$\pi$, and let $n$ be an integer such that $n \geq 0$.
In the first section, we prove that if $W \in V_s (n)$,
then $W\neq 0$ if and only if $W$ does not vanish on the diagonal
subgroup of $\GSp(4,F)$. 
To describe our second result, assume that $n \geq \max(N_{\pi,s},2)$,
and recall the 
surjective level lowering operator $\sigma_{n-1}:V_s(n) \to V_s(n-1)$
from Sect.~\ref{sigmasec}. We prove that 
$$
V_s(n) = (V(n-1) \oplus V(n)) + \ker(\sigma_{n-1})
$$
for integers $n \geq \max(N_{\pi,s},2)$. Consequently, $\sigma_{n-1}$ induces
an isomorphism
$$
(V(n-1) \oplus V(n))/\mathcal{K}_{n} \stackrel{\sim}{\longrightarrow} V_s(n-1)
$$
for integers $n \geq \max(N_{\pi,s},2)$, where $\mathcal{K}_n$ is the intersection
of $V(n-1) \oplus V(n)$ with $\ker(\sigma_{n-1})$. 
We characterize $\mathcal{K}_n$ for $n \geq \max(N_{\pi,s},2)$; in particular,
this subspace is at most two-dimensional. To prove these results we
use the nonvanishing statement from the first section. 
 Finally, in the last section,
we consider $\pi$ such that ${L(s,\pi)=1}$; this includes all $\pi$ that are supercuspidal.
Assuming that $L(s,\pi)=1$, we are able to define another, graphical, model
for $V_s(n)$ for integers $n \geq N_{\pi,s}$. The existence of this alternative
model also uses the result about the non-vanishing of stable Klingen vectors
on the diagonal of $\GSp(4,F)$. In this model, our level changing operators have simple
and visual interpretations. Finally, still under the hypothesis that $L(s,\pi)=1$,
we prove that $N_\pi \geq 4$. 

Before beginning, 
we remind the reader that $\lambda_\pi$, $\mu_\pi$, and $\varepsilon_\pi$ are the 
eigenvalues of $T_{0,1}$, $T_{1,0}$, and $u_{N_\pi}$, respectively, on the 
one-dimensional space $V(N_\pi)$. Also, by Theorem \ref{linindparatheorem} the subspace $V(n-1)+V(n)$ of $V_s(n)$
is a direct sum $V(n-1) \oplus V(n)$ for $n \geq 1$. Finally, we will use the 
equivalences of Proposition \ref{shadowsigmavnewprop}.

\section{Non-vanishing on the diagonal}

Let $(\pi,V)$ be a generic, irreducible, admissible representation of $\GSp(4,F)$
with trivial central character, and let $n$ be an integer such that $n \geq 0$.
Let $W \in V_s(n)$. In this section we prove that $W \neq 0$ if and only if $W$
does not vanish on the diagonal subgroup of $\GSp(4,F)$. This result generalizes
Corollary 4.3.8 of \cite{NF}. We begin by proving this statement for the elements
of the subspace $V(n-1)\oplus V(n)$ of $V_s(n)$. 

\begin{proposition}\label{paraetaprop}
Let $(\pi,V)$ be a generic, irreducible, admissible representation of $\GSp(4,F)$ with trivial central character,
and let $n$ be an integer such that $n \geq 0$. Let $V=\mathcal{W}(\pi,\psi_{c_1,c_2})$. Assume that $W\in V(n-1)\oplus V(n)$ is such that $W(\Delta_{i,j})=0$ for all $i,j\in\Z$; if $n=0$, then we define $V(-1)=0$. Then $W=0$.
\end{proposition}
\begin{proof}
The statement is trivial if $n\leq N_\pi-1$. For $n\geq N_\pi$ we use induction on~$n$. The case $n=N_\pi$ follows from Corollary 4.3.8 of \cite{NF}. Assume that $n>N_\pi$, and that the statement holds for $n-1$. Let $W_1\in V(n-1)$ and $W_2\in V(n)$, and assume that
\begin{equation}\label{paraetapropeq1}
 (W_1+W_2)(\Delta_{i,j})=0
\end{equation}
for all $i,j\in\Z$; we need to prove that $W_1=W_2=0$. Define $W'=q^{-1} \theta' W_1+W_2$. Then $W' \in V(n)$ and:
\begin{align*}
Z(s,W')&= q^{-1} Z(s,\theta' W_1) + Z(s, W_2)\\
&=Z(s,W_1)+Z(s,W_2)\qquad \text{(by \eqref{zetanewspanlemmaeq31})}\\
&=Z(s,W_1+W_2)\\
&=0 \qquad \text{(by \eqref{paraetapropeq1})}.
\end{align*}
Assume that $n\geq2$. By the $\eta$-principle\index{eta@$\eta$-principle}, Theorem 4.3.7 of \cite{NF}, there exists $W_3\in V(n-2)$ such that
\begin{equation}\label{paraetapropeq2}
W'=\eta W_3.
\end{equation}
If $n=1$, then we set $W_3=0$, and equation \eqref{paraetapropeq2} still holds by Theorem 4.3.7 of \cite{NF}. Evaluating at $\Delta_{i,j}$, we get by \eqref{thetap1defeq}
$$
 q^{-1}W_1(\Delta_{i,j}\eta)+W_1(\Delta_{i,j})+W_2(\Delta_{i,j})=W_3(\Delta_{i,j}\eta).
$$
By \eqref{paraetapropeq1},
$$
 q^{-1}W_1(\Delta_{i,j}\eta)=W_3(\Delta_{i,j}\eta),
$$
and hence
$$
 q^{-1}W_1(\Delta_{i,j})=W_3(\Delta_{i,j})
$$
for all $i,j\in\Z$. This implies that $(W_3-q^{-1}W_1)(\Delta_{i,j})=0$ for all $i,j \in \Z$.
Since $W_3 - q^{-1} W_1 \in V(n-2)+V(n-1)$ we have 
$W_1=W_3=0$ by the induction hypothesis. Then also $W_2=0$. This concludes the proof.
\end{proof}

\begin{proposition}\label{parazerodiagprop}
Let $(\pi,V)$ be a generic, irreducible, admissible representation of $\GSp(4,F)$ with trivial central character, and 
let $n$ be an integer such that $n \geq 0$. Let $V=\mathcal{W}(\pi,\psi_{c_1,c_2})$. Assume that $W\in V_s(n)$ is such that $W(\Delta_{i,j})=0$ for all $i,j\in\Z$. Then $W=0$.
\end{proposition}
\begin{proof}
Write $W=W_1+W_2$ with $W_1\in V(n)$ and $W_2\in V_s(n)$ such that $p_n(W_2)=0$ (see \eqref{pnVsneq4}. By \eqref{rhopoplemmaeq3},
\begin{equation}\label{parazerodiaglemmaeq1}
 \rho'_nW_2=q^{-1}\eta W_2-q\tau_{n+1}W_2\in V(n+1).
\end{equation}
Hence
\begin{align}\label{parazerodiaglemmaeq5}
 (\rho'_nW_2)(\Delta_{i,j})&=q^{-1}(\eta W_2)(\Delta_{i,j})-q(\tau_{n+1}W_2)(\Delta_{i,j})\nonumber\\
 &=q^{-1}W_2(\Delta_{i,j}\eta)-qW_2(\Delta_{i,j})\nonumber\\
 &=-q^{-1}W_1(\Delta_{i,j}\eta)+qW_1(\Delta_{i,j}).
\end{align}
It follows that
\begin{align*}
Z(s,\rho'_nW_2-\theta'W_1)
&=-q^{-1}Z(s,\eta W_1)+qZ(s,W_1)-Z(s,\theta' W_1)\\
&=0\qquad\text{(by \eqref{zetanewspanlemmaeq31} and \eqref{zetanewspanlemmaeq32})}.
\end{align*}
The vector $\rho'_nW_2-\theta'W_1$ is in $V(n+1)$.
If $n\geq1$, then by the $\eta$-principle, Theorem~4.3.7 of \cite{NF}, there exists $W_3\in V(n-1)$ such that
\begin{equation}\label{parazerodiaglemmaeq3}
 \rho'_nW_2-\theta'W_1=\eta W_3.
\end{equation}
If $n=0$, then we set $W_3=0$, and \eqref{parazerodiaglemmaeq3} still holds by Theorem 4.3.7 of \cite{NF}. Hence, by \eqref{thetap1defeq},
\begin{align}\label{parazerodiaglemmaeq4}
 (\rho'_nW_2)(\Delta_{i,j})&=(\theta'W_1)(\Delta_{i,j})+(\eta W_3)(\Delta_{i,j})\nonumber\\
 &=W_1(\Delta_{i,j}\eta)+qW_1(\Delta_{i,j})+W_3(\Delta_{i,j}\eta)
\end{align}
for all $i,j \in \Z$.
Substituting into \eqref{parazerodiaglemmaeq5}, we get
\begin{align*}
 W_1(\Delta_{i,j}\eta)+qW_1(\Delta_{i,j})+W_3(\Delta_{i,j}\eta)&=-q^{-1}W_1(\Delta_{i,j}\eta)+qW_1(\Delta_{i,j}),\\
 W_1(\Delta_{i,j}\eta)+W_3(\Delta_{i,j}\eta)&=-q^{-1}W_1(\Delta_{i,j}\eta),\\
 (1+q^{-1})W_1(\Delta_{i,j}\eta)+W_3(\Delta_{i,j}\eta)&=0
\end{align*}
for all $i,j \in \Z$.
It follows that
$$
(1+q^{-1})W_1(\Delta_{i,j})+W_3(\Delta_{i,j})=0
$$
for all $i,j\in\Z$. By Proposition \ref{paraetaprop} and Theorem \ref{linindparatheorem} we obtain $W_1=W_3=0$. Hence $W=W_2$, so that $p_n(W)=0$. By \eqref{parazerodiaglemmaeq3}, $\rho'_nW=0$. By \eqref{rhopoplemmaeq3},
\begin{equation}\label{parazerodiaglemmaeq7}
 0=\rho'_nW=q^{-1}\eta W-q\tau_{n+1}W.
\end{equation}
Hence $\tau_{n+1}W=q^{-2}\eta W$. Therefore $W=0$ by Proposition \ref{funnyidentityprop}.  This completes the proof.
\end{proof}
\section{The kernel of a level lowering operator}
\label{sigmakernelsec}
Let $(\pi,V)$ be a generic, irreducible, admissible representation of $\GSp(4,F)$
with trivial central character, and let $n$ be an integer such that $n \geq \max(N_{\pi,s},2)$.
In Theorem \ref{W0theorem} we proved the structural result \eqref{W0theoremeq10} for the spaces $V_s(n)$.
In this section we will prove that also $V_s(n)=(V(n-1)\oplus V(n))+\ker(\sigma_{n-1})$, 
where $\sigma_{n-1}:V_s(n)\to V_s(n-1)$ is the surjective level lowering operator 
defined in Sect.~\ref{sigmasec}. We will precisely determine the 
intersection $(V(n-1)\oplus V(n))\cap\ker(\sigma_{n-1})$; it is at most two-dimensional.

\begin{lemma}\label{operatorWDeltalemma}
Let $(\pi,V)$ be a generic, irreducible, admissible representation of $\GSp(4,F)$ with trivial central character,
and let $n$ be an integer such that $n\geq 0$. Let $V=\mathcal{W}(\psi_{c_1,c_2})$. Assume that $i,j \in \Z$ with $i,j \geq 0$. Then, for $W\in V_s(n)$, 
\begin{align*}
  (\sigma_{n+1}\theta W)(\Delta_{i,j})&=
    W(\Delta_{i+1,j-1})+qW(\Delta_{i,j+1}),\\
  (\sigma_{n+1}\eta W)(\Delta_{i,j})&=W(\Delta_{i,j}),\\
  (\sigma_{n+1}\theta^2W)(\Delta_{i,j})&=
   \begin{cases}
    W(\Delta_{i+1,j-2})+2qW(\Delta_{i,j})+q^2W(\Delta_{i-1,j+2})&\text{if }j\geq1,\\
    qW(\Delta_{i,0})+q^2W(\Delta_{i-1,2})&\text{if }j=0.
   \end{cases}
 \end{align*}
 Furthermore, for $W\in V(n)$,
 \begin{align*}
  (\sigma_{n+1}\theta' W)(\Delta_{i,j})&=W(\Delta_{i,j})+qW(\Delta_{i+1,j}),\\
  (\sigma_{n+1}\theta\theta'W)(\Delta_{i,j})&=W(\Delta_{i,j-1})+qW(\Delta_{i-1,j+1})+qW(\Delta_{i+1,j-1})+q^2W(\Delta_{i,j+1}),\\
  (\sigma_{n+1}{\theta'}^2W)(\Delta_{i,j})&=W(\Delta_{i-1,j})+2qW(\Delta_{i,j})+q^2W(\Delta_{i+1,j}).
 \end{align*}
 Here, $\Delta_{i,j}$ is the diagonal matrix defined in \eqref{Deltaijdefeq}.
\end{lemma}
\begin{proof}
By \eqref{thetadefeq3}, for $W\in V_s(n)$,
\begin{equation}\label{thetaWdiageq}
 (\theta W)(\Delta_{i,j})=
  \begin{cases}
   W(\Delta_{i,j-1})+qW(\Delta_{i-1,j+1})&\text{if }j\geq0,\\
   0&\text{if }j<0,
  \end{cases}
\end{equation}
and
\begin{align}\label{theta2Wdiageq}
 (\theta^2 W)(\Delta_{i,j})
  &=
  \begin{cases}
   W(\Delta_{i,j-2})+2qW(\Delta_{i-1,j})+q^2W(\Delta_{i-2,j+2})&\text{if }j\geq1,\\
   qW(\Delta_{i-1,0})+q^2W(\Delta_{i-2,2})&\text{if }j=0,\\
   0&\text{if }j<0.
  \end{cases}
\end{align}
By \eqref{thetap1defeq}, for $W\in V(n)$,
\begin{equation}\label{thetapWdiageq}
 (\theta' W)(\Delta_{i,j})=W(\Delta_{i-1,j})+qW(\Delta_{i,j})
\end{equation}
and
\begin{equation}\label{thetap2Wdiageq}
 ({\theta'}^2 W)(\Delta_{i,j})=W(\Delta_{i-2,j})+2qW(\Delta_{i-1,j})+q^2W(\Delta_{i,j}).
\end{equation}
By \eqref{sigmaopseq}, for $n\geq1$ and $W\in V_s(n+1)$,
\begin{equation}\label{sigmaWdiageq}
 (\sigma_n W)(\Delta_{i,j})=
  \begin{cases}
   W(\Delta_{i+1,j})&\text{if }i\geq0,\\
   0&\text{if }i<0.
  \end{cases}
\end{equation}
Our asserted formulas follow easily from \eqref{thetaWdiageq}--\eqref{sigmaWdiageq}.
\end{proof}

\begin{lemma}\label{typeIlevel2kersigmalemma}
 Let $(\pi,V)$ be a generic, irreducible, admissible representation of $\GSp(4,F)$ with trivial central character and $N_\pi=0$. Let $W_{\rm new}\in V(0)$ be non-zero. Then the subspace
 $$
  (V(1)+V(2))\cap\ker(\sigma_1)
 $$
 of $V_s(2)$ is two-dimensional, and is spanned by the vectors
 \begin{align}
  \label{typeIlevel2kersigmalemmaeq1}W_1&=((q^2-q)\theta+\theta\theta'-\lambda_\pi\eta)W_{\rm new},\\
  \label{typeIlevel2kersigmalemmaeq2}W_2&=(q(q^2-1)\theta'+q\theta^2+{\theta'}^2-(\mu_\pi+q^3+q^2+q+1)\eta)W_{\rm new}.
 \end{align}
 The map
 \begin{equation}\label{typeIlevel2kersigmalemmaeq3}
  \sigma_1:\:V(1)+V(2)\longrightarrow V_s(1)
 \end{equation}
 is surjective.
\end{lemma}
\begin{proof}
We have 
$$
V(1)=\langle\theta W_{\rm new},\theta'W_{\rm new}\rangle \quad \text{and}\quad 
V(2)=\langle\theta^2 W_{\rm new},\theta\theta'W_{\rm new},{\theta'}^2 W_{\rm new},\eta W_{\rm new}\rangle
$$
by Theorem 7.5.6 of \cite{NF}. Evidently, $W_1,W_2\in V(1)+ V(2)\subset  V_s(2)$. 
We will show that $\sigma_1W_1=\sigma_1W_2=0$. As in \eqref{cijdefeq}, let $c_{i,j}=W_{\rm new}(\Delta_{i,j})$ for $i,j \in \Z$. 
By Lemma~\ref{operatorWDeltalemma}, for $i,j \in \Z$ with $i,j\geq0$,
\begin{align*}
 (\sigma_1\theta W_{\rm new})(\Delta_{i,j})&=c_{i+1,j-1}+qc_{i,j+1},\\
 (\sigma_1\theta' W_{\rm new})(\Delta_{i,j})&=c_{i,j}+qc_{i+1,j},\\
 (\sigma_1\eta W_{\rm new})(\Delta_{i,j})&=c_{i,j},\\
 (\sigma_1\theta^2W_{\rm new})(\Delta_{i,j})&=
  \begin{cases}
   c_{i+1,j-2}+2qc_{i,j}+q^2c_{i-1,j+2}&\text{if }j\geq1,\\
   qc_{i,0}+q^2c_{i-1,2}&\text{if }j=0,
  \end{cases}\\
 (\sigma_1\theta\theta'W_{\rm new})(\Delta_{i,j})&=c_{i,j-1}+qc_{i-1,j+1}+qc_{i+1,j-1}+q^2c_{i,j+1},\\
 (\sigma_1{\theta'}^2W_{\rm new})(\Delta_{i,j})&=c_{i-1,j}+2qc_{i,j}+q^2c_{i+1,j}.
\end{align*}
Using Lemma 7.1.2 of \cite{NF}, we get for $i,j \in \Z$ with $i,j \geq 0$,
\begin{align}
 (\sigma_1\theta W_{\rm new})(\Delta_{i,j})&=q^{-2}\lambda_\pi  c_{i,j}-q^{-1}c_{i-1,j+1}-q^{-2}c_{i,j-1},\nonumber\\
 (\sigma_1\theta' W_{\rm new})(\Delta_{i,j})&=c_{i,j}+qc_{i+1,j},\label{typeIlevel2kersigmalemmaeq21}\\
 (\sigma_1\eta W_{\rm new})(\Delta_{i,j})&=c_{i,j},\nonumber\\
 (\sigma_1\theta^2W_{\rm new})(\Delta_{i,j})&=(q^{-1}\mu_\pi+q+q^{-1})c_{i,j}-q^3c_{i+1,j}-q^{-1}c_{i-1,j},\nonumber\\
 (\sigma_1\theta\theta'W_{\rm new})(\Delta_{i,j})&=q^{-1}\lambda_\pi c_{i,j}+(1-q^{-1})(c_{i,j-1}+qc_{i-1,j+1}),\nonumber\\
 (\sigma_1{\theta'}^2W_{\rm new})(\Delta_{i,j})&=2qc_{i,j}+q^2c_{i+1,j}+c_{i-1,j}.\nonumber
\end{align}
It follows that $(\sigma_1W_1)(\Delta_{i,j})=(\sigma_1W_2)(\Delta_{i,j})=0$ for $i,j \in \Z$ with $i,j \geq 0$. By \eqref{delltavanisheq} and Proposition \ref{parazerodiagprop}, $\sigma_1W_1=\sigma_1W_2=0$.
So far we have proven that $\langle W_1,W_2 \rangle \subset (V(1)+V(2))\cap \ker(\sigma_1)$. From above, and by Table \ref{dimensionstable} we have $\dim (V(1)+V(2))=6$ and $\dim V_s(1)=4$; also, by \eqref{sigmasurjectivepropeq2}
of Proposition \ref{sigmasurjectiveprop}, the three-dimensional subspace $V(0)+V(1)$ of $V_s(1)$ is contained
in $\sigma_1(V(1)+V(2))$. Thus, to prove that $\sigma_1:V(1)+V(2) \to V_s(1)$ is surjective and
that $\langle W_1,W_2 \rangle = (V(1)+V(2)) \cap \ker(\sigma_1)$, it will suffice to prove that, say, $\sigma_1 \theta' W_{\mathrm{new}} \notin V(0)+V(1)$. We have
\begin{align*}
&Z(s,\sigma_1 \theta' W_{\mathrm{new}})
=(1-q^{-1})\sum_{j=0}^\infty (\sigma_1 \theta' W_{\mathrm{new}})(
\begin{bsmallmatrix}
\varpi^j&&&\\
&\varpi^j&&\\
&&1&\\
&&&1
\end{bsmallmatrix}) (q^{-s+\frac{3}{2}})^j \qquad \text{(by \eqref{Vszetaintposdef})}\\
&\qquad=(1-q^{-1})\sum_{j=0}^\infty 
(c_{0,j}+q c_{1,j})
 (q^{-s+\frac{3}{2}})^j \qquad \text{(by \eqref{typeIlevel2kersigmalemmaeq21})}\\
&\qquad=(1+
q^{-3}\mu_\pi+q^{-3}-\lambda_\pi q^{-\frac{5}{2}}q^{-s}+q^{-1} q^{-2s}) Z(s,W_{\mathrm{new}})
\qquad \text{(by \eqref{c1jZsWneweq})}.
\end{align*}
On the other hand, $Z(s,V(0)+V(1)) = (\C + \C q^{-s}) Z(s, W_{\mathrm{new}})$ since $V(0)+V(1)$ is spanned
by $W_{\mathrm{new}}, \theta W_{\mathrm{new}}$, and $\theta' W_{\mathrm{new}}$ and since
\eqref{zetanewspanlemmaeq31} holds. It follows that 
$\sigma_1 \theta' W_{\mathrm{new}} \notin V(0)+V(1)$, as desired.
\end{proof}

\begin{lemma}\label{typeIIalevel2kersigmalemma}
 Let $(\pi,V)$ be a generic, irreducible, admissible representation of 
 $\GSp(4,F)$ with trivial central character and $N_\pi=1$. Let $W_{\rm new}\in V(1)$ be non-zero. Then the subspace
 $$
  (V(1)+V(2))\cap\ker(\sigma_1)
 $$
 of $V_s(2)$ is one-dimensional, and is spanned by the vector
 \begin{equation}\label{typeIIalevel2kersigmalemmaeq1}
  W_3=(q(q^2-1)+\varepsilon_\pi q\theta+\theta')W_{\rm new}.
 \end{equation}
 The map
 \begin{equation}\label{typeIIalevel2kersigmalemmaeq3}
  \sigma_1:\:V(1)+V(2)\longrightarrow V_s(1)
 \end{equation}
 is surjective.
\end{lemma}
\begin{proof}
We have $V(1)=\langle W_{\rm new}\rangle$ and $V(2)=\langle\theta W_{\rm new},\theta'W_{\rm new}\rangle$ by Theorem~7.5.6 of \cite{NF}. Evidently, $W_3\in V(1)+ V(2)$. We will show that $\sigma_1W_3=0$. As in \eqref{cijdefeq}, let $c_{i,j}=W_{\rm new}(\Delta_{i,j})$ for $i,j \in \Z$. By \eqref{thetaWdiageq}, \eqref{thetapWdiageq} and \eqref{sigmaWdiageq}, for $i,j \in \Z$ with $i,j\geq0$,
\begin{align*}
 (\sigma_1 W_{\rm new})(\Delta_{i,j})&=c_{i+1,j},\\
 (\sigma_1\theta W_{\rm new})(\Delta_{i,j})&=c_{i+1,j-1}+qc_{i,j+1},\\
 (\sigma_1\theta' W_{\rm new})(\Delta_{i,j})&=c_{i,j}+qc_{i+1,j}.
\end{align*}
By Lemmas 7.2.1 and 7.2.2 of \cite{NF},
\begin{align}
 (\sigma_1 W_{\rm new})(\Delta_{i,j})&=c_{i+1,j},\nonumber \\
 (\sigma_1\theta W_{\rm new})(\Delta_{i,j})&=-\varepsilon q^{-1}(c_{i,j}+q^3c_{i+1,j}),\nonumber\\
 (\sigma_1\theta' W_{\rm new})(\Delta_{i,j})&=c_{i,j}+qc_{i+1,j} \label{sigma1thetapeq}
\end{align}
for $i,j \in \Z$ with $i,j \geq 0$. 
It follows that $(\sigma_1 W_3)(\Delta_{i,j})=0$ for $i,j \in \Z$ with $i,j \geq 0$. By \eqref{delltavanisheq} and Proposition \ref{parazerodiagprop}, $\sigma_1W_3=0$.
So far we have proven that $\langle W_3\rangle\subset (V(1)+V(2))\cap\ker(\sigma_1)$. 
From above and by Table \ref{dimensionstable} we have $\dim (V(1)+V(2))=3$
and $\dim V_s(1)=2$; also, by \eqref{sigmasurjectivepropeq2}
of Proposition \ref{sigmasurjectiveprop}, the one-dimensional subspace $V(1)$ of $V_s(1)$ is contained
in $\sigma_1(V(1)+V(2))$. Thus, to prove that $\sigma_1:V(1)+V(2) \to V_s(1)$ is surjective and
that $\langle W_3 \rangle = (V(1)+V(2)) \cap \ker(\sigma_1)$, it will suffice to prove that, say, $\sigma_1 \theta' W_{\mathrm{new}} \notin V(1)$. We have
\begin{align*}
Z(s,\sigma_1 \theta' W_{\mathrm{new}})
&=(1-q^{-1})\sum_{j=0}^\infty (\sigma_1 \theta' W_{\mathrm{new}})(
\begin{bsmallmatrix}
\varpi^j&&&\\
&\varpi^j&&\\
&&1&\\
&&&1
\end{bsmallmatrix}) (q^{-s+\frac{3}{2}})^j \qquad \text{(by \eqref{Vszetaintposdef})}\\
&=(1-q^{-1})\sum_{j=0}^\infty 
(c_{0,j}+q c_{1,j})
 (q^{-s+\frac{3}{2}})^j \qquad \text{(by \eqref{sigma1thetapeq})}\\
&=(1+
q^{-3}\mu_\pi+\varepsilon_\pi q^{-\frac{3}{2}}q^{-s}) Z(s,W_{\mathrm{new}})
\qquad \text{(by \eqref{c1jZsWneweq})}.
\end{align*}
On the other hand, $Z(s,V(1)) = \C  Z(s, W_{\mathrm{new}})$ since $V(1)$ is spanned
by $W_{\mathrm{new}}$. It follows that 
$\sigma_1 \theta' W_{\mathrm{new}} \notin V(1)$, as desired.
\end{proof}

\begin{lemma}\label{sigmakerneldimensionlemma}
Let $\pi$ be a generic, irreducible, admissible representation of the group $\GSp(4,F)$ with trivial central character, and
let $n$ be an integer such that 
$n\geq \max(N_{\pi,s},2)=\max(\bar N_{\pi,s},2)$. We consider the homomorphism 
$$
\sigma_{n-1}:V_s(n)\to V_s(n-1).
$$
 \begin{enumerate}
  \item \label{sigmakerneldimensionlemmaitem1} We have
   \begin{equation}\label{sigmakerneldimensionlemmaeq1}
    \dim\ker(\sigma_{n-1})=n-N_\pi+2.
   \end{equation}
  \item \label{sigmakerneldimensionlemmaitem2} The map 
   \begin{equation}\label{sigmakerneldimensionlemmaeq2}
    \ker(\sigma_{n-1})\stackrel{\sim}{\longrightarrow}Z(s,W_{\rm new})(\C+\C q^{-s}+\ldots+\C(q^{-s})^{n-N_\pi+1})
   \end{equation}
   defined by $W \mapsto Z(s,W)$ is an isomorphism of vector spaces. 
 \end{enumerate}
\end{lemma}
\begin{proof}
\ref{sigmakerneldimensionlemmaitem1}. By Proposition \ref{sigmasurjectiveprop}, the map $\sigma_{n-1}:V_s(n)\to V_s(n-1)$ is surjective. Hence $\dim\ker(\sigma_{n-1})=\dim V_s(n)-\dim V_s(n-1)$. Now \eqref{sigmakerneldimensionlemmaeq1} follows from \eqref{W0theoremeq15}.

\ref{sigmakerneldimensionlemmaitem2}. 
By \eqref{W0theoremeq20}, the map \eqref{sigmakerneldimensionlemmaeq2} is well-defined. By \eqref{sigmakerneldimensionlemmaeq1}, it suffices to show that the map is injective. Hence assume that $W\in\ker(\sigma_{n-1})$ and $Z(s,W)=0$. The latter condition implies that
\begin{equation}\label{sigmakerneldimensionlemmaeq6}
 W(\Delta_{0,j})=0\qquad\text{for all $j \in \Z$ with $j \geq 0$}.
\end{equation}
It follows from \eqref{sigmaopseq} that
\begin{equation}\label{sigmakerneldimensionlemmaeq7}
 0=(\sigma_{n-1}W)(\Delta_{i,j})=W(\Delta_{i+1,j})\qquad\text{for all $i,j \in \Z$ with $i,j \geq 0$}.
\end{equation}
Hence $W(\Delta_{i,j})=0$ for all $i,j \in \Z$ with $i,j\geq0$, and then also for all $i,j\in\Z$ by~\eqref{delltavanisheq}. By Proposition \ref{parazerodiagprop} it follows that $W=0$. This concludes the proof.
\end{proof}

\begin{lemma}\label{sigmakerlemma}
 Let $\pi$ be a generic, irreducible, admissible representation of the group 
$\GSp(4,F)$ with trivial central character. Assume that $N_\pi\geq2$. Let $W_0$ be as in Theorem \ref{W0theorem}. For any $n\geq\max(\bar N_{\pi,s},2)=\max(N_{\pi,s},2)$, and any $j\in\{0,\ldots,n-\bar N_{\pi,s}\}$, there exists a vector of the form
 \begin{equation}\label{sigmakerlemmaeq1}
  W_1+W_2+\tau^i\theta^jW_0,\qquad W_1\in V(n-1),\;W_2\in V(n),\;i+j=n-\bar N_{\pi,s},
 \end{equation}
 which is in the kernel of $\sigma_{n-1}:V_s(n)\to V_s(n-1)$.
\end{lemma}
\begin{proof}
We will use induction on $n$, distinguishing cases for the beginning of the induction.

Assume that $N_\pi=0$. Then $\bar N_{\pi,s}=1$ and the beginning of the induction is at $n=2$. By Theorem \ref{W0theorem},
$$
 V_s(2)=V(1)\oplus V(2)\oplus\C\tau W_0\oplus\C\theta W_0.
$$
By Lemma \ref{typeIlevel2kersigmalemma} and its proof, $\dim((V(1)\oplus V(2))\cap\ker(\sigma_1))=2$, $\dim V_s(1)=4$, and the map $\sigma_1:\:V(1)\oplus V(2)\longrightarrow V_s(1)$ is surjective. It follows that the map
$$
 \sigma_1:\:V(1)\oplus V(2)\oplus\C\tau W_0\longrightarrow V_s(1)
$$
is surjective and $\dim((V(1)\oplus V(2)\oplus\C\tau W_0)\cap\ker(\sigma_1))=3$. In particular, there exists an element in $(V(1)\oplus V(2)\oplus\C\tau W_0)\cap\ker(\sigma_1)$ with non-trivial $\tau W_0$ component. Similarly we see that there exists an element in $(V(1)\oplus V(2)\oplus\C\theta W_0)\cap\ker(\sigma_1)$ with non-trivial $\theta W_0$ component. This proves the assertion for $N_\pi=0$ and $n=2$.

Assume that $N_\pi=1$. Then $\bar N_{\pi,s}=1$ and the beginning of the induction is at $n=2$. By Theorem \ref{W0theorem},
$$
 V_s(2)=V(1)\oplus V(2)\oplus\C\tau W_0\oplus\C\theta W_0.
$$
By Lemma \ref{typeIIalevel2kersigmalemma} and its proof, $\dim((V(1)\oplus V(2))\cap\ker(\sigma_1))=1$, $\dim V_s(1)=2$, and the map $\sigma_1:\:V(1)\oplus V(2)\longrightarrow V_s(1)$ is surjective. It follows that the map
$$
 \sigma_1:\:V(1)\oplus V(2)\oplus\C\tau W_0\longrightarrow V_s(1)
$$
is surjective and $\dim((V(1)\oplus V(2)\oplus\C\tau W_0)\cap\ker(\sigma_1))=2$. In particular, there exists an element in $(V(1)\oplus V(2)\oplus\C\tau W_0)\cap\ker(\sigma_1)$ with non-trivial $\tau W_0$ component. Similarly we see that there exists an element in $(V(1)\oplus V(2)\oplus\C\theta W_0)\cap\ker(\sigma_1)$ with non-trivial $\theta W_0$ component. This proves the assertion for $N_\pi=1$ and $n=2$.

Assume that $N_\pi=2$ and $\mu_\pi \neq0$. Then $\bar N_{\pi,s}=N_{\pi,s}=1$ and the beginning of the induction is at $n=2$. We have $\dim V_s(1)=1$ and $\dim V_s(2)=3$. By Theorem~\ref{W0theorem},
$$
 V_s(2)=V(2)\oplus\C\tau W_0\oplus\C\theta W_0.
$$
By Proposition \ref{shadowsigmavnewprop}, $\sigma_1:V(2)\to V_s(1)$ is an isomorphism. It follows that $\dim((V(2)\oplus\C\tau W_0)\cap\ker(\sigma_1))=\dim((V(2)\oplus\C\theta W_0)\cap\ker(\sigma_1))=1$. This proves the assertion for $N_\pi=2$, $\mu_\pi \neq0$ and $n=2$.

Assume that $N_\pi\geq3$ and $\mu_\pi \neq0$. Then $\bar N_{\pi,s}=N_{\pi,s}=N_\pi-1$ and the beginning of the induction is at $n=N_\pi-1$. The assertion is true for $n=N_\pi-1$, since $\sigma_{(N_\pi-1)-1}(W_0)=0$. 

Assume that $N_\pi\geq2$ and $\mu_\pi=0$. Then $\bar N_{\pi,s}=N_{\pi,s}=N_\pi$ and the beginning of the induction is at $n=N_\pi$. The assertion is true for $n=N_\pi$, since $\sigma_{N_\pi-1}(W_0)=0$.

We completed the beginning of the induction and proceed to the induction step.
Assume that $n>\max(\bar N_{\pi,s},2)$, and that the assertion has been proven for $n-1$. Hence, for any $j\in\{0,\ldots,n-1-\bar N_{\pi,s}\}$, there exists a vector of the form
\begin{equation}\label{sigmakerlemmaeq2}
 W_1+W_2+\tau^i\theta^jW_0,\qquad W_1\in V(n-2),\;W_2\in V(n-1),\;i+j=n-1-\bar N_{\pi,s},
\end{equation}
which is in the kernel of $\sigma_{n-2}:V_s(n-1)\to V_s(n-2)$; then the vectors
\begin{equation}\label{sigmakerlemmaeq3}
 \tau(W_1+W_2+\tau^i\theta^jW_0)=\tau(W_1+W_2)+\tau^{i+1}\theta^jW_0,\qquad i+1+j=n-\bar N_{\pi,s},
\end{equation}
lie in the kernel of $\sigma_{n-1}:V_s(n)\to V_s(n-1)$ by \ref{sigmaopslemmaitem6} of Lemma \ref{sigmaopslemma}, and also $\tau(W_1+W_2)\in V(n-1)\oplus V(n)$ by the proof of Lemma \ref{taunlevelraisinglemma}. It remains to find a vector of the form
\begin{equation}\label{sigmakerlemmaeq4}
 W_1'+W_2'+\theta^{n-\bar N_{\pi,s}}W_0,\qquad W_1'\in V(n-1),\;W_2'\in V(n),
\end{equation}
which is in the kernel of $\sigma_{n-1}:V_s(n)\to V_s(n-1)$. Suppose that no such vector would exist. Then, by Theorem \ref{W0theorem}, every element of ${\rm ker}(\sigma_{n-1})$ would be of the form
\begin{equation}\label{sigmakerlemmaeq5}
 W_1'+W_2'+\sum_{\substack{i+j=n-\bar N_{\pi,s}\\j<n-\bar N_{\pi,s}}}c_i\tau^i\theta^j W_0,\qquad W_1'\in V(n-1),\;W_2'\in V(n),\;c_i\in\C.
\end{equation}
By Theorem \ref{W0theorem}, in particular \eqref{W0theoremeq3}, it would follow that
\begin{equation}\label{sigmakerlemmaeq6}
 Z(s,\ker(\sigma_{n-1}))\subset Z(s,W_{\rm new})(\C+\C q^{-s}+\ldots+\C(q^{-s})^{n-N_\pi}).
\end{equation}
This contradicts \eqref{sigmakerneldimensionlemmaeq2}.
\end{proof}

\begin{theorem}\label{genericsigmakertheorem}
 Let $\pi$ be a generic\index{generic representation}\index{representation!generic}, irreducible, admissible representation of the group 
 $\GSp(4,F)$ with trivial central character. Let
$n$ be an integer such that $n\geq\max(N_{\pi,s},2)=\max(\bar N_{\pi,s},2)$.
 \begin{enumerate}
  \item \label{genericsigmakertheoremitem1} We have
   \begin{equation}\label{genericsigmakertheoremeq1}
    V_s(n)=(V(n-1)\oplus V(n))+\ker(\sigma_{n-1}).
   \end{equation}
  \item \label{genericsigmakertheoremitem2} Let $\mathcal{K}_n=(V(n-1)\oplus V(n))\,\cap\,\ker(\sigma_{n-1})$. Then
   \begin{equation}\label{genericsigmakertheoremeq2}
    \dim(\mathcal{K}_n)=
           \begin{cases}
            2&\text{if }N_\pi=0,\\
            1&\text{if }N_\pi=1,\\
            1&\text{if $N_\pi\geq2$ and $\mu_\pi=0$},\\
            0&\text{if $N_\pi\geq2$ and $\mu_\pi\neq0$}.
           \end{cases}
   \end{equation}
   If $N_\pi=0$, then $\mathcal{K}_n=\langle\tau_{n-1}W_1,\tau_{n-1}W_2\rangle$, where $W_1,W_2$ are the vectors defined in \eqref{typeIlevel2kersigmalemmaeq1}, \eqref{typeIlevel2kersigmalemmaeq2}; if $N_\pi=1$, then $\mathcal{K}_n=\langle\tau_{n-1}W_3\rangle$, where $W_3$ is the vector defined in \eqref{typeIIalevel2kersigmalemmaeq1}; if $N_\pi\geq2$ and $\mu_\pi=0$, then $\mathcal{K}_n=\langle\tau_{n-1}W_{\rm new}\rangle$.
  \item \label{genericsigmakertheoremitem3} The operator $\sigma_{n-1}$ induces an isomorphism
   \begin{equation}\label{genericsigmakertheoremeq4}
    (V(n-1)\oplus V(n))/\mathcal{K}_n\stackrel{\sim}{\longrightarrow}V_s(n-1).
   \end{equation}
 \end{enumerate}
\end{theorem}
\begin{proof}
By \ref{sigmaopslemmaitem6} of Lemma \ref{sigmaopslemma}, if $W\in V_s(m)$ lies in the kernel of $\sigma_{m-1}$, then $\tau_mW$ lies in the kernel of $\sigma_m$.
It follows that
\begin{equation}\label{genericsigmakertheoremeq5}
  \mathcal{K}_n\supset
           \begin{cases}
            \langle\tau_{n-1}W_1,\tau_{n-1}W_2\rangle&\text{if }N_\pi=0,\\
            \langle\tau_{n-1}W_3\rangle&\text{if }N_\pi=1,\\
            \langle\tau_{n-1}W_{\rm new}\rangle&\text{if $N_\pi\geq2$ and $\mu_\pi=0$}.
           \end{cases}
\end{equation}
Since $\tau_{n-1}$ is injective, it follows that the numbers on the right hand side of \eqref{genericsigmakertheoremeq2} are lower bounds for the dimension of $\mathcal{K}_n$. Hence
\begin{equation}\label{genericsigmakertheoremeq7}
 \dim(\mathcal{K}_n)\geq\deg(P_0),
\end{equation}
where $P_0$ is the polynomial defined in Theorem \ref{W0theorem}.

Next, using Lemma \ref{sigmakerlemma}, we see that there exist vectors $Y_0,\ldots,Y_{n-\bar N_{\pi,s}}$ in $V(n-1)\oplus V(n)$ such that
\begin{equation}\label{genericsigmakertheoremeq6}
Y_i+\tau^i \theta^j W_0 \qquad \text{$0\leq i \leq n-\bar N_{\pi,s}$, $i+j=n-\bar N_{\pi,s}$}
\end{equation}
lie in $\ker(\sigma_{n-1})$. Now
\begin{align*}
n-N_\pi+2&=\dim(\ker(\sigma_{n-1}))\qquad\text{(by Lemma \ref{sigmakerneldimensionlemma})}\\
&\geq \dim(\mathcal{K}_n)+n-\bar N_{\pi,s}+1  
\qquad \text{(by \eqref{W0theoremeq10} and \eqref{genericsigmakertheoremeq6})} \\
&\geq\deg(P_0)+n-\bar N_{\pi,s}+1 \qquad\text{(by \eqref{genericsigmakertheoremeq7})}\\
&=n-N_\pi+2 \qquad \text{(by \eqref{W0theoremeq3})}.
\end{align*}
It follows that $\dim(\mathcal{K}_n)=\deg(P_0)$. This proves \ref{genericsigmakertheoremitem2}.

To prove \ref{genericsigmakertheoremitem1}, we calculate as follows:
\begin{align*}
&\dim\big((V(n-1)\oplus V(n))+\ker(\sigma_{n-1})\big)\\
&\qquad=\dim(V(n-1)\oplus V(n))+\dim(\ker(\sigma_{n-1}))-\dim(\mathcal{K}_n)  \\
&\qquad=\dim(V(n-1)\oplus V(n))+n-N_\pi+2-\dim(\mathcal{K}_n) \qquad\text{(by \eqref{sigmakerneldimensionlemmaeq1})} \\
&\qquad=\Big\lfloor\frac{(n-N_\pi+1)^2}4\Big\rfloor+\Big\lfloor\frac{(n-N_\pi+2)^2}4\Big\rfloor\\
&\qquad\quad+n-N_\pi+2-\dim(\mathcal{K}_n) \qquad \text{(by Theorem 7.5.6 of \cite{NF})}\\
&\qquad=\frac{(n-N_\pi+1)(n-N_\pi+2)}2+n-N_\pi+2-\dim(\mathcal{K}_n)\\
&\qquad=\frac{(n-N_\pi+2)(n-N_\pi+3)}2-\dim(\mathcal{K}_n)\\
&\qquad=\dim V_s(n) \qquad \text{(by \eqref{W0theoremeq15})}.
\end{align*}
This proves \eqref{genericsigmakertheoremeq1}.

Finally, \ref{genericsigmakertheoremitem3} follows from \ref{genericsigmakertheoremitem1} and Proposition \ref{sigmasurjectiveprop}.
\end{proof}

\section{The alternative model}\label{Lspi1sec}
Let $(\pi,V)$ be a generic, irreducible, admissible representation of $\GSp(4,F)$ with trivial central character.
In this and the next section we will consider $\pi$ with $L(s,\pi)=1$; this is a substantial
family of representations, and includes all $\pi$ that are supercuspidal.
We will prove two  results about such $\pi$.
First, in this section, for $\pi$ with $L(s,\pi)=1$, we will prove the existence of an alternative model for the spaces $V_s(n)$ of stable Klingen vectors in $V$ for integers $n \geq N_{\pi,s}$.
In this model, our level raising and lowering operators have simple visual interpretations.
Second, in the next section, still for $\pi$ with $L(s,\pi)=1$, we will prove that $N_\pi \geq 4$ using zeta integrals.

We begin with some preliminary observations. Assume that $L(s,\pi)=1$. Then  Theorem 7.5.3 of \cite{NF} implies that $N_\pi\geq2$, and that the Hecke eigenvalues $\lambda_\pi$ and $\mu_\pi$ are given by $\lambda_\pi=0$ and $\mu_\pi=-q^2$;
consequently, we have $N_{\pi,s}=N_\pi-1$ by Proposition \ref{shadowsigmavnewprop}, and $\pi$ is a \catone\index{category 1 representation}\index{representation!category 1}representation.

Let $(\pi,V)$ be a generic, irreducible, admissible representation of $\GSp(4,F)$ with trivial central character.
We will work in the Whittaker model $V=\mathcal{W}(\pi,\psi_{c_1,c_2})$ of~$\pi$. For $W \in V$ and
integers $i,j \in \Z$ with $i,j \geq 0$ we define 
$$
 m(W)_{ij} =W( \Delta_{i,j} )
$$
and let $m(W)$ be the matrix
\begin{equation}\label{mWeq}
m(W)
= (m(W)_{ij})_{0 \leq i,j < \infty}
 =
\begin{bsmallmatrix}
W(\Delta_{0,0})&W(\Delta_{0,1})&W(\Delta_{0,2})&\cdots\\
W(\Delta_{1,0})&W(\Delta_{1,1})&W(\Delta_{1,2})&\cdots\\
W(\Delta_{2,0})&W(\Delta_{2,1})&W(\Delta_{2,2})&\cdots\\
\vdots&\vdots&\vdots&
\end{bsmallmatrix}.
\end{equation}
Evidently, $m(W)$ is an element of the complex vector space $\Mat_{\infty \times \infty} (\mathbb C)$ consisting of all matrices $(m_{ij})_{0 \leq i,j <\infty}$ with $m_{ij} \in \C$ for $i,j \in \Z$ with $i,j \geq 0$. 
Now let $n$ be an integer such that $n \geq 0$. By \eqref{delltavanisheq} and Proposition \ref{parazerodiagprop},
the map $V_s(n) \to \Mat_{\infty \times \infty} (\mathbb C)$ defined by $W \mapsto m(W)$
is injective. We denote by $\Mat_s(n)$
\label{Msndef} the $\mathbb C$ vector space of all $m(W)$ for $W \in V_s(n)$, so that 
the map 
\begin{equation}
\label{altmodelisoeq}
V_s(n) \stackrel{\sim}{\longrightarrow} \Mat_s(n)
\end{equation} 
defined by $W \mapsto m(W)$
is an isomorphism. 
We refer to $\Mat_s(n)$ as the \emph{alternative model}\index{alternative model}
\index{stable Klingen vectors!alternative model}
for the space $V_s(n)$ of stable Klingen vectors. 

Our main result about the alternative model is an explicit description of the elements
of $\Mat_s(n)$. Before we present and prove this result it will be convenient to realize our usual 
level changing operators in the alternative model. 
We will write the elements $A$ of $\Mat_{\infty \times \infty} (\mathbb C)$ as a column of rows,
$$
 A = 
\begin{bsmallmatrix}
r_0 \vphantom{0}\\ r_1 \vphantom{0}\\ r_2 \vphantom{0}\\ \vdots \vphantom{0}
\end{bsmallmatrix}.
$$
We define two shift operations $\Left$ and $\Right$ on row vectors,
\begin{align*}
 \Left[a_0, a_1, a_2, \dots]&=[a_1,a_2,a_3, \dots],\\
 \Right[a_0,a_1,a_2, \dots]&=[0,a_0,a_1, \dots].
\end{align*}
Using this notation we can describe the level changing operators $\theta$, $\tau$, $\eta$, and $\sigma$ in the 
alternative model. 

\begin{proposition} \label{Ksdiagramsproposition}
 Let $\pi$ be a generic, irreducible, admissible representation of the group $\GSp(4,F)$ with 
 trivial central character, and let $V=\mathcal{W}(\pi,\psi_{c_1,c_2})$. Define
 $$
  \theta, \tau, \eta, \sigma: \Mat_{\infty \times \infty} (\mathbb C) \to \Mat_{\infty \times \infty} (\mathbb C) \label{altlevelraisdef}
  $$
 by
$$
\theta ( \begin{bsmallmatrix} 
r_0 \vphantom{\Right(r_0)}\\ 
r_1\vphantom{\Right(r_0)} \\ 
r_2\vphantom{\Right(r_0)} \\ 
\vdots\vphantom{\Right(r_0)}
\end{bsmallmatrix}  ) =
q \begin{bsmallmatrix} 0\vphantom{\Right(r_0)} \\ 
\Left(r_0)\vphantom{\Right(r_0)} \\  
\Left(r_1)\vphantom{\Right(r_0)} \\ 
\vdots \vphantom{\Right(r_0)}
\end{bsmallmatrix} +
\begin{bsmallmatrix} \Right(r_0) \\ \Right(r_1) \\ \Right(r_2) \\ \vdots\vphantom{\Right(r_0)} \end{bsmallmatrix}, \qquad 
\tau(\begin{bsmallmatrix} r_0\vphantom{0} \\ r_1\vphantom{0} \\ r_2\vphantom{0} \\ \vdots\vphantom{0} \end{bsmallmatrix} ) 
=(\begin{bsmallmatrix} r_0\vphantom{0} \\ r_1\vphantom{0} \\ r_2\vphantom{0} \\ \vdots\vphantom{0} \end{bsmallmatrix}),
$$
and
$$
\eta(\begin{bsmallmatrix} r_0\vphantom{0} \\ r_1\vphantom{0} \\ r_2\vphantom{0} \\ \vdots\vphantom{0} \end{bsmallmatrix})
=\begin{bsmallmatrix} 0\vphantom{0} \\ r_0\vphantom{0} \\ r_1\vphantom{0} \\ \vdots\vphantom{0} \end{bsmallmatrix},\qquad
\sigma(\begin{bsmallmatrix} r_0\vphantom{0} \\ r_1 \vphantom{0} \\ r_2\vphantom{0} \\ \vdots\vphantom{0} \end{bsmallmatrix})
=\begin{bsmallmatrix} r_1\vphantom{0} \\ r_2\vphantom{0} \\ r_3\vphantom{0} \\ \vdots\vphantom{0} \end{bsmallmatrix}.
 $$
 The diagrams
 \begin{equation}\label{Ksdiagramspropositioneq1}
  \begin{CD}
   V_s(n+1) @>\sim>> \Mat_s(n+1) \\
   @A\theta AA @AA\theta A\\
   V_s(n) @>\sim>> \Mat_s(n)
   \end{CD}, \qquad
   \begin{CD}
   V_s(n+1) @>\sim>> \Mat_s(n+1) \\
   @A\tau AA @AA\tau A\\
   V_s(n) @>\sim>> \Mat_s(n)
  \end{CD}
 \end{equation}
 and
 \begin{equation}\label{Ksdiagramspropositioneq2}
  \begin{CD}
   V_s(n+2) @>\sim>> \Mat_s(n+2) \\
   @A\eta AA @AA\eta A\\
   V_s(n) @>\sim>> \Mat_s(n)
   \end{CD},\qquad
  \begin{CD}
   V_s(n+1) @>\sim>> \Mat_s(n+1) \\
   @V\sigma VV @VV\sigma V\\
   V_s(n) @>\sim>> \Mat_s(n)
   \end{CD}
 \end{equation}
 commute.
\end{proposition}
\begin{proof}
This follows from \eqref{taundefeq}, \eqref{thetadefeq3} and \eqref{sigmaopseq}.
\end{proof}

We will also need the translation of the paramodular level raising operators\index{level raising operator} to the $\Mat_{\infty \times \infty}(\C)$
setting. 
Again let $(\pi,V)$ be a generic, irreducible, admissible representation of $\GSp(4,F)$ with trivial
central character, and assume that $V=\mathcal{W}(\pi,\psi_{c_1,c_2})$. Let $n$ be an integer such 
that $n \geq 0$. We denote by $\Mat(n)$
\label{localMnmodeldef} 
the $\C$ vector space of all $m(W)$ for $W \in V(n)$. We have
$\Mat(n) \subset \Mat_s(n)$, and the map $V(n) \to \Mat(n)$ defined by $W \mapsto m(W)$ for $W \in V(n)$
is an isomorphism of $\C$ vector spaces by Proposition~\ref{parazerodiagprop}. Define
$$
\theta': \Mat_{\infty \times \infty}(\C) \longrightarrow \Mat_{\infty \times \infty}(\C)
$$
by 
$$ 
\theta' (\begin{bsmallmatrix} r_0\vphantom{0} \\ r_1 \\ r_2 \\ \vdots \end{bsmallmatrix} ) = q 
\begin{bsmallmatrix} r_0 \vphantom{0}\\ r_1 \\ r_2 \\ \vdots \end{bsmallmatrix} +
\begin{bsmallmatrix} 0 \\ r_0 \\ r_1 \\ \vdots \end{bsmallmatrix}.
$$
Calculations using \eqref{eta1tdefeq}, \eqref{thetadefeq} and \eqref{thetap1defeq} show that the diagrams
 \begin{equation}\label{diagramspropositioneq1}
  \begin{CD}
   V(n+1) @>\sim >> \Mat(n+1)\\
   @A\theta AA @AA\theta A\\
   V(n) @> \sim >> \Mat(n)
   \end{CD}, 
   \qquad
   \begin{CD}
   V(n+1) @>\sim >> \Mat(n+1)\\
   @A\theta' AA @AA\theta' A\\
   V(n) @> \sim >> \Mat(n)
  \end{CD}
 \end{equation}
 and
 \begin{equation}\label{diagramspropositioneq2}
  \begin{CD}
   V(n+2) @>\sim >> \Mat(n+2)\\
   @A\eta AA @AA\eta A\\
   V(n) @> \sim >> \Mat(n)
   \end{CD}
 \end{equation}
 commute.

We turn now to the main result of this subsection. 
Let $n$ be an integer such that $n \geq 0$. 
We define
\begin{equation}\label{triangspaceeq}
 \Mat_\triangle(n)=\{A=(A_{ij})_{i,j\geq0}\in \Mat_{\infty\times\infty}(\C)\:|\:A_{ij}=0\text{ for }i+j\geq n\}.
\end{equation}
Hence $\Mat_\triangle(0)$ is the zero matrix, and
$$
 \Mat_\triangle(1)=\begin{bsmallmatrix}*&0&0&\cdots\\
                 0&0&0&\cdots\\
                 0&0&0&\cdots\\
                 \vdots&\vdots&\vdots\vphantom{0_A}
                \end{bsmallmatrix},\qquad
 \Mat_\triangle(2)=\begin{bsmallmatrix}*&*&0&\cdots\\
                 *&0&0&\cdots\\
                 0&0&0&\cdots\\
                 \vdots&\vdots&\vdots\vphantom{0_A}
                \end{bsmallmatrix},\qquad\cdots
$$
We have $\dim \Mat_\triangle(n)=\frac{n(n+1)}2$.

\begin{proposition}\label{triangprop}
Let $\pi$ be a generic, irreducible, admissible representation of the group $\GSp(4,F)$ with trivial central character, and let $V=\mathcal{W}(\pi,\psi_{c_1,c_2})$. Assume that $L(s,\pi)=1$. 
Then 
  \begin{equation}\label{triangpropeq1}
   \Mat_s(n)=\Mat_\triangle(n-(N_\pi-2))
  \end{equation}
  for all integers $n$ such that $n \geq N_{\pi,s}=N_{\pi}-1$.
\end{proposition}
\begin{proof}
By Theorem \ref{genericdimensionstheorem},
  \begin{equation}\label{triangpropeq4}
   \dim V_s(n)=
   \begin{cases}
    \displaystyle\frac{(n-N_\pi+2)(n-N_\pi+3)}2&\text{for }n\geq N_\pi-1,\\
    0&\text{for }n<N_\pi-1.
   \end{cases}
  \end{equation}
Hence $\dim V_s(n)=\dim \Mat_\triangle(n-(N_\pi-2))$ for $n\geq N_{\pi}-1$.
It is therefore enough to prove that $\Mat_s(n) \subset  \Mat_\triangle(n-(N_\pi-2))$ for $n\geq N_\pi-1$. For this we will use induction on $n$. Let $W_{\mathrm{new}}$ be a non-zero element of the one-dimensional space $V(N_\pi)$, 
and let $W_s$ be the shadow of $W_{\mathrm{new}}$.
By Theorem \ref{genericdimensionstheorem}, the space $V_s(N_\pi-1)$ is one-dimensional and spanned by  $W_s$. Moreover, it follows from Lemma 7.4.1 and Corollary 7.4.6 of \cite{NF} that
$$
m(W_s)=
\begin{bsmallmatrix}W_{\mathrm{new}}(\Delta_{0,0})&0&0&\cdots\\
0&0&0&\cdots\\
0&0&0&\cdots\\
\vdots&\vdots&\vdots\vphantom{0_A}
\end{bsmallmatrix}.
$$
Hence, $m(W_s)$ is a non-zero element of $M_\triangle(1)$. From this it follows that $\Mat_s(n) \subset  \Mat_\triangle(n-(N_\pi-2))$ for $n=N_\pi-1$.
Now assume that $n>N_\pi-1$, and that $M_s(k) \subset  \Mat_\triangle(k-(N_\pi-2))$  for $k=n-1$. The second diagram in \eqref{Ksdiagramspropositioneq1} shows that
\begin{equation}\label{triangpropeq5}
 \Mat_\triangle(n-1-(N_\pi-2))\subset \Mat_s(n).
\end{equation}
Let $E_{ij}\in \Mat_{\infty\times\infty}(\C)$ be the matrix with $(i,j)$-coefficient $1$ and all other coefficients~$0$. For a positive integer $k$ let $\Mat_\triangle'(k)$ be the subspace of $\Mat_\triangle(k)$ spanned by $E_{i,k-i-1}$ with $0\leq i\leq k-2$. Hence,
$$
 \Mat_\triangle'(1)=\begin{bsmallmatrix}0&0&\cdots\\
                 0&0&\cdots\\
                 \vdots&\vdots\vphantom{0_A}
                \end{bsmallmatrix},\quad
 \Mat_\triangle'(2)=\begin{bsmallmatrix}0&*&0&\cdots\\
                 0&0&0&\cdots\\
                 0&0&0&\cdots\\
                 \vdots&\vdots&\vdots\vphantom{0_A}
                \end{bsmallmatrix},\quad
 \Mat_\triangle'(3)=\begin{bsmallmatrix}0&0&*&0&\cdots\\
                 0&*&0&0&\cdots\\
                 0&0&0&0&\cdots\\
                 0&0&0&0&\cdots\\
                 \vdots&\vdots&\vdots&\vdots\vphantom{0_A}
                \end{bsmallmatrix},\qquad\cdots
$$
The first diagram in \eqref{Ksdiagramspropositioneq1} shows that
\begin{equation}\label{triangpropeq6}
 \Mat'_\triangle(n-(N_\pi-2))\subset \Mat_s(n).
\end{equation}
The second diagram in \eqref{diagramspropositioneq1} shows that
\begin{equation}\label{triangpropeq7}
 E_{n-N_\pi+1,0}\in(\theta')^{n-N_\pi}(W_{\rm new})+\Mat_\triangle(n-1-(N_\pi-2))\subset m(V_s(n)).
\end{equation}
The inclusions \eqref{triangpropeq5}, \eqref{triangpropeq6} and \eqref{triangpropeq7} show that
\begin{equation}\label{triangpropeq8}
 \Mat_\triangle(n-(N_\pi-2))\subset \Mat_s(n).
\end{equation}
This completes the proof. 
\end{proof}

Using the alternative model we can provide an explicit description
of the kernel of the level lowering operator $\sigma_{n-1}$ when
$\pi$ is generic and $L(s,\pi)=1$. 
For a positive integer $k$ let $\Mat_0(k)$
\label{M0kdef} be the subspace of $\Mat_\triangle(k)$ spanned by $E_{0,i}$ with $0\leq i\leq k-1$. Hence,
$$
 \Mat_0(1)=\begin{bsmallmatrix}*&0&\cdots\\
                 0&0&\cdots\\
                 \vdots&\vdots
                \end{bsmallmatrix},\quad
 \Mat_0(2)=\begin{bsmallmatrix}*&*&0&\cdots\\
                 0&0&0&\cdots\\
                 0&0&0&\cdots\\
                 \vdots&\vdots&\vdots
                \end{bsmallmatrix},\quad
 \Mat_0(3)=\begin{bsmallmatrix}*&*&*&0&\cdots\\
                 0&0&0&0&\cdots\\
                 0&0&0&0&\cdots\\
                 0&0&0&0&\cdots\\
                 \vdots&\vdots&\vdots&\vdots
                \end{bsmallmatrix},\qquad\cdots
$$

\begin{corollary}\label{triangpropcor}
Let $\pi$ be a generic, irreducible, admissible representation of $\GSp(4,F)$ with trivial central character, and let $V=\mathcal{W}(\pi,\psi_{c_1,c_2})$. Assume that $L(s,\pi)=1$. 
 Under the isomorphism \eqref{altmodelisoeq}, we have
  \begin{equation}\label{triangpropcoreq1}
   \ker(\sigma_{n-1})\stackrel{\sim}{\longrightarrow}\Mat_0(n-(N_\pi-2))
  \end{equation}
  for all $n\geq N_\pi$.
\end{corollary}
\begin{proof}
This follows  from Proposition \ref{triangprop} and Proposition \ref{Ksdiagramsproposition}.
\end{proof}

\section{A lower bound on the paramodular level} 
Let $(\pi,V)$ be a generic, irreducible, admissible representation of $\GSp(4,F)$
with trivial central character. In this final section we prove that if $L(s,\pi)=1$,
then $N_\pi \geq 4$. We begin with a result about representations of $\GL(2,F)$. In the
following lemma, $\Gamma_0(\p)$ is the subgroup of $\begin{bsmallmatrix}a&b\\c&d \end{bsmallmatrix}
$ in $\GL(2,\OF)$ such that $c \in \p$. 

\begin{lemma}
\label{Iwahorihyperlemma}
Let $(\tau,W)$ be a smooth representation of $\GL(2,F)$.
Let $w \in W$ be such that $\tau(k)w=w$ for $k \in \Gamma_0(\p)$.
Assume that there exists an integer $j$ such that $j \geq 1$ and
$$
\int\limits_{\p^{-j}} \tau(
\begin{bsmallmatrix}
1&y\\
&1
\end{bsmallmatrix})w\, dy
=0.
$$
Then $w=0$.
\end{lemma}
\begin{proof}
Let $X$ be the subset of $w \in W^{\Gamma_0(\p)}$ such that
there exists an integer $j$ such that $j \geq 1$ and 
\begin{equation}
\label{Iwahorihypereq1}
\int\limits_{\p^{-j}}\tau(
\begin{bsmallmatrix}
1&y\\
&1
\end{bsmallmatrix})w\, dy =0.
\end{equation}
The set $X$ is a subspace of $W^{\Gamma_0(\p)}$; we need to prove
that $X=0$. Assume that $X \neq 0$; we will obtain a contradiction.
Let $j_0$ be the smallest integer such that $j_0 \geq 1$ and there
exists a non-zero element $w_0 \in X$ such that \eqref{Iwahorihypereq1}
holds with $j$ and $w$ replaced by $j_0$ and $w_0$, respectively. We 
will first prove that $j_0=1$. Suppose that $j_0>1$; we will obtain a 
contradiction.

We introduce an operator on $W$. Define $U: W \to W$ by
$$
Uv = \tau(
\begin{bsmallmatrix}
1&\\
&\varpi^{-1}
\end{bsmallmatrix})
\int\limits_{\p^{-1}}
\tau(
\begin{bsmallmatrix}
 1&x\\
 &1
\end{bsmallmatrix})v\, dx
$$
for $v \in W$. We claim that if $v \in W^{\Gamma_0(\p)}$, then 
$Uv \in W^{\Gamma_0(\p)}$. To see this, let $v \in W^{\Gamma_0(\p)}$.
It is clear that 
$\tau(\begin{bsmallmatrix} a&\\&d\end{bsmallmatrix})v =v$
and $\tau(\begin{bsmallmatrix}1&b\\&1\end{bsmallmatrix})v=v$
for $a,d \in \OF^\times$ and $b \in \OF$. Let $c \in \p$. Then
\begin{align*}
\tau(\begin{bsmallmatrix} 1& \\ c&1 \end{bsmallmatrix}) Uv
& = 
\tau(\begin{bsmallmatrix} 1& \\ & \varpi^{-1} \end{bsmallmatrix})
\int\limits_{\p^{-1}}
\tau(\begin{bsmallmatrix} (1+xc\varpi)^{-1} &x\\ & 1+ xc\varpi \end{bsmallmatrix}
\begin{bsmallmatrix} 1& \\ (1+xc\varpi)^{-1} c\varpi & 1 \end{bsmallmatrix}) v \, dx\\
&= 
\tau(\begin{bsmallmatrix} 1& \\ & \varpi^{-1} \end{bsmallmatrix})
\int\limits_{\p^{-1}}
\tau(\begin{bsmallmatrix} 1 &(1+xc\varpi)^{-1}x \\ & 1 \end{bsmallmatrix}
) v \, dx\\
&=Uv.
\end{align*}
It follows that $Uv \in W^{\Gamma_0(\p)}$. 

By the last paragraph, $U w_0 \in W^{\Gamma_0(\p)}$.  
Since $j_0 \geq 2$, and by the definition of $U$ and the minimality of $j_0$, we must have $Uw_0 \neq 0$.
Define $w_1=Uw_0$. Then $w_1 \in W^{\Gamma_0(\p)}$. We have
$$
\int\limits_{\p^{-(j_0-1)}} 
\tau(\begin{bsmallmatrix} 1&y \\ &1 \end{bsmallmatrix})w_1 \, dy
=\tau(\begin{bsmallmatrix}1& \\ & \varpi^{-1} \end{bsmallmatrix})
\int\limits_{\p^{-j_0} }\tau(
\begin{bsmallmatrix} 1&y \\ &1 \end{bsmallmatrix})\, w_0 \, dy =0.
$$
This implies that $w_1 \in X$, and  
 contradicts the minimality of $j_0$; it follows that $j_0=1$.

We now have
\begin{align*}
0& = \int\limits_{\p^{-1}}\tau(
\begin{bsmallmatrix} 1&y \\ &1 \end{bsmallmatrix})w_0\, dy\\
& = \int\limits_{\OF}\tau(
\begin{bsmallmatrix} 1&y \\ &1 \end{bsmallmatrix})w_0\, dy
+ \int\limits_{\OF^\times \varpi^{-1}}\tau(
\begin{bsmallmatrix} 1&y \\ &1 \end{bsmallmatrix})w_0\, dy\\
& = w_0
+q \int\limits_{\OF^\times}\tau(
\begin{bsmallmatrix} 1&u\varpi^{-1}\\ &1 \end{bsmallmatrix})w_0\, du.
\end{align*}
Hence, since $w_0 \in W^{\Gamma_0(\p)}$, 
\begin{align}
w_0
&=-q \int\limits_{\OF^\times}\tau(
\begin{bsmallmatrix} \vphantom{\varpi^{-1}} 1&\\ u^{-1} \varpi &1 \end{bsmallmatrix}
\begin{bsmallmatrix} u\varpi^{-1} &\\ &u^{-1} \varpi \end{bsmallmatrix}
\begin{bsmallmatrix}\vphantom{\varpi^{-1}} & 1 \\ -1 & \vphantom{\varpi^{-1}}\end{bsmallmatrix}
\begin{bsmallmatrix}\vphantom{\varpi^{-1}} 1&\\ u^{-1} \varpi &1 \end{bsmallmatrix}
)w_0\, du\nonumber\\
&=-q \int\limits_{\OF^\times}\tau(
\begin{bsmallmatrix}  1&\\ u^{-1} \varpi &1 \end{bsmallmatrix}
\begin{bsmallmatrix} \varpi^{-1} &\\ & \varpi \end{bsmallmatrix}
\begin{bsmallmatrix} & 1 \\ -1 & \end{bsmallmatrix}
)w_0\, du\nonumber \\
w_0&=-q
\tau(
\begin{bsmallmatrix} \varpi^{-1} &\\ & \varpi \end{bsmallmatrix}
\begin{bsmallmatrix} & 1 \\ -1 & \end{bsmallmatrix})
\int\limits_{\OF^\times}\tau(
\begin{bsmallmatrix}  1 &u^{-1} \varpi^{-1}\\ &1 \end{bsmallmatrix}
)w_0\, du.\label{Iwahorihypereq2}
\end{align}
Using \eqref{Iwahorihypereq2} we find that
\begin{align*}
w_0 & = q \int\limits_{\p} 
\tau(
\begin{bsmallmatrix}
1&\\
c&1 
\end{bsmallmatrix})w_0\, dc\\
& = -q^2 \int\limits_{\p} 
\tau(
\begin{bsmallmatrix}
1&\\
c&1 
\end{bsmallmatrix})
\tau(
\begin{bsmallmatrix} \varpi^{-1} &\\ & \varpi \end{bsmallmatrix}
\begin{bsmallmatrix} & 1 \\ -1 & \end{bsmallmatrix})
\int\limits_{\OF^\times}\tau(
\begin{bsmallmatrix}  1 &u^{-1} \varpi^{-1}\\ &1 \end{bsmallmatrix}
)w_0\, du\, dc\\
& = -q^2 \tau(
\begin{bsmallmatrix} \varpi^{-1} &\\ & \varpi \end{bsmallmatrix}
\begin{bsmallmatrix} & 1 \\ -1 & \end{bsmallmatrix})
\int\limits_{\OF^\times} 
\int\limits_{\p}\tau(
\begin{bsmallmatrix}  1 &(-c + u^{-1} \varpi)\varpi^{-2}\\ &1 \end{bsmallmatrix}
)w_0\, dc\, du\\
& = -q^2(1-q^{-1}) \tau(
\begin{bsmallmatrix} \varpi^{-1} &\\ & \varpi \end{bsmallmatrix}
\begin{bsmallmatrix} & 1 \\ -1 & \end{bsmallmatrix})
\int\limits_{\p}\tau(
\begin{bsmallmatrix}  1 &c\varpi^{-2}\\ &1 \end{bsmallmatrix}
)w_0\, dc\\
& = -(1-q^{-1}) \tau(
\begin{bsmallmatrix} \varpi^{-1} &\\ & \varpi \end{bsmallmatrix}
\begin{bsmallmatrix} & 1 \\ -1 & \end{bsmallmatrix})
\int\limits_{\p^{-1}}\tau(
\begin{bsmallmatrix}  1 &c\\ &1 \end{bsmallmatrix}
)w_0\, dc\\
&=0,
\end{align*}
where the last step follows because $j_0=1$. This contradicts $w_0 \neq 0$,
and completes the proof.
\end{proof}

Let $n$ be an integer such that $n \geq 2$. We define
$$
\mathrm{Kl}_{s,1}(\p^n)
=
\{ g \in \GSp(4,F)\mid \lambda (g) \in \OF^\times\}
\cap
\begin{bsmallmatrix}
\OF&\OF&\OF&\p^{-n+1}\\
\p^{n-1}&\OF&\OF&\OF\\
\p^n&\p&\OF&\OF\\
\p^n&\p^n&\p^{n-1}&\OF
\end{bsmallmatrix}. \label{Ks1pndef}
$$
Using \eqref{ginveq} it is easy to verify that $\mathrm{Kl}_{s,1}(\p^n)$ is 
a subgroup of $\GSp(4,F)$. 

\begin{lemma}
\label{Ks1iwahoridecomplemma}
Let $n$ be an integer such that $n \geq 2$. Then $\mathrm{Kl}_{s,1}(\p^n)$ is equal to
$$
(\begin{bsmallmatrix}
1&\OF&\OF&\p^{-n+1}\\
&1&&\OF\\
&&1&\OF\\
&&&1
\end{bsmallmatrix}
\cap \mathrm{Kl}_{s,1}(\p^n) )
(\begin{bsmallmatrix}
\OF^\times&&&\\
&\OF&\OF&\\
&\p&\OF&\\
&&&\OF^\times
\end{bsmallmatrix}
\cap \mathrm{Kl}_{s,1}(\p^n))
(\begin{bsmallmatrix}
1&&&\\
\p^{n-1}&1&&\\
\p^n&&1&\\
\p^n&\p^n&\p^{n-1}&1
\end{bsmallmatrix}
\cap 
\mathrm{Kl}_{s,1}(\p^n) ).
$$
\end{lemma}
\begin{proof}
Let $k \in \mathrm{Kl}_{s,1}(\p^n)$; we need to prove that $k$ is in the product.
Let $k=(k_{ij})_{1 \leq i,j \leq 4}$. A calculation shows that 
$\det(k) \equiv k_{11}k_{22}k_{33}k_{44}\ \text{(mod $\p$)}$;
since $\det(k)^2 = \lambda(k)^4$ and $\lambda(k) \in \OF^\times$, we obtain
$k_{11},k_{22},k_{33},k_{44} \in \OF^\times$.
It follows that 
$$
\begin{bsmallmatrix}
1&&&-k_{14}k_{44}^{-1}\\
&1&&\\
&&1&\\
&&&1
\end{bsmallmatrix}k
\in \Kl{n-1}.
$$
Since $\Kl{n-1}$ has an Iwahori decomposition\index{Iwahori decomposition}, there exist $x,y,z \in \OF$,
$t,\lambda \in \OF^\times$, $\begin{bsmallmatrix}a&b\\c&d\end{bsmallmatrix} 
\in \GL(2,\OF)$ with $ad-bc = \lambda$, and $x',y',z' \in \p^{n-1}$ such that
$$
\begin{bsmallmatrix}
1&&&-k_{14}k_{44}^{-1}\\
&1&&\\
&&1&\\
&&&1
\end{bsmallmatrix}k
=
\begin{bsmallmatrix}
1&x&y&z\vphantom{x'}\\
&1&&y\vphantom{x'}\\
&&1&-x\vphantom{x'}\\
&&&1\vphantom{x'}
\end{bsmallmatrix}
\begin{bsmallmatrix}
t&&&\vphantom{x'}\\
&a&b&\vphantom{x'}\\
&c&d&\vphantom{x'}\\
&&&\lambda t^{-1}\vphantom{x'}
\end{bsmallmatrix}
\begin{bsmallmatrix}
1&&&\\
x'&1&&\\
y'&&1&\\
z'&y'&-x'&1
\end{bsmallmatrix}.
$$
From this equation we see that $y',z' \in \p^n$, $x' \in \p^{n-1}$, and $c \in \p$. The 
lemma follows.
\end{proof}

For later use, we note that if $n$ is an integer such that $n \geq 2$, then 
\begin{equation}
\label{Kls1uninveq}
u_n \mathrm{Kl}_{s,1}(\p^n) u_n^{-1} =
\mathrm{Kl}_{s,1}(\p^n)
\end{equation}
Here $u_n$ is as in \eqref{ALeq}.

In the following we will use the Iwahori subgroup $I$ defined in \eqref{iwahoridefeq}.
\begin{lemma}
\label{VIzerozjlemma}
Let $(\pi,V)$ be a smooth representation of $\GSp(4,F)$
for the which the center acts trivially. Assume that $V^I=0$.
Then 
\begin{equation}
\label{VIzerozjeq1}
\int\limits_{\OF}\int\limits_{\OF}\int\limits_{\OF}
\pi(
\begin{bsmallmatrix}
1&x\varpi^{-1}&y\varpi^{-1}&z\varpi^{-2}\\
&1&&y\varpi^{-1}\\
&&1&-x\varpi^{-1}\\
&&&1
\end{bsmallmatrix}
)w\, dx\,dy\,dz=0
\end{equation}
for all $w \in V^{\mathrm{Kl}_{s,1}(\p^2)}$. 
\end{lemma}

\begin{proof}
Let $w \in V^{\mathrm{Kl}_{s,1}(\p^2)}$, and let $u$ be the integral in 
\eqref{VIzerozjeq1}; we need to prove that $u=0$. Define $u_1=\pi(\eta^{-1})u$.
Given the assumption $V^I=0$, to prove that $u=0$ it will suffice to prove
that $u_1 \in V^I$. Since $I=I_+T(\OF)I_-$ (see \eqref{iwahoriplusminuseq}),
it suffices to show that $\pi(k)u_1=u_1$ for $k \in I_+$, $k \in T(\OF)$,
and $k \in I_-$; these statements can be verified by direct
computations using $w \in V^{\mathrm{Kl}_{s,1}(\p^2)}$. 
\end{proof}

Let $(\pi,V)$ be a smooth representation of $\GSp(4,F)$ for which
the center acts trivially, and let $n$ be an integer such that $n \geq 2$.
As in Sect.~7.3 of \cite{NF}, we define a linear map $R_{n-1}:V \to V$  by 
\begin{equation}
\label{Rdefeq}
R_{n-1}v =
q 
\int\limits_{\OF} \pi(
\begin{bsmallmatrix}
1&&&\\
x\varpi^{n-1}&1&&\\
&&1&\\
&&-x\varpi^{n-1}&1 
\end{bsmallmatrix})v\, dx
\end{equation}
for $v \in V$.
Evidently, if $v \in V$ is invariant under the elements
$\begin{bsmallmatrix} 1&&&\\ x&1&& \\ &&1& \\ &&-x&1 \end{bsmallmatrix}$
for $x\in\p^n$, then 
\begin{equation}
\label{Rnminus1sumeq}
R_{n-1} v
=
\sum_{x \in \OF/\p} \pi(
\begin{bsmallmatrix}
1&&&\\
x\varpi^{n-1}&1&&\\
&&1&\\
&&-x\varpi^{n-1}&1
\end{bsmallmatrix})v.
\end{equation}
We also define a linear map 
$S:V \to V$ by \label{Soperdef}
\begin{equation}
\label{Sdefeq}
Sv =
\vl(\Gamma_0(\p))^{-1} 
\int\limits_{\GL(2,\OF)}
\pi(
\begin{bsmallmatrix}
1&&\\
&g&\\
&&\det(g)
\end{bsmallmatrix})v\, dg
\end{equation}
for $v \in V$. 
Here, $dg$ is a Haar measure on $\GL(2,F)$. 
We note that if $v \in V$ is invariant under 
the elements $\begin{bsmallmatrix} 1&&\\&k&\\&& \det(k) \end{bsmallmatrix}$
for $k \in \Gamma_0(\p)$, then 
\begin{equation}
\label{simpSdefeq}
Sv =
\pi(s_2) v
+ q \int\limits_{\OF} \pi(
\begin{bsmallmatrix}
1&&&\\
&1&&\\
&y&1&\\
&&&1
\end{bsmallmatrix}) v\, dy.
\end{equation}
The following result has some overlap with Lemma~3 of \cite{Yi}.

\begin{lemma}
\label{RSinjlemma}
Let $(\pi,V)$ be a smooth representation of $\GSp(4,F)$
for which the center acts trivially, and let $n$ be an
integer such that $n \geq 2$. If $v \in V_s(n)$, then
$R_{n-1}v \in V^{\mathrm{Kl}_{s,1}(\p^n)}$, and if $v \in 
V^{\mathrm{Kl}_{s,1}(\p^n)}$, then $Sv \in V_s(n)$. The 
linear map $R_{n-1}:V_s(n) \to V^{\mathrm{Kl}_{s,1}(\p^n)}$ is injective.
If $V^I=0$, then the map $S:V^{\mathrm{Kl}_{s,1}(\p^2)} \to V_s(2)$
is injective. 
\end{lemma}
\begin{proof}
If $v \in V_s(n)$, then calculations using the Iwahori decomposition from Lemma~\ref{Ks1iwahoridecomplemma}
show that $R_{n-1}v \in V^{\mathrm{Kl}_{s,1}(\p^n)}$; if $v \in V^{\mathrm{Kl}_{s,1}(\p^n)}$,
then calculations using the Iwahori decomposition \eqref{Ksniwahorieq1} show
that $Sv \in V_s (n)$. Next, let $v \in V_s(n)$. Then by \eqref{Rnminus1sumeq} and \eqref{simpSdefeq}, 
\begin{align*}
&S(R_{n-1}(v))\\
&\qquad=\sum_{y\in\OF/\p}\pi(
\begin{bsmallmatrix}1\\&1\\&y&1\\&&&1\end{bsmallmatrix})R_{n-1}(v)+\pi(s_2)R_{n-1}(v)\\
&\qquad=\sum_{x,y\in\OF/\p}\pi(
\begin{bsmallmatrix}
1\vphantom{\varpi^{n-1}}\\
&1\vphantom{\varpi^{n-1}}\\
&y&1\vphantom{\varpi^{n-1}}\\
&&&1\vphantom{\varpi^{n-1}}
\end{bsmallmatrix}
\begin{bsmallmatrix}1\\x\varpi^{n-1}&1\\&&1\\&&-x\varpi^{n-1}&1\end{bsmallmatrix})v\\
&\qquad\quad+\sum_{x\in\OF/\p}\pi(s_2
\begin{bsmallmatrix}1\\x\varpi^{n-1}&1\\&&1\\&&-x\varpi^{n-1}&1\end{bsmallmatrix})v\\
&\qquad=\sum_{x,y\in\OF/\p}\pi(\begin{bsmallmatrix}1\\x\varpi^{n-1}&1\\&&1\\&&-x\varpi^{n-1}&1\end{bsmallmatrix}
\begin{bsmallmatrix}1\\&1\\yx\varpi^{n-1}&y&1\\yx^2\varpi^{2n-2}&yx\varpi^{n-1}&&1\end{bsmallmatrix})v\\
&\qquad\quad+\sum_{x\in\OF/\p}\pi(
\begin{bsmallmatrix}1\\&1\\x\varpi^{n-1}&&1\\&x\varpi^{n-1}&&1\end{bsmallmatrix}s_2)v\\
&\qquad=qv+\sum_{x\in(\OF/\p)^\times}\sum_{y\in\OF/\p}\pi(
\begin{bsmallmatrix}1\\x\varpi^{n-1}&1\\&&1\\&&-x\varpi^{n-1}&1\end{bsmallmatrix}
\begin{bsmallmatrix}1\\&1\\yx\varpi^{n-1}&&1\\&yx\varpi^{n-1}&&1\end{bsmallmatrix})v\\
&\qquad\quad+\sum_{x\in\OF/\p}\pi(
\begin{bsmallmatrix}1\\&1\\x\varpi^{n-1}&&1\\&x\varpi^{n-1}&&1\end{bsmallmatrix})v\\
&\qquad=qv+\sum_{x\in(\OF/\p)^\times}\sum_{y\in\OF/\p}\pi(
\begin{bsmallmatrix}1\\x\varpi^{n-1}&1\\&&1\\&&-x\varpi^{n-1}&1\end{bsmallmatrix}
\begin{bsmallmatrix}1\\&1\\y\varpi^{n-1}&&1\\&y\varpi^{n-1}&&1\end{bsmallmatrix})v\\
&\qquad\quad+\sum_{y\in\OF/\p}\pi(
\begin{bsmallmatrix}1\\&1\\y\varpi^{n-1}&&1\\&y\varpi^{n-1}&&1\end{bsmallmatrix})v\\
&\qquad=qv+\sum_{x,y\in\OF/\p}\pi(
\begin{bsmallmatrix}1\\x\varpi^{n-1}&1\\&&1\\&&-x\varpi^{n-1}&1\end{bsmallmatrix}
\begin{bsmallmatrix}1\\&1\\y\varpi^{n-1}&&1\\&y\varpi^{n-1}&&1\end{bsmallmatrix})v\\
&\qquad=qv+\sum_{x,y\in\OF/\p}\pi(
\begin{bsmallmatrix}1\\&1\\y\varpi^{n-1}&&1\\&y\varpi^{n-1}&&1\end{bsmallmatrix}
\begin{bsmallmatrix}1\\x\varpi^{n-1}&1\\&&1\\&&-x\varpi^{n-1}&1\end{bsmallmatrix})v\\
&\qquad=qv+\sum_{y\in\OF/\p}\pi(
\begin{bsmallmatrix}1\\&1\\y\varpi^{n-1}&&1\\&y\varpi^{n-1}&&1\end{bsmallmatrix})R_{n-1}(v).
\end{align*}
By the last equation, if $R_{n-1}(v)=0$, then $v=0$, so that $R_{n-1}:V_s(n) \to V^{\mathrm{Kl}_{s,1}(\p^n)}$
is injective.

Finally, assume that $V^I=0$. 
Let $v \in V^{\mathrm{Kl}_{s,1}(\p^2)}$, and assume that
$Sv=0$. Since $Sv=0$, we have by \eqref{simpSdefeq},
$$
\pi(s_2) v =
- q \int\limits_{\OF} \pi(
\begin{bsmallmatrix}
1&&&\\
&1&&\\
&y&1&\\
&&&1
\end{bsmallmatrix}) v\, dy.
$$
Applying $u_2$ to this equation, and using $\pi(u_2 s_2 u_2^{-1})\pi(u_2)v = \pi(t_2)\pi(u_2)v$, we obtain
$$
\pi(t_2) v_1
=
- q \int\limits_{\OF} \pi(
\begin{bsmallmatrix}
1&&&y\varpi^{-2}\\
&1&&\\
&&1&\\
&&&1
\end{bsmallmatrix}) v_1\, dy
$$
where $v_1= \pi(u_2) v$. Since $u_2$ normalizes
$\mathrm{Kl}_{s,1}(\p^2)$, it follows that 
$v_1 \in V^{\mathrm{Kl}_{s,1}(\p^2)}$. 
If $x \in \OF$, then 
$$
t_2^{-1} 
\begin{bsmallmatrix}
1&&x\varpi^{-1}&\\
&1&&x\varpi^{-1}\\
&&1&\\
&&&1
\end{bsmallmatrix}
t_2 =
\begin{bsmallmatrix}
1&&&\\
x\varpi&1&&\\
&&1&\\
&&-x\varpi&1
\end{bsmallmatrix}.
$$
It follows that $\pi(t_2)v_1$ is invariant under the 
elements 
$$
\begin{bsmallmatrix}
1&&x\varpi^{-1}&\\
&1&&x\varpi^{-1}\\
&&1&\\
&&&1
\end{bsmallmatrix}
$$
for $x \in \OF$. Hence,
\begin{equation}
\label{RSinjeq21}
\pi(t_2) v_1
=
- q \int\limits_{\OF}\int\limits_{\OF} \pi(
\begin{bsmallmatrix}
1&&x\varpi^{-1}&y\varpi^{-2}\\
&1&&x\varpi^{-1}\\
&&1&\\
&&&1
\end{bsmallmatrix}) v_1\, dy\, dx.
\end{equation}
We now define a smooth representation $(\tau,W)$
of $\GL(2,F)$ by letting $W=V$ and setting
$\tau(g)w = \pi(\begin{bsmallmatrix} g&\\ & g' \end{bsmallmatrix})w$
for $g \in \GL(2,F)$ and $w \in W$. 
Let $w = \pi(t_2) v_1$. Calculations 
using \eqref{RSinjeq21} and 
$v_1 \in V^{\mathrm{Kl}_{s,1}(\p^2)}$ show that 
$\tau(k)w=w$ for $k \in \Gamma_0(\p)$. By Lemma \ref{VIzerozjlemma}
we have 
$$
\int\limits_{\p^{-1}} \tau(
\begin{bsmallmatrix}
1&y\\
&1
\end{bsmallmatrix})w\, dy=0.
$$
Lemma \ref{Iwahorihyperlemma} now implies that $w=0$, so that $v=0$.
\end{proof}

Let $(\pi,V)$ be a generic, irreducible, admissible representation
with trivial central character. As in (7.14) of \cite{NF}, if 
$W \in \mathcal{W}(\pi,\psi_{c_1,c_2})$, then we define
$$
Z_N(s,W)
=
\int\limits_{F^\times}
W(
\begin{bsmallmatrix}
a&&&\\
&a&&\\
&&1&\\
&&&1
\end{bsmallmatrix})
|a|^{s-\frac{3}{2}}\, d^\times a. \label{simpzetadef}
$$

\begin{lemma}
\label{ZNs2Runlemma}
Let $(\pi,V)$ be a generic, irreducible, admissible
representation of the group $\GSp(4,F)$ with trivial central 
character. Let $n$ be an integer such that $n \geq 2$.
Let $W \in V_s(n)$. Then
$$
Z(s, \pi(s_2)R_{n-1} W)=
\gamma(s,\pi)^{-1}q^{n/2-ns} Z_N(1-s,\pi(u_n)W).
$$
\end{lemma}
\begin{proof}
The proof of this assertion is the same as the proof of 
Proposition~7.3.2 of \cite{NF}, except that the last
line of the last display in that proof is omitted.
\end{proof}

\begin{lemma}
\label{ZNs2RNminus2zerolemma}
Let $(\pi,V)$ be a generic, irreducible, admissible
representation of the group $\GSp(4,F)$ with trivial central
character. Assume that $L(s,\pi)=1$, so that $N_\pi \geq 2$ and  
$\pi$ 
is a \index{category 1 representation}\index{representation!category 1}\catone representation. Let $W_{\mathrm{new}}$
be a non-zero vector in $V(N_\pi)$, and let $W_s
\in V_s(N_\pi-1)$ be the \index{shadow of a newform}shadow of $W_{\mathrm{new}}$. 
Then
$$
Z(s,\pi(s_2)R_{N_\pi-2} W_s)=0.
$$
\end{lemma}
\begin{proof}
We will abbreviate $N=N_\pi$ in this proof. 
By Lemma \ref{ZNs2Runlemma} it suffices to prove that
$Z_N(s,\pi(u_{N-1})W_s)=0$. By \eqref{shadowdefeq200},
$$
W_s
=
q^3
\int\limits_{\OF^3}
\pi(
\begin{bsmallmatrix}
1&&&\\
x\varpi^{N-1}&1&&\\
y\varpi^{N-1}&&1&\\
z\varpi^{N-1}&y\varpi^{N-1}&-x\varpi^{N-1}&1
\end{bsmallmatrix}) W_{\mathrm{new}}\, dx\, dy\, dz.
$$
Hence,
$$
\pi(u_{N-1})W_s
=
q^3
\int\limits_{\OF^3}
\pi(
\begin{bsmallmatrix}
1&&y&\\
x\varpi^{N-1}&1&z&y\\
&&1&\\
&&-x\varpi^{N-1}&1
\end{bsmallmatrix}) \pi(u_{N-1}) W_{\mathrm{new}}\, dx\, dy\, dz.
$$
Therefore, 
\begin{align*}
&Z_N(s,\pi(u_{N-1})W_s) \\
&\qquad = q^3 \int\limits_{F^\times}
\int\limits_{\OF^3}
W_{\mathrm{new}}(
\begin{bsmallmatrix}
a\vphantom{\varpi^{N-1}}&&&\\
&a\vphantom{\varpi^{N-1}}&&\\
&&1\vphantom{\varpi^{N-1}}&\\
&&&1\vphantom{\varpi^{N-1}}
\end{bsmallmatrix}
\begin{bsmallmatrix}
1&&y&\\
x\varpi^{N-1}&1&z&y\\
&&1&\\
&&-x\varpi^{N-1}&1
\end{bsmallmatrix} u_{N-1})|a|^{s-\frac{3}{2}}\, dx\, dy\, dz\, d^\times a \\
&\qquad = q^3 \int\limits_{F^\times}
\int\limits_{\OF^3}
W_{\mathrm{new}}(
\begin{bsmallmatrix}
a\vphantom{\varpi^{N-1}}&&&\\
&a\vphantom{\varpi^{N-1}}&&\\
&&1\vphantom{\varpi^{N-1}}&\\
&&&1\vphantom{\varpi^{N-1}}
\end{bsmallmatrix}
\begin{bsmallmatrix}
1&&y&\\
x\varpi^{N-1}&1&z&y\\
&&1&\\
&&-x\varpi^{N-1}&1
\end{bsmallmatrix}\\
&\qquad\quad\times
\begin{bsmallmatrix}
1\vphantom{\varpi^{N-1}}&&&\\
&1\vphantom{\varpi^{N-1}}&&\\
&&\varpi^{-1}\vphantom{\varpi^{N-1}}&\\
&&&\varpi^{-1}\vphantom{\varpi^{N-1}}
\end{bsmallmatrix}
u_{N})|a|^{s-\frac{3}{2}}\, dx\, dy\, dz\, d^\times a \\
&\qquad = q^3 \int\limits_{F^\times}
\int\limits_{\OF^3}
W_{\mathrm{new}}(
\begin{bsmallmatrix}
a\varpi\vphantom{\varpi^{N-1}}&&&\\
&a\varpi\vphantom{\varpi^{N-1}}&&\\
&&1\vphantom{\varpi^{N-1}}&\\
&&&1\vphantom{\varpi^{N-1}}
\end{bsmallmatrix}
\begin{bsmallmatrix}
1&&y\varpi^{-1}&\\
x\varpi^{N-1}&1&z\varpi^{-1}&y\varpi^{-1}\\
&&1&\\
&&-x\varpi^{N-1}&1
\end{bsmallmatrix}
u_{N})\\
&\qquad\quad\times |a|^{s-\frac{3}{2}}\, dx\, dy\, dz\, d^\times a \\
&\qquad = q^{s+\frac{3}{2}}\varepsilon_\pi \int\limits_{F^\times}
\int\limits_{\OF^3}
W_{\mathrm{new}}(
\begin{bsmallmatrix}
a\vphantom{\varpi^{N-1}}&&&\\
&a\vphantom{\varpi^{N-1}}&&\\
&&1\vphantom{\varpi^{N-1}}&\\
&&&1\vphantom{\varpi^{N-1}}
\end{bsmallmatrix}
\begin{bsmallmatrix}
1&&y\varpi^{-1}&\\
&1&z\varpi^{-1}&y\varpi^{-1}\\
&&1&\\
&&&1
\end{bsmallmatrix}\\
&\qquad\quad \times
\begin{bsmallmatrix}
1&&&\\
x\varpi^{N-1}&1&&\\
&&1&\\
&&-x\varpi^{N-1}&1
\end{bsmallmatrix}
)|a|^{s-\frac{3}{2}}\, dx\, dy\, dz\, d^\times a \\
&\qquad = q^{s+\frac{1}{2}}\varepsilon_\pi \int\limits_{F^\times}
(\int\limits_{\OF} \psi(c_2 az\varpi^{-1})\, dz)
(R_{N-1}W_{\mathrm{new}})(
\begin{bsmallmatrix}
a&&&\\
&a&&\\
&&1&\\
&&&1
\end{bsmallmatrix})\\
&\qquad\quad\times |a|^{s-\frac{3}{2}}\, dx\, dy\, dz\, d^\times a \\
&\qquad = q^{s+\frac{1}{2}}\varepsilon_\pi \int\limits_{\{a \in F^\times \mid v(a) \geq 1 \} }
(R_{N-1}W_{\mathrm{new}})(
\begin{bsmallmatrix}
a&&&\\
&a&&\\
&&1&\\
&&&1
\end{bsmallmatrix}
)|a|^{s-\frac{3}{2}}\, dx\, dy\, dz\, d^\times a. 
\end{align*}
By Lemma 2.4 of \cite{HNF} we have $Z_N(s,R_{N-1}W_{\mathrm{new}})=
-q^3Z_N(s,\pi(\eta^{-1})W_{\mathrm{new}})$. 
This implies that if $k \in \Z$, then
$$
(R_{N-1}W_{\mathrm{new}})
(
\begin{bsmallmatrix}
\varpi^k&&&\\
&\varpi^k&&\\
&&1&\\
&&&1
\end{bsmallmatrix}
)
=
-q^3
W_{\mathrm{new}}
(
\begin{bsmallmatrix}
\varpi^k&&&\\
&\varpi^k&&\\
&&1&\\
&&&1
\end{bsmallmatrix}
\eta^{-1})\\
=
-q^3 
W_{\mathrm{new}}
(\Delta_{1,k}).
$$
By assumption, $L(s,\pi)=1$. It follows from Theorem 7.5.1 and Corollary
7.4.6 of \cite{NF} that $W_{\mathrm{new}}(\Delta_{i,j})=0$ if $i,j \in \Z$, $(i,j) \neq (0,0)$ and $(i,j) \neq (1,0)$.
In particular, we have $W_{\mathrm{new}}(\Delta_{1,k})=0$ if $k \in \Z$ and $k \neq 0$. 
It follows that 
$$
\int\limits_{\{a \in F^\times \mid v(a) \geq 1 \} }
(R_{N-1}W_{\mathrm{new}})(
\begin{bsmallmatrix}
a&&&\\
&a&&\\
&&1&\\
&&&1
\end{bsmallmatrix}
)|a|^{s-\frac{3}{2}}\, dx\, dy\, dz\, d^\times a =0,
$$
so that $Z_N(s,\pi(u_{N-1})W_s)=0$.
\end{proof}

\begin{theorem}
\label{Npi4theorem}
Let $(\pi,V)$ be a generic, irreducible, admissible
representation of $\GSp(4,F)$ with trivial central
character. Assume that $L(s,\pi)=1$. Then $N_\pi \geq 4$. 
\end{theorem}
\begin{proof}
Since $L(s,\pi)=1$ we have $N_\pi \geq 2$ and $\pi$
is a \catone representation. In particular, 
by iii) of Theorem 7.5.3
of \cite{NF} we must have $\lambda_\pi=0$ and $\mu_\pi=-q^2$. 

We first claim that $V^I=0$. Suppose that $V^I \neq 0$;
we will obtain a contradiction. Since $V^I \neq 0$,
and since $N_\pi \geq 2$, by Table A.13 of \cite{NF}
we see that $\pi$ belongs to Group IIIa, IVa, Va, or VIa
with unramified inducing data. 
This contradicts the values of $\lambda_\pi$ and $\mu_\pi$
from Table A.14 of \cite{NF}; hence, $V^I=0$. 

Now assume that $N_\pi \leq 3$; we will obtain a contradiction.
If $N_\pi=2$, then since $\pi$ is a \catone
representation we have $V_s(1)\neq 0$; this contradicts $V^I=0$.
Therefore, $N_\pi=3$ and $N_{\pi,s}=2$. By Theorem \ref{genericdimensionstheorem}, the space $V_s(2)$ is one-dimensional
and is spanned by the shadow $W_s$ of a newform $W_{\mathrm{new}} 
\in V(3)$. By Lemma \ref{RSinjlemma} the maps 
$R_1: V_s(2) \to V^{\mathrm{Kl}_{s,1}(\p^2)}$ 
and $S: V^{\mathrm{Kl}_{s,1}(\p^2)} \to V_s(2)$ are injective
and hence isomorphisms. Let $c \in \C^\times$ be such that 
$S R_1 W_s = c W_s$. By  the definition of $S$, this means
that
$$
c W_s = \pi(s_2) R_1 W_s + q 
\int\limits_{\OF}
\pi(
\begin{bsmallmatrix}
1&&&\\
&1&&\\
&z&1&\\
&&&1
\end{bsmallmatrix}) R_1 W_s \, dy.
$$
It follows that
$$
c Z(s,W_s)
=Z(s,\pi(s_2)R_1  W_s) +q Z(s,R_1 W_s).
$$
By Lemma \ref{ZNs2RNminus2zerolemma} we have 
$Z(s,\pi(s_2)R_1  W_s)=0$. Also, by \eqref{WsWnewzetaeq}, $Z(s,W_s)$
is a non-zero constant. It follows that $Z(s,R_1 W_s) = C$
for some $C \in \C^\times$. Since $\dim V^{\mathrm{Kl}_{s,1}(\p^2)}
=1$ and $\pi(u_2) V^{\mathrm{Kl}_{s,1}(\p^2)} =
V^{\mathrm{Kl}_{s,1}(\p^2)}$, there exists $\varepsilon\in  \{ \pm 1\}$
such that $\pi(u_2) R_1 W_s = \varepsilon R_1 W_s$. Because $L(s,\pi)=1$
we have $\gamma (s,\pi) = \varepsilon (s,\pi)$. By Corollary 7.5.5
of \cite{NF}, $\varepsilon(s,\pi) = \varepsilon_\pi q^{-3(s-\frac{1}{2})}$. 
By (2.61) of \cite{NF} we have
\begin{align*}
Z(1-s,\pi(u_2) R_1 W_s)
&= q^{2(s-\frac{1}{2})} \gamma(s,\pi) Z(s,R_1 W_s).
\end{align*}
Substituting, we obtain:
\begin{align*}
\varepsilon C &= q^{2(s-\frac{1}{2})}  \varepsilon_\pi q^{-3(s-\frac{1}{2})} C\\
\varepsilon&=q^{-(s-\frac{1}{2})}\varepsilon_\pi.
\end{align*}
This contradiction completes the proof. 
\end{proof}

%% file: SKMS_chapter9.tex
\chapter{Iwahori-spherical Representations}
\label{iwahorichap}
Let $(\pi,V)$ be an irreducible admissible
representation of $\GSp(4,F)$ with trivial 
central character. We will say that $\pi$
is 
\emph{Iwahori-spherical}\index{Iwahori-spherical representation}\index{representation!Iwahori-spherical}
if the space
$V^I$ is non-zero, where $I \subset \Kl{}$ is the 
Iwahori subgroup of $\GSp(4,F)$ as defined below
in \eqref{iwahoridefeq}. 
The goal of this chapter
is to  describe
the actions of the stable Hecke operators $T_{0,1}^s$
and $T_{1,0}^s$ on $V_s(1)$ when $\pi$ is
Iwahori-spherical. 
Since $\Ks{}=\Kl{}$, $V_s(1)$ is simply the space
$V^{\Kl{}}$ of vectors fixed by the Klingen congruence
subgroup $\Kl{}$ of level $\p$. 
The dimensions of the spaces $V^K$ for
all standard parahoric subgroups $K$
of $\GSp(4,F)$ and 
 for all Iwahori-spherical
$\pi$ was calculated in \cite{SI}.
We will use this information, along with
the theory of Iwahori-spherical representations, 
to calculate the action of the stable Hecke operators\index{stable Hecke operators}.

\section{Some background}
\label{backtheorysec}

We begin by recalling some group theory about
$\GSp(4,F)$ centered around the Iwahori subgroup. 
This will allow the introduction of the Iwahori-Hecke algebra,
which will be used in the analysis of  the stable Hecke operators $T_{0,1}^s$ and $T_{1,0}^s$
as mentioned above.

\subsection*{An extended Tits system}
\begin{figure}
\caption{The Coxeter graphs of the Weyl group $W$ and
affine Weyl group $W^a$, respectively}
\label{coxeterfig}
\bigskip
\hspace{20ex}\begin{tikzpicture}[x=1cm,y=1cm]
\draw (0.5,0) -- (-0.5,0);
\draw (0.5,0) -- (0,0.866);
\draw (-0.5,0) -- (0,0.866);
\draw[fill=white] (0.5,0) circle (0.05);
\draw[fill=white] (-0.5,0) circle (0.05);
\draw[fill=white] (0,0.866) circle (0.05);
\node at (0.5,0.4){\scriptsize 2};
\node at (-0.5,0.4){\scriptsize 4};
\node at (0,-0.25){\scriptsize 4};
\node at (0,1.1){\scriptsize $s_0$};
\node at (0.8,-0.1){\scriptsize $s_2$};
\node at (-0.75,-0.1){\scriptsize $s_1$};
\draw (-3.5,0.4) -- (-2.5,0.4);
\draw[fill=white] (-3.5,0.4) circle (0.05);
\draw[fill=white] (-2.5,0.4) circle (0.05);
\node at (-3,0.6){\scriptsize 4};
\node at (-3.5,0.1){\scriptsize $s_1$};
\node at (-2.5,0.1){\scriptsize $s_2$};
\end{tikzpicture}
\end{figure}
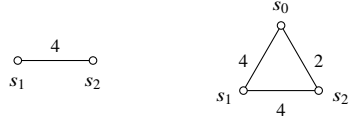

As in \eqref{Udefeq}, let $T$ be the subgroup of $\GSp(4,F)$
consisting of diagonal matrices, and let $\Normalizer(T)$ be
the normalizer of $T$ in $\GSp(4,F)$. Also, let 
\begin{equation}
\label{iwahoridefeq}
 I=\GSp(4,\OF)\cap\begin{bsmallmatrix}\OF&\OF&\OF&\OF\\\p&\OF&\OF&\OF\\\p&\p&\OF&\OF\\\p&\p&\p&\OF\end{bsmallmatrix}.
\end{equation}
Then $I$ is a compact, open subgroup of $\GSp(4,\OF)$, and we
refer to $I$ as the \emph{Iwahori subgroup}\index{Iwahori subgroup}
of $\GSp(4,F)$. The intersection $I \cap \Normalizer(T)$
is equal to $T(\OF)$, where $T(\OF)=T\cap \GSp(4,\OF)$ is the subgroup of $T$ of elements
with diagonal entries in $\OF^\times$, and $T(\OF) = I \cap \Normalizer(T)$
is a normal subgroup of $\Normalizer(T)$. We let
\begin{equation}
\label{iwahoriweylgroupeq}
W^{e} = \Normalizer(T)/(I \cap \Normalizer(T)) = \Normalizer(T) /T(\OF),
\end{equation}
and call this group the \emph{Iwahori-Weyl group}
\index{Iwahori-Weyl group}, or
the \emph{extended affine Weyl group}\index{extended affine Weyl group}. \label{Wedef}
The group $W^e$ is generated by the images of
\begin{equation}
\label{t1defeq}
s_0=t_1 = 
\begin{bsmallmatrix}
&&&-\varpi^{-1}\\
&1&&\\
&&1&\\
\varpi&&&
\end{bsmallmatrix},\quad
s_1 = 
\begin{bsmallmatrix}
&1&&\\
1&&&\\
&&&1\\
&&1&
\end{bsmallmatrix},\quad
s_2 = 
\begin{bsmallmatrix}
1&&&\\
&&1&\\
&-1&&\\
&&&1
\end{bsmallmatrix}
\end{equation}
and
\begin{equation}
u_1 = 
\begin{bsmallmatrix}
&&1&\\
&&&-1\\
\varpi&&&\\
&-\varpi&&
\end{bsmallmatrix}. \label{u1defeq}
\end{equation}
We note that $s_0=t_1$ and $u_1$ were already defined in \eqref{tndefeq}
and \eqref{ALeq}, respectively.
We let $W^a$ be the subgroup of $W^e$ generated by the images 
of $s_0$, $s_1$, and $s_2$; the subgroup $W^a$ of $W^e$ is called
the \emph{affine Weyl group}
\index{affine Weyl group}. 
The affine Weyl group $W^a$ is a Coxeter group with Coxeter \label{Wadef}
generators $s_0,s_1$ and $s_2$, and Coxeter graph as in Fig.~\ref{coxeterfig} (see \cite{BB}).
The group $W^a$ is isomorphic to the $(2,4,4)$ triangle group; since $1/2+1/4+1/4=1$,
$W^a$ is infinite.
The subgroup $W$ of $W^e$ generated by $s_1$ and $s_2$ is called the 
\emph{Weyl group}.\index{Weyl group}
We have \label{Weyldef}
$W = (\Normalizer(T) \cap \GSp(4,\OF))/T(\OF)$. The inclusion
$$
W =(\Normalizer(T) \cap \GSp(4,\OF))/T(\OF) \stackrel{\sim}{\longrightarrow} \Normalizer(T)/T
$$
is surjective, so that $W$ is naturally isomorphic to $\Normalizer(T)/T$. 
The Weyl group is a Coxeter group with Coxeter generators $s_1$ and $s_2$ and
Coxeter graph as in Fig.~\ref{coxeterfig}.
Finally, we
let $\Omega$
be the subgroup of $W^e$ generated by the image of $u_1$. \label{omegadef}
The group $\Omega$ is an infinite cyclic group.
It can be proven that the triple $(\GSp(4,F),I,\Normalizer(T))$ is  a generalized 
Tits system (see \cite{I}) with respect to $W^a$, the 
generators $s_0,s_1$ and $s_2$ of $W^a$, and $\Omega$. 
This means that: $W^e$ is the semi-direct
product of $\Omega$ and the normal subgroup $W^a$;
$\sigma I w \subset I\sigma wI \cup I\sigma I$ for $\sigma \in W^e$
and $w \in \{s_0,s_1,s_2\}$; $wIw^{-1} \neq I$ for $w \in \{s_0,s_1,s_2\}$;
$u_1\{s_0,s_1,s_2\}u_1^{-1} = \{s_0,s_1,s_2\}$ in $W^e$;
$u_1 I u_1^{-1}=I$; $I \rho \neq I$ for $\rho \in \Omega-\{1\}$;
and $\GSp(4,F)$ is generated by $I$ and $\Normalizer(T)$. 

\subsection*{Parahoric subgroups} A \emph{standard parahoric subgroup}
\index{standard parahoric subgroup} of $\GSp(4,F)$ is by definition
a compact subgroup of $\GSp(4,F)$ containing the Iwahori subgroup $I$.
Note that such a subgroup is necessarily open since it contains $I$.
The standard parahoric subgroups of $\GSp(4,F)$ are in bijection with the
proper subsets of $\{s_0,s_1,s_2\}$, the set of Coxeter generators of the 
affine Weyl group $W^a$ (see Sect.~1 of  \cite{I}). More precisely,
let $J$ be a proper subset of $\{s_0,s_1,s_2\}$. Define $W^a_J$
to be the subgroup of $W^a$ generated by $J$. Since $J$ is proper, \label{WaJdef}
$W^a_J$ is finite. The standard parahoric subgroup corresponding
to $J$ is 
$$
K_J = \bigsqcup_{w \in W^a_J} IwI. \label{KJdef}
$$
The union is disjoint because $\sqcup_{w \in W^e} IwI$ is disjoint (this
is a property of generalized Tits systems). There are seven standard parahoric
subgroups. The parahoric subgroup corresponding to the empty subset
is $K_\emptyset = I$. The parahoric subgroup corresponding to $\{s_1\}$
is called the Siegel congruence subgroup of level $\p$, and is often
denoted by $\Gamma_0(\p)$. We have 
$$
K_{\{s_1\}}=\Gamma_0(\p) = \GSp(4,\OF) \cap   
\begin{bsmallmatrix}
\OF&\OF&\OF&\OF\\
\OF&\OF&\OF&\OF\\
\p&\p&\OF&\OF\\
\p&\p&\OF&\OF
\end{bsmallmatrix}.
$$
The parahoric subgroup corresponding to $\{s_2\}$ is the Klingen
congruence subgroup $\Kl{}$ of level $\p$ which was defined in \eqref{Klndefeq}.
Thus, 
$$
K_{\{s_2\}} = \Kl{} 
 = \GSp(4,\OF) \cap   
\begin{bsmallmatrix}
\OF&\OF&\OF&\OF\\
\p&\OF&\OF&\OF\\
\p&\OF&\OF&\OF\\
\p&\p&\p&\OF
\end{bsmallmatrix}.
$$
For the parahoric subgroup corresponding to $\{s_0\}$ we have 
$W^a_{\{s_0\}}=\{1,s_0\}$ so that
$K_{\{s_0\}}=I \sqcup I s_0I$. Since $u_1 s_2 u_1^{-1} = s_0$ as elements
of $W^a$, there is an equality $K_{\{s_0\}} = u_1K_{s_2}u_1^{-1}$. 
The parahoric subgroup corresponding to $\{s_1,s_2\}$ is $K_{\{s_1,s_2\}}=\GSp(4,\OF)$. 
Note that $W_{\{s_1,s_2\}}=W$, the Weyl group.
Since $u_1 s_1 u_1^{-1} = s_1$ in $W^a$, the parahoric subgroup
corresponding to $\{s_0,s_1\}$ is $K_{\{s_0,s_1\}}=u_1\GSp(4,\OF)u_1^{-1}$. 
Finally, the parahoric subgroup corresponding to $\{s_0,s_2\}$ is the 
paramodular group $\K{}$ of level $\p$. 
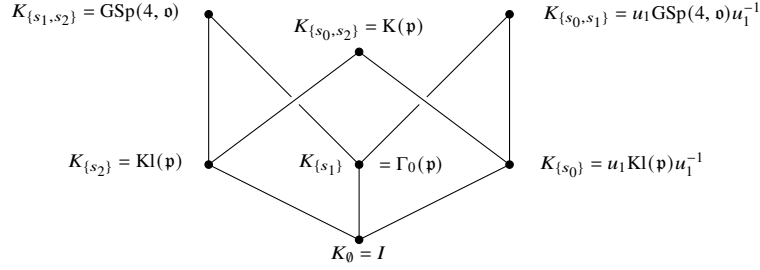
\begin{figure}
\caption{Inclusion relations between the standard parahoric subgroups
of $\GSp(4,F)$}
\label{parahoricfig}
\bigskip
\begin{tikzpicture}[x=1cm,y=1cm]
\draw[fill=black] (2,2.5) circle (0.05);
\node at (2,2.8){\scriptsize $K_{\{s_0,s_2\}}=\K{}$};
\draw[fill=black] (2,1) circle (0.05);
\node at (1.5,1){\scriptsize $K_{\{s_1\}}$};
\node at (2.7,1){\scriptsize $=\Gamma_0(\p)$};
\draw[fill=black] (2,0) circle (0.05);
\node at (2,-0.2){\scriptsize $K_{\emptyset}=I$};
\draw[fill=black] (0,1) circle (0.05);
\node at (-1.1,1){\scriptsize $K_{\{s_2\}}=\Kl{}$};
\draw[fill=black] (4,1) circle (0.05);
\node at (5.5,1){\scriptsize $K_{\{s_0\}}=u_1\Kl{}u_1^{-1}$};
\draw[fill=black] (0,3) circle (0.05);
\node at (-1.5,3){\scriptsize $K_{\{s_1,s_2\}}=\GSp(4,\OF)$};
\draw[fill=black] (4,3) circle (0.05);
\node at (5.9,3){\scriptsize $K_{\{s_0,s_1\}}=u_1\GSp(4,\OF)u_1^{-1}$};
\draw (2,0) -- (0,1);
\draw (2,0) -- (4,1);
\draw (0,3) -- (0,1);
\draw (4,3) -- (4,1);
\draw (2,2.5) -- (0,1);
\draw (2,2.5) -- (4,1);
\draw (2,0) -- (2,1);
\draw (1.2,1.8) -- (2,1);
\draw (1.1,1.9) -- (0,3);
\draw (2,1) -- (2.8,1.8);
\draw (4,3) -- (2.9  ,1.9);
\end{tikzpicture}
\end{figure}
The inclusion relations between the standard parahoric subgroups
of $\GSp(4,F)$ are indicated in Fig.~\ref{parahoricfig}.

\subsection*{The Iwahori-Hecke algebra}
Next, we recall  the Iwahori-Hecke algebra \label{IwahoriHeckedef}
of $\GSp(4,F)$ associated to the Iwahori subgroup $I$. 
As a complex vector space, the \emph{Iwahori-Hecke algebra}
\index{Iwahori-Hecke algebra}
 $\mathcal{H}(\GSp(4,F),I)$
is the set of all compactly supported left and right $I$-invariant
complex-valued functions on $\GSp(4,F)$ 
with addition defined by pointwise addition of functions. 
If $T,T' \in  \mathcal{H}(\GSp(4,F),I)$, then we define the product $T\cdot T'$ 
by
$$
 (T\cdot T')(x)=\int\limits_{\GSp(4,F)}T(xy^{-1})T'(y)\,d^Iy
$$
for $x \in \GSp(4,F)$. 
Here, $d^Iy$
is the Haar measure on $\GSp(4,F)$ which gives $I$ volume~$1$. \label{Imeasure}
If $T,T' \in \mathcal{H}(\GSp(4,F),I)$, then $T\cdot T' \in  \mathcal{H}(\GSp(4,F),I)$, and 
with this definition  $\mathcal{H}(\GSp(4,F),I)$ is a $\C$-algebra. For $g \in \GSp(4,F)$,
we let $T_g \in  \mathcal{H}(\GSp(4,F),I)$ be the characteristic function of $IgI$. 
The identity element of  $\mathcal{H}(\GSp(4,F),I)$ is $T_1$, the characteristic
function of $I\cdot 1 \cdot I=I$. We will write $e=T_1$.
The element $T_{u_1}$
is the characteristic function of $Iu_1 I=u_1I =I u_1$; recall that $u_1$
normalizes $I$. We will write $u_1=T_{u_1}$.
\label{Tu1def} For $i=0,1,$ and $2$, let $e_i=T_{s_i}$,
the characteristic function of $Is_iI$. \label{eielement} By Sect.~5 of \cite{I} (see
also Ex.~25 c) of \S 2 of Chap.~IV of \cite{Bour}), the Iwahori-Hecke
algebra $\mathcal{H}(\GSp(4,F),I)$ is generated by $u_1,e_0,e_1$, and~$e_2$,
and the following relations hold:
\begin{align}
&e_0u_1=u_1e_2,\qquad e_1 u_1=u_1 e_1,\qquad e_2u_1=u_1e_0,\label{IHgeneq1}\\
&e_i^2=(q-1)e_i+qe\qquad \text{for $i=0,1,2$},\label{IHgeneq2}\\
&e_0e_1e_0e_1=e_1e_0e_1e_0,\qquad e_1e_2e_1e_2=e_2e_1e_2e_1,\qquad e_0e_2=e_2e_0.\label{IHgeneq3}
\end{align}
Equation \eqref{IHgeneq2} implies that $e_0$, $e_1$, and $e_2$ are 
invertible. The element $u_1$ is also invertible, with inverse $T_{\varpi^{-1}u_1}$. 
For $w_1,w_2\in W^a$ we  have
\begin{equation}\label{Tw1w2eq}
 T_{w_1w_2}=T_{w_1}\cdot T_{w_2}\qquad\text{if $\ell(w_1w_2)=\ell(w_1)+\ell(w_2)$}.
\end{equation}
Here, $\ell$ is the length function for the Coxeter group $W^a$. 
Also, if $x \in \GSp(4,F)$ normalizes $I$, and $g \in \GSp(4,F)$, then
\begin{equation}
\label{xnormeq}
T_x\cdot T_g=T_{xg}\quad\text{and}\quad\text T_g\cdot T_x = T_{gx}.
\end{equation}
Finally, let $w \in W^a$, and let $IwI=\sqcup_{i\in X} x_i I$
be a disjoint decomposition; we define 
$q_w=\# X$.
\label{qwdef} Since $\GSp(4,F)$
is unimodular, $q_w$ is also the number $\# Y$, where $IwI=\sqcup_{j \in Y} Iy_i$
is a disjoint decomposition. We have
\begin{equation}\label{qsieq}
 q_{s_i}=q\qquad\text{for $i=0,1,2$};
\end{equation}
for this, see Lemma \ref{sIdecomplemma}. 
For $w_1,w_2 \in W^a$ one has
\begin{equation}\label{qw1w2eq}
 q_{w_1w_2}=q_{w_1}q_{w_2}\qquad\text{if $\ell(w_1w_2)=\ell(w_1)+\ell(w_2)$}.
\end{equation} 
\subsection*{Representations}
Now let  $(\pi,V)$ be a smooth representation of $\GSp(4,F)$.
For $T \in  \mathcal{H}(\GSp(4,F),I)$ and $v \in V$ we define
\begin{equation}
\label{IHactioneq}
Tv=\int\limits_{\GSp(4,F)}T(g)\pi(g)v\,d^Ig. 
\end{equation}
Calculations show that with this definition, $V$ is an  $\mathcal{H}(\GSp(4,F),I)$-module.
Furthermore, if $v \in V^I$, and $T \in  \mathcal{H}(\GSp(4,F),I)$, then $Tv \in V^I$.
Here, $V^I$ is the subspace of $V$ of vectors fixed by $I$. 
\subsection*{Volumes} Using \eqref{qw1w2eq} and the definitions
of the parahoric subgroups\index{standard parahoric subgroup}, we have
\begin{align}
&\vl(\GSp(4,\OF))  = \vl(K_{\{s_1,s_2\}} ) = (1+q)^2(1+q^2),\label{gsp4ovoleq}\\
&\vl(\K{})  = \vl(K_{\{s_0,s_2\}}) = (1+q)^2,\label{paravoleq}\\
&\vl(\Gamma_0(\p))=\vl(\Kl{})=1+q,\label{siegelklingenvoleq}\\
&\vl(I)=1.\label{iwahorivoleq}
\end{align}
Again, we use the Haar measure on $\GSp(4,F)$ that assigns $I$ volume $1$. Differences
in volumes are reflected in vertical positioning in Fig.~\ref{parahoricfig}.

\section{Action of the Iwahori-Hecke algebra}

Let $(\pi,V)$ be an irreducible admissible representation of
$\GSp(4,F)$ with trivial central character,
and assume that $\pi$ is Iwahori-spherical\index{Iwahori-spherical representation}\index{representation!Iwahori-spherical}, i.e., $V^I \neq 0$. 
As mentioned
at the beginning of this chapter, our goal is to calculate
the action of the stable Hecke operators on $V_s(1)$. 
As a consequence of a theorem of Casselman (see Sect.~3.6 of
\cite{Car} for a summary), $\pi$ is an irreducible
subquotient of $\chi_1 \times \chi_2 \rtimes \sigma$ for
some unramified characters $\chi_1$, $\chi_2$, and $\sigma$
of $F^\times$ such that $\chi_1 \chi_2 \sigma^2=1$. Consequently,
as a first step toward our goal, in this section we will
explicitly calculate the action of the Iwahori-Hecke algebra
on $(\chi_1 \times \chi_2 \rtimes \sigma)^I$ for 
unramified characters $\chi_1$, $\chi_2$, and $\sigma$ of $F^\times$.

For the remainder of this section, we fix 
unramified characters $\chi_1,\chi_2,\sigma$ of $F^\times$.
We define $\chi_1 \times \chi_2 \rtimes \sigma$, as in Sect.~\ref{repsec}.
The standard space of $\chi_1 \times \chi_2 \rtimes \sigma$ consists of the locally constant
functions $f: \GSp(4,F) \to \C$ such that 
\begin{equation}
\label{indIeq}
 f(\begin{bsmallmatrix}a&*&*&*\\&b&*&*\\&&cb^{-1}&*\\&&&ca^{-1}\end{bsmallmatrix}g)=\chi_1(a)\chi_2(b)\sigma(c)|a^2b|\,|c|^{-3/2}f(g)
\end{equation}
for all $g\in\GSp(4,F)$ and $a,b,c \in F^\times$. 
The action of $\GSp(4,F)$ on $\chi_1 \times \chi_2 \rtimes \sigma$ will be denoted by $\pi$.
We let $(\chi_1 \times \chi_2 \rtimes \sigma)^I$ be the subspace of $I$-invariant vectors
in $\chi_1 \times \chi_2 \rtimes \sigma$. 

For the  next lemma,
we recall that the  elements of the  Weyl group $W$ are the images of $1,s_1,s_2,s_2s_1, s_1s_2s_1,
s_1s_2,s_1s_2s_1s_2$ and $s_2s_1s_2$. 

\begin{lemma}
\label{Ibasislemma}
The function $(\chi_1 \times \chi_2 \rtimes \sigma)^I \to \C^8$ 
defined by 
$$
f \mapsto 
(f(1),f(s_1),f(s_2),f(s_2s_1),f(s_1s_2s_1),f(s_1s_2),f(s_1s_2s_1s_2),f(s_2s_1s_2))
$$
is an isomorphism of $\C$-vector spaces.
\end{lemma}
\begin{proof}
Let $f \in (\chi_1 \times \chi_2 \rtimes \sigma)^I$, and assume that $f$ maps to zero. Let $g \in \GSp(4,F)$.
The Iwasawa decomposition\index{Iwasawa decomposition} asserts that $\GSp(4,F)=B \GSp(4,\OF)$, where $B$ is the Borel subgroup defined in 
Sect.~\ref{repsec}. By Sect.~\ref{backtheorysec}, we have $\GSp(4,\OF)=\sqcup_{w \in W} IwI$. 
The \index{Iwahori subgroup}Iwahori subgroup $I$ admits a decomposition $I=I_+ T(\OF) I_-$, where $T(\OF)$
is as in Sect.~\ref{backtheorysec},  
\begin{equation}
\label{iwahoriplusminuseq}
I_+ = I \cap 
\begin{bsmallmatrix}
1&\OF&\OF&\OF\\
&1&\OF&\OF\\
&&1&\OF\\
&&&1
\end{bsmallmatrix} \quad\text{and}\quad
I_- = I \cap 
\begin{bsmallmatrix}
1&&&\\
\p&1&&\\
\p&\p&1&\\
\p&\p&\p&1
\end{bsmallmatrix}.
\end{equation}
We may thus write $g=butu'wk$ for some $b \in B$, $u\in I_+$, $t \in T(\OF)$,
$u' \in I_-$, $w \in W$, and $k \in I$. By \eqref{indIeq}, and the assumption
that $f \in (\chi_1 \times \chi_2 \rtimes \sigma)^I $, there exists $c \in \C^\times$
such that $f(g) = f(butu'wk) = cf(u' w)$. Since $w^{-1} u' w \in I$, we have 
$f(g) = c f(u'w) = cf(w w^{-1} u' w) = cf(w)$. By assumption $f(w)=0$; hence, $f(g)=0$.
It follows that our map is injective. To see that the map is surjective, let $w \in W$. 
Define $f_w:\GSp(4,F) \to \C$
in the following way. If $g \notin BIwI$, define $f_w(g)=0$.
If $g \in B IwI$, define
$$
f_w(g)=\chi_1(a)\chi_2(b)\sigma(c)|a^2b|\,|c|^{-3/2} \label{fwdef}
$$
where $g = pk$, with
$$
p= \begin{bsmallmatrix}a&*&*&*\\&b&*&*\\&&cb^{-1}&*\\&&&ca^{-1}\end{bsmallmatrix} \in B
$$
for some $a,b,c \in F^\times$, and $k \in IwI$. Using that $\chi_1,\chi_2$, and $\sigma$
are unramified, it is straightforward to verify that $f_w$ is a well-defined element of
$(\chi_1 \times \chi_2 \rtimes \sigma)^I$. Moreover, it is clear that $f_w$ maps to 
the vector with entry $1$ at the $w$-th position, and zeros elsewhere. It follows that our
map is surjective. 
\end{proof}

\begin{lemma}
\label{sIdecomplemma}
There are disjoint decompositions
\begin{align}
Is_0I&=
\bigsqcup_{x \in \OF/\p} s_0
\begin{bsmallmatrix}
1&&&x\varpi^{-1}\\
&1&&\\
&&1&\\
&&&1
\end{bsmallmatrix}
I,\label{sIdecompeq1}\\
Is_1I&=
\bigsqcup_{x \in \OF/\p} s_1
\begin{bsmallmatrix}
1&&&\\
x&1&&\\
&&1&\\
&&-x&1
\end{bsmallmatrix}
I,\label{sIdecompeq2}\\
Is_2I&=
\bigsqcup_{x \in \OF/\p} s_2
\begin{bsmallmatrix}
1&&&\\
&1&&\\
&x&1&\\
&&&1
\end{bsmallmatrix}
I.\label{sIdecompeq3}
\end{align}
\end{lemma}
\begin{proof}
We first prove \eqref{sIdecompeq2}. It is easy to see
that cosets on the right-hand side of \eqref{sIdecompeq2}
are disjoint. To see that the left-hand side of \eqref{sIdecompeq2}
is contained in the right-hand side, let $k \in I$; 
we need to prove that $ks_1$ is contained in 
 the right-hand side of \eqref{sIdecompeq2}. We have
$$
k=
\begin{bsmallmatrix}
a_1&a_2&*&*\\
a_3&a_4&*&*\\
*&*&*&*\\
*&*&*&*
\end{bsmallmatrix}
$$
for some $a_1,a_2,a_3,a_4 \in \OF$ with $a_1,a_4 \in \OF^\times$ and $a_3 \in \p$. Hence,
$$
s_1ks_1 =
\begin{bsmallmatrix}
a_4&a_3&*&*\\
a_2&a_1&*&*\\
*&*&*&*\\
*&*&*&*
\end{bsmallmatrix}
\in 
\begin{bsmallmatrix}
\OF&\p&\OF&\OF\\
\OF&\OF&\OF&\OF\\
\p&\p&\OF&\p\\
\p&\p&\OF&\OF
\end{bsmallmatrix}.
$$
Since $a_4 \in \OF^\times$, there exists $x \in \OF$ such that $a_2+xa_4=0$. 
It follows that
$$
k'=
\begin{bsmallmatrix}
1&&&\\
x&1&&\\
&&1&\\
&&-x&1
\end{bsmallmatrix}
s_1ks_1 
\in 
\begin{bsmallmatrix}
\OF&\p&\OF&\OF\\
\p&\OF&\OF&\OF\\
\p&\p&\OF&\p\\
\p&\p&\OF&\OF
\end{bsmallmatrix}.
$$
This inclusion, along with $k' \in \GSp(4,\OF)$, also implies that $k'_{43} \in \p$; hence,
$k' \in I$. We now see that $ks_1$ is contained 
in 
 the right-hand side of \eqref{sIdecompeq2}. The decomposition \eqref{sIdecompeq3}
is proven in a similar fashion. Finally, \eqref{sIdecompeq1} follows from 
\eqref{sIdecompeq3} by conjugation by $u_1$. 
\end{proof}

By Lemma \ref{Ibasislemma} and its proof, the vector space 
$(\chi_1 \times \chi_2 \rtimes \sigma)^I$ has as basis
the functions $f_w$, $w\in W$, where $f_w$ is the unique $I$-invariant function in 
$\chi_1 \times \chi_2 \rtimes \sigma$ with $f_w(w)=1$ and $f_w(w')=0$ for $w'\in W$, $w'\neq w$. It is convenient to order the basis as follows:
\begin{equation}\label{VIbasiseq2}
 f_e,\quad f_1,\quad f_2,\quad f_{21},\quad f_{121},\quad
 f_{12},\quad f_{1212},\quad f_{212},
\end{equation}
where $e$ is the identity element of $W$, and where we have abbreviated $f_1=f_{s_1}$ and so on. Having fixed
this basis, the operators $e_0,e_1,e_2$ and $u$ on $(\chi_1 \times \chi_2 \rtimes \sigma)^I$ become
$8\times8$ matrices. These are given in the following lemma. We use the notation
\begin{equation}\label{satakeeq}
 \alpha=\chi_1(\varpi),\qquad\beta=\chi_2(\varpi),\qquad\gamma=\sigma(\varpi)
\end{equation}
for the Satake parameters.

\begin{lemma}\label{IHbasisopmatprop}
Let $\chi_1$, $\chi_2$, and $\sigma$
be unramified characters of $F^\times$. With respect to the basis (\ref{VIbasiseq2})
 of $(\chi_1\times \chi_2 \rtimes \sigma)^I$, the action of the elements $e_1$ and $e_2$ 
 is given by the following matrices:
 $$
e_1=\begin{bsmallmatrix}
   0& q\\ 1& q-1\\
  && 0& q\\&& 1& q-1\\
  &&&& q-1& 1\\&&&& q& 0\\
  &&&&&& q-1& 1\\&&&&&& q& 0
  \end{bsmallmatrix},\;
e_2=\begin{bsmallmatrix}
   0& 0& q& 0&
    0& 0& 0& 0\\
   0& 0& 0& 0&
    0& q& 0& 0\\
   1& 0& q-1& 0&
    0& 0& 0& 0\\
   0& 0& 0& 0&
    0& 0& 0& q\\
   0& 0& 0& 0&
    0& 0& q& 0\\
   0& 1& 0& 0&
    0& q-1& 0& 0\\
   0& 0& 0& 0&
    1& 0& q-1& 0\\
   0& 0& 0& 1&
    0& 0& 0& q-1
  \end{bsmallmatrix}
 $$
 The action of $u_1$ is given by
 $$
u_1=\begin{bsmallmatrix}
   &&&&&&& \gamma q^{3/2}\\&&&&&& \gamma q^{3/2}\\
   &&&&& \beta\gamma q^{1/2}\\&&&& \beta\gamma q^{1/2}\\
   &&& \alpha\gamma q^{-1/2}\\&& \alpha\gamma q^{-1/2}\\
   & \alpha\beta\gamma q^{-3/2}\\ \alpha\beta\gamma q^{-3/2}
  \end{bsmallmatrix}.
 $$
 The matrix of $e_0$ is given by the matrix of $u_1e_2u_1^{-1}$.
\end{lemma}
\begin{proof}
Let $w,w' \in \{1,s_1,s_2,s_2s_1,s_1s_2s_1,s_1s_2,s_1s_2s_1s_2,s_2s_1s_2\}$. 
By \eqref{IHactioneq} and \eqref{sIdecompeq2}, we have
\begin{align*}
(e_1f_w)(w')
&=\int\limits_{Is_1I}f_w(w'g)\,d^Ig\\
&=\sum_{x \in \OF/\p} f_w(w' s_1 
\begin{bsmallmatrix}
1&&&\\
x&1&&\\
&&1&\\
&&-x&1
\end{bsmallmatrix}).
\end{align*}
Calculations now verify that the matrix of $e_1$ is as stated;
for this, it is useful to employ \eqref{usefuleq}. The formulas for 
$e_2$ and $u_1$ have similar proofs. The final claim follows
from $e_0=u_1e_2u_1^{-1}$ in the Iwahori-Hecke algebra.
\end{proof}
\subsection*{Projections and bases}
Let $K$ be a parahoric subgroup\index{standard parahoric subgroup} of $\GSp(4,F)$
and let $(\pi,V)$ be a smooth representation\index{representation!smooth}\index{smooth representation} of $\GSp(4,F)$. 
We define 
$$
d_K: V \longrightarrow V^K \label{dKdef}
$$
by
$$
d_K v = \vl(K)^{-1} \int\limits_K \pi(g) v\, d^Ig
$$
for $v \in V$. Again we use the Haar measure on $\GSp(4,F)$
that gives $I$ measure $1$. Evidently, $d_K$ is a projection onto the subspace $V^K$
of $V$. 

\begin{lemma}
\label{basesprojlemma}
Let $\chi_1,\chi_2,$ and $\sigma$ be unramified characters of $F^\times$. 
\begin{enumerate}
\item\label{basesprojlemmaitem1} The projection $d_{\Gamma_0(\p)}:\chi_1\times\chi_2 \rtimes \sigma \to (\chi_1\times\chi_2 \rtimes \sigma)^{\Gamma_0(\p)}$ is given by
\begin{equation}
\label{dgamma0eq}
d_{\Gamma_0(\p)} = (1+q)^{-1}(e+e_1).
\end{equation}
The subspace $(\chi_1 \times \chi_2 \rtimes \sigma)^{\Gamma_0(\p)}$ is four-dimensional and has basis
\begin{equation}\label{VP1basiseq}
 f_e+f_1,\qquad f_2+f_{21},\qquad f_{121}+f_{12},\qquad f_{1212}+f_{212}.
\end{equation}
\item\label{basesprojlemmaitem2} The projection $d_{\Kl{}}:\chi_1\times\chi_2 \rtimes \sigma \to (\chi_1\times\chi_2 \rtimes \sigma)^{\Kl{}}$ is given by
\begin{equation}
\label{dKleq}
d_{\Kl{}} = (1+q)^{-1}(e+e_2).
\end{equation}
The subspace $(\chi_1 \times \chi_2 \rtimes \sigma)^{\Kl{}}$ is four-dimensional and has basis
\begin{equation}\label{VP2basiseq}
 f_e+f_2,\qquad f_1+f_{12},\qquad f_{21}+f_{212},\qquad f_{121}+f_{1212}.
\end{equation}
\item \label{basesprojlemmaitem3} The projection $d_{\K{}}:\chi_1\times\chi_2 \rtimes \sigma \to (\chi_1\times\chi_2 \rtimes \sigma)^{\K{}}$ is given by 
\begin{equation}
\label{dKeq}
d_{\K{}} = (1+q)^{-2}(e+e_0+e_2+e_0e_2).
\end{equation}
The subspace $(\chi_1 \times \chi_2 \rtimes \sigma)^{\K{}}$ is two-dimensional and has basis
\begin{align}\label{VP02basiseq}
 f^{\rm para}_1:&=f_e+f_2+\alpha q^{-2}(f_{121}+f_{1212}),\\
 f^{\rm para}_2:&=f_1+f_{12}+\beta q^{-1}(f_{21}+f_{212}).
\end{align}
\item \label{basesprojlemmaitem4} The projection $d_{\GSp(4,\OF)}:\chi_1\times\chi_2 \rtimes \sigma \to 
(\chi_1\times\chi_2 \rtimes \sigma)^{\GSp(4,\OF)}$ is given by
\begin{align}
d_{\GSp(4,\OF)}&=(1+q)^{-1}(1+q^2)^{-1}\big( e+e_1+e_2+e_1e_2\nonumber \\
&\quad+e_2e_1+e_1e_2e_1+e_2e_1e_2+e_1e_2e_1e_2\big).
\end{align}
The subspace $(\chi_1 \times \chi_2 \rtimes \sigma)^{\GSp(4,\OF)}$ is one-dimensional and has basis
\begin{equation}\label{VP12basiseq}
f_0=f_e+f_1+f_2+f_{21}+f_{121}+f_{12}+f_{1212}+f_{212}.
\end{equation}
\end{enumerate}
\end{lemma}
\begin{proof}
Let $K$ be one of $\Gamma_0(\p)$, $\Kl{}$, $\K{}$, or $\GSp(4,\OF)$. 
The formula for $d_K$ follows from the definition of $d_K$,  \eqref{Tw1w2eq},
and the formula for $\vl(K)$ in Sect.~\ref{backtheorysec}. The subspace 
$(\chi_1 \times \chi_2 \rtimes \sigma)^K$ is exactly the $1$-eigenspace of $d_K|_{(\chi_1 \times \chi_2 \rtimes \sigma)^I}$.
It is straightforward to calculate this eigenspace using the formula for $d_K$,
 the matrices from Lemma \ref{IHbasisopmatprop}, and a computer algebra program. 
\end{proof}

\begin{lemma}
\label{thetabasislemma}
Let $\chi_1$, $\chi_2$, and $\sigma$ be 
unramified characters of $F^\times$. Let $f_0$ be 
the basis element of the one-dimensional space
$(\chi_1 \times \chi_2 \rtimes \sigma)^{\GSp(4,\OF)}$
defined in \eqref{VP12basiseq}. 
Let $\theta,\theta':(\chi_1 \times \chi_2 \rtimes \sigma)^{\GSp(4,\OF)}
\to (\chi_1 \times \chi_2 \rtimes \sigma)^{\K{}}$
be as in \eqref{thetadefeq} and \eqref{thetap1defeq}.
We have 
\begin{align*}
\theta f_0 &= \sigma(\varpi)q^{3/2}(1+\chi_2(\varpi))f_1^{\mathrm{para}}
+\sigma(\varpi)q^{3/2} (1+\chi_1(\varpi)q^{-1})f_2^{\mathrm{para}},\\
\theta' f_0 &=  (\chi_1(\varpi)^{-1}q^2+ q)f_1^{\mathrm{para}}+(\chi_2(\varpi)^{-1}q+ q)f_2^{\mathrm{para}}.
\end{align*}
Assume further that $\chi_1\chi_2 \sigma^2=1$, so that $\chi_1 \times \chi_2 \rtimes \sigma$
has trivial central character. Then $\theta f_0$ and $\theta' f_0$ are linearly
independent if and only if $\chi_1 \chi_2 \neq \nu^{-1}$ and $\chi_1 \chi_2^{-1} \neq \nu^{-1}$. 
\end{lemma}
\begin{proof}
The vector $\theta f_0$ is contained in $(\chi_1 \times \chi_2 \rtimes \sigma)^{\K{}}$.
Since $f_1^{\mathrm{para}}$ and $f_2^{\mathrm{para}}$ form a basis for 
$(\chi_1 \times \chi_2 \rtimes \sigma)^{\K{}}$ by Lemma \ref{basesprojlemma}, 
there exist $a,b \in \C$ such that
$\theta f_0 = a f_1^{\mathrm{para}}+bf_2^{\mathrm{para}}$. Evaluating at $1,s_1 \in \GSp(4,F)$,
we have $a=(\theta f_0)(1)$ and $b=(\theta f_0)(s_1)$. Calculations using the 
definition of $\theta$, along with \eqref{indIeq} show that $a$ and $b$ are
as in the statement of the lemma. The expression for $\theta' f_0$ is similarly
verified. The final assertion follows from the fact that  $f_1^{\mathrm{para}}$
and $f_2^{\mathrm{para}}$ are always linearly independent.
\end{proof}

\section{Stable Hecke operators and the Iwahori-Hecke algebra}
Let $(\pi,V)$ be a smooth representation of $\GSp(4,F)$ for which the center acts trivially. 
In this section we express the endomorphisms $T_{0,1}^s$ and $T_{1,0}^s$ of $V_s(1)$ defined 
in Sect.~\ref{stableheckesec} in terms of certain elements of the Iwahori-Hecke algebra.

\begin{lemma}\label{IgIlemma}
Let $(\pi,V)$ be a smooth representation of $\GSp(4,F)$. If $g \in \GSp(4,F)$ and $v \in V^I$, then 
$(e\circ\pi(g))(v)=q_g^{-1} T_g(v)$, with $q_g$ as in Sect.~\ref{backtheorysec}.
\end{lemma}
\begin{proof}
Let $g \in \GSp(4,F)$. 
We have a bijection
\begin{equation}\label{IgIlemmaeq1}
 I/(I\cap gIg^{-1})\stackrel{\sim}{\longrightarrow}IgI/I
\end{equation}
given by $x(I\cap gIg^{-1})\mapsto xgI$. Let $v\in V^I$. Then
\begin{align*}
(e\circ\pi(g))(v)&=\int\limits_I\pi(h)\pi(g)v\,d^Ih\\
&=\frac1{\#I/(I\cap gIg^{-1})}\sum_{x\in I/(I\cap gIg^{-1})}\pi(x)\pi(g)v\\
&=q_g^{-1}\sum_{x\in I/(I\cap gIg^{-1})}\pi(x)\pi(g)v.
\end{align*}
On the other hand,
\begin{equation}\label{IgIlemmaeq3}
 T_g(v)=\int\limits_{IgI}\pi(h)v\,d^Ih=\sum_{x\in IgI/I}\pi(x)v=\sum_{x\in I/(I\cap gIg^{-1})}\pi(xg)v.
\end{equation}
The assertion follows.
\end{proof}

\begin{lemma}\label{T01Iwahorilemma}
Let $(\pi,V)$ be a smooth representation of $\GSp(4,F)$ for which the center acts trivially. Then, for $v\in V_s(1)$,
\begin{align}
T^s_{0,1}(v)&=(1+q^{-1})d_{\Kl{}}e_2e_1e_2u_1(v)\label{T01Iwahorieq6}\\
&=(1+q^{-1})(1+q)^{-1}(e+e_2)e_2e_1e_2u_1(v).\label{T01Iwahorieq65}
\end{align}
and
 \begin{align}
T^s_{1,0}(v)&=d_{\Kl{}}e_1e_2e_1e_0(v)\label{T01Iwahorieq16}\\
&=(1+q)^{-1}(e+e_2)e_1e_2e_1e_0(v).\label{T01Iwahorieq165}
\end{align}
\end{lemma}
\begin{proof}
Let 
$$
h=\begin{bsmallmatrix}\varpi\\&\varpi\\&&1\\&&&1\end{bsmallmatrix}.
$$
Then
\begin{equation}\label{T01Iwahorieq7}
 h=s_2s_1s_2u_1\begin{bsmallmatrix}1\\&-1\\&&-1\\&&&1\end{bsmallmatrix}.
\end{equation}
Hence
\begin{align}\label{T01Iwahorieq4}
T_h &=T_{s_2s_1s_2u_1}\nonumber\\
&=T_{s_2s_1s_2}T_{u_1}\qquad \text{(see \eqref{xnormeq})} \nonumber\\
&=T_{s_2}T_{s_1}T_{s_2}T_{u_1}\qquad \text{(see \eqref{Tw1w2eq})} \nonumber\\
&=e_2e_1e_2u_1.
\end{align}
Let $v\in V_s(1)=V^{\Kl{}}$. 
By \eqref{Heckecommutelemmaeq20},
\begin{align}
 T^s_{0,1}(v)&=\frac{q^2+q^3}{{\rm vol}(\Kl{})} \int\limits_{\Kl{}}\pi(k)\pi(h)v\,d^I k\nonumber\\
&=(q^2+q^3)d_{\Kl{}}(\pi(h)v).\label{T01Iwahorieq2}
\end{align}
Also, by Lemma \ref{IgIlemma},
$$
e(\pi(h)v)=q_h^{-1}T_h v.
$$
Applying $d_{\Kl{}}$ to this equation, we obtain
$$
d_{\Kl{}}(\pi(h)v)=q_h^{-1}d_{\Kl{}}(T_h v).
$$
Hence, by \eqref{T01Iwahorieq2},
$$
 T^s_{0,1}(v)=(q^2+q^3)q_h^{-1}d_{\Kl{}}(T_h v).
$$
By \eqref{T01Iwahorieq4} this is
\begin{equation}\label{T01Iwahorieq5}
 T^s_{0,1}(v)=(q^2+q^3)q_h^{-1}d_{\Kl{}}e_2e_1e_2u_1(v).
\end{equation}
We have
\begin{align*}
q_h
&=q_{s_2s_1s_2u_1} \qquad \text{(by \eqref{T01Iwahorieq7})}\\
&=q_{s_2s_1s_2}\\
&=q_{s_2}q_{s_1}q_{s_2} \qquad \text{(by \eqref{qw1w2eq})}\\
&=q^3. \qquad\text{(by \eqref{qsieq})}
\end{align*}
Substituting into \eqref{T01Iwahorieq5} now proves \eqref{T01Iwahorieq6};
\eqref{T01Iwahorieq65} follows from \eqref{T01Iwahorieq6} and \eqref{dKleq}.

Similarly, let
\begin{equation}\label{T01Iwahorieq8}
 h'=\begin{bsmallmatrix}\varpi\\&1\\&&1\\&&&\varpi^{-1}\end{bsmallmatrix}=s_1s_2s_1s_0.
\end{equation}
Then
\begin{align}\label{T01Iwahorieq10}
T_{h'}&=T_{s_1s_2s_1s_0}\nonumber\\
 &=T_{s_1}T_{s_2}T_{s_1}T_{s_0}\qquad\text{(see \eqref{Tw1w2eq})}\nonumber\\
 &=e_1e_2e_1e_0.
\end{align}
Let $v\in V_s(1)=V^{P_2}$. 
By \eqref{Heckecommutelemmaeq21}, 
\begin{align}
 T^s_{1,0}(v)&=\frac{q^4}{\vl(\Kl{})}\int\limits_{\Kl{}}\pi(k)\pi(h')v\,d^Ik\nonumber\\
&=q^4 d_{\Kl{}}(\pi(h')v).
 \label{T01Iwahorieq11}
\end{align}
By Lemma \ref{IgIlemma},
$$
e(\pi(h')v)=q_{h'}^{-1}T_{h'} v.
$$
Applying $d_{\Kl{}}$ to this equation yields
\begin{equation}\label{T01Iwahorieq13}
d_{\Kl{}}(\pi(h')v)=q_{h'}^{-1}d_{\Kl{}}(T_{h'} v).
\end{equation}
Therefore, by \eqref{T01Iwahorieq11},
$$
T^s_{1,0}(v)
=q^4 q_{h'}^{-1}d_{\Kl{}}(T_{h'} v).
$$
By \eqref{T01Iwahorieq10} we now have
\begin{equation}
\label{T01Iwahorieq78}
T^s_{1,0}(v)
=q^4 q_{h'}^{-1}d_{\Kl{}}e_1e_2e_1e_0(v).
\end{equation}
Now
\begin{align*}
q_{h'}
 &=q_{s_1s_2s_1s_0} \qquad \text{(by \eqref{T01Iwahorieq7})}\\
 &=q_{s_1}q_{s_2}q_{s_1}q_{s_0}\qquad\text{(by \eqref{qw1w2eq})}\\
 &=q^4.\qquad\text{(by \eqref{qsieq})}
\end{align*}
Substituting into \eqref{T01Iwahorieq78} proves \eqref{T01Iwahorieq16}; also, 
\eqref{T01Iwahorieq165} follows from \eqref{T01Iwahorieq16} and \eqref{dKleq}.
\end{proof}

\section{Characteristic polynomials}
\label{charpolysec}
Let $(\pi,V)$ be an irreducible, admissible representation of 
$\GSp(4,F)$ with trivial central character. 
In this section we calculate the characteristic polynomials of the 
stable Hecke operators\index{stable Hecke operators} $T_{0,1}^s$ and $T_{1,0}^s$ acting on the 
space $V_s(1)$ when $\pi$ is Iwahori-spherical\index{Iwahori-spherical representation}\index{representation!Iwahori-spherical} 
and $V_s(1)$ is non-zero.
By a theorem of Casselman (see Sect.~3.6 of
\cite{Car} for a summary), $\pi$ is Iwahori-spherical if 
and only if $\pi$ is an irreducible
subquotient of $\chi_1 \times \chi_2 \rtimes \sigma$ for
some unramified characters $\chi_1$, $\chi_2$, and $\sigma$
of $F^\times$ such that $\chi_1 \chi_2 \sigma^2=1$. 
It follows that $\pi$ is Iwahori-spherical
if and only if $\pi$ is a group I, II, III, IV, V, or VI
representation with unramified inducing data. 
We begin this section by calculating the characteristic
polynomials of $T_{0,1}^s$ and $T_{1,0}^s$ for three families of representations formed
from unramified characters. These representations,
which may be reducible, occur in the exact sequences
involving group I-VI representations with unramified
inducing data that are described in Sect.~\ref{repsec}. 
With these lemmas in place, we then prove the main result of 
this chapter, Theorem \ref{charpolyiwahoritheorem}.

In this section, if $W$ is a finite-dimensional complex vector space,
and $T:W\to W$ is a linear operator, then the characteristic polynomial 
of $T$ is by definition
$$
p(T,W,X)=\det(X\cdot 1_W -T).
$$

\begin{lemma}
\label{exactcharlemma}
Let 
$$
\begin{CD}
0 @>>> W_1 @>>> W_2 @>>> W_3 @>>> 0\\
@. @VT_1VV @VT_2VV @VT_3VV @.\\
0 @>>> W_1 @>>> W_2 @>>> W_3 @>>> 0\\
\end{CD}
$$
be a commutative diagram with
exact rows, where $W_1$, $W_2$, and
$W_3$ are finite-dimensional complex vector spaces,
and all maps are linear. Then
$$
p(T_2,W_2,X)=p(T_1,W_1,X)p(T_3, W_3,X).
$$
\end{lemma}
\begin{proof}
The proof is left to the reader. 
\end{proof}

\begin{lemma}
\label{genTTlemma}
Let $\chi_1$, $\chi_2$, and $\sigma$ be unramified characters of $F^\times$,
and assume that $\chi_1\chi_2\sigma^2=1$, so that the center of 
$\GSp(4,F)$ acts trivially on $\chi_1 \times \chi_2 \rtimes \sigma$. 
Recall that $(\chi_1 \times \chi_2 \rtimes \sigma)^{\Kl{}}$ is four-dimensional
with basis \eqref{VP2basiseq}. We have

\begin{align*}
&p(T^s_{0,1}, (\chi_1\times \chi_2 \rtimes \sigma)^{\Kl{}},X)=
\left(X-\chi_1(\varpi)\left(1+\chi_2(\varpi))\sigma(\varpi\right) q^{\frac{3}{2}}\right)\\
&\qquad\times  \left(X-\chi_2(\varpi)\left(1+\chi_1(\varpi)\right)\sigma(\varpi) q^{\frac{3}{2}}\right)
\left(X-\left(1+\chi_1(\varpi)\right)\sigma(\varpi) q^{\frac{3}{2}}\right)\\
&\qquad\times \left(X-\left(1+\chi_2(\varpi)\right)\sigma(\varpi)q^{\frac{3}{2}}\right)
\end{align*}
and
\begin{align*}
&p(T^s_{1,0}, (\chi_1\times \chi_2 \rtimes \sigma)^{\Kl{}},X)=
 \left( X - \chi_1(\varpi) q^2\right)
 \left( X - \chi_2(\varpi) q^2 \right)\\
&\qquad\times\left( X - \chi_2(\varpi)^{-1} q^2 \right)
\left( X- \chi_1(\varpi)^{-1} q^2 \right).
\end{align*}
\end{lemma}
\begin{proof}
We work in $(\chi_1 \times \chi_2 \rtimes \sigma)^I$, using the basis \eqref{VIbasiseq2} and the notation
\eqref{satakeeq}. 
By Lemma \ref{basesprojlemma}, the space $(\chi_1 \times \chi_2 \rtimes \sigma)^{\Kl{}}$ is a 
four-dimensional subspace of $(\chi_1 \times \chi_2 \rtimes \sigma)^I$,
with basis \eqref{VP2basiseq}. 
The elements $e_0,e_1,e_2,u_1 \in 
\mathcal{H}(\GSp(4,F),I)$ act on $(\chi_1 \times \chi_2 \rtimes \sigma)^I$,
and have matrices as in Lemma \ref{IHbasisopmatprop} with respect
to the basis \eqref{VIbasiseq2}. 
Using the formulas \eqref{T01Iwahorieq65}
and
\eqref{T01Iwahorieq165}
for 
$T_{0,1}^s$ and $T_{1,0}^s$, respectively, 
in terms of 
of $e_0,e_1,e_2,$ and $u_1$,
it is straightforward to verify that the 
matrix of $T^s_{0,1}$ on $(\chi_1 \times \chi_2 \rtimes \sigma)^{\Kl{}}$ 
with respect to the basis \eqref{VP2basiseq} is

 $$
  \begin{bsmallmatrix}\alpha(1+\beta)\gamma q^{3/2}&0&0&0\\
  \alpha\gamma(q-1)q^{1/2}&\beta(1+\alpha)\gamma q^{3/2}&0&0\\
  \alpha\beta\gamma(q-1)q^{-1/2}&(1+\alpha)\beta\gamma(q-1)q^{1/2}&(1+\alpha)\gamma q^{3/2}&0\\
  \alpha(1+\beta)\gamma(q-1)q^{-1/2}&\alpha\beta\gamma(q-1)q^{1/2}&\alpha\gamma(q-1)q^{1/2}&(1+\beta)\gamma q^{3/2}\end{bsmallmatrix},
 $$
and that the matrix of $T^s_{1,0}$ on $(\chi_1 \times \chi_2 \rtimes \sigma)^{\Kl{}}$
with respect to the basis \eqref{VP2basiseq} is
 $$
  \begin{bsmallmatrix}
  \alpha q^2&0&0&0\\
  \alpha(q-1)q&\beta q^2&0&0\\
  \alpha(q-1)q&(1+\beta)(q-1)q&\beta^{-1}q^2&0\\
  (q-1)(1+\alpha q)&(q-1)(q-1+\beta q)&(q-1)q\beta^{-1}&\alpha^{-1}q^2
  \end{bsmallmatrix}.
 $$
The formulas for the characteristic polynomials are now immediate.
\end{proof}

\begin{lemma}
\label{Sp4stablelemma}
Let $\chi$ and $\sigma$ be characters of $F^\times$. 
\begin{enumerate}
\item \label{Sp4stablelemmaitem1} Assume that $\chi^2 \sigma^2 =1$.  The representation
$\chi 1_{\GSp(4,F)} \rtimes \sigma$ contains a non-zero
vector fixed by $\SSp(4,F)$ if and only if $\chi=\nu^{-\frac{3}{2}}$.
\item \label{Sp4stablelemmaitem2} Assume that $\chi \sigma^2=1$. The representation
$\chi \rtimes \sigma 1_{\GSp(2)}$ contains a non-zero vector
fixed by $\SSp(4,F)$ if and only if $\chi=\nu^{-2}$. 
\end{enumerate}
\end{lemma}
\begin{proof}
\ref{Sp4stablelemmaitem1}. Assume that $\chi 1_{\GL(2)} \rtimes \sigma$ contains a
non-zero vector $f$ fixed by $\SSp(4,F)$. For $g \in \GSp(4,F)$,
we define $g_1=\begin{bsmallmatrix}1&\\ &\lambda(g)^{-1} \end{bsmallmatrix}
g \in \SSp(4,F)$. If $g \in \GSp(4,F)$, then 
\begin{align}
f(g)&=f( \begin{bsmallmatrix}1&\\ &\lambda(g) \end{bsmallmatrix} g_1)\nonumber \\
&=f( \begin{bsmallmatrix}1&\\ &\lambda(g) \end{bsmallmatrix})\nonumber \\
f(g)&=|\lambda(g)|^{-\frac{3}{2}}\sigma (\lambda(g)) f(1). \label{Sp4stablelemmaeq1}
\end{align}
Since $f \neq 0$, it follows that $f(1) \neq 0$. Now let $a \in F^\times$. Then
\begin{align*}
f(
\begin{bsmallmatrix}
a&&&\\
&1&&\\
&&1&\\
&&&a^{-1}
\end{bsmallmatrix})
&=|\det(\begin{bsmallmatrix} a&\\&1 \end{bsmallmatrix})|^{\frac{3}{2}}\chi(\det(\begin{bsmallmatrix} a&\\&1 \end{bsmallmatrix}))f(1)\\
f(1)&=|a|^{\frac{3}{2}}\chi(a) f(1).
\end{align*}
We conclude that $\chi = \nu^{-\frac{3}{2}}$. Now assume  $\chi = \nu^{-\frac{3}{2}}$.
Define $f:\GSp(4,F) \to \C$ by $f(g)=|\lambda(g)|^{-\frac{3}{2}}\sigma (\lambda(g))$.
Calculations shows that $f$ is a 
non-zero $\SSp(4,F)$ invariant element of $\chi 1_{\GL(2)} \rtimes \sigma$. 

The proof of \ref{Sp4stablelemmaitem2} is similar to the proof of \ref{Sp4stablelemmaitem1}.
\end{proof}

\begin{lemma}
\label{typeIIlemma}
Let $\chi$ and $\sigma$ be unramified characters of $F^\times$.
Assume that $\chi^2 \sigma^2=1$, $\chi^2 \neq \nu^{-1}$, and $\chi \neq \nu^{-3/2}$.
The vector spaces $(\chi \St_{\GL(2)} \rtimes \sigma)^{\Kl{}}$ and 
$(\chi 1_{\GL(2)} \rtimes \sigma)^{\Kl{}}$ are both two-dimensional.
We have
\begin{align*}
&p(T^s_{0,1}, (\chi\St_{\GL(2)} \rtimes \sigma)^{\Kl{}},X)
=X^2-\Big(2(\chi\sigma)(\varpi) q+\left(\sigma(\varpi)+\sigma(\varpi)^{-1}\right)q^{\frac{3}{2}}\Big)X\\
&\qquad+q^2+q^3+\left(\chi(\varpi)+\chi(\varpi)^{-1}\right)q^{\frac{5}{2}},\\
&p(T^s_{1,0}, (\chi\St_{\GL(2)} \rtimes \sigma)^{\Kl{}},X)
=X^2-\left(\chi(\varpi)+\chi(\varpi)^{-1}\right) q^{\frac{3}{2}}X+q^3,\\
&p(T^s_{0,1}, (\chi 1_{\GL(2)} \rtimes \sigma)^{\Kl{}},X)
=X^2- \left(2 (\chi\sigma)(\varpi) q^2+\left(\sigma(\varpi)+\sigma(\varpi)^{-1}\right) q^{\frac{3}{2}}\right)X\\
&\qquad+q^3+q^4+\left(\chi(\varpi)+\chi(\varpi)^{-1}\right) q^{\frac{7}{2}},\\
&p(T^s_{1,0}, (\chi 1_{\GL(2)} \rtimes \sigma)^{\Kl{}},X)
=X^2-\left(\chi(\varpi)+\chi(\varpi)^{-1}\right) q^{5/2} X+ q^5.
\end{align*}
\end{lemma}
\begin{proof}
Since there is an exact sequence
 $$
0 \longrightarrow \chi \St_{\GL(2)} \longrightarrow \nu^{1/2} \chi \times \nu^{-1/2} \chi \longrightarrow
\chi 1_{\GL(2)} \longrightarrow 0
$$
of $\GL(2,F)$ representations, 
there is an exact sequence
$$
0 \longrightarrow \chi \St_{\GL(2)}\rtimes \sigma  \longrightarrow (\nu^{1/2} \chi \times \nu^{-1/2} \chi)\rtimes \sigma \longrightarrow
\chi 1_{\GL(2)}\rtimes \sigma \longrightarrow 0
$$
of $\GSp(4,F)$ representations. 
There is an isomorphism $(\nu^{1/2} \chi \times \nu^{-1/2} \chi)\rtimes \sigma \stackrel{\sim}{\longrightarrow}
\nu^{1/2} \chi \times \nu^{-1/2} \chi\rtimes \sigma$ of $\GSp(4,F)$ representations given by $f \mapsto F_f$, where 
$F_f:\GSp(4,F) \to \C$ is defined by
$F_f(g)=\big( f(g) \big)( \begin{bsmallmatrix} 1& \\ &1 \end{bsmallmatrix})$ for $f$ in the standard
model of $(\nu^{1/2} \chi \times \nu^{-1/2} \chi)\rtimes \sigma$ and $g \in \GSp(4,F)$. We thus have an 
exact sequence
\begin{equation}
\label{IIabproofeq1}
0 \longrightarrow \chi \St_{\GL(2)}\rtimes \sigma  \longrightarrow \nu^{1/2} \chi \times \nu^{-1/2} \chi \rtimes \sigma \longrightarrow
\chi 1_{\GL(2)}\rtimes \sigma \longrightarrow 0
\end{equation}
of $\GSp(4,F)$ representations. 
Let $v_0$ be a newform
in the space of $\chi \St_{\GL(2)} \subset \nu^{1/2} \chi \times \nu^{-1/2} \chi$ (see
\cite{RNF} for an account of newforms for representations of $\GL(2,F)$). By (26) of \cite{RNF} we
may assume that, as an element of $\nu^{1/2} \chi \times \nu^{-1/2} \chi$, $v_0(\begin{bsmallmatrix}1&\\&1\end{bsmallmatrix})=q$ and $v_0(\begin{bsmallmatrix}&-1\\1&\end{bsmallmatrix})=-1$.
By the proof of Theorem 5.2.2 of \cite{NF}, the function $F^a_{\mathrm{para}}:\GSp(4,F) \to \chi \St_{\GL(2)}$ 
defined by 
$$
F^a_{\mathrm{para}}(g)=|c^{-1}\det(A) |^{3/2}\sigma(c)(\chi \St_{\GL(2)})(A)v_0
$$
for $g \in \GSp(4,F)$, $g=pk$, with $p=\begin{bsmallmatrix}A&*\\&cA' \end{bsmallmatrix} \in P$, $A \in \GL(2,F)$, $c \in F^\times$,  and $k \in \K{}$, is a non-zero element of 
$(\chi \St_{\GL(2)}\rtimes \sigma)^{\K{}}$. We also write $F^a_{\mathrm{para}}$ for the image of $F^a_{\mathrm{para}}$ in
$\nu^{1/2} \chi \times \nu^{-1/2} \chi \rtimes \sigma$. 
The element $F^a_{\mathrm{para}}$ is contained in $(\nu^{1/2} \chi \times \nu^{-1/2} \chi\rtimes \sigma)^{\K{}}$, which is two-dimensional and spanned by the vectors $f_1^{\mathrm{para}}$ and $f_2^{\mathrm{para}}$
from Lemma \ref{basesprojlemma} (applied to $\nu^{1/2}\chi\times \nu^{-1/2} \chi \rtimes \sigma$).
Calculations show that $F^a_{\mathrm{para}}(1)=q$ and $F^a_{\mathrm{para}}(s_1)=-1$. It follows that 
$F^a_{\mathrm{para}}=qf_1^{\mathrm{para}}-f_2^{\mathrm{para}}$. Further calculations using
Lemma \ref{IHbasisopmatprop}, 
Lemma \ref{T01Iwahorilemma}, and the  assumption $\chi^2 \neq \nu^{-1}$
show that $F^a_{\mathrm{para}}$ and $T_{0,1}^s F^a_{\mathrm{para}}$
are linearly independent. It follows that $(\chi \St_{\GL(2)}\rtimes \sigma)^{\Kl{}}$
is at least two-dimensional. 

Next, there is also an exact sequence
\begin{equation}
\label{IIabproofeq2}
0 \longrightarrow \chi 1_{\GL(2)}\rtimes \sigma   \longrightarrow \nu^{-1/2} \chi \times \nu^{1/2} \chi \rtimes \sigma \longrightarrow
 \chi \St_{\GL(2)}\rtimes \sigma \longrightarrow 0.
\end{equation}
Let $F^b_0$ be a non-zero element of 
$(\chi 1_{\GL(2)}\rtimes \sigma)^{\GSp(4,\OF)}$. We regard $F^b_0$
as an element of $(\nu^{-1/2} \chi \times \nu^{1/2} \chi \rtimes \sigma)^{\GSp(4,\OF)}$
via \eqref{IIabproofeq2}.
The vector $\theta F^b_0$ is 
non-zero by Lemma~\ref{thetabasislemma}
and the assumption $\chi^2 \neq \nu^{-1}$. 
Also, $F^b_0$
and $\theta F^b_0$
are linearly independent by Theorem~\ref{linindparatheorem}; this uses \ref{Sp4stablelemmaitem1} of Lemma \ref{Sp4stablelemma} and
the assumption $\chi \neq \nu^{-\frac{3}{2}}$. Since $F^b_0$
and $\theta F^b_0$ are elements
of  $(\chi 1_{\GL(2)}\rtimes \sigma)^{\Kl{}}$,
it follows that $(\chi 1_{\GL(2)}\rtimes \sigma)^{\Kl{}}$
is at least two-dimensional. Since
$(\chi \St_{\GL(2)}\rtimes \sigma)^{\Kl{}}$
and $(\chi 1_{\GL(2)}\rtimes \sigma)^{\Kl{}}$
are at least two-dimensional, and since
$(\nu^{1/2} \chi \times \nu^{-1/2} \chi \rtimes \sigma)^{\Kl{}}$
is four-dimensional, we conclude that 
$(\chi \St_{\GL(2)}\rtimes \sigma)^{\Kl{}}$
and $(\chi 1_{\GL(2)}\rtimes \sigma)^{\Kl{}}$
are two-dimensional. 

Using Lemma \ref{IHbasisopmatprop} and 
Lemma \ref{T01Iwahorilemma}, it now follows that
the matrix of
$T^s_{0,1}$ on $(\chi \St_{\GL(2)}\rtimes \sigma)^{\Kl{}}$ with respect to the basis 
$F^a_{\mathrm{para}}$, $T^s_{0,1}F^a_{\mathrm{para}}$ is
$$
   \begin{bsmallmatrix}
   0&-(1+q^{1/2}(\chi(\varpi)+\chi(\varpi)^{-1})+q)q^2\\
   1&2(\chi\sigma)(\varpi)q+(\sigma(\varpi)+\sigma(\varpi)^{-1})q^{3/2}
   \end{bsmallmatrix}
  $$
and the matrix of $T^s_{1,0}$ on $(\chi \St_{\GL(2)}\rtimes \sigma)^{\Kl{}}$ with respect to the basis $F^a_{\mathrm{para}}$, $T^s_{0,1}F^a_{\mathrm{para}}$  is
   $$
   \begin{bsmallmatrix}
    -q&\;-q^2(q+1+q^{1/2}(\chi(\varpi)+\chi(\varpi)^{-1}))(\chi\sigma)(\varpi)\\
    (\chi\sigma)(\varpi)&q(1+q^{1/2}(\chi(\varpi)+\chi(\varpi)^{-1}))
    \end{bsmallmatrix}.
   $$
Similarly, 
the matrix of $T^s_{0,1}$ on $(\chi 1_{\GL(2)}\rtimes \sigma)^{\Kl{}}$ with respect to the basis
$F^b_0$, $\theta F^b_0$ is
$$
\begin{bsmallmatrix}
(\chi\sigma)(\varpi)(q+q^2)+(\sigma(\varpi)+\sigma(\varpi)^{-1})q^{3/2}
& q^2+q^3+(\chi\sigma)(\varpi)(\sigma(\varpi)+\sigma(\varpi)^{-1})q^{5/2}\\
-1& (\chi\sigma)(\varpi) (q^2-q)
\end{bsmallmatrix},
$$
and the matrix of $T^s_{1,0}$ on $(\chi 1_{\GL(2)}\rtimes \sigma)^{\Kl{}}$ with respect to the basis
$F^b_0$, $\theta F^b_0$ is
$$
\begin{bsmallmatrix}
q^2+(\chi\sigma)(\varpi)(\sigma(\varpi)+\sigma(\varpi)^{-1})q^{5/2}&
(\chi\sigma)(\varpi)(q^3+q^4)+(\sigma(\varpi)+\sigma(\varpi)^{-1})q^{7/2}\\
-(\chi\sigma)(\varpi)q& -q^2
\end{bsmallmatrix}.
$$
The assertions about characteristic polynomials follow by direct calculations.
\end{proof}

\begin{lemma}
\label{typeIIIlemma}
Let $\chi$ and $\sigma$ be unramified characters of $F^\times$. Assume  that
 ${\chi \sigma^2=1}$, $\chi \neq 1$, and $\chi \neq \nu^{-2}$. 
The vector space $(\chi \rtimes \sigma\St_{\GSp(2)})^{\Kl{}}$ 
is one-dimensional, and the vector space 
$(\chi  \rtimes \sigma 1_{\GSp(2)})^{\Kl{}}$ is  three-dimensional.
We have
\begin{align*}
&p(T^s_{0,1}, (\chi \rtimes \sigma \St_{\GSp(2)})^{\Kl{}},X)=X-q\left(\sigma(\varpi)+\sigma(\varpi)^{-1}\right),\\
&p(T^s_{1,0}, (\chi \rtimes \sigma \St_{\GSp(2)})^{\Kl{}},X)=X-q,\\
&p(T^s_{0,1}, (\chi \rtimes \sigma 1_{\GSp(2)})^{\Kl{}},X)=X^3-\left(\sigma(\varpi)+\sigma(\varpi)^{-1}\right)q(1+2q)X^2\\
&\qquad+\left(1+3q+(\chi(\varpi)+\chi(\varpi)^{-1})q\right)q^2(q+1)X\\
&\qquad-\left(\sigma(\varpi)+\sigma(\varpi)^{-1}\right)q^4(q+1)^2,\\
&p(T^s_{1,0}, (\chi \rtimes \sigma 1_{\GSp(2)})^{\Kl{}},X)=X^3-\left(\left(\chi(\varpi)+\chi(\varpi)^{-1}\right)q^2+q^3\right)X^2\\
&\qquad+\left(\left(\chi(\varpi)+\chi(\varpi)^{-1}\right)q^5+q^4\right)X-q^7.
\end{align*}
\end{lemma}
\begin{proof}
Since there is an exact sequence
$$
0 \longrightarrow \chi \St_{\GSp(2)} \longrightarrow \nu \rtimes \nu^{-1/2} \chi \longrightarrow \chi 1_{\GSp(2)} 
\longrightarrow 0
$$
of representations of $\GSp(2,F) = \GL(2,F)$, there is an exact sequence
$$
0 \longrightarrow \chi \rtimes \sigma \St_{\GSp(2)} \longrightarrow \chi \rtimes (\nu \rtimes \nu^{-1/2}\sigma)
\longrightarrow \chi \rtimes \sigma 1_{\GSp(2)} \longrightarrow 0
$$
of representations of $\GSp(4,F)$. There is an isomorphism 
$\chi \rtimes (\nu \rtimes \nu^{-1/2}\sigma) \stackrel{\sim}{\to} \chi \times \nu \rtimes \nu^{-1/2}\sigma$
of $\GSp(4,F)$ representations given by $f \mapsto F_f$, where $F_f:\GSp(4,F) \to \C$ is defined
by $F_f(g) = \big(f(g)\big)(\begin{bsmallmatrix} 1& \\ &1 \end{bsmallmatrix})$ for $f$ in 
the standard model of $\chi \rtimes (\nu \rtimes \nu^{-1/2}\sigma)$ and $g \in \GSp(4,F)$. We thus
have an exact sequence
\begin{equation}
\label{IIIabeq1}
0 \longrightarrow \chi \rtimes \sigma \St_{\GSp(2)} \longrightarrow \chi \times \nu \rtimes \nu^{-1/2}\sigma
\longrightarrow \chi \rtimes \sigma 1_{\GSp(2)} \longrightarrow 0
\end{equation}
of $\GSp(4,F)$ representations. Let $v_0$ be a newform in the space of 
$\sigma \St_{\GSp(2)} \subset \nu \rtimes \nu^{-1/2}\sigma$
(again, see \cite{RNF} for an account of newforms for representations of $\GL(2,F)$). By (26) of \cite{RNF} we
may assume that, as an element of $\nu \rtimes \nu^{-1/2}\sigma = \nu^{1/2}\sigma \times \nu^{1/2} \sigma$, 
$v_0(\begin{bsmallmatrix}1&\\&1\end{bsmallmatrix})=q$ and $v_0(\begin{bsmallmatrix}&-1\\1&\end{bsmallmatrix})=-1$.
Let 
$$
L_1 = 
\begin{bsmallmatrix}
1&&&\\
\varpi&1&&\\
&&1&\\
&&-\varpi&1
\end{bsmallmatrix},
$$
as on p.~153 of \cite{NF}. By the proof of Theorem 5.4.2 of \cite{NF}, the function $F^a_{\mathrm{para}}:
\GSp(4,F) \to \sigma \St_{\GSp(2)}$ defined by $F^a_{\mathrm{para}}(g)=0$ for $g \notin PL_1\K{2}$ and by 
$$
F^a_{\mathrm{para}}(g) = |y^2 \det(A)^{-1}| \chi(y) (\sigma \St_{\GSp(2)})(A) v_0
$$
for $g = p L_1 k$, $p=\begin{bsmallmatrix} y&*&*\\ &A&* \\ && y^{-1} \det(A) \end{bsmallmatrix}$, $y \in F^\times$,
$A \in \GL(2,F)$, and $k \in \K{2}$, is a well-defined non-zero element in $ (\chi \rtimes \sigma \St_{\GSp(2)})^{\K{2}}$. By definition,
the vector $\sigma_1 F^a_{\mathrm{para}}$ is contained in $ (\chi \rtimes \sigma \St_{\GSp(2)})^{\Kl{}}$. Moreover, a calculation using \eqref{sigmaopseq} and \eqref{usefuleq}
shows that $(\sigma_1 F^a_{\mathrm{para}})(s_1)
\neq 0$, so that $ (\chi \rtimes \sigma \St_{\GSp(2)})^{\Kl{}}$ is at least one-dimensional. 

Next, there is also an exact sequence
\begin{equation}
\label{IIIabproofeq2}
0 \longrightarrow  \chi \rtimes \sigma 1_{\GSp(2)}   \longrightarrow \chi \times \nu^{-1} \rtimes \nu^{1/2}\sigma
\longrightarrow  \chi \rtimes \sigma \St_{\GSp(2)} \longrightarrow 0
\end{equation}
of $\GSp(4,F)$ representations. The representation $\chi \rtimes \sigma 1_{\GSp(2)}$ is unramified
by Table~A.12 of \cite{NF}; let $F_0^b$ be a non-zero element of 
$(\chi \rtimes \sigma 1_{\GSp(2)})^{\GSp(4,\OF)}$. 
We regard $F_0^b$ as an element of $\chi \times \nu^{-1} \rtimes \nu^{1/2}\sigma$
via the inclusion from \eqref{IIIabproofeq2}. By Lemma~\ref{thetabasislemma}
and Theorem \ref{linindparatheorem}, the vectors $F_0^b$, $\theta F_0^b$, and $\theta' F_0^b$
are linearly independent; the application of Theorem \ref{linindparatheorem}
uses \ref{Sp4stablelemmaitem2} of Lemma \ref{Sp4stablelemma} and
the assumption ${\chi \neq \nu^{-2}}$.
It follows that 
$(\chi \rtimes \sigma 1_{\GSp(2)})^{\Kl{}}$  is at least three-dimensional.
Since $(\chi \times \nu^{-1} \rtimes \nu^{1/2}\sigma)^{\Kl{}}$ is four-dimensional by 
Lemma \ref{basesprojlemma}, we see from \eqref{IIIabproofeq2}
that $(\chi \rtimes \sigma \St_{\GSp(2)})^{\Kl{}}$
is one-dimensional and that $(\chi \rtimes \sigma \St_{\GSp(2)})^{\Kl{}}$ is three-dimensional. 
The vectors $F_0^b$, $\theta F_0^b$, and $\theta' F_0^b$ are thus a basis for
$(\chi \rtimes \sigma \St_{\GSp(2)})^{\Kl{}}$; let $W$ be the subspace of 
$\chi \times \nu^{-1} \rtimes \nu^{1/2}\sigma$ spanned by these vectors,
regarded as elements of $\chi \times \nu^{-1} \rtimes \nu^{1/2}\sigma$.
Using Lemma \ref{IHbasisopmatprop} and 
Lemma \ref{T01Iwahorilemma}, it now follows that
the matrix of $T_{0,1}^s$ in the basis $F_0^b$, $\theta F_0^b$, and $\theta' F_0^b$ for
$(\chi \rtimes \sigma \St_{\GSp(2)})^{\Kl{}}$ is
$$
\begin{bsmallmatrix}
 q (q+1) \left(\sigma(\varpi)+\sigma(\varpi)^{-1}\right) & 2 q^2 (q+1) & q^2 (q+1) \left(\sigma(\varpi)+\sigma(\varpi)^{-1}\right) \\
 -1 & q^2 \left(\sigma(\varpi)+\sigma(\varpi)^{-1}\right) & (q-1) q \\
 0 & -q (q+1) & 0 \\
\end{bsmallmatrix},
$$
and the matrix of $T_{1,0}^s$ in the same basis is
$$
\begin{bsmallmatrix}
 q^2 \left(q+\chi(\varpi)+\chi(\varpi)^{-1}\right) & q^3 (q+1) \left(\sigma(\varpi)+\sigma(\varpi)^{-1}\right) & q^3 \left(2 q+\chi(\varpi)+\chi(\varpi)^{-1}\right) \\
 -q \left(\sigma(\varpi)+\sigma(\varpi)^{-1}\right) & -q^2 & -q^2 \left(\sigma(\varpi)+\sigma(\varpi)^{-1}\right) \\
 q & 0 & q^2 \\
\end{bsmallmatrix}.
$$
The formulas for  $p(T^s_{0,1}, (\chi 1_{\GL(2)} \rtimes \sigma)^{\Kl{}},X)$
and $p(T^s_{1,0}, (\chi 1_{\GL(2)} \rtimes \sigma)^{\Kl{}},X)$ are calculated using these matrices.
In view of Lemma \ref{exactcharlemma} and \eqref{IIIabproofeq2},
the formulas for  $p(T^s_{0,1}, (\St_{\GL(2)} \rtimes \sigma)^{\Kl{}},X)$
and $p(T^s_{1,0}, (\St_{\GL(2)} \rtimes \sigma)^{\Kl{}},X)$ are obtained by dividing 
$$
p(T^s_{0,1}, (\chi \times \nu^{-1} \rtimes \nu^{1/2}\sigma)^{\Kl{}},X)
\quad\text{and}\quad
p(T^s_{1,0}, (\chi \times \nu^{-1} \rtimes \nu^{1/2}\sigma)^{\Kl{}},X)
$$
from Lemma \ref{genTTlemma}
by
$$
p(T^s_{0,1}, (\chi 1_{\GL(2)} \rtimes \sigma)^{\Kl{}},X)
\quad\text{and}\quad
p(T^s_{1,0}, (\chi 1_{\GL(2)} \rtimes \sigma)^{\Kl{}},X)
$$
respectively.
\end{proof}

\begin{theorem}
\label{charpolyiwahoritheorem}
Let $(\pi,V)$ be an irreducible, admissible
representation of the group $\GSp(4,F)$ with trivial central
character. The space $V_s(1)$ is non-zero if and only
if $\pi$ is one of the representations in Table \ref{Iwahoritable}.
If $\pi$ is one of the representations in Table \ref{Iwahoritable},
then $\pi$ is paramodular, and the characteristic 
polynomials of $T_{0,1}^s$ and $T_{1,0}^s$ acting
on $V_s(1)$ are as in Table \ref{Iwahoritable}
and Table \ref{eigenIwahoritable}; in the latter table,
$\lambda_\pi$, $\mu_\pi$, and $\varepsilon_\pi$ are the eigenvalues of the paramodular
Hecke operators\index{paramodular Hecke operators!local}\index{Hecke operator!paramodular} $T_{0,1}$ and $T_{1,0}$, and the Atkin-Lehner
operator $\pi(u_{N_\pi})$, respectively,
on a paramodular newform of $\pi$. 
\end{theorem}
\begin{proof}
Assume that $V_s(1)=V^{\Kl{}}$ is non-zero. Then $V^I$
is non-zero, i.e., $V$ is Iwahori-spherical. 
By a theorem of Casselman (see Sect.~3.6 of 
\cite{Car} for a summary), $\pi$ is an irreducible
subquotient of $\chi_1 \times \chi_2 \rtimes \sigma$
for some unramified characters $\chi_1$, $\chi_2$, and
$\sigma$ of $F^\times$ such that $\chi_1\chi_2\sigma^2=1$.
Therefore, $\pi$ is an entry in Table
3 of \cite{SI}; examining this table, and using
the assumption $V_s(1) \neq 0$, we now see that $\pi$ is one of the
representations in Table \ref{Iwahoritable}. 

Conversely,
assume that $\pi$ is one of the representations in 
Table \ref{Iwahoritable}. Then $\pi$ is an entry in 
Table 3 of \cite{SI}, and by inspection of this table
we see that $V_s(1) \neq 0$. 

Assume that $\pi$ is an entry in Table \ref{Iwahoritable}.
Then $\pi$ is paramodular by Table A.12 of~\cite{NF}.

Now let $\pi$ be one of representations in Table \ref{Iwahoritable};
we will prove the characteristic polynomials $p(T_{0,1}^s,V_s(1),X)$
and $p(T_{1,0}^s,V_s(1),X)$ are as in Table \ref{Iwahoritable}. If $\pi$ belongs
to Group I, II or III, then the formulas for the characteristic
polynomials follow from 
Lemma \ref{genTTlemma}, Lemma \ref{typeIIlemma},
and Lemma \ref{typeIIIlemma}.

Assume that $\pi$ belongs to group IVb, so that
$\pi=L(\nu^2,\nu^{-1}\sigma\St_{\GSp(2)})$ for some unramified
character $\sigma$ of $F^\times$ with $\sigma^2=1$.
By \eqref{groupIVtableeq}, there is an exact sequence
$$
0\to
(\sigma \St_{\GSp(4)})^{\Kl{}}
\to
(\nu^2 \rtimes \nu^{-1}\sigma \St_{\GSp(2)})^{\Kl{}}
\to
L(\nu^2, \nu^{-1}\St_{\GSp(2)})^{\Kl{}}
\to 0.
$$
By Table 3 of \cite{SI}, $(\sigma \St_{\GSp(4)})^{\Kl{}}=0$; hence
\begin{align*}
p(T_{0,1}^s,V_s(1),X)=p(T_{0,1}^s,(\nu^2 \rtimes \nu^{-1}\sigma \St_{\GSp(2)})^{\Kl{}},X),\\
p(T_{1,0}^s,V_s(1),X)=p(T_{1,0}^s,(\nu^2 \rtimes \nu^{-1}\sigma \St_{\GSp(2)})^{\Kl{}},X).
\end{align*}
The formulas for $p(T_{0,1}^s,V_s(1),X)$ and $p(T_{1,0}^s,V_s(1),X)$ follow  from 
Lemma~\ref{typeIIIlemma}.

If $\pi$ belongs to group IVc,  then the argument is similar to the IVb case.

If $\pi$ belongs to group IVd, so that $\pi=\sigma 1_{\GSp(4)}$ for some unramified
character $\sigma$ of $F^\times$ with $\sigma^2=1$, then the formulas for 
the characteristic polynomials follow by direct calculations from the 
involved definitions. 

Assume that $\pi$ belongs to group Vb, so that $\pi=L(\nu^{\frac{1}{2}}\xi\St_{\GL(2)},\nu^{-\frac{1}{2}}\sigma)$
with $\sigma$ and $\xi$ unramified characters of $F^\times$, $\xi \neq 1$, $\xi^2=1$, and $\sigma^2=1$. 
By Table 3 of \cite{SI}, $V(1)$ and $V_s(1)$ are both one-dimensional, and $V(0)=0$. It follows that
$V(1)=V_s(1)$ and that the paramodular level $N_\pi$ of 
$\pi$ is $1$; let $v$ be a non-zero element of $V(1)=V_s(1)$. Let $\lambda, \mu \in \C$ be
such that $T_{0,1}v=\lambda v$ and $T_{1,0}v=\mu v$. By Table A.9 of \cite{NF} we have 
$\lambda = \sigma(\varpi)(q^2-1)$ and $\mu=-q(q+1)$. Let $c_{0,1},c_{1,0} \in \C$ be
such that $T_{0,1}^sv = c_{0,1} v$ and $T_{1,0}^sv=c_{1,0}v$. By Lemma \ref{T01Ts01lemma} we
obtain $c_{0,1}=\sigma(\varpi)(q^2-q)$ and $c_{1,0}=-q^2$. The formulas for 
$p(T_{0,1}^s,V_s(1),X)$ and $p(T_{1,0}^s,V_s(1),X)$ now follow immediately.

If $\pi$ belongs to group Vc, then  the computation is similar to the Vb case.

If $\pi$ belongs to group Va or Vd, then  the computation is similar to the IVb case,
using \eqref{groupVtableeq}, Lemma \ref{exactcharlemma},  Lemma \ref{typeIIlemma}, and 
the Vb or Vc case.

If $\pi$ belongs to group VId, then the computation is similar to the IVb case, using
\eqref{groupVItableeq}, Lemma \ref{exactcharlemma}, and Lemma \ref{typeIIlemma} (note
that $\tau(T,\nu^{-\frac{1}{2}}\sigma)^{\Kl{}}=0$ by Table 3 of \cite{SI}).

If $\pi$ belongs to group VIc, then the computation is similar to the Vb case (again,
$V(1)$ and $V_s(1)$ are both one-dimensional by Table 3 of \cite{SI}).

If $\pi$ belongs to group VIa, then the computation is similar to the IVb case, using
\eqref{groupVItableeq}, Lemma \ref{exactcharlemma},  Lemma \ref{typeIIlemma}, and the VIc case;
this completes the verification of the characteristic polynomials in Table \ref{Iwahoritable}.

Finally, using Table A.9 of \cite{NF}, it is
straightforward to verify that the characteristic polynomials
of Table \ref{Iwahoritable} can be written in terms of the 
paramodular eigenvalues $\lambda_\pi$, $\mu_\pi$, and $\varepsilon_\pi$ as
in Table \ref{eigenIwahoritable}.
\end{proof}

%% file: SKMS_part2.tex
\begin{partbacktext}
\part{Siegel Modular Forms}
\end{partbacktext}

%% file: SKMS_chapter10.tex
\chapter{Background on Siegel Modular Forms}
\label{backSMFchap}

The remainder of this text explores applications of the local theory developed in the first part of this work to Siegel modular forms of degree two. 
In this chapter, we recall some essential definitions.  In the next chapter we translate the operators on stable Klingen vectors defined in
previous chapters to modular forms.

\section{Basic definitions}
\label{basicdefsec}

\subsection*{The symplectic similitude group} Let $R$ be a commutative ring with identity $1$.
In this second part of this work we let
\begin{equation}
\label{newJeq}
J=
\begin{bsmallmatrix}
&&1&\\
&&&1\\
-1&&&\\
&-1&&
\end{bsmallmatrix},
\end{equation}
and we define $\GSp(4,R)$
\label{GSp4Rdef} and $\SSp(4,R)$
using $J$. Thus, $\GSp(4,R)$ is defined to be the 
set of $g$ in $\GL(4,R)$ such that $\transpose{g}Jg= \lambda J$ for some $\lambda \in R^\times$.
If $g \in \GSp(4,R)$, then the unit $\lambda \in R^\times$ such that $\transpose{g}Jg=\lambda J$
is unique, and will be denoted by $\lambda(g)$.
\label{sim2def}
The set $\GSp(4,R)$ is a subgroup of
$\GL(4,R)$. If $\begin{bsmallmatrix} A&B \\ C&D \end{bsmallmatrix} \in \GSp(4,R)$,
then 
\begin{equation}
\label{gnewinveq}
g^{-1} = \lambda(g)^{-1}\begin{bsmallmatrix} \transpose{D} & -\transpose{B} \\ -\transpose{C}& \transpose{A} \end{bsmallmatrix}.
\end{equation}
We define $\SSp(4,R)$ to be the subgroup of $g \in \GSp(4,R)$ such that $\lambda (g) =1$, \label{Sp4Rdef}
i.e., $\transpose{g} J g = J$. In Part~1 of this work we defined $\GSp(4)$ and $\SSp(4)$
with respect to 
\begin{equation}
\label{oldJeq}
\begin{bsmallmatrix}
&&&1\\
&&1&\\
&-1&&\\
-1&&&
\end{bsmallmatrix}.
\end{equation}
We will convert between  Part~1 and Part~2 of this work by conjugating by the matrix
\begin{equation}
\label{convmatrixeq}
\begin{bsmallmatrix}
&1&&\\
1&&&\\
&&1&\\
&&&1
\end{bsmallmatrix}.
\end{equation}
This matrix is its own inverse, and
the conjugate of $\GSp(4,R)$ as defined with respect to \eqref{newJeq}
is $\GSp(4,R)$ as defined with respect to \eqref{oldJeq}.
When $R=\R$, we also define $\GSp(4,\R)^+$ as the subgroup of $g \in \GSp(4,\R)$
such that $\lambda(g) >0$. 

\subsection*{The Siegel upper half-space} 
We define $\mathcal{H}_2$
to be the subset of $\Mat(2,\C)$ consisting of the matrices
$Z=X+\I Y$ with $X,Y \in \Mat(2,\R)$ such that $\transpose{X}=X$, $\transpose{Y}=Y$,
and $Y$ is positive-definite. We refer to $\mathcal{H}_2$ as the 
\emph{Siegel upper half-space of degree $2$}\index{Siegel upper half-space}. \label{siegelupperdef}
The set $\mathcal{H}_2$
is a simply connected, open subset of $\Sym(2,\C)$, the $\C$
vector space of $2 \times 2$  symmetric matrices with 
entries from $\C$. The group $\GSp(4,\R)^+$ acts on $\mathcal{H}_2$
via the formula
$$
g\langle Z \rangle =(AZ+B)(CZ+D)^{-1} \label{GH2action}
$$
for $g = \begin{bsmallmatrix} A&B \\ C&D \end{bsmallmatrix} \in \GSp(4,\R)^+$
and $Z \in \mathcal{H}_2$; in particular, $CZ+D \in \GL(2,\C)$.
The action of $\GSp(4,\R)^+$ on $\mathcal{H}_2$ is transitive. 
Define $j:\GSp(4,\R)^+ \times \mathcal{H}_2 \to \C^\times$ by 
$
j(g,Z)=\det(CZ+D)
$
for $g = \begin{bsmallmatrix} A&B \\ C&D \end{bsmallmatrix} \in \GSp(4,\R)^+$
and $Z \in \mathcal{H}_2$. 
We have $j(g_1g_2,Z) = j(g_1,g_2\langle Z \rangle) j(g_2,Z)$ for $Z\in \mathcal{H}_2$ and $g_1,g_2 \in \GSp(4,\R)^+$. \label{autofactordef}
Next, let ${F:\mathcal{H}_2 \to \C}$ be a function,
and let $k$ be an integer such that $k>0$. For ${g \in \GSp(4,\R)^+}$ we 
define $F\big |_k g :\mathcal{H}_2 \to \C$
by 
$$
(F\big |_k g)(Z) = \lambda(g)^k j(g,Z)^{-k} F(g\langle Z \rangle) \label{slashactiondef}
$$
for $Z \in \mathcal{H}_2$. If $g_1,g_2 \in \GSp(4,\R)^+$, then 
$(F\big |_k g_1)\big |_k g_2 = F\big |_k g_1g_2$. We will abbreviate
\begin{equation}
\label{Idefeq}
I=\begin{bsmallmatrix}
\I&\\
&\I
\end{bsmallmatrix}.
\end{equation}

\subsection*{Additional notation} If $A$ and $B$ are square matrices of the same 
size, then we define $A[B]=\transpose{B}AB$, \label{transposeBAB} 
and we define $\mathrm{Tr}(A)$ 
to be the trace of $A$. Let $N$ be an integer such that $N >0$
and let $p$ be a prime number. There exist unique integers $n$ and $M$ such 
that $n \geq 0$, $M>0$, $p$ and $M$ are relatively prime, and $N=Mp^n$; we define $v_p(N)=n$.
\label{vpNdef}
Given a prime number in $\Z$ denoted by a letter of the italic roman font, we will denote the  prime ideal
in the ring of integers of the corresponding local field by the same letter in the Fraktur font.
For example, if $p$ is a prime number in $\Z$, then the prime ideal $p \Z_p$ of $\Z_p$ will be 
denoted by $\mathfrak{p}$. 
Let $p$ be a prime of $\Z$, let $n$ be an integer such that $n \geq 0$, and let $(\pi,V)$ be a smooth representation
of $\GSp(4,\Q_p)$ for which the center acts trivially. In this second part, we will write
\begin{equation}\label{pexpKeq}
V(\p^n) = \{ v \in V \mid \pi(k) v = v \text{ for $k \in \mathrm{K}(\p^{n})$} \}
\end{equation}
and 
\begin{equation}\label{pexpKseq}
V_s(\p^n) = \{ v \in V \mid \pi(k) v = v \text{ for $k \in \mathrm{K}_s(\p^{n})$} \}.
\end{equation}
Here, $\mathrm{K}(\p^{n})$
is the local paramodular group defined in \eqref{paradefeq},
and $\mathrm{K}_s(\p^{n})$
is the local stable Klingen subgroup defined in \eqref{Ksndefeq} (with 
the change in the definition of $\GSp(4)$ mentioned at the beginning of this section). 
Previously, in Part~1, the subspaces $V(\p^n)$ and $V_s(\p^n)$ were denoted by 
$V(n)$ and $V_s(n)$, respectively; in this second part we need the notation to also reflect the choice
of prime. 
We will denote the adeles \label{adelesofQ} of $\Q$ by $\A$,
and we denote the finite
adeles \label{finiteadelesofQ} of $\Q$ by $\A_{\mathrm{fin}}$.
\label{additionalnotation}
Let $g \in \GSp(4,\Q)$ and let $v$ be a place of $\Q$. Then $g_v$
will denote the element of $\GSp(4,\A)$ that is $g$ at the place $v$ 
and $1$ at all other places.

\section{Modular forms}
\label{modformsdefsec}

\subsection*{Congruence subgroups}
Let $N$ be an integer such that $N>0$. We define
\begin{align}
\Gamma_0'(N)
&=\SSp(4,\Q)\cap
\begin{bsmallmatrix}
\Z&N\Z&\Z&\Z \vphantom{N_s^{-1}}\\
\Z&\Z&\Z&\Z\vphantom{N_s^{-1}}\\
\Z&N\Z&\Z&\Z\vphantom{N_s^{-1}}\\
N\Z&N\Z&N\Z&\Z\vphantom{N_s^{-1}}
\end{bsmallmatrix},
\label{classKlNeq2}\\
\mathrm{K}(N)
&=\SSp(4,\Q)\cap
\begin{bsmallmatrix}
\Z&N\Z&\Z&\Z\vphantom{N_s^{-1}} \\
\Z&\Z&\Z&N^{-1}\Z\vphantom{N_s^{-1}} \\
\Z&N\Z&\Z&\Z\vphantom{N_s^{-1}} \\
N\Z&N\Z&N\Z&\Z \vphantom{N_s^{-1}}
\end{bsmallmatrix},
\label{classKNeq}
\end{align}
and 
\begin{equation}
\mathrm{K}_s(N)
=\SSp(4,\Q)\cap
\begin{bsmallmatrix}
\Z&N\Z&\Z&\Z\vphantom{N_s^{-1}}\\
\Z&\Z&\Z&N_s^{-1}\Z\vphantom{N_s^{-1}} \\
\Z&N\Z&\Z&\Z\vphantom{N_s^{-1}} \\
N\Z&N\Z&N\Z&\Z\vphantom{N_s^{-1}}
\end{bsmallmatrix},
\label{classKsNeq2}
\end{equation}
where 
\begin{equation}\label{tildeNdefeq}
N_s=N\prod_{p\mid N}\frac1p.
\end{equation}
In \eqref{tildeNdefeq} $p \mid N$ means that $p$ runs over the primes dividing $N$. 
Using \eqref{gnewinveq} it is easy to verify that $\Gamma_0'(N)$, $\mathrm{K}(N)$, 
and $\mathrm{K}_s(N)$ are subgroups of $\SSp(4,\Q)$. The group $\Gamma_0'(N)$ is the  
\emph{Klingen congruence subgroup}\index{Klingen congruence subgroup} of level $N$, the group $\mathrm{K}(N)$ is the
\emph{paramodular congruence subgroup}\index{paramodular congruence subgroup} of level $N$, and the group $\mathrm{K}_s(N)$
is the 
\emph{stable Klingen congruence subgroup}\index{stable Klingen congruence subgroup} of level $N$.
We note that if $N$ is a square-free integer, then $\mathrm{K}_s(N)=\Gamma_0'(N)$. 
The groups $\Gamma_0'(N)$, $\mathrm{K}(N)$, and $\mathrm{K}_s(N)$ are
commensurable with $\SSp(4,\Z)$ (see 6.3 of Chap.~II on p.~126 of \cite{Fr} for the definition
of commensurable).

\subsection*{Siegel modular forms} Let $k$ be an integer such that $k>0$, 
and let $\Gamma$ be a subgroup of $\SSp(4,\Q)$ commensurable with $\SSp(4,\Z)$. 
A \emph{Siegel modular form of weight $k$ with respect to $\Gamma$}\index{Siegel modular form} \label{SMFdef} is a 
holomorphic function $F: \mathcal{H}_2 \to \C$ such that $F\big |_k \gamma = F$ 
for $\gamma \in \Gamma$. If $\Gamma=\mathrm{K}(N)$ for some positive integer $N$, then $F$ is 
often referred to as a \emph{paramodular form}\index{paramodular form};
if $\Gamma=\mathrm{K}_s(N)$ for some positive integer $N$, then we call $F$ a \emph{stable Klingen form}\index{stable Klingen form}.
We denote by $M_k(\Gamma)$
the $\C$ vector space 
of all Siegel modular forms of weight $k$ with respect to $\Gamma$. If $F \in M_k(\Gamma)$,
then we say that $F$ is a \emph{cusp form} if $F$ satisfies 6.9 of Chap.~II on p.~129 of \cite{Fr};
we let $S_k(\Gamma)$
be the subspace of $F \in M_k(\Gamma)$ such that $F$ is a
cusp form. We define two sets of positive semi-definite symmetric matrices,
\begin{align}
 A(N)&=\{\begin{bsmallmatrix}\alpha&\beta\\\beta&\gamma\end{bsmallmatrix} \mid  \alpha, 2\beta, \gamma\in\Z, \, N\mid\gamma,\, \alpha\gamma-\beta^2 \geq 0,\, \alpha \geq 0,\, \gamma\geq 0\},\label{ANdefeq}\\
 B(N)&=\{\begin{bsmallmatrix}\alpha&\beta\\\beta&\gamma\end{bsmallmatrix} \mid \alpha, 2\beta, \gamma\in \Z,\, N_s\mid\gamma,\,\alpha\gamma-\beta^2 \geq 0,\,\alpha \geq 0,\, \gamma \geq 0\}.\label{BNdefeq}
\end{align}
We denote the subset of positive definite elements of $A(N)$ and $B(N)$ by $A(N)^+$ and $B(N)^+$, respectively.
We note that $A(N)\subset B(N)\subset A(1)$.
Define $\Gamma_0(N)_{\pm}$ to be the subgroup of $\begin{bsmallmatrix} a&b \\ c&d \end{bsmallmatrix} 
\in \GL(2,\Z)$ such that $c \equiv 0\Mod{N}$. The group $\Gamma_0(N)_{\pm}$
acts on $A(N)$ and $B(N)$ via the definition 
\begin{equation}
\label{gamma0acteq}
g \cdot S = g S \transpose{g} =S[\transpose{g}]
\end{equation}
for $g \in \Gamma_0(N)_{\pm}$ and $S \in A(N)$ or $S \in B(N)$. This action preserves the subsets $A(N)^+$ and $B(N)^+$. 
If $F\in M_k(\mathrm{K}(N))$, then $F$ has a \emph{Fourier expansion}\index{Fourier expansion} of the form
\begin{equation}\label{Fourierexpeq}
F(Z)=\sum_{S\in A(N)}a(S)\E^{2\pi \I \mathrm{Tr}(SZ)},
\end{equation}
and if $F\in M_k(\mathrm{K}_s(N))$, then $F$ has a Fourier expansion of the form
\begin{equation}\label{NSFourierexpeq}
F(Z)=\sum_{S\in B(N)}a(S)\E^{2\pi \I \mathrm{Tr}(SZ)}.
\end{equation} 
If $F$ is a cusp form, then the sum is over only $A(N)^+$ in \eqref{Fourierexpeq} and over only $B(N)^+$ in \eqref{NSFourierexpeq}. 
Let $F$ be in $M_k(\mathrm{K}(N))$ or $M_k(\mathrm{K}_s(N))$ with Fourier expansion as in  \eqref{Fourierexpeq} or \eqref{NSFourierexpeq}.
Since $F\big |_k \begin{bsmallmatrix} \transpose{g}^{-1}& \\ & g \end{bsmallmatrix} = F$ for 
$g \in \Gamma_0(N)_{\pm}$, it follows that if $F \in M_k(\mathrm{K}(N))$ (respectively, 
$F \in M_k(\mathrm{K}_s(N))$), then for all $S\in A(N)$ (respectively, $S\in B(N)$) we have  
\begin{equation}
\label{gamma0transruleeq}
a(g \cdot S)  = a ( g S \transpose{g} ) = 
\det(g)^{k} a(S)
\end{equation}
for all $g\in\Gamma_0(N)_{\pm}$.

\subsection*{Fourier-Jacobi expansions}
Let $k$ be an integer such that $k>0$. The elements of 
$M_k(\mathrm{K}(N))$ and $M_k(\mathrm{K}_s(N))$ have Fourier-Jacobi
expansions. To explain this, let $m$ be an integer 
such that $m\geq 0$, let $\mathcal{H}_1$
be the complex
upper half-plane, \label{complexupperhalfplane} and let $f: \mathcal{H}_1 \times \C \to \C$
be a function. For $\lambda, \mu,\kappa \in \R$ define
$$
[\lambda,\mu,\kappa]
=
\begin{bsmallmatrix}
1&&&\mu\\
\lambda&1&\mu&\kappa\\
&&1&-\lambda\\
&&&1
\end{bsmallmatrix}. \label{heisenelementdef}
$$
The elements $[\lambda,\mu,\kappa]$ for $\lambda,\mu,\kappa \in \R$
form a group $H(\R)$ under multiplication of matrices called the 
\emph{Heisenberg group}.\index{Heisenberg group}
\label{Heisenberggroupdef}
The group law is given by 
$$
[\lambda,\mu,\kappa] \cdot [\lambda',\mu',\kappa']
=
[\lambda+\lambda',\mu+\mu',\kappa+\kappa'+\lambda \mu'-\mu \lambda']
$$
for $\lambda,\lambda',\mu,\mu',\kappa,\kappa' \in \R$. 
For $\begin{bsmallmatrix}a&b\\c&d\end{bsmallmatrix}\in\SL(2,\R)$, 
$[\lambda,\mu,\kappa] \in H(\R)$, and $(\tau,z) \in \mathcal{H}_1 \times \C$
we define 
\begin{equation}\label{jacobiformsl2defeq}
(f\big  |_{k,m}\begin{bsmallmatrix}a&b\\c&d\end{bsmallmatrix})(\tau,z)
=(c\tau+d)^{-k}\E^{2\pi \I m (\frac{-cz^2}{c\tau+d})}f(\frac{a\tau+b}{c\tau+d},\frac{z}{c\tau+d})
\end{equation}
and 
\begin{equation}\label{jacobiformhgpdefeq}
(f \big |_{m}[\lambda,\mu,\kappa])(\tau,z)=\E^{2\pi \I m(\lambda^2\tau+2\lambda z+\lambda\mu+\kappa)}f(\tau,z+\lambda\tau+\mu).
\end{equation}
Calculations show that $(f\big |_{k,m} g_1)\big |_{k,m} g_2 = f\big |_{k,m} g_1g_2$ for $g_1,g_2 \in \SL(2,\R)$
and $(f\big |_m h_1 )\big |_m h_2 = f\big |_m h_1h_2$ for $h_1,h_2 \in H(\R)$.  
Also, if $g=\begin{bsmallmatrix} a&b \\ c&d \end{bsmallmatrix} \in \SL(2,\R)$ and $[\lambda,\mu,\kappa] \in H(\R)$,
then we have the following commutation rule\index{Jacobi form!commutation rule}
\begin{equation}
\label{JFcommruleeq}
(f\big |_m [\lambda,\mu,\kappa])\big |_{k,m} \begin{bsmallmatrix} a&b \\ c&d \end{bsmallmatrix}=
(f\big |_{k,m} \begin{bsmallmatrix} a&b \\ c&d \end{bsmallmatrix}) \big |_m [\lambda a + \mu c, \lambda b + \mu d, \kappa].
\end{equation}
We say that $f$ is a \emph{Jacobi form of 
weight $k$ and index $m$ on $\SL(2,\Z)$}\index{Jacobi form} if $f$ is holomorphic, if 
$f\big |_{k,m} g = f$ for $g \in \SL(2,\Z)$, $f \big |_m h = f$ for 
$h \in H(\Z)$, and $f$ has an absolutely convergent Fourier expansion\index{Jacobi form!Fourier expansion} of the form
$$
f(\tau,z) = \sum_{n=0}^\infty \sum_{\substack{r \in \Z\\ r^2 \leq 4nm}} c(n,r) \E^{2\pi \I (n\tau + rz)}. \label{Jacobiformexp}
$$
We say that $f$ is a \emph{cusp form}\index{Jacobi form!cusp form} if $c(n,r)=0$ for $r^2=4nm$.
See \cite{EZ}, p.~1 and p.~9. We denote the space of Jacobi forms of weight $k$ and index $m$ by $J_{k,m}$, and its subspace of cusp forms by $J_{k,m}^{\mathrm{cusp}}$.
If $f$ is a Jacobi form of 
weight $k$ and index $m$ on $\SL(2,\Z)$, then 
\begin{equation}
\label{FJbasiceq}
f(\tau+n_1,z+n_2) = f(\tau,z)\quad\text{and}\quad f(\tau,-z) = (-1)^k f(\tau,z)
\end{equation}
for $n_1,n_2 \in \Z$, $\tau \in \mathcal{H}_1$, and $z \in \C$. 
Let $F\in M_k(\mathrm{K}_s(N))$ or $F\in M_k(\mathrm{K}(N))$. Then $F$ has a \emph{Fourier-Jacobi expansion}\index{Fourier-Jacobi expansion}
\begin{equation}\label{FourierJacobieq}
F(Z)=\sum_{m=0}^\infty f_m(\tau,z)\E^{2\pi \I m\tau'},\qquad Z=\begin{bsmallmatrix}\tau&z\\z&\tau'\end{bsmallmatrix} \in \mathcal{H}_2,
\end{equation}
where $f_m$ is a Jacobi form of weight $k$ and index $m$ on $\SL(2,\Z)$; if $F$ is a cusp form, then $f_0=0$. Suppose that $F\in M_k(\mathrm{K}(N))$. Then it is easy to see that
\begin{equation}\label{SkKNFJvanisheq}
f_m\neq 0\quad\Longrightarrow\quad N \mid m.
\end{equation}
Suppose that $M\in M_k(\mathrm{K}_s(N))$. Then, similarly,
\begin{equation}\label{SkKsNFJvanisheq}
f_m\neq 0\quad\Longrightarrow\quad N_s \mid m.
\end{equation}
Here, $N_s$ is as in \eqref{tildeNdefeq}.

\subsection*{Adelic automorphic forms}
To connect the representation
theory of Part~1 to Siegel modular forms we need to define a certain representation of 
$\GSp(4,\A_{\mathrm{fin}})$. 
Let
$$
K_\infty = \{ \begin{bsmallmatrix} A&B \\ -B&A \end{bsmallmatrix} \in \GL(4,\R) \mid 
\transpose{A} A + \transpose{B} B = 1, \transpose{A} B = \transpose{B} A \}. \label{Kinfdef}
$$
Then $K_\infty$ is a maximal compact subgroup of $\SSp(4,\R)$ and is the stabilizer in $\GSp(4,\R)^+$ of $I$.
Let 
$k$ be an integer such that $k>0$. We define $\mathcal{A}_k$
to be the set of all continuous functions \label{adelicweightk}
$\Phi: \GSp(4,\A) \to \C$ such that 
\begin{enumerate}
\item $\Phi (\rho g) = \Phi (g)$ for all $\rho \in \GSp(4,\Q)$ and $g \in \GSp(4,\A)$;
\item $\Phi (gz) = \Phi (g)$ for all $z \in \A^\times$ and $g \in \GSp(4,\A)$;
\item for some compact, open subgroup $K$ of $\GSp(4,\A_{\mathrm{fin}})$ we have $\Phi (g \kappa) = \Phi(g)$ for all $\kappa \in K$ and $g \in \GSp(4,\A)$;
\item $\Phi (g \kappa) = j(\kappa, I )^{-k} \Phi(g)$
for all $\kappa \in K_\infty$ and $g \in \GSp(4,\A)$;
\item For any $g_{\mathrm{fin}} \in \GSp(4,\A_{\mathrm{fin}})$, the function $\GSp(4,\R)^+ \to \C$ 
defined by $g \mapsto \Phi(g_\mathrm{fin} g)$ is smooth and is annihilated by $\mathfrak{p}_\C^-$ 
(we refer to Sect.~3.5 of \cite{AS} for the definition of $\mathfrak{p}_\C^-$). 
\end{enumerate}
Then $\mathcal{A}_k$ is a complex vector space under addition of functions.
Moreover, it is evident that $\mathcal{A}_k$ is a smooth 
representation of $\GSp(4,\A_{\mathrm{fin}})$ under
the right translation action. 
We let $\mathcal{A}_k^\circ$
be the subspace of  $\Phi \in \mathcal{A}_k$  \label{cuspidaladelicweightk}
that satisfy the following  additional condition:
\begin{enumerate}
\item[(6)] For any proper parabolic subgroup $P$ of $\GSp(4)$ and $g \in \GSp(4,\A)$ we have
$$
\int\limits_{N_P(\Q) \backslash N_P(\A)} \Phi (ng)\, dn = 0;
$$
here, $N_P$ is the unipotent radical of $P$.
\end{enumerate}
We refer to elements of $\mathcal{A}_k^\circ$ as cusp forms. The subspace
$\mathcal{A}_k^\circ$ is closed under the action of $\GSp(4,\A_{\mathrm{fin}})$. 

Next, let $\{K_p\}_p$, where $p$ runs over the primes of $\Z$,
be a family of compact, open subgroups of $\GSp(4,\Q_p)$. 
We will say that $\{K_p\}_p$ is an \emph{admissible family}\index{admissible family of compact, open subgroups}
if  $K_p = \GSp(4,\Z_p)$ for all
but finitely many primes $p$ of $\Z$ and $\lambda (K_p) = \Z_p^\times$ for all primes $p$ of $\Z$. 
Assume that $\{K_p\}_p$ is admissible. 
Define 
\begin{equation}
\label{Kprodefeq}
\mathcal{K} = \prod_{p<\infty} K_p.
\end{equation}
Then $\mathcal{K}$ is a compact, open subgroup of $\GSp(4,\A_{\mathrm{fin}})$. We define
\begin{equation}
\label{gammfromKdefeq}
\Gamma = \GSp(4,\Q) \cap \GSp(4,\R)^+ \mathcal{K}.
\end{equation}
Then $\Gamma$ is a subgroup of $\SSp(4,\Q)$ commensurable with $\SSp(4,\Z)$. 
Also, since strong approximation holds for $\SSp(4)$ (see Satz 2 of \cite{Kn}), since $\Q$
has class number one, and since $\lambda (K_p) = \Z_p^\times$
for all primes $p$ of $\Z$, we have 
$$
\GSp(4,\A) = \GSp(4,\Q) \GSp(4,\R)^+ \mathcal{K}.
$$
We now define 
$$
\mathcal{A}_k(\mathcal{K})
=
\{\Phi \in \mathcal{A}_k\mid \kappa \cdot \Phi = \Phi\text{ for $\kappa\in \mathcal{K}$} \};
$$
here, the action of $\mathcal{K}$ on $\mathcal{A}_k(\mathcal{K})$ is defined by $(\kappa\cdot\Phi)(g)=\Phi(g\kappa)$ for $\kappa\in\mathcal{K}$ and $g\in\GSp(4,\A)$.
We let $\mathcal{A}^\circ_k(\mathcal{K})$
be the subspace of $\mathcal{A}_k(\mathcal{K})$
consisting of cusp forms. 
An element of $\mathcal{A}_k(\mathcal{K})$ is called an \emph{adelic  automorphic form of weight $k$
with respect to $\mathcal{K}$}, \label{adelicweightkwrtK} and an element of  $\mathcal{A}^\circ_k(\mathcal{K})$ is called a \emph{cuspidal adelic  automorphic form of weight $k$
with respect to $\mathcal{K}$}.\index{adelic automorphic form}
\index{cuspidal adelic automorphic form}\index{adelic automorphic form!cuspidal}
The next lemma proves that the vector space $M_k(\Gamma)$ of Siegel modular forms and the vector
space $\mathcal{A}_k(\mathcal{K})$ of adelic automorphic forms are naturally isomorphic. \label{cuspidaladelicweightkwrtK}

\begin{lemma}
\label{twospacesisolemma}
Let $k$ be an integer such that $k>0$, and let $\{K_p\}_p$, where $p$ runs over the primes of $\Z$,
be an admissible family of compact, open subgroups of $\GSp(4,\Q_p)$.
Define $\mathcal{K}$ as in \eqref{Kprodefeq} and define $\Gamma$ as in \eqref{gammfromKdefeq}. 
For $F \in M_k(\Gamma)$, define
$\Phi_F: \GSp(4,\A) \to \C$ by 
\begin{equation}
\label{PhiFdefeq}
\Phi_F (\rho g \kappa) = \lambda(g)^k j (g,I)^{-k}
F( g \langle I \rangle )
\end{equation}
for $\rho \in \GSp(4,\Q)$, $g \in \GSp(4,\R)^+$, and $\kappa \in \mathcal{K}$. 
Then $\Phi_F$ is a well-defined element of $\mathcal{A}_k(\mathcal{K})$, so that there is a linear map
\begin{equation}
\label{twospacesisolemmaeq1}
M_k(\Gamma) \longrightarrow \mathcal{A}_k(\mathcal{K})
\end{equation}
defined by $F \mapsto \Phi_F$. This map sends $S_k(\Gamma)$ into $\mathcal{A}^\circ_k(\mathcal{K})$. Conversely, for $\Phi \in \mathcal{A}_k(\mathcal{K})$ define 
$F_\Phi: \mathcal{H}_2 \to \C$ by 
\begin{equation}
\label{FPhidefeq}
F_\Phi(Z) = \lambda (g)^{-k} j(g, I)^k \Phi(g)
\end{equation}
for $Z \in \mathcal{H}_2$ and $g \in \GSp(4,\R)^+$ such that $g \langle I\rangle
=Z$. Then $F_\Phi$ is a well-defined element of $M_k(\Gamma)$, so that there is a linear
map
\begin{equation}
\label{twospacesisolemmaeq2}
\mathcal{A}_k(\mathcal{K}) \longrightarrow   M_k(\Gamma)
\end{equation}
defined by $\Phi \mapsto F_\Phi$. This map sends $\mathcal{A}^\circ_k(\mathcal{K})$ into $S_k (\Gamma)$. Moreover, the maps in \eqref{twospacesisolemmaeq1} and \eqref{twospacesisolemmaeq2}
are inverses of each other. 
\end{lemma}
\begin{proof}
See the proof of Lemma 4.1 of \cite{JLR} and also \cite{AS}.
\end{proof}

We may apply this theory to the paramodular and stable Klingen settings. 
Let $N$  be an integer such that $N>0$.
Then
$$
\{ \mathrm{K}(\p^{v_p(N)})\}_p \quad \text{and} \quad
\{\mathrm{K}_s(\p^{v_p(N)})\}_p
$$
are admissible families of compact, open subgroups of $\GSp(4,\Q_p)$.
We define
$$
\mathcal{K}(N)=\prod_p \mathrm{K}(\p^{v_p(N)})\quad \text{and} \quad
\mathcal{K}_s(N) = \prod_p \mathrm{K}_s(\p^{v_p(N)}). \label{adelicparadef}
$$
Here, $\mathrm{K}(\p^{v_p(N)})$ is defined as in \eqref{paradefeq}
and $\mathrm{K}_s(\p^{v_p(N)})$ is defined as in \eqref{Ksndefeq};
we note again that in this second part of this work we use
\eqref{newJeq}, while in the first part we used \eqref{oldJeq};
the conversion between formulas from Parts~1 and  2 is achieved by conjugation
by the matrix \eqref{convmatrixeq}. 
It is straightforward to verify that 
$$
\mathrm{K}(N) = \GSp(4,\Q) \cap \GSp(4,\R)^+ \prod_p \mathrm{K}(\p^{v_p(N)})
$$
and 
$$
\mathrm{K}_s(N) = \GSp(4,\Q) \cap \GSp(4,\R)^+ \prod_p \mathrm{K}_s(\p^{v_p(N)}).
$$
Applying Lemma \ref{twospacesisolemma}, we see that 
there are natural isomorphisms
$$
M_k(\mathrm{K}(N)) \stackrel{\sim}{\longleftrightarrow} \mathcal{A}_k(\mathcal{K}(N))
\quad\text{and}\quad
M_k(\mathrm{K}_s(N))  \stackrel{\sim}{\longleftrightarrow}\mathcal{A}_k(\mathcal{K}_s(N)).
$$
Under these isomorphisms cusp forms
are mapped onto cusp forms.

\subsection*{Paramodular old- and newforms} To define paramodular old- and 
newforms we first recall  three level raising operators from \cite{RSA}.
These level raising operators\index{level raising operator} are obtained from the corresponding local operators from 
Sect.~\ref{parasec}.

\begin{lemma}
\label{SFlevelraiselemma}
Let $N$ and $k$ be integers such that $N>0$ and $k>0$, and 
let $p$ be a prime of $\Z$. 
Let $F \in M_k(\mathrm{K}(N))$.  Define 
\begin{align}
\eta_p F&= F\big |_k \begin{bsmallmatrix} 1&&&\\ &p&&\\ &&1& \\ &&&p^{-1} \end{bsmallmatrix},\label{Fetapdefeq}\\
\theta_p F&=F\big |_k 
\begin{bsmallmatrix}
1&&&\\
&1&&\\
&&p^{-1}&\\
&&&p^{-1}
\end{bsmallmatrix}
+
\sum_{x \in \Z/p\Z}
F\big |_k 
\begin{bsmallmatrix}
p^{-1}&&&\\
&1&&\\
&&1&\\
&&&p^{-1}
\end{bsmallmatrix}
\begin{bsmallmatrix}
\vphantom{p^{-1}}1&&x& \\
&1&&\\
&&1&\\
&&&1\vphantom{p^{-1}}
\end{bsmallmatrix},\label{Fthetapdefeq}\\
\theta_p' F&=\eta_p F
+\sum_{c \in \Z/ p\Z}
F\big |_k
\begin{bsmallmatrix}
1&&&\\
&1&&cp^{-1} N^{-1}\\
&&1&\\
&&&1
\end{bsmallmatrix}. \label{Fthetaprimepdefeq}
\end{align}
Then $\eta_p F \in M_k(\mathrm{K}(Np^2))$ and $\theta_p F, \theta_p' F \in M_k(\mathrm{K}(Np))$.
If $F$ is a cusp form, then $\eta_p F$, $\theta_p F$, and $\theta_p' F$ are cusp forms. \label{etathetathetapglobaldef}
\end{lemma}
\begin{proof}
We will prove the lemma for $\theta_pF$; the proofs for $\eta_pF$ and $\theta_p' F$ are similar.
Let $W$ be the subspace of $\mathcal{A}_k$ of all $\Phi \in \mathcal{A}_k$ such that 
$\kappa \cdot \Phi = \Phi$ for all $\kappa \in \prod_{q \neq p} \mathrm{K}(\mathfrak{q}^{v_q(N)})$.
Then $W$ is a smooth representation of $\GSp(4,\Q_p)$. 
Let $\theta: W(\p^{v_p(N)}) \to W (\p^{v_p(N)+1})$ be the level raising operator from 
Sect.~\ref{parasec}; note that we are also using the adjusted notation from \eqref{pexpKeq}.
Since $W(\p^{v_p(N)})=\mathcal{A}_k(\mathcal{K}(N))$ and 
$W(\p^{v_p(N)+1})=\mathcal{A}_k(\mathcal{K}(Np))$, we obtain 
an operator $\theta: \mathcal{A}_k(\mathcal{K}(N)) \to
\mathcal{A}_k(\mathcal{K}(Np))$.
To prove the lemma, it will suffice prove that 
under the composition
$$
M_k(\mathrm{K}(N)) \stackrel{\sim}{\longrightarrow} 
\mathcal{A}_k(\mathcal{K}(N)) \stackrel{\theta}{\longrightarrow} \mathcal{A}_k(\mathcal{K}(Np))
\stackrel{\sim}{\longrightarrow} M_k(\mathrm{K}(Np))
$$
the Siegel modular form $F$ is mapped to $\theta_p F$; here the first and last maps are
the isomorphisms from Lemma \ref{twospacesisolemma}. Define $\Phi=\Phi_F$;
the image of $F$ under the composition is $F_{\theta \Phi}$. 
Let $Z \in \mathcal{H}_2$, and let $g \in \GSp(4,\R)^+$
be such that $g \langle I\rangle = Z$.
Then 
\begin{align}
F_{\theta \Phi}(Z)
&=\lambda(g)^{-k}  j(g, I)^k (\theta \Phi)(g)\nonumber \\
&=\lambda(g)^{-k}  j(g, I)^k \big( \Phi(g 
\begin{bsmallmatrix}
1&&&\\
&1&&\\
&&p&\\
&&&p
\end{bsmallmatrix}_p
)\nonumber \\
&\quad +
\sum_{x \in \Z/p \Z} \Phi (g
\begin{bsmallmatrix}
1&&x&\\
&1&&\\
&&1&\\
&&&1
\end{bsmallmatrix}_p
\begin{bsmallmatrix}
p&&&\\
&1&&\\
&&1&\\
&&&p
\end{bsmallmatrix}_p) \big)\label{SFlevelraiselemmaeq1} \\
&=\lambda(g)^{-k}  j(g, I)^k \big( \Phi(
\underbrace{\begin{bsmallmatrix}
1&&&\\
&1&&\\
&&p&\\
&&&p
\end{bsmallmatrix}^{-1}}_{\text{all places}}
g
\begin{bsmallmatrix}
1&&&\\
&1&&\\
&&p&\\
&&&p
\end{bsmallmatrix}_p
)\nonumber \\
&\quad+
\sum_{x \in \Z/p \Z} \Phi (
\underbrace{\begin{bsmallmatrix}
p&&&\\
&1&&\\
&&1&\\
&&&p
\end{bsmallmatrix}^{-1}
\begin{bsmallmatrix}
1&&x&\\
&1&&\\
&&1&\\
&&&1
\end{bsmallmatrix}^{-1}}_{\text{all places}}
g
\begin{bsmallmatrix}
1&&x&\\
&1&&\\
&&1&\\
&&&1
\end{bsmallmatrix}_p
\begin{bsmallmatrix}
p&&&\\
&1&&\\
&&1&\\
&&&p
\end{bsmallmatrix}_p) \big)\nonumber \\
&=\lambda(g)^{-k}  j(g, I)^k \big( \Phi(
\begin{bsmallmatrix}
1&&&\\
&1&&\\
&&p^{-1}&\\
&&&p^{-1}
\end{bsmallmatrix}_\infty
g
)\nonumber \\
&\quad+
\sum_{x \in \Z/p \Z} \Phi (
\begin{bsmallmatrix}
p^{-1}&&&\\
&1&&\\
&&1&\\
&&&p^{-1}
\end{bsmallmatrix}_\infty
\begin{bsmallmatrix}
1\vphantom{p^{-1}}&&x&\\
&1&&\\
&&1&\\
&&&1\vphantom{p^{-1}}
\end{bsmallmatrix}_\infty
g ) \big)\nonumber \\
&=
p^{k}  F(
\begin{bsmallmatrix}
1&&&\\
&1&&\\
&&p^{-1}&\\
&&&p^{-1}
\end{bsmallmatrix}
 \langle Z \rangle) 
+
\sum_{x \in \Z/p \Z} 
F ((
\begin{bsmallmatrix}
p^{-1}&&&\\
&1&&\\
&&1&\\
&&&p^{-1}
\end{bsmallmatrix}
\begin{bsmallmatrix}
1\vphantom{p^{-1}}&&x&\\
&1&&\\
&&1&\\
&&&1\vphantom{p^{-1}}
\end{bsmallmatrix}
 )\langle Z\rangle)\nonumber  \\
&=(\theta_p F)(Z). \nonumber
\end{align}
We have proved that $\theta_p F \in M_k (\mathrm{K}(Np))$. If $F$ is a cusp form 
then the same argument, with $\mathcal{A}_k^\circ$ replacing $\mathcal{A}_k$, shows
that $\theta_p F$ is a cusp form.
\end{proof}

Next, let $M$ be an integer such that $M>0$. \label{peterssoninnerdef}
Let $\langle \cdot ,\cdot \rangle$
be the 
Petersson inner product\index{Petersson inner product}\index{inner product} on $S_k(\mathrm{K}(M))$ 
from Chap.~2, Theorem 5.3 on p.~89 of \cite{AZ} or Chap.~IV, Hilfsatz~4.10 on p.~271 of
\cite{Fr}. We define the subspace $S_k(\mathrm{K}(M))_{\mathrm{old}}$ of \emph{oldforms}\index{oldforms}
in $S_k(\mathrm{K}(M))$ to be the linear span of the functions $\eta_p F$ where $p$ is a prime such that $p^2 \mid M$ and 
$F \in S_k(\mathrm{K}(Mp^{-2}))$ and the functions $\theta_p F$ and $\theta_p' F$ where $p$ is a prime such that $p \mid M$
and $F \in S_k(\mathrm{K}(Mp^{-1}))$. We define the subspace $S_k(\mathrm{K}(M))_{\mathrm{new}}$ of \emph{newforms}
\index{newforms}
in $S_k(\mathrm{K}(M))$ to be the orthogonal complement
in $S_k(\mathrm{K}(M))$ of the subspace $S_k(\mathrm{K}(M))_{\mathrm{old}}$, so that $S_k(\mathrm{K}(M))_{\mathrm{new}} = 
S_k(\mathrm{K}(M))_{\mathrm{old}}^\perp$. \label{oldnewformsdef}

\subsection*{Paramodular Hecke and Atkin-Lehner operators} Let $N$ and $k$ be integers such that $N>0$ and $k>0$. For each
prime of $\Z$, there exist three natural and important operators  on $M_k(\mathrm{K}(N))$ and 
$S_k(\mathrm{K}(N))$. Let $p$ be a prime of $\Z$ and let $n=v_p(N)$. There are finite disjoint decompositions
$$
\mathrm{K}(N) \begin{bsmallmatrix} 1&&& \\ &1&& \\ &&p& \\ &&&p \end{bsmallmatrix} \mathrm{K}(N)
= \sqcup_i \mathrm{K}(N) h_i, \qquad
\mathrm{K}(N) \begin{bsmallmatrix} p&&& \\ &1&& \\ &&p& \\ &&&p^2 \end{bsmallmatrix} \mathrm{K}(N)
= \sqcup_j \mathrm{K}(N) h_j'.
$$
Define 
\begin{equation}
\label{classicalHeckeoperatoreq}
T(1,1,p,p), T(1,p,p,p^2): M_k(\mathrm{K}(N)) \longrightarrow M_k(\mathrm{K}(N))
\end{equation}
by 
\begin{equation}
\label{classicalHeckeoperatordefeq}
T(1,1,p,p)F = p^{k-3} \sum_i F\big |_k h_i \qquad  T(1,p,p,p^2)F =p^{2(k-3)} \sum_j F\big |_k h_j' 
\end{equation}
for $F \in M_k(\mathrm{K}(N))$. 
It is straightforward to verify that $T(1,1,p,p)$ and 
$T(1,p,p,p^2)$ are well defined and map cusp forms to cusp forms. 
Note that we use the same normalization as in the case $N=1$ (see, for example, (1.3.3) of \cite{AE}).
We refer to $T(1,1,p,p)$ and 
$T(1,p,p,p^2)$ as the \emph{paramodular Hecke operators}.\index{paramodular Hecke operators!global}
Next, since the 
canonical map $\SSp(4,\Z) \to \SSp(4,\Z/N\Z) \cong \SSp(4,\Z/p^{n}\Z) \times \SSp(4,\Z/Np^{-n}\Z)$ is surjective,
there exists $\gamma \in \SSp(4,\Z)$ such that 
$$
\gamma \equiv \begin{bsmallmatrix} &&&-1\\ &&1& \\ &1&& \\ -1&&& \end{bsmallmatrix} \Mod{p^{n}}, \qquad 
\gamma \equiv \begin{bsmallmatrix} 1&&&\vphantom{-1}\\ &1&& \\ &&1& \\\vphantom{-1}  &&&1 \end{bsmallmatrix} \Mod{N p^{-n}}.
$$
Define 
$$
u = \begin{bsmallmatrix} 1&&& \\ &1&& \\ &&p^{n} & \\ &&& p^{n} \end{bsmallmatrix} \gamma.
$$
Then arguments show that $u$ normalizes $\mathrm{K}(N)$ and that $u^2 \in p^{n} \mathrm{K}(N)$. We now 
define
\begin{equation}
\label{ATglobaldefeq}
w_p: M_k(\mathrm{K}(N)) \longrightarrow M_k(\mathrm{K}(N))
\end{equation}
by 
\begin{equation}
\label{ATglobaldefeq2}
w_p F = F\big |_k u.
\end{equation}
We see that $w_p $ is a well-defined involution of $M_k(\mathrm{K}(N))$; it is also evident that 
$w_p $ maps cusp forms to cusp forms. We refer to $w_p $ as a \emph{paramodular Atkin-Lehner operator}.
\index{paramodular Atkin-Lehner operator!global}\index{Atkin-Lehner operator, paramodular!global}

These operators are essentially the operators induced by the corresponding local operators from Sect.~\ref{parasec}. 
To see this, let $W$ be the smooth representation of $\GSp(4,\Q_p)$  from the proof of Lemma \ref{SFlevelraiselemma}.
Let 
$$
T_{0,1},T_{1,0},u_n: W(\p^{n}) \longrightarrow W(\p^{n})
$$
be the operators from \eqref{Heckeopdefeq} and \eqref{ALeq}. Since
$W(\p^{n}) = \mathcal{A}_k(\mathcal{K}(N))$, we obtain operators 
$$
T_{0,1},T_{1,0},u_n: \mathcal{A}_k(\mathcal{K}(N)) \longrightarrow \mathcal{A}_k(\mathcal{K}(N)).
$$
Let $T_{0,1}(p), T_{1,0}(p),$ and $u_{n,p}$  be the respective compositions
\begin{equation}
\label{T01T10unppeq}
M_k(\mathrm{K}(N)) \stackrel{\sim}{\longrightarrow} 
\mathcal{A}_k(\mathcal{K}(N)) \stackrel{T_{0,1},T_{1,0},u_n}{\longrightarrow} \mathcal{A}_k(\mathcal{K}(N))
\stackrel{\sim}{\longrightarrow} M_k(\mathrm{K}(N)),
\end{equation}
where the first and last maps  are
the isomorphisms from Lemma \ref{twospacesisolemma}.

\begin{lemma}
\label{classicallocallemma}
Let $T(1,1,p,p)$, $T(1,p,p,p^2)$, and $w_p $ be as in \eqref{classicalHeckeoperatoreq} and \eqref{ATglobaldefeq}, and let 
$T_{0,1}(p)$, $T_{1,0}(p)$, and $u_{n,p}$ be as in \eqref{T01T10unppeq}. 
Then 
$$
T(1,1,p,p)= p^{k-3} T_{0,1}(p), \qquad T(1,p,p,p^2) = p^{2(k-3)}  T_{1,0}(p), \qquad w_p  = u_{n,p}.
$$
\end{lemma}
\begin{proof}
Using the decompositions in Sect.~6.1 of \cite{NF}, along with some further refinements, one 
can prove that there exist finite disjoint decompositions
$$
\K{n} 
\begin{bsmallmatrix}
p&&&\\
&p&&\\
&&1&\\
&&&1
\end{bsmallmatrix}
\K{n} =\sqcup_i g_i \K{n}, \qquad
\K{n} 
\begin{bsmallmatrix}
p&&&\\
&p^2&&\\
&&p&\\
&&&1
\end{bsmallmatrix}
\K{n} =\sqcup_j g'_j \K{n}
$$
such that $g_i,g_j' \in \GSp(4,\Q)$ and 
\begin{align*}
\mathrm{K}(N)
\begin{bsmallmatrix}
1&&&\\
&1&&\\
&&p&\\
&&&p
\end{bsmallmatrix}\mathrm{K}(N)
&= \sqcup_i \mathrm{K}(N) p g_i^{-1}, \\
\mathrm{K}(N)
\begin{bsmallmatrix}
p&&&\\
&1&&\\
&&p&\\
&&&p^2
\end{bsmallmatrix}\mathrm{K}(N)
&= \sqcup_j \mathrm{K}(N) p^2 g_j'^{-1}.
\end{align*}
Since $\mathrm{K}(N) \subset \mathrm{K}(\mathfrak{q}^{v_q(N)})$ for all primes $q$
of $\Z$, we have $g_i,g_j'  \in \mathrm{K}(\mathfrak{q}^{v_q(N)})$ for all primes $q$ 
of $\Z$ such that $q \neq p$. 
Let $F \in M_k(\mathrm{K}(N))$ and $Z \in \mathcal{H}_2$. Let $h \in \GSp(4,\R)^+$
be such that $Z=h\langle I \rangle$. Then 
\begin{align}
p^{k-3}(T_{0,1}(p) F)(Z)
&=p^{k-3}\lambda(h)^{-k} j(h,I)^{k} T_{0,1} \Phi_F (h)\nonumber \\
&=p^{k-3}\lambda(h)^{-k} j(h,I)^{k} \sum_i \Phi_F (hg_{i,p})\label{classicallocallemmaeq1} \\
&=p^{k-3}\lambda(h)^{-k} j(h,I)^{k} \sum_i \Phi_F (\underbrace{g_i^{-1}}_{\text{all places}}hg_{i,p})\nonumber \\
&=p^{k-3}\lambda(h)^{-k} j(h,I)^{k} \sum_i \Phi_F (g_{i,\infty}^{-1}h)\nonumber \\
&=p^{k-3}(\sum_i F\big |_k g_i^{-1})(Z)\nonumber \\
&=p^{k-3}(\sum_i F\big |_k pg_i^{-1})(Z)\nonumber \\
&= (T(1,1,p,p)F)(Z).\nonumber 
\end{align}
The proofs of the remaining statements are similar. 
\end{proof}

\begin{lemma}
\label{TTUpreservelemma}
Let $p$ be a prime of $\Z$.
Let $T(1,1,p,p)$, $T(1,p,p,p^2)$, and $w_p $ be the endomorphisms 
of $S_k(\mathrm{K}(N))$ as in \eqref{classicalHeckeoperatoreq} and \eqref{ATglobaldefeq}.
Then $T(1,1,p,p)$, $T(1,p,p,p^2)$, and $w_p $ map $S_k(\mathrm{K}(N))_{\mathrm{old}}$ to $S_k(\mathrm{K}(N))_{\mathrm{old}}$
and map $S_k(\mathrm{K}(N))_{\mathrm{new}}$ to $S_k(\mathrm{K}(N))_{\mathrm{new}}$. 
\end{lemma}
\begin{proof}
Let $A$ be one of $T(1,1,p,p)$, $T(1,p,p,p^2)$, and $w_p $. 
We first prove that $S_k(\mathrm{K}(N))_{\mathrm{old}}$ is preserved by $A$.
Let $F \in S_k(\mathrm{K}(N))_{\mathrm{old}}$.
We may assume that there exists a prime $q$ of $\Z$ such that $q \mid N$ and $F$ is of the form $\theta_q F_0$ or
$\theta'_q F_0$ for some $F_0 \in S_k(\mathrm{K}(Nq^{-1}))$, or 
that there exists a prime $q$ of $\Z$ such that $q^2 \mid N$ and $F$ is of the form  $\eta_q F_0$
for some $F_0 \in S_k(\mathrm{K}(Nq^{-2}))$. Assume that $q \neq p$. If $F=\theta_q F_0$,
then $AF=A\theta_q F_0 = \theta_q A' F_0 \in S_k(\mathrm{K}(N))$. Here, in $\theta_q A' F_0$, the map 
$A'$ is the endomorphism of $S_k(\mathrm{K}(Nq^{-1}))$ defined in the same way as $A$; also, to
verify $A\theta_q F_0 = \theta_q A' F_0$ we
use the formulas \eqref{SFlevelraiselemmaeq1} and \eqref{classicallocallemmaeq1}. If
$F=\theta_q' F_0$ or $F=\eta_q F_0$ then a similar argument proves that $AF \in S_k(\mathrm{K}(N))_{\mathrm{old}}$.
Now assume that $q=p$. Then $AF \in S_k(\mathrm{K}(N))_{\mathrm{old}}$ by Lemma \ref{localoldpreservinglemma}; for the connection
to representation theory we use the space $W$ from the proof of Lemma \ref{SFlevelraiselemma}. 
Next we prove that $S_k(\mathrm{K}(N))_{\mathrm{new}}$ is preserved by $A$. Let $\langle \cdot,\cdot \rangle'$
be the inner product on $L^2(\A^\times \GSp(4,\Q) \backslash \GSp(4,\A))$ defined by 
$$
\langle \Phi_1, \Phi_2 \rangle' = \int\limits_{\A^\times \GSp(4,\Q) \backslash \GSp(4,\A) } \Phi_1(g) \overline{\Phi_2 (g)}\, dh
$$
for $\Phi_1,\Phi_2 \in L^2(\A^\times \GSp(4,\Q) \backslash \GSp(4,\A))$
where $dg$ is a fixed right $\GSp(4,\A)$ invariant measure on $\A^\times \GSp(4,\Q) \backslash \GSp(4,\A)$. 
There exists a positive real number $c$ such that if $F_1,F_2 \in S_k(\mathrm{K}(N))$, then 
$\langle F_1, F_2 \rangle = c \langle \Phi_{F_1},\Phi_{F_2} \rangle'$; see Lemma~6 of \cite{AS}. Using this formula, \eqref{Tgformeq}, 
and Lemma \ref{classicallocallemma} one can verify that $\langle AF_1,F_2 \rangle = \langle F_1,AF_2 \rangle$
for $F_1, F_2 \in S_k(\mathrm{K}(N))$, i.e., the operator $A$ is self-adjoint. Since
$A S_k(\mathrm{K}(N))_{\mathrm{old}} \subset S_k(\mathrm{K}(N))_{\mathrm{old}}$, the definition
of $S_k(\mathrm{K}(N))_{\mathrm{new}}$ now implies that $A$ preserves $S_k(\mathrm{K}(N))_{\mathrm{new}}$. 
\end{proof}

%% file: SKMS_chapter11.tex
\chapter{Operators on Siegel Modular Forms}
\label{operchap}

In this chapter we translate the upper block level changing operators\index{upper block operator}
from Chap.~\ref{basicfactschapter} to operators on Siegel modular
forms defined with respect to the stable Klingen congruence subgroups.
We  give a slash formula for each such operator; since this formula
involves only upper block matrices, we are able to calculate the Fourier
and Fourier-Jacobi expansions of the resulting Siegel modular forms.

\section{Overview}\label{overviewsec}
The local level changing operators  from Chap.~\ref{basicfactschapter} 
naturally induce operators on Siegel modular forms. To describe this, 
let $p$ be a prime of $\Z$, and 
let $N'$, $k$, $n_1$, and $n_2$
be integers such that $N'>0$ and $N'$ is relatively prime to $p$, 
${k>0}$, and  ${n_1,n_2 \geq 0}$. Define $N_1=N' p^{n_1}$ and $N_2 = N' p^{n_2}$.
Define $\mathcal{A}_k$ as in Sect.~\ref{modformsdefsec}, 
and let $W$ be the subspace of all $\Phi \in \mathcal{A}_k$ such that $\kappa \cdot \Phi = \Phi$
for all $\kappa \in \prod_{q \neq p} \mathrm{K}(\mathfrak{q}^{v_q(N')})$. Then $W$ is a smooth
representation of $\GSp(4,\Q_p)$. We may consider the subspaces $W_s(\p^{n_1})$
and $W_s(\p^{n_2})$ of stable Klingen vectors\index{stable Klingen vector} as defined in \eqref{Vsndefeq} (with the 
notational adjustment of \eqref{pexpKseq}). 
Let
$$
\delta: W_s(\p^{n_1}) \longrightarrow W_s(\p^{n_2})
$$
be a linear map; a typical example for us will be one of the level changing operators from Chap.~\ref{basicfactschapter}. 
From the involved definitions we have 
$$
W_s(\p^{n_1}) = \mathcal{A}_k(\mathcal{K}_s(N_1))\quad\text{and}\quad W_s(\p^{n_2}) = \mathcal{A}_k(\mathcal{K}_s(N_2)),
$$
so that we may view $\delta$ as map from $\mathcal{A}_k(\mathcal{K}_s(N_1))$ to $\mathcal{A}_k(\mathcal{K}_s(N_2))$:
$$
\delta: \mathcal{A}_k(\mathcal{K}_s(N_1)) \longrightarrow \mathcal{A}_k(\mathcal{K}_s(N_2)). 
$$
By Lemma \ref{twospacesisolemma} we have natural isomorphisms
$$
M_k(\mathrm{K}_s(N_1)) \stackrel{\sim}{\longrightarrow}  \mathcal{A}_k(\mathcal{K}_s(N_1))
\quad \text{and} \quad 
 \mathcal{A}_k(\mathcal{K}_s(N_2))  \stackrel{\sim}{\longrightarrow} M_k(\mathrm{K}_s(N_2))
$$
Therefore, we may compose and obtain a linear map 
\begin{equation}
\label{overviewTddefAeq}
\delta_p: M_k(\mathrm{K}_s(N_1)) \longrightarrow M_k(\mathrm{K}_s(N_2))
\end{equation}
so that the diagram
\begin{equation}
\begin{CD}
\label{overvieweq}
\mathcal{A}_k(\mathcal{K}_s(N_1))  @<\sim<<   M_k(\mathrm{K}_s(N_1)) \\
 @V\delta VV   @VV \delta_p V \\
 \mathcal{A}_k(\mathcal{K}_s(N_2)) @>\sim>>   M_k(\mathrm{K}_s(N_2))
\end{CD}
\end{equation}
commutes. If $F \in M_k(\mathrm{K}_s(N_1))$, then in the notation of Lemma \ref{twospacesisolemma},
\begin{equation}
\label{overviewTddefeq}
( \delta_p F)(Z) = \lambda(h)^{-k} j(h,I)^k (\delta \Phi_F)(h)
\end{equation}
where $Z \in \mathcal{H}_2$ and $h \in \GSp(4,\R)^+$ is such that $h\langle I \rangle =Z$. 
In this construction, if the subspace of cusp forms 
$\mathcal{A}_k^\circ$ is used in place of $\mathcal{A}_k$, 
then the result is a linear map $S_k(\mathrm{K}_s(N_1))  \to S_k(\mathrm{K}_s(N_2))$
that is the restriction of $\delta_p$ to $S_k(\mathrm{K}_s(N_1))$.
Thus, $\delta_p$ maps cusp forms to cusp forms.

In the following sections we will calculate $\delta_p$ for each 
of the upper block level changing  operators $\delta$ from Chap.~\ref{basicfactschapter}.
Before beginning, we note a general feature of the calculation of Fourier expansions.
Let $F \in S_k(\mathrm{K}_s(N_1))$, and let 
$$
F(Z)=\sum_{S\in B(N_1)} a(S)\E^{2\pi \I \mathrm{Tr}(SZ)}
$$
be the Fourier expansion of $F$. In calculating the Fourier expansion of 
$\delta_p F$, we will typically arrive at an expression of the form
$$
(\delta_p F)(Z) = \sum_{S \in X} a(S) \E^{2 \pi \I \mathrm{Tr}( t(S) Z)}.
$$
Here, $X$ is a subset of $B(N_1)$ and $t: X \stackrel{\sim}{\to} Y$
is a bijection, where $Y$ is a subset of $B(N_2)$. We may rewrite this expression as
$$
(\delta_p F)(Z) = \sum_{S \in Y} a(t^{-1}(S)) \E^{2 \pi \I \mathrm{Tr}( S Z)}
$$
and thus obtain the Fourier expansion of $\delta_p F$. We point out that it may not 
be immediately obvious from the expression for $t^{-1}(S)$ 
appearing in the statements of the lemmas below that  this matrix lies in the 
domain $B(N_1)$ for $a(\cdot,F): B(N_1) \to \C$. Nevertheless, it is apparent from 
the method that this is the case.

Let $N$ and $k$ be integers such that $N>0$ and $k>0$, and let $F \in M_k(\mathrm{K}_s(N))$. In this
chapter we will write the Fourier expansion and Fourier-Jacobi expansions of $F$ as in \eqref{NSFourierexpeq}
and \eqref{FourierJacobieq}, respectively, so that 
\begin{align}
F(Z)&=\sum_{S\in B(N)} a(S)\E^{2\pi \I \mathrm{Tr}(SZ)}, \label{NSFourierexpagaineq}\\
F(Z)&=\sum_{\substack{m=0\\N_s \mid m}}^\infty f_m(\tau,z)\E^{2\pi \I m\tau'},\qquad Z
=\begin{bsmallmatrix}\tau&z\\z&\tau'\end{bsmallmatrix}.\label{FourierJacobiagaineq}
\end{align}

\section{Level raising operators}
\index{level raising operator} 
We begin by calculating formulas for the level raising operators $\tau_p$, $\theta_p$, and $\eta_p$ 
obtained from the corresponding local operators defined in Sect.~\ref{levelraisingsec} via the procedure
discussed in Sect.~\ref{overviewsec}. For $\theta_p$ and $\eta_p$ we find that the resulting formulas in
terms of the Fourier-Jacobi expansion involve the well-known operators $V_p$ and $U_p$ on Jacobi forms defined in \cite{EZ}.

\begin{lemma}\label{tauformulalemma}
Let $N$ and $k$ be integers such that $N>0$ and $k>0$, and let $p$ be 
a prime of $\Z$. Define $n = v_p(N)$. Let $\tau: \mathcal{A}_k(\mathcal{K}_s(N)) \to \mathcal{A}_k(\mathcal{K}_s(Np))$
be the linear map obtained as in Sect.~\ref{overviewsec} from the level raising operator \eqref{taundefeq},
and let $\tau_p$ be as in \eqref{overviewTddefAeq}, so that $\tau_p$ is a linear map
$$
\tau_p : M_k(\mathrm{K}_s(N)) \longrightarrow M_k(\mathrm{K}_s(Np)).
$$
Let $F\in M_k(\mathrm{K}_s(N))$ with Fourier expansion as in \eqref{NSFourierexpagaineq}.
Then 
\begin{align}\tau_{p}F
&=p^{-1}\sum_{x\in\Z/p\Z}F\big |_k\begin{bsmallmatrix}1&&&\\&1&&xp^{-n}\\&&1&\\&&&1\end{bsmallmatrix}.\label{taunslashalteq}
\end{align}
If $F$ is a cusp form, then $\tau_{p}F$ is a cusp form.  We have 
\begin{equation}
(\tau_{p}F)(Z)= \sum_{S  \in B(Np) }    a(S) \E^{2\pi \I \mathrm{Tr}(SZ)}.\label{taunfouriereq}
\end{equation} 
The Fourier-Jacobi expansion of $\tau_{p}F$ is given by 
\begin{equation}
(\tau_{p}F)(Z)=\sum_{\substack{m=0\\ N_s\mid m,\ p^n \mid m}}^\infty f_m(\tau,z)\E^{2\pi \I m\tau'},
\qquad Z=\begin{bsmallmatrix}\tau&z\\z&\tau'\end{bsmallmatrix} \in \mathcal{H}_2.\label{taunfourierjacobieq}
\end{equation}
\end{lemma}
\begin{proof}
Define $\Phi=\Phi_F$. Let $Z\in\mathcal{H}_2$ 
and let $h\in \GSp(4,\R)^+$ be such that $Z=h\langle I \rangle$. Then by \eqref{overviewTddefeq}, 
\begin{align*}
(\tau_pF)(Z)
&=\lambda(h)^{-k}j(h,I)^k\tau\Phi(h)\\
&=\lambda(h)^{-k}j(h,I)^k\int\limits_{\Z_p} \Phi(h \begin{bsmallmatrix} 1&&& \\ &1&&xp^{-n} \\ &&1& \\ &&&1 \end{bsmallmatrix} )\, dx\\
&=\lambda(h)^{-k}j(h,I)^k\int\limits_{\Z_p} \Phi( \begin{bsmallmatrix} 1&&& \\ &1&& xp^{-n}\\ &&1& \\ &&&1 \end{bsmallmatrix} h )\, dx\\
&=\lambda(h)^{-k}j(h,I)^k\int\limits_{\Z_p/p^n \Z_p} \int\limits_{p^n \Z_p} \Phi( \begin{bsmallmatrix} 1&&& \\ &1&& (x_1+x_2)p^{-n}\\ &&1& \\ &&&1 \end{bsmallmatrix} h )\, dx_1 \, dx_2\\
&=\lambda(h)^{-k}j(h,I)^kp^{-n}\sum_{x \in \Z/p^n \Z}  \Phi( \begin{bsmallmatrix} 1&&& \\ &1&& xp^{-n}\\ &&1& \\ &&&1 \end{bsmallmatrix}_p h )\\
&=\lambda(h)^{-k}j(h,I)^kp^{-n}\sum_{x \in \Z/p^n \Z}  \Phi( \underbrace{\begin{bsmallmatrix} 1&&& \\ &1&&-x p^{-n} \\ &&1& \\ &&&1 \end{bsmallmatrix}}_{\text{all places}} \begin{bsmallmatrix} 1&&& \\ &1&& xp^{-n}\\ &&1& \\ &&&1 \end{bsmallmatrix}_p h )\\
&=\lambda(h)^{-k}j(h,I)^kp^{-n}\sum_{x \in \Z/p^n \Z}  \Phi( \underbrace{\begin{bsmallmatrix} 1&&& \\ &1&&xp^{-n} \\ &&1& \\ &&&1 \end{bsmallmatrix}}_{\text{all places but $p$}}  h )\\
&=\lambda(h)^{-k}j(h,I)^kp^{-n}\sum_{x \in \Z/p^n \Z}  \Phi( \begin{bsmallmatrix} 1&&& \\ &1&& xp^{-n}\\ &&1& \\ &&&1 \end{bsmallmatrix}_{\infty}  h \underbrace{\begin{bsmallmatrix} 1&&& \\ &1&&xp^{-n} \\ &&1& \\ &&&1 \end{bsmallmatrix}}_{\text{all places but $p, \infty$}}  )\\
&=\lambda(h)^{-k}j(h,I)^kp^{-n}\sum_{x \in \Z/p^n \Z}  \Phi(\begin{bsmallmatrix} 1&&& \\ &1&& xp^{-n}\\ &&1& \\ &&&1 \end{bsmallmatrix}_\infty h   )\\
&=  \big( p^{-n}\sum_{x \in \Z/p^n \Z}  F \big |_k \begin{bsmallmatrix} 1&&& \\ &1&&xp^{-n} \\ &&1& \\ &&&1 \end{bsmallmatrix} \big)  (Z)\\
&=  \big( p^{-1}\sum_{x \in \Z/p \Z}  F \big |_k \begin{bsmallmatrix} 1&&& \\ &1&&xp^{-n} \\ &&1& \\ &&&1 \end{bsmallmatrix} \big)  (Z).
\end{align*}
This proves \eqref{taunslashalteq}. To compute the Fourier expansion, we calculate as follows:
\begin{align*}
(\tau_{p}F)(Z)
&= p^{-n}\sum_{x \in \Z/p^n \Z}  \sum_{S \in B(N)}  a(S) \E^{2\pi \I \mathrm{Tr}\big(S(Z+\begin{bsmallmatrix} & \\ &xp^{-n} \end{bsmallmatrix})\big)} \\
&=   \sum_{S=\begin{bsmallmatrix} \alpha & \beta \\ \beta & \gamma \end{bsmallmatrix} \in B(N)} \big( p^{-n}\sum_{x \in \Z/p^n \Z}  \E^{2\pi \I\gamma x p^{-n}}  \big)   a(S) \E^{2\pi \I \mathrm{Tr}(SZ)} \\
&=   \sum_{\substack{S=\begin{bsmallmatrix} \alpha & \beta \\ \beta & \gamma \end{bsmallmatrix}  \in B(N) \\p^n\mid \gamma}}    a(S) \E^{2\pi \I \mathrm{Tr}(SZ)}\\
&= \sum_{S  \in B(Np)}    a(S) \E^{2\pi \I \mathrm{Tr}(SZ)}.
\end{align*}
This proves  \eqref{taunfouriereq}.
Finally, 
\begin{align*}
(\tau_{p}F)(Z)&=p^{-n}\sum_{x\in\Z/p^n\Z} (F\big |_k\begin{bsmallmatrix}1&&&\\&1&&xp^{-n}\\&&1&\\&&&1\end{bsmallmatrix})(Z)\\
&=p^{-n}\sum_{x\in\Z/p^n\Z}F(Z+\begin{bsmallmatrix}&\\&xp^{-n}\end{bsmallmatrix})\\
&=p^{-n}\sum_{x\in\Z/p^{n}\Z}\sum_{\substack{m=0\\N_s\mid m}}^\infty f_m(\tau,z)\E^{2\pi \I m (\tau'+xp^{-n})}\\
&=p^{-n}\sum_{\substack{m=0\\N_s\mid m}}^\infty f_m(\tau,z)\E^{2\pi \I m \tau'}\sum_{x\in\Z/p^{n}\Z}\E^{2\pi \I m (xp^{-n})}\\
&=\sum_{\substack{m=0\\  N_s\mid m,\ p^n \mid m}}^\infty f_m(\tau,z)\E^{2\pi \I m\tau'}.
\end{align*}
This completes the proof.
\end{proof}	
\begin{lemma}\label{thetaformulalemma}
Let $N$ and $k$ be integers such that $N>0$ and $k>0$, and let $p$ be 
a prime of $\Z$.  Let $\theta: \mathcal{A}_k(\mathcal{K}_s(N)) \to \mathcal{A}_k(\mathcal{K}_s(Np))$
be the linear map obtained as in Sect.~\ref{overviewsec} from the level raising operator \eqref{thetadefeq1},
and let $\theta_p$ be as in \eqref{overviewTddefAeq}, so that $\theta_p$ is a linear map
$$
\theta_p : M_k(\mathrm{K}_s(N)) \longrightarrow M_k(\mathrm{K}_s(Np)).
$$
Let $F \in M_k(\mathrm{K}_s(N))$ with Fourier expansion as in \eqref{NSFourierexpagaineq}.
Then
\begin{equation}\theta_p F=F\big |_k\begin{bsmallmatrix}1&&&\\&1&&\\&&p^{-1}&\\&&&p^{-1}\end{bsmallmatrix}+\sum_{x\in\Z/p\Z}F|_k(\begin{bsmallmatrix}p^{-1}&&&\\&1&&\\&&1&\\&&&p^{-1}\end{bsmallmatrix}
\begin{bsmallmatrix}1&&x&\vphantom{p^{-1}}\\&1&&\\&&1&\\&&&1\vphantom{p^{-1}}\end{bsmallmatrix}).\label{thetaslasheq}
\end{equation} 
If $F$ is a cusp form, then $\theta_pF$ is a cusp form. 
We have  
\begin{align}
(\theta_pF)(Z)
&=\sum_{\substack{S\in B(N)\\p\mid \alpha,\ p \mid 2\beta,\  N_s p\mid\gamma}}p^k  a(p^{-1}S)\E^{2\pi \I \mathrm{Tr} (SZ)}\nonumber\\
&\quad+\sum_{\substack{S\in B(N)\\  N_s p\mid\gamma}} p\,  a(p^{-1}S[\begin{bsmallmatrix}p&\\&1\end{bsmallmatrix}])\E^{2\pi \I \mathrm{Tr}(SZ)}.
\label{thetafouriereq}
\end{align}
The 
Fourier-Jacobi expansion of $\theta_p F$ is  
\begin{equation}
\label{thetafourierjacobieq}
(\theta_{p}F)(Z)=p\sum_{\substack{m=0\\ N_s p\mid m}}^\infty (f_{mp^{-1}}|_{k,mp^{-1}} V_p) (\tau,z)\E^{2\pi \I m\tau'}, 
\qquad Z=\begin{bsmallmatrix}\tau&z\\z&\tau'\end{bsmallmatrix} \in \mathcal{H}_2.
\end{equation}
Here, $V_p$ is the operator from (2) of Chap.~I, Sect.~4 on p.~41 of \cite{EZ}. Explicitly, 
$$
(f_{mp^{-1}}\big |_{k,mp^{-1}} V_p) (\tau,z)
=
p^{k-1} f_{mp^{-1}}(p\tau,pz)+p^{-1}\sum_{x \in \Z/p\Z} f_{mp^{-1}}(p^{-1}(\tau+x),z)
$$
for $\tau \in \mathcal{H}_1$ and $z \in \C$.
\end{lemma}
\begin{proof}
The proof of \eqref{thetaslasheq} is very similar to the proof of Lemma \ref{SFlevelraiselemma} and will be omitted. 
To compute the Fourier expansion we use \eqref{thetaslasheq}:
\begin{align}
&(\theta_p F)(Z)=(F\big |_k\begin{bsmallmatrix}1&&&\\&1&&\\&&p^{-1}&\\&&&p^{-1}\end{bsmallmatrix})(Z)+\sum_{x\in\Z/p\Z}F|_k(\begin{bsmallmatrix}p^{-1}&&&\\&1&&\\&&1&\\&&&p^{-1}\end{bsmallmatrix}
\begin{bsmallmatrix}1&&x&\vphantom{p^{-1}}\\&1&&\\&&1&\\&&&1\vphantom{p^{-1}}\end{bsmallmatrix})(Z)\nonumber \\
&\qquad=p^kF(pZ)+\sum_{x\in\Z/p\Z}F(\begin{bsmallmatrix}p^{-1}&\\&1\end{bsmallmatrix}Z\begin{bsmallmatrix}1\vphantom{p^{-1}}&\\&p\end{bsmallmatrix}
+\begin{bsmallmatrix}p^{-1}x&\vphantom{p^{-1}}\\& \end{bsmallmatrix})\nonumber \\
&\qquad=p^k\sum_{S\in B(N)} a(S)\E^{2\pi \I \mathrm{Tr} (pSZ)}\nonumber \\
&\quad\qquad+\sum_{x\in\Z/p\Z}\sum_{S\in B(N)} a(S)\E^{2\pi \I \mathrm{Tr}(S(\begin{bsmallmatrix}p^{-1}&\\&1\end{bsmallmatrix}Z\begin{bsmallmatrix}1&\\&p\vphantom{p^{-1}}\end{bsmallmatrix}
+\begin{bsmallmatrix}p^{-1}x&\\&\end{bsmallmatrix}))}\nonumber  \\
&\qquad=p^k\sum_{S=\begin{bsmallmatrix} \alpha&\beta\\ \beta&\gamma \end{bsmallmatrix}\in B(N)} a(S)\E^{2\pi \I \mathrm{Tr} (pSZ)}\nonumber \\
&\quad\qquad+\sum_{S=\begin{bsmallmatrix} \alpha&\beta\\ \beta&\gamma \end{bsmallmatrix} \in B(N)} a(S)(\sum_{x\in\Z/p\Z}\E^{2\pi \I p^{-1}x\alpha})
\E^{2\pi \I \mathrm{Tr}(S\begin{bsmallmatrix}p^{-1}&\\&1\end{bsmallmatrix}Z\begin{bsmallmatrix}1&\\&p\end{bsmallmatrix})}\nonumber \\
&\qquad=p^k\sum_{S=\begin{bsmallmatrix} \alpha&\beta\\ \beta&\gamma \end{bsmallmatrix}\in B(N)} a(S)\E^{2\pi \I \mathrm{Tr} (pSZ)}\nonumber \\
&\qquad\quad+p\sum_{\substack{S=\begin{bsmallmatrix} \alpha&\beta\\ \beta&\gamma \end{bsmallmatrix}\in B(N)\\p\mid\alpha}} 
a(S)\E^{2\pi \I \mathrm{Tr}(\begin{bsmallmatrix}1&\\&p\end{bsmallmatrix}S\begin{bsmallmatrix}p^{-1}&\\&1\end{bsmallmatrix}Z)}\nonumber \\
&\qquad=\sum_{\substack{S=\begin{bsmallmatrix} \alpha&\beta\\ \beta&\gamma \end{bsmallmatrix}\in B(N)\\ p\mid \alpha,\ p \mid 2\beta, \  N_s p\mid\gamma}}p^k  
a(p^{-1}S)\E^{2\pi \I \mathrm{Tr} (SZ)}\nonumber \\
&\qquad\quad+\sum_{\substack{S=\begin{bsmallmatrix} \alpha&\beta\\ \beta&\gamma \end{bsmallmatrix}\in B(N)\\  N_s p \mid\gamma}}
p  a(\begin{bsmallmatrix}1&\\&p^{-1}\end{bsmallmatrix}S
\begin{bsmallmatrix}p&\\&1\vphantom{p^{-1}}\end{bsmallmatrix})\E^{2\pi \I \mathrm{Tr}(SZ)}.\nonumber 
\end{align}
This proves \eqref{thetafouriereq}. 

To obtain the Fourier-Jacobi expansion of $\theta_pF$ we proceed from the calculation for \eqref{thetafouriereq}
and use \eqref{FourierJacobiagaineq}:
\begin{align*}
&(\theta_pF)(Z)
=p^kF(pZ)+\sum_{x\in\Z/p\Z}F(\begin{bsmallmatrix}p^{-1}&\\&1\end{bsmallmatrix}Z\begin{bsmallmatrix}1\vphantom{p^{-1}}&\\&p\end{bsmallmatrix}
+\begin{bsmallmatrix}p^{-1}x&\vphantom{p^{-1}}\\& \vphantom{1}\end{bsmallmatrix}) \\
&\qquad=p^k\sum_{\substack{m=0\\ N_s \mid m}}^\infty f_m(p\tau,pz) \E^{2\pi \I mp\tau'} + \sum_{x \in \Z/p\Z} \sum_{\substack{m=0\\ N_s \mid m}}^\infty f_m(p^{-1}\tau+p^{-1}x,z)
\E^{2\pi \I mp\tau'}\\
&\qquad=\sum_{\substack{m=0\\  N_s p \mid m}}^\infty ( p^k f_{mp^{-1}}(p\tau,pz) + 
\sum_{x \in \Z/p\Z}  f_{mp^{-1}}(p^{-1}\tau+p^{-1}x,z) )
\E^{2\pi \I m\tau'}\\
&\qquad=p\sum_{\substack{m=0\\ N_s p \mid m}}^\infty (f_{mp^{-1}}\big |_{k} V_p) (\tau,z)\E^{2\pi \I m\tau'}.
\end{align*}
This completes the proof. 
\end{proof}

\begin{lemma}\label{etaformulalemma}
Let $N$ and $k$ be integers such that $N>0$ and $k>0$, and let $p$ be a prime of $\Z$. 
Let $\eta: \mathcal{A}_k(\mathcal{K}_s(N)) \to \mathcal{A}_k(\mathcal{K}_s(Np^2))$ 
be the linear map obtained as in Sect.~\ref{overviewsec} from the level raising operator \eqref{etadefeq},
and let $\eta_p$ be as in \eqref{overviewTddefAeq}, so that $\eta_p$ is a linear map 
$$
\eta_p: M_k(\mathrm{K}_s(N)) \longrightarrow M_k(\mathrm{K}_s(Np^2)). 
$$
Let $F \in M_k(\mathrm{K}_s(N))$ with Fourier expansion as in \eqref{NSFourierexpagaineq}. Then
\begin{equation}
\eta_{p}F=F\big |_k\begin{bsmallmatrix}1&&&\\&p&&\\&&1&\\&&&p^{-1}\end{bsmallmatrix}.\label{etaslasheq}
\end{equation}
If $F$ is a cusp form, then $\eta_p F$ is a cusp form. We have 
\begin{equation}
(\eta_{p}F)(Z)= \sum_{\substack{S=\begin{bsmallmatrix} \alpha & \beta \\ \beta & \gamma \end{bsmallmatrix}
\in B(Np^2)\\ p\mid 2\beta,\ p^2 \mid \gamma}}
p^k  a(S[\begin{bsmallmatrix}1&\\&p^{-1}\end{bsmallmatrix}])\E^{2\pi \I \mathrm{Tr}(SZ)}.\label{etafouriereq}
\end{equation}
The Fourier-Jacobi expansion of $\eta_p F$ is 
\begin{equation}
(\eta_{p}F)(Z)=\sum_{\substack{m=0\\p^2\mid m,\ N_s\mid m}}^\infty p^k (f_{mp^{-2}}\big |_{k,mp^{-2}} U_p)(\tau,z)\E^{2\pi \I m\tau'}\label{etafourierjacobieq}
\end{equation}
for $Z=\begin{bsmallmatrix}\tau&z\\z&\tau'\end{bsmallmatrix} \in \mathcal{H}_2$. 
Here, $U_p$
is the operator from (1) of Chap.~I, Sect.~4 on p.~41 of \cite{EZ}. Explicitly, we have 
$$
(f_{mp^{-2}}\big |_{k,mp^{-2}} U_p)(\tau,z)=f_{mp^{-2}}(\tau,pz) \label{Updef}
$$
for integers $m$ such that $m \geq 0$, $p^2 \mid m$, $\tau \in \mathcal{H}_1$, and $z \in \C$. 
\end{lemma}	
\begin{proof}
Define $\Phi=\Phi_F$. Let $Z\in\mathcal{H}_2$ 
and let $h\in \GSp(4,\R)^+$ be such that $Z=h\langle I \rangle$. Then by \eqref{overviewTddefeq}, 
\begin{align*}
(\eta_pF)(Z)&=\lambda(h)^{-k}j(h,I)^k\eta\Phi(h)\\
&=\lambda(h)^{-k}j(h,I)^k\Phi(h \begin{bsmallmatrix}1&&&\\&p^{-1}&&\\&&1&\\&&&p\end{bsmallmatrix}_p)\\
&=\lambda(h)^{-k}j(h,I)^k\Phi(\begin{bsmallmatrix}1&&&\\&p^{-1}&&\\&&1&\\&&&p\end{bsmallmatrix}_p h)\\	
&=\lambda(h)^{-k}j(h,I)^k\Phi(\underbrace{\begin{bsmallmatrix}1&&&\\&p&&\\&&1&\\&&&p^{-1}\end{bsmallmatrix}}_{\text{all places}}\begin{bsmallmatrix}1&&&\\&p^{-1}&&\\&&1&\\&&&p\end{bsmallmatrix}_ph)\\
&=\lambda(h)^{-k}j(h,I)^k\Phi(\begin{bsmallmatrix}1&&&\\&p&&\\&&1&\\&&&p^{-1}\end{bsmallmatrix}_{\infty}h)\\
&=(F\big |_k\begin{bsmallmatrix}1&&&\\&p&&\\&&1&\\&&&p^{-1}\end{bsmallmatrix})(Z).	
\end{align*}
To prove  \eqref{etafouriereq} we use \eqref{etaslasheq}:
\begin{align*}
(\eta_p F)(Z)&=(F\big |_k\begin{bsmallmatrix}1&&&\\&p&&\\&&1&\\&&&p^{-1}\end{bsmallmatrix})(Z)\\
&=p^{k}F(\begin{bsmallmatrix}1&\\&p\end{bsmallmatrix}Z\begin{bsmallmatrix}1&\\&p\end{bsmallmatrix})\\
&=\sum_{S\in B(N)}p^k a(S)\E^{2\pi \I \mathrm{Tr}(S\begin{bsmallmatrix}1&\\&p\end{bsmallmatrix}Z\begin{bsmallmatrix}1&\\&p\end{bsmallmatrix})}\\
&=\sum_{S\in B(N)}p^k a(S)\E^{2\pi \I \mathrm{Tr}(\begin{bsmallmatrix}1&\\&p\end{bsmallmatrix}S\begin{bsmallmatrix}1&\\&p\end{bsmallmatrix}Z)}\\
&=\sum_{\substack{S=\begin{bsmallmatrix} \alpha & \beta \\ \beta & \gamma \end{bsmallmatrix} \in B(Np^2)\\p\mid 2\beta,\ p^2 \mid \gamma}}
p^k  a(\begin{bsmallmatrix}1&\\&p^{-1}\end{bsmallmatrix}S\begin{bsmallmatrix}1&\\&p^{-1}\end{bsmallmatrix})\E^{2\pi \I \mathrm{Tr}(SZ)}\\
&=\sum_{\substack{S=\begin{bsmallmatrix} \alpha & \beta \\ \beta & \gamma \end{bsmallmatrix} \in B(Np^2)\\ p\mid 2\beta,\ p^2 \mid \gamma}}
p^k a(S[\begin{bsmallmatrix}1&\\&p^{-1}\end{bsmallmatrix}])\E^{2\pi \I \mathrm{Tr}(SZ)}.
\end{align*}
This proves \eqref{etafouriereq}.

To obtain the Fourier-Jacobi expansion, we calculate as follows:
\begin{align*}
(\eta_p F)(Z)&=(F\big |_k\begin{bsmallmatrix}1&&&\\&p&&\\&&1&\\&&&p^{-1}\end{bsmallmatrix})(Z)\\
&=p^{k}F(\begin{bsmallmatrix}1&\\&p\end{bsmallmatrix}\begin{bsmallmatrix}\tau&z\vphantom{p}\\z&\tau'\vphantom{p}\end{bsmallmatrix}\begin{bsmallmatrix}1&\\&p\end{bsmallmatrix})\\
&=p^k\sum_{\substack{m=0\\N_s\mid m}}^\infty f_m(\tau,pz)\E^{2\pi \I mp^2\tau'}\\
&=p^k\sum_{\substack{m=0\\p^2\mid m,\ N_s \mid m}}^\infty f_{mp^{-2}}(\tau,pz)\E^{2\pi \I m\tau'}\\
&=p^k\sum_{\substack{m=0\\p^2 \mid m,\ N_s\mid m}}^\infty (f_{mp^{-2}}\big |_k U_p)(\tau,z)\E^{2\pi \I m\tau'}.
\end{align*}
This completes the proof.
\end{proof}

\section{A level lowering operator}
\index{level lowering operator}
In this section we calculate the formula for the level lowering operator $\sigma_p$ obtained 
from the corresponding local operator defined in Sect.~\ref{sigmasec}. The formula for the 
Fourier-Jacobi expansion of $\sigma_pF$ involves a certain index lowering operator $L_{c^2}$ on Jacobi forms. 
This operator is described in the next lemma.
Using the Petersson inner product as described in \cite{EZ} and \cite{KS}, one can define the adjoint of a
linear operator between spaces of Jacobi cusp forms.
One can show that, up to a constant, $L_{c^2}:J^{\mathrm{cusp}}_{k,mc^2}\to J_{k,m}^{\mathrm{cusp}}$ 
is the adjoint of the operator $U_c:J_{k,m}^{\mathrm{cusp}}\to J_{k,mc^2}^{\mathrm{cusp}}$ from Chap.~I, Sect.~4 on p.~41 of \cite{EZ}.

\begin{lemma}\label{Lpoperatorlemma}
Let $k$, $m$, and $c$ be integers such that $k>0$, $m\geq 0$ and $c >0$.  
Let $f$ be a Jacobi form of weight $k$ and index $mc^2$ on $\SL(2,\Z)$. Define 
$L_{c^2}f: \mathcal{H}_1 \times \C \to \C$
by 
\begin{equation}\label{Lcoperatoreq}
(L_{c^2}f)(\tau,z)=c^{-k}\sum_{\substack{a,b\in\Z/c\Z}}f(\tau,(b+a\tau+z)c^{-1}) \E^{2\pi \I m (2az+a^2\tau)}
\end{equation}
for $\tau \in \mathcal{H}_1$ and $z \in \C$. Then $L_{c^2} f$ is a well-defined
Jacobi form of weight $k$ and index $m$ on $\SL(2,\Z)$. If $f$ is a cusp form, then $L_{c^2}f$ is a cusp form.
\end{lemma}
\begin{proof}
If $g:\mathcal{H}_1 \times \C \to \C$ is a function, then we define $g\big |U:\mathcal{H}_1 \times \C \to \C$
by $g(\tau,z) = g(\tau,c^{-1} z)$ for $\tau \in \mathcal{H}_1$ and $z \in \C$. If $\gamma \in \SL(2,\R)$,
$[\lambda,\mu,\kappa] \in H(\R)$, 
and $g:\mathcal{H}_1 \times \C \to \C$ is a function, then 
\begin{align}
(g\big |U)\big |_{k,m} \gamma & = (g\big |_{k,mc^2} \gamma)|U,\label{Lpoperatorlemmaeq1}\\
(g\big |U)\big |_m [\lambda,\mu,\kappa] &=
(g\big |_{mc^2} [c^{-1} \lambda, c^{-1} \mu, c^{-2} \kappa] )\big | U.\label{Lpoperatorlemmaeq2}
\end{align}
Using \eqref{Lpoperatorlemmaeq2} we see that the sum 
$$
c^{-k}\sum_{a,b \in \Z/c \Z} (f\big |U)\big |_m [a,b,0]
=c^{-k}\sum_{\substack{a,b\in\Z/c\Z}}f(\tau,(b+a\tau+z)c^{-1}) \E^{2\pi \I m (2az+a^2\tau)}
$$
is well-defined. It follows that $L_{c^2}f$ is well-defined. 
It is straightforward  to verify that  $(L_{c^2}f)\big |_m[\lambda,\mu,\kappa]=L_{c^2}f$ 
for $[\lambda,\mu,\kappa] \in H(\Z)$. 
Let $\gamma=\begin{bsmallmatrix} \gamma_1&\gamma_2 \\ \gamma_3&\gamma_4 \end{bsmallmatrix} \in\SL(2,\Z)$. Then
\begin{align*}
(L_{c^2}f)\big |_{k,m}\gamma&=c^{-k}\sum_{a,b\in\Z/c\Z}((f\big |U)\big |_m[a,b,0])\big |_{k,m}\gamma\\
&=c^{-k}\sum_{a,b\in\Z/c\Z}((f\big |U)\big |_{k,m}\gamma)\big |_m [a\gamma_1+b\gamma_3,a\gamma_2+b\gamma_4,0]\\
&=c^{-k}\sum_{a,b\in\Z/c\Z}((f\big |_{k,mc^2}\gamma)|U)\big |_m[a,b]\\
&=c^{-k}\sum_{a,b\in\Z/c\Z}(f\big |U)\big |_m[a,b]\\
&=L_{c^2}f.
\end{align*}
Here, the second equality follows from the commutation rule \eqref{JFcommruleeq} and the third 
equality is verified by a calculation.
We next show that $L_{c^2}f$ has an absolutely convergent Fourier expansion of the form 
$$
L_{c^2}f(\tau,z)=\sum_{n'=0}^\infty \sum_{\substack{r \in \Z\\ r'^2 \leq 4n'm}} c'(n',r')\E^{2\pi \I (n'\tau+ r'z)}.
$$
Let 
$$
f(\tau,z)=\sum_{n=0}^\infty\sum_{\substack{r\in\Z \\ r^2 \leq 4nmc^2 }}c(n,r)\E^{2\pi \I (n\tau+ r z)}
$$
be the Fourier expansion of $f$. Then
\begin{align*}
&(L_{c^2}f)(\tau,z)=c^{-k}\sum_{a,b\in\Z/c\Z}f(\tau,(b+a\tau+z)c^{-1})\E^{2\pi \I m(2az+a^2\tau)}\\
&\qquad=c^{-k}\sum_{a=0}^{c-1}\sum_{n=0}^\infty \sum_{\substack{r\in\Z \\ r^2 \leq 4nmc^2 }} (\sum_{b\in \Z/c\Z}\E^{2\pi \I rbc^{-1}})
c(n,r)\\
&\qquad \quad \times \E^{2\pi \I (n+arc^{-1}+ma^2)\tau }\E^{2\pi \I (rc^{-1}+2ma)z}\\
&\qquad =c^{1-k}\sum_{a=0}^{c-1}\sum_{n=0}^\infty \sum_{\substack{r\in\Z \\ c\mid r,\ r^2 \leq 4nmc^2 }} 
c(n,r)\\
&\qquad\quad \times \E^{2\pi \I (n+arc^{-1}+ma^2)\tau }\E^{2\pi \I (rc^{-1}+2ma)z}\\
&\qquad=c^{1-k}\sum_{a=0}^{c-1}\sum_{n=0}^\infty \sum_{\substack{r\in\Z \\ r^2 \leq 4nm }} 
c(n,rc)\E^{2\pi \I (n+ar+ma^2)\tau }\E^{2\pi \I (r+2ma)z}.
\end{align*}
Let $X=\{(n,r) \in \Z \times \Z \mid n\geq 0,\ r^2 \leq 4nm\}$. 
Let $a \in \{0,\dots,c-1\}$. 
Define $i: \Z \times \Z \to \Z \times \Z$ by $i(n,r)=(n+ra+ma^2,r+2ma)$ for $(n,r) \in \Z \times \Z$.
Then $i$ is injective. 
Set $X_a = i(X)$.  Calculations show that if $(n',r') \in X_a$, then  $n' \geq 0$ and 
$r'^2 \leq 4 n' m$.
It follows that $X_a \subset X$. We now have
$$
(L_{c^2}f)(\tau,z)
=c^{1-k}\sum_{a=0}^{c-1}\sum_{(n',r') \in X_a} 
c(n,rc) \E^{2\pi \I (n' \tau + r' z)}
$$
where $(n,r)= i^{-1}(n',r')$. This proves that $L_{c^2}f$ has a Fourier expansion with the 
required properties.
It is easy to see from the Fourier expansion that if $f$ is a cusp form, then so is $L_{c^2}f$.
\end{proof}

\begin{lemma}\label{sigmaformulalemma}
Let $N$ and $k$ be integers such that $N>0$ and $k>0$, and let $p$ be a prime of $\Z$.
Assume that $v_p(N)\geq 2$. Let $\sigma: \mathcal{A}_k(\mathcal{K}_s(N)) \to
\mathcal{A}_k (\mathcal{K}_s(Np^{-1}))$ be the linear map obtained as in 
Sect.~\ref{overviewsec} from the level lowering operator \eqref{sigmaopseq},
and let $\sigma_p$ be as in \eqref{overviewTddefAeq}, so that $\sigma_p$ is a linear map 
$$
\sigma_p: M_k(\mathrm{K}_s(N)) \longrightarrow M_k(\mathrm{K}_s(Np^{-1})). 
$$
Let $F \in M_k(\mathrm{K}_s(N))$with Fourier expansion as in \eqref{NSFourierexpagaineq}. Then 
\begin{equation}
\label{sigmaslasheq}
\sigma_p F =p^{-3}\sum_{a,b,c\in\Z/p\Z}F\big |_k\begin{bsmallmatrix}p\\&1\\&&p\\&&&p^2\end{bsmallmatrix}
\begin{bsmallmatrix}1&&&b\\a&1&b&cp^{-n+2}\\&&1&-a\\&&&1\end{bsmallmatrix}.
\end{equation}
If $F$ is a cusp form, then $\sigma_p F$ is a cusp form. 
We have 
\begin{align}
(\sigma_{p}F)(Z)&= \sum_{\substack{S\in B(Np^{-1})}}\sum_{\substack{a\in\Z/p\Z}} p^{-k-1}
 a(S[\begin{bsmallmatrix}1&\\a&p\end{bsmallmatrix}]) \E^{2\pi \I \mathrm{Tr}{(SZ)}}.\label{sigmafouriereq}
\end{align}
The Fourier-Jacobi expansion of $\sigma_p F$ is given by 
\begin{equation}
(\sigma_{p}F)(Z)=p^{-2}\sum_{\substack{m=0\\p^{-1}N_s\mid m}}^{\infty} (L_{p^2}f_{mp^2})(\tau,z)\E^{2\pi \I m\tau'},
\qquad Z=\begin{bsmallmatrix}\tau&z\\z&\tau'\end{bsmallmatrix} \in \mathcal{H}_2.\label{sigmafourierjacobieq}
\end{equation}
Here $L_{p^2}$ is as in Lemma \ref{Lpoperatorlemma}. 
\end{lemma}	
\begin{proof}
Define $\Phi=\Phi_F$. Let $Z\in\mathcal{H}_2$ 
and let $h\in \GSp(4,\R)^+$ be such that $Z=h\langle I \rangle$. Let $n=v_p(N)$. Then
by \eqref{overviewTddefeq} and \eqref{sigmaalteq},
\begin{align*}
&(\sigma_p F)(Z)=\lambda(h)^{-k}j(h,I)^k\sigma\Phi(h)\\
&\qquad=\lambda(h)^{-k}j(h,I)^k p^{-3} \sum_{a,b,c \in \Z/ p\Z}
\Phi(h
\begin{bsmallmatrix}1&&&b\\a&1&b&cp^{-n+2}\\&&1&-a\\&&&1\end{bsmallmatrix}_{p}
\begin{bsmallmatrix}1\\&p\\&&1\\&&&p^{-1}\end{bsmallmatrix}_p)\\
&\qquad=\lambda(h)^{-k}j(h,I)^k p^{-3} \sum_{a,b,c\in\Z/p\Z}\Phi(
\underbrace{\begin{bsmallmatrix}1\\&p^{-1}\\&&1\\&&&p\end{bsmallmatrix}}_{\text{all places}}\\
&\qquad\quad\times\underbrace{\begin{bsmallmatrix}1&&&-b\\-a&1&-b&-cp^{-n+2}\\&&1&a\\&&&1\end{bsmallmatrix}}_{\text{all places}}h
\begin{bsmallmatrix}1&&&b\\a&1&b&cp^{-n+2}\\&&1&-a\\&&&1\end{bsmallmatrix}_{p}
\begin{bsmallmatrix}1\\&p\\&&1\\&&&p^{-1}\end{bsmallmatrix}_p)\\
&\qquad=\lambda(h)^{-k}j(h,I)^kp^{-3}\sum_{a,b,c\in\Z/p\Z}\Phi(
\begin{bsmallmatrix}p\\&1\\&&p\\&&&p^2\end{bsmallmatrix}_{\infty}
\begin{bsmallmatrix}1&&&b\\a&1&b&cp^{-n+2}\\&&1&-a\\&&&1\end{bsmallmatrix}_{\infty}h)\\
&\qquad=p^{-3}\sum_{a,b,c\in\Z/p\Z}(F\big |_k\begin{bsmallmatrix}p\\&1\\&&p\\&&&p^2\end{bsmallmatrix}
\begin{bsmallmatrix}1&&&b\\a&1&b&cp^{-n+2}\\&&1&-a\\&&&1\end{bsmallmatrix})(Z).
\end{align*}
To prove \eqref{sigmafouriereq} we use \eqref{sigmaslasheq}:
\begin{align*}
&(\sigma_pF)(Z)=p^{-3}\sum_{a,b,c\in\Z/p\Z}
(F\big |_k\begin{bsmallmatrix}p\\&1\\&&p\\&&&p^2\end{bsmallmatrix}
\begin{bsmallmatrix}1&&&b\\a&1&b&cp^{-n+2}\\&&1&-a\\&&&1\end{bsmallmatrix})(Z)\\
&\qquad=p^{-k-3}\sum_{a,b,c\in\Z/p\Z}
F(\begin{bsmallmatrix}p&\vphantom{p^{-1}}\\a&1\vphantom{p^{-1}}\end{bsmallmatrix}Z\begin{bsmallmatrix} p^{-1}&ap^{-2}\\&p^{-2}\end{bsmallmatrix}
+\begin{bsmallmatrix} & bp^{-1}\\ bp^{-1} & abp^{-2} + cp^{-n}\end{bsmallmatrix})\\
&\qquad=p^{-k-3}\sum_{S=\begin{bsmallmatrix} \alpha&\beta\\ \beta&\gamma \end{bsmallmatrix} \in B(N)}
(\sum_{a,b\in\Z/p\Z}\E^{2\pi \I (\gamma abp^{-2}+2\beta b p^{-1})})
(\sum_{c\in\Z/p\Z}\E^{2\pi \I \gamma c p^{-n}})\\
&\qquad \quad \times  a(S)\E^{2\pi \I \mathrm{Tr}(\begin{bsmallmatrix}p^{-1}&ap^{-2}\\&p^{-2}\end{bsmallmatrix}S\begin{bsmallmatrix}p&\vphantom{p^{-1}}\\a&1\vphantom{p^{-1}}\end{bsmallmatrix}Z)}\\
&\qquad=p^{-k-2}
\sum_{\substack{S=\begin{bsmallmatrix} \alpha&\beta\\ \beta&\gamma \end{bsmallmatrix} \in B(N)\\
p^n\mid\gamma}}\sum_{a\in\Z/p\Z}(\sum_{b\in\Z/p\Z}\E^{2\pi \I (2\beta b p^{-1})} )\\
&\qquad \quad \times  a(S)\E^{2\pi \I \mathrm{Tr}(\begin{bsmallmatrix}p^{-1}&ap^{-2}\\&p^{-2}\end{bsmallmatrix}S\begin{bsmallmatrix}p&\vphantom{p^{-1}}\\a&1\vphantom{p^{-1}}\end{bsmallmatrix}Z)}\\
&\qquad=p^{-k-1}\sum_{a\in\Z/p\Z}
\sum_{\substack{S=\begin{bsmallmatrix} \alpha&\beta\\ \beta&\gamma \end{bsmallmatrix} \in B(N)\\
p^n\mid\gamma,\ p\mid 2\beta}}
 a(S) \E^{2\pi \I \mathrm{Tr}(\begin{bsmallmatrix}p^{-1}&ap^{-2}\\&p^{-2}\end{bsmallmatrix}S\begin{bsmallmatrix}p&\vphantom{p^{-1}}\\a&1\vphantom{p^{-1}}\end{bsmallmatrix}Z)}\\
&\qquad=p^{-k-1}\sum_{a\in\Z/p\Z}\sum_{\substack{S\in B(Np^{-1})}}
 a(\begin{bsmallmatrix}1&a\\&p\end{bsmallmatrix}S\begin{bsmallmatrix}1&\\a&p\end{bsmallmatrix}) \E^{2\pi \I \mathrm{Tr}{(SZ)}}\\
&\qquad=\sum_{\substack{S\in B(Np^{-1})}}\sum_{\substack{a\in\Z/p\Z}}p^{-k-1}
 a(S[\begin{bsmallmatrix}1&\\a&p\end{bsmallmatrix}]) \E^{2\pi \I \mathrm{Tr}{(SZ)}}.
\end{align*}
This proves \eqref{sigmafouriereq}.
 
Finally, beginning from \eqref{sigmaslasheq}, we have:
\begin{align*}
&(\sigma_pF)(Z)=p^{-k-3}\sum_{a,b,c\in\Z/p\Z}
F(\begin{bsmallmatrix}p&\vphantom{p^{-1}}\\a&1\vphantom{p^{-1}}\end{bsmallmatrix}Z\begin{bsmallmatrix} p^{-1}&ap^{-2}\\&p^{-2}\end{bsmallmatrix}
+\begin{bsmallmatrix} & bp^{-1}\\ bp^{-1} & abp^{-2} + cp^{-n}\end{bsmallmatrix})\\
&\qquad=p^{-k-3}\sum_{a,b,c\in\Z/p\Z}
\sum_{\substack{m=0\\N_s\mid m}}^{\infty}f_m(\tau,(b+a\tau+z)p^{-1})\\
&\qquad\quad\times \E^{2\pi \I m(2ap^{-2}z+a^2p^{-2}\tau+cp^{-n}+abp^{-2}+p^{-2}\tau')}\\
&\qquad=p^{-k-3}\sum_{\substack{a,b\in\Z/p\Z}}
\sum_{\substack{m=0\\N_s\mid m}}^{\infty}f_m(\tau,(b+a\tau+z)p^{-1})\\
&\quad\qquad\times (\sum_{c\in\Z/p\Z}e^{2\pi i mcp^{-n}}) \E^{2\pi \I m(2ap^{-2}z+a^2p^{-2}\tau+abp^{-2}+p^{-2}\tau')}\\
&\qquad=p^{-k-2}\sum_{\substack{a,b\in\Z/p\Z}}
\sum_{\substack{m=0\\pN_s\mid m}}^{\infty}f_m(\tau,(b+a\tau+z)p^{-1}) \E^{2\pi \I mp^{-2}(2az+a^2\tau+\tau')}\\
&\qquad=p^{-k-2}\sum_{\substack{a,b\in\Z/p\Z}}
\sum_{\substack{m=0\\pN_s\mid m}}^{\infty}f_m(\tau,(b+a\tau+z)p^{-1}) \E^{2\pi \I mp^{-2}(2az+a^2\tau)}\E^{2\pi \I \tau'mp^{-2}}\\
&\qquad=p^{-k-2}\sum_{\substack{m=0\\p^{-1}N_s\mid m}}^{\infty}
\sum_{\substack{a,b\in\Z/p\Z}}f_{mp^2}(\tau,(b+a\tau+z)p^{-1}) \E^{2\pi \I m (2az+a^2\tau)}\E^{2\pi \I m\tau'}\\
&\qquad=p^{-2}\sum_{\substack{m=0\\p^{-1}N_s\mid m}}^{\infty} (L_{p^2}f_{mp^2})(\tau,z)\E^{2\pi \I m\tau'}.
\end{align*}
This completes the proof.
\end{proof}

\section{Hecke operators}
In this section we calculate the formulas 
for the Hecke operators $T_{0,1}^s(p)$ and $T_{1,0}^s(p)$ obtained 
from the corresponding local operators defined in Sect.~\ref{stableheckesec}.
The formula for the Fourier-Jacobi expansion of $T_{1,0}^s(p)F$ involves the 
index lowering operator $L_{c^2}$ on Jacobi forms discussed in the previous section.
For the formula for the Fourier-Jacobi expansion of $T_{0,1}^s(p)F$ we introduce in 
Lemma~\ref{Lppoperlemma} below another index lowering operator $L'_p$ on Jacobi forms. 
One can show that, up to a constant, $L'_p:J^{\mathrm{cusp}}_{k,mp}\to J_{k,m}^{\mathrm{cusp}}$ is 
the adjoint of the operator $V_p:J_{k,m}^{\mathrm{cusp}}\to J_{k,mp}^{\mathrm{cusp}}$ from Chap.~I, Sect.~4 on p.~41 of \cite{EZ}. 
The adjoint of $V_N:J_{k,1}^{\mathrm{cusp}}\to J_{k,N}^{\mathrm{cusp}}$ was also considered in~\cite{KS}.

\begin{lemma}
\label{Lppoperlemma}
Let $k$ and $m$ be integers such that $k >0$ and $m \geq 0$.
Let $p$ be a prime of $\Z$.
Let $f$ be a Jacobi form of weight $k$ and index $pm$ on $\SL(2,\Z)$. \label{Lpprimedef}
Define $L_p' f:\mathcal{H}_1 \times \C \to \C$
by 
\begin{align*}
L'_pf(\tau,z)&=p\sum_{a\in\Z/p\Z}f(p\tau, a\tau+z)\E^{2\pi \I m(a^2\tau+2az)}\\
&\quad+p^{1-k}\sum_{a,b\in\Z/p\Z} f(p^{-1}(\tau+b),p^{-1}(z+a)) 
\end{align*}
for $\tau \in \mathcal{H}_1$ and $z \in \C$. 
Then $L'_pf$ is a well-defined Jacobi form of weight $k$ and index $m$ on $\SL(2,\Z)$. 
If $f$ is a cusp form, then $L'_pf$ is a cusp form.
\end{lemma}
\begin{proof}
To prove that $L_p'f$ is well-defined it will suffice to prove that 
\begin{equation}
\label{Lppoperlemmaeq1}
f(p\tau,(a+\ell p)+z)\E^{2\pi \I m((a+\ell p)^2 \tau + 2(a+\ell p) z)}
= f(p\tau, a\tau+z)\E^{2\pi \I m (a^2 \tau + 2 a z)}
\end{equation}
and
\begin{equation}
\label{Lppoperlemmaeq2}
f(p^{-1}(\tau+b+jp),p^{-1}(z+a+\ell p))
=f(p^{-1}(\tau+b),p^{-1}(z+a))
\end{equation}
for $a,b,\ell,j \in \Z$, $\tau \in \mathcal{H}_1$, and $z \in \C$. 
The equation \eqref{Lppoperlemmaeq2} follows from the fact that $f(\tau+n_1,z+n_2)
=f(\tau,z)$ for $\tau\in \mathcal{H}_1$, $z \in \C$, and $n_1,n_2 \in \Z$. 
For \eqref{Lppoperlemmaeq1}, we note first that 
$$
f(p\tau, a\tau+z)\E^{2\pi \I m (a^2 \tau + 2 a z)}
= (f\big |_{pm} [ap^{-1},0,0])(p\tau,z)
$$
for $a \in \Z$, $\tau \in \mathcal{H}_1$, and $z \in \C$; \eqref{Lppoperlemmaeq1}
now follows from $f\big |_{pm}[n,0,0] = f$ for $n \in \Z$. 

To prove that $L_p'f\big |_{k,m} g = L_p' f$ for $g \in \SL(2,\Z)$
it will suffice to prove that $(L'_pf)\big |_{k,m}\begin{bsmallmatrix}&1\\-1&\end{bsmallmatrix} = L_p' f$
and $(L'_pf)\big |_{k,m}\begin{bsmallmatrix}1&1\\&1\end{bsmallmatrix} = L_p'f$. 
For each $b \in \Z$ relatively prime to $p$ fix $e_b,f_b \in \Z$ such that 
$bf_b-e_b p=1$. In the following calculation we will use that the map $(\Z/p\Z)^\times \to (\Z/p\Z)^\times$
defined by $b \mapsto -f_b$ is a bijection; we will also use \eqref{FJbasiceq}.
Let $\tau \in \mathcal{H}_1$ and $z \in \C$. Then 
\begin{align*}
&(L'_pf)\big |_{k,m}\begin{bsmallmatrix}&1\\-1&\end{bsmallmatrix}(\tau,z)\\
&\qquad=(-\tau)^{-k}\E^{2\pi \I m(-z^2\tau^{-1})}L'_pf(-\tau^{-1},-z\tau^{-1})\\
&\qquad=(-\tau)^{-k}e^{2\pi i m(-z^2\tau^{-1})}\\
&\qquad\quad\times\Big(p\sum_{a\in\Z/p\Z} f(-p\tau^{-1}, -a\tau^{-1}-z\tau^{-1})\E^{2\pi \I m(-a^2\tau^{-1}-2az\tau^{-1})}\\
&\qquad\quad+p^{1-k}\sum_{a,b\in\Z/p\Z} f(p^{-1}(-\tau^{-1}+b),p^{-1}(-z\tau^{-1}+a))\Big)\\
&\qquad=(-\tau)^{-k}e^{2\pi i m(-z^2\tau^{-1})}\\
&\qquad\quad\times\Big(p\sum_{a\in\Z/p\Z} (f\big |_{k,pm}\begin{bsmallmatrix}&-1\\1&\end{bsmallmatrix})(-p\tau^{-1}, -a\tau^{-1}-z\tau^{-1})\E^{2\pi \I m(-a^2\tau^{-1}-2az\tau^{-1})}\\
&\qquad\quad+p^{1-k}\sum_{a,b\in\Z/p\Z} (f\big |_{k,pm}\begin{bsmallmatrix}&-1\\1&\end{bsmallmatrix})(p^{-1}(-\tau^{-1}+b),p^{-1}(-z\tau^{-1}+a))\Big)\\
&\qquad=p^{1-k}(-\tau)^{-k} \E^{2\pi \I m(-z^2\tau^{-1})}\\
&\qquad\quad\times\Big(\sum_{a\in\Z/p\Z}(-\tau)^k f(p^{-1}\tau, p^{-1}(z+a)) \E^{2\pi \I m(a+z)^2\tau^{-1}}\E^{2\pi \I m(-a^2\tau^{-1}-2az\tau^{-1})}\\
&\qquad\quad+\sum_{a,b\in\Z/p\Z}(\frac{p\tau}{b\tau-1})^k \E^{2\pi \I m\frac{-(a\tau-z)^2}{\tau(b\tau-1)}}f(\frac{-p\tau}{b\tau-1},\frac{a\tau-z}{b\tau-1})\Big)\\
&\qquad=p^{1-k}\sum_{a\in\Z/p\Z} f(p^{-1}\tau, p^{-1}(z+a))\\
&\qquad\quad+p\sum_{a,b\in\Z/p\Z}(\frac{-1}{b\tau-1})^k \E^{2\pi \I m(\frac{-a^2\tau+2az-z^2b}{b\tau-1})}f(\frac{-p\tau}{b\tau-1},\frac{a\tau-z}{b\tau-1})\\
&\qquad=p\sum_{a\in\Z/p\Z} \E^{2\pi \I m(a^2\tau-2az)}f(p\tau,-a\tau+z)+p^{1-k}\sum_{a\in\Z/p\Z} f(p^{-1}\tau, p^{-1}(z+a))\\
&\qquad\quad+p\sum_{\substack{a,b\in\Z/p\Z\\b\neq 0}}(\frac{-1}{b\tau-1})^k  \E^{2\pi \I m(\frac{-a^2\tau+2az-z^2b}{b\tau-1})}\\
&\qquad\quad\times(f\big |_{k,pm}\begin{bsmallmatrix}e_b&f_b\\-b&-p\end{bsmallmatrix})(\frac{-p\tau}{b\tau-1},\frac{a\tau-z}{b\tau-1})\\
&\qquad=p\sum_{a\in\Z/p\Z} \E^{2\pi \I m(a^2\tau+2az)}f(p\tau,a\tau+z)+p^{1-k}\sum_{a\in\Z/p\Z} f(p^{-1}\tau, p^{-1}(z+a))\\
&\qquad\quad+p^{1-k}(-1)^k\sum_{\substack{a,b\in\Z/p\Z\\b\neq 0}} \E^{2\pi \I m(a^2\tau-2az)}f(p^{-1}(\tau-f_b),p^{-1}(a\tau-z))\\
&\qquad=p\sum_{a\in\Z/p\Z} \E^{2\pi \I m(a^2\tau+2az)}f(p\tau,a\tau+z)+p^{1-k}\sum_{a\in\Z/p\Z} f(p^{-1}\tau, p^{-1}(z+a))\\
&\qquad\quad+p^{1-k}(-1)^k\sum_{\substack{a,b\in\Z/p\Z\\b\neq 0}} \E^{2\pi \I m(a^2\tau-2az)}f(p^{-1}(\tau+b),p^{-1}(a\tau-z))\\
&\qquad=p\sum_{a\in\Z/p\Z} \E^{2\pi \I m(a^2\tau+2az)}f(p\tau,a\tau+z)+p^{1-k}\sum_{a\in\Z/p\Z} f(p^{-1}\tau, p^{-1}(z+a))\\
&\qquad\quad+p^{1-k}(-1)^k\sum_{\substack{a,b\in\Z/p\Z\\b\neq 0}} \E^{2\pi \I m(a^2\tau-2az)}\\
&\qquad\quad\times(f\big |_{pm}[-a,0,0])(p^{-1}(\tau+b),p^{-1}(a\tau-z))\\
&\qquad=p\sum_{a\in\Z/p\Z} \E^{2\pi \I m(a^2\tau+2az)}f(p\tau,a\tau+z)+p^{1-k}\sum_{a\in\Z/p\Z} f(p^{-1}\tau, p^{-1}(z+a))\\
&\qquad\quad+p^{1-k}(-1)^k\sum_{\substack{a,b\in\Z/p\Z\\b\neq 0}}f(p^{-1}(\tau+b),-p^{-1}(z+ab))\\
&\qquad=p\sum_{a\in\Z/p\Z} \E^{2\pi \I m(a^2\tau+2az)}f(p\tau,a\tau+z)+p^{1-k}\sum_{a\in\Z/p\Z} f(p^{-1}\tau, p^{-1}(z+a))\\
&\qquad\quad+p^{1-k}\sum_{\substack{a,b\in\Z/p\Z\\b\neq 0}}f(p^{-1}(\tau+b),p^{-1}(z+a))\\
&\qquad=p\sum_{a\in\Z/p\Z} \E^{2\pi \I m(a^2\tau+2az)}f(p\tau,a\tau+z)\\
&\qquad\quad+p^{1-k}\sum_{\substack{a,b\in\Z/p\Z}}f(p^{-1}(\tau+b),p^{-1}(z+a))\\
&\qquad=L'_p(\tau,z).	
\end{align*}
And:
\begin{align*}
(L'_pf)\big |_{k,m}\begin{bsmallmatrix}1&1\\&1\end{bsmallmatrix}(\tau,z)&=(L'_pf)(\tau+1,z)\\
&=p\sum_{a\in\Z/p\Z}f(p\tau+p, a\tau+z+a) \E^{2\pi \I m(a^2(\tau+1)+2az)}\\
&\quad+p^{1-k}\sum_{a,b\in\Z/p\Z} f(p^{-1}(\tau+b+1),p^{-1}(z+a))\\
&=p\sum_{a\in\Z/p\Z}f(p\tau, a\tau+z) \E^{2\pi \I m(a^2\tau+2az)}\\
&\quad+p^{1-k}\sum_{a,b\in\Z/p\Z} f(p^{-1}(\tau+b),p^{-1}(z+a))\\
&=L'_p(\tau,z),
\end{align*}
where the penultimate formula follows from \eqref{FJbasiceq}. 
Next, to prove that $L_p'f\big |_m h = L_p'f$ for $h \in H(\Z)$, it suffices to prove that
$L_p'f\big |_m h = L_p'f$ for $h=[1,0,0]$ and $h=[0,1,0]$. We have
\begin{align*}
&(L'_pf\big |_{m}[1,0,0])(\tau,z)=\E^{2\pi \I m(\tau+2z)}L'_p(\tau,z+\tau)\\
&\qquad=\E^{2\pi \I m(\tau+2z)}\Big(p\sum_{a\in\Z/p\Z}f(p\tau, a\tau+(z+\tau))\E^{2\pi \I m(a^2\tau+2a(z+\tau))}\\
&\qquad\quad+p^{1-k}\sum_{a,b\in\Z/p\Z} f(p^{-1}(\tau+b),p^{-1}(z+\tau+a))\Big)\\
&\qquad=\E^{2\pi \I m(\tau+2z)}\Big(p\sum_{a\in\Z/p\Z}f(p\tau, a\tau+z)\E^{2\pi \I m((a-1)^2\tau+2(a-1)(z+\tau))}\\
&\qquad\quad+p^{1-k}\sum_{a,b\in\Z/p\Z} (f\big |_{pm}[-1,0,0])(p^{-1}(\tau+b),p^{-1}(z+\tau+a))\Big)\\
&\qquad=p\sum_{a\in\Z/p\Z}f(p\tau, a\tau+z) \E^{2\pi \I m(a^2\tau+2az)}\\
&\qquad\quad+p^{1-k}\sum_{a,b\in\Z/p\Z} \E^{2\pi \I m(\tau+2z)} \E^{2\pi \I m(-\tau-2z+b-2a)}\\
&\qquad\quad\times f(p^{-1}(\tau+b),p^{-1}(z+a-b))\\
&\qquad=p\sum_{a\in\Z/p\Z}f(p\tau, a\tau+z) \E^{2\pi \I m(a^2\tau+2az)}\\
&\qquad\quad+p^{1-k}\sum_{a,b\in\Z/p\Z} f(p^{-1}(\tau+b),p^{-1}(z+a))\\
&\qquad=L'_p(\tau,z).
\end{align*}
And
\begin{align*}
(L'_pf\big |_{m}[0,1,0])(\tau,z)&=L'_p(\tau,z+1)\\
&=p\sum_{a\in\Z/p\Z}f(p\tau, a\tau+z+1) \E^{2\pi \I m(a^2\tau+2a(z+1))}\\
&\quad+p^{1-k}\sum_{a,b\in\Z/p\Z} f(p^{-1}(\tau+b),p^{-1}(z+1+a))\\
&=p\sum_{a\in\Z/p\Z}f(p\tau, a\tau+z) \E^{2\pi \I m(a^2\tau+2az)}\\
&\quad+p^{1-k}\sum_{a,b\in\Z/p\Z} f(p^{-1}(\tau+b),p^{-1}(z+a))\\
&=L'_p(\tau,z).
\end{align*}

We next show that $L'_pf$ has an absolutely convergent Fourier expansion of the form 
$$
L'_pf(\tau,z)=\sum_{n'=0}^\infty \sum_{\substack{r'\in\Z\\r'^2 \leq 4n'm}}c'(n',r') \E^{2\pi \I (n'\tau+r'z)}.
$$
We will show that each of the summands of $L_p'f$ has this property.  Let 
$$
f(\tau,z) = \sum_{n=0}^\infty \sum_{\substack{r\in\Z\\r^2 \leq 4npm}}c(n,r) \E^{2\pi \I (n\tau+rz)}
$$
be the Fourier expansion of $f$. 
Let $a$ be an integer such that $a \geq 0$.
We have
$$
\E^{2\pi \I m(a^2\tau+2az)}f(p\tau,a\tau+z)
=\sum_{n=0}^\infty \sum_{\substack{r\in\Z\\ r^2 \leq 4npm}}c(n,r) \E^{2\pi \I ((np+ra+ma^2)\tau+(r+2am)z)}.
$$
Let $X=\{(n,r) \in \Z \times \Z: n\geq 0, r^2 \leq 4npm\}$. Define $i:\Z \times \Z \to \Z \times \Z$
by $i(n,r)=(np+ra+ma^2,r+2am)$, and let $X'=i(X)$. The map $i$ is an injection. Calculations show
that if $(n',r') \in X'$, then $n' \geq 0$ and $r'^2 \leq 4n'm$. We now have 
$$
\E^{2\pi \I m(a^2 \tau +2 a z)} f(p\tau,a\tau +z)
=
\sum_{(n',r') \in X'} c(i^{-1}(n',r')) \E^{2 \pi \I (n' \tau+ r' z)}.
$$
It follows that the first summand 
$$
p\sum_{a\in\Z/p\Z}f(p\tau, a\tau+z)\E^{2\pi \I m(a^2\tau+2az)}
$$
of $L_p'f$ has a Fourier expansion with the required properties.
We may calculate the second summand of $L_p'f$ as follows:
\begin{align*}
&p^{1-k}\sum_{a,b\in\Z/p\Z} f(p^{-1}(\tau+b),p^{-1}(z+a))\\
&\qquad=p^{1-k}\sum_{a,b\in\Z/p\Z}\sum_{n=0}^\infty\sum_{\substack{r\in \Z\\ r^2\leq 4npm}} c(n,r) \E^{2\pi \I ( np^{-1}(\tau+b)+ rp^{-1}(z+a))}\\
&\qquad=p^{1-k}\sum_{n=0}^\infty\sum_{\substack{r\in \Z\\ r^2\leq 4npm}} 
(\sum_{a,b\in\Z/p\Z} \E^{2\pi \I (np^{-1}b+rp^{-1}a)})c(n,r) \E^{2\pi \I (np^{-1}\tau+ rp^{-1}z)}\\
&\qquad=p^{3-k}\sum_{\substack{n=0\\ p\mid n}}^\infty\sum_{\substack{r\in \Z\\ r^2\leq 4npm,\ p \mid r}} c(n,r) \E^{2\pi \I (np^{-1}\tau+ rp^{-1}z)}\\
&\qquad=p^{3-k}\sum_{n=0}^\infty\sum_{\substack{r\in \Z\\ r^2\leq 4nm}} c(np,rp) \E^{2\pi \I (n\tau+ rz)}.
\end{align*}
This shows that the second summand also has a Fourier expansion with the required properties.
It is easy to see from the Fourier expansion that if $f$ is a cusp form, then so is $L'_pf$.
\end{proof}
\begin{lemma}\label{T01formulalemma}
Let $N$ and $k$ be integers such that $N>0$ and $k>0$, and let $p$ be a prime of $\Z$.
Assume that $v_p(N) \geq 1$. 
Let $T^s_{0,1}:\mathcal{A}_k(\mathcal{K}_s(N)) \to \mathcal{A}_k(\mathcal{K}_s(N))$ be the linear
map obtained as in 
Sect.~\ref{overviewsec} from the Hecke operator \eqref{Ts01eq},and let $T_{0,1}^{s}(p)=(T_{0,1}^s)_p$ 
be as in \eqref{overviewTddefAeq}, so that $T_{0,1}^{s}(p)$ is a linear map
$$
T_{0,1}^{s}(p): M_k(\mathrm{K}_s(N)) \longrightarrow M_k(\mathrm{K}_s(N)). 
$$
Let $F\in M_k(\mathcal{K}_s(N))$ with Fourier expansion as in \eqref{NSFourierexpagaineq}. Then 
\begin{align}
T_{0,1}^{s}(p)F&=\sum_{a,b\in\Z/p\Z}F\big |_k\begin{bsmallmatrix}p&&&\vphantom{bp^{-n+1}}\\&1&&\\&&1&\\&&&p\end{bsmallmatrix}
\begin{bsmallmatrix}1&&&\\a&1&&bp^{-n+1}\\&&1&-a\\&&&1\end{bsmallmatrix}\nonumber\\
&\quad+\sum_{a,b,c\in\Z/p\Z}F\big |_k \begin{bsmallmatrix}1&&&\vphantom{bp^{-n+1}}\\&1&&\\&&p&\\&&&p\end{bsmallmatrix}
\begin{bsmallmatrix}1&&c&a\\&1&a&bp^{-n+1}\\&&1&\\&&&1\end{bsmallmatrix}.\label{T01slasheq}
\end{align}
If $F$ is a cusp form, then $T_{0,1}^{s}(p)F$ is a cusp form. 
We have 
\begin{align}
(T_{0,1}^{s}(p)F)(Z)
&=\sum_{S=\begin{bsmallmatrix}\alpha&\beta\\ \beta&\gamma \end{bsmallmatrix}\in B(N)} 
\Big( p^{3-k}  a(pS)\nonumber\\
&\quad+  \sum_{\substack{y \in \Z/p\Z \\p\mid(\alpha+2\beta y+\gamma y^2)}} p  a(p^{-1}S[\begin{bsmallmatrix}1&\\y&p\end{bsmallmatrix}]) \Big) \E^{2\pi \I \mathrm{Tr}(SZ)}.\label{T01fouriereq}
\end{align}
The Fourier-Jacobi expansion of $T_{0,1}^{s}(p) F$ is given by 
\begin{equation}
(T_{0,1}^{s}(p)F)(Z)=\sum_{\substack{m=0\\N_s\mid m}}^\infty(L'_pf_{mp})(\tau,z)\E^{2\pi \I m\tau'},
\qquad Z=\begin{bsmallmatrix}\tau&z\\z&\tau'\end{bsmallmatrix} \in \mathcal{H}_2.\label{T01fourierjacobieq}
\end{equation}
Here $L_p'$ is the operator on Jacobi forms from Lemma \ref{Lppoperlemma}.
\end{lemma}	
\begin{proof}
Define $\Phi=\Phi_F$. Let $Z\in\mathcal{H}_2$ 
and let $h\in \GSp(4,\R)^+$ be such that $Z=h\langle I \rangle$. Let $n=v_p(N)$. Then
by \eqref{overviewTddefeq} and \eqref{Ts01eq} we have
\begin{align*}
&(T_{0,1}^{s}(p)F)(Z)=\lambda(h)^{-k}j(h,I)^kT^{s}_{0,1}\Phi(h)\\
&\qquad=\lambda(h)^{-k}j(h,I)^k\big(\sum_{y,z\in\Z/p\Z}
\Phi(h
\begin{bsmallmatrix}1&&&\\y&1&&zp^{-n+1}\\&&1&-y\\&&&1\end{bsmallmatrix}_p
\begin{bsmallmatrix}1\vphantom{p^{-n+1}}&&&\\ &p&&\vphantom{-y}\\&&p&\\&&&1\end{bsmallmatrix}_p)\\
&\qquad\quad+\sum_{c,y,z\in\Z/p\Z}\Phi(h 
\begin{bsmallmatrix}1&&c&y\\&1&y&zp^{-n+1}\\&&1&\\&&&1\end{bsmallmatrix}_p
\begin{bsmallmatrix}p\vphantom{p^{-n+1}}&&&\\&p&&\vphantom{-y}\\&&1&\\&&&1\end{bsmallmatrix}_p)\big)\\	
&\qquad=\lambda(h)^{-k}j(h,I)^k\big(\sum_{y,z\in\Z/p\Z}
\Phi(\underbrace{\begin{bsmallmatrix}1&&&\\&p^{-1}&&\\&&p^{-1}&\\&&&1\end{bsmallmatrix}
\begin{bsmallmatrix}1&&&\\-y&1&&-zp^{-n+1}\\&&1&y\\&&&1\end{bsmallmatrix}}_{\text{all places}}h\\
&\qquad\quad\times\begin{bsmallmatrix}1&&&\\y&1&&zp^{-n+1}\\&&1&-y\\&&&1\end{bsmallmatrix}_p
\begin{bsmallmatrix}1\vphantom{p^{-n+1}}&&&\\&p&&\vphantom{-y}\\&&p&\\&&&1\end{bsmallmatrix}_p)\\
&\qquad\quad+\sum_{c,y,z\in\Z/p\Z}\Phi(\underbrace{\begin{bsmallmatrix}p^{-1}&&&\\&p^{-1}&&\\&&1&\\&&&1\end{bsmallmatrix}
\begin{bsmallmatrix}1&&-c&-y\\&1&-y&-zp^{-n+1}\\&&1&\\&&&1\end{bsmallmatrix}}_{\text{all places}}h\\
&\qquad\quad\times\begin{bsmallmatrix}1&&c&y\\&1&y&zp^{-n+1}\\&&1&\\&&&1\end{bsmallmatrix}_p
\begin{bsmallmatrix}p\vphantom{p^{-n+1}}&&&\\&p&&\vphantom{-y}\\&&1&\\&&&1\end{bsmallmatrix}_p) \big)\\	
&\qquad=\lambda(h)^{-k}j(h,I)^k\big(\sum_{y,z\in\Z/p\Z}\Phi(\begin{bsmallmatrix}1&&&\\&p^{-1}&&\\&&p^{-1}&\\&&&1\end{bsmallmatrix}_\infty
\begin{bsmallmatrix}1&&&\\-y&1&&-zp^{-n+1}\\&&1&y\\&&&1\end{bsmallmatrix}_{\infty}h)\\
&\qquad\quad+\sum_{c,y,z\in\Z/p\Z}\Phi(\begin{bsmallmatrix}p^{-1}&&&\\&p^{-1}&&\\&&1&\\&&&1\end{bsmallmatrix}_\infty
\begin{bsmallmatrix}1&&-c&-y\\&1&-y&-zp^{-n+1}\\&&1&\\&&&1\end{bsmallmatrix}_{\infty}h)\big)\\
&\qquad=\lambda(h)^{-k}j(h,I)^k\big(\sum_{y,z\in\Z/p\Z}
\Phi(\begin{bsmallmatrix}p\vphantom{p^{-1}}&&&\\&1&&\\&&1&\\&&&p\vphantom{p^{-1}}\end{bsmallmatrix}_\infty
\begin{bsmallmatrix}1&&&\\y&1&&zp^{-n+1}\\&&1&-y\\&&&1\end{bsmallmatrix}_{\infty}h)\\
&\qquad\quad+\sum_{c,y,z\in\Z/p\Z}\Phi(\begin{bsmallmatrix}1&&&\\&1&&\\&&p\vphantom{p^{-1}}&\\&&&p\vphantom{p^{-1}}\end{bsmallmatrix}_\infty
\begin{bsmallmatrix}1&&c&y\\&1&y&zp^{-n+1}\\&&1&\\&&&1\end{bsmallmatrix}_{\infty}h)\big)\\
&\qquad=\sum_{y,z\in\Z/p\Z}F\big |_k(\begin{bsmallmatrix}p\vphantom{p^{-1}}&&&\\&1&&\\&&1&\\&&&p\vphantom{p^{-1}}\end{bsmallmatrix}
\begin{bsmallmatrix}1&&&\\y&1&&zp^{-n+1}\\&&1&-y\\&&&1\end{bsmallmatrix})(Z)\\
&\qquad\quad+\sum_{c,y,z\in\Z/p\Z}F\big |_k(\begin{bsmallmatrix}1&&&\\&1&&\\&&p\vphantom{p^{-1}}&\\&&&p\vphantom{p^{-1}}\end{bsmallmatrix}
\begin{bsmallmatrix}1&&c&y\\&1&y&zp^{-n+1}\\&&1&\\&&&1\end{bsmallmatrix})(Z).
\end{align*}
This proves \eqref{T01slasheq}. 
To prove \eqref{T01fouriereq} we  use \eqref{T01slasheq}:
\begin{align*}
&(T_{0,1}^{s}(p)F)(Z)\\
&\qquad=\sum_{y,z\in\Z/p\Z}F\big |_k(\begin{bsmallmatrix}p&&&\vphantom{bp^{-n+1}}\\&1&&\\&&1&\\&&&p\end{bsmallmatrix}
\begin{bsmallmatrix}1&&&\\y&1&&zp^{-n+1}\\&&1&-y\\&&&1\end{bsmallmatrix})(Z)\\
&\qquad\quad+\sum_{c,y,z\in\Z/p\Z}F\big |_k(\begin{bsmallmatrix}1&&&\vphantom{bp^{-n+1}}\\&1&&\\&&p&\\&&&p\end{bsmallmatrix}
\begin{bsmallmatrix}1&&c&y\\&1&y&zp^{-n+1}\\&&1&\\&&&1\end{bsmallmatrix})(Z)\\
&\qquad=\sum_{y,z\in\Z/p\Z}F\big |_k(\begin{bsmallmatrix}p&&&\\y&1&&zp^{-n+1}\\&&1&-y\\&&&p\end{bsmallmatrix})(Z)
+\sum_{c,y,z\in\Z/p\Z}F\big |_k(\begin{bsmallmatrix}1&&c&y\\&1&y&zp^{-n+1}\\&&p&\\&&&p\end{bsmallmatrix})(Z)\\
&\qquad=\sum_{y,z\in\Z/p\Z}F(\begin{bsmallmatrix}p&&&\\y&1&&zp^{-n+1}\\&&1&-y\\&&&p\end{bsmallmatrix}\langle Z\rangle)
+p^{-k}\sum_{c,y,z\in\Z/p\Z}F(\begin{bsmallmatrix}1&&c&y\\&1&y&zp^{-n+1}\\&&p&\\&&&p\end{bsmallmatrix}\langle Z\rangle )\\
&\qquad=\sum_{y,z\in\Z/p\Z}F(\begin{bsmallmatrix}p\vphantom{p^{-1}}&\\y&1\vphantom{p^{-1}}\end{bsmallmatrix}Z
\begin{bsmallmatrix}1&yp^{-1}\\&p^{-1}\end{bsmallmatrix}
+\begin{bsmallmatrix}\vphantom{p^{-1}}&\\ \vphantom{p^{-1}}&zp^{-n}\end{bsmallmatrix})\\
&\qquad\quad+p^{-k}\sum_{y,z,c\in\Z/p\Z}F(p^{-1}Z+\begin{bsmallmatrix}cp^{-1}&yp^{-1}\\yp^{-1}&zp^{-n}\end{bsmallmatrix})\\
&\qquad=\sum_{y,z\in\Z/p\Z}\sum_{S\in B(N)} a(S) \E^{2\pi \I \mathrm{Tr}(S
\begin{bsmallmatrix}p\vphantom{p^{-1}}&\\y&1\vphantom{p^{-1}}\end{bsmallmatrix}Z
\begin{bsmallmatrix}1&yp^{-1}\\&p^{-1}\end{bsmallmatrix}
+S\begin{bsmallmatrix}\vphantom{p^{-1}}&\\ \vphantom{p^{-1}} &zp^{-n}\end{bsmallmatrix})}\\
&\qquad\quad+p^{-k}\sum_{y,z,c\in\Z/p\Z}\sum_{S\in B(N)} a(S)
\E^{2\pi \I \mathrm{Tr} (p^{-1}SZ+S\begin{bsmallmatrix}cp^{-1}&yp^{-1}\\ yp^{-1}&zp^{-n}\end{bsmallmatrix})}\\
&\qquad=\sum_{y\in\Z/p\Z}\sum_{S=\begin{bsmallmatrix} \alpha&\beta\\ \beta&\gamma \end{bsmallmatrix}\in B(N)}
\sum_{z\in\Z/p\Z} a(S)\E^{2\pi \I \gamma z p^{-n}} \E^{2\pi \I \mathrm{Tr}(\begin{bsmallmatrix}1&yp^{-1}\\&p^{-1}\end{bsmallmatrix}S
\begin{bsmallmatrix}p \vphantom{p^{-1}}&\\y&1 \vphantom{p^{-1}}\end{bsmallmatrix}Z)}\\
&\qquad\quad+p^{-k}\sum_{S=\begin{bsmallmatrix} \alpha&\beta\\ \beta&\gamma \end{bsmallmatrix}\in B(N)}
\sum_{y,z,c\in\Z/p\Z} a(S)
\E^{2\pi \I (\alpha c p^{-1}+2\beta yp^{-1}+\gamma z p^{-n})} \E^{2\pi \I \mathrm{Tr}(p^{-1}SZ)}\\
&\qquad=p\sum_{y\in\Z/p\Z}\sum_{\substack{S=\begin{bsmallmatrix} \alpha&\beta\\ \beta&\gamma \end{bsmallmatrix}\in B(N)\\p^n\mid\gamma}}
 a(S) \E^{2\pi \I 
\mathrm{Tr}(\begin{bsmallmatrix}1&yp^{-1}\\&p^{-1}\end{bsmallmatrix}S
\begin{bsmallmatrix}p\vphantom{p^{-1}}&\\y&1\vphantom{p^{-1}}\end{bsmallmatrix}Z)}\\
&\qquad\quad+p^{3-k}\sum_{\substack{S=\begin{bsmallmatrix} \alpha&\beta\\ \beta&\gamma \end{bsmallmatrix}\in B(N)\\p^n\mid\gamma,\ p\mid\alpha,\ p\mid2\beta}} a(S)
\E^{2\pi \I \mathrm{Tr}(p^{-1}SZ)}\\
&\qquad=p\sum_{y\in\Z/p\Z}
\sum_{\substack{S=\begin{bsmallmatrix} \alpha&\beta\\ \beta&\gamma \end{bsmallmatrix}\in B(N)\\p\mid(\alpha-2\beta y+\gamma y^2)}}
 a(p^{-1}S[\begin{bsmallmatrix}1&\\-y&p\end{bsmallmatrix}]) \E^{2\pi \I \mathrm{Tr}(SZ)}\\
&\quad\qquad+p^{3-k}\sum_{\substack{S\in B(N)}}  a(pS) \E^{2\pi \I \mathrm{Tr}(SZ)}\\
&\qquad=\sum_{S=\begin{bsmallmatrix} \alpha&\beta\\ \beta&\gamma \end{bsmallmatrix}\in B(N)} 
\Big( p^{3-k}  a(pS) 
+ p \sum_{\substack{y \in \Z/p\Z \\p\mid(\alpha+2\beta y+\gamma y^2)}}  a(p^{-1}S[\begin{bsmallmatrix}1&\\y&p\end{bsmallmatrix}]) \Big) 
\E^{2\pi \I \mathrm{Tr}(SZ)}.
\end{align*}
This is \eqref{T01fouriereq}. 

Finally, to prove \eqref{T01fourierjacobieq}, we proceed from an equation from the 
proof of \eqref{T01fouriereq}:
\begin{align*}
&(T_{0,1}^{s}(p)F)(Z)\\
&\qquad=\sum_{a,b\in\Z/p\Z}F(\begin{bsmallmatrix}p\vphantom{p^{-1}}&\\a&1\vphantom{p^{-1}}\end{bsmallmatrix}Z
\begin{bsmallmatrix}1&ap^{-1}\\&p^{-1}\end{bsmallmatrix}
+\begin{bsmallmatrix}\vphantom{p^{-1}}&\\ \vphantom{p^{-1}}&bp^{-n}\end{bsmallmatrix})\\
&\qquad\quad+p^{-k}\sum_{a,b,c\in\Z/p\Z}F(p^{-1}Z+\begin{bsmallmatrix}cp^{-1}&ap^{-1}\\ ap^{-1}&bp^{-n}\end{bsmallmatrix})\\
&\qquad=\sum_{a,b\in\Z/p\Z}F(\begin{bsmallmatrix}p\tau&a\tau+z\\a\tau+z&(a^2\tau+2az+\tau'+bp^{-n+1})p^{-1}\end{bsmallmatrix})\\
&\qquad\quad+p^{-k}\sum_{a,b,c\in\Z/p\Z}F(\begin{bsmallmatrix}(\tau+c)p^{-1}&(z+a)p^{-1}\\ (z+a)p^{-1}&(\tau'+bp^{-n+1})p^{-1}\end{bsmallmatrix})\\
&\qquad=\sum_{a,b\in\Z/p\Z}\sum_{\substack{m=0\\ N_s\mid m}}^\infty f_m(p\tau, a\tau+z) \E^{2\pi \I mp^{-1}(a^2\tau+2az+\tau'+bp^{-n+1})}\\
&\qquad\quad+p^{-k}\sum_{a,b,c\in\Z/p\Z}\sum_{\substack{m=0\\ N_s\mid m}}^\infty f_m(p^{-1}(\tau+c),p^{-1}(z+a)) \E^{2\pi \I mp^{-1}(\tau'+bp^{-n+1})}\\
&\qquad=p\sum_{a\in\Z/p\Z}\sum_{\substack{m=0\\ N_s p\mid m}}^\infty f_m(p\tau, a\tau+z) \E^{2\pi  \I mp^{-1}(a^2\tau+2az+\tau')}\\
&\qquad\quad+p^{1-k}\sum_{a,c\in\Z/p\Z}\sum_{\substack{m=0\\ N_s p\mid m}}^\infty f_m(p^{-1}(\tau+c),p^{-1}(z+a)) \E^{2\pi \I mp^{-1}\tau'}\\
&\qquad=\sum_{\substack{m=0\\ N_s p\mid m}}^\infty\Big(p\sum_{a\in\Z/p\Z} f_m(p\tau, a\tau+z)\E^{2\pi \I mp^{-1}(a^2\tau+2az)}\\
&\qquad\quad+p^{1-k}\sum_{a,c\in\Z/p\Z} f_m(p^{-1}(\tau+c),p^{-1}(z+a))\Big) \E^{2\pi \I mp^{-1}\tau'}\\
&\qquad=\sum_{\substack{m=0\\N_s\mid m}}^\infty\Big(p\sum_{a\in\Z/p\Z} f_{mp}(p\tau, a\tau+z) \E^{2\pi \I m(a^2\tau+2az)}\\
&\qquad\quad+p^{1-k}\sum_{a,c\in\Z/p\Z} f_{mp}(p^{-1}(\tau+c),p^{-1}(z+a))\Big) \E^{2\pi \I  m\tau'}\\
&\qquad=\sum_{\substack{m=0\\N_s\mid m}}^\infty(L'_pf_{mp})(\tau,z) \E^{2\pi \I m\tau'}.
\end{align*}
This completes the proof.
\end{proof}

\begin{lemma}\label{T10formulalemma}
Let $N$ and $k$ be integers such that $N>0$ and $k>0$, and let
$p$ be a prime of $\Z$. Assume that $v_p(N) \geq 1$. Let
$T_{1,0}^s: \mathcal{A}_k(\mathcal{K}_s(N)) \to \mathcal{A}_k(\mathcal{K}_s(N))$
be the linear map obtained as in Sect.~\ref{overviewsec} from the Hecke operator \eqref{Ts10eq},
and let $T_{1,0}^{s}(p)=(T_{1,0}^s)_p$
be as in \eqref{overviewTddefAeq}, so that $T_{1,0}^{s}(p)$ is a linear map 
$$
T_{1,0}^{s}(p): M_k(\mathrm{K}_s(N)) \longrightarrow M_k(\mathrm{K}_s(N)). 
$$
Let $F\in M_k(\mathrm{K}_s(N))$ with Fourier expansion as in \eqref{NSFourierexpagaineq}. Then
\begin{align}
T_{1,0}^{s}(p)F&=\sum_{\substack{x,y\in\Z/p\Z\\z\in\Z/p^2\Z}}F\big |_k
\begin{bsmallmatrix}p&&&\\&1&&\\&&p&\\&&&p^2\end{bsmallmatrix}
\begin{bsmallmatrix}1&&&y\\x&1&y&zp^{-n+1}\\&&1&-x\\&&&1\end{bsmallmatrix}.\label{T10slashopeq}
\end{align}
If $F$ is a cusp form, then $T_{1,0}^{s}(p) F$ is a cusp form. We have 
\begin{equation}
(T_{1,0}^{s}(p)F)(Z)= \sum_{\substack{S\in B(N)}} \sum_{\substack{a\in\Z/p\Z}} p^{3-k}
 a(S[\begin{bsmallmatrix}1&\\a&p\end{bsmallmatrix}]) \E^{2\pi \I \mathrm{Tr} (SZ)}.\label{T10fouriereq}
\end{equation} 
The Fourier-Jacobi expansion of $T_{1,0}^{s}(p)F$ is given by 
\begin{equation}
(T_{1,0}^{s}(p)F)(Z)=p^{2}\sum_{\substack{m=0\\N_s\mid m}}^\infty (L_{p^2}f_{mp^2})(\tau,z) \E^{2\pi \I m\tau'},
\qquad Z=\begin{bsmallmatrix} \tau & z \\ z & \tau' \end{bsmallmatrix} \in \mathcal{H}_2. 
\label{T10fourierjacobieq}
\end{equation}
Here $L_{p^2}$ is as in Lemma \ref{Lpoperatorlemma}. 
\end{lemma}	
\begin{proof}
Define $\Phi=\Phi_F$. Let $Z\in\mathcal{H}_2$ 
and let $h\in \GSp(4,\R)^+$ be such that $Z=h\langle I \rangle$. Let $n=v_p(N)$. 
By \eqref{overviewTddefeq} and \eqref{Ts10eq} we have 
\begin{align*}
T_{1,0}^s(p) F(Z)&=\lambda(h)^{-k}j(h,I)^kT^{s}_{1,0}\Phi(h)\\
&=\lambda(h)^{-k}j(h,I)^k\sum_{\substack{x,y\in\Z/p\Z\\z\in\Z/p^2\Z}}
\Phi(h\begin{bsmallmatrix}1&&&y\\x&1&y&zp^{-n+1}\\&&1&-x\\&&&1\end{bsmallmatrix}_p
\begin{bsmallmatrix}p&&&\\&p^2\vphantom{p^{-1}}&&\\&&p&\\&&&1\vphantom{p^{-1}}\end{bsmallmatrix}_p)\\
&=\lambda(h)^{-k}j(h,I)^k\sum_{\substack{x,y\in\Z/p\Z\\z\in\Z/p^2\Z}}
\Phi(
\underbrace{
\begin{bsmallmatrix}
p&&&\\
&p^2\vphantom{p^{-1}}&&\\
&&p&\\
&&&1\vphantom{p^{-1}}
\end{bsmallmatrix}^{-1}}_{\text{all places}}
\underbrace{\begin{bsmallmatrix}1&&&y\\x&1&y&zp^{-n+1}\\&&1&-x\\&&&1\end{bsmallmatrix}^{-1}}_{\text{all places}}\\
&\quad\times h\begin{bsmallmatrix}1&&&y\\x&1&y&zp^{-n+1}\\&&1&-x\\&&&1\end{bsmallmatrix}_p
\begin{bsmallmatrix}p&&&\\&p^2\vphantom{p^{-1}}&&\\&&p&\\&&&1\vphantom{p^{-1}}\end{bsmallmatrix}_p)\\
&=\lambda(h)^{-k}j(h,I)^k\sum_{\substack{x,y\in\Z/p\Z\\z\in\Z/p^2\Z}}
\Phi(\begin{bsmallmatrix}p^{-1}&&&\\&p^{-2}&&\\&&p^{-1}&\\&&&1\end{bsmallmatrix}_\infty
\begin{bsmallmatrix}1&&&-y\\-x&1&-y&-zp^{-n+1}\\&&1&x\\&&&1\end{bsmallmatrix}_\infty h)\\
&=\lambda(h)^{-k}j(h,I)^k\sum_{\substack{x,y\in\Z/p\Z\\z\in\Z/p^2\Z}}
\Phi(\begin{bsmallmatrix}p&&&\\&1&&\\&&p&\\&&&p^2\end{bsmallmatrix}_\infty
\begin{bsmallmatrix}1&&&y\\x&1&y&zp^{-n+1}\\&&1&-x\\&&&1\end{bsmallmatrix}_\infty h)\\
&=\sum_{\substack{x,y\in\Z/p\Z\\z\in\Z/p^2\Z}}(F\big |_k\begin{bsmallmatrix}p&&&\\&1&&\\&&p&\\&&&p^2\end{bsmallmatrix}
\begin{bsmallmatrix}1&&&y\\x&1&y&zp^{-n+1}\\&&1&-x\\&&&1\end{bsmallmatrix})(Z)
\end{align*}
This proves \eqref{T10slashopeq}. To prove \eqref{T10fouriereq} we use \eqref{T10slashopeq}:
\begin{align*}
&(T_{1,0}^{s}(p)F)(Z)=\sum_{\substack{a,b\in\Z/p\Z\\c\in\Z/p^2\Z}}(F\big |_k\begin{bsmallmatrix}p&&&\\&1&&\\&&p&\\&&&p^2\end{bsmallmatrix}
\begin{bsmallmatrix}1&&&b\\a&1&b&cp^{-n+1}\\&&1&-a\\&&&1\end{bsmallmatrix})(Z)\\
&\qquad =p^{-k}\sum_{\substack{a,b\in\Z/p\Z\\c\in\Z/p^2\Z}}F(\begin{bsmallmatrix}p\vphantom{p^{-1}}&\\a&1\vphantom{p^{-1}}\end{bsmallmatrix}Z
\begin{bsmallmatrix}p^{-1}&ap^{-2}\\&p^{-2}\end{bsmallmatrix}+\begin{bsmallmatrix}&p^{-1}b\\ bp^{-1}&abp^{-2}+cp^{-n-1}\end{bsmallmatrix})\\
&\qquad =p^{-k}\sum_{\substack{a,b\in\Z/p\Z\\c\in\Z/p^2\Z}}\sum_{S\in B(N)} a(S) 
\E^{2\pi \I \mathrm{Tr} (S(\begin{bsmallmatrix}p\vphantom{p^{-1}}&\\a&1\vphantom{p^{-1}}\end{bsmallmatrix}Z
\begin{bsmallmatrix}p^{-1}&ap^{-2}\\&p^{-2}\end{bsmallmatrix}
+\begin{bsmallmatrix}&p^{-1}b\\p^{-1}b&abp^{-2}+cp^{-n-1}\end{bsmallmatrix}))}\\
&\qquad =p^{-k}\sum_{S=\begin{bsmallmatrix} \alpha & \beta \\ \beta & \gamma \end{bsmallmatrix} \in B(N)}
\sum_{\substack{a,b\in\Z/p\Z\\c\in\Z/p^2\Z}} a(S) \E^{2\pi \I (2\beta b p^{-1}
+\gamma abp^{-2}+\gamma cp^{-n-1})}\\
&\qquad \quad\times \E^{2\pi \I \mathrm{Tr} (\begin{bsmallmatrix}p^{-1}&ap^{-2}\\&p^{-2}\end{bsmallmatrix}
S\begin{bsmallmatrix}p\vphantom{p^{-1}}&\\a&1\vphantom{p^{-1}}\end{bsmallmatrix}   Z)}\\
&\qquad =p^{2-k}\sum_{\substack{S=\begin{bsmallmatrix} \alpha & \beta \\ \beta & \gamma \end{bsmallmatrix}\in B(N)\\p^{n+1}\mid\gamma}}
\sum_{\substack{a,b\in\Z/p\Z}} a(S) \E^{2\pi \I (2\beta b  p^{-1})}
\E^{2\pi \I \mathrm{Tr} (\begin{bsmallmatrix}p^{-1}&ap^{-2}\\&p^{-2}\end{bsmallmatrix}
S\begin{bsmallmatrix}p\vphantom{p^{-1}}&\\a&1\vphantom{p^{-1}}\end{bsmallmatrix}   Z)}\\
&\qquad =p^{3-k}\sum_{\substack{a\in\Z/p\Z}}
\sum_{\substack{S=\begin{bsmallmatrix} \alpha & \beta \\ \beta & \gamma \end{bsmallmatrix}\in B(N)\\
p^{n+1}\mid\gamma,\  p\mid 2\beta}} a(S) \E^{2\pi \I \mathrm{Tr} (\begin{bsmallmatrix}p^{-1}&ap^{-2}\\&p^{-2}\end{bsmallmatrix}
S\begin{bsmallmatrix}p\vphantom{p^{-1}}&\\a&1\vphantom{p^{-1}}\end{bsmallmatrix}   Z)}\\
&\qquad =\sum_{\substack{S\in B(N)}}p^{3-k}\sum_{\substack{a\in\Z/p\Z}} a(S[\begin{bsmallmatrix}1&\\a&p\end{bsmallmatrix}]) \E^{2\pi \I \mathrm{Tr} (SZ)}.
\end{align*}
This is \eqref{T10fouriereq}. 

Finally, to prove \eqref{T10fourierjacobieq}, we proceed from an equation from the 
proof of \eqref{T10fouriereq}:
\begin{align*}
&(T_{1,0}^{s}(p)F)(Z)= p^{-k}\sum_{\substack{a,b\in\Z/p\Z\\c\in\Z/p^2\Z}}F(\begin{bsmallmatrix}p\vphantom{p^{-1}}&\\a&1\vphantom{p^{-1}}\end{bsmallmatrix}Z
\begin{bsmallmatrix}p^{-1}&ap^{-2}\\&p^{-2}\end{bsmallmatrix}+\begin{bsmallmatrix}&p^{-1}b\\ bp^{-1}&abp^{-2}+cp^{-n-1}\end{bsmallmatrix})\\
&\qquad =p^{-k}\sum_{\substack{a,b\in\Z/p\Z\\c\in\Z/p^2\Z}}
F(\begin{bsmallmatrix}\tau&(a\tau+b+z)p^{-1}\\(a\tau+b+z)p^{-1}&(a^2\tau+2az+\tau'+ab+cp^{-n+1})p^{-2}\end{bsmallmatrix})\\
&\qquad =p^{-k}\sum_{\substack{a,b\in\Z/p\Z\\c\in\Z/p^2\Z}}\sum_{\substack{m=0\\N_s\mid m}}^\infty f_m(\tau,(a\tau+b+z)p^{-1})
\E^{2\pi \I m(a^2\tau+2az+\tau'+ab+cp^{-n+1})p^{-2}}\\
&\qquad =p^{2-k}\sum_{\substack{a,b\in\Z/p\Z}}\sum_{\substack{m=0\\p^2N_s\mid m}}^\infty f_m(\tau,(a\tau+b+z)p^{-1})
\E^{2\pi \I mp^{-2}(a^2\tau+2az)} \E^{2\pi \I mp^{-2}\tau'}\\
&\qquad =p^{2-k}\sum_{\substack{m=0\\N_s\mid m}}^\infty\sum_{\substack{a,b\in\Z/p\Z}} f_{mp^2}(\tau,(b+ a\tau+z)p^{-1})e^{2\pi i m(2az+a^2\tau)}\E^{2\pi \I m\tau'}\\
&\qquad =p^{2}\sum_{\substack{m=0\\N_s\mid m}}^\infty (L_{p^2}f_{mp^2})(\tau,z) \E^{2\pi \I m\tau'},
\end{align*}
where $L_{p^2}$ is defined in \eqref{Lcoperatoreq}. This completes the proof.
\end{proof}

\section{Some relations between operators}
We conclude this chapter by stating some relations between the operators on stable Klingen forms introduced above; 
these relations follow from local statements proved in Chap.~\ref{basicfactschapter}.
\begin{proposition}\label{globalrelationsprop}
 Let $N$ and $k$ be integers such that $N>0$ and $k>0$, and let $p$ be a prime of $\Z$. Then, for $F\in M_k(\mathrm{K}_s(N))$,
 \begin{align}
  \label{globalrelationspropeq1}\theta_p\tau_pF&=\tau_p\theta_pF,\\
  \label{globalrelationspropeq2}T_{0,1}^s(p)T_{1,0}^s(p)F&=T_{1,0}^s(p)T_{0,1}^s(p)F\qquad\text{if }p\mid N,\\
  \label{globalrelationspropeq3}T_{0,1}^s(p)F&= p^2\sigma_p \theta_pF\qquad\text{if }p\mid N,\\
  \label{globalrelationspropeq4}T_{1,0}^s(p)F&= p^4 \sigma_p \tau_pF\qquad\text{if }p\mid N,\\
  \label{globalrelationspropeq5}T_{1,0}^s(p)F&= p^4 \tau_p \sigma_pF\qquad\text{if }p^2\mid N,\\
  \label{globalrelationspropeq6}T^s_{0,1}(p)\tau_pF&=\tau_pT^s_{0,1}(p)F\qquad\text{if }p\mid N,\\
  \label{globalrelationspropeq7}T^s_{1,0}(p)\tau_pF&=\tau_pT^s_{1,0}(p)F\qquad\text{if }p\mid N,\\
  \label{globalrelationspropeq8}T^s_{1,0}(p)\theta_pF&=p^2T^s_{0,1}(p)\tau_pF,\\
  \label{globalrelationspropeq9}T^s_{0,1}(p)\theta_pF&=\theta_p T^s_{0,1}(p)F+p^3\tau_pF-p^3\eta_p\sigma_pF\qquad\text{if }p^2\mid N.
 \end{align}
\end{proposition}
\begin{proof}
These statements follow from the corresponding local statements contained in 
Lemma~\ref{Heckecommutelemma}, Lemma~\ref{heckeupdownlemma}, and Lemma \ref{Ts01thetalemma}. 
For the connection to representation theory we use the space $W$ from Sect.~\ref{overviewsec}.
\end{proof}

%% file: SKMS_chapter12.tex
\chapter{Hecke Eigenvalues and Fourier Coefficients}
\label{heckechap}

In this final chapter we present some applications of the local 
theory of Part~1 to the Hecke eigenvalues and Fourier coefficients
of Siegel modular newforms $F$ in $S_k(\mathrm{K}(N))_{\mathrm{new}}$ of degree
two with paramodular level $N$. Assuming that $F$ is an eigenform
for the Hecke operators $T(1,1,p,p)$ and $T(1,p,p,p^2)$
for all primes $p$ of $\Z$,\footnote{In fact, assuming that $F$ is
an eigenform for $T(1,1,p,p)$ and $T(1,p,p,p^2)$ for all but
finitely many $p$ is sufficient: see Theorem \ref{Ralftheorem} below.}
we begin by proving in Sect.~\ref{appsec} that the local results
from Part~1 imply identities involving $F$ and its images under
the upper block operators from Chap.~\ref{operchap} at $p$ for $p^2 \mid N$. 
We then show in Corollary \ref{fourierheckethm} that these identities
yield relations between Fourier coefficients and Hecke eigenvalues as well as
conditions which determine properties of the attached local representations at $p$. 
Corollary \ref{fourierheckethm} may be regarded as a solution to the problem
mentioned in the introduction to this work. In the second section we show that 
the formulas of Corollary \ref{fourierheckethm} can be rewritten in terms
of the action of the Hecke ring of $\Gamma_0(N)$ on the vector space of $\C$
valued functions on the set of positive semi-definite $2 \times 2$ matrices
with rational entries; these results represent a more conceptual
and compact statement of the formulas of Corollary~\ref{fourierheckethm}. 
We conclude this chapter with two applications of Corollary~\ref{fourierheckethm}.
First, in Sect.~\ref{examplessec} and Sect.~\ref{eigenvaluescompsec} we show that 
the equations of Corollary \ref{fourierheckethm} do indeed hold for the examples of
\cite{PSY} and  \cite{PSYW}, and we indicate how the equations could be used to 
calculate Hecke eigenvalues from Fourier coefficients in other instances.  Finally, in Sect.~\ref{recurrencesec}
we apply Corollary \ref{fourierheckethm} to prove that the radial Fourier
coefficients $a(p^t S)$ for $t \geq 0$ and $p^2 \mid N$ satisfy a recurrence
relation determined by the spin $L$-factor of $F$ at $p$. This extends results
known in other cases (see, for example,  Sect.~4.3.2 of \cite{MR884891}).

\section{Applications}
\label{appsec}

In this section we apply some of the local results of Part~1 to the Siegel modular newforms 
with paramodular level. To connect such newforms to representation theory we will use 
the following theorem. In this statement
 $\mathfrak{g}'$ is the Lie algebra of $\GSp(4,\R)$ and $K^{\pm}$ is the maximal compact subgroup of $\GSp(4,\R)$
as in \cite{S1}. See Chap.~\ref{backSMFchap} for general definitions about Siegel modular forms.

\begin{theorem}
\label{Ralftheorem}
Let $N$ and $k$ be integers such that $N>0$ and $k>0$, and let $F \in S_k(\mathrm{K}(N))$. 
Assume that $F$ is a newform and is an eigenvector for the Hecke operators $T(1,1,q,q)$ and $T(1,q,q,q^2)$
for all but finitely many of the primes $q$ of $\Z$ such that $q \nmid N$. Then $F$ is an 
eigenvector for $T(1,1,q,q)$, $T(1,q,q,q^2)$, and $w_q$  for all primes $q$ of $\Z$. Moreover, 
let $V$ be the $\GSp(4,\A_{\mathrm{fin}}) \times (\mathfrak{g}',K^{\pm})$ submodule  of
the space of cuspidal automorphic forms on $\GSp(4,\A)$ generated by $\Phi_F$, with
$\Phi_F$ as in Lemma \ref{twospacesisolemma}. Then  $V$ is irreducible, so that $V\simeq\otimes_{v\leq\infty}\pi_v$, and 
where, for each prime $q$ of $\Z$, $\pi_q$ is an irreducible, admissible representation of $\GSp(4,\Q_q)$ with trivial central character,
and $\pi_\infty$ is an irreducible, admissible $(\mathfrak{g}',K^{\pm})$ module.  Under this isomorphism, $\Phi_F$ corresponds to a pure tensor $\otimes_{v\leq\infty}x_v$, where $x_q$ is a newform in $\pi_q$ for each prime $q$ of $\Z$. Let $q$ be a prime of 
$\Z$, and let $\lambda_q,\mu_q \in \C$ and $\varepsilon_q \in \{ \pm 1\}$ be such that 
$T(1,1,q,q)F= q^{k-3}\lambda_q F$, $T(1,q,q,q^2) F = q^{2(k-3)} \mu_q F$, and $w_q F = \varepsilon_q F$. 
Then $T_{0,1} x_q = \lambda_q x_q$, $T_{1,0} x_q = \mu_q x_q$, and $u_{v_q(N)} x_q = \varepsilon_q x_q$. 
The representation $\pi_q$ is non-generic for some prime $q$ of $\Z$ if and only if 
$\pi_q$ is a Saito-Kurokawa representation for all primes $q$ of $\Z$. 
\end{theorem}
\begin{proof}
See Theorem 5.5 and Corollary 5.4 of \cite{S3}.
\end{proof}

Let $k$ and $N$ be positive integers. Using this theorem and further results from \cite{S2} and \cite{S3}, one may describe the space $S_k(\mathrm{K}(N))_{\rm new}$ as follows. The Hecke operators $T(1,1,q,q)$ and $T(1,q,q,q^2)$ for primes $q\nmid N$ act on this space and commute; see Lemma~\ref{TTUpreservelemma}. We thus may diagonalize $S_k(\mathrm{K}(N))_{\rm new}$ according to this action. Each eigenspace is $1$-dimensional, and any eigenform is also an eigenform for the Hecke operators at places $p\mid N$. Each eigenform corresponds to an automorphic representation as in Theorem~\ref{Ralftheorem}. Such eigenforms come in two types, known as lifts\index{lifts} and non-lifts\index{non-lifts}. The lifts are modular forms of Saito-Kurokawa type (type {\bf(P)}), and the non-lifts are modular forms of general type (type {\bf(G)}) in the sense of \cite{S2}. The lifts can be characterized in a number of ways; in particular, $F$ is a lift if and only if $\pi_q$ is non-generic for one (then all) primes $q$ of~$\Z$. (See \cite{FPRS} for a survey of further characterizations.)

We now apply the local results from Part~1 to Siegel modular forms with paramodular level.

\begin{theorem}\label{globalalgthm}
Let $N$ and $k$ be integers such that $N>0$ and $k>0$, and let $F \in S_k(\mathrm{K}(N))$. 
Assume that $F$ is a newform and is an eigenvector for the Hecke operators $T(1,1,q,q)$ and $T(1,q,q,q^2)$
for all but finitely many of the primes $q$ of $\Z$ such that $q \nmid N$; then by Theorem
\ref{Ralftheorem}, $F$ is an eigenvector for $T(1,1,q,q)$ and $T(1,q,q,q^2)$ for all primes $q$ of $\Z$.
Let $\otimes_{v \leq \infty} \pi_v$ be as in Theorem \ref{Ralftheorem}. 
For every prime $q$ of $\Z$ let $\lambda_q, \mu_q \in \C$ be such that 
\begin{align}
T(1,1,q,q)F&=q^{k-3}\lambda_qF,\label{globallambdaeq}\\
T(1,q,q,q^2)F&=q^{2(k-3)}\mu_qF\label{globalmueq}.
\end{align}
Let $p$ be a prime of $\Z$ with $v_p(N)\geq2$. Then 
\begin{equation}
\label{globalalgthmeq01}
\mu_p = 0 \Longleftrightarrow \sigma_p F = 0 \Longleftrightarrow T_{1,0}^s(p) F = 0.
\end{equation}
Moreover:
\begin{enumerate}
\item \label{globalalgthmitem15} If $v_p(N) \geq 3$, then $\sigma_p^2 F =0$.
\item \label{globalalgthmitem16} We have
\begin{equation}
\label{globalnewformreleq}
\mu_p F = p^4 \tau_p^2 \sigma_p F - p^2 \eta_p \sigma_p F.
\end{equation}
\item \label{globalalgthmitem2} Assume that $\mu_p\neq 0$. Then  $\sigma_p F \neq 0$, and 
\begin{align}
T^{s}_{0,1}(p)(\sigma_p F)&=\lambda_p(\sigma_pF),\label{stablelambdacase1eq}\\
T^{s}_{1,0}(p)(\sigma_p F)&=(\mu_p+p^2)(\sigma_pF)\label{stablemucase1eq},
\end{align}
and the representation $\pi_p$ is generic.
\item \label{globalalgthmitem3} Assume that $\mu_p=0$. Then
\begin{equation}
\label{sigmazeromuzeroeq}
T_{1,0}^s(p) F =0.
\end{equation}
If $F$ is not an eigenvector for $T^{s}_{0,1}(p)$, then $\pi_p$ is generic and 
\begin{align}
T^{s}_{0,1}(p)^2F=-p^3F+\lambda_p T^{s}_{0,1}(p)F. \label{stablelambdacase3eq}
\end{align}
The newform $F$ is an eigenvector for $T^{s}_{0,1}(p)$ if and only if 
\begin{align}
T^{s}_{0,1}(p)F=(1+p^{-1})^{-1}\lambda_pF;\label{stablelambdacase2eq}
\end{align}
in this case $\pi_p$ is non-generic. 
\end{enumerate}
\end{theorem}
\begin{proof}
Let $\Phi_F$ and $V \cong \otimes_{v \leq \infty} \pi_v$ be as in Theorem \ref{Ralftheorem}. For each place $v$
of $\Q$, let $V_v$ be the space of $\pi_v$. By Theorem \ref{Ralftheorem}, $\Phi_F$ corresponds to a 
pure tensor $\otimes_v x_v$ under the fixed isomorphism $i: \otimes_v V_v \stackrel{\sim}{\to} V$. 
Define $V_p \hookrightarrow V$  by $x\mapsto i(x\otimes \otimes_{v \neq p} x_v)$
for $x \in V_p$; then $V_p \hookrightarrow V$ is a well-defined injective $\GSp(4,\Q_p)$ map. 
Moreover, the properties of $x_\infty$, and of $x_q$ for $q$ a  prime of $\Z$ such that $q \neq p$,
imply that the image of this map is contained in 
\begin{equation}
\label{Npseteq}
W = \{\, \Phi \in \mathcal{A}_k^\circ\mid \kappa \cdot \Phi = \Phi \ \text{for $\kappa \in \prod_{q \neq p} \mathrm{K}_s(\mathfrak{q}^{v_q(N)})$}\, \},
\end{equation}
which is also a $\GSp(4,\Q_p)$ space.
We thus obtain an injective $\GSp(4,\Q_p)$  map 
\begin{equation}
\label{Npinceq}
t: V_p \hookrightarrow W.
\end{equation}
Let $N'=Np^{-v_p(N)}$ so that $p$ and $N'$ are relatively prime.
Let $n_1$ and $n_2$ be integers such that $n_1,n_2 \geq 0$, and define $N_1 = N' p^{n_1}$ and $N_2= N' p^{n_2}$.  
For $i=1$ or $2$, since the subspace of $\mathrm{K}_s(\p^{n_i})$-invariant  vectors in $W$ is $\mathcal{A}_k^\circ(\mathcal{K}_s(N_i))$
and since $t$ is a $\GSp(4,\Q_p)$ map, the map $t$ restricts to give inclusions
$$
t: V_{p,s}(\p^{n_1}) \hookrightarrow \mathcal{A}_k^\circ (\mathcal{K}_s(N_1)) \quad \text{and} \quad 
t: V_{p,s}(\p^{n_2}) \hookrightarrow \mathcal{A}_k^\circ (\mathcal{K}_s(N_2))
$$
Next, let $g \in \GSp(4,\Q_p)$, and let 
$$
\delta=T_g: V_{p,s}(\p^{n_1}) \longrightarrow V_{p,s}(\p^{n_2})\quad \text{and} \quad
\delta =T_g: W_{p,s}(\p^{n_1}) \longrightarrow W_{p,s}(\p^{n_2}) 
$$
be the  operators corresponding to $g$ 
as in Sect.~\ref{opersksec}.
Since $t$ is a $\GSp(4,\Q_p)$ map,  there is a commutative diagram
\begin{equation}
\label{firstTCDeq}
\begin{CD}
V_{p,s}(\p^{n_1})@>t>>   \mathcal{A}_k^\circ (\mathcal{K}_s(N_1))  \\
@V \delta VV @VV  \delta V \\
V_{p,s}(\p^{n_2})@>t>>   \mathcal{A}_k^\circ (\mathcal{K}_s(N_2))  
\end{CD} 
\end{equation}
Also, from Sect.~\ref{overviewsec} there is a commutative diagram
\begin{equation}
\label{secondTCDeq}
\begin{CD}
\mathcal{A}_k^\circ (\mathcal{K}_s(N_1))@>r>\sim>   S_k (\mathrm{K}_s(N_1))  \\
@V  \delta VV @VV  \delta_p V \\
\mathcal{A}_k^\circ (\mathcal{K}_s(N_2))@>r>\sim>   S_k (\mathrm{K}_s(N_2))  
\end{CD}
\end{equation}
Combining \eqref{firstTCDeq} with \eqref{secondTCDeq}, we  obtain a commutative diagram
\begin{equation}
\label{thirdTCDeq}
\begin{CD}
V_{p,s}(\p^{n_1})@>r\circ t>>   S_k (\mathrm{K}_s(N_1))  \\
@V  \delta VV @VV  \delta_p V \\
V_{p,s}(\p^{n_2})@>r\circ t>>   S_k (\mathrm{K}_s(N_2))  
\end{CD}
\end{equation}
We note that the map $r \circ t$ is injective. 
Also, for the following arguments we remind the reader
that $v_p(N)$ is the paramodular level $N_{\pi_p}$ of $\pi_p$, 
$x_p$ is a newform in $V_p(\p^{v_p(N)})\subset V_{p,s}(\p^{v_p(N)})$, 
$\lambda_{\pi_p} =\lambda_p$, $\mu_{\pi_p}=\mu_p$, 
and  $F=(r \circ t) x_p$ (these statements follow from Theorem \ref{Ralftheorem}).  

To prove that $\mu_p = 0$  if and only if  $\sigma_p F= 0$
we use the diagram \eqref{thirdTCDeq} with $n_1=v_p(N)$, $n_2=v_p(N)-1$ and 
$$
\delta=\sigma=\sigma_{v_p(N)-1}: V_{p,s}(\p^{v_p(N)}) \longrightarrow V_{p,s}(\p^{v_p(N)-1})
$$
to obtain the commutative diagram
\begin{equation}
\label{thirdsigmaCDeq}
\begin{CD}
V_{p,s}(\p^{v_p(N)})@>r\circ t>>   S_k (\mathrm{K}_s(N))  \\
@V  \sigma VV @VV  \sigma_p V \\
V_{p,s}(\p^{v_p(N)-1})@>r\circ t>>   S_k (\mathrm{K}_s(Np^{-1}))  
\end{CD}
\end{equation}
We have $x_p \in V_{p,s}(\p^{v_p(N)})$ and $(r\circ t)(x_p)=F$. 
The commutativity of  \eqref{thirdsigmaCDeq}  implies that $\sigma_p F = 
(r \circ t)(\sigma x_p)$.  Since $r \circ t$ is injective it follows that 
$\sigma_p F = 0$ if and only if $\sigma x_p = 0$. 
By the equivalence of \ref{shadowsigmavnewpropitem2} and \ref{shadowsigmavnewpropitem3} of Proposition 
\ref{shadowsigmavnewprop} we have $\sigma x_p = 0$ if and only if $\mu_{\pi_p} =\mu_p = 0$, 
so that $\mu_p = 0$  if and only if $\sigma_p F = 0$. 
Next, by \eqref{thirdTCDeq} with $n_1=n_2=v_p(N)$ and $\delta=T_{1,0}^s$ we have
the commutative diagram
\begin{equation}
\label{T22CDeq}
\begin{CD}
V_{p,s}(\p^{v_p(N)})@>r\circ t>>   S_k (\mathrm{K}_s(N))  \\
@V T_{1,0}^s VV @VV  T_{1,0}^s(p) V \\
V_{p,s}(\p^{v_p(N)})@>r\circ t>>   S_k (\mathrm{K}_s(N))  
\end{CD}
\end{equation}
The commutativity of \eqref{T22CDeq} implies that $T_{1,0}^s(p)F = (r \circ t)(T_{1,0}^s x_p)$. 
By the injectivity of $r\circ t$ we have $T_{1,0}^s(p)F = 0$ if and only if $T_{1,0}^s x_p = 0$. 
By the equivalence of \ref{shadowsigmavnewpropitem2} and \ref{shadowsigmavnewpropitem5} of Proposition 
\ref{shadowsigmavnewprop} we have $\mu_{\pi_p} =\mu_p= 0$ if and only if $T_{1,0}^s x_p = 0$. 
Therefore, $\mu_p = 0$ if and only if $T_{1,0}^s(p)F = 0$. 

\ref{globalalgthmitem15} Assume that $v_p(N) \geq 3$. We combine diagram \eqref{thirdTCDeq}
with $n_1=v_p(N)$, $n_2=v_p(N)-1$, and $\delta=\sigma=\sigma_{v_p(N)-1}$ and diagram \eqref{thirdTCDeq} 
with $n_1=v_p(N)-1$, $n_2=v_p(N)-2$, and $\delta=\sigma=\sigma_{v_p(N)-2}$ to obtain the commutative 
diagram 
$$
\begin{CD}
V_{p,s}(\p^{v_p(N)})@>r \circ t>>   S_k (\mathrm{K}_s(N))  \\
@V  \sigma^2 VV @VV  \sigma_p^2 V \\
V_{p,s}(\p^{v_p(N)-2})@>r \circ t>>   S_k (\mathrm{K}_s(Np^{-2}))  
\end{CD}
$$
Since $v_p(N)=N_{\pi_p}$ we have $V_{p,s}(\p^{v_p(N)-2})=0$ by Corollary \ref{upperboundpropcor}. 
This implies that $\sigma_p^2 ( (r\circ t)(x_p))=0$; since $F=(r \circ t)(x_p)$, we obtain
$\sigma_p^2 F=0$. 

\ref{globalalgthmitem16}. Letting $n_1=v_p(N)-1$, $n_2=v_p(N)$, and 
$$
\delta=\tau=\tau_{v_p(N)-1}: V_{p,s}(\p^{v_p(N)-1}) \longrightarrow V_{p,s}(\p^{v_p(N)})
$$
in \eqref{thirdTCDeq}, we obtain the commutative diagram
\begin{equation}
\label{16CDeq1}
\begin{CD}
V_{p,s}(\p^{v_p(N)-1})@>r \circ t>>   S_k (\mathrm{K}_s(Np^{-1}))  \\
@V  \tau VV @VV  \tau_p V \\
V_{p,s}(\p^{v_p(N)})@>r \circ t>>   S_k (\mathrm{K}_s(N))  
\end{CD}
\end{equation}
Similarly, there are commutative diagrams
\begin{equation}
\label{16CDeq2}
\begin{CD}
V_{p,s}(\p^{v_p(N)})@>r \circ t>>   S_k (\mathrm{K}_s(N))  \\
@V  \tau VV @VV  \tau_p V \\
V_{p,s}(\p^{v_p(N)+1})@>r \circ t>>   S_k (\mathrm{K}_s(Np))  
\end{CD}
\end{equation}
and
\begin{equation}
\label{16CDeq3}
\begin{CD}
V_{p,s}(\p^{v_p(N)-1})@>r \circ t>>   S_k (\mathrm{K}_s(Np^{-1}))  \\
@V  \eta VV @VV  \eta_p V \\
V_{p,s}(\p^{v_p(N)+1})@>r \circ t>>   S_k (\mathrm{K}_s(Np))  
\end{CD}
\end{equation}
By \eqref{vnewalluppereq} of Lemma \ref{Wsgenericlemma1} we have 
$q^{-1} \mu_p x_p = q^3 \tau^2 \sigma x_p - q \eta\sigma x_p$. Applying 
${r\circ t}$ to this equation and using the commutativity of
\eqref{16CDeq1}, \eqref{16CDeq2}, and \eqref{16CDeq3}
yields~\eqref{globalnewformreleq}.

\ref{globalalgthmitem2}. Assume that $\mu_p  \neq 0$. Then $\sigma_p F \neq 0$
by \eqref{globalalgthmeq01}. 
By \eqref{thirdTCDeq} with $n_1=n_2=v_p(N)-1$ and $\delta=T_{0,1}^s,T_{1,0}^s$ we have
the commutative diagrams
\begin{equation}
\label{TCDeq}
\begin{CD}
V_{p,s}(\p^{v_p(N)-1})@>r\circ t>>   S_k (\mathrm{K}_s(Np^{-1}))  \\
@V T_{0,1}^s,\ T_{1,0}^s VV @VV  T_{0,1}^s(p),\ T_{1,0}^s(p) V \\
V_{p,s}(\p^{v_p(N)-1})@>r\circ t>>   S_k (\mathrm{K}_s(Np^{-1}))  
\end{CD}
\end{equation}
Since $\mu_{\pi_p}=\mu_p \neq 0$ we have that $\pi_p$ is a \catone\index{category 1 representation}\index{representation!category 1}representation by Proposition \ref{shadowsigmavnewprop}.
By \eqref{Wsgenericlemma1eq1} 
and \ref{TSKleveltheoremitem1} of Theorem \ref{TSKleveltheorem} we have that $T_{0,1}^s \sigma x_p = \lambda_p \sigma x_p$
and $T_{1,0}^s \sigma x_p = (\mu_p + p^2) \sigma x_p$. 
Applying $r \circ t$ to these equations and using 
the commutativity of  \eqref{thirdsigmaCDeq}
and \eqref{TCDeq}, we obtain 
\eqref{stablelambdacase1eq} and \eqref{stablemucase1eq}.  The remaining claim follows
from Lemma \ref{nogenorIVblemma} and the fact that representations of type IVb are non-unitary while
$\pi_p$ is necessarily unitary.

\ref{globalalgthmitem3}.  Assume that $\mu_p F = 0$. By \eqref{globalalgthmeq01} we have $T_{1,0}^s(p) F=0$ which is \eqref{sigmazeromuzeroeq}. 
Since $\mu_{\pi_p}=\mu_p = 0$ we have that $\pi_p$ is a \cattwo\index{category 2 representation}\index{representation!category 2}representation by Proposition~\ref{shadowsigmavnewprop}. 
By \eqref{thirdTCDeq} with $n_1=n_2=v_p(N)$ and $\delta=T_{0,1}^s$ we have
the commutative diagram
\begin{equation}
\label{T2CDeq}
\begin{CD}
V_{p,s}(\p^{v_p(N)})@>r\circ t>>   S_k (\mathrm{K}_s(N))  \\
@V T_{0,1}^s VV @VV  T_{0,1}^s(p) V \\
V_{p,s}(\p^{v_p(N)})@>r\circ t>>   S_k (\mathrm{K}_s(N))  
\end{CD}
\end{equation}
Assume that $F$ is not an eigenform for $T_{0,1}^s(p)$.
Then by \ref{TSKleveltheoremitem2} of Theorem \ref{TSKleveltheorem} the representation $\pi_p$ is generic and 
 $(T_{0,1}^s)^2 x_p = -p^3 x_p + \lambda_p T_{0,1}^s x_p$. 
Applying $r \circ t$ to this equation and using the 
 commutativity of \eqref{T2CDeq}
yields \eqref{stablelambdacase3eq}.
Assume that $F$ is an eigenform for $T_{0,1}^s(p)$. Then by \ref{TSKleveltheoremitem2} of Theorem \ref{TSKleveltheorem}
the representation $\pi_p$ is non-generic and 
 $T_{0,1}^s x_p = (1+p^{-1})^{-1} \lambda_p x_p$. 
Again applying $r \circ t$ to this equation  and using 
the commutativity of \eqref{T2CDeq},
we obtain \eqref{stablelambdacase2eq}.
\end{proof}

We now combine Theorem~\ref{globalalgthm} and the formulas from Chap.~\ref{operchap} 
to obtain relations between the Fourier coefficients and the Hecke eigenvalues of a paramodular newform.

\begin{corollary}\label{fourierheckethm}
Let $N$ and $k$ be integers such that $N>0$ and $k>0$, and let $F \in S_k(\mathrm{K}(N))$. 
Assume that $F$ is a newform and is an eigenvector for the Hecke operators $T(1,1,q,q)$ and $T(1,q,q,q^2)$
for all but finitely many of the primes $q$ of $\Z$ such that $q \nmid N$;  by Theorem
\ref{Ralftheorem}, $F$ is an eigenvector for $T(1,1,q,q)$, $T(1,q,q,q^2)$ and $w_q$ for all primes $q$ of $\Z$.
Let $\lambda_q, \mu_q \in \C$ and $\varepsilon_q \in \{\pm 1\}$ be such that 
\begin{align*}
T(1,1,q,q)F&=q^{k-3}\lambda_qF,\\
T(1,q,q,q^2)F&=q^{2(k-3)}\mu_qF,\\
w_qF&=\varepsilon_q F
\end{align*}
for all primes $q$ of $\Z$. Regard $F$ as an element of $S_k(\mathrm{K}_s(N))$, and let 
\begin{align*}
F(Z)=\sum_{S\in B(N)^+}a(S) \E^{2\pi \I \mathrm{Tr}(SZ)}
\end{align*}
be the Fourier expansion of $F$. Let $\pi\cong \otimes_{v \leq \infty} \pi_v$ be as in Theorem \ref{globalalgthm}. 
Let $p$ be a prime of $\Z$ with $v_p(N)\geq2$. Then
\begin{gather}
\sum_{a \in \Z/p\Z} a(S[ \begin{bsmallmatrix} 1& \\ a&p \end{bsmallmatrix} ] ) =0 \qquad
\text{for $S \in B(Np^{-1})^+$}  \nonumber\\
\Updownarrow\label{fourierheckethmeq1000}\\
\mu_p = 0 \nonumber\\
\Updownarrow\label{fourierheckethmeq1001}\\
\sum_{a \in \Z/p\Z} a(S[ \begin{bsmallmatrix} 1& \\ a&p \end{bsmallmatrix} ] ) =0 \qquad 
\text{for $S \in B(N)^+$}. \nonumber
\end{gather}
Moreover:
\begin{enumerate}
\item \label{fourierheckethmsigmassquareitem}
If $v_p(N) \geq 3$ and $S\in B(Np^{-2})^+$, then 
\begin{equation}
\label{fourierheckethmsigmassquareitemeq}
\sum_{z \in \Z/ p^2 \Z}  a(S[\begin{bsmallmatrix} 1& \\ z&p^2 \end{bsmallmatrix}]) = 0.
\end{equation}
\item \label{fourierheckethmmufinderitem}
If $S= \begin{bsmallmatrix} \alpha&\beta\\ \beta&\gamma \end{bsmallmatrix} \in B(Np)^+$, then 
\begin{equation}\label{mufouriereq}
\mu_p  a(S)=
\begin{cases}
\sum\limits_{x\in\Z/p\Z} p^{3-k}  a(S[\begin{bsmallmatrix}1&\\x&p\end{bsmallmatrix}])&\text{if $p\nmid 2\beta$,}\\
\sum\limits_{x\in\Z/p\Z} p^{3-k}  a(S[\begin{bsmallmatrix}1&\\x&p\end{bsmallmatrix}])
-\sum\limits_{x\in\Z/p\Z} p  a(S[\begin{bsmallmatrix}1&\\xp^{-1}&1\end{bsmallmatrix}])
&\text{if $p\mid 2\beta$}.
\end{cases}
\end{equation}
\item \label{fourierheckethmsigmaneqzeroitem}
Assume that $\mu_p \neq 0$. 
If $S= \begin{bsmallmatrix} \alpha&\beta\\ \beta&\gamma \end{bsmallmatrix}\in B(Np^{-1})^+$, then 
\begin{align}
\lambda_p \sum_{x\in\Z/p\Z} a(S[\begin{bsmallmatrix}1&\\x&p\end{bsmallmatrix}])
&=\sum_{x\in\Z/p\Z}p^{3-k}   a(pS[\begin{bsmallmatrix}1&\\x&p\end{bsmallmatrix}])\nonumber \\
&\quad+\sum_{\substack{z\in\Z/p^2\Z\\p\mid(\alpha+2\beta z+\gamma z^2)}}p  a(p^{-1}S[\begin{bsmallmatrix}1&\\z&p^2\end{bsmallmatrix}])
\label{munozerolambdaeq}
\end{align}
and 
\begin{equation}
\sum_{y\in\Z/p^2\Z}  a(S[\begin{bsmallmatrix}1&\\ y&p^2\end{bsmallmatrix}])
=\begin{cases}\varepsilon_p \sum\limits_{x\in\Z/p\Z}p^{k-2}  a(S[\begin{bsmallmatrix}1&\\x&p\end{bsmallmatrix}])&\text{if } v_p(N)=2,\\
0&\text{if } v_p(N)>2.
\end{cases}\label{epsilonfourierformulaeq}
\end{equation}
The representation $\pi_p$ is generic.
\item Assume that $\mu_p=0$. \label{fourierheckethmsigmazeroitem}
Then 
\begin{equation}
\sum_{x \in \Z/p\Z} a(S[ \begin{bsmallmatrix} 1& \\ x&p \end{bsmallmatrix} ] ) =0\qquad \text{for $S \in B(Np^{-1})^+$.}
\label{fourierheckethmeq1002}
\end{equation}
We have $T_{0,1}^s(p)F = \sum_{S \in B(N)^+} c(S) \E^{2\pi \I \mathrm{Tr}(SZ)}$ where
\begin{equation}
\label{fourierheckethmeq10021}
c(S)= p^{3-k}  a(pS)+ \sum_{\substack{x\in\Z/p\Z\\ p\mid (\alpha+2\beta x)}}
p  a(p^{-1}S[\begin{bsmallmatrix}1&\\ x&p\end{bsmallmatrix}])
\end{equation}
for $S=\begin{bsmallmatrix} \alpha&\beta\\ \beta&\gamma \end{bsmallmatrix} \in B(N)^+$. 
If $F$ is not an eigenvector for  $T^{s}_{0,1}(p)$, then $\pi_p$ is generic, 
\begin{align}
\lambda_p  c(S)
&=p^3 a(S)+p^{6-2k}  a(p^2S)+\sum_{\substack{y\in\Z/p\Z\\ p\mid(\alpha+2\beta y)}}
p^{4-k}  a(S[\begin{bsmallmatrix}1&\\y&p\end{bsmallmatrix}])\nonumber \\
&\quad+\sum_{\substack{z\in\Z/p^2\Z\\ p^2 \mid (\alpha+ 2\beta z+\gamma z^2)}}
p^{2}  a(p^{-2}S[\begin{bsmallmatrix}1&\\z&p^2\end{bsmallmatrix}])\label{muzerosecondcaseeq}
\end{align}
for $S= \begin{bsmallmatrix} \alpha&\beta\\ \beta&\gamma \end{bsmallmatrix} \in B(N)^+$,
and $c(S) \neq 0$ for some $S \in B(N)^+$. 
The newform $F$ is an eigenvector for $T^{s}_{0,1}(p)$ if and only if
\begin{equation} 
\label{muzerofirstcaseeq}
\lambda_p  a(S) = (1+p^{-1}) c(S) 
=(1+p) p^{2-k}  a(pS)+\sum_{\substack{x\in\Z/p\Z\\ p \mid (\alpha +  2\beta x)}}
(1+p)  a(p^{-1}S[\begin{bsmallmatrix}1&\\x&p\end{bsmallmatrix}])
\end{equation}
for $S= \begin{bsmallmatrix} \alpha&\beta\\ \beta&\gamma \end{bsmallmatrix} \in B(N)^+$; in this case $\pi_p$ is non-generic.
\end{enumerate}
\end{corollary}

\begin{proof}
The equivalences \eqref{fourierheckethmeq1000} and \eqref{fourierheckethmeq1001} follow from \eqref{globalalgthmeq01},
\eqref{sigmafouriereq}, and \eqref{T10fouriereq}.

\ref{fourierheckethmsigmassquareitem}. Assume that $v_p(N) \geq 3$. Then $\sigma_p^2 F =0$ by \ref{globalalgthmitem15}
of Theorem \ref{globalalgthm}. Two applications of \eqref{sigmafouriereq} yield
\begin{align*}
(\sigma_p^2 F) (Z)
&= \sum_{S \in B(Np^{-2})^+} \sum_{y \in \Z/p\Z} \sum_{x \in \Z/p\Z} p^{-2k-2} 
 a\big((S[\begin{bsmallmatrix} 1& \\ x& p \end{bsmallmatrix} ])[\begin{bsmallmatrix} 1& \\ y& p \end{bsmallmatrix} ]\big)
\E^{2 \pi \I \mathrm{Tr}(SZ)}\\
&= \sum_{S \in B(Np^{-2})^+} \sum_{z\in \Z/p^2\Z} p^{-2k-2} 
 a( S[\begin{bsmallmatrix} 1& \\ z& p^2 \end{bsmallmatrix} ])
\E^{2 \pi \I \mathrm{Tr}(SZ)}.
\end{align*}
Since $\sigma_p^2 F =0$ we obtain \eqref{fourierheckethmsigmassquareitemeq}.

\ref{fourierheckethmmufinderitem}. By \ref{globalalgthmitem16} of Theorem \ref{globalalgthm} we have 
$\mu_p F = p^4 \tau_p^2 \sigma_p F - p^2 \eta_p \sigma_p F$. 
Each term of this equality is contained in $S_k(\mathrm{K}_s(Np))$. 
Let
\begin{equation}
\label{sigmabeq}
(\sigma_pF)(Z)=\sum_{S\in B(Np^{-1})^+}b(S) \E^{2\pi \I \mathrm{Tr}(SZ)}
\end{equation}
be the Fourier expansion of $\sigma_p F$. 
By \eqref{sigmafouriereq}, for $S\in B(Np^{-1})^+$,
\begin{equation}
\label{bformeq}
b(S)=\sum_{x\in\Z/p\Z} p^{-k-1}  a(S[\begin{bsmallmatrix}1&\\x&p\end{bsmallmatrix}]).
\end{equation}
Therefore, by \eqref{taunfouriereq},
\begin{align*}
p^4 \tau_p^2 \sigma_p F (Z)
&= \sum_{S \in B(Np)^+}p^4  b(S) \E^{2 \pi \I \mathrm{Tr}(SZ)}\\
&= \sum_{S \in B(Np)^+} \sum_{x \in \Z/ p\Z} p^{3-k} a(S[\begin{bsmallmatrix}1&\\x&p\end{bsmallmatrix}]) \E^{2 \pi \I \mathrm{Tr}(SZ)}.
\end{align*}
And by \eqref{etafouriereq}, 
\begin{align*}
(\eta_p \sigma_p F)(Z)
& = \sum_{\substack{S=\begin{bsmallmatrix} \alpha&\beta \\ \beta&\gamma \end{bsmallmatrix} \in B(Np)^+\\ p \mid 2\beta,\ p^2 \mid \gamma}}
p^k b( S[ \begin{bsmallmatrix} 1& \\ & p^{-1} \end{bsmallmatrix} ]) \E^{2 \pi \I \mathrm{Tr}(SZ)} \\
& = \sum_{\substack{S=\begin{bsmallmatrix} \alpha&\beta \\ \beta&\gamma \end{bsmallmatrix} \in B(Np)^+\\ p \mid 2\beta}}
p^k b( S[ \begin{bsmallmatrix} 1& \\ & p^{-1} \end{bsmallmatrix} ]) \E^{2 \pi \I \mathrm{Tr}(SZ)} \\
& = \sum_{\substack{S=\begin{bsmallmatrix} \alpha&\beta \\ \beta&\gamma \end{bsmallmatrix} \in B(Np)^+\\ p \mid 2\beta}}
\sum_{x\in\Z/p\Z} p^{-1}  
a\big( (S[ \begin{bsmallmatrix} 1& \\ & p^{-1} \end{bsmallmatrix} ])[\begin{bsmallmatrix}1&\\ x&p\vphantom{p^{-1}}\end{bsmallmatrix}] \big)
\E^{2 \pi \I \mathrm{Tr}(SZ)} \\
& = \sum_{\substack{S=\begin{bsmallmatrix} \alpha&\beta \\ \beta&\gamma \end{bsmallmatrix} \in B(Np)^+\\ p \mid 2\beta}}
\sum_{x\in\Z/p\Z} p^{-1}  
a( S[ \begin{bsmallmatrix} 1& \\ xp^{-1}  & 1 \end{bsmallmatrix} ] )
\E^{2 \pi \I \mathrm{Tr}(SZ)}.
\end{align*}
The assertion \eqref{mufouriereq} follows now from $\mu_p F = p^4 \tau_p^2 \sigma_p F - p^2 \eta_p \sigma_p F$.

\ref{fourierheckethmsigmaneqzeroitem}. Assume that $\mu_p \neq 0$. 
Again let the Fourier expansion of $\sigma_p F$ be as \eqref{sigmabeq} with $b$ as in \eqref{bformeq}.
By \eqref{T01fouriereq}, 
\begin{align*}
&(T^{s}_{0,1}(p)\sigma_pF)(Z)\\
&\qquad=\sum_{S=\begin{bsmallmatrix} \alpha&\beta\\ \beta&\gamma \end{bsmallmatrix}\in B(Np^{-1})^+}\Big(p^{3-k}b(pS)
+\!\!\sum_{\substack{y\in\Z/p\Z\\p\mid (\alpha+2\beta y+\gamma y^2)}}\!p 
b(p^{-1}S[\begin{bsmallmatrix}1&\\x&p\end{bsmallmatrix}])\Big) \E^{2\pi \I \mathrm{Tr}(SZ)}.
\end{align*}
By \eqref{stablelambdacase1eq}, $T^{s}_{0,1}(p)\sigma_pF=\lambda_p\sigma_pF$. Hence, for 
$S=\begin{bsmallmatrix} \alpha&\beta\\ \beta&\gamma \end{bsmallmatrix}\in B(Np^{-1})^+$,
\begin{align*}
\lambda_p b(S)=p^{3-k}b(pS)
+\sum_{\substack{y\in\Z/p\Z\\p\mid(\alpha+2\beta y+\gamma y^2)}}pb(p^{-1}S[\begin{bsmallmatrix}1&\\ y&p\end{bsmallmatrix}]).
\end{align*}
Substituting  \eqref{bformeq} we obtain, for $S=\begin{bsmallmatrix} \alpha&\beta\\ \beta&\gamma \end{bsmallmatrix}\in B(Np^{-1})^+$,	 
\begin{align*}
&\lambda_p \sum_{x\in\Z/p\Z}a(S[\begin{bsmallmatrix}1&\\x&p\end{bsmallmatrix}])\\
&\qquad=\sum_{x\in\Z/pZ}p^{3-k} a(pS[\begin{bsmallmatrix}1&\\x&p\end{bsmallmatrix}])
+\sum_{\substack{z\in\Z/p^2\Z\\p\mid(\alpha+2\beta z+\gamma z^2)}} p a(p^{-1}S[\begin{bsmallmatrix}1&\\z&p^2\end{bsmallmatrix}]).
\end{align*}
This is \eqref{munozerolambdaeq}.
Next we prove \eqref{epsilonfourierformulaeq}. By \eqref{T10fouriereq} and  \eqref{sigmafouriereq} we have  
$$
(T^{s}_{1,0}(p) \sigma_pF)(Z)= \sum_{S\in B(Np^{-1})^+} \sum_{y\in\Z/p^2\Z}
p^{2-2k} a(S[\begin{bsmallmatrix}1&\\y&p^2\end{bsmallmatrix}]) \E^{2\pi \I \mathrm{Tr}(SZ)}.
$$
By \eqref{sigmafouriereq} again, 
$$
(\sigma_pF)(Z)=\sum_{S\in B(Np^{-1})^+} \sum_{x\in\Z/p\Z}p^{-k-1} a(S[\begin{bsmallmatrix}1&\\x&p\end{bsmallmatrix}]) \E^{2\pi \I \mathrm{Tr}(SZ)}.
$$
Hence, \eqref{stablemucase1eq} implies that if $S\in B(Np^{-1})^+$, then
\begin{equation}
\label{epsilonfourierformulapreeq}
\sum_{y\in\Z/p^2\Z}a(S[\begin{bsmallmatrix}1&\\y&p^2\end{bsmallmatrix}])
=(\mu_p+p^2)p^{k-3}\sum_{x\in\Z/p\Z}a(S[\begin{bsmallmatrix}1&\\x&p\end{bsmallmatrix}]).
\end{equation}
Let $V \cong \otimes_v \pi_v$ be as in Theorem \ref{globalalgthm}. 
Applying  Corollary \ref{muvalueslemma} to $\pi_p$, we find that 
if $v_p(N)=2$, then $\mu_p=-p^2+\varepsilon_pp$, and if $v_p(N)>2$, then $\mu_p=-p^2$. This implies \eqref{epsilonfourierformulaeq}.
The remaining claim follows from \ref{globalalgthmitem2} of Theorem \ref{globalalgthm}. 

\ref{fourierheckethmsigmazeroitem}. Assume that $\mu_p = 0$.  The assertion \eqref{fourierheckethmeq1002}
follows immediately from \eqref{fourierheckethmeq1001}.
Let
$$
T^{s}_{0,1}(p)F=\sum_{S\in B(N)^+}c(S) \E^{2\pi \I \mathrm{Tr}(SZ)}
$$
be the Fourier expansion of $T^{s}_{0,1}(p)F$. 
By \eqref{T01fouriereq}, for $S=\begin{bsmallmatrix} \alpha & \beta \\ \beta & \gamma \end{bsmallmatrix} \in B(N)^+$,
\begin{align*}
c(S)=p^{3-k}a(pS)+\sum_{\substack{x\in\Z/p\Z\\p\mid(\alpha+2\beta x)}}pa(p^{-1}S[\begin{bsmallmatrix}1&\\ x&p\end{bsmallmatrix}]).
\end{align*}
For this, we note that $p \mid \gamma$ because $p \mid N_s$ and $N_s \mid \gamma$; hence $\alpha+2\beta x + \gamma x^2  \equiv \alpha+2\beta x
\Mod{p}$ for $x \in \Z$, which simplifies \eqref{T01fouriereq}. 
This proves \eqref{fourierheckethmeq10021}. 
Assume that $F$ is not an eigenvector for $T^s_{0,1}(p)$. Then $\pi_p$ is generic by
\ref{globalalgthmitem3} of Theorem \ref{globalalgthm}.  
Since $F$ is not an eigenvector for $T_{0,1}^s(p)$, $c(S)$ is non-zero for some $S \in B(N)^+$. 
Applying  \eqref{T01fouriereq} again, and also using \eqref{fourierheckethmeq1001}, we obtain
\begin{align*}
&(T_{0,1}^{s}(p)^2F)(Z)\\
&\qquad=\sum_{S=\begin{bsmallmatrix} \alpha & \beta \\ \beta & \gamma \end{bsmallmatrix}\in B(N)^+}
\big( p^{3-k}c(pS)+\sum_{\substack{y\in\Z/p\Z\\p\mid(\alpha+2\beta y)}}
p c(p^{-1}S[\begin{bsmallmatrix}1&\\y&p\end{bsmallmatrix}])\big) \E^{2\pi \I \mathrm{Tr}(SZ)}\\
&\qquad=\sum_{S=\begin{bsmallmatrix} \alpha & \beta \\ \beta & \gamma \end{bsmallmatrix}\in B(N)^+}
\Big(p^{3-k}\big(p^{3-k}a(p^2 S)+\sum_{\substack{x\in\Z/p\Z}}p a(S[\begin{bsmallmatrix}1&\\x&p\end{bsmallmatrix}])\big)\\
&\qquad\quad+\sum_{\substack{y\in\Z/p\Z\\p\mid(\alpha+2\beta y)}}p\big(p^{3-k}
a(S[\begin{bsmallmatrix}1&\\y&p\end{bsmallmatrix}])\\
&\qquad\quad+\sum_{\substack{x\in\Z/p\Z\\p\mid(\alpha'+2\beta x)\\\alpha'=p^{-1}(\alpha +2\beta y+\gamma y^2)}}
p a(p^{-2}S[\begin{bsmallmatrix}1&\\y+px&p^2\end{bsmallmatrix}])\big)\Big) \E^{2\pi \I \mathrm{Tr}(SZ)}\\
&\qquad=\sum_{S=\begin{bsmallmatrix} \alpha & \beta \\ \beta & \gamma \end{bsmallmatrix}\in B(N)^+}\big(  
p^{6-2k}a(p^2S)
+\sum_{\substack{y\in\Z/p\Z\\p\mid(\alpha+2\beta y)}}p^{4-k} a(S[\begin{bsmallmatrix}1&\\y&p\end{bsmallmatrix}])\\
&\qquad\quad+\sum_{\substack{y\in\Z/p\Z\\p\mid(\alpha+2\beta y)}} 
\sum_{\substack{x\in\Z/p\Z\\  p \mid (\alpha'+2\beta x) \\ \alpha'=p^{-1}(\alpha +2\beta y+\gamma y^2)}}
p^2 a(p^{-2}S[\begin{bsmallmatrix}1&\\y+px&p^2\end{bsmallmatrix}])\big) \E^{2\pi \I \mathrm{Tr}(SZ)}\\
&\qquad=\sum_{S=\begin{bsmallmatrix} \alpha & \beta \\ \beta & \gamma \end{bsmallmatrix} \in B(N)^+}\big( p^{6-2k}a(p^2S)
+\sum_{\substack{y\in\Z/p\Z\\ p\mid(\alpha+2\beta y)}} p^{4-k} a(S[\begin{bsmallmatrix}1&\\ y&p\end{bsmallmatrix}])\\
&\qquad\quad+\sum_{\substack{z\in\Z/p^2\Z\\ p^2 \mid (\alpha +2\beta z+ \gamma z^2)}}
p^2 a(p^{-2}S[\begin{bsmallmatrix}1&\\z&p^2\end{bsmallmatrix}])\big) \E^{2\pi \I \mathrm{Tr}(SZ)}.
\end{align*}
By \eqref{stablelambdacase3eq},  $F=p^{-3}\lambda_p T^{s}_{0,1}(p)F-p^{-3}T^{s}_{0,1}(p)^2 F$. 
Hence, for $S=\begin{bsmallmatrix} \alpha & \beta \\ \beta & \gamma \end{bsmallmatrix}\in B(N)^+$,
\begin{align*}
a(S)&=p^{-3}\lambda_p \big(p^{3-k}a(pS)+\sum_{\substack{x\in\Z/p\Z\\p\mid(\alpha+2\beta x)}}
p a(p^{-1}S[\begin{bsmallmatrix}1&\\ x&p\end{bsmallmatrix}])\big)\\
&\quad-p^{-3}\big( p^{6-2k}a(p^2S)+\sum_{\substack{y\in\Z/p\Z\\ p\mid(\alpha+2\beta y)}} 
p^{4-k} a(S[\begin{bsmallmatrix}1&\\ y&p\end{bsmallmatrix}])\\
&\quad+\sum_{\substack{z\in\Z/p^2\Z\\ p^2 \mid (\alpha +2\beta z+\gamma z^2)}}
p^2 a(p^{-2}S[\begin{bsmallmatrix}1&\\z&p^2\end{bsmallmatrix}])\big)\\
&=\lambda_pp^{-k}a(pS)+\lambda_p \sum_{\substack{x\in\Z/p\Z\\p\mid(\alpha+2\beta x)}}
p^{-2} a(p^{-1}S[\begin{bsmallmatrix}1&\\ x&p\end{bsmallmatrix}])\\
&\quad-p^{3-2k}a(p^2S)-\sum_{\substack{y\in\Z/p\Z\\ p\mid(\alpha+2\beta y)}}
p^{1-k} a(S[\begin{bsmallmatrix}1&\\y&p\end{bsmallmatrix}])\\
&\quad-\sum_{\substack{z\in\Z/p^2\Z\\ p^2 \mid (\alpha + 2\beta z+ \gamma z^2)}}
p^{-1} a(p^{-2}S[\begin{bsmallmatrix}1&\\z&p^2\end{bsmallmatrix}]).
\end{align*}
This proves \eqref{muzerosecondcaseeq}.
Finally, if $F$ is an eigenvector for  $T_{0,1}^{s}(p)$, then 
$\pi_p$ is non-generic by
\ref{globalalgthmitem3} of Theorem \ref{globalalgthm} and 
\eqref{muzerofirstcaseeq}
follows from \eqref{stablelambdacase2eq} and \eqref{fourierheckethmeq10021}. 
\end{proof}

We mention that the identities of Theorem~\ref{globalalgthm} can also be 
expressed in terms of the Fourier-Jacobi coefficients of the newform $F$, using the formulas computed in Chap.~\ref{operchap}.

\section{Another formulation}
\label{altformsec}

In this section we present another formulation of the consequences of Theorem~\ref{globalalgthm} for Fourier coefficients. This 
approach is based on the observation that the formulas from Corollary \ref{fourierheckethm}
involve right coset representatives for elements of the classical Hecke ring
of $\Gamma_0(N) \subset \SL(2,\Z)$, suggesting a more compact and abstract description.

To describe the setting, 
let $A(\Q)$
\label{AQdef} be the set of all positive semi-definite symmetric $2 \times 2$ matrices
with entries from $\Q$; if $X$ is a subset of $A(\Q)$, then 
let $X^+$ be the subset of elements of $X$ that are positive definite. 
We define an action of $\GL(2,\Q)$ on $A(\Q)$ by 
\begin{equation}
\label{dotactionexteq}
g\cdot S = g S \transpose{g}= S[\transpose{g}]
\end{equation}
for $g \in \GL(2,\Q)$ and $S \in A(\Q)$; this definition extends the action \eqref{gamma0acteq}. 
Let $\mathcal{F}$ \label{Fspace} be the set  of all functions $a:A(\Q) \to \C$. Then $\mathcal{F}$
is a $\C$-algebra under pointwise addition and multiplication of functions. 
Characteristic functions lying in $\mathcal{F}$ will be useful. 
Let $X$ be a subset of $A(\Q)$. We define $\setchar_X:A(\Q) \to \C$ to be the characteristic function
of $X$, so that
$$
\setchar_X(S) = 
\begin{cases}
1&\text{ if $S \in X$,}\\
0&\text{ if $S \notin X$.}
\end{cases}
$$
We define a right action of $\GL(2,\Q)$ on $\mathcal{F}$ by
$$
(a\big\vert g)(S)  = a(g \cdot S) = a(g S \transpose{g}) = a (S[\transpose{g}])
$$
for $a \in \mathcal{F}$, $g \in \GL(2,\Q)$ and $S \in A(\Q)$. 
Let $N$ be an integer such that $N>0$, and let $\Gamma_0(N)$ \label{gamma0Ndef}
be the subgroup of $\begin{bsmallmatrix}a&b \\ c&d \end{bsmallmatrix} \in \SL(2,\Z)$
such that $c \equiv 0 \Mod{N}$.
We define $\mathcal{F}^{\Gamma_0(N)}$ to be the subspace 
of $a \in \mathcal{F}$ such that $a\big \vert k = a$ for $k \in \Gamma_0(N)$. 
We note that  $\setchar_{B(N)}$ and $\setchar_{B(N)^+}$ are  in $\mathcal{F}^{\Gamma_0(N)}$. 
Fourier coefficients of Siegel modular forms are also contained in $\mathcal{F}^{\Gamma_0(N)}$. 
More precisely, let $k$ be an integer such that  $k>0$, and let 
$F \in M_k(\mathrm{K}_s(N))$. Let 
$$
F(Z) = \sum_{S \in B(N)} a(S) \E^{2 \pi \I \mathrm{Tr}(SZ)}
$$
be the Fourier expansion of $F$. To $F$ we  associate the function 
\begin{equation}
\label{hataFeq}
\hat a(F): A(\Q) \longrightarrow \C
\end{equation}
defined by 
$$
\hat a(F)(S) =
\begin{cases}
a(S)&\text{if $S \in B(N)$,}\\
0&\text{if $S \notin B(N)$.}
\end{cases}
$$
By \eqref{gamma0transruleeq} the function $\hat a(F)$ is contained in $\mathcal{F}^{\Gamma_0(N)}$.
Next, let $\Delta_0(N)$ be the subset of $\begin{bsmallmatrix} a&b \\ c&d \end{bsmallmatrix}
\in \Mat(2,\Z)$ such that $c \equiv 0 \Mod{N}$, $(a,N)=1$, and $ad-bc>0$. Then $\Delta_0(N)$
is a semigroup under multiplication of matrices. We let 
$\mathcal{R}(N) = \mathcal{R}(\Gamma_0(N),\Delta_0(N))$ be the  Hecke ring
with respect to $\Gamma_0(N)$ and $\Delta_0(N)$;  see, for example,
(4.5.4) of \cite{Mi}. The Hecke ring $\mathcal{R}(N)$ is
commutative (see Theorem~4.5.3
of \cite{Mi}), and we define an action of this ring
on $\mathcal{F}^{\Gamma_0(N)}$ via the formula 
$$
a \big \vert \Gamma_0(N) g \Gamma_0(N) = \sum_{i} a \big | g_i
$$
where $a \in\mathcal{F}^{\Gamma_0(N)}$, $g \in \Delta_0(N)$ and $\Gamma_0(N) g \Gamma_0(N) = \sqcup_i \Gamma_0(N) g_i$
is a disjoint decomposition. Following \cite{Mi}, we define
\begin{align}
T(l,m) &= \Gamma_0(N) \begin{bsmallmatrix}l&\\&m\end{bsmallmatrix} \Gamma_0(N),\label{Tlmdefeq}\\
T(n) & = \sum_{\det(g) = n} \Gamma_0(N) g \Gamma_0(N),\label{Tndefeq}
\end{align}
where $l \mid m$ and $(l,N)=1$, and 
in the second equality the summation is over all double cosets $\Gamma_0(N)g\Gamma_0(N)$
where  $g \in \Delta_0(N)$ with $\det(g) = n$. 
Let $p$ be a prime of $\Z$ dividing $N$. By (4.5.17) and Lemma 4.5.7 of \cite{Mi} we have
\begin{gather}
T(p) = T(1,p), \label{Tpformeq}\\
T(p)^2=T(p^2)=T(1,p^2). \label{Tp2formeq}
\end{gather}
By Lemma 4.5.6 of \cite{Mi} we have
\begin{align}
\Gamma_0(N) \begin{bsmallmatrix}1&\\&p\end{bsmallmatrix} \Gamma_0(N)
& = \bigsqcup_{a \in \Z/p\Z} \Gamma_0(N) \begin{bsmallmatrix}1&a\\&p\end{bsmallmatrix}, \label{Tpdecompeq}\\
\Gamma_0(N) \begin{bsmallmatrix}1&\\&p^2\end{bsmallmatrix} \Gamma_0(N)
& = \bigsqcup_{a \in \Z/p^2\Z} \Gamma_0(N) \begin{bsmallmatrix}1&a\\&p^2\end{bsmallmatrix}. \label{Tp2decompeq}
\end{align}
Finally, we will need three other operators on $\mathcal{F}$. 
Let $t \in \Z$ with $t>0$.
We define
$$
\Delta_t^+,\ \Delta_t^-,\ \nabla_t: \mathcal{F} \longrightarrow \mathcal{F} \label{deltatypedef}
$$
by 
\begin{align*}
(\Delta_t^+ a) (S) & = a(tS),\\
(\Delta_t^- a) (S) & = a(t^{-1}S),\\
(\nabla_t a) (S) & = (a\big \vert \begin{bsmallmatrix}1&\\&t^{-1}\end{bsmallmatrix})(S)=a(\begin{bsmallmatrix}1&\\&t^{-1}\end{bsmallmatrix}S\begin{bsmallmatrix}1&\\&t^{-1}\end{bsmallmatrix})
\end{align*}
for $a \in \mathcal{F}$ and $S \in A(\Q)$. The operators $\Delta_t^+$, $\Delta_t^-$, 
and $\nabla_t$ mutually commute, and the action of the Hecke algebra $\mathcal{R}(\Gamma_0(N),\Delta_0(N))$ commutes with $\Delta_t^+$ and $\Delta_t^-$. 
Also, if $a \in \mathcal{F}$ and $g \in \GL(2,\Q)$, then 
\begin{align}
(\mathrm{char}_{A(\Q)^+} \cdot a) \big | g &= \mathrm{char}_{A(\Q)^+} \cdot a \big | g,\label{AQpluscompeq1}\\
\Delta_t^+ (\mathrm{char}_{A(\Q)^+} \cdot a)& = \mathrm{char}_{A(\Q)^+} \cdot \Delta_t^+ a,\label{AQpluscompeq2}\\
\Delta_t^- (\mathrm{char}_{A(\Q)^+} \cdot a)& = \mathrm{char}_{A(\Q)^+} \cdot \Delta_t^- a,\label{AQpluscompeq3}\\
\nabla_t (\mathrm{char}_{A(\Q)^+} \cdot a) &= \mathrm{char}_{A(\Q)^+} \cdot \nabla_t a.\label{AQpluscompeq4}
\end{align}

The next lemma translates the Fourier coefficient calculations of Chap.~\ref{operchap}
into the language of the last paragraph. 

\begin{lemma}
\label{blockoperreformlemma}
Let $N$ and $k$ be integers such that $N>0$ and $k>0$, and let $p$
be a prime of $\Z$. Let $F \in M_k(\mathrm{K}_s(N))$, and set $\hat a = \hat a (F)$. 
Then:
\begin{align}
\hat a(\tau_{p} F) &= \setchar_{B(Np)}\cdot  \hat a, \label{taunfourieralteq}\\
\hat a(\theta_p F)
&=p^k\Delta_p^-  \hat a + p\cdot 
\left\{\begin{array}{ll}
\mathrm{char}_{B(Np^2)}&\text{if $v_p(N)=0$}\\
\mathrm{char}_{B(Np)}&\text{if $v_p(N)>0$}\\
\end{array}\right\} \cdot \Delta_p^+ \nabla_p  \hat a, \label{thetahateq} \\ 
\hat a(\eta_p F) &= p^k \nabla_p  \hat a,\label{etafourieralteq}\\
\hat a(\sigma_pF)
&=p^{-k-1}\setchar_{B(Np^{-1})}\cdot (\hat a \big | T(p)) \qquad\text{if $v_p(N) \geq 2$},\label{sigmafourieralteq}\\
\hat a(\sigma_p^2 F)
&=p^{-2k-2}\setchar_{B(Np^{-2})}\cdot ( \hat a \big | T(p)^2) \qquad\text{if $v_p(N) \geq 3$}, \label{sigmasquarefourieralteq}\\
\hat a (T_{0,1}^{s}(p) F) 
&=p^{3-k} \setchar_{B(N)} \cdot  \Delta_p^+ \hat a \nonumber \\
&\quad+ p\, \setchar_{B(N)}\cdot ( \Delta_p^-  \hat a \big | T(p)) \qquad\text{if $v_p(N) \geq 1$}, \label{T01fourieralteq}\\
\hat a(T_{1,0}^{s}(p)F) &=  p^{3-k} \setchar_{B(N)} \cdot  (\hat a \big | T(p))  \qquad\text{if $v_p(N) \geq 1$},\label{T10fourieralteq}\\
\hat a(T_{1,0}^s(p) \sigma_p F) & = p^{2-2k} \mathrm{char}_{B(Np^{-1})} \cdot (\hat a \big | T(p)^2) \qquad\text{if $v_p(N) \geq 2$.}
\label{T10sigmafourieralteq}
\end{align}
If $F \in M_k(\mathrm{K}(N))$ and $v_p(N) \geq 2$, then
\begin{align}
\hat a(T_{0,1}^s(p) \sigma_p F) & = p^{2-2k} \mathrm{char}_{B(Np^{-1})} \cdot ( \Delta_p^+ \hat a \big \vert T(p) )\nonumber \\
&\quad + p^{-k} \mathrm{char}_{B(Np^{-1})} \cdot ( \Delta_p^- \hat a \big \vert T(p)^2 ), \label{T01sigmafouriereq}\\
\hat a (T_{0,1}^s(p)^2 F) & = \mathrm{char}_{B(N)} \cdot \big( p^{6-2k} \Delta_{p^2}^+ \hat a + p^{4-k} \hat a \big \vert T(p) \nonumber \\
&\quad +p^{4-k}   \Delta_p^- (\mathrm{char}_{B(N)} \cdot \hat a)  \big \vert T(p) 
 + p^2  \Delta_{p^2}^- \hat a \big \vert T(p)^2\big). \label{T01T01fouriereq}
\end{align}
If $F$ is a cusp form, then the same equations with all appearing sets of the form $B(M)$ replaced
by $B(M)^+$, where $M$ an integer such that $M >0$, hold.
\end{lemma}
\begin{proof}
\eqref{taunfourieralteq}.  This follows immediately from \eqref{taunfouriereq}.

\eqref{thetahateq}. Let 
\begin{align*}
X_1&=\{\, \begin{bsmallmatrix} \alpha&\beta\\ \beta&\gamma \end{bsmallmatrix} \in B(N) \mid\  p\mid \alpha,\ p\mid 2\beta,\  N_s p \mid \gamma\,\},\\
X_2 &= \{\, \begin{bsmallmatrix} \alpha&\beta\\ \beta&\gamma \end{bsmallmatrix} \in B(N) \mid\  N_s p \mid \gamma\, \}.
\end{align*}
By   \eqref{thetafouriereq}, since $\hat a (\theta_p F)$ is supported in $B(Np)$ and $B(Np) \subset B(N)$, 
\begin{align*}
\hat a (\theta_p F) (S)
&=  p^k  \mathrm{char}_{B(Np)}(S)\,  \mathrm{char}_{X_1}(S)\hat a (p^{-1} S) \\
&\quad + p\,\mathrm{char}_{B(Np)}(S)\,  \setchar_{X_2} (S)  \hat a(p^{-1} S[ \begin{bsmallmatrix} p& \\ & 1 \end{bsmallmatrix}])
\end{align*}
for $S \in A(\Q)$. 
Since $X_1,X_2 \subset B(Np)$ we obtain 
\begin{align*}
 \hat a (\theta_p F) (S)
&=  p^k  \mathrm{char}_{X_1}(S)  \hat a (p^{-1} S) 
 + p\,\setchar_{X_2} (S)  \hat a(p^{-1} S[ \begin{bsmallmatrix} p& \\ & 1 \end{bsmallmatrix}])
\end{align*}
for $S \in A(\Q)$. Next, we note that 
$$
X_1 = p B(N) \quad \text{and} \quad
X_2 = 
\begin{cases} 
B(Np^2) &\text{if $v_p(N)=0$,}\\
B(Np) &\text{if $v_p(N) > 0$}.                            
\end{cases}
$$
Therefore, if $S \in A(\Q)$, then
\begin{align*}
p^k \mathrm{char}_{X_1}(S) \hat a (p^{-1} S)
& = p^k  \mathrm{char}_{pB(N)}(S) \hat a (p^{-1} S)\\
& = p^k  \mathrm{char}_{B(N)}(p^{-1}S) \hat a (p^{-1} S)\\
& = p^k   \hat a (p^{-1} S)\\
& = p^k   (\Delta_p^-\hat a) (S).
\end{align*}
For the third equality we used that the support of $\hat a$ lies in $B(N)$. Also, if $S \in A(\Q)$, then
\begin{align*}
p \,  \setchar_{X_2} (S)  \hat a(p^{-1} S[ \begin{bsmallmatrix} p& \\ & 1 \end{bsmallmatrix}])
&= p   \left\{\begin{array}{ll}
\mathrm{char}_{B(Np^2)}(S)&\text{if $v_p(N)=0$}\\
\mathrm{char}_{B(Np)}(S)&\text{if $v_p(N)>0$}\\
\end{array}\right\} (\Delta_p^+ \nabla_p \hat a)(S).
\end{align*}
This proves \eqref{thetahateq}. 

\eqref{etafourieralteq}. Let 
$$
X = \{\, \begin{bsmallmatrix} \alpha&\beta\\ \beta&\gamma \end{bsmallmatrix} \in B(Np^2) \mid \ p \mid 2\beta,\ p^2 \mid \gamma \, \}.
$$
Then an argument shows that $X=\begin{bsmallmatrix}1&\\&p\end{bsmallmatrix}B(N)\begin{bsmallmatrix}1&\\&p\end{bsmallmatrix}$. 
By \eqref{etafouriereq}, for $S \in A(\Q)$, 
\begin{align*}
\hat a (\eta_p F)(S)
& = p^k \mathrm{char}_X(S) (\nabla_p \hat a)(S) \\
& = p^k \mathrm{char}_{\begin{bsmallmatrix}1&\\&p\end{bsmallmatrix}B(N)\begin{bsmallmatrix}1&\\&p\end{bsmallmatrix}}(S) 
(\nabla_p \hat a)(S) \\
& = p^k (\nabla_p \mathrm{char}_{B(N)})(S)  (\nabla_p \hat a)(S)\\
& = p^k \nabla_p (\mathrm{char}_{B(N)}\cdot  \hat a)(S)\\
& = p^k  (\nabla_p \hat a)(S). 
\end{align*}
This is \eqref{etafourieralteq}. 

\eqref{sigmafourieralteq}. This follows
immediately from \eqref{sigmafouriereq}. 

\eqref{sigmasquarefourieralteq}. To prove this, assume that $v_p(N) \geq 3$. Let $S \in A(\Q)$.
If $S \notin B(Np^{-2})$, then both sides of \eqref{sigmasquarefourieralteq} evaluated at $S$ are zero. Assume
that $S=\begin{bsmallmatrix} \alpha&\beta \\ \beta &\gamma \end{bsmallmatrix} \in B(Np^{-2})$. By \eqref{sigmafourieralteq},
\begin{align*}
\hat a(\sigma_p^2 F)(S)
&= p^{-k-1} (\hat a(\sigma_p F)\big | T(p))(S)\\
&= p^{-k-1} \sum_{x \in \Z/ p\Z} \hat a(\sigma_p F)(S[\begin{bsmallmatrix} 1&\\ x&p \end{bsmallmatrix} ]).
\end{align*}
If $x \in \Z$ then $S[\begin{bsmallmatrix} 1&\\ x&p \end{bsmallmatrix} ] 
=\begin{bsmallmatrix} \alpha + 2\beta x + \gamma x^2 & \beta p + \gamma p x \\
 \beta p + \gamma p x  & \gamma p^2 \end{bsmallmatrix}
\in B(Np^{-1})$. 
Using \eqref{sigmafourieralteq} again, we have:
\begin{align*}
&p^{-k-1} \sum_{x \in \Z/ p\Z} \hat a(\sigma_p F)(S[\begin{bsmallmatrix} 1&\\ x&p \end{bsmallmatrix} ])\\
&\qquad= p^{-2k-2} \sum_{x \in \Z/ p\Z} \sum_{y \in \Z/ p\Z} \hat a( F)((S[\begin{bsmallmatrix} 1&\\ x&p \end{bsmallmatrix} ]) 
[\begin{bsmallmatrix} 1&\\ y&p \end{bsmallmatrix} ])\\
&\qquad=p^{-2k-2}\mathrm{char}_{B(Np^{-2})}(S) \cdot (\hat a \big | T(p)^2 )(S).
\end{align*}
This proves \eqref{sigmasquarefourieralteq}. 

\eqref{T01fourieralteq}. Assume that $v_p(N) \geq 1$, and 
let $S = \begin{bsmallmatrix} \alpha&\beta\\ \beta&\gamma \end{bsmallmatrix} \in A(\Q)$. 
If $S \notin B(N)$,
then both sides of \eqref{T01fourieralteq} evaluated at $S$ are trivially zero. Assume that $S \in B(N)$. Then
by \eqref{Tpdecompeq}, 
\begin{align*}
&p^{3-k} \setchar_{B(N)} (S)\,  (\Delta_p^+ \hat a) (S)
+ p\, \setchar_{B(N)}(S)\, (\Delta_p^-  \hat a \big | T(p))(S) \\
&\qquad = p^{3-k} \setchar_{B(N)} (S)\,  \hat a (pS)
+p\, \setchar_{B(N)}(S) \sum_{y \in \Z/p\Z}  \hat a(p^{-1} S[\begin{bsmallmatrix}1&\\y&p\end{bsmallmatrix}])\\
&\qquad=p^{3-k}   \hat a (pS)
+p \sum_{y \in \Z/p\Z}  \hat a(p^{-1} S[\begin{bsmallmatrix}1&\\y&p\end{bsmallmatrix}])\\
&\qquad=p^{3-k}    a (pS)
+p \sum_{y \in \Z/p\Z} \setchar_{B(N)} (p^{-1} S[\begin{bsmallmatrix}1&\\y&p\end{bsmallmatrix}])  \hat a(p^{-1} S[\begin{bsmallmatrix}1&\\y&p\end{bsmallmatrix}])\\
&\qquad=p^{3-k}   a (pS)
+p \sum_{y \in \Z/p\Z} \setchar_{B(N)} (
\begin{bsmallmatrix} 
(\alpha +2\beta y +\gamma y^2)p^{-1}& \beta +y \gamma \\
\beta +y \gamma & p\gamma             
\end{bsmallmatrix}) 
 \hat a(p^{-1} S[\begin{bsmallmatrix}1&\\y&p\end{bsmallmatrix}])\\
&\qquad=p^{3-k}   a (pS)
+p \sum_{\substack{y \in \Z/p\Z\\ p \mid (\alpha + 2\beta y+\gamma y^2)}} 
 a(p^{-1} S[\begin{bsmallmatrix}1&\\y&p\end{bsmallmatrix}])\\
&\qquad =  a(T_{0,1}^{s}(p)F)(S)\qquad \text{(by \eqref{T01fouriereq})}\\
&\qquad =   \hat a(T_{0,1}^{s}(p)F)(S).
\end{align*}
This proves \eqref{T01fourieralteq}. 

\eqref{T10fourieralteq}. This follows immediately from \eqref{T10fouriereq}.

\eqref{T10sigmafourieralteq}. The proof of this is similar to the proof of \eqref{sigmasquarefourieralteq}.

\eqref{T01sigmafouriereq}. 
Assume that $v_p(N) \geq 2$ and $F \in M_k(\mathrm{K}(N))$, and 
let $S \in A(\Q)$. Then using \eqref{T01fourieralteq} and \eqref{sigmafourieralteq}, we have:
\begin{align*}
&\hat a (T^{s}_{0,1}(p)\sigma_p F)(S)\\
&\qquad=p^{3-k} \setchar_{B(Np^{-1})} (S) \hat a(\sigma_p F) (pS)\\
&\quad\qquad + p\, \setchar_{B(Np^{-1})}(S) \sum_{x \in \Z/p\Z} \hat a(\sigma_pF)(p^{-1} S[\begin{bsmallmatrix}1&\\x&p\end{bsmallmatrix}])\\
&\qquad=p^{3-k} \setchar_{B(Np^{-1})} (S) \Big( p^{-k-1}\setchar_{B(Np^{-1})}(pS) 
\sum_{\substack{x\in\Z/p\Z}}\hat a(pS[\begin{bsmallmatrix}1&\\x&p\end{bsmallmatrix}])\Big) \\
&\quad\qquad + p\, \setchar_{B(Np^{-1})}(S) \sum_{x \in \Z/p\Z} \Big( 
p^{-k-1}\setchar_{B(Np^{-1})}(p^{-1} S[\begin{bsmallmatrix}1&\\x&p\end{bsmallmatrix}])\\
&\quad\qquad \times \sum_{\substack{y\in\Z/p\Z}}
\hat a(p^{-1} S[\begin{bsmallmatrix}1&\\x&p\end{bsmallmatrix}\begin{bsmallmatrix}1&\\y&p\end{bsmallmatrix}]) \Big) \\
&\qquad=p^{2-2k} \setchar_{B(Np^{-1})} (S)
\sum_{\substack{x\in\Z/p\Z}}\hat a(pS[\begin{bsmallmatrix}1&\\x&p\end{bsmallmatrix}]) \\
&\quad\qquad + p^{-k}\, \setchar_{B(Np^{-1})}(S) \sum_{x \in \Z/p\Z} \Big( 
\setchar_{B(Np^{-1})}(p^{-1} S[\begin{bsmallmatrix}1&\\x&p\end{bsmallmatrix}])\\
&\quad\qquad \times  \sum_{\substack{y\in\Z/p\Z}}
\hat a(p^{-1} S[\begin{bsmallmatrix}1&\\x&p\end{bsmallmatrix}\begin{bsmallmatrix}1&\\y&p\end{bsmallmatrix}]) \Big).
\end{align*}
We now observe the following. Let $x,y \in \Z$, and assume that
\begin{equation}
\label{doublecondeq}
\setchar_{B(Np^{-1})}(S) \neq 0
\quad \text{and} \quad \hat a ( p^{-1} S[ \begin{bsmallmatrix}1&\\x&p\end{bsmallmatrix} \begin{bsmallmatrix}1&\\y&p\end{bsmallmatrix}]) \neq 0.
\end{equation}
Let $S = \begin{bsmallmatrix} \alpha & \beta \\ \beta & \gamma \end{bsmallmatrix}$, and write $N = M p^{v_p(N)}$ with $M \in \Z$
and $p$ and $M$ relatively prime. Then $(Np^{-1})_s = M_s p^{v_p(N)-2}$. 
Since $\setchar_{B(Np^{-1})}(S) \neq 0$ we have $\alpha,2 \beta, \gamma \in \Z$ and $M_s p^{v_p(N)-2} \mid \gamma$. 
Since $F \in M_k(\mathrm{K}(N))$ the second statement in \eqref{doublecondeq} implies that
$$
p^{-1} S[ \begin{bsmallmatrix}1&\\x&p\end{bsmallmatrix} \begin{bsmallmatrix}1&\\y&p\end{bsmallmatrix}]
=
\begin{bsmallmatrix}
\gamma   y^2 p+2 \gamma  x y+ 2 \beta  y+( \alpha+2 \beta  x +\gamma  x^2)p^{-1} & \gamma   y p^2+\beta  p+\gamma   x p \\
 \gamma  y  p^2+\beta  p+\gamma   x p& \gamma  p^3 \\
\end{bsmallmatrix}
\in A(N).
$$
This implies that $p \mid (  \alpha+2 \beta  x +\gamma  x^2 )$ and $N \mid \gamma p^3$. 
Now 
\begin{equation}
\label{Sxmatrixeq}
p^{-1} S[ \begin{bsmallmatrix}1&\\x&p\end{bsmallmatrix}]= \begin{bsmallmatrix}
( \alpha+2 \beta  x +\gamma  x^2)p^{-1} & \beta +\gamma  x \\
 \beta +\gamma  x & \gamma  p \\
\end{bsmallmatrix}.
\end{equation}
Since  $N \mid \gamma p^3$, the integer $(Np^{-1})_s = M_s p^{v_p(N)-2}$ divides $ \gamma p$; 
therefore, since also $p \mid (  \alpha+2 \beta  x +\gamma  x^2 )$, 
the matrix in \eqref{Sxmatrixeq} is in 
$B(Np^{-1})$. 
Thus,  
$$
\setchar_{B(Np^{-1})} (p^{-1} S[ \begin{bsmallmatrix}1&\\x&p\end{bsmallmatrix}]) =1. 
$$
It follows that:
\begin{align*}
\hat a (T^{s}_{0,1}(p)\sigma_p F)(S)
&=p^{2-2k} \setchar_{B(Np^{-1})} (S)
\sum_{\substack{x\in\Z/p\Z}}\hat a(pS[\begin{bsmallmatrix}1&\\x&p\end{bsmallmatrix}]) \\
&\quad + p^{-k}\, \setchar_{B(Np^{-1})}(S) \sum_{x \in \Z/p\Z} 
\sum_{\substack{y\in\Z/p\Z}}\hat a(p^{-1} S[\begin{bsmallmatrix}1&\\x&p\end{bsmallmatrix}\begin{bsmallmatrix}1&\\y&p\end{bsmallmatrix}]) \\
&=p^{2-2k} \setchar_{B(Np^{-1})} (S) \cdot (\Delta_p  \hat a \big \vert T(p) )(S)\\
&\quad + p^{-k}\, \setchar_{B(Np^{-1})}(S) \cdot  (\Delta_p^-  \hat a \big \vert T(p)^2 )(S).
\end{align*}
This proves \eqref{T01sigmafouriereq}. 

\eqref{T01T01fouriereq}. Assume again that $v_p(N) \geq 2$ and $F \in M_k(\mathrm{K}(N))$, and 
let $S \in A(\Q)$.
Then:
\begin{align*}
&\hat a (T_{0,1}^{s}(p) T_{0,1}^{s}(p) F)(S)\\
&\qquad  = p^{3-k} \setchar_{B(N)} (S)\, \hat a (T_{0,1}^{s}(p) F)(pS) \\
&\quad\qquad  +p \,\setchar_{B(N)} (S) \sum_{x \in \Z/p\Z} \hat a(T_{0,1}^{s}(p) F) ( p^{-1} S[ \begin{bsmallmatrix}1&\\x&p\end{bsmallmatrix}])\\
&\qquad  = p^{3-k} \setchar_{B(N)} (S) \Big(p^{3-k} \setchar_{B(N)} (pS)\, \hat a (p^2 S)\\
&\quad\qquad + p \,\setchar_{B(N)} (pS) \sum_{x \in \Z/p\Z} \hat a (  S[ \begin{bsmallmatrix}1&\\x&p\end{bsmallmatrix}]) \Big) \\
&\quad\qquad  +p \,\setchar_{B(N)} (S) \sum_{x \in \Z/p\Z} 
\Big(  p^{3-k} \setchar_{B(N)} (p^{-1} S[ \begin{bsmallmatrix}1&\\x&p\end{bsmallmatrix}])\, \hat a (S[ \begin{bsmallmatrix}1&\\x &p\end{bsmallmatrix}]) \\
&\quad\qquad  +p \,\setchar_{B(N)} (p^{-1} S[\begin{bsmallmatrix}1&\\x&p\end{bsmallmatrix}]) \sum_{y \in \Z/p\Z} \hat 
a ( p^{-2} S[ \begin{bsmallmatrix}1&\\ x&p\end{bsmallmatrix} \begin{bsmallmatrix}1&\\y&p\end{bsmallmatrix}])  \Big) \\
&\qquad  = p^{3-k} \setchar_{B(N)} (S) \Big(p^{3-k} \setchar_{B(N)} (pS)\, \hat a (p^2 S)\\
&\quad\qquad+ p \,\setchar_{B(N)} (pS) \sum_{x \in \Z/p\Z} \hat a (  S[ \begin{bsmallmatrix}1&\\ x&p\end{bsmallmatrix}]) \Big) \\
&\quad\qquad +p^{4-k} \,\setchar_{B(N)} (S) \sum_{x \in \Z/p\Z} 
\setchar_{B(N)} (p^{-1} S[ \begin{bsmallmatrix}1&\\ x&p\end{bsmallmatrix}])\, \hat a (S[ \begin{bsmallmatrix}1&\\ x&p\end{bsmallmatrix}]) \\
&\quad\qquad +p^2 \,\setchar_{B(N)} (S)\\
&\quad  \qquad \times \sum_{x \in \Z/p\Z}  \sum_{y \in \Z/p\Z}\setchar_{B(N)} (p^{-1} S[ \begin{bsmallmatrix}1&\\x&p\end{bsmallmatrix}]) \hat 
a ( p^{-2} S[ \begin{bsmallmatrix}1&\\x&p\end{bsmallmatrix} \begin{bsmallmatrix}1&\\y&p\end{bsmallmatrix}]).
\end{align*}
Let $x,y \in \Z$, and assume that
\begin{equation}
\label{doublecondeq2}
\mathrm{char}_{B(N)} (S) \neq 0 \quad \text{and} \quad 
\hat a ( p^{-2} S[ \begin{bsmallmatrix}1&\\x&p\end{bsmallmatrix} \begin{bsmallmatrix}1&\\y&p\end{bsmallmatrix}], F) \neq 0.
\end{equation}
Let $S = \begin{bsmallmatrix} \alpha & \beta \\ \beta & \gamma \end{bsmallmatrix}$, and write $N = M p^{v_p(N)}$ with $M \in \Z$
and $p$ and $M$ relatively prime. Then $N_s = M_s p^{v_p(N)-1}$. 
Since $\setchar_{B(N)}(S) \neq 0$ we have $\alpha,2 \beta, \gamma \in \Z$ and $M_s p^{v_p(N)-1} \mid \gamma$. 
Since $F \in M_k(\mathrm{K}(N))$ the second equation in \eqref{doublecondeq2} implies that
$$
p^{-2} S[ \begin{bsmallmatrix}1&\\x&p\end{bsmallmatrix} \begin{bsmallmatrix}1&\\y&p\end{bsmallmatrix}]
=
\begin{bsmallmatrix}
(\alpha  +2\beta x  + \gamma x^2  )p^{-2}+(2\beta y+2  \gamma xy)p^{-1}+ \gamma y^2  &\gamma  x  +\gamma y   p+\beta  \\
 x \gamma + \gamma  y  p+\beta  & \gamma  p^2 \\
\end{bsmallmatrix}
\in A(N).
$$
This implies that $p \mid (\alpha +2 \beta x + \gamma x^2)$ and $N \mid \gamma p^2$. 
Now 
\begin{equation}
\label{Sxnextmatrixeq}
p^{-1} S[ \begin{bsmallmatrix}1&\\x&p\end{bsmallmatrix}]= \begin{bsmallmatrix}
( \alpha+2 \beta  x +\gamma x^2)p^{-1} & \beta +\gamma  x \\
 \beta +\gamma  x & \gamma  p \\
\end{bsmallmatrix}.
\end{equation}
Since  $N \mid \gamma p^2$, the integer $N_s = M_s p^{v_p(N)-1}$ divides $ \gamma p$; 
therefore, since  we also have $p \mid (  \alpha+2 \beta  x  +\gamma  x^2)$, 
the matrix in \eqref{Sxnextmatrixeq} is in 
$B(N)$. 
Thus,  
$$
\setchar_{B(N)} (p^{-1} S[ \begin{bsmallmatrix}1&\\x&p\end{bsmallmatrix}]) =1. 
$$
It follows that
\begin{align*}
&\hat a (T_{0,1}^{s}(p) T_{0,1}^{s}(p) F)(S)\\
&\qquad = p^{3-k} \setchar_{B(N)} (S) \big(p^{3-k} \setchar_{B(N)} (pS)\, \hat a (p^2 S)\\
&\quad\qquad+ p \,\setchar_{B(N)} (pS) \sum_{x \in \Z/p\Z} \hat a (  S[ \begin{bsmallmatrix}1&\\x&p\end{bsmallmatrix}]) \big) \\
&\quad\qquad +p^{4-k} \,\setchar_{B(N)} (S) \sum_{x \in \Z/p\Z} 
\setchar_{B(N)} (p^{-1} S[ \begin{bsmallmatrix}1&\\x&p\end{bsmallmatrix}])\, \hat a (S[ \begin{bsmallmatrix}1&\\x&p\end{bsmallmatrix}]) \\
&\quad\qquad +p^2 \,\setchar_{B(N)} (S) \sum_{x \in \Z/p\Z}  \sum_{y \in \Z/p\Z} \hat 
a ( p^{-2} S[ \begin{bsmallmatrix}1&\\x&p\end{bsmallmatrix} \begin{bsmallmatrix}1&\\y&p\end{bsmallmatrix}])\\
&\qquad = p^{3-k} \setchar_{B(N)} (S) \big(p^{3-k} \setchar_{B(N)} (pS)\, \hat a (p^2 S)\\
&\quad\qquad+ p \,\setchar_{B(N)} (pS) \sum_{x \in \Z/p\Z} \hat a (  S[ \begin{bsmallmatrix}1&\\x&p\end{bsmallmatrix}]) \big) \\
&\quad\qquad +p^{4-k} \,\setchar_{B(N)} (S) \sum_{x \in \Z/p\Z} 
\Delta^-_{p}(\setchar_{B(N)} )( S[ \begin{bsmallmatrix}1&\\x&p\end{bsmallmatrix}])\, \hat a (S[ \begin{bsmallmatrix}1&\\x&p\end{bsmallmatrix}]) \\
&\quad\qquad +p^2 \,\setchar_{B(N)} (S) \cdot ( \Delta_{p^2}^- \hat a \big\vert T(p)^2  ) (S)\\
&\qquad = p^{3-k} \setchar_{B(N)} (S) \big(p^{3-k} \, (\Delta_{p^2}^+\hat a)(S)
+ p \, (\hat a\big\vert T(p))(S) \big) \\
&\quad\qquad +p^{4-k} \,\setchar_{B(N)} (S) \cdot \big( (\Delta_p^-(\setchar_{B(N)}) \cdot \hat a)\big\vert T(p) \big)(S)  \\
&\quad\qquad +p^2 \,\setchar_{B(N)} (S) \cdot \big( \Delta_{p^2}^-  \hat a \big\vert T(p)^2 ) \big) (S)\\
&\qquad =  \setchar_{B(N)} (S) \big(  p^{6-2k} \, (\Delta_{p^2}^+\hat a)(S)
+ p^{4-k} \, (\hat a\big\vert T(p))(S)  \\
&\quad\qquad +p^{4-k} ( \Delta_p^-(\setchar_{B(N)} )\cdot \hat a)\big\vert T(p) )(S)  
+p^2 ( \Delta_{p^2}^-  \hat a \big\vert T(p^2) )  (S) \big).
\end{align*}
This is \eqref{T01T01fouriereq}.

Finally, assume that $F$ is a cusp form. Straightforward calculations using   
$\hat a = \mathrm{char}_{A(\Q)} \cdot \hat a$, \eqref{AQpluscompeq1}, \eqref{AQpluscompeq2}, \eqref{AQpluscompeq3},
and \eqref{AQpluscompeq4} then prove that the same equations hold with all sets of the form $B(M)$ replaced
by $B(M)$, where $M$ is an integer such that $M>0$. 
\end{proof}

We can now present  another formulation of the consequences of Theorem \ref{globalalgthm}
for Fourier coefficients.

\begin{theorem}
\label{abstractformulathm}
Let $N$ and $k$ be integers such that $N>0$ and $k>0$, and let $F \in S_k(\mathrm{K}(N))$. 
Assume that $F$ is a newform and is an eigenvector for the Hecke operators $T(1,1,q,q)$ and $T(1,q,q,q^2)$
for all but finitely many of the primes $q$ of $\Z$ such that $q \nmid N$;  by Theorem
\ref{Ralftheorem}, $F$ is an eigenvector for $T(1,1,q,q)$, $T(1,q,q,q^2)$ and $w_q$ for all primes $q$ of $\Z$.
Let $\lambda_q, \mu_q \in \C$ be such that 
\begin{align*}
T(1,1,q,q)F&=q^{k-3}\lambda_qF,\\
T(1,q,q,q^2)F&=q^{2(k-3)}\mu_qF.
\end{align*}
for all primes $q$ of $\Z$. 
Regard $F$ as an element of $S_k(\mathrm{K}_s(N))$ and 
define $\hat a = \hat a(F)$ as in \eqref{hataFeq}. 
Let $p$ be a prime of $\Z$ with $v_p(N)\geq2$. 
Then:
\begin{enumerate}
\item \label{abstractformulathmitem1} If $v_p(N) \geq 3$, then 
\begin{equation}
\setchar_{B(Np^{-2})}\cdot (\hat a \big | T(p^2) )=0. 
\end{equation}
\item \label{abstractformulathmitem2} We have
\begin{equation}
\label{abstractmueq}
\mu_p\, \hat a = p^{-k+3}\, \setchar_{B(Np)^+} \cdot  (\hat a \big\vert T(p))
-p\, \nabla_p\big(\setchar_{B(Np^{-1})^+} \cdot (\hat a  \big \vert T(p))  \big).
\end{equation}
\item \label{abstractformulathmitem3} Assume that $\mu_p \neq 0$. Then
\begin{align}
0&= \setchar_{B(Np^{-1})^+} \cdot \big( \lambda_p \, \hat a\big\vert T(p)
-p^{3-k} \, \Delta_p^+ \hat a \big \vert T(p) - p\, \Delta_p^-  \hat a \big \vert T(p^2) \big),\label{abstractmuneq0lambdaeq}\\
0&=\setchar_{B(Np^{-1})^+} \cdot \big( \hat a \big \vert T(p^2)
-(\mu_p+p^2)p^{k-3}  \,  \hat a \big \vert T(p)  \big).\label{abstractmuneq0eq}
\end{align}
\item \label{abstractformulathmitem4} Assume that $\mu_p=0$. Then
\begin{equation}
0=\setchar_{B(Np^{-1})^+} \cdot (\hat a \big \vert T(p)). \label{muzerosigmazeroeq}
\end{equation}
If $F$ is not an eigenform for $T_{0,1}^{s}(p)$, then
\begin{align}
0&=
\setchar_{B(N)^+} \cdot  \Big(p^2 \hat a+ p^{5-2k} \, \Delta_{p^2}^+ \hat a + p^{3-k} \, \hat a\big\vert T(p) +p \, \Delta_{p^2}^-  \hat a \big\vert T(p)^2\nonumber  \\
&\quad+p^{3-k} \, \big( \Delta_p^-(\setchar_{B(N)^+}) \cdot \hat a\big)\big\vert T(p)  
  - \lambda_p\big( p^{2-k} \, \Delta_p^+ \hat a +    \Delta_p^- \hat a \big\vert T(p)  \big) \Big). \label{abstractmuzeronoteigeneq}
\end{align}
If $F$ is an eigenform for $T_{0,1}^{s}(p)$, then
\begin{equation}
\label{abstractmuzeroeigeneq}
0=\setchar_{B(N)^+} \cdot \big( (1+p^{-1})^{-1} \lambda_p\, \hat a
-p^{3-k}  \ \Delta_p^+ \hat a 
-p \,  \Delta_p^-  \hat a \big \vert T(p) \big).
\end{equation}
\end{enumerate}
\end{theorem}				
\begin{proof}
Throughout the proof we will use the formulas from Lemma \ref{blockoperreformlemma}, taking 
into account the last statement of this lemma. 

\ref{abstractformulathmitem1}. This follows from \eqref{sigmasquarefourieralteq}. 

\ref{abstractformulathmitem2}. By \eqref{globalnewformreleq} we have 
$$
\mu_p F = p^4 \tau_p^2 \sigma_p F - p^2 \eta_p \sigma_p F.
$$
Therefore,
\begin{align*}
\mu_p \hat a(F)
& = p^4 \hat a( \tau_p^2 \sigma_p F) - p^2 \hat a (\eta_p \sigma_p F)\\
& = p^4 \mathrm{char}_{B(N)^+} \cdot \hat a( \tau_p \sigma_p F) - p^{k+2} \nabla_p \hat a (\sigma_p F)\qquad\!\text{(by \eqref{taunfourieralteq}
and \eqref{etafourieralteq})}\\
& = p^4 \mathrm{char}_{B(N)^+} \cdot \mathrm{char}_{B(Np)^+}\cdot \hat a(\sigma_p F) \\
&\quad- p \nabla_p (\mathrm{char}_{B(Np^{-1})^+} \cdot\hat a \big | T(p))\qquad\text{(by \eqref{taunfourieralteq}
and \eqref{sigmafourieralteq})}\\
& = p^{-k+3} \mathrm{char}_{B(N)^+} \cdot \mathrm{char}_{B(Np)^+}\cdot \mathrm{char}_{B(Np^{-1})^+}\cdot \hat a \big | T(p) \\
&\quad- p \nabla_p (\mathrm{char}_{B(Np^{-1})^+} \cdot\hat a \big | T(p))\qquad\text{(by  \eqref{sigmafourieralteq})}\\
\mu_p \hat a & = p^{-k+3}  \mathrm{char}_{B(Np)^+}\cdot \hat a \big | T(p) - p \nabla_p (\mathrm{char}_{B(Np^{-1})^+} \cdot\hat a \big | T(p)).
\end{align*}
This proves \eqref{abstractmueq}.

\ref{abstractformulathmitem3}.  Assume that $\mu_p \neq 0$. 
By \eqref{stablelambdacase1eq} we have
$T^{s}_{0,1}(p)(\sigma_p F)=\lambda_p(\sigma_pF)$. 
Hence, by \eqref{sigmafourieralteq} and \eqref{T01sigmafouriereq},
\begin{align*}
\lambda_p \hat a (\sigma_p F)
& = \hat a ( T_{0,1}^s(p) \sigma_p F) \\
\lambda_p p^{-k-1} \mathrm{char}_{B(Np^{-1})^+} \cdot ( \hat a \big \vert T(p) ) 
& = p^{2-2k} \mathrm{char}_{B(Np^{-1})^+} \cdot ( \Delta_p^+ \hat a \big \vert T(p))\\
&\quad + p^{-k} \mathrm{char}_{B(Np^{-1})^+} \cdot (\Delta_p^- \hat a \big \vert T(p)^2) \\
\lambda_p  \mathrm{char}_{B(Np^{-1})^+} \cdot ( \hat a \big \vert T(p) ) 
& = p^{3-k} \mathrm{char}_{B(Np^{-1})^+} \cdot ( \Delta_p^+ \hat a \big \vert T(p))\\
&\quad + p\, \mathrm{char}_{B(Np^{-1})^+} \cdot (\Delta_p^- \hat a \big \vert T(p)^2).
\end{align*}
This is \eqref{abstractmuneq0lambdaeq}.
Next, by \eqref{stablemucase1eq} we have
$T^{s}_{1,0}(p)(\sigma_p F)=(\mu_p+p^2)(\sigma_pF)$. Therefore, by \eqref{sigmafourieralteq}
and \eqref{T10sigmafourieralteq}, 
\begin{align*}
\hat a ( T_{1,0}^s(p) \sigma_p F)
& = (\mu_p+p^2) \hat a (\sigma_p F) \\
p^{2-2k}\mathrm{char}_{B(Np^{-1})^+} \cdot \hat a \big\vert T(p)^2
& = (\mu_p+p^2) p^{-k-1} \mathrm{char}_{B(Np^{-1})^+} \cdot \hat a \big\vert T(p).
\end{align*}
This is \eqref{abstractmuneq0eq}. 

\ref{abstractformulathmitem4}. Assume that $\mu_p=0$. The assertion \eqref{muzerosigmazeroeq} follows directly from \eqref{globalalgthmeq01} and  \eqref{sigmafourieralteq}.
Assume that $F$ is not an eigenform for $T_{0,1}^s (p)$. 
By \eqref{stablelambdacase3eq} we have $T^{s}_{0,1}(p)^2F=-p^3F+\lambda_p T^{s}_{0,1}(p)F$. Hence, 
\begin{equation}
\label{Ts01twiceeq}
\hat a (T_{0,1}^s(p)^2 F) = - p^3 \hat a + \lambda_p \hat a (T^{s}_{0,1}(p)F).
\end{equation}
By \eqref{T01T01fouriereq}
\begin{align*}
\hat a (T_{0,1}^s(p)^2 F) & = \mathrm{char}_{B(N)^+} \cdot \big( p^{6-2k} \Delta_{p^2}^+ \hat a + p^{4-k} \hat a \big \vert T(p) \nonumber \\
&\quad +p^{4-k}   \Delta_p^- (\mathrm{char}_{B(N)^+} \cdot \hat a)  \big \vert T(p) 
 + p^2  \Delta_{p^2}^- \hat a \big \vert T(p)^2\big)\\
& = p\, \mathrm{char}_{B(N)^+} \cdot  \big( p^{5-2k} \Delta_{p^2}^+ \hat a + p^{3-k} \hat a \big \vert T(p) \nonumber \\
&\quad +p^{3-k}   \Delta_p^- (\mathrm{char}_{B(N)^+} \cdot \hat a)  \big \vert T(p) 
 + p  \Delta_{p^2}^- \hat a \big \vert T(p)^2\big)\\
\end{align*}
and by \eqref{T01fourieralteq}
\begin{align*}
&-p^3\hat a(F) + \lambda_p \hat a (T_{0,1}^{s}(p) F) \\
&\qquad =-p^3\hat a 
+  \lambda_p \big(p^{3-k}  \setchar_{B(N)^+} \cdot  \Delta_p^+ \hat a 
+ p\,   \setchar_{B(N)^+}\cdot  (\Delta_p^-  \hat a \big | T(p)) \big)\\
&\qquad = p\, \mathrm{char}_{B(N)^+} \cdot  \big( -p^2\hat a+  \lambda_p (  p^{2-k} \Delta_p^+ \hat a 
+   \Delta_p^-  \hat a \big | T(p) ) \big).
\end{align*}
Substituting into \eqref{Ts01twiceeq} now proves \eqref{abstractmuzeronoteigeneq}; note that we have used that 
$\hat a(F) =  \setchar_{B(N)^+}\cdot \hat a (F)$. 

Finally, assume  that  $F$ is an eigenform for $T_{0,1}^{s}(p)$. By \eqref{stablelambdacase2eq} we have
$T^{s}_{0,1}(p)F=(1+p^{-1})^{-1}\lambda_pF$. Therefore, by \eqref{T01fourieralteq},
\begin{align*}
(1+p^{-1})^{-1}\lambda_p \hat a 
& = \hat a ( T^{s}_{0,1}(p)F ) \\
(1+p^{-1})^{-1}\lambda_p  \setchar_{B(N)^+} \cdot  \hat a 
&=p^{3-k} \setchar_{B(N)^+} \cdot  \Delta_p^+ \hat a 
+ p\, \setchar_{B(N)^+}\cdot ( \Delta_p^-  \hat a \big | T(p)).
\end{align*}
This is \eqref{abstractmuzeroeigeneq}. 
\end{proof}

\section{Examples}
\label{examplessec}

The work \cite{PSY} and the website \cite{PSYW} present examples of 126 distinct newforms in $S_k(\mathrm{K}(16))$ for $k=6,7,8,9,10,11,12,13$ and $14$ that are eigenvectors for $T(1,1,q,q)$, $T(1,q,q,q^2)$, and $w_q$ for all primes $q$ of $\Z$. 
These newforms are non-lifts and form 58 Galois orbits.
Altogether, about 67,500 Fourier coefficients appear in \cite{PSYW}. These works also determine the Hecke eigenvalues 
$\lambda_2$ and $\mu_2$ and provide information about $\pi_2$; in particular, $\pi_2$ is always generic.
In this section we will describe how we verified that the seven equations \eqref{fourierheckethmsigmassquareitemeq}--\eqref{muzerofirstcaseeq} from Corollary \ref{fourierheckethm} hold
for this data. This verification is a welcome extra check of Corollary~\ref{fourierheckethm} (in addition to the proof!).

\subsection*{Indices and Fourier coefficients} 
We begin by reviewing some background
concerning Fourier coefficients and their
index set. 
Let $N,k \in \Z$ with $N>0$ and $k>0$, and let $F \in S_k(\mathrm{K}(N))$. Let 
$$
F(Z) = \sum\limits_{S \in A(N)^+} a(S) \E^{2 \pi \I \mathrm{Tr}(SZ)}
$$
be the Fourier expansion of $F$. As defined in \eqref{gamma0acteq}, the group $\Gamma_0(N)_{\pm}$ acts on $A(N)^+$.
Thus, the set $A(N)^+$ is partitioned into $\Gamma_0(N)_{\pm}$ orbits. 
By \eqref{gamma0transruleeq},  if $S \in A(N)^+$ and $a(S)$ is known, then $a(T)$ is known for all $T$
in the orbit $\Gamma_0(N)_{\pm} \cdot S$. 
Also, if $d \in \Z$ with $d>0$,
then the number of orbits $\Gamma_0(N)_{\pm} \cdot S$ for $S \in A(N)^+$ such that $4 \det(S) =d$ is finite.
This is a consequence of a result of Gauss; see also Lemma \ref{Ydequivlemma} below. 
We may thus effectively list all the Fourier coefficients of $F$ by first making a list of orbit representatives $S$ for $\Gamma_0(N)_{\pm} \backslash A(N)^+$ with $4\det(S)$ in increasing order, and then specifying for each such orbit representative $S$ the Fourier coefficient $a(S)$.
Below, we will describe a method for obtaining such a list of orbit representatives for $\Gamma_0(N)_{\pm} \backslash A(N)^+$; as we will see, this technique can also be used to generate other lists of indices needed for the verification. 
To obtain  orbit representatives we need three lemmas.
\begin{lemma}
\label{Ydsetlemma}
Let $d\in \Z$ be such that $d>0$, and let $Y(d)$ be the set of all $(a,b,c) \in \Z^3$
such that 
\begin{enumerate}
\item\label{Ydsetlemmaitem1} $d = 4ac-b^2$;
\item\label{Ydsetlemmaitem2} $0<a \leq c$;
\item\label{Ydsetlemmaitem3} $|b| \leq a$. 
\end{enumerate}
If $(a,b,c) \in Y(d)$, then 
$$
1 \leq a \leq \sqrt{d/3},  \quad -a \leq b \leq a \leq c, \quad \text{and} \quad c = (d+b^2)/(4a).
$$
In particular, the set $Y(d)$ is finite. 
\end{lemma}
\begin{proof}
Let $(a,b,c) \in Y(d)$. Then
\begin{align*}
d & = 4ac-b^2 \\
& \geq 4a^2 -a^2 \qquad \text{(by \ref{Ydsetlemmaitem2} and \ref{Ydsetlemmaitem3})}\\
& = 3a^2.
\end{align*}
This implies that $\sqrt{d/3} \geq a$. It is immediate that $-a \leq b \leq a \leq c$ and that $c=(d+b^2)/(4a)$. 
Since the function
$$
Y(d) \longrightarrow \{ 1,\dots,\lfloor \sqrt{d/3} \rfloor \} \times
\{ -\lfloor \sqrt{d/3} \rfloor ,\dots,\lfloor \sqrt{d/3} \rfloor \}
$$
defined by $(a,b,c) \mapsto (a,b)$ for $(a,b,c) \in Y(d)$ is well-defined and injective, the set $Y(d)$ is finite. 
\end{proof}

Let $N$ and $d$ be integers such that $N>0$ and $d>0$. We define
\begin{equation}
\label{ANddefeq}
A(N)^+(d) = \{ S \in A(N)^+\mid 4 \det(S) = d\}. 
\end{equation}
We are interested in obtaining orbit representatives for $\Gamma_0(N)_{\pm} \backslash A(N)^+(d)$. We 
note that the set $A(1)^+(d)$ contains $A(N)^+(d)$. 

\begin{lemma}
\label{Ydequivlemma}
Let $d\in \Z$ be such that $d>0$. If $S \in A(1)^+(d)$, then 
there exists an element $k \in \SL(2,\Z)$ such that 
$$
k \cdot S 
\in
\{ \begin{bsmallmatrix}   a & b/2 \\ b/2 & c \end{bsmallmatrix} \in A(1)^+(d)\mid(a,b,c) \in Y(d) \}. 
$$
\end{lemma}
\begin{proof}
We follow the argument in \cite{Z}, p.~59.
Let $S \in A(1)^+(d)$, and write
$$
S = \begin{bsmallmatrix} a & b/2 \\ b/2 & c \end{bsmallmatrix}
$$
with $a,b,c \in \Z$. Since $S$ is positive-definite, we have $ax^2 + bxy +c y^2 >0$ for $(x,y) \in \Z^2$ with $(x,y) \neq 0$. 
Let $a'$ be the smallest integer in the set $\{ax^2 + bxy +c y^2\mid (x,y) \in \Z^2, (x,y) \neq 0 \}$. The integer $a'$ is
positive; let $(x,y) \in \Z^2$ be such that $a'=ax^2 + bxy +c y^2$. The definition of $a'$ implies that $x$ and $y$ are relatively
prime. Let $w,z \in \Z$ be such that $xz-yw=1$. Then
$$
k = \begin{bsmallmatrix}    x & y \\ w & z \end{bsmallmatrix} \in \SL(2,\Z).
$$
Moreover,
$$
k \cdot S = 
\begin{bsmallmatrix}  ax^2 + bxy +c y^2 & b''/2 \\ b''/2 & aw^2 + bwz +c z^2 \end{bsmallmatrix}
= \begin{bsmallmatrix}  a' & b''/2 \\ b''/2 & c'' \end{bsmallmatrix}
$$
where $b'' \in \Z$ and $c''\in \Z$. Next, there exists an integer
$n$ such that 
$$
2n-1 \leq b''/ a' \leq 2n+1.
$$
This implies that
$$
-a' \leq b'' -2 a'n  \leq a'.
$$
Let 
$$
h =\begin{bsmallmatrix} 1 &  \\ -n & 1 \end{bsmallmatrix} \quad \text{and} \quad b'=b''-2a'n.
$$
Then $h \in \SL(2,\Z)$ and 
$$
 h \cdot (k \cdot S)  = \begin{bsmallmatrix}  a' & b'/2 \\ b'/2 & c' \end{bsmallmatrix}
$$
for some $c' \in \Z$. Since $-a' \leq b'' -2 a'n  \leq a'$ we have $|b'| \leq a'$. Concerning $c'$, we
note that $h \cdot (k \cdot S) = (hk) \cdot S$; since
$$
hk = \begin{bsmallmatrix} x & y \\ w - n x & -n y + z \end{bsmallmatrix}
$$
we obtain
$$
c' = a (w-nx)^2+b(w-nx)(-ny+z)+c(-ny+z)^2.
$$
By the definition of $a'$, we must have $a' \leq c'$. 
\end{proof}

\begin{lemma}
\label{cosetrepgammalemma}
Let $N$ and $d$ be integers such that $N>0$ and $d>0$, and let $\Gamma$ be 
a subgroup of $\Gamma_0(N)$. Let $Y(d) = \{ (a_1,b_1,c_1),\dots,(a_t,b_t,c_t)\}$,
and let $\SL(2,\Z) = \sqcup_{i \in I} \Gamma g_i$ be a disjoint decomposition.
Define $X(\Gamma,N,d)$ to be the set of all the $g_i \cdot \begin{bsmallmatrix} a_j & b_j/2 \\ b_j/2 & c_j \end{bsmallmatrix}$
that are contained in $A(N)^+(d)$, where $i \in I$ and $j \in \{1,\dots,t\}$. Then 
$
A(N)^+(d) = \cup_{S \in X(\Gamma,N,d)} \Gamma \cdot S.
$
\end{lemma}
\begin{proof}
Clearly, $\cup_{S \in X(\Gamma,N,d)} \Gamma \cdot S \subset A(N)^+(d)$. Let $T \in A(N)^+(d)$. By Lemma~\ref{Ydequivlemma}
there exists $g \in \SL(2,\Z)$ such that $g \cdot T =\begin{bsmallmatrix} a_j & b_j/2 \\ b_j/2 & c_j \end{bsmallmatrix}$
for some $j \in \{1,\dots,t\}$. Let $i \in I$ and $h \in \Gamma$ be such that $g^{-1} = h g_i$. Then
$T= h \cdot (g_i \cdot \begin{bsmallmatrix} a_j & b_j/2 \\ b_j/2 & c_j \end{bsmallmatrix})$, so that $T \in \cup_{S \in X(\Gamma,N,d)} \Gamma \cdot S$. 
\end{proof}

We can now state the steps for determining a set of orbit representatives for $\Gamma_0(N)_{\pm} \backslash A(N)^+(d)$
when $N, d\in \Z$ with $N>0$ and $d>0$:
\begin{enumerate}
\item Determine the finite set $Y(d)$.
\item Find a disjoint decomposition $\SL(2,\Z)=\sqcup_{i \in I} \Gamma_0(N) g_i$.
\item Determine the set $X(\Gamma_0(N),N,d)$.
\item By Lemma \ref{cosetrepgammalemma} the set $X(\Gamma_0(N),N,d)$ contains a set of orbit representatives for 
$\Gamma_0(N) \backslash A(N)^+(d)$, and hence a set of orbit representatives for $\Gamma_0(N)_{\pm} \backslash A(N)^+(d)$.
Refine the set $X(\Gamma_0(N),N,d)$ to obtain a set of orbit representatives for $\Gamma_0(N)_{\pm} \backslash A(N)^+(d)$.
\end{enumerate}
We will not describe an explicit implementation of this algorithm, but do
 mention two points. First, if $N=p^t$ is a power of a prime $p$, e.g., $16$, then there is 
a convenient disjoint decomposition
$$
\SL(2,\Z) = \bigsqcup_{x \in \Z/ p^t \Z} \Gamma_0(p^t) \begin{bsmallmatrix} 1& \\ x&1 \end{bsmallmatrix} \sqcup 
\bigsqcup_{y \in p\Z / p^t \Z} \Gamma_0(p^t) \begin{bsmallmatrix} 1& \\ y&1 \end{bsmallmatrix} \begin{bsmallmatrix} &1 \\ -1 & \end{bsmallmatrix}.
$$
Second, for the implementation of (4), given $S_1,S_2 \in A(N)^+(d)$, one needs a method for determining whether or not there exists $k \in \Gamma_0(N)_{\pm}$ such that $k \cdot S_1 = S_2$. One such method is described in Lemma 2.6 of \cite{Zac}. Alternatively, one may  use the meta-algorithm function
{\tt Solve} from Mathematica. For this, writing $S_1=\begin{bsmallmatrix} x_1&y_1\\y_1&z_1\end{bsmallmatrix}$ and 
$S_2=\begin{bsmallmatrix} x_2&y_2\\y_2&z_2\end{bsmallmatrix}$, it is convenient to note that 
there exists $k \in \Gamma_0(N)_{\pm}$ such that $k \cdot S_1 = S_2$ if and only if there exist $a,b,c,d \in \Z$ such that 
$$
\begin{bsmallmatrix} a&b\\ Nc&d \end{bsmallmatrix} \begin{bsmallmatrix} x_1&y_1\\y_1\vphantom{N}&z_1\end{bsmallmatrix}
=\begin{bsmallmatrix} x_2&y_2\vphantom{N}\\y_2&z_2\end{bsmallmatrix} \begin{bsmallmatrix} d&-Nc\\ -b&a \end{bsmallmatrix}\quad\text{and}\quad ad-bNc=1
$$
or there exist $a,b,c,d \in \Z$ such that 
$$
\begin{bsmallmatrix} a&b\\ Nc&d \end{bsmallmatrix} \begin{bsmallmatrix} x_1&y_1\\y_1\vphantom{N}&z_1\end{bsmallmatrix}
=\begin{bsmallmatrix} x_2&y_2\vphantom{N}\\y_2&z_2\end{bsmallmatrix} \begin{bsmallmatrix} -d&Nc\\ b&-a \end{bsmallmatrix}\quad\text{and}\quad ad-bNc=-1;
$$
moreover, the number of $k \in \Gamma_0(N)_{\pm}$ such that $k \cdot S_1 = S_2$ is finite.

\subsection*{The data from \cite{PSYW}}
Next, we review the data from \cite{PSYW}.
Let $\mathcal{O}$ be one of the 58 Galois orbits from \cite{PSYW}. Let $p_{\mathcal{O}}(X)\in \Z[X]$ be the 
irreducible monic polynomial corresponding to $\mathcal{O}$; also, let $k$ be the weight corresponding to 
 $\mathcal{O}$. Fix a root $\alpha$ of $p_{\mathcal{O}}(X)$. Then corresponding to $\alpha$
there exists a newform $F_\alpha \in S_k(\mathrm{K}(16))$ that is an eigenform for the Hecke operators
$T(1,1,p,p)$ and $T(1,p,p,p^2)$ for all primes $p$ of $\Z$. Let
$$
F_\alpha(Z) = \sum_{S \in A(16)^+} a(S)\E^{2 \pi \I \mathrm{Tr}(SZ)}
$$
be the Fourier expansion of $F_\alpha$. 
Some of the Fourier coefficients of $F_\alpha$
are given in \cite{PSYW}
via a three-column table. A row of the table has the form
\begin{equation}
\label{tablelineeq}
\def\arraystretch{1.5}
\begin{array}{|c|c|c|}
\hline
(\alpha,\beta,\gamma) & d & z\\
\hline
\end{array}
\end{equation}
Here, $\alpha,2\beta$ and $\gamma$ are integers, and the triple $(\alpha,\beta,\gamma) $ is an abbreviation for  the element $S$ of $A(16)^+$
given by 
\begin{equation}
\label{Sfromtableeq}
S  =\begin{bsmallmatrix} \alpha & \beta \\ \beta & \gamma \end{bsmallmatrix}.
\end{equation} 
Also, 
$$
d = 4\det(S) = 4\det(\begin{bsmallmatrix} \alpha & \beta \\ \beta & \gamma \end{bsmallmatrix} )= 4(\alpha\gamma -\beta^2).
$$
The integer $d$ is necessarily positive and satisfies $d \equiv \text{$0$ or $3$} \Mod{4}$.
Finally, $z$ is an element of $\Q(\alpha)$ such that 
$$
z = a(S) = a(\begin{bsmallmatrix} \alpha & \beta \\ \beta & \gamma \end{bsmallmatrix} ), 
$$
The tables from \cite{PSYW} have several further properties. First, the rows are ordered
by the second entry, $d=4\det(S)$, and $d$ is such that $d<500$;
most tables have approximately $1200$ entries. Second, if $S_1 \in A(16)^+$ and $S_2 \in A(16)^+$ 
correspond to the first
entries in two distinct rows, then $S_1$ and $S_2$ are not in the same $\Gamma_0(16)_{\pm}$ orbit
(the action is defined in \eqref{gamma0acteq}).  Thus,
each row corresponds to a unique orbit of $\Gamma_0(16)_{\pm}$ acting on $A(16)^+$. 
Note that if $S \in A(16)^+$ corresponds to the first entry of a table, so that $a(S)$ is known, then 
$a(S')$ is also known for all $S'$ in the $\Gamma_0(16)_{\pm}$ orbit of $S$ by 
\eqref{gamma0transruleeq}. 
Finally, not every orbit of $\Gamma_0(16)_{\pm}$ acting on $A(16)^+(d)$ with $d<500$ is
represented in every table. (In fact, our summary of the tables from \cite{PSYW} is not completely
accurate: in the tables from \cite{PSYW} the first entry of each row is given as $(2\alpha,2\beta,2\gamma)$ rather
than $(\alpha,\beta,\gamma)$.) As an illustration, Table \ref{Ftable} shows the
first few Fourier coefficients of $F_{7{\text -}16{\text -}2}$ and $F_{10{\text -}16{\text -}2}$.
 Here, $7$ and $10$ are the weights of $F_{7{\text -}16{\text -}2}$ and $F_{10{\text -}16{\text -}2}$, respectively, 16 refers
to the common level of $F_{7{\text -}16{\text -}2}$ and $F_{10{\text -}16{\text -}2}$, and $2$ refers to the positions of $F_{7{\text -}16{\text -}2}$ and $F_{10{\text -}16{\text -}2}$ in the ordering from \cite{PSYW}.

\begin{table}
\caption{Some Fourier coefficients of $F_{7{\text -}16{\text -}2}$ and $F_{10{\text -}16{\text -}2}$}
\label{Ftable}
$
\begin{array}{lrr}
\toprule
\multicolumn{3}{c}{F_{7{\text -}16{\text -}2}}\\
\midrule
S & d & a(S)\\
\midrule
(2, -53/2, 352)&7& 1\\
(4, -89/2, 496)&15& -45\\
(2, -25/2, 80)& 15& 27\\
(6, -109/2, 496)& 23& -131\\
(3, -13/2, 16)& 23& -131\\ 
(3, -19/2, 32)& 23& -229\\
(7, -21, 64)& 28& 112\\ 
(7, -35, 176)& 28& 112\\ 
(4, -21, 112)& 28& 112\\ 
(4, -53, 704)& 28& -112\\
\bottomrule
\end{array}
$
\qquad
$
\begin{array}{lrr}
\toprule
\multicolumn{3}{c}{F_{10{\text -}16{\text -}2}}\\
\midrule
S & d & a(S)\\
\midrule
(2,-53/2,352)& 7&-1\\
(4,-89/2,496)& 15&145\\
(2,-25/2,80)& 15&-95\\
(6,-109/2,496)& 23&-133\\
(3,-13/2,16)& 23&-133\\
(3,-19/2,32)&23&779\\
(7,-21,64)&28&-256\\
(7,-35,176)&28&-1280\\
(4,-21,112)&28&-256\\
(4,-53,704)&28&768\\
\bottomrule
\end{array}
$ 

\smallskip
{For $S=\begin{bsmallmatrix} \alpha&\beta\\ \beta&\gamma \end{bsmallmatrix}=(\alpha,\beta,\gamma) \in A(16)^+$ we write $d=4\det(S)$}.
\end{table}
\subsection*{Verification} In this subsection we explain how we verified 
 that the seven equations of Corollary~\ref{fourierheckethm}
hold for the newforms from \cite{PSY} and  \cite{PSYW} when the data applies. 

Our first step was to calculate appropriate finite sets of indices at which to evaluate the seven equations
of Corollary \ref{fourierheckethm}.  More precisely, let $i \in \{1,\dots,7\}$. 
By the summary of the data from   \cite{PSYW} and 
the statement of Corollary \ref{fourierheckethm}, we see that, for the $i$-th equation of  Corollary \ref{fourierheckethm}, $S$ must lie in the set $Z_i$ from Table \ref{Scondtable}.
We calculated a finite subset $X_i \subset Z_i$ such that  the $i$-th equation
holds for $S \in X_i$ if and only if the $i$-th equation holds for all $S \in Z_i$. The idea is similar to the method
for obtaining orbit representatives for $\Gamma_0(16)_{\pm} \backslash A(16)^+$ described above.
\begin{table}
\caption{Conditions on $S$ as in Corollary \ref{fourierheckethm} needed for application of the data from \cite{PSYW}}
\label{Scondtable}
\begin{tabular}{cl}
\toprule
\multicolumn{1}{c}{for} & \multicolumn{1}{c}{$S$ must lie in} \\
\toprule
\eqref{fourierheckethmsigmassquareitemeq} & $Z_1=\{\, S \in B(Np^{-2})^+ = A(2)^+ \:|\: 4 \det(S) < \lfloor 500/16 \rfloor =31\, \}$\\
\eqref{mufouriereq} & $Z_2=\{\, S \in B(Np)^+ = A(16)^+ \:|\: 4 \det(S) <  \lfloor 500/4 \rfloor =125\, \}$\\
\eqref{munozerolambdaeq} & $Z_3=\{\, S \in B(Np^{-1})^+ = A(4)^+ \:|\: 4 \det (S) < \lfloor 500/16 \rfloor =31\, \}$\\
\eqref{epsilonfourierformulaeq} & $Z_4=\{\, S \in B(Np^{-1})^+ = A(4)^+ \:|\: 4 \det(S) <  \lfloor 500/16 \rfloor =31\, \}$\\
\eqref{fourierheckethmeq1002} & $Z_5=\{\, S \in B(Np^{-1})^+ = A(4)^+ \:|\: 4 \det (S) < \lfloor 500/4 \rfloor =125\, \}$\\
\eqref{muzerosecondcaseeq}& $Z_6=\{\, S \in B(N)^+ = A(8)^+ \:|\: 4 \det(S) <  \lfloor 500/16 \rfloor =31\, \}$\\
\eqref{muzerofirstcaseeq}&$Z_7=\{\, S \in B(N)^+ = A(8)^+ \:|\: 4 \det(S) <  \lfloor 500/4 \rfloor =125\, \}$\\
\bottomrule
\end{tabular}
\end{table}

To explain this, let the notation be as in Corollary \ref{fourierheckethm},
and consider, say, \eqref{fourierheckethmsigmassquareitemeq} for the examples from \cite{PSYW}.
This is the assertion that 
\begin{equation}
\label{part1Seq}
\sum_{z \in \Z/ 4 \Z}  a(S[\begin{bsmallmatrix} 1& \\ z& 4 \end{bsmallmatrix}]) = 0
\end{equation}
for $S \in Z_1$. To arrive at the set $X_1$ we will take advantage of the 
invariance property of the Fourier coefficient function $a: A(16)^+ \to \C$. 
We recall that the group $\GL(2,\Q)^+$ acts on $A(\Q)^+$ via the formula
$g \cdot S = g S \transpose{g}$ for $g \in \GL(2,\Q)$ and $S \in A(\Q)^+$,
and that with this action,  $\Gamma_0(16)$ preserves the subset $A(16)^+$. By \eqref{gamma0transruleeq} we have
\begin{equation} 
\label{aequivkeq}
a(k \cdot S) = a(S)
\end{equation}
for $k \in \Gamma_0(16)$ and $S \in A(16)^+$. 
The action also satisfies the identity
\begin{equation}
\label{gSBeq}
(g \cdot S) [B] =  (\transpose{B} g \transpose{B}^{-1}) \cdot S[B]
\end{equation}
for $g, B \in \GL(2,\Q)$ and $S \in A(\Q)^+$.
Since \eqref{gSBeq} holds, we see that 
if  $H$ is a  subgroup of $\Gamma_0(2)$ 
(so that $H$ acts on $B(Np^{-2})^+ = A(2)^+$)  and 
$$
g \in H \implies \transpose{\begin{bsmallmatrix} 1& \\ z&4 \end{bsmallmatrix}} g \transpose{\begin{bsmallmatrix} 1& \\ z&4 \end{bsmallmatrix}}^{-1} \in 
\Gamma_0(16)
$$
for $z \in \Z$, 
then \eqref{part1Seq} will hold for $S \in Z_1$ if and only if 
\eqref{part1Seq} holds for all the elements in the orbit $H \cdot S$. Similar
observations apply to the other equations of Corollary \ref{fourierheckethm}. We will define such a group
$H$ for each of the equations \eqref{fourierheckethmsigmassquareitemeq}--\eqref{muzerofirstcaseeq}.
We first define
\begin{align*}
G_1&=\{ \begin{bsmallmatrix} a&b \\ c&d \end{bsmallmatrix} \in \SL(2,\Z)\mid b \equiv c \equiv 0 \Mod{4},\ a \equiv d \Mod{4}\}\\
&=\{ \begin{bsmallmatrix} a&b \\ c&d \end{bsmallmatrix} \in \SL(2,\Z)\mid b \equiv c \equiv 0 \Mod{4}\},\\
G_2&=\{ \begin{bsmallmatrix} a&b \\ c&d \end{bsmallmatrix} \in \SL(2,\Z)\mid b \equiv  0 \Mod{2},\ c \equiv 0 \Mod{16},\  a \equiv d \Mod{2}\}\\
&=\{ \begin{bsmallmatrix} a&b \\ c&d \end{bsmallmatrix} \in \SL(2,\Z)\mid b \equiv  0 \Mod{2},\ c \equiv 0 \Mod{16} \},\\
G_{3}&=\{ \begin{bsmallmatrix} a&b \\ c&d \end{bsmallmatrix} \in \SL(2,\Z)\mid b \equiv  0 \Mod{4},\ c \equiv 0 \Mod{8},\  a \equiv d \Mod{4}\},\\
&=\{ \begin{bsmallmatrix} a&b \\ c&d \end{bsmallmatrix} \in \SL(2,\Z)\mid b \equiv  0 \Mod{4},\ c \equiv 0 \Mod{8}\}.
\end{align*}
We have $G_1 \subset \Gamma_0(4)\subset \Gamma_0(2)$, $G_2 \subset \Gamma_0(16)$, and $G_3 \subset \Gamma_0(8)\subset \Gamma_0(4)$. 
Also, let $\Gamma(16)$ be the subgroup of $g \in \SL(2,\Z)$ such that $g \equiv 1 \Mod{16}$. 

\begin{lemma}
\label{gammactlemma}
If $x \in \Z$, then
\begin{align*}
g \in G_1 &\implies \transpose{\begin{bsmallmatrix} 1& \\ x&4 \end{bsmallmatrix}} g \transpose{\begin{bsmallmatrix} 1& \\ x&4 \end{bsmallmatrix}}^{-1} \in 
\Gamma_0(16),\\
g \in G_2 &\implies \transpose{\begin{bsmallmatrix} 1& \\ x&2 \end{bsmallmatrix}} g \transpose{\begin{bsmallmatrix} 1& \\ x&2 \end{bsmallmatrix}}^{-1},\ 
\transpose{\begin{bsmallmatrix} 1& \\ x/2&1 \end{bsmallmatrix}} g \transpose{\begin{bsmallmatrix} 1& \\ x/2&1 \end{bsmallmatrix}}^{-1} \in 
\Gamma_0(16),\\
g \in G_3 &\implies \transpose{\begin{bsmallmatrix} 1& \\ x&2 \end{bsmallmatrix}} g \transpose{\begin{bsmallmatrix} 1& \\ x&2 \end{bsmallmatrix}}^{-1},\ 
\transpose{\begin{bsmallmatrix} 1& \\ x&4 \end{bsmallmatrix}} g \transpose{\begin{bsmallmatrix} 1& \\ x&4 \end{bsmallmatrix}}^{-1}\in 
\Gamma_0(16).
\end{align*}
We have $[\SL(2,\Z): \Gamma(16)]=3072$, $[\SL(2,\Z):G_1]=24$, $[\SL(2,\Z):G_2]=48$, and $[\SL(2,\Z): G_3]=48$. 
\end{lemma}
\begin{proof}
Let $g = \begin{bsmallmatrix} a&b \\ c&d \end{bsmallmatrix}$. Then 
\begin{align*}
\transpose{\begin{bsmallmatrix} 1& \\ x&4 \end{bsmallmatrix}} g \transpose{\begin{bsmallmatrix} 1& \\ x&4 \end{bsmallmatrix}}^{-1} 
&= 
\begin{bsmallmatrix}
 a+c x & b/4-c x^2/4+(d-a) x/4 \\
 4 c & d-c x 
\end{bsmallmatrix},\\
\transpose{\begin{bsmallmatrix} 1& \\ x&2 \end{bsmallmatrix}} g \transpose{\begin{bsmallmatrix} 1& \\ x&2 \end{bsmallmatrix}}^{-1} &=
\begin{bsmallmatrix}
 a+c x & b/2-c x^2/2+(d -a)x/2 \\
 2 c & d-c x 
\end{bsmallmatrix},\\
\transpose{\begin{bsmallmatrix} 1& \\ x/2&1 \end{bsmallmatrix}} g \transpose{\begin{bsmallmatrix} 1& \\ x/2&1 \end{bsmallmatrix}}^{-1}
&=
\begin{bsmallmatrix}
 a+c x/2 & b-c x^2/4+(d-a) x/2 \\
 c & d-c x/2 
\end{bsmallmatrix}.
\end{align*}
The inclusion statements of the lemma follow from these identities.
The equality $[\SL(2,\Z): \Gamma(16)]=3072$ follows from p.~22 of 
\cite{Shimura}. To prove $[\SL(2,\Z):G_1]=24$, let $u=\begin{bsmallmatrix} &1 \\ -4 & \end{bsmallmatrix}$. Then 
$u$ normalizes $\Gamma_0(4)$, and $u G_1 u^{-1} = \Gamma_0(16)$. It follows that
$[\Gamma_0(4):G_1]=[\Gamma_0(4):\Gamma_0(16)]$. Hence,
\begin{align*}
[\SL(2,\Z): G_1]
&= [\SL(2,\Z): \Gamma_0(4)] [\Gamma_0(4):G_1]\\
&= [\SL(2,\Z): \Gamma_0(4)] [\Gamma_0(4):\Gamma_0(16)]\\
&=[\SL(2,\Z): \Gamma_0(16)]\\
&=24 \qquad \text{(see p.~24 of \cite{Shimura})}.
\end{align*}
Next, let $g=\begin{bsmallmatrix}1& \\ &2 \end{bsmallmatrix}$. We have 
$G_2 \subset \transpose{\Gamma_0(2)}$. Conjugating this inclusion by $g$,
we obtain $\Gamma_0(32)=gG_2 g^{-1} \subset \Gamma_0(2)$. Hence, 
$[\transpose{\Gamma_0(2)} : G_2] = [\Gamma_0(2): \Gamma_0(32)]$ and:
\begin{align*}
[\SL(2,\Z): G_2]
&= [\SL(2,\Z): \transpose{\Gamma_0(2)}] [\transpose{\Gamma_0(2)}: G_2]\\
&= [\SL(2,\Z): \Gamma_0(2)] [\Gamma_0(2): \Gamma_0(32)]\\
&= [\SL(2,\Z):  \Gamma_0(32)]\\
&= 48 \qquad \text{(see p.~24 of \cite{Shimura})}.
\end{align*}
Finally, let $u$ be as before; conjugating the inclusion $G_{3} \subset \Gamma_0(4)$
we have $G_2 = uG_{3}u^{-1} \subset \Gamma_0(4)$. Hence $[\SL(2,\Z):G_3]=[\SL(2,\Z):G_2] = 48$. 
\end{proof}

We define the desired subgroups as follows.

\begin{lemma}
\label{Hdeflemma}
Define $H_1=G_1$, $H_2=G_2$, $H_3=G_3$, $H_4=G_3$, $H_5=G_2$, $H_6=G_3$, and 
$H_7=G_2$. Let $i \in \{1,\dots,7\}$. The group $H_i$ acts on $Z_i$. If $S \in Z_i$,
then the corresponding statement from Corollary \ref{fourierheckethm} as in Table \ref{Scondtable}
holds for $S$ if and only if it holds for all the elements of the orbit $H_i \cdot S$. 
\end{lemma}
\begin{proof}
This follows from the definitions of the groups $H_i$ for $i \in \{1,\dots,7\}$ and 
Lemma \ref{gammactlemma}. 
\end{proof}

By Lemma \ref{cosetrepgammalemma} the  sets
$\cup_{d=1}^{30} X(H_1,2,d)$, $\cup_{d=1}^{124} X(H_2,16,d)$, 
$\cup_{d=1}^{30} X(H_3,4,d)$, $\cup_{d=1}^{30} X(H_4,4,d)$,
$\cup_{d=1}^{124} X(H_5,4,d)$, $\cup_{d=1}^{30} X(H_6,8,d)$, and $\cup_{d=1}^{124} X(H_7,8,d)$
contain orbit representatives for $H_1\backslash Z_1$, $H_2\backslash Z_2$, $H_3\backslash Z_3$,  $H_4\backslash Z_4$, 
$H_5\backslash Z_5$, $H_6\backslash Z_6$, and $H_7\backslash Z_7$, respectively.
However, each of the aforementioned sets contains multiple pairs of elements that define the same orbit. To avoid 
repetitive calculations, for each $i \in \{1,\dots,7\}$ we refined each set to obtain a set $X_i$ of orbit
representatives for $H_i \backslash Z_i$.  
The orders of these sets appear in Table \ref{orbitreptable}. By our discussion,
verifying the $i$-th equation for $i \in \{1,\dots,7\}$ for $S \in Z_i$ is now reduced to verifying this equation for $S \in X_i$. 
\begin{table}
\caption{The number of orbit representatives}
\label{orbitreptable}
\hspace{10ex}
\begin{tabular}{ccc}
\toprule
\multicolumn{1}{c}{$X_i$} &$H_i$ &\multicolumn{1}{c}{order of $X_i \cong H_i \backslash Z_i$} \\
\toprule
$X_1$ &$G_1$& \hphantom{8}208\\
$X_2$ &$G_2$& \hphantom{8}536\\
$X_3$ &$G_3$& \hphantom{8}184\\
$X_4$ &$G_3$& \hphantom{8}\hphantom{8}88\\
$X_5$ &$G_2$& 2320\\
$X_6$ &$G_3$& \hphantom{8}\hphantom{8}88\\
$X_7$ &$G_2$& 1120\\
\bottomrule
\end{tabular}
\end{table}

For every newform $F$ from \cite{PSYW} and $S$ in the appropriate
set $X_i$ we verified the applicable equations from the Corollary \ref{fourierheckethm}. We note that for some equations and some $S \in X_i$ it may happen that not all the Fourier coefficients
were available; as pointed out above, the data from \cite{PSYW} is not quite complete. The calculations had several 
steps. First, we determined the list $L$ of all indices $S$ that appear in the tables defining the newforms from \cite{PSYW}
as discussed above. As mentioned, each element $S$ of $L$ is a distinct orbit representative for 
$\Gamma_0(16)_{\pm}\backslash A(16)^+$ with $4\det(S)<500$. 
To explain the second step, consider, for example, the first equation \eqref{fourierheckethmsigmassquareitemeq}. In this equation there appear the indices
$S[\begin{bsmallmatrix} 1& \\ z&4\end{bsmallmatrix} ]$ for $S \in X_1$ and $z=0,1,2,3$; evidently, each of these 
indices is contained in $ A(16)^+$. For  each $z =0,1,2,3$ we explicitly computed a function $r_z:X_1 \to \{\pm 1\} \times L$
such that $r_z(S) =(\varepsilon_z(S),T_z(S))$ if and only if there exists $k \in \Gamma_0(16)_{\pm}$ such that $k \cdot T_z(S) = S[\begin{bsmallmatrix} 1& \\ z&4\end{bsmallmatrix} ]$ and $\det(k)=\varepsilon_z(S)$. Note
that the functions $r_z$ for $z \in \{0,1,2,3\}$ do not involve any particular newform: this is algebra that applies
to all the examples. We computed similar functions for the remaining six equations of Corollary \ref{fourierheckethm}. Finally, we
calculated the equations. For example, for a particular newform $F$ with weight $k$ and $S \in X_1$,  the first equation \eqref{fourierheckethmsigmassquareitemeq} 
now has the form
$$
\varepsilon_0(S)^k a(T_0(S)) + \varepsilon_1(S)^k a(T_1(S)) + \varepsilon_2(S)^k a(T_2(S)) + \varepsilon_3(S)^k a(T_3(S))=0
$$
where we have used \eqref{gamma0transruleeq}. 
\begin{table}
\caption{Examples of \eqref{fourierheckethmsigmassquareitemeq}  for $F_{7{\text -}16{\text -}2}$ and $F_{10{\text -}16{\text -}2}$}
\label{eq1table}
$
\begin{array}{lrrrrrc}
\toprule
\multicolumn{7}{c}{F_{7{\text -}16{\text -}2}:}\\
\midrule
S \in X_1 & d & A_0& A_1 & A_2 & A_3& \sum_{i=0}^3 A_i\\
\midrule
(257, 241, 226) &4 & 2304 & 2304 & -2304 & -2304 & 0 \\
(242, 453/2,  212)&7& 4352& 9216& -4352& -9216& 0\\
(1 ,  0 , 2)&8& 0 & 4608& 0 & -4608 &  0\\
(259, 243, 228)&12& -27648& -55296& 27648& 55296& 0\\
(244,   457/2, 214 ) &15& 25344& -27648& -25344& 27648& 0 \\
( 1, 0, 4  )&16&  0& -55296 & 0 & 55296& 0\\
(261 , 245  , 230 ) & 20 & 129024& 27648& -129024& -27648& 0\\
(246, 461/2  , 216 ) &23& -109312& 193536& 109312& -193536& 0 \\
(1, 0, 6) &24 & 0&248832& 0& -248832& 0\\
(263, 247, 232) &28& -258048& 442368& 258048& -442368& 0\\
\midrule
\multicolumn{7}{c}{F_{10{\text -}16{\text -}2}:}\\
\midrule
S \in X_1 & d & A_0& A_1 & A_2 & A_3& \sum_{i=0}^3 A_i\\
\midrule
(257, 241, 226) &4& -24576& 24576& 24576& -24576& 0\\
(242, 453/2,  212)&7& 0& 131072& -262144& 131072& 0\\
(1 ,  0 , 2)&8& -81920& 524288& -966656& 524288& 0\\
(259, 243, 228)&12& 5799936& -2949120& -5799936& 2949120& 0\\
(244,   457/2, 214 ) &15& 0& 3276800& -6553600& 3276800& 0\\
( 1, 0, 4  )&16&  -5898240& 6291456& -6684672& 6291456& 0\\
(261 , 245  , 230 ) & 20 & -48168960& 22118400& 48168960& -22118400& 0\\
(246, 461/2  , 216 ) &23& 0& -42336256& 84672512& -42336256& 0\\
(1, 0, 6) &24 & 12419072& -30408704& 48398336& -30408704& 0\\
(263, 247, 232) &28& -168689664& -53477376& 168689664& 53477376& 0\\
\bottomrule
\end{array}
$

\smallskip
For $S=\begin{bsmallmatrix} \alpha&\beta\\ \beta&\gamma \end{bsmallmatrix}=(\alpha,\beta,\gamma)\in X_1$, $d=4\det(S)$; $A_0= a(S[\begin{bsmallmatrix} 1&\\ &4 \end{bsmallmatrix}])$, $A_1=a(S[\begin{bsmallmatrix} 1&\\ 1&4 \end{bsmallmatrix}])$,
$A_2=a(S[\begin{bsmallmatrix} 1&\\ 2&4 \end{bsmallmatrix}])$, $A_3=a(S[\begin{bsmallmatrix} 1&\\ 3&4 \end{bsmallmatrix}])$.
\end{table}
As an illustration, in Table \ref{eq1table} 
we present the data for some of the verifications of \eqref{fourierheckethmsigmassquareitemeq} for the earlier presented newforms $F_{7{\text -}16{\text -}2}$ and $F_{10{\text -}16{\text -}2}$.
By  \cite{PSYW}, the $\mu_2$ eigenvalue of $F_{7{\text -}16{\text -}2}$ is $0$, the $\mu_2$ eigenvalue
of $F_{10{\text -}16{\text -}2}$ is $-4$, the $\lambda_2$ eigenvalue of $F_{7{\text -}16{\text -}2}$ is $-3$, and the $\lambda_2$ eigenvalue of $F_{10{\text -}16{\text -}2}$ is $-2$; also, $\pi_2$ is generic for both of these newforms (as is the case for all
of newforms from \cite{PSYW} described earlier in this section).  Interestingly, for $F_{7{\text -}16{\text -}2}$, the first and third terms already
sum to zero, as do the second and fourth terms. This is a consequence of $\mu_2=0$ for $F_{7{\text -}16{\text -}2}$; for this, use \eqref{fourierheckethmeq1002}. 
\begin{table}
 \caption{Examples of \eqref{fourierheckethmeq1002} for  $F_{7{\text -}16{\text -}2}$}
\label{eq5table}
$
\begin{array}{lrrrc}
\toprule
\multicolumn{5}{c}{F_{7{\text -}16{\text -}2}}\\
\midrule
S \in X_5 & d & A_0 & A_1 & A_0+A_1\\
\midrule
(242, 453/2, 212)& 7& -112& 112& 0\\
(259, 243, 228)& 12& -192& 192& 0\\
(199156, 373529/2, 175144) &15& -1584& 1584& 0\\
(1, 0, 4) & 16& 0& 0& 0\\
(246, 461/2, 216)& 23& -3152& 3152& 0\\
(263, 247, 232)& 28& 9216& -9216& 0\\
(401656, 753329/2, 353228)& 31& 800& -800& 0\\
(1, 0, 8)& 32& 0& 0& 0\\
(250, 469/2, 220)& 39& -3504& 3504& 0\\
(267, 251, 236)& 44& -8896& 8896& 0\\
\bottomrule
\end{array}
$

\smallskip
For $S=\begin{bsmallmatrix} \alpha&\beta\\ \beta&\gamma \end{bsmallmatrix}=(\alpha,\beta,\gamma) \in X_5$, 
$d=4\det(S)$; $A_i= a(S[\begin{bsmallmatrix}1& \\ i & 2 \end{bsmallmatrix}])$.
\end{table}
Some examples of \eqref{fourierheckethmeq1002} appear in Table~\ref{eq5table}. 
The example $F_{10{\text -}16{\text -}2}$ shows that this canceling pairs phenomenon does not hold if $\mu_2 \neq 0$.

For all of these calculations, it is important to note that in Corollary \ref{fourierheckethm} the Siegel modular form $F$
is contained in $S_k(\mathrm{K}(16))$ but is regarded as an element of the larger space $S_k(\mathrm{K}_s(16))$;
as such, $F$ has the Fourier expansion
$$
F(Z) = \sum_{S \in B(16)^+} a(S) e^{2 \pi i \mathrm{Tr}(SZ)}.
$$
Thus, we trivially have $a(S)=0$ for $S \in B(16)^+$ with $S \notin A(16)^+$.

\section{Computing eigenvalues} 
\label{eigenvaluescompsec}

A natural application of Corollary \ref{fourierheckethm} is to the computation of 
the Hecke eigenvalues $\lambda_p$ and $\mu_p$ of a paramodular newform $F$. As mentioned in the introduction, knowledge of $\lambda_p$ and $\mu_p$ is equivalent to knowing the $L$-factor at $p$ of the spin $L$-function attached to $F$. This new method of calculating $\lambda_p$ and $\mu_p$ is easier than the method of ``restriction to a modular curve'' used in \cite{PSY}.

Let the notation be
as in Corollary \ref{fourierheckethm}. The following algorithm then determines the 
Hecke eigenvalues $\lambda_p$ and $\mu_p$ and determines whether or not $\pi_p$ is generic for $p^2 \mid N$ in terms of the Fourier
coefficients of $F$. 

\begin{enumerate}
\item\label{globalline1} Find $S \in A(N)^+ \subset B(Np)^+$ such that $a(S) \neq 0$. Solve \eqref{mufouriereq}
for $\mu_p$. If $\mu_p \neq 0$, then go to \ref{globalline2}. If $\mu_p=0$, then go to \ref{globalline3}.
\item\label{globalline2} Find $S \in B(Np^{-1})^+$ such that $\sum_{x \in \Z/p\Z} a(S[\begin{bsmallmatrix}1& \\ x&p \end{bsmallmatrix}])
\neq 0$; such an $S$ exists by \eqref{fourierheckethmeq1000}. Solve \eqref{munozerolambdaeq} for $\lambda_p$. 
By Corollary \ref{fourierheckethm} the representation $\pi_p$ is generic.
\item\label{globalline3} Using \eqref{fourierheckethmeq10021} determine whether or not $F$ is an eigenvector for $T_{0,1}^s(p)$.
If $F$ is an eigenvector for $T_{0,1}^s(p)$, then find $S \in A(N)^+ \subset B(N)^+$ such that $a(S) \neq 0$, and  
use \eqref{muzerofirstcaseeq} to solve for $\lambda_p$; in this case $\pi_p$ is non-generic.
If $F$ is not an eigenvector for $T_{0,1}^s(p)$, then find $S \in B(N)$ such that $c(S)\neq 0$ (such an $S$ exists by \ref{fourierheckethmsigmazeroitem} of Corollary \ref{fourierheckethm}), and use \eqref{muzerosecondcaseeq} to solve for $\lambda_p$; in this case $\pi_p$
is generic. 
\end{enumerate}

We note that, in case it is known that $F$ is a non-lift, then $\pi_p$ is generic, 
and hence $F$ is not an eigenform for $T_{0,1}^s(p)$. Thus, in this situation, \ref{globalline3} of the algorithm simplifies.

For the examples discussed in Sect.~\ref{examplessec} we found that the above algorithm quickly determines $\lambda_2$ and $\mu_2$
from a few Fourier coefficients.
As an illustration, we consider the newforms $F_{7{\text -}16{\text -}2}$ and $F_{10{\text -}16{\text -}2}$ of \cite{PSYW} that were recalled in 
Sect.~\ref{examplessec}. In Table \ref{eq2table} we show some examples of the determination of $\mu_2$ for 
$F_{7{\text -}16{\text -}2}$ and $F_{10{\text -}16{\text -}2}$ as in \ref{globalline1} of the algorithm. 
We see that $\mu_2(F_{10{\text -}16{\text -}2})=-4$ and 
$\mu_2(F_{7{\text -}16{\text -}2})=0$. Since $\mu_2(F_{10{\text -}16{\text -}2}) \neq 0$, we proceed to \ref{globalline2}
for $F_{10{\text -}16{\text -}2}$ to determine $\lambda_2(F_{10{\text -}16{\text -}2})$. In Table \ref{eq3table}
we show some examples of the determination of $\lambda_2(F_{10{\text -}16{\text -}2})$ as in \ref{globalline2} of the algorithm. 
Since $\mu_2(F_{7{\text -}16{\text -}2})=0$, we proceed to \ref{globalline3} for $F_{7{\text -}16{\text -}2}$ to determine
$\lambda_2(F_{7{\text -}16{\text -}2})$. In Table \ref{eq7table} we present some examples showing that $F_{7{\text -}16{\text -}2}$
cannot be an eigenvector for $T_{0,1}^s(2)$. The argument is by contradiction. If $F_{7{\text -}16{\text -}2}$
were an eigenvector for $T_{0,1}^s(2)$, then \eqref{muzerofirstcaseeq} would have to hold  for all $S \in B(16)^+ =
A(8)^+$; however, in Table \ref{eq7table} we see that this impossible, either because it leads to equations of the form
$0 = C$ where $C \neq 0$, or because it leads to multiple distinct values for $\lambda_2(F_{7{\text -}16{\text -}2})$. Since
$F_{7{\text -}16{\text -}2}$ is not an eigenvector for $T_{0,1}^s(2)$, we may use \eqref{muzerosecondcaseeq} to solve
for $\lambda_2$. This is carried out for some examples in Table \ref{eq6table}.

\begin{table}
 \caption{Determining $\mu_2$  for $F_{7{\text -}16{\text -}2}$ and $F_{10{\text -}16{\text -}2}$}
\label{eq2table}
$
\begin{array}{lrrrrrrr}
\toprule
\multicolumn{8}{c}{F_{7{\text -}16{\text -}2}:}\\
\midrule
S \in X_2 & d & a(S)& A_0 & A_1 & B_0 & B_1 & \mu_2 \\
\midrule
(7666, 14389/2, 6752)& 7& -1& 112& -112&-&-& 0\\
(36340, 68185/2, 31984)& 15& 45& -1584& 1584&-&-& 0\\
(98038, 183917/2, 86256) &23& 131& -3152& 3152&-&-& 0\\
(354631, 332565, 311872) &28& -112& 9216& -9216& -112& 112& 0\\
(205048, 384625/2, 180368) &31& -178& -800& 800&-&-& 0\\
(369658, 693349/2, 325120) &39& 123& 3504& -3504&-&-& 0\\
(604156, 1133129/2, 531312)& 47& 290& 3232& -3232&-&-& 0\\
(1, 2, 16) &48& -192& -55296& 55296& -192& 192& 0\\
(254, 477/2, 224) &55& -1871& 44816& -44816&-&-& 0\\
(271, 255, 240) &60& 1584& -27648& 27648& 1584& -1584& 0\\
\midrule
\multicolumn{8}{c}{F_{10{\text -}16{\text -}2}:}\\
\midrule
S \in X_2 & d & a(S)& A_0 & A_1 & B_0 & B_1 & \mu_2 \\
\midrule
(7666, 14389/2, 6752)& 7& -1& 768& -256&-&-& -4\\
(36340, 68185/2, 31984)& 15& 145& -49920& -24320&-&-& -4\\
(98038, 183917/2, 86256) & 23& -133& -131328& 199424&-&-& -4\\
(354631, 332565, 311872) & 28& -256& 131072& 131072& -256& 768& -4\\
(205048, 384625/2, 180368) & 31& -1402& 327168& 390656&-&-& -4\\
(369658, 693349/2, 325120) & 39& -1341& 1226496& -539904&-&-& -4\\
(604156, 1133129/2, 531312)& 47& -2558& 1853952& -544256&-&-& -4\\
(1, 2, 16) & 48& 0& 2949120& -2949120& 0& 0& *\\
(254, 477/2, 224) & 55& -14575& 3452160& 4010240&-&-& -4\\
(271, 255, 240) & 60& 1280& -3276800& -3276800& 1280& -24320& -4\\
\bottomrule
\end{array}
$

\smallskip
For $S=\begin{bsmallmatrix} \alpha&\beta\\ \beta&\gamma \end{bsmallmatrix}=(\alpha,\beta,\gamma)\in X_2$, $d=4\det(S)$;  
$A_i=a(S[\begin{bsmallmatrix} 1& \\ i&p \end{bsmallmatrix}])$; $B_i = a(S[\begin{bsmallmatrix} 1& \\ i/2 & 1 \end{bsmallmatrix}])$ 
if $2\mid 2 \beta$; the $\mu_2$ column lists the numbers obtained by solving \eqref{mufouriereq} for $\mu_2$; a * 
indicates an equation with both sides equal to zero.
\end{table}

\begin{table}
 \caption{Determining $\lambda_2$   for  $F_{10{\text -}16{\text -}2}$}
\label{eq3table}

\noindent\resizebox{0.98\textwidth}{!}{
$
\begin{array}{lrrrrrrrrrr}
\toprule
\multicolumn{11}{c}{F_{10{\text -}16{\text -}2}}\\
\midrule
S \in X_3 & d & A_0 & A_1 & B_0 & B_1 & C_0 & C_1 & C_2 & C_3 & \lambda_2 \\
\midrule
(242, 453/2, 212) &7& -256& -1280& 0& 262144 & -256& -&768&- &-2\\
(259, 243, 228)& 12& 0& 0& 0& 0& -&-&-&-&*\\
(199156, 373529/2, 175144)& 15& -49920& -24320& 6553600& 0&11520& - & 37120&-& -2\\
(1, 0, 4)& 16& 49152& 0& -12582912& 0& -&-&-&-& -2\\
(246, 461/2, 216)& 23& -131328& 199424& 24903680& -59768832& 199424& -&-131328&-& -2\\
(263, 247, 232) &28& -131072& -131072& 33554432& 33554432& -&-&-&-& -2\\
\bottomrule
\end{array}
$
}

\smallskip
For $S=\begin{bsmallmatrix} \alpha&\beta\\ \beta&\gamma \end{bsmallmatrix}=(\alpha,\beta,\gamma)\in X_3$, $d=4\det(S)$;
$A_i= a(S[\begin{bsmallmatrix} 1&\\ i&2 \end{bsmallmatrix}])$;  $B_i= a(2S[\begin{bsmallmatrix} 1&\\ i&2 \end{bsmallmatrix}])$
$C_i= a(2^{-1} S[\begin{bsmallmatrix} 1&\\ i&4 \end{bsmallmatrix}])$ if $p\mid (\alpha +2\beta i +\gamma i^2)$ and 
is undefined if $p\nmid (\alpha +2\beta i +\gamma i^2)$;  the $\lambda_2$ column lists the numbers obtained by 
solving \eqref{munozerolambdaeq} for $\lambda_2$; a * indicates an equation with both sides equal to zero.
\end{table}

\begin{table}

\caption{Showing that $F_{7{\text -}16{\text -}2}$ is not an eigenvector for $T_{0,1}^s(2)$}
\label{eq7table}
$
\begin{array}{lrrrrrr}
\toprule
\multicolumn{7}{c}{F_{10{\text -}16{\text -}2}}\\
\midrule
S \in X_7 & d & A(S) &A(2S) & A_0 & A_1 & A\\
\midrule
(71666, 134421/2, 63032)& 7& 0& -112& -1&-& \Rightarrow\Leftarrow \\
(199156, 373529/2, 175144)& 15& 0& 432& -27& -& \Rightarrow\Leftarrow \\
(1, 2, 8)& 16& 0& 2304& -& -& \Rightarrow\Leftarrow \\
(246, 461/2, 216)& 23& 0& 3152& 131& -& \Rightarrow\Leftarrow \\
( 263, 247, 232) & 28& 0& -9216& -& -& \Rightarrow\Leftarrow  \\
(205048, 384625/2, 180368)& 31& -178& -800& -178& -&  609/178\\
(1, 0, 8) &32& 0& 0& -& -& 0=0\\
(369658, 693349/2, 325120)& 39& 123& -18384& 195& -& -759/82\\
(604156, 1133129/2, 531312)& 47& 290& -11744& -358& -& -15/2\\
(1, 2, 16)& 48& -192& 33792& -& -& -33/2\\
\bottomrule
\end{array}
$

\smallskip
For $S=\begin{bsmallmatrix} \alpha&\beta\\ \beta&\gamma \end{bsmallmatrix}=(\alpha,\beta,\gamma)\in X_7$, 
$d=4\det(S)$; $A_i= a(2^{-1}S[\begin{bsmallmatrix} 1&\\ i&2 \end{bsmallmatrix}])$ if $2 \mid (\alpha +2 \beta i)$ 
and is undefined if $2 \nmid (\alpha +2 \beta i)$;  $A$ is $\Rightarrow\Leftarrow$ if this instance 
of \eqref{muzerofirstcaseeq} produces a contradiction of the form $0=C$ where $C \neq 0$, $A$ is $0=0$ 
if both sides of \eqref{muzerofirstcaseeq} are zero, and $A$ is a number, putatively equal to $\lambda_2$, 
if \eqref{muzerofirstcaseeq} can be solved for $\lambda_2$.
\end{table}

\begin{table}
 \caption{Determining $\lambda_2$ for  $F_{7{\text -}16{\text -}2}$}
\label{eq6table}
\noindent\resizebox{0.98\textwidth}{!}{
$
\begin{array}{lrrrrrrrrrrrrr}
\toprule
\multicolumn{14}{c}{F_{7{\text -}16{\text -}2}}\\
\midrule
S \in X_6 & d & A(S) &A(2S) & \mc{A_0} & \mc{A_1} &A(4S)& \mc{B_0} & \mc{B_1} & \mc{C_0} & \mc{C_1} & \mc{C_2} & \mc{C_3} & \lambda_2 \\
\midrule
(71666, 134421/2, 63032)& 7 & 0& -112&-1& - & 4352 & 112& 0&  -& - &-1& - &   -3\\
(199156, 373529/2, 175144)& 15 & 0& 432& -27&-& 25344&-1584& 0&45&-&-&-&-3\\
( 1, 2, 8) &16 & 0& 2304& -& -& -110592& -& 0& -& -& -& -& -3\\
(246, 461/2, 216)&23& 0& 3152& 131& -& -17152& -3152& 0& -& -& -229& -&  -3\\
(263, 247, 232) & 28& 0& -9216& -& -& 442368& -& 0& -& -& -& -& -3\\
\bottomrule
\end{array}
$}

\smallskip
For $S=\begin{bsmallmatrix} \alpha&\beta\\ \beta&\gamma \end{bsmallmatrix}=(\alpha,\beta,\gamma) \in X_6$, 
$d=4\det(S)$; $A_i= a(2^{-1}S[\begin{bsmallmatrix} 1&\\ i&2 \end{bsmallmatrix}])$ if $2 \mid (\alpha +2 \beta i)$ 
and is undefined if $2 \nmid (\alpha +2 \beta i)$;  $B_i= a(S[\begin{bsmallmatrix} 1&\\ i&2 \end{bsmallmatrix}])$ 
if $2 \mid (\alpha +2 \beta i)$ and is undefined if $2 \nmid (\alpha +2 \beta i)$; 
$C_i= a(2^{-2} S[\begin{bsmallmatrix} 1&\\ i&4 \end{bsmallmatrix}])$ if $4\mid (\alpha +2\beta i +\gamma i^2)$ 
and is undefined if $p\nmid (\alpha +2\beta i +\gamma i^2)$; the $\lambda_2$ column lists the numbers 
obtained by solving \eqref{muzerosecondcaseeq} for $\lambda_2$.
\end{table}

\section{A recurrence relation}
\label{recurrencesec}

In this final section, assuming that $F$ is as in Corollary \ref{fourierheckethm}, we
prove that the radial Fourier coefficients $a(p^t S)$ of $F$ for $t \geq 0$, $p^2 \mid N$, satisfy
a recurrence relation determined by the spin $L$-factor $L_p(s,F)$ of $F$ at $p$. This 
theorem extends results known in other situations (e.g., Sect.~4.3.2 of \cite{MR884891}).
We begin with the following lemma.

\begin{lemma}
\label{psformulacor}
Let $N$ and $k$ be integers such that $N>0$ and $k>0$, and let $F \in S_k(\mathrm{K}(N))$. 
Assume that $F$ is a newform and is an eigenvector for the Hecke operators $T(1,1,q,q)$ and $T(1,q,q,q^2)$
for all but finitely many of the primes $q$ of $\Z$ such that $q \nmid N$;  by Theorem
\ref{Ralftheorem}, $F$ is an eigenvector for $T(1,1,q,q)$, $T(1,q,q,q^2)$ and $w_q$ for all primes $q$ of $\Z$.
Let $\lambda_q, \mu_q \in \C$  be such that 
$
T(1,1,q,q)F=q^{k-3}\lambda_qF
$
and 
$
T(1,q,q,q^2)F=q^{2(k-3)}\mu_qF
$
for all primes $q$ of $\Z$. Regard $F$ as an element of $S_k(\mathrm{K}_s(N))$, and let 
\begin{align*}
F(Z)=\sum_{S\in B(N)^+}a(S) \E^{2\pi \I \mathrm{Tr}(SZ)}
\end{align*}
be the Fourier expansion of $F$. Let $\pi\cong \otimes_{v \leq \infty} \pi_v$ be as in Theorem \ref{globalalgthm}. 
Let $p$ be a prime of $\Z$ with $v_p(N)\geq2$. 
\begin{enumerate}
\item \label{psformulacori1} Assume that $\mu_p \neq 0$. If $S \in B(N)^+$, then there is 
a formal identity
\begin{equation}
\label{psformulacoreq15}
\sum_{t=0}^\infty a(p^t S) X^t
=
\dfrac{N_1(X,S)}{D_1(X)}
\end{equation}
where
\begin{align}
N_1(X,S)
&=a(S) +\big( a(pS) -p^{k-3} \lambda_p a(S) \big)X \nonumber \\
&\quad+\big( a(p^2 S) - p^{k-3} \lambda_p a(pS) + p^{2k-5} (\mu_p +p^2) a(S) \big) X^2,
\label{psformulacoreq16}\\
D_1(X)&= 1-p^{k-3} \lambda_p X + p^{2k-5} (\mu_p +p^2) X^2.\label{psformulacoreq17}
\end{align}
\item \label{psformulacori2} Assume that $\mu_p =0$. If $S \in B(N)^+$ and $F$ is not an eigenvector for $T_{0,1}^s(p)$ (so that 
$\pi_p$ is generic by Theorem \ref{globalalgthm}), then there is a formal identity
\begin{equation}
\label{psformulacoreq5}
\sum_{t=0}^\infty a(p^t S) X^t
=
\dfrac{N_2(X,S)}{D_2(X)}
\end{equation}
where
\begin{align}
N_2(X,S)
&=a(S) + \big(a(pS)- p^{k-3}\lambda_p  a(S) \big) X\nonumber \\
&\quad+ \big( a(p^2 S) - p^{k-3} \lambda_p a(pS) + p^{2k-3} a(S) \big) X^2,\label{psformulacoreq6}\\
D_2(X)&= 1-p^{k-3} \lambda_p X + p^{2k-3} X^2.\label{psformulacoreq7}
\end{align}
If  $S \in B(N)^+$ and $F$ is  an eigenvector for $T_{0,1}^s(p)$ (so that 
$\pi_p$ is non-generic by Theorem \ref{globalalgthm}), then there is a formal identity
\begin{equation}
\label{psformulacoreq50}
\sum_{t=0}^\infty a(p^t S) X^t
=
\dfrac{N_3(X,S)}{D_3(X)}
\end{equation}
where
\begin{align}
N_3(X,S)
&=a(S) + \big(a(pS)- p^{k-2}(1+p)^{-1}  \lambda_p a(S) \big)X,\label{psformulacoreq60}\\
D_3(X)&= 1-p^{k-2}(1+p)^{-1}  \lambda_p X.\label{psformulacoreq70}
\end{align}
\end{enumerate}
\end{lemma}
\begin{proof}\ref{psformulacori1}.  For $S \in B(Np^{-1})^+$ define
\begin{equation}
\label{psformulacoreq1}
a'(S) = \sum_{x \in \Z/p\Z} a(S[\begin{bsmallmatrix} 1& \\ x&p \end{bsmallmatrix} ]).
\end{equation}
We will first prove that if $S \in B(Np^{-1})^+$, then there is a formal identity
\begin{equation}
\label{psformulacoreq2}
\sum_{t=0}^\infty a'(p^t S) X^t = \dfrac{N_0(X,S)}{D_1(X)}
\end{equation}
where
\begin{align}
N_0(X,S) & = a'(S) + \big( a'(pS) -p^{k-3} \lambda_p a'(S) \big) X, \label{psformulacoreq3}\\
D_1(X) & = 1-p^{k-3} \lambda_p X + p^{2k-5} (\mu_p +p^2) X^2.\label{psformulacoreq4}
\end{align}
Let $S \in B(Np^{-1})^+$. Let $i \in \Z$ be such that $i \geq 1$ and define
$S' = p^i S$. Then $S' \in B(Np^{-1})^+$. Write $S' = \begin{bsmallmatrix} \alpha & \beta \\ \beta & \gamma \end{bsmallmatrix}$. 
By \eqref{munozerolambdaeq} we have 
\begin{align*}
\lambda_p \sum_{x\in\Z/p\Z} a(S'[\begin{bsmallmatrix}1&\\x&p\end{bsmallmatrix}])
&=\sum_{x\in\Z/p\Z}p^{3-k}   a(pS'[\begin{bsmallmatrix}1&\\x&p\end{bsmallmatrix}])\nonumber \\
&\quad+\sum_{\substack{z\in\Z/p^2\Z\\p\mid(\alpha+2\beta z+\gamma z^2)}}p  a(p^{-1}S'[\begin{bsmallmatrix}1&\\z&p^2\end{bsmallmatrix}]).
\end{align*}
Since $p$ divides $\alpha$, $2\beta$, and $\gamma$, by \eqref{psformulacoreq1}  this is 
\begin{equation}
\label{aprimeeq1}
\lambda_p a'(p^i S) 
=p^{3-k}   a'(p^{i+1}S)+\sum_{z\in\Z/p^2\Z}p a(p^{-1} S'[\begin{bsmallmatrix}1&\\z&p^2\end{bsmallmatrix}]).
\end{equation}
Applying \eqref{epsilonfourierformulapreeq} to $p^{-1}S'\in B(Np^{-1})^+$, we also have 
\begin{align}
\sum_{z\in\Z/p^2\Z} pa(p^{-1} S'[\begin{bsmallmatrix}1&\\z&p^2\end{bsmallmatrix}]) 
&= (\mu_p+p^2)p^{k-2}\sum_{x\in\Z/p\Z}a(p^{-1}S'[\begin{bsmallmatrix}1&\\x&p\end{bsmallmatrix}]) \nonumber \\
&= (\mu_p+p^2)p^{k-2} a'(p^{i-1} S). \label{aprimeeq2}
\end{align}
Substituting \eqref{aprimeeq2} into \eqref{aprimeeq1},
we now have
$$
\lambda_p a'(p^i S) 
=p^{3-k}   a'(p^{i+1}S)+ p^{k-2}(\mu_p+p^2) a'(p^{i-1} S). 
$$
Rewriting, we have proven that for $t \in \Z$, $t \geq 2$, and $S \in B(Np^{-1})^+$, 
\begin{equation}
a'(p^t S) -p^{k-3} \lambda_p a'(p^{t-1} S)  +p^{2k-5} (\mu_p+p^2) a'(p^{t-2} S) =0. \label{part1recurreleq}
\end{equation}
It is straightforward to derive \eqref{psformulacoreq2} from this recurrence relation. Let 
$S \in B(N)^+$. Then
\begin{align*}
&\mu_p \sum_{t=0}^\infty a(p^t S) X^t 
= \mu_p a(S) +\mu_p X\sum_{t=1}^\infty a(p^{t} S) X^{t-1} \\
&\qquad= \mu_p a(S) + X\sum_{t=0}^\infty\mu_p a(p^{t+1} S) X^{t} \\
&\qquad= \mu_p a(S) 
+X\sum_{t=0}^\infty  \sum_{x \in \Z/p\Z} p^{3-k} a(p^{t+1} S[\begin{bsmallmatrix} 1& \\ x&p \end{bsmallmatrix} ])X^t\\
&\qquad\quad-X\sum_{t=0}^\infty \sum_{x \in \Z/p\Z} p a(p^{t+1} S[\begin{bsmallmatrix} 1& \\ xp^{-1}&1 \end{bsmallmatrix} ])X^t
\qquad\text{(by \eqref{mufouriereq})}\\
&\qquad= \mu_p a(S) 
+ X\sum_{t=0}^\infty  \sum_{x \in \Z/p\Z} p^{3-k} a(p^{t+1}S[\begin{bsmallmatrix} 1& \\ x&p \end{bsmallmatrix} ])X^t\\
&\qquad\quad-X\sum_{t=0}^\infty \sum_{x \in \Z/p\Z} p 
a(p^{t+1}S[\begin{bsmallmatrix} 1& \\ &p^{-1} \end{bsmallmatrix} ][\begin{bsmallmatrix} 1& \\ x&p\vphantom{p^{-1}} \end{bsmallmatrix} ])X^t\\
&\qquad= \mu_p a(S) 
+ p^{3-k}X\sum_{t=0}^\infty   a'(p^{t+1} S)X^t \\
&\qquad\quad- pX\sum_{t=0}^\infty 
a'(p^{t+1}S[\begin{bsmallmatrix} 1& \\ &p^{-1} \end{bsmallmatrix} ])X^t \\
&\qquad= \mu_p a(S) 
+ p^{3-k}X N_0(X,pS) D_1(X)^{-1} \\
&\qquad\quad- pX 
N_0(X,pS[\begin{bsmallmatrix} 1& \\ &p^{-1} \end{bsmallmatrix} ]) D_1(X)^{-1} \qquad\text{(use \eqref{psformulacoreq2})}\\
&\qquad= \big( \mu_p a(S) D_1(X)
+ p^{3-k} a'(pS)X+ p^{3-k}a'(p^2 S)X^2 - \lambda_p a'(pS)X^2 \\
&\qquad\quad-  
p a'(pS[\begin{bsmallmatrix} 1& \\ &p^{-1} \end{bsmallmatrix} ])X
-pa'(p^2S[\begin{bsmallmatrix} 1& \\ &p^{-1} \end{bsmallmatrix} ])X^2\\
&\qquad\quad+p^{k-2}\lambda_p a'(pS[\begin{bsmallmatrix} 1& \\ &p^{-1} \end{bsmallmatrix} ]))X^2 \big) D_1(X)^{-1}
\qquad \text{(by \eqref{psformulacoreq3})}\\
&\qquad = ( \mu_p a(S) D_1(X) 
+( p^{3-k}a'(pS) -pa'(pS[\begin{bsmallmatrix} 1& \\ & p^{-1} \end{bsmallmatrix} ])) X\\
&\qquad\quad+(p^{3-k} a'(p^2S) -p a'(p^2 S[\begin{bsmallmatrix} 1& \\ & p^{-1} \end{bsmallmatrix}]))X^2\\
&\qquad\quad-p^{k-3} \lambda_p (p^{3-k} a'(pS) -p a'(pS[\begin{bsmallmatrix} 1& \\ & p^{-1} \end{bsmallmatrix}])) X^2 ) 
D_1(X)^{-1}\\
&\qquad = ( \mu_p a(S) D_1(X) 
+\mu_p a(pS) X
+\mu_p a(p^2 S) X^2\\
&\qquad\quad-p^{k-3} \lambda_p\mu_p a(pS) X^2 ) 
D_1(X)^{-1}\qquad\text{(by \eqref{mufouriereq})}\\
&\qquad = \mu_p ( a(S) +( a(pS) -p^{k-3} \lambda_p a(S) )X \\
&\qquad\quad+( a(p^2 S) - p^{k-3} \lambda_p a(pS) + p^{2k-5} (\mu_p +p^2) a(S) ) X^2 ) D_1(X)^{-1}.
\end{align*}
For the last equality we used \eqref{psformulacoreq4}. 
Canceling $\mu_p$ from both sides of the last equation now yields \eqref{psformulacoreq15}. 

\ref{psformulacori2}. Assume first that $F$ is not an eigenvector for $T_{0,1}^s(p)$. 
Let $S \in B(N)^+$. Let $i \in \Z$ be such that $i \geq 1$ and define $S'=p^i S$. 
Then $S' \in B(N)^+$; let $S' = \begin{bsmallmatrix} \alpha & \beta \\ \beta & \gamma \end{bsmallmatrix}$. 
By \eqref{muzerosecondcaseeq} we have  
\begin{align*}  
&p^{3-k} \lambda_p a(pS')+ \lambda_p \sum_{\substack{x\in\Z/p\Z\\ p\mid (\alpha+2\beta x)}}
p  a(p^{-1}S'[\begin{bsmallmatrix}1&\\ x&p\end{bsmallmatrix}])\\
&\qquad=p^3 a(S')+p^{6-2k}  a(p^2S')+\sum_{\substack{y\in\Z/p\Z\\ p\mid(\alpha+2\beta y)}}
p^{4-k}  a(S'[\begin{bsmallmatrix}1&\\y&p\end{bsmallmatrix}])\nonumber \\
&\qquad\quad+\sum_{\substack{z\in\Z/p^2\Z\\ p^2 \mid (\alpha+ 2\beta z+\gamma z^2)}}
p^{2}  a(p^{-2}S'[\begin{bsmallmatrix}1&\\z&p^2\end{bsmallmatrix}]).
\end{align*}
Since $p$ divides $\alpha$ and $2\beta$, and since $p^2$ divides $\gamma$, this equation is
\begin{align}  
&p^{3-k} \lambda_p a(pS')+ \lambda_p \sum_{\substack{x\in\Z/p\Z}}
p  a(p^{-1}S'[\begin{bsmallmatrix}1&\\ x&p\end{bsmallmatrix}])\nonumber \\
&\qquad=p^3 a(S')+p^{6-2k}  a(p^2S')+\sum_{\substack{y\in\Z/p\Z}}
p^{4-k}  a(S'[\begin{bsmallmatrix}1&\\y&p\end{bsmallmatrix}])\nonumber \\
&\qquad\quad+\sum_{\substack{z\in\Z/p^2\Z\\ p^2 \mid (\alpha+ 2\beta z)}}
p^{2}  a(p^{-2}S'[\begin{bsmallmatrix}1&\\z&p^2\end{bsmallmatrix}]). \label{psformulacoreq100}
\end{align}
Since $ p^{-1} S', S' \in B(N)^+ \subset B(Np^{-1})^+$ 
the identity
\eqref{fourierheckethmeq1002} 
implies that 
\begin{equation}
 \sum_{\substack{x\in\Z/p\Z}}
  a(p^{-1}S'[\begin{bsmallmatrix}1&\\ x&p\end{bsmallmatrix}]) = 0
\quad\text{and}\quad
\sum_{\substack{y\in\Z/p\Z}}
  a(S'[\begin{bsmallmatrix}1&\\y&p\end{bsmallmatrix}])=0. \label{psformulacoreq101}
\end{equation}
We also have
\begin{align}
&\sum_{\substack{z\in\Z/p^2\Z\\ p^2 \mid (\alpha+ 2\beta z)}} a(p^{-2} S' [\begin{bsmallmatrix} 1& \\ z & p^2 \end{bsmallmatrix} ])
 = \sum_{x \in \Z/ p\Z} \sum_{\substack{y \in \Z/ p \Z\\ p^2 \mid (\alpha +2 \beta (x+yp))} }
a(p^{-2} S'[ \begin{bsmallmatrix} 1& \\ x+yp & p^2 \end{bsmallmatrix} ]) \nonumber \\
&\qquad = \sum_{x \in \Z/ p\Z} \sum_{\substack{y \in \Z/ p \Z\\ p^2 \mid (\alpha +2 \beta x)} } 
a( (p^{-2} S'[ \begin{bsmallmatrix} 1& \\ x & p \end{bsmallmatrix} ])
[ \begin{bsmallmatrix} 1& \\ y & p \end{bsmallmatrix} ]) \qquad \text{(note that $p \mid 2\beta$)}\nonumber \\
&\qquad = \sum_{\substack{x \in \Z/ p \Z\\ p^2 \mid (\alpha +2 \beta x)} }  \sum_{y \in \Z/ p \Z } 
a( (p^{-2} S'[ \begin{bsmallmatrix} 1& \\ x & p \end{bsmallmatrix} ])
[ \begin{bsmallmatrix} 1& \\ y & p \end{bsmallmatrix} ]) \nonumber \\
&\qquad = \sum_{\substack{x \in \Z/ p \Z\\ p^2 \mid (\alpha +2 \beta x)} } 0 \nonumber \qquad \text{(by \eqref{fourierheckethmeq1002})}\\
&\qquad= 0. \label{psformulacoreq102}
\end{align}
Substituting \eqref{psformulacoreq101} and \eqref{psformulacoreq102} into \eqref{psformulacoreq100} now gives
\begin{align*}
p^{3-k} \lambda_p a(pS')
&=p^3 a(S')+p^{6-2k}  a(p^2S')\\
p^{3-k} \lambda_p a(p^{i+1} S)
&=p^3 a(p^i S)+p^{6-2k}  a(p^{i+2} S).
\end{align*}
Rewriting, we have proven that if $t \in \Z$ with $t \geq 3$, then 
$$
a(p^t S) - p^{k-3} \lambda_p a(p^{t-1} S) + p^{2k-3} a(p^{t-2} S) =0.
$$
This recurrence relation implies \eqref{psformulacoreq5}. Finally, assume that 
$F$ is an eigenvector for $T_{0,1}^s(p)$. 
Let $S \in B(N)^+$. Let $i \in \Z$ be such that $i \geq 1$ and define $S'=p^i S$. 
Then $S' \in B(N)^+$; let $S' = \begin{bsmallmatrix} \alpha & \beta \\ \beta & \gamma \end{bsmallmatrix}$. 
By \eqref{muzerofirstcaseeq} we have 
$$
\lambda_p  a(S') 
=(1+p) p^{2-k}  a(pS')+\sum_{\substack{x\in\Z/p\Z\\ p \mid (\alpha +  2\beta x)}}
(1+p)  a(p^{-1}S'[\begin{bsmallmatrix}1&\\x&p\end{bsmallmatrix}]).
$$
Since $p$ divides $\alpha$ and $2\beta$, this is
$$
\lambda_p  a(S') 
=(1+p) p^{2-k}  a(pS')+\sum_{\substack{x\in\Z/p\Z}}
(1+p)  a(p^{-1}S'[\begin{bsmallmatrix}1&\\x&p\end{bsmallmatrix}]).
$$
Also, since $p^{-1} S' \in B(N)^+ \subset B(Np^{-1})^+$, \eqref{fourierheckethmeq1002}
implies that 
$$
\sum_{\substack{x\in\Z/p\Z}}
  a(p^{-1}S'[\begin{bsmallmatrix}1&\\x&p\end{bsmallmatrix}])=0.
$$
Therefore, 
\begin{align*}
\lambda_p  a(S') 
&=(1+p) p^{2-k}  a(pS')\\
\lambda_p  a(p^i S) 
&=(1+p) p^{2-k}  a(p^{i+1} S). 
\end{align*}
Rewriting, we have proven that if $t \in \Z$ and $t \geq 2$, then 
\begin{equation}
\label{SKrecureq}
a(p^t S) -(1+p)^{-1} p^{k-2} \lambda_p a(p^{t-1} S) =0. 
\end{equation}
This recurrence relation implies \eqref{psformulacoreq50}. 
\end{proof}

We can now prove the main result of this section.

\begin{theorem}
\label{recurrencetheorem}
Let $N$ and $k$ be integers such that $N>0$ and $k>0$, and let $F \in S_k(\mathrm{K}(N))$. 
Assume that $F$ is a newform and is an eigenvector for the Hecke operators $T(1,1,q,q)$ and $T(1,q,q,q^2)$
for all but finitely many of the primes $q$ of $\Z$ such that $q \nmid N$;  by Theorem
\ref{Ralftheorem}, $F$ is an eigenvector for $T(1,1,q,q)$, $T(1,q,q,q^2)$ and $w_q$ for all primes $q$ of $\Z$.
Let $\lambda_q, \mu_q \in \C$  be such that 
\begin{align*}
T(1,1,q,q)F&=q^{k-3}\lambda_qF,\\
T(1,q,q,q^2)F&=q^{2(k-3)}\mu_qF
\end{align*}
for all primes $q$ of $\Z$. Regard $F$ as an element of $S_k(\mathrm{K}_s(N))$, and let 
\begin{align*}
F(Z)=\sum_{S\in B(N)^+}a(S) \E^{2\pi \I \mathrm{Tr}(SZ)}
\end{align*}
be the Fourier expansion of $F$. 
Let $p$ be a prime of $\Z$ with $v_p(N)\geq2$. If $S \in B(N)^+$, then there is a formal identity
of power series in $p^{-s}$
\begin{equation}
\label{recurrencetheoremeq1}
\sum\limits_{t=0}^\infty \dfrac{a(p^t S)}{p^{ts}} = N(p^{-s},S) L_p(s,F)
\end{equation}
where 
\begin{align}
N(p^{-s},S) &= a(S) + \big( a(pS) - p^{k-3} \lambda_p a(S) \big) p^{-s}\nonumber \\
&\quad + \big( a(p^2 S) - p^{k-3} \lambda_p a(pS) + p^{2k-5} (\mu_p +p^2) a(S) \big) p^{-2s}
\label{recurrencetheoremeq2}
\end{align}
and 
\begin{equation}
\label{recurrencetheoremeq3}
L_p(s,F) = \dfrac{1}{1 - p^{k-3} \lambda_p p^{-s} +p^{2k-5} (\mu_p + p^2 ) p^{-2s} }
\end{equation}
is the spin $L$-factor of $F$ at $p$ (e.g., see \cite{MR2887605}, p.~547).
\end{theorem}
\begin{proof}
The identity \eqref{recurrencetheoremeq1} follows from Lemma \ref{psformulacor} if $\mu_p \neq 0$ or 
if $\mu_p =0$ and $F$ is not an eigenvector for $T_{0,1}^s(p)$. Assume that $\mu_p=0$ and $F$
is an eigenvector for~$T_{0,1}^s(p)$. Let $\pi \cong \otimes_{v \leq \infty} \pi_v$ be as in 
Theorem \ref{Ralftheorem}. By \ref{globalalgthmitem3} of Theorem \ref{globalalgthm} $\pi_p$
is non-generic and by Theorem \ref{Ralftheorem} $\pi_p$ is a Saito-Kurokawa representation. 
An inspection of Table \ref{levelsandeigenvaluestable} shows that $\lambda_p^2 = p^2 (p+1)^2$. 
Let $S \in B(N)^+$. 
By \eqref{psformulacoreq50} there is a formal identity
\begin{equation}
\label{recurrencetheoremeq4}
\sum_{t=0}^\infty \dfrac{a(p^t S) }{p^{ts}}
=
\dfrac{a(S) +\big( a(pS) -p^{k-2}(1+p)^{-1} \lambda_p a(S) \big) p^{-s}}{1 - p^{k-2} (1+p)^{-1} \lambda_p p^{-s}}.
\end{equation}
To prove \eqref{recurrencetheoremeq1} we will multiply the numerator and denominator in \eqref{recurrencetheoremeq4}
by   $1-p^{k-2}(1+p)^{-1} \lambda_p p^{-(s+1)}$. Using $\lambda_p^2 = p^2 (p+1)^2$  and $\mu_p=0$ it is easy
to verify that 
$$
(1 - p^{k-2} (1+p)^{-1} \lambda_p p^{-s})(1-p^{k-2}(1+p)^{-1} \lambda_p p^{-(s+1)}) = L_p(s,F)^{-1}.
$$
And:
\begin{align*}
&(a(S) +( a(pS) -p^{k-2}(1+p)^{-1} \lambda_p a(S) ) p^{-s})(1-p^{k-2}(1+p)^{-1} \lambda_p p^{-(s+1)})\\
&\qquad = a(S) +(a(pS)-p^{k-3}\lambda_p a(S))p^{-s} \\
&\qquad\quad+ (-p^{k-3} (1+p)^{-1} \lambda_p a(pS) +p^{2k-3} a(S)) p^{-2s} \\
&\qquad = a(S) +(a(pS)-p^{k-3}\lambda_p a(S))p^{-s} \\
&\qquad\quad+ (a(p^2S) -a(p^2 S) -p^{k-3} (1+p)^{-1} \lambda_p a(pS) +p^{2k-3} a(S)) p^{-2s} \\
&\qquad = a(S) +(a(pS)-p^{k-3}\lambda_p a(S))p^{-s} \\
&\qquad\quad+ (a(p^2S) -p^{k-2}(1+p)^{-1}\lambda_p a(pS)\\
&\qquad\quad-p^{k-3} (1+p)^{-1} \lambda_p a(pS) +p^{2k-3} a(S)) p^{-2s} \qquad\text{(by \eqref{SKrecureq})}\\
&\qquad = a(S) +(a(pS)-p^{k-3}\lambda_p a(S))p^{-s} \\
&\qquad\quad+ (a(p^2S) -p^{k-3}\lambda_pa(pS)+p^{2k-3} a(S)) p^{-2s} \\
&\qquad=N (p^{-s},S).
\end{align*}
This completes the proof.
\end{proof}

Finally, we remark  that the examples of \cite{PSY} and  \cite{PSYW}
can be used to exhibit $F$ and $S$ such that the coefficients of $p^{-ts}$ for $t=0,1$ and $2$ in $N(p^{-s},S)$ are non-zero. 

%% file: SKMS_appendix.tex
\appendix
\chapter{Tables}
In this appendix we present a number of tables summarizing 
some of the results of the text. The following is a guide
to these tables. 

\subsection*{Table \ref{maintable}.\ Non-supercuspidal 
representations of $\GSp(4,F)$} This table lists all the non-supercuspidal,
irreducible, admissible representations $\pi$ of $\GSp(4,F)$; see
Sect.~\ref{repsec} for more details. In addition, this table
shows in the ``tempered'' column the conditions required for $\pi$ to be tempered;
indicates with a $\bullet$ in the ``ess.~$L^2$'' column
when $\pi$ is essentially square-integrable, i.e., square-integrable
after an appropriate twisting; and indicates with a $\bullet$
in the ``generic'' column if $\pi$ is generic.

\subsection*{Table \ref{noklingentable}.\ Non-paramodular representations}
This table lists all the irreducible, admissible
representations $\pi$ of $\GSp(4,F)$ with trivial central
character that are not paramodular. In this table,
$\pi$ is not paramodular if and only if the listed condition on the 
defining data for $\pi$ is satisfied. This table follows from
Theorem 3.4.3 of \cite{NF} and Table A.12 of \cite{NF}.

\subsection*{Table \ref{dimensionstable}.\ Stable Klingen dimensions}
This table lists, 
for every irreducible, admissible representation $(\pi,V)$ of $\GSp(4,F)$ 
with trivial central character, the paramodular level $N_\pi$,
the stable Klingen level $N_{\pi,s}$, 
the dimensions of the spaces $V_s(n)$ for integers $n \geq N_{\pi,s}$,
the quotient stable Klingen level $\bar N_{\pi,s}$, 
the dimensions of the spaces $\bar V_s(n)$ for integers $n \geq \bar N_{\pi,s}$,
the category of $\pi$, and
some additional information for some $\pi$. 
The symbols and abbreviations \catonebox, \cattwobox, ``SK'', non-unit., and one-dim., 
stand for \catone and \cattwo (see p.~\pageref{catdefpageanchor}), Saito-Kurokawa (see p.~\pageref{SKsubsec}), 
non-unitary, and one-dimensional, respectively.
The entry for $\bar N_{\pi,s}$ is --  if $V_s(n)=V(n)$
for all integers $n\geq 0$, so that $\bar N_{\pi,s}$ is not defined. See Theorem \ref{dimensionstheorem}.

\subsection*{Table \ref{levelsandeigenvaluestable}.\ Hecke eigenvalues}
For every irreducible, admissible representation $(\pi,V)$ of $\GSp(4,F)$ 
with trivial central character, this table lists 
the paramodular level $N_\pi$,
the Atkin-Lehner eigenvalue~$\varepsilon_\pi$, and
the paramodular Hecke eigenvalues $\lambda_\pi$ and $\mu_\pi$ from Table A.14 of~\cite{NF};
also listed are 
the category of $\pi$,  
and some additional information for some~$\pi$. 
The symbols and abbreviations \catonebox, \cattwobox, ``SK'', non-unit., and one-dim., 
stand for \catone and \cattwo (see p.~\pageref{catdefpageanchor}), Saito-Kurokawa (see p.~\pageref{SKsubsec}), 
non-unitary, and one-dimensional, respectively.
For typesetting reasons, some of the eigenvalues are given as below.
\begin{align*}
 \text{(A1)}&=q^{3/2}\sigma(\varpi)\big(1+\chi_1(\varpi)+\chi_2(\varpi)+\chi_1(\varpi)\chi_2(\varpi)\big),\\
 \text{(A2)}&=q^2\big(\chi_1(\varpi)+\chi_2(\varpi)+\chi_1(\varpi)^{-1}+\chi_2(\varpi)^{-1}+1-q^{-2}\big),\\
 \text{(B1)}&=q^{3/2}(\sigma(\varpi)+\sigma(\varpi)^{-1})+(q+1)(\sigma\chi)(\varpi),\\
 \text{(B2)}&=q^{3/2}(\chi(\varpi)+\chi(\varpi)^{-1}),\\
 \text{(C1)}&=q^{3/2}(\sigma(\varpi)+\sigma(\varpi)^{-1})+q(q+1)(\sigma\chi)(\varpi),\\
 \text{(C2)}&=q^{3/2}(q+1)(\chi(\varpi)+\chi^{-1}(\varpi))+q^2-1,\\
 \text{(D)}&=q^2(\chi(\varpi)+\chi^{-1}(\varpi)+q+1)+q-1.
\end{align*}
The last column provides additional information about the representation. 

\subsection*{Table \ref{Iwahoritable}.\ Characteristic polynomials of $T_{0,1}^s$ and $T_{1,0}^s$ on $V_s(1)$  in terms of inducing data} 
This table lists, for every irreducible, admissible representation $(\pi,V)$ of $\GSp(4,F)$
with trivial central character such that $V_s(1) \neq 0$, the characteristic polynomials of 
$T_{0,1}^s$ and $T_{1,0}^s$ in terms of the inducing data for $\pi$. 
See Theorem \ref{charpolyiwahoritheorem}.

\subsection*{Table \ref{eigenIwahoritable}.\ Characteristic polynomials of $T_{0,1}^s$ and $T_{1,0}^s$ on $V_s(1)$  in terms of paramodular eigenvalues} 
This table lists, for every irreducible, admissible representation $(\pi,V)$ of $\GSp(4,F)$
with trivial central character such that ${V_s(1) \neq 0}$, the characteristic polynomials of 
$T_{0,1}^s$ and $T_{1,0}^s$ in terms of the paramodular eigenvalues $\varepsilon_\pi$, 
$\lambda_\pi$, and $\mu_\pi$ for $\pi$. 
See Theorem \ref{charpolyiwahoritheorem}.

\clearpage

\begin{table}
\caption{Non-supercuspidal representations of $\GSp(4,F)$}
\label{maintable}

\noindent\resizebox{\textwidth}{!}{
$
\renewcommand{\arraystretch}{1.37}
\begin{array}{ccccccc}
\toprule
\mbox{constituent of}&\multicolumn{2}{c}{\mbox{group}}&\mbox{representation}
   &{\rm tempered}& \mbox{ess.~$L^2$}  &\mbox{generic}\\
\bottomrule
\chi_1\times\chi_2\rtimes\sigma\ \text{(irreducible)}  &\mathrm{I}&&\chi_1\times\chi_2\rtimes\sigma&\mbox{$\chi_i,\sigma$ unit.}
   &&\bullet\\
\midrule
  \nu^{1/2}\chi\times\nu^{-1/2}\chi\rtimes\sigma&\mathrm{II}&\mbox{a}
   &\chi\St_{\GL(2)}\rtimes\sigma
   &\mbox{$\chi,\sigma$ unit.}&&\bullet\\
\cmidrule{3-7}
   (\chi^2\neq\nu^{\pm1},\chi\neq\nu^{\pm3/2})&&\mbox{b}
   &\chi\triv_{\GL(2)}\rtimes\sigma
   &&&\\
\midrule
  \chi\times\nu\rtimes\nu^{-1/2}\sigma&\mathrm{III}&\mbox{a}&\chi\rtimes\sigma\St_{\GSp(2)}&\mbox{$\chi,\sigma$ unit.}&
   &\bullet\\
\cmidrule{3-7}
  (\chi\notin\{1,\nu^{\pm2}\})&&\mbox{b}&\chi\rtimes\sigma\triv_{\GSp(2)}&&&\\
\midrule
 \nu^2\times\nu\rtimes\nu^{-3/2}\sigma&\mathrm{IV}&\mbox{a}&\sigma\St_{\GSp(4)}&\mbox{$\sigma$ unit.}&\bullet&\bullet\\
\cmidrule{3-7}
  &&\mbox{b}&L(\nu^2,\nu^{-1}\sigma\St_{\GSp(2)})&&&\\
\cmidrule{3-7}
  &&\mbox{c}&L(\nu^{3/2}\St_{\GL(2)},\nu^{-3/2}\sigma)&&&\\
\cmidrule{3-7}
  &&\mbox{d}&\sigma\triv_{\GSp(4)}&&&\\
\midrule
  \nu\xi \times\xi \rtimes\nu^{-1/2}\sigma&\mathrm{V}&\mbox{a}&\delta([\xi ,\nu\xi ],\nu^{-1/2}\sigma)
   &\mbox{$\sigma$ unit.}&\bullet&\bullet\\
\cmidrule{3-7}
  (\xi ^2=1,\:\xi \neq1)&&\mbox{b}&L(\nu^{1/2}\xi \St_{\GL(2)},\nu^{-1/2}\sigma)&&&\\
\cmidrule{3-7}
  &&\mbox{c}&L(\nu^{1/2}\xi \St_{\GL(2)},\xi \nu^{-1/2}\sigma)&&&\\
\cmidrule{3-7}
  &&\mbox{d}&L(\nu\xi ,\xi \rtimes\nu^{-1/2}\sigma)&&&\\
\midrule
  \nu\times1_{F^\times}\rtimes\nu^{-1/2}\sigma&\mathrm{VI}&\mbox{a}&\tau(S,\nu^{-1/2}\sigma)&\mbox{$\sigma$ unit.}&&\bullet\\
\cmidrule{3-7}
  &&\mbox{b}&\tau(T,\nu^{-1/2}\sigma)&\mbox{$\sigma$ unit.}&&\\
\cmidrule{3-7}
  &&\mbox{c}&L(\nu^{1/2}\St_{\GL(2)},\nu^{-1/2}\sigma)&&&\\
\cmidrule{3-7}
  &&\mbox{d}&L(\nu,1_{F^\times}\rtimes\nu^{-1/2}\sigma)&&&\\
\midrule
\chi\rtimes\pi\ \text{(irreducible)}  &{\rm VII}&&\chi\rtimes \pi&\mbox{$\chi,\pi$ unit.}
   &&\bullet\\
\midrule
1_{F^\times}\rtimes\pi&\mathrm{VIII}&\mbox{a}&\tau(S,\pi)&\mbox{$\pi$ unit.}&&\bullet\\
\cmidrule{3-7}
&&\mbox{b}&\tau(T,\pi)&\mbox{$\pi$ unit.}&&\\
\midrule
\nu\xi \rtimes\nu^{-1/2}\pi &\mbox{IX}&\mbox{a}
   &\delta(\nu\xi ,\nu^{-1/2}\pi)&\mbox{$\pi$ unit.}&\bullet&\bullet\\
\cmidrule{3-7}
(\xi \neq1,\:\xi \pi=\pi)&&\mbox{b}&L(\nu\xi ,\nu^{-1/2}\pi)&&&\\
\midrule
\pi\rtimes \sigma\   \text{(irreducible)} &{\rm X}&&\pi\rtimes\sigma&\mbox{$\pi,\sigma$ unit.}
   &&\bullet\\
\midrule
\nu^{1/2}\pi\rtimes\nu^{-1/2}\sigma&\mbox{XI}&\mbox{a}&\delta(\nu^{1/2}\pi,\nu^{-1/2}\sigma)
&\mbox{$\pi,\sigma$ unit.}&\bullet&\bullet\\
\cmidrule{3-7}
(\omega_\pi=1)&&\mbox{b}&L(\nu^{1/2}\pi,\nu^{-1/2}\sigma)&&&\\
\bottomrule
 \end{array}
$}
\end{table}

\clearpage

\begin{table}
\caption{Non-paramodular representations}
\label{noklingentable}
\renewcommand{\arraystretch}{1.5}
\begin{tabular}{cccc}
\toprule
&&\mbox{\rm representation}&\mbox{\rm condition on defining data} \\ 
\bottomrule
 II&b&$\chi\triv_{\GL(2)}\rtimes\sigma$
   & $\chi \sigma$ ramified
   \\
\midrule
 III&b
   &$\chi\rtimes\sigma\triv_{\GSp(2)}$
   & $\sigma$ ramified \\ 
\midrule
&b&$L(\nu^2,\nu^{-1}\sigma\St_{\GSp(2)})$& $\sigma$ ramified
   \\ 
   \cmidrule{3-4}
IV &c&$L(\nu^{3/2}\St_{\GL(2)},\nu^{-3/2}\sigma)$
   &$\sigma$ ramified \\ 
   \cmidrule{3-4}
  &d&$\sigma\triv_{\GSp(4)}$& $\sigma$ ramified\\
\midrule
&b
   &$L(\nu^{1/2}\xi\St_{\GL(2)},\nu^{-1/2}\sigma)$&$\sigma$ ramified\\ 
   \cmidrule{3-4}
V &c
   &$L(\nu^{1/2}\xi\St_{\GL(2)},\xi\nu^{-1/2}\sigma)$& $\xi \sigma$ ramified
   \\ 
   \cmidrule{3-4}
 &d
   &$L(\nu\xi,\xi\rtimes\nu^{-1/2}\sigma)$&$\sigma$  or $\xi$ ramified \\ 
   \midrule
  &b&$\tau(T,\nu^{-1/2}\sigma)$&none\\
  \cmidrule{3-4}
  VI &c
   &$L(\nu^{1/2}\St_{\GL(2)},\nu^{-1/2}\sigma)$&$\sigma$ ramified\\
   \cmidrule{3-4}
  &d&$L(\nu,1_{F^\times}\rtimes\nu^{-1/2}\sigma)$
   &$\sigma$  ramified\\ 
   \midrule
   VIII&b&$\tau(T,\pi)$& none\\ 
\midrule
  IX
   &b&$L(\nu\xi,\nu^{-1/2}\pi)$&none\\ 
   \midrule
  XI &b
   &$L(\nu^{1/2}\pi,\nu^{-1/2}\sigma)$&$\sigma$ ramified\\ 
\midrule
 &&$\pi$ supercuspidal
   &  non-generic\\ 
\bottomrule
\end{tabular}
\end{table}

\clearpage

\begin{table}
\caption{Stable Klingen dimensions}
\label{dimensionstable}
\noindent\resizebox{\textwidth}{!}{
$
\renewcommand{\arraystretch}{1.3}\renewcommand{\arraycolsep}{1ex}
\begin{array}{clcccccccc}
\toprule
&&\text{inducing data}&N_\pi&N_{\pi,s}&\dim V_s(n)&\bar N_{\pi,s}&\dim \bar V_s(n)&\text{Cat.}&\text{Cmt.}\\
\bottomrule
\mathrm{I}&&\chi_1,\chi_2,\sigma\text{ unr.}&0&0&\frac{n^2+5n+2}2&1&n&\text{\cattwobox}&\\
\cmidrule{3-10}
&&\chi_1,\chi_2\text{ ram.},\sigma\text{ unr}.&\cellcolor{white}a:=a(\chi_1)+a(\chi_2)&a&\frac{(n-a+1)(n-a+4)}2&a&n-a+1&\text{\cattwobox}&\\
\cmidrule{3-10}
&&\chi_i\sigma\text{ unr}.,\sigma\text{ ram}.&\cellcolor{white}a:=2a(\sigma)&a&\frac{(n-a+1)(n-a+4)}2&a&n-a+1&\text{\cattwobox}&\\
\cmidrule{3-10}
&&\chi_i\sigma\text{ ram}.,\sigma\text{ ram}.&\cellcolor{white}\scriptstyle a:=a(\chi_1\sigma)+a(\chi_2\sigma)+2a(\sigma)&a-1&\frac{(n-a+2)(n-a+3)}2
&a-1&n-a+2&\text{\catonebox}&\\
\midrule
\mathrm{II}&\mathrm{a}&\sigma,\chi\text{ unr.}&1&1&\frac{n(n+3)}2&1&n&\text{\cattwobox}&\\
\cmidrule{3-10}
&&\sigma\text{ ram}.,\chi\sigma\text{ unr.}&\cellcolor{white}a:=2a(\sigma)+1&a-1&\frac{(n-a+2)(n-a+3)}2&a-1&n-a+2&\text{\catonebox}&\\
\cmidrule{3-10}
&&\sigma\text{ unr}.,\chi\sigma\text{ ram.}&\cellcolor{white}a:=2a(\sigma\chi)&a&\frac{(n-a+1)(n-a+4)}2&a&n-a+1&\text{\cattwobox}&\\
\cmidrule{3-10}
&&\sigma,\chi\sigma\text{ ram.}&\cellcolor{white}a:=2a(\chi\sigma)+2a(\sigma)&a-1&\frac{(n-a+2)(n-a+3)}2&a-1&n-a+2&\text{\catonebox}&\\
\cmidrule{2-10}
&\text{b}&\chi\sigma\text{ unr.}, \sigma\text{ unr.}&0&0&n+1&\text{---}&0&\text{\cattwobox}&\text{SK}\\
\cmidrule{3-10}
&&\chi\sigma\text{ unr.}, \sigma\text{ ram.}&\cellcolor{white}a:=2a(\sigma)&a&n-a+1&\text{---}&0&\text{\cattwobox}&\text{SK}\\
\cmidrule{3-10}
&&\chi\sigma\text{ ram.}&\multicolumn{5}{c}{\text{not paramodular}}&&\text{SK}\\
\midrule
\mathrm{III}&\mathrm{a}&\sigma\text{ unr.}&\cellcolor{white}2&1&\frac{n(n+1)}2&1&n&\text{\catonebox}&\\
\cmidrule{3-10}
&&\sigma\text{ ram.}&\cellcolor{white}a:=4a(\sigma)&a-1&\frac{(n-a+2)(n-a+3)}2&a-1&n-a+2&\text{\catonebox}& \\
\cmidrule{2-10}
&\mathrm{b}&\sigma\;\text{unr.}&0&0&2n+1&\text{---}&0&\text{\cattwobox}&\\
\cmidrule{3-10}
&&\sigma\;\text{ram.}&\multicolumn{5}{c}{\text{not paramodular}}&&\\
\midrule
\mathrm{IV}&\mathrm{a}&\sigma\;\text{unr.}&\cellcolor{white}3&2&\frac{(n-1)n}2&2&n-1&\text{\catonebox}&\\
\cmidrule{3-10}
&&\sigma\;\text{ram.}&\cellcolor{white}a:=4a(\sigma)&a-1&\frac{(n-a+2)(n-a+3)}2&a-1&n-a+2&\text{\catonebox}&\\
\cmidrule{2-10}
&\mathrm{b}&\sigma\;\text{unr.}&\cellcolor{white}2&1&n&1&1&\text{\catonebox}&\text{non-unit.}\\
\cmidrule{3-10}
&&\sigma\;\text{ram.}&\multicolumn{5}{c}{\text{not paramodular}}&&\text{non-unit.}\\
\cmidrule{2-10}
&\mathrm{c}&\sigma\;\text{unr.}&1&1&2n&1&1&\text{\cattwobox}&\text{non-unit.}\\
\cmidrule{3-10}
&&\sigma\;\text{ram.}&\multicolumn{5}{c}{\text{not paramodular}}&&\text{non-unit.}\\
\cmidrule{2-10}
&\mathrm{d}&\sigma\;\text{unr.}&0&0&1&\text{---}&0&\text{\cattwobox}&\text{one-dim.}\\
\cmidrule{3-10}
&&\sigma\;\text{ram.}&\multicolumn{5}{c}{\text{not paramodular}}&&\text{one-dim.}\\
\midrule
\mathrm{V}&\mathrm{a}&\text{$\sigma$, $\xi$ unr.}&\cellcolor{white}2&1&\frac{n(n+1)}2&1&n&\text{\catonebox}&\\
\cmidrule{3-10}
&&\text{$\sigma$ unr., $\xi$ ram.}&\cellcolor{white}a:=2a(\xi)+1&a-1&\frac{(n-a+2)(n-a+3)}2&a-1&n-a+2&\text{\catonebox}&\\
\cmidrule{3-10}
&&\text{$\sigma$ ram., $\sigma\xi$ unr.}&\cellcolor{white}a:=2a(\sigma)+1&a-1&\frac{(n-a+2)(n-a+3)}2&a-1&n-a+2&\text{\catonebox}&\\
\cmidrule{3-10}
&&\text{$\sigma$, $\sigma\xi$ ram.}&\cellcolor{white}a:=2a(\xi\sigma)+2a(\sigma)&a-1&\frac{(n-a+2)(n-a+3)}2&a-1&n-a+2&\text{\catonebox}&\\
\cmidrule{2-10}
&\mathrm{b}&\text{$\sigma$, $\xi$ unr.}&1&1&n&\text{---}&0&\text{\cattwobox}&\text{SK}\\
\cmidrule{3-10}
&&\text{$\sigma$ unr., $\xi$ ram.}&\cellcolor{white}2a(\xi)&2a(\xi)&n-2a(\xi)+1&\text{---}&0&\text{\cattwobox}&\text{SK}\\
\cmidrule{3-10}
&&\text{$\sigma$ ram., $\sigma\xi$ unr.}&\multicolumn{5}{c}{\text{not paramodular}}&&\text{SK}\\
\cmidrule{3-10}
&&\text{$\sigma$, $\sigma\xi$ ram.}&\multicolumn{5}{c}{\text{not paramodular}}&&\text{SK}\\
\bottomrule
\end{array}
$
}
\end{table}
 \clearpage

\noindent\resizebox{\textwidth}{!}{
$
\renewcommand{\arraystretch}{1.3}
\begin{array}{cccccccccc}
\toprule
&&\text{inducing data}&N_\pi&N_{\pi,s}&\dim V_s(n)&\bar N_{\pi,s}&\dim \bar V_s(n)&\text{Cat.}&\text{Cmt.}\\
\bottomrule
\mathrm{V}&\mathrm{c}&\text{$\sigma$, $\xi$ unr.}&1&1&n&\text{---}&0&\text{\cattwobox}&\text{SK}\\
\cmidrule{3-10}
&&\text{$\sigma$ unr., $\xi$ ram.}&\multicolumn{5}{c}{\text{not paramodular}}&&\text{SK}\\
\cmidrule{3-10}
&&\text{$\sigma$ ram., $\sigma\xi$ unr.}&\cellcolor{white}2a(\sigma)&2a(\sigma)&n-2a(\sigma)+1&\text{---}&0&\text{\cattwobox}&\text{SK}\\
\cmidrule{3-10}
&&\text{$\sigma$, $\sigma\xi$ ram.}&\multicolumn{5}{c}{\text{not paramodular}}&&\text{SK}\\
\cmidrule{2-10}
&\mathrm{d}&\sigma,\xi\;\text{unr.}&0&0&1&\text{---}&0&\text{\cattwobox}&\\
\cmidrule{3-10}
&&\text{$\sigma$ or $\xi$ ram.}&\multicolumn{5}{c}{\text{not paramodular}}&&\\
\midrule
\mathrm{VI}&\mathrm{a}&\text{$\sigma$ unr.}&\cellcolor{white}2&1&\frac{n(n+1)}2&1&n&\text{\catonebox}&\\
\cmidrule{3-10}
&&\text{$\sigma$ ram.}&\cellcolor{white}a:=4a(\sigma)&a-1&\frac{(n-a+2)(n-a+3)}2&a-1&n-a+2&\text{\catonebox}&\\
\cmidrule{2-10}
&\mathrm{b}&\text{$\sigma$ unr.}&\multicolumn{5}{c}{\text{not paramodular}}&&\\
\cmidrule{3-10}
&&\text{$\sigma$ ram.}&\multicolumn{5}{c}{\text{not paramodular}}&&\\
\cmidrule{2-10}
&\mathrm{c}&\sigma\;\text{unr.}&1&1&n&\text{---}&0&\text{\cattwobox}&\text{SK}\\
\cmidrule{3-10}
&&\sigma\;\text{ram.}&\multicolumn{5}{c}{\text{not paramodular}}&&\text{SK}\\
\cmidrule{2-10}
&\mathrm{d}&\sigma\;\text{unr.}&0&0&n+1&\text{---}&0&\text{\cattwobox}&\\
\cmidrule{3-10}
&&\sigma\;\text{ram.}&\multicolumn{5}{c}{\text{not paramodular}}&&\\
\midrule
\mathrm{VII}&&&\cellcolor{white}a:=2a(\pi)&a-1&\frac{(n-a+2)(n-a+3)}2&a-1&n-a+2&\text{\catonebox}&\\
\midrule
\mathrm{VIII}&\mathrm{a}&&\cellcolor{white}a:=2a(\pi)&a-1&\frac{(n-a+2)(n-a+3)}2&a-1&n-a+2&\text{\catonebox}&\\
\cmidrule{2-10}
&\mathrm{b}&&\multicolumn{5}{c}{\text{not paramodular}}&&\\
\midrule
\mathrm{IX}&\mathrm{a}&&\cellcolor{white}a:=2a(\pi)&a-1&\frac{(n-a+2)(n-a+3)}2&a-1&n-a+2&\text{\catonebox}\\
\cmidrule{2-10}
&\mathrm{b}&&\multicolumn{5}{c}{\text{not paramodular}}&&\\
\midrule
\mathrm{X}&&\sigma\text{ unr}.&\cellcolor{white}a:=a(\pi)&a&\frac{(n-a+1)(n-a+4)}2&a&n-a+1&\text{\cattwobox}&\\
\cmidrule{3-10}
&&\sigma\text{ ram}.&\cellcolor{white}a:=a(\sigma\pi)+2a(\sigma)&a-1&\frac{(n-a+2)(n-a+3)}2&a-1&n-a+2&\text{\catonebox}&\\
\midrule
\mathrm{XI}&\mathrm{a}&\sigma\;\text{unr.}&\cellcolor{white}a:=a(\sigma\pi)+1&a-1&\frac{(n-a+2)(n-a+3)}2&a-1&n-a+2&\text{\catonebox}&\\
\cmidrule{3-10}
&&\sigma\text{ ram.}&\cellcolor{white}a:=a(\sigma\pi)+2a(\sigma)&a-1&\frac{(n-a+2)(n-a+3)}2&a-1&n-a+2&\text{\catonebox}&\\
\cmidrule{2-10}
&\mathrm{b}&\sigma\text{ unr.}&\cellcolor{white}a:=a(\pi)&a&n-a+1&\text{---}&0&\text{\cattwobox}&\text{SK}\\
\cmidrule{3-10}
&&\sigma\text{ ram.}\quad&\multicolumn{5}{c}{\text{not paramodular}}&&\text{SK}\\
\midrule
\text{s.c.}&&\text{generic}&\cellcolor{white}a\geq4&a-1&\frac{(n-a+2)(n-a+3)}2&a-1&n-a+2&\text{\catonebox}&\\
\cmidrule{3-10}
&&\text{non-generic}&\multicolumn{5}{c}{\text{not paramodular}}&& \\
\bottomrule
\end{array}
$
}

\clearpage

\begin{table}
\caption{Hecke eigenvalues}
\label{levelsandeigenvaluestable}

\noindent\resizebox{\textwidth}{!}{
$
\renewcommand{\arraystretch}{1.4}\renewcommand{\arraycolsep}{0.7ex}
\begin{array}{ccccccccc}
\toprule
&&\text{inducing data}&N_{\pi}&\varepsilon_\pi&\lambda_\pi&\mu_\pi&\text{Cat.}&\text{Cmt.}\\
\bottomrule
\mathrm{I}&&\chi_1,\chi_2,\sigma\text{ unr.}&0&1&\text{(A1)}&\text{(A2)}&\text{\cattwobox}& \\
\cmidrule{3-9}
&&\chi_1,\chi_2\text{ ram.},\sigma\text{ unr}.&\cellcolor{white}a(\chi_1)+a(\chi_2)&\chi_1(-1)&q^{3/2}(\sigma(\varpi)+\sigma(\varpi^{-1}))&0
&\text{\cattwobox}&\\
\cmidrule{3-9}
&&\chi_i\sigma\text{ unr}.,\sigma\text{ ram}.&\cellcolor{white}2a(\sigma)&\chi_1(-1)&q^{3/2}((\chi_1\sigma)(\varpi)+(\chi_2\sigma)(\varpi))&0
&\text{\cattwobox}&\\
\cmidrule{3-9}
&&\chi_i\sigma\text{ ram}.,\sigma\text{ ram}.&\cellcolor{white}\scriptstyle a(\chi_1\sigma)+a(\chi_2\sigma)+2a(\sigma)&\chi_1(-1)&0&-q^2
&\text{\catonebox}&\\
\midrule
\mathrm{II}&\mathrm{a}&\sigma,\chi\text{ unr.}&1&-(\chi\sigma)(\varpi)&\text{(B1)}&\text{(B2)}&\text{\cattwobox}&\\
\cmidrule{3-9}
&&\sigma\text{ ram}.,\chi\sigma\text{ unr.}&\cellcolor{white}2a(\sigma)+1&\scriptstyle-\sigma(-1)(\chi\sigma)(\varpi)&q(\chi\sigma)(\varpi)&-q^2
&\text{\catonebox}&\\
\cmidrule{3-9}
&&\sigma\text{ unr}.,\chi\sigma\text{ ram.}&\cellcolor{white}2a(\sigma\chi)&\chi(-1)&q^{3/2}(\sigma(\varpi)+\sigma(\varpi)^{-1})&0&\text{\cattwobox}\\
\cmidrule{3-9}
&&\sigma,\chi\sigma\text{ ram.}&\cellcolor{white}2a(\chi\sigma)+2a(\sigma)&\chi(-1)&0&-q^2&\text{\catonebox}&\\
\cmidrule{2-9}
&\mathrm{b}&\chi\sigma\text{ unr.}, \sigma\text{ unr.}&0&1 &\text{(C1)}&\text{(C2)}&\text{\cattwobox} &\text{SK}\\
\cmidrule{3-9}
&&\chi\sigma\text{ unr.}, \sigma\text{ ram.}&\cellcolor{white}2a(\sigma)&\chi(-1)&q(q+1)(\sigma\chi)(\varpi)&0&\text{\cattwobox}&\text{SK}\\
\cmidrule{3-9}
&&\chi\sigma\text{ ram.}&\multicolumn{4}{c}{\text{not paramodular}}& &\text{SK}\\
\midrule
\mathrm{III}&\mathrm{a}&\sigma\text{ unr.}&\cellcolor{white}2&1&q(\sigma(\varpi)+\sigma(\varpi)^{-1})&-q^2+q&\text{\catonebox}&\\
\cmidrule{3-9}
&&\sigma\text{ ram.}&\cellcolor{white}4a(\sigma)&1&0&-q^2&\text{\catonebox}&\\
\cmidrule{3-9}
&\mathrm{b}&\sigma\;\text{unr.}&0&1&q(q+1)\sigma(\varpi)(1+\chi(\varpi))&\text{(D)}&\text{\cattwobox}&\\
\cmidrule{3-9}
&&\sigma\;\text{ram.}&\multicolumn{4}{c}{\text{not paramodular}}&& \\
\midrule
\mathrm{IV}&\mathrm{a}&\sigma\;\text{unr.}&\cellcolor{white}3&-\sigma(\varpi)&\sigma(\varpi)&-q^2&\text{\catonebox}&\\
\cmidrule{3-9}
&&\sigma\;\text{ram.}&\cellcolor{white}4a(\sigma)&1&0&-q^2&\text{\catonebox}&\\
\cmidrule{2-9}
&\mathrm{b}&\sigma\;\text{unr.}&\cellcolor{white}2&1&\sigma(\varpi)(1+q^2)&-q^2+q&\text{\catonebox}&\text{non-unit.}\\
\cmidrule{3-9}
&&\sigma\;\text{ram.}&\multicolumn{4}{c}{\text{not paramodular}}& &\text{non-unit.}\\
\cmidrule{2-9}
&\mathrm{c}&\sigma\;\text{unr.}&1&-\sigma(\varpi)&\sigma(\varpi)(q^3+q+2)&q^3+1&\text{\cattwobox}&\text{non-unit.}\\
\cmidrule{3-9}
&&\sigma\;\text{ram.}&\multicolumn{4}{c}{\text{not paramodular}}&&\text{non-unit.}\\
\cmidrule{2-9}
&\mathrm{d}&\sigma\;\text{unr.}&0&1&\sigma(\varpi)(q+1)(q^2+1)&q(q+1)(q^2+1)&\text{\cattwobox}&\text{one-dim.}\\
\cmidrule{3-9}
&&\sigma\;\text{ram.}&\multicolumn{4}{c}{\text{not paramodular}}&&\text{one-dim.}\\
\midrule
\mathrm{V}&\mathrm{a}&\text{$\sigma$, $\xi$ unr.}&\cellcolor{white}2&-1&0&-q^2-q&\text{\catonebox}&\\
\cmidrule{3-9}
&&\text{$\sigma$ unr., $\xi$ ram.}&\cellcolor{white}2a(\xi)+1&-\sigma(\varpi)\xi(-1)&\sigma(\varpi)q&-q^2&\text{\catonebox}&\\
\cmidrule{3-9}
&&\text{$\sigma$ ram., $\sigma\xi$ unr.}&\cellcolor{white}2a(\sigma)+1&\scriptstyle-\sigma(-1)(\xi\sigma)(\varpi)&(\xi\sigma)(\varpi)q&-q^2&
\text{\catonebox}&\\
\cmidrule{3-9}
&&\text{$\sigma$, $\sigma\xi$ ram.}&\cellcolor{white}2a(\xi\sigma)+2a(\sigma)&\xi(-1)&0&-q^2&\text{\catonebox}&\\
\cmidrule{2-9}
&\mathrm{b}&\text{$\sigma$, $\xi$ unr.}&1&\sigma(\varpi)&\sigma(\varpi)(q^2-1)&-q^2-q&\text{\cattwobox}&\text{SK}\\
\cmidrule{3-9}
&&\text{$\sigma$ unr., $\xi$ ram.}&\cellcolor{white}2a(\xi)&\xi(-1)&\sigma(\varpi)(q^2+q)&0&\text{\cattwobox}&\text{SK} \\
\cmidrule{3-9}
  &&\text{$\sigma$ ram., $\sigma\xi$ unr.}&\multicolumn{4}{c}{\text{not paramodular}}&&\text{SK}\\
\cmidrule{3-9}
    &&\text{$\sigma$, $\sigma\xi$ ram.}&\multicolumn{4}{c}{\text{not paramodular}}&&\text{SK}\\
\bottomrule
 \end{array}
$
}
\end{table}

\clearpage

\noindent\resizebox{\textwidth}{!}{
$
\renewcommand{\arraystretch}{1.4}\renewcommand{\arraycolsep}{0.7ex}
\begin{array}{ccccccccc}
\toprule
&&\text{inducing data}&N_{\pi}&\varepsilon_\pi&\lambda_\pi&\mu_\pi&\text{Cat.}&\text{Cmt.}\\
\bottomrule
\mathrm{V}&\mathrm{c}&\text{$\sigma$, $\xi$ unr.}&1&-\sigma(\varpi)&-\sigma(\varpi)(q^2-1)&-q^2-q&\text{\cattwobox}&\text{SK}\\
\cmidrule{3-9}
&&\text{$\sigma$ unr., $\xi$ ram.}&\multicolumn{4}{c}{\text{not paramodular}}&&\text{SK}\\
\cmidrule{3-9}
&&\text{$\sigma$ ram., $\sigma\xi$ unr.}&\cellcolor{white}2a(\sigma)&\xi(-1)&(\xi\sigma)(\varpi)(q^2+q)&0&\text{\cattwobox}&\text{SK}\\
\cmidrule{3-9}
&&\text{$\sigma$, $\sigma\xi$ ram.}&\multicolumn{4}{c}{\text{not paramodular}}&&\text{SK}\\
\cmidrule{2-9}
&\mathrm{d}&\sigma,\xi\;\text{unr.}&0&1&0&-(q+1)(q^2+1)&\text{\cattwobox}&\\
\cmidrule{3-9}
    &&\text{$\sigma$ or $\xi$ ram.}&\multicolumn{4}{c}{\text{not paramodular}}&&\\
\midrule
\mathrm{VI}&\mathrm{a}&\text{$\sigma$ unr.}&\cellcolor{white}2&1&2q\sigma(\varpi)&-q^2+q&\text{\catonebox}&\\
\cmidrule{3-9}
&&\text{$\sigma$ ram.}&\cellcolor{white}4a(\sigma)&1&0&-q^2&\text{\catonebox}&\\
\cmidrule{2-9}
&\mathrm{b}&\text{$\sigma$ unr.}&\multicolumn{4}{c}{\text{not paramodular}}&&\\
\cmidrule{3-9}
&&\text{$\sigma$ ram.}&\multicolumn{4}{c}{\text{not paramodular}}&&\\
\cmidrule{2-9}
&\mathrm{c}&\sigma\;\text{unr.}&1&-\sigma(\varpi)&\sigma(\varpi)(q+1)^2&q(q+1)&\text{\cattwobox}&\text{SK}\\
\cmidrule{3-9}
&&\sigma\;\text{ram.}&\multicolumn{4}{c}{\text{not paramodular}}& &\text{SK}\\
\cmidrule{2-9}
&\mathrm{d}&\sigma\;\text{unr.}&0&1&2q(q+1)\sigma(\varpi)&(q+1)(q^2+2q-1)&\text{\cattwobox}& \\
\cmidrule{3-9}
&&\sigma\;\text{ram.}&\multicolumn{4}{c}{\text{not paramodular}}& \\
\midrule
\mathrm{VII}&&&\cellcolor{white}2a(\pi)&\chi(-1)&0&-q^2&\text{\catonebox}&\\
\midrule
\mathrm{VIII}&\mathrm{a}&&\cellcolor{white}2a(\pi)&1&0&-q^2&\text{\catonebox}\\
\cmidrule{3-9}
&\mathrm{b}&&\multicolumn{4}{c}{\text{not paramodular}}& \\
\midrule
\mathrm{IX}&\mathrm{a}&&\cellcolor{white}2a(\pi)&\xi(-1)&0&-q^2&\text{\catonebox}&\\
\cmidrule{3-9}
&\mathrm{b}&&\multicolumn{4}{c}{\text{not paramodular}}&& \\
\midrule
\mathrm{X}&&\sigma\text{ unr}.&\cellcolor{white}a(\pi)&\varepsilon(1/2,\sigma\pi)&q^{3/2}(\sigma(\varpi)+\sigma(\varpi)^{-1})&0&\text{\cattwobox}&\\
\cmidrule{3-9}
&&\sigma\text{ ram}.&\cellcolor{white}a(\sigma\pi)+2a(\sigma)&\scriptstyle\sigma(-1)\varepsilon(1/2,\sigma\pi)&0&-q^2&\text{\catonebox}&\\
\midrule
\mathrm{XI}&\mathrm{a}&\sigma\;\text{unr.}&\cellcolor{white}a(\sigma\pi)+1&\scriptstyle-\sigma(\varpi)\varepsilon(1/2,\sigma\pi)&q\sigma(\varpi)&-q^2
&\text{\catonebox}&\\
\cmidrule{3-9}
&&\sigma\text{ ram.}&\cellcolor{white}a(\sigma\pi)+2a(\sigma)&\scriptstyle\sigma(-1)\varepsilon(1/2,\sigma\pi)&0&-q^2&\text{\catonebox}&\\
\cmidrule{2-9}
&\mathrm{b}&\sigma\text{ unr.}&\cellcolor{white}a(\pi)&\varepsilon(1/2,\sigma\pi)&(q^2+q)\sigma(\varpi)&0&\text{\cattwobox}&\text{SK}\\
\cmidrule{3-9}
&&\sigma\text{ ram.}&\multicolumn{4}{c}{\text{not paramodular}}&&\text{SK}\\
\midrule
  \text{s.c.}&&\text{generic}&\cellcolor{white}N_\pi\geq 4&\varepsilon_\pi&0&-q^2&\text{\catonebox}&\\
\cmidrule{3-9}
  &&\text{non-generic}&\multicolumn{4}{c}{\text{not paramodular}}&&\\
\bottomrule
 \end{array}
$
}

\clearpage

\begin{table}
\caption{Characteristic polynomials of $T_{0,1}^s$ and $T_{1,0}^s$
on $V_s(1)$  in terms
of inducing data}
\label{Iwahoritable}
\noindent\resizebox{\textwidth}{!}{
$
\renewcommand{\arraystretch}{1.6}
 \begin{array}{lccll}
\toprule
\multicolumn{2}{c}{\mathrm{type}}&\pi&\multicolumn{1}{c}{p(T^s_{0,1},V_s(1),X)}&\multicolumn{1}{c}{p(T^s_{1,0},V_s(1),X)}\\
\bottomrule
\mathrm{I}&&\chi_1\times\chi_2\rtimes \sigma&\left(X-\chi_1(\varpi)\left(1+\chi_2(\varpi))\sigma(\varpi\right) q^{\frac{3}{2}}\right)& \left( X - \chi_1(\varpi) q^2\right)\\
\chi_1,\chi_2,\sigma\ \mathrm{unram.}&&&\times\left(X-\chi_2(\varpi)\left(1+\chi_1(\varpi))\sigma(\varpi\right) q^{\frac{3}{2}}\right)&\times \left( X - \chi_2(\varpi) q^2 \right)\\
\chi_1\chi_2\sigma^2=1&&&\times\left(X-\left(1+\chi_1(\varpi)\right)\sigma(\varpi) q^{\frac{3}{2}}\right)&\times\left( X - \chi_2(\varpi)^{-1} q^2 \right)\\
&&&\times\left(X-\left(1+\chi_2(\varpi)\right)\sigma(\varpi)q^{\frac{3}{2}}\right)&\times\left( X- \chi_1(\varpi)^{-1} q^2 \right)\\
\midrule
 \mathrm{II}&\mathrm{a}&\chi\St_{\GL(2)}\rtimes\sigma&X^2
&X^2-\left(\chi(\varpi)+\chi(\varpi)^{-1}\right) q^{\frac{3}{2}}X+q^3\\
\chi,\sigma\ \mathrm{unram.}&&&-\Big(2(\chi\sigma)(\varpi) q+\left(\sigma(\varpi)+\sigma(\varpi)^{-1}\right)q^{\frac{3}{2}}\Big)X
&\\
\chi^2\sigma^2=1&&& +q^2+q^3+\left(\chi(\varpi)+\chi(\varpi)^{-1}\right)q^{\frac{5}{2}}&\\ 
\cmidrule{2-5}
 &\mathrm{b}&\chi 1_{\GL(2)}\rtimes\sigma&X^2&X^2\\
&&&- \left(2 (\chi\sigma)(\varpi) q^2+\left(\sigma(\varpi)+\sigma(\varpi)^{-1}\right) q^{\frac{3}{2}}\right)X&- \left(\chi(\varpi)+\chi(\varpi)^{-1}\right) q^{5/2} X+q^5\\
&&&+q^3+q^4+\left(\chi(\varpi)+\chi(\varpi)^{-1}\right) q^{\frac{7}{2}}&\\
\midrule
 \mathrm{III}&\mathrm{a}&\chi \rtimes \sigma\St_{\GSp(2)}&X-q(\sigma(\varpi)+\sigma(\varpi)^{-1})&X-q\\
\chi,\sigma\ \mathrm{unram.}&&&&\\
\chi\sigma^2=1&&&&\\
\cmidrule{2-5}
&\mathrm{b}&\chi\rtimes \sigma 1_{\GSp(2)} &X^3&X^3\\
&&&-\left(\sigma(\varpi)+\sigma(\varpi)^{-1}\right)q(1+2q)X^2&
-\left(\left(\chi(\varpi)+\chi(\varpi)^{-1}\right)q^2+q^3\right)X^2\\
&&&+\left(1+3q+(\chi(\varpi)+\chi(\varpi)^{-1})q\right)q^2(q+1)X&
+\left(\left(\chi(\varpi)+\chi(\varpi)^{-1}\right)q^5+q^4\right)X\\
&&&-\left(\sigma(\varpi)+\sigma(\varpi)^{-1}\right)q^4(q+1)^2&-q^7\\
\midrule
  \mathrm{IV}&\mathrm{b}&L(\nu^2,\nu^{-1}\sigma \St_{\GSp(2)})&X-\sigma(\varpi)(1+q^2)&X-q\\
\cmidrule{2-5}
\sigma\ \mathrm{unram.}&\mathrm{c}&L(\nu^{\frac{3}{2}}\St_{\GL(2)},\nu^{-\frac{3}{2}}\sigma)&X^2-\sigma(\varpi)(1+2q+q^3)X&X^2-(1+q^3)X+q^3\\
\sigma^2=1&&&+q(1+q)(1+q^2)&\\
\cmidrule{2-5}
 &\mathrm{d}&\sigma1_{\GSp(4)}&X-\sigma(\varpi)(q^2+q^3)&X-q^4\\
\midrule
 \mathrm{V}&\mathrm{a}&\delta([\xi,\nu \xi],\nu^{-\frac{1}{2}}\sigma)&X&X+q\\
 \cmidrule{2-5}
\sigma,\xi\ \mathrm{unram.}&\mathrm{b} &L(\nu^{\frac{1}{2}}\xi\St_{\GL(2)},\nu^{-\frac{1}{2}}\sigma)&X-\sigma(\varpi)(q^2-q)&X+q^2\\
\cmidrule{2-5}
\sigma^2=\xi^2=1&\mathrm{c} &L(\nu^{\frac{1}{2}}\xi\St_{\GL(2)},\xi\nu^{-\frac{1}{2}}\sigma)&X+\sigma(\varpi)(q^2-q)&X+q^2\\
\cmidrule{2-5}
\xi\neq 1&\mathrm{d}&L(\nu\xi,\xi \rtimes \nu^{-\frac{1}{2}}\sigma)&X&X+q^3\\
\midrule
 \mathrm{VI}&\mathrm{a}&\tau(S,\nu^{-\frac{1}{2}}\sigma)&X-2q\sigma(\varpi)&X-q\\
\cmidrule{2-5}
\sigma\ \mathrm{unram.}&\mathrm{c}&L(\nu^{\frac{1}{2}}\St_{\GL(2)},\nu^{-\frac{1}{2}}\sigma)&X-\sigma(\varpi)(q^2+q)&X-q^2\\
\cmidrule{2-5}
\sigma^2=1&\mathrm{d}&L(\nu,1_{F^\times}\rtimes \nu^{-\frac{1}{2}}\sigma)&X^2-\sigma(\varpi)(q+3q^2)X&(X-q^2)(X-q^3)\\
 &&&+2q^3(q+1)&\\
\bottomrule
\end{array}
$
}
\end{table}

\clearpage

\begin{table}
\caption{Characteristic polynomials of $T_{0,1}^s$ and $T_{1,0}^s$
on $V_s(1)$  in terms
of paramodular eigenvalues}
\label{eigenIwahoritable}
\noindent\resizebox{\textwidth}{!}{
$
\renewcommand{\arraystretch}{1.5}
 \begin{array}{lccll}
\toprule
\multicolumn{2}{c}{\mathrm{type}}&\pi&\multicolumn{1}{c}{p(T^s_{0,1},V_s(1),X)}&\multicolumn{1}{c}{p(T^s_{1,0},V_s(1),X)}\\
\bottomrule
 \mathrm{I}&&\chi_1\times\chi_2\rtimes \sigma&X^4-2\lambda_\pi X^3&X^4-(\mu_\pi-q^2+1)X^3\\
\chi_1,\chi_2,\sigma\ \mathrm{unram.} &&&+(\lambda_\pi^2+q\mu_\pi+3q^3+q)X^2&+(q\lambda_\pi^2-2q^2\mu_\pi-2q^2)X^2\\
\chi_1\chi_2\sigma^2=1&&&-q\lambda_\pi(\mu_\pi+3q^2+1)X+q^3\lambda_\pi^2&-q^4(\mu_\pi-q^2+1)X+q^8\\
\midrule
 \mathrm{II}&\mathrm{a}&\chi\St_{\GL(2)}\rtimes\sigma&X^2&X^2-\mu_\pi X+q^3\\
\chi,\sigma\ \mathrm{unram.}&&&+((q-1)\varepsilon_\pi-\lambda_\pi)X&\\
\chi^2\sigma^2=1&&&+q(q^2+q+\mu)&\\ 
\cmidrule{2-5}
&\mathrm{b}&\chi 1_{\GL(2)}\rtimes\sigma&X^2&X^2\\
&&&-\lambda_\pi\dfrac{\mu_\pi+1+q^2+2q^3}{\mu_\pi+1+q+q^2+q^3}X&-q(\mu_\pi(q+1)^{-1}-q+1)X\\
&&&+q^2\mu_\pi(q+1)^{-1}+q^2(q^2+1)&+q^5\\
\midrule
 \mathrm{III}&\mathrm{a}&\chi \rtimes \sigma\St_{\GSp(2)}&X-\lambda_\pi&X-q\\
\cmidrule{2-5}
\chi,\sigma\ \mathrm{unram.}&\mathrm{b}&\chi\rtimes \sigma 1_{\GSp(2)} &X^3&X^3\\
\chi\sigma^2=1&&&-\lambda_\pi(1+q(1+q)^{-1})X^2&-(\mu_\pi+1-q-q^2)X^2\\
&&&+(1+2q^2-q^3+\mu_\pi)q(q+1)X&+q^3(\mu_\pi+1-q^2-q^3)X\\
&&&-\lambda_\pi q^3(q+1)&-q^7\\
\midrule
  \mathrm{IV}&\mathrm{b}&L(\nu^2,\nu^{-1}\sigma \St_{\GSp(2)})&X-\lambda_\pi&X-q\\
\cmidrule{2-5}
\sigma\ \mathrm{unram.}&\mathrm{c}&L(\nu^{\frac{3}{2}}\St_{\GL(2)},\nu^{-\frac{3}{2}}\sigma)&X^2&X^2-(1+q^3)X+q^3\\
\sigma^2=1&&&-\lambda_\pi\dfrac{1+2q+q^3}{2+q+q^3}X&\\
&&&+q(1+q)(1+q^2)&\\
\cmidrule{2-5}
 &\mathrm{d}&\sigma1_{\GSp(4)}&X-\lambda_\pi q^2(q^2+1)^{-1}&X-q^4\\
\midrule
 \mathrm{V}&\mathrm{a}&\delta([\xi,\nu \xi],\nu^{-\frac{1}{2}}\sigma)&X&X+q\\
 \cmidrule{2-5}
\sigma,\xi\ \mathrm{unram.}&\mathrm{b} &L(\nu^{\frac{1}{2}}\xi\St_{\GL(2)},\nu^{-\frac{1}{2}}\sigma)&X-\lambda_\pi q(q+1)^{-1}&X+q^2\\
\cmidrule{2-5}
\sigma^2=\xi^2=1&\mathrm{c} &L(\nu^{\frac{1}{2}}\xi\St_{\GL(2)},\xi\nu^{-\frac{1}{2}}\sigma)&X-\lambda_\pi q(q+1)^{-1}&X+q^2\\
\cmidrule{2-5}
\xi \neq 1&\mathrm{d}&L(\nu\xi,\xi \rtimes \nu^{-\frac{1}{2}}\sigma)&X&X+q^3\\
\midrule
 \mathrm{VI}&\mathrm{a}&\tau(S,\nu^{-\frac{1}{2}}\sigma)&X-\lambda_\pi&X-q\\
\cmidrule{2-5}
\sigma\ \mathrm{unram.}&\mathrm{c}&L(\nu^{\frac{1}{2}}\St_{\GL(2)},\nu^{-\frac{1}{2}}\sigma)&X-\lambda_\pi q (q+1)^{-1}&X-q^2\\
\cmidrule{2-5}
\sigma^2=1&\mathrm{d}&L(\nu,1_{F^\times}\rtimes \nu^{-\frac{1}{2}}\sigma)&X^2-\frac{\lambda_\pi}{ 2}(1+2q(1+q)^{-1})X&(X-q^2)(X-q^3)\\
 &&&+2q^3(q+1)&\\
\bottomrule
\end{array}
$
}
\end{table}

%% file: SKMS_symbols.tex
\Extrachap{Glossary of Notations}
\label{notationchap}

\subsection*{Roman symbols}
\begin{tabbing}
Symbolxxxxxxxxx\=Descriptionxxxxxxxxxxxxxxxxxxxxxxxxxxxxxxxxxxxxxxxxxxxx\kill\\
$\hat{a}(F)$\>Fourier coefficient function of $F$,  \pageref{hataFeq}\\
$a(\chi)$\>conductor of the character $\chi$,  \pageref{conductorpageref}\\
$a(S)$\>Fourier coefficient of Siegel modular form, \pageref{Fourierexpeq}\\
$A'$\>conjugate-inverse-transpose of the $2\times 2$ matrix $A$,  \pageref{Aprimedefeq}\\
$\A$\>adeles of $\Q$,  \pageref{adelesofQ}\\
$\A_f$\>finite adeles of $\Q$,  \pageref{finiteadelesofQ}\\
$A[B]$\>$\transp{B}AB$ for $2\times 2$ matrices $A$ and $B$,  \pageref{transposeBAB}\\
$\mathcal{A}_k$\>adelic automorphic forms of weight $k$,  \pageref{adelicweightk}\\
$\mathcal{A}_k^\circ$\>cuspidal adelic automorphic forms of weight $k$,  \pageref{cuspidaladelicweightk}\\
$\mathcal{A}_k(\mathcal{K})$\>weight $k$ adelic automorphic forms with respect to $\mathcal{K}$,  \pageref{adelicweightkwrtK}\\
$\mathcal{A}_k^\circ(\mathcal{K})$\>cuspidal elements of $\mathcal{A}_k(\mathcal{K})$,   \pageref{cuspidaladelicweightkwrtK}\\
$A(N)$\>paramodular Fourier coefficient indices,  \pageref{ANdefeq}\\
$A(N)^+$\>cuspidal paramodular Fourier coefficient indices,  \pageref{ANdefeq}\\
$A(\Q)$\>positive semi-definite $2\times 2$ matrices with entries from $\Q$,  \pageref{AQdef}\\
$A(\Q)^+$\>positive definite $2\times 2$ matrices with entries from $\Q$,  \pageref{AQdef}\\
$B$\>Borel subgroup of $\GSp(4,F)$,  \pageref{BorelsubgroupGSp4}\\
$B(N)$\>stable Klingen Fourier coefficient indices,  \pageref{BNdefeq}\\
$B(N)^+$\>cuspidal stable Klingen Fourier coefficient indices,  \pageref{BNdefeq}\\
$c_{i,j}$\>value of a local newform evaluated at $\Delta_{i,j}$,  \pageref{cijdefeq}\\
$c(n,r)$\>Fourier coefficient of a Jacobi form,  \pageref{Jacobiformexp}\\
$d^I y$\>Haar measure on $\GSp(4,F)$ with $m(I)=1$,  \pageref{Imeasure}\\
$d_K$\>projection using the parahoric subgroup $K$,  \pageref{dKdef}\\
$e_i=T_{s_i}$\>element of the Iwahori-Hecke algebra defined by $s_i$,  \pageref{eielement}\\
$e=T_1$\>identity element of the Iwahori-Hecke algebra,  \pageref{eielement}\\
$\mathcal{F}$\>vector space of $\C$ valued functions on $A(\Q)$,  \pageref{Fspace}\\
$F\big |_k g$\>slash action on functions on $\mathcal{H}_2$,  \pageref{slashactiondef}\\
$f\big |_{k,m} g$\>slash action on functions on $\mathcal{H}_1 \times \C$,  \pageref{jacobiformsl2defeq}\\
$f_1^{\mathrm{para}},f_2^{\mathrm{para}}$\> basis elements of $(\chi_1 \times \chi_2 \rtimes \sigma)^{\K{}}$,  \pageref{VP02basiseq}\\
$f_w$\>certain element of $(\chi_1 \times \chi_2 \rtimes \sigma)^I$ supported on $BIwI$,  \pageref{fwdef}\\
$g\langle Z \rangle$\> action of $\GSp(4,\R)^+$ on Siegel upper half-space $\mathcal{H}_2$,   \pageref{GH2action}\\
$\transp{g}$\>transpose of the matrix $g$,  \pageref{transposeg}\\
$G^J$\>Jacobi group,  \pageref{Jacobigroupdef}\\
$\GSp(4,F)$\>$4\times 4$ symplectic similitudes with entries from $F$,  \pageref{GSp4Fdef}\\
$\GSp(4,R)$\>$4\times 4$ symplectic similitudes with entries from $R$,  \pageref{GSp4Rdef}\\
$\mathcal{H}_1$\>complex upper half-plane,  \pageref{complexupperhalfplane}\\
$\mathcal{H}_2$\>Siegel upper half-space of degree $2$,  \pageref{siegelupperdef}\\
$\mathcal{H}(\GSp(4,F),I)$ \> Iwahori-Hecke algebra  ,  \pageref{IwahoriHeckedef}\\
$H(\R)$\> Heisenberg group ,  \pageref{Heisenberggroupdef}\\
$I$\>Iwahori subgroup of $\GSp(4,F)$,  \pageref{iwahoridefeq}\\
$I=\begin{bsmallmatrix} i&\\&i \end{bsmallmatrix}$\>element of the Siegel upper half-space $\mathcal{H}_2$,  \pageref{Idefeq}\\
$J$\> matrix defining $\GSp(4)$ (particular to Parts~1, 2),   \pageref{Jdefeq}, \pageref{newJeq}\\
$j(g,Z)$\>degree $2$ factor of automorphy,  \pageref{autofactordef}\\
$K_\infty$\>maximal compact subgroup of $\SSp(4,\R)$,  \pageref{Kinfdef}\\
$\mathcal{K}$\>compact, open subgroup of $\GSp(4,\A_{\mathrm{fin}})$,  \pageref{Kprodefeq}\\
$\K{n}$\>local paramodular subgroup of level $\p^n$,  \pageref{paradefeq}\\
$\Kl{n}$\>local Klingen congruence subgroup of level $\p^n$ ,  \pageref{Klndefeq} \\
$\Ks{n}$\>local stable Klingen congruence subgroup of level $\p^n$,   \pageref{Ksndefeq}\\
$\mathrm{K}_{s,1}(\p^n)$\>auxiliary  congruence subgroup ,  \pageref{Ks1pndef}\\
$\mathrm{K}(N)$\>paramodular congruence subgroup of level $N$,  \pageref{classKNeq}\\
$\mathrm{K}_s(N)$\>stable Klingen congruence subgroup of level $N$,  \pageref{classKsNeq2}\\
$\mathcal{K}(N)$\>adelic paramodular congruence subgroup of level $N$,  \pageref{adelicparadef}\\
$\mathcal{K}_s(N)$\>adelic stable Klingen congruence subgroup of level $N$,  \pageref{adelicparadef}\\
$K_J$\>standard parahoric subgroup corresponding to $J$,  \pageref{KJdef}\\
$L_{c^2}$\>index lowering operator on Jacobi forms,  \pageref{Lcoperatoreq}\\
$L_p'$\>index lowering operator on Jacobi forms,  \pageref{Lpprimedef}\\
$L(s,\pi)$\>$L$-function of $\pi$, \pageref{Lfunctiondefeq}\\
$\mathrm{M}(n)$\>alternative model for local paramodular vectors of level $\p^n$,  \pageref{localMnmodeldef}\\
$\mathrm{M}_{\Delta}(n)$\>$n$-th triangular subspace of $\mathrm{M}_{\infty\times\infty}(\C)$,  \pageref{triangspaceeq}\\
$\mathrm{M}_s(n)$\>alternative model for  stable Klingen vectors of level $\p^n$,  \pageref{Msndef}\\
$\mathrm{M}_0(k)$\>subspace of $\mathrm{M}_{\Delta}(k)$,  \pageref{M0kdef}\\
$M_i$\>certain element of $\GSp(4,\OF)$,  \pageref{Midefeq}\\
$M_k(\Gamma)$\>Siegel modular forms of weight $k$ with respect to $\Gamma$,  \pageref{SMFdef}\\
$m(W)$\>paramodular vector $W$ in the alternative model,  \pageref{mWeq}\\
$N_\pi$\>paramodular level of $\pi$,  \pageref{paraleveldef}\\
$N_{\pi,s}$\>stable Klingen level of $\pi$,  \pageref{stableKlingenleveldef}\\
$\bar N_{\pi,s}$\>quotient stable Klingen level of $\pi$,  \pageref{stableKlingenleveldef}\\
$N_s$\>the integer $N \prod_{p \mid N} p^{-1}$ ,  \pageref{tildeNdefeq}\\
$N_\tau$\>level of the local $\GL(2)$ representation $\tau$,  \pageref{gl2leveldef}\\
$P$\>Siegel parabolic subgroup of $\GSp(4,F)$,  \pageref{Siegelparadef}\\
$P_0$\>certain polynomial,  \pageref{P0def}\\
$P_3$\>subgroup of $\GL(3,F)$,  \pageref{P3defeq}\\
$p_n$\>paramodularization map,  \pageref{pndefeq}\\
$Q$\>Klingen parabolic subgroup of $\GSp(4,F)$,  \pageref{KlingenQdef}\\
$\bar Q$\>subgroup of Klingen parabolic subgroup $Q$ of $\GSp(4,F)$,  \pageref{Qbardef}\\
$Q(\nu^{\frac{1}{2}}\pi,\nu^{-\frac{1}{2}}\sigma)$\> Saito-Kurokawa representation,  \pageref{SKsubsec}\\
$q_w$\>number of left $I$ cosets in a disjoint decomposition of $IwI$,  \pageref{qwdef}\\
$R_{n-1}$\>certain auxiliary operator,  \pageref{Rdefeq}\\
$s_0=t_1$\>Weyl group element for $\GSp(4,F)$,  \pageref{t1defeq}\\
$s_1,s_2$\>Weyl group elements for $\GSp(4,F)$,  \pageref{s1s2defeq}\\
$S_k(\Gamma)$\>Siegel modular cusp forms of weight $k$ with respect to $\Gamma$,  \pageref{SMFdef}\\
$S_k(\mathrm{K}(M))_{\mathrm{old}}$\> subspace of oldforms in $S_k(\mathrm{K}(M))$,  \pageref{oldnewformsdef}\\
$S_k(\mathrm{K}(M))_{\mathrm{new}}$\> subspace of newforms in $S_k(\mathrm{K}(M))$,  \pageref{oldnewformsdef}\\
$s_n$\>trace operator on stable Klingen vectors,  \pageref{snsetupdefeq}\\
$\SSp(4,F)$\>$4 \times 4$ symplectic isometries with entries from $F$,  \pageref{Sp4Fdef}\\
$\SSp(4,R)$\> symplectic isometries with entries from $R$,  \pageref{Sp4Rdef}\\
$S:V\to V$\> auxiliary operator,  \pageref{Soperdef}\\
$\mathcal{S}(X)$\>Schwartz functions on $X$,  \pageref{schwartzfndef}\\
$T_{0,1},T_{1,0}$\> local paramodular Hecke operators,  \pageref{Heckeopdefeq}\\
$T^s_{0,1},T^s_{1,0}$\> local stable Klingen Hecke operators,  \pageref{Ts01eq}\\
$T^s_{0,1}(p)$\>global stable Klingen Hecke operator at $p$,  \pageref{T01formulalemma}\\
$T^s_{1,0}(p)$\>global stable Klingen Hecke operator at $p$,  \pageref{T10formulalemma}\\
$T(1,1,p,p)$\> global paramodular Hecke operator at $p$,  \pageref{classicalHeckeoperatoreq}\\
$T(1,p,p,p^2)$\> global paramodular Hecke operator at $p$,  \pageref{classicalHeckeoperatoreq}\\
$T(l,m)$\>$\GL(2)$ Hecke operator,  \pageref{Tlmdefeq}\\
$T(n)$\>$\GL(2)$ Hecke operator,  \pageref{Tndefeq}\\
$t_n$\>element of $\K{n}$, level raising operator,   \pageref{tndefeq}, \pageref{tnlevelraisinglemmaeq1}\\
$U$\>subgroup of the Borel subgroup of $\GSp(4,F)$,  \pageref{Udefeq}\\
$U_p$\>index raising operator on Jacobi forms,  \pageref{Updef}\\
$U_3$\>subgroup of $P_3$,  \pageref{U3def}\\
$u_1$\>extended affine Weyl group element,  \pageref{u1defeq}\\
$u_1=T_{u_1}$\>element of the Iwahori-Hecke algebra,  \pageref{Tu1def}\\
$u_n$\>Atkin-Lehner element of $\GSp(4,F)$,  \pageref{ALeq}\\
$V^\vee$\>contragredient of $V$,  \pageref{contradef}\\
$V_{N,\chi}$\>Jacquet module,  \pageref{Jacquetmoddef}\\
$v_p(N)$\>exponent of the power of $p$ that divides $N$ exactly,  \pageref{vpNdef}\\
$V(n)$\>vectors in $V$ fixed by $\K{n}$,  \pageref{nparadef}\\
$V_s(n)$\>vectors in $V$ fixed by $\Ks{n}$,  \pageref{Vsndefeq}\\
$\bar V_s(n)$\>quotient of $V_s(n)$ by $V(n-1)+V(n)$,  \pageref{Vsbardefeq}\\
$V(\p^n)$\>vectors in $V$ fixed by $\K{n}$ (global setting notation),  \pageref{pexpKeq}\\
$V_s(\p^n)$\>vectors in $V$ fixed by $\Ks{n}$ (global setting notation) ,  \pageref{pexpKseq}\\
$v_s$\>shadow of a given newform $v_{\mathrm{new}}$,  \pageref{shadowdefeq2}\\
$V_{Z^J}$\>$P_3$-quotient of $V$,  \pageref{VZJdef}\\
$W$\>Weyl group,  \pageref{Weyldef}\\
$W^a$\>affine Weyl group,  \pageref{Wadef}\\
$W^a_J$\>subgroup of $W^a$ generated by $J$,  \pageref{WaJdef}\\
$W^e$\>extended affine Weyl group (Iwahori-Weyl group),  \pageref{Wedef}\\
$w_p$\>global paramodular Atkin-Lehner operator at $p$,  \pageref{ATglobaldefeq}\\
$\mathcal{W}(\pi,\psi_{c_1,c_2})$\>Whittaker model of $\pi$, \pageref{whitmodeldef}\\
$Z$\>center of $\GSp(4,F)$, \pageref{Zcenterdef}\\
$Z^J$\>center of the Jacobi group, \pageref{ZJdefeq}\\
$Z(s,W)$\>zeta integral, \pageref{zetaidef}\\
$Z_N(s,W)$\>simplified zeta integral, \pageref{simpzetadef}
\end{tabbing}

\subsection*{Greek and other symbols}

\begin{tabbing}
Symbolxxxxxxxxx\=Descriptionxxxxxxxxxxxxxxxxxxxxxxxxxxxxxxxxxxxxxxxxxxxx\kill\\
$\Gamma_0(\p^n)$\>Hecke congruence subgroup of $\GL(2,F)$, \pageref{K1defeq}\\
$\Gamma_0(N)$\>Hecke congruence subgroup of $\SL(2,\Z)$ of level $N$, \pageref{gamma0Ndef}\\
$\Gamma_0'(N)$\>Klingen congruence subgroup of level $N$, \pageref{classKlNeq2}\\
$\gamma(s,\pi)$\>gamma factor of $\pi$, \pageref{gammafactordef}\\
$\Delta_t^+,\Delta_t^-$\>operators on $\mathcal{F}$, \pageref{deltatypedef}\\
$\Delta_{i,j}$\>certain diagonal element of $\GSp(4,F)$, \pageref{deltaijdefeq}\\
$\delta_P$\>modular character, \pageref{modchardef}\\
$\varepsilon_\pi$\>eigenvalue of $\pi(u_{N_\pi})$ on newform in $\pi$, \pageref{epsilonlambdamudefeq}\\
$\varepsilon(s,\pi)$\>epsilon factor of $\pi$, \pageref{epsilonfactordef}\\
$\eta$\>local level raising operator, \pageref{levelraisdef}, \pageref{etadefeq}, \pageref{altlevelraisdef}\\
$\eta_p$\>global level raising operator at $p$ induced by $\eta$, \pageref{etathetathetapglobaldef}, \pageref{etaformulalemma}\\
$\theta$\>local level raising operator, \pageref{levelraisdef}, \pageref{thetadefeq1}, \pageref{altlevelraisdef}\\
$\theta_p$\> global level raising operator at $p$ induced by $\theta$, \pageref{etathetathetapglobaldef}, \pageref{thetaformulalemma}\\
$\theta'$\>local level raising operator, \pageref{levelraisdef}\\
$\theta'_p$\>global level raising operator at $p$ induced by $\theta'$, \pageref{etathetathetapglobaldef}\\
$[\lambda,\mu,\kappa]$\>element of the Heisenberg group, \pageref{heisenelementdef}\\
$\lambda(g)$\> similitude factor, \pageref{simdef}, \pageref{sim2def}\\
$\lambda_\pi$\>eigenvalue of $T_{0,1}$ on newform in $\pi$, \pageref{epsilonlambdamudefeq}\\
$\mu_\pi$\>eigenvalue of $T_{1,0}$ on newform in $\pi$, \pageref{epsilonlambdamudefeq}\\
$\nu$\> absolute value character, \pageref{absolutevalchardef}\\
$\pi^\vee$\> contragredient of $\pi$, \pageref{contradef}\\
$\sigma_n,\sigma$\> local level lowering operator, \pageref{sigmandefformeq}, \pageref{Ksdiagramsproposition}\\
$\sigma_p$\> global level lowering operator at $p$ induced by $\sigma$, \pageref{sigmaformulalemma}\\
$\tau_n,\tau$\> local level raising operator, \pageref{taunlevelraisinglemmaeq1}, \pageref{taurefdef}, \pageref{Ksdiagramsproposition}\\
$\tau_p$\> global level raising operator at $p$ induced by $\tau$, \pageref{tauformulalemma}\\
$\psi_{c_1,c_2}$\> character of $U$, \pageref{whitmodeldef}\\
$\Omega$\>subgroup of $W^e$, \pageref{omegadef}\\
$\omega_\pi$\> central character of $\pi$, \pageref{centralchardef}\\
$\langle \cdot,\cdot \rangle$\> Petersson inner product, \pageref{peterssoninnerdef}\\
$\nabla_t$\>operator on $\mathcal{F}$, \pageref{deltatypedef}
\end{tabbing}
 
 